\documentclass[openany]{amsbook}
\usepackage{amssymb, amsfonts}
\usepackage{enumerate}
\usepackage[all,arc]{xy}

\makeatletter

\def\makenoteindex{%
  \newwrite\@noteindexfile
  \immediate\openout\@noteindexfile=notation.idx
  \def\noteindex{\@bsphack\begingroup
             \@sanitize
             \@wrnoteindex}\typeout
    {Writing index file notation.idx}%
  \let\makenoteindex\@empty
}
\@onlypreamble\makeindex 

\def\@wrnoteindex#1{%
   \protected@write\@noteindexfile{}%
      {\string\indexentry{#1}{\thepage}}%
 \endgroup
 \@esphack}

\def\noteindex{\@bsphack\begingroup \@sanitize\@noteindex}
\def\@noteindex#1{\endgroup\@esphack}

\makeatother

\newtheorem{thm}{Theorem}[section]
\newtheorem{cor}[thm]{Corollary}
\newtheorem{prop}[thm]{Proposition}
\newtheorem{lem}[thm]{Lemma}

\theoremstyle{definition}
\newtheorem{hyp}[thm]{Hypothesis}
\newtheorem{defn}[thm]{Definition}

\newtheorem{con}[thm]{Construction}
\newtheorem{ouch0}[thm]{Counterexample}
\newtheorem{exmp}[thm]{Example}

\newtheorem{notn}[thm]{Notation}
\newtheorem{notns}[thm]{Notations}

\theoremstyle{remark}
\newtheorem{rem}[thm]{Remark}

\newtheorem{warn}[thm]{Warning}
\newtheorem{sch}[thm]{Scholium}

\DeclareFontFamily{OMS}{rsfs}{\skewchar\font'60}
\DeclareFontShape{OMS}{rsfs}{m}{n}{<-5>rsfs5 <5-7>rsfs7 <7->rsfs10 }{}
\DeclareSymbolFont{rsfs}{OMS}{rsfs}{m}{n}
\DeclareSymbolFontAlphabet{\scr}{rsfs}

\newcommand{\sA}{\scr{A}}
\newcommand{\sB}{\scr{B}}
\newcommand{\sC}{\scr{C}}
\newcommand{\sD}{\scr{D}}
\newcommand{\sE}{\scr{E}}
\newcommand{\sF}{\scr{F}}
\newcommand{\sG}{\scr{G}}
\newcommand{\sH}{\scr{H}}
\newcommand{\sI}{\scr{I}}
\newcommand{\sJ}{\scr{J}}
\newcommand{\sK}{\scr{K}}

\newcommand{\sM}{\scr{M}}

\newcommand{\sO}{\scr{O}}
\newcommand{\sP}{\scr{P}}
\newcommand{\sQ}{\scr{Q}}

\newcommand{\sS}{\scr{S}}
\newcommand{\sT}{\scr{T}}
\newcommand{\sU}{\scr{U}}
\newcommand{\sV}{\scr{V}}
\newcommand{\sW}{\scr{W}}

\newcommand{\bF}{\mathbb{F}}

\newcommand{\bH}{\mathbb{H}}

\newcommand{\bN}{\mathbb{N}}

\newcommand{\bP}{\mathbb{P}}
\newcommand{\bQ}{\mathbb{Q}}
\newcommand{\bR}{\mathbb{R}}

\newcommand{\bU}{\mathbb{U}}

\newcommand{\bZ}{\mathbb{Z}}

\newcommand{\al}{\alpha}
\newcommand{\be}{\beta}
\newcommand{\ga}{\gamma}
\newcommand{\de}{\delta}
\newcommand{\epz}{\varepsilon}
\newcommand{\ph}{\phi}

\newcommand{\et}{\eta}
\newcommand{\io}{\iota}
\newcommand{\ka}{\kappa}
\newcommand{\la}{\lambda}
\newcommand{\tha}{\theta}

\newcommand{\rh}{\rho}
\newcommand{\si}{\sigma}
\newcommand{\ta}{\tau}
\newcommand{\ch}{\chi}
\newcommand{\ps}{\psi}
\newcommand{\ze}{\zeta}
\newcommand{\om}{\omega}
\newcommand{\GA}{\Gamma}
\newcommand{\LA}{\Lambda}
\newcommand{\DE}{\Delta}
\newcommand{\SI}{\Sigma}

\newcommand{\OM}{\Omega}

\newcommand{\PI}{\Pi}

\newcommand{\PH}{\Phi}

\newcommand{\com}{\circ}     
\newcommand{\iso}{\cong}     
\newcommand{\htp}{\simeq}    
\newcommand{\sma}{\wedge}    
\newcommand{\wed}{\vee}      
\newcommand{\rtarr}{\longrightarrow}
\newcommand{\ltarr}{\longleftarrow}

\newcommand{\U}[1]{\bar{P}{#1}}
\newcommand{\T}[1]{P\times_{\PI}{#1}}
\newcommand{\rB}{\bar{r}}
\newcommand{\rE}{r}

\newcommand{\Ho}{\text{Ho}\,}

\makeatletter
\let\c@equation\c@thm
\makeatother
\numberwithin{equation}{section}

\makeatletter
\ifx\SK@label\undefined\let\SK@label\label\fi
 \let\your@thm\@thm
 \def\@thm#1#2#3{\gdef\currthmtype{#3}\your@thm{#1}{#2}{#3}}
 \def\mylabel#1{{\let\your@currentlabel\@currentlabel\def\@currentlabel
  {\currthmtype~\your@currentlabel}
 \SK@label{#1@}}\label{#1}}
 \def\myref#1{\ref{#1@}}
\makeatother
  
\begin{document}

\frontmatter

\title{Parametrized Homotopy Theory}
\author{J.\,P. May}
\author{J. Sigurdsson}
\address{Department of Mathematics, The University of Chicago, Chicago, IL 60637}
\address{Department of Mathematics, The University of Notre Dame, Notre Dame,
IN, 46556-4618}
\email{may@math.uchicago.edu}
\email{jsigurds@nd.edu}

\thanks{May was partially supported by the NSF}




\begin{abstract}
Part I: We set the stage for our homotopical work with
preliminary chapters on the point-set topology necessary 
to parametrized homotopy theory, the base change and other
functors that appear in over and under categories, and
generalizations of several classical results about 
equivariant bundles and fibrations to the context of
proper actions of non-compact Lie groups.


Part II:
Despite its long history, the homotopy theory of ex-spaces requires
further development before it can serve as the starting point
for a rigorous modern treatment of parametrized stable homotopy theory.
We give a leisurely account that emphasizes several issues that
are of independent interest in the theory and applications of
topological model categories. The essential point is to resolve
problems about the homotopy theory of ex-spaces that are 
absent from the homotopy theory of spaces. In contrast to 
previously encountered situations, model theoretic techniques 
are intrinsically insufficient to a full development of the basic 
foundational properties of the homotopy category of ex-spaces.  
Instead, a rather intricate blend of model theory and classical 
homotopy theory is required.  However, considerable new material 
on the general theory of topologically enriched model categories 
is also required. 

Part III:
We give a systematic treatment of the foundations of parametrized stable 
homotopy theory, working equivariantly and with highly structured smash 
products and function spectra. The treatment is based on equivariant 
orthogonal spectra, which are simpler for the purpose than 
alternative kinds of spectra.  Again, the parametrized context introduces many
difficulties that have no nonparametrized counterparts and cannot be dealt 
with using standard model theoretic techniques. The space level techniques
of Part II only partially extend to the spectrum level, and many new twists
are encountered. Most of the difficulties are already present in the
nonequivariant special case. Equivariantly, we show how change of universe,
passage to fixed points, and passage to orbits behave in the parametrized
setting.

Part IV: 
We give a fiberwise duality 
theorem that allows fiberwise recognition of dualizable and 
invertible parametrized 
spectra.  This allows direct application of the formal theory of duality in 
symmetric monoidal categories to the construction and analysis of transfer maps.
The relationship between transfer for general Hurewicz fibrations and for fiber
bundles is illuminated by the construction of fiberwise bundles of spectra, which
are like bundles of tangents along fibers, but with spectra replacing spaces as fibers. Using this construction, we obtain a simple conceptual proof of a generalized Wirthm\"uller isomorphism theorem that calculates the right adjoint to base change along an equivariant bundle with manifold fibers in terms of a shift of the left adjoint. Due to the generality of our bundle theoretic context, the Adams isomorphism theorem relating orbit and fixed point spectra is a direct consequence.
\end{abstract}


\maketitle

\tableofcontents

\mainmatter

\chapter*{Prologue}

What is this book about and why is it so long?  Parametrized homotopy theory concerns systems of spaces and spectra that are parametrized as fibers over points of a given base space $B$. Parametrized spaces, or ``ex-spaces'', are just spaces over and under $B$, with a projection, often a fibration, and a section. Parametrized spectra are analogous but considerably more sophisticated objects.  They provide a world in which one can apply the methods of stable homotopy theory without losing track of fundamental groups and other unstable information.  Parametrized homotopy theory is a natural and important part of homotopy theory that is implicitly central to all of bundle and fibration theory. Results that make essential use of it are widely scattered throughout the literature.  For classical examples, the theory of transfer maps is
intrinsically about parametrized homotopy theory, and Eilenberg-Moore type spectral sequences are parametrized K\"unneth theorems.  Several new and 
current directions, such as ``twisted'' cohomology theories 
and parametrized fixed point theory cry out for the rigorous foundations that
we shall develop.

On the foundational level, homotopy theory, and especially stable homotopy theory, has undergone a thorough reanalysis in recent years. Systematic use of Quillen's theory of model categories has illuminated the structure of the subject and has done so in a way that makes the general methodology widely applicable to other branches of mathematics. The discovery of categories of spectra with associative and commutative smash products has revolutionized stable homotopy theory. The systematic study and application of equivariant algebraic topology has greatly enriched the subject.

There has not been a thorough and rigorous study of parametrized homotopy theory that takes these developments into account.  It is the purpose of this monograph to provide such a study, although we shall leave many interesting loose ends.  We shall also give some direct applications, especially to equivariant stable homotopy theory where the new theory is particularly essential. The reason this study is so lengthy is that, rather unexpectedly, many seemingly trivial nonparametrized results fail to generalize, and many of the conceptual and technical obstacles to a rigorous treatment have no nonparametrized counterparts. For this reason, the resulting theory is considerably more subtle than its nonparametrized precursors. We indicate some of these problems here.

The central conceptual subtlety in the theory enters when we try to prove that structure enjoyed by the point-set level categories of parametrized spaces descends to their homotopy categories. Many of our basic functors occur in Quillen adjoint pairs, and such structure
descends directly to homotopy categories. Recall that an adjoint pair of functors $(T,U)$ between model categories is a Quillen adjoint pair, or a Quillen adjunction, if the left adjoint $T$ preserves cofibrations and acyclic cofibrations or, equivalently, the right adjoint $U$ preserves fibrations and acyclic fibrations. It is a Quillen equivalence if, further, the induced adjunction on homotopy categories is an adjoint equivalence. We originally
hoped to find a model structure on para\-me\-trized spaces in which all of the relevant adjunctions are Quillen  adjunctions.  It eventually became clear that there can be no such model structure, for altogether trivial reasons. Therefore, 
it is intrinsically impossible to lay down the basic foundations of parametrized homotopy theory using only the standard methodology of model category theory.  

The force of parametrized theory largely comes from base change functors associated to maps $f\colon A\longrightarrow B$. The existing literature on fiberwise homotopy theory says surprisingly little about such functors. 
This is especially strange since they are the most important feature that makes parametrized homotopy theory useful for the study of ordinary homotopy theory:
such functors are used to transport information from the parametrized context
to the nonparametrized context.  One of the goals of our work is to 
fill this gap. 

On the point-set level,
there is a pullback functor $f^*$ from ex-spaces (or spectra) over $B$ to ex-spaces (or spectra) over $A$. That functor has a left adjoint $f_!$ and a right adjoint $f_*$. We would like both  of these to be Quillen adjunctions, but that is not possible unless the model structures lead to trivial homotopy categories. We mean literally trivial: one object and one morphism. We explain why. It will be clear that the explanation is generic and applies equally well to a variety of sheaf theoretic situations where one encounters analogous base change functors. 
\begin{ouch0}\mylabel{noway}  Consider the following diagram.
$$\xymatrix{
\emptyset \ar[r]^-{\ph} \ar[d]_{\ph}  &  B \ar[d]^{i_{0}} \\
B \ar[r]_-{i_{1}} & B\times I}$$
Here $\emptyset$ is the empty set and $\phi$ is the initial (empty) map into $B$. This diagram is a pullback since $B\times\{0\}\cap B\times\{1\} = \emptyset$. The category of ex-spaces over $\emptyset$ is the trivial category with one object, and it admits a unique model structure. Let $*_B$ denote the ex-space $B$ over $B$, with section and projection the identity map. Both $(\ph_!,\ph^*)$ and $(\ph^*,\ph_*)$ are Quillen adjoint pairs for any model structure on the category of ex-spaces over $B$. Indeed, $\ph_!$ and $\ph_*$ preserve weak equivalences, fibrations, and cofibrations since both take $*_{\emptyset}$ to $*_B$. We have $(i_0)^*\com (i_1)_! \iso \ph_!\com\ph^*$ since both composites take any ex-space over $B$ to $*_B$.  If $(i_1)_!$ and $(i_0)^*$ were both Quillen left adjoints, it would follow that this isomorphism descends to homotopy categories. If, further, the functors $(i_1)_!$ and $(i_0)^*$ on homotopy categories were equivalences of categories, this would imply that the homotopy category of ex-spaces over $B$ with respect to the given model structure is equivalent to the trivial category.
\end{ouch0}

Information in ordinary homotopy theory is derived from results in parametrized homotopy theory by use of the base change functor $r_!$ associated to the trivial map $r\colon B\rtarr *$. For this and other reasons, we choose our basic model structure to be one such that $(f_!,f^*)$ is a Quillen adjoint pair for every map $f\colon A\rtarr B$ and is a Quillen equivalence when $f$ is a homotopy equivalence. Then $(f^*,f_*)$ cannot be a Quillen adjoint pair in general.  However, it is essential that we still have the adjunction $(f^*,f_*)$ after passage to homotopy categories. For example, taking $f$ to be the diagonal map on $B$, this adjunction is used to obtain the adjunction on homotopy categories that relates the fiberwise smash product functor $\sma_B$ on ex-spaces over $B$ to the function ex-space functor $F_B$.  To construct the homotopy category level right adjoints $f_*$, we shall have to revert to more classical methods, using Brown's representability theorem.  However, it is not clear how to verify the hypotheses of Brown's theorem in the model theoretic framework. 

\myref{noway} also illustrates the familiar fact that a commutative diagram of functors on the point-set level need not induce a commutative diagram of functors on homotopy categories.  When commuting left and right adjoints,
this is a problem even when all functors in sight are parts of Quillen
adjunctions.  Therefore, proving that compatibility relations that hold on the point-set level descend to the homotopy category level is far from automatic. 
In fact, proving such ``compatibility relations'' is often a highly 
non-trivial problem, but one which is essential to the applications. 
We do not know how to prove the most interesting compatibility relations 
working only model theoretically.

Even in the part of the theory in which model theory works, it does not work as expected. There is an obvious naive model structure on ex-spaces over $B$ in which the weak equivalences, fibrations, and cofibrations are the ex-maps whose maps of total spaces are weak equivalences, fibrations, and cofibrations of spaces in the usual Quillen model structure.  This ``$q$-model structure'' is the natural starting point for the theory, but it turns out to have severe drawbacks that limit its space level utility and bar it from serving as the starting point for the development of a useful spectrum level stable model structure. In fact, it has two opposite drawbacks. First, it has too many cofibrations.  In particular, the model theoretic cofibrations need not be cofibrations in the intrinsic homotopical sense.  That is, they fail to satisfy the fiberwise homotopy extension property (HEP) defined in terms of parametrized mapping cylinders.  This already fails for the sections of cofibrant objects and for the inclusions of cofibrant objects in their cones.  Therefore the classical theory of cofiber sequences fails to mesh with the model category structure. 

Second, it also has too many fibrations.  The fibrant ex-spaces are Serre fibrations, and Serre fibrations are not preserved by fiberwise colimits. Such colimits are preserved by a more restrictive class of fibrations, namely the well-sectioned Hurewicz fibrations, which we call ex-fibrations. Such preservation properties are crucial to resolving the problems with base change functors that we have indicated.

In model category theory, decreasing the number of cofibrations increases the number of fibrations, so that these two problems cannot admit a solution in common. Rather, we require two different equivalent descriptions of our homotopy categories of ex-spaces. First, we have another model structure, the ``$qf$-model structure'', which has the same weak equivalences as the $q$-model structure but has fewer cofibrations, all of which satisfy the fiberwise HEP. Second, we have a description in terms of the classical theory of ex-fibrations, which does not fit naturally into a model theoretic framework. The former is vital to the development of the stable model structure on parametrized spectra. The latter is vital to the solution of the intrinsic problems with base change functors.  

Before getting to the issues just discussed, we shall have to resolve various others that also have no nonparametrized analogues. Even the point set topology requires care since function ex-spaces take us out of the category of compactly generated spaces. Equivariance raises further problems, although most of our new foundational work is already necessary nonequivariantly. Passage to the spectrum level raises more serious problems. The main source of difficulty is that the underlying total space functor is too poorly behaved, especially with respect to smash products and fibrations, to give good control of homotopy groups as one passes from parametrized spaces to parametrized spectra. Moreover, the resolution of base change problems requires a different set of details on the spectrum level than on the space level.

The end result may seem intricate, but it gives a very powerful framework in which to study homotopy theory. We illustrate by showing how fiberwise duality and transfer maps work out and by showing that the basic change of groups isomorphisms of equivariant stable homotopy theory, namely the generalized Wirthm\"uller and Adams isomorphisms, drop out directly from the foundations. Costenoble and Waner \cite{CWNew} have already given other applications in equivariant stable homotopy theory, using our foundations to study Poincar\'e duality in ordinary $RO(G)$-graded cohomology.  Further applications are work 
in progress.

The theory here gives perhaps the first worked example in which a model theoretic approach to derived homotopy categories is intrinsically insufficient and must be blended with a quite different approach even to establish the essential structural features of the derived category.  Such a blending of techniques seems essential in analogous sheaf theoretic contexts that have 
not yet received a modern model theoretic treatment.
Even nonequivariantly, the basic results on base change, smash products, and function ex-spaces that we obtain do not appear in the literature. Such results are essential to serious work in parametrized homotopy theory. 

Much of our work should have applications beyond the new parametrized theory.
The model theory of topological enriched categories has received much less attention in the literature than the model theory of simplicially enriched
categories.  Despite the seemingly equivalent nature of these variants,
the topological situation is actually quite different from the simplicial one,
as our applications make clear.  In particular, the interweaving of $h$-type 
and $q$-type model structures that pervades our work seems to have no simplicial 
counterpart. This interweaving does also appear in algebraic contexts 
of model categories enriched over chain complexes, where foundations analogous
to ours can be developed.  One of our goals is to give a thorough analysis and axiomatization of how this interweaving works in general in topologically enriched model categories. 

\vspace{.7mm}

{\em History.} This project began with unpublished notes, dating from the summer of 2000, of the first author \cite{May}. He put the project aside and returned to it in the fall of 2002, when he was joined by the second author. Some of Parts I and II was originally in a draft of the first author that was submitted and accepted for publication, but was later withdrawn.  That draft was correct, but it did not include the ``$qf$-model structure'', which comes from the second author's 2004 PhD thesis \cite{Sig}. The first author's notes \cite{May} claimed to construct the stable model structure on parametrized spectra starting from the $q$-model structure on ex-spaces.  Following \cite{May}, the monograph \cite{Hu} of Po Hu also takes that starting point and makes that claim.   The second author realized that, with the obvious definitions, the axioms for the stable model structure cannot be proven from that starting point and that any
naive variant would be disconnected with cofiber sequences and other essential
needs of a fully worked out theory.  His $qf$-model structure is the crucial new ingredient that is used to solve this problem. Although implemented quite differently, the applications of Chapter 16 were inspired by Hu's work.

\pagebreak

{\em Thanks.}  We thank the referee of the partial first version for several helpful suggestions,  Gaunce Lewis and Peter Booth for help with the point set topology, Mike Cole for sharing his remarkable insights about model
categories, and Mike Mandell for much technical help. We thank Kathleen
Lewis for working out the counterexample in Theorem 1.1.1 and Victor Ginzburg for giving us the striking Counterexample 11.6.2. 
We are especially grateful to Kate Ponto for a meticulously
careful reading that uncovered many obscurities and infelicities.
Needless to say, she is not to blame for those that remain.


\part{Point-set topology, change functors, and proper actions}

\chapter{The point-set topology of parametrized spaces}

\section*{Introduction}

We develop the basic point-set level properties of the category of ex-spaces over a fixed base space $B$ in this chapter. In \S1.1, we discuss convenient categories of topological spaces. The usual category of compactly generated spaces is not adequate for our study of ex-spaces, and we shall see later that the interplay between model structures and the relevant convenient categories is quite subtle. In \S1.2, we give basic facts about based and unbased topologically bicomplete categories. This gives the language that is needed to describe the good formal properties of the various categories in 
which we shall work. We discuss convenient categories of ex-spaces in \S1.3, 
and we discuss convenient categories of ex-$G$-spaces in \S1.4.

As a matter of recovery of lost 
folklore, \S\ref{sec:topass} is an appendix, the substance of 
which is due to Kathleen Lewis. It is only at her insistence that she is not 
named as its author. It documents the nonassociativity of the smash product 
in the ordinary category of based spaces, as opposed to the category 
of based $k$-spaces. When writing the historical paper \cite{History}, 
the first author came across several 1950's references to this 
nonassociativity, 
including an explicit, but unproven, counterexample in a 1958 paper of 
Puppe \cite{Puppe}. However, we know of no reference that gives details, and we
feel that this nonassociativity should be documented in the modern literature.

We are very grateful to Gaunce Lewis for an extended correspondence 
and many details about the material of this chapter, but he is not to be blamed for the point of view that we have taken.  We are also much indebted to Peter Booth. He is the main pioneer of the theory of fibered mapping spaces (see \cite{B1, B2, B3}) and function ex-spaces, and he sent us several detailed proofs about them.  

\section{Convenient categories of topological spaces}\label{sec:pt}

We recall the following by now standard definitions.

\begin{defn} 
Let $B$ be a space and $A$ a subset. Let $f\colon K\rtarr B$ run 
over all continuous maps from compact Hausdorff spaces $K$ into $B$.
\begin{enumerate}[(i)]
\item $A$ is \emph{compactly closed}\index{compactly closed} if each $f^{-1}(A)$ is closed.
\item $B$ is \emph{weak Hausdorff}\index{weak Hausdorff} if each $f(K)$ is closed.
\item $B$ is a \emph{$k$-space}\index{kspace@$k$-space} if each compactly closed subset is closed.
\item $B$ is \emph{compactly generated}\index{compactly generated!space} if it is a weak Hausdorff $k$-space.
\end{enumerate}
Let ${\sT}\!op$\noteindex{Top@${\sT}op$}\index{category!of spaces} be the category of all topological spaces and let $\sK$,\noteindex{K@$\sK$} $w\sH$, and $\sU = \sK\cap w\sH$\noteindex{U@$\sU$}
be its full subcategories of $k$-spaces, weak Hausdorff spaces, and compactly generated 
spaces. The \emph{$k$-ification}\index{kification@$k$-ification functor} functor $k\colon {\sT}\!op\rtarr \sK$ assigns to a space $X$ the
same set with the finer topology that is obtained by requiring all compactly closed 
subsets to be closed. It is right adjoint to the inclusion $\sK\rtarr {\sT}op$. The 
\emph{weak Hausdorffication}\index{weak Hausdorffication functor} functor $w\colon {\sT}op\rtarr w\sH$ assigns to a space $X$ its 
maximal weak Hausdorff quotient. It is left adjoint to the inclusion $w\sH\rtarr {\sT}\!op$.
\end{defn}

>From now on, we work in $\sK$, implicitly $k$-ifying any space that
is not a $k$-space to begin with. In particular, products and function spaces are understood to be $k$-ified.  With this convention, $B$ is weak Hausdorff 
if and only if the diagonal map embeds it as a closed subspace of $B\times B$.
Let $A\times_c B$ denote the classical cartesian product in ${\sT}op$ and 
recall that $B$ is Hausdorff if and only if the diagonal embeds it as a 
closed subspace of $B\times_c B$. The following result is proven in 
\cite[App.\S2]{Lewis0}. 

\begin{prop}\mylabel{HauswHaus} 
Let $A$ and $B$ be $k$-spaces. If one of them is locally
compact or if both of them are first countable, then
$$A\times B = A\times_c B.$$
Therefore, if $B$ is either locally compact or first countable, then
$B$ is Hausdorff if and only if it is weak Hausdorff.
\end{prop}

We would have preferred to work in $\sU$ rather than $\sK$, since there are
many counterexamples which reveal the pitfalls of working without a separation
property.  However, as we will explain in \S1.3, several inescapable facts 
about ex-spaces force us out of that convenient category. 
Like $\sU$, the category $\sK$ is closed cartesian monoidal. This means
that it has function spaces $\text{Map}(X,Y)$\noteindex{MapXY@$\text{Map}(X,Y)$} with homeomorphisms
$$ \text{Map}(X\times Y,Z)\iso \text{Map}(X,\text{Map}(Y,Z)).$$
This was proven by Vogt \cite{Vogt}, who uses the term compactly 
generated for our $k$-spaces. See also \cite{Wylie}. An early unpublished preprint by Clark \cite{Clark} also showed 
this, and an exposition of ex-spaces based on \cite{Clark} was given by Booth \cite{B2}.

Philosophically, we can justify a preference for $\sK$ over $\sU$ by remarking
that the functor $w$ is so poorly behaved that we prefer to minimize its use.
In $\sU$, colimits must be constructed by first constructing them in $\sK$ and then applying the functor $w$, which changes the underlying point set and loses 
homotopical control.  However, this justification would be more persuasive were it not that colimits in $\sK$ that are not colimits in $\sU$ can already be quite badly behaved topologically. For example, $w$ itself is a colimit construction in $\sK$.  We describe a relevant situation in which colimits behave better in $\sU$ than in $\sK$ in \myref{Umod} below.

More persuasively, $w$ is a formal construction that only retains formal 
control because both colimits and the functor $w$ are left adjoints.  We will encounter right adjoints constructed in $\sK$ that do not preserve the weak Hausdorff property when restricted to $\sU$, and in such situations we cannot apply $w$ without losing the adjunction.  In fact, when restricted to $\sU$, 
the relevant left adjoints do not commute with colimits and so cannot be left adjoints there.  We shall encounter other reasons for working in $\sK$ later.  An obvious advantage of $\sK$ is that $\sU$ sits inside it, so that we can use $\sK$ when it is needed, but can restrict to the better behaved category $\sU$ whenever possible. Actually, the situation is more subtle than a simple dichotomy. In our study of ex-spaces, it is essential to combine use of the 
two categories, requiring base spaces to be in $\sU$ but allowing 
total spaces to be in $\sK$.

We have concomitant categories $\sK_*$\noteindex{Kpt@$\sK_*$}\index{category!of based spaces} and $\sU_*$\noteindex{Upt@$\sU_*$} of based spaces in $\sK$ and in $\sU$. We generally write $\sT$\noteindex{T@$\sT$} for $\sU_*$ to mesh with 
a number of relevant earlier papers.  Using duplicative notations, we write $\text{Map}(X,Y)$ for the space $\sK(X,Y)$ of maps $X\rtarr Y$ 
and $F(X,Y)$\noteindex{FXY@$F(X,Y)$} for the based space $\sK_*(X,Y)$ of based maps $X\rtarr Y$ between based spaces. Both $\sK_*$ and $\sT$ are closed symmetric monoidal categories under $\sma$ and $F$ \cite{Lewis0, Vogt, Wylie}. This means that the smash
product is associative, commutative, and unital up to coherent natural isomorphism and that $\sma$ and $F$ are related by the usual adjunction homeomorphism
$$ F(X\sma Y, Z) \iso F(X,F(Y,Z)).$$
The need for $k$-ification is illustrated by the nonassociativity of the smash product in ${\sT}\!op_*$; see \S1.5.

We need a few observations about inclusions and colimits.
Recall that a map is an inclusion if it is a homeomorphism onto its image.
Of course, inclusions need not have closed image.  As noted by Str{\o}m \cite{Strom1}, 
the simplest example of a non-closed inclusion in $\sK$ is the inclusion $i\colon \{a\}\subset \{a,b\}$, where $\{a,b\}$ has the indiscrete topology. 
Here $i$ is both the inclusion of a retract and a Hurewicz cofibration (satisfies the homotopy extension property, or HEP). As is well-known, 
such pathology cannot occur in $\sU$.

\begin{lem}\mylabel{coflemma} Let $i\colon A\rtarr X$ be a map in $\sK$.
\begin{enumerate}[(i)]
\item If there is a map $r\colon X\rtarr A$ such that $r\com i = \text{id}$,
then $i$ is an inclusion.  If, further, $X$ is in $\sU$, then $i$ is a closed inclusion.
\item If $i$ is a Hurewicz cofibration, then $i$ is an inclusion. If, further, $X$ is in $\sU$, then $i$ is a closed inclusion.
\end{enumerate}
\end{lem}

\begin{proof}
Inclusions $i\colon A\rtarr X$ are characterized by the property that a function $j\colon Y\rtarr A$ is continuous if and only if $i\com j$ is continuous. This implies the first statement in (i).  
Alternatively, one can note that a map in $\sK$ is an inclusion if and only if it is an equalizer in $\sK$, and a map in $\sU$ is a closed inclusion if and only if it is an equalizer in $\sU$ \cite[7.6]{Lewis0}. Since $i$ is the equalizer of $i\com r$ and the identity map of $X$, this implies both statements in (i). 
 For (ii), let $Mi$ be the mapping cylinder of $i$. The canonical map $j\colon Mi\rtarr X\times I$ has a left inverse $r$ and is thus an inclusion or closed inclusion in the respective cases. The evident closed inclusions $i_1\colon A\rtarr Mi$ and $i_1\colon X\rtarr X\times I$ satisfy $j\com i_1 = i_1\com i$, and the conclusions of (ii) follow.
\end{proof}

The following remark, which we learned from Mike Cole \cite{Cole3}
and Gaunce Lewis, compares certain colimits in $\sK$ and $\sU$.
It illuminates the difference between these categories and will 
be needed in our later discussion of $h$-type model structures.
\begin{rem}\mylabel{Umod}
Suppose given a sequence of inclusions $g_n\colon X_n\rtarr X_{n+1}$ and maps $f_n\colon X_n\rtarr Y$ in $\sK$ such that $f_{n+1} g_n = f_{n}$. Let $X=\text{colim}\, X_n$ and let 
$f\colon X\rtarr Y$ be obtained by passage to colimits.  Fix a map $p\colon Z\rtarr Y$.  
The maps $Z\times_Y X_n\rtarr Z\times_Y X$ induce a map 
$$\al\colon \text{colim}\, (Z\times_Y X_n) \rtarr Z\times_Y X.$$ 
Lewis has provided counterexamples showing that $\al$ need not 
be a homeomorphism in general. However, if $Y\in \sU$, then a 
result of his \cite[App.\,10.3]{Lewis0} shows that $\al$ is a 
homeomorphism for any $p$ and any maps $g_n$. In fact, 
as in \myref{second} below, if $Y\in\sU$, then the pullback 
functor $p^*\colon \sK/Y\rtarr \sK/Z$ is a left adjoint and 
therefore commutes with all colimits.  To see what goes wrong 
when $Y$ is not in $\sU$, consider the diagram
$$\xymatrix{
\text{colim}\,(Z\times_Y X_n)\ar[r]^-{\al} \ar[d]_{\io} & Z\times_Y X \ar[d] \\
\text{colim}\,(Z\times X_n) \ar[r] & Z\times X.}
$$
Products commute with colimits, so the bottom arrow is a homeomorphism,
and the top arrow $\al$ is a continuous bijection.  The right 
vertical arrow is an inclusion by the construction of pullbacks. If the
left vertical arrow $\io$ is an inclusion, then the diagram implies that $\al$
is a homeomorphism.  The problem is that $\io$ need not be an inclusion.  
One point is that the maps $Z\times_Y X_n\rtarr  Z\times X_n$ are closed inclusions if $Y$ is weak Hausdorff, but not in general otherwise. 
Now assume that all spaces in sight are in $\sU$. Since
the $g_n$ are inclusions, the relevant colimits, when computed in $\sK$, are 
weak Hausdorff and thus give colimits in $\sU$.  Therefore the commutation of
$p^*$ with colimits (which is a result about colimits in $\sK$) applies to 
these particular colimits in $\sU$ to show that $\al$ is a homeomorphism.
\end{rem}

The following related observation will be needed for applications of
Quillen's small object argument to $q$-type model structures in \S4.5
and elsewhere.

\begin{lem}\mylabel{little}  
Let $X_n\rtarr X_{n+1}$, $n\geq 0$, be a 
sequence of inclusions in $\sK$ with colimit $X$.  Suppose that $X/X_0$ 
is in $\sU$. Then, for a compact Hausdorff space $C$, the natural map
$$\text{colim}\, \sK(C,X_n)\rtarr \sK(C,X)$$
is a bijection.
\end{lem}\begin{proof}
The point is that $X_0$ need not be in $\sU$. Let $f\colon C\rtarr X$ be a map. 
Then the composite of $f$ with the quotient 
map $X\rtarr X/X_0$ takes image in some $X_n/X_0$, hence $f$ takes image in $X_n$. 
The conclusion follows.\end{proof}

\begin{sch} 
One might expect the conclusion to hold for colimits of
sequences of closed inclusions $X_{n-1}\rtarr X_n$ such that $X_n-X_{n-1}$ 
is a $T_1$ space. This is stated as \cite[4.2]{IJ}, whose authors
got the statement from May. However, Lewis has shown us a counterexample.
\end{sch}

\section{Topologically bicomplete categories and ex-objects}\label{Sbicat}

We need some standard and some not quite so standard categorical language. 
All of our categories $\sC$ will be topologically enriched, with the enrichment given by a topology on the underlying set of morphisms.  We therefore agree to write $\sC(X,Y)$\noteindex{CXY@$\sC(X,Y)$} for the space of morphisms $X\rtarr Y$ in $\sC$. Enriched category theory would have us distinguish notationally between morphism spaces and morphism sets, but we shall not do that.  A topological category $\sC$ is said to be \emph{topologically bicomplete}\index{topologically bicomplete category} if, in addition to being bicomplete in the usual sense of having all limits and colimits, it is bitensored in the sense that it is tensored and cotensored over $\sK$. We shall denote the tensors\index{tensor!with spaces} and cotensors\index{cotensor!with spaces} by $X\times K$\noteindex{XK@$X\times K$} and $\text{Map}(K,X)$\noteindex{MapXY@$\text{Map}(K,X)$} for a space $K$ and an object $X$ of $\sC$. The defining adjunction homeomorphisms are
\begin{equation}\label{tencoten1}
\sC(X\times K,Y)\iso \sK(K,\sC(X,Y))\iso \sC(X,\text{Map}(K,Y)).
\end{equation}
By the Yoneda lemma, these have many standard implications.
For example,
\begin{equation}\label{tenunit}
X\times * \iso X \quad\text{and}\quad  \text{Map}(*,Y) \iso Y,
\end{equation}
\begin{equation}\label{tenass}
X\times (K\times L)\iso (X\times K)\times L \ \ \text{and}\ \  
\text{Map}(K,\text{Map}(L,X))\iso \text{Map}(K\times L, X).
\end{equation}

We say that a bicomplete topological category $\sC$ is \emph{based}\index{based bicomplete category}\index{category!based bicomplete} if the unique map from the initial object $\emptyset$ to the terminal object $*$ is an isomorphism.  In that case, $\sC$ is enriched in the category $\sK_*$ of based $k$-spaces, the basepoint of $\sC(X,Y)$ being the unique map that factors through $*$.  We then say that $\sC$ is \emph{based topologically bicomplete}\index{based topologically bicomplete category}\index{category!based topologically bicomplete} if it is tensored and cotensored over $\sK_*$. We denote the tensors\index{tensor!with based spaces} and cotensors\index{cotensor!with based spaces} by $X\sma K$\noteindex{XK@$X\sma K$} and $F(K,X)$\noteindex{FKX@$F(K,X)$} for a based space $K$ and an object $X$ of $\sC$. The defining adjunction homeomorphisms are
\begin{equation}\label{tencoten20}
\sC(X\sma K,Y)\iso \sK_*(K,\sC(X,Y))\iso \sC(X,F(K,Y)).
\end{equation}
The based versions of (\ref{tenunit}) and (\ref{tenass}) are 
\begin{equation}\label{tenunit0}
X\sma S^0 \iso X \quad \text{and}\quad  F(S^0,Y) \iso Y,
\end{equation}
\begin{equation}\label{tenass0}
X\sma (K\sma L) \iso (X\sma K)\sma L \quad \text{and}\quad 
F(K,F(L,X))\iso F(K\sma L,X).
\end{equation}

Although not essential to our work, a formal comparison between 
the based and unbased notions of bicompleteness is illuminating. 
The following result allows us to interpret topologically bicomplete 
to mean based topologically bicomplete whenever $\sC$ is based, a
convention that we will follow throughout.

\begin{prop}\mylabel{ped1}
Let $\sC$ be a based and bicomplete topological category. Then $\sC$ is topologically bicomplete if and only if it is based topologically bicomplete.\index{based topologically bicomplete category}\index{category!based topologically bicomplete}
\end{prop}

\begin{proof}
Suppose given tensors and cotensors for unbased spaces $K$ and write them as $X\ltimes K$ and $\text{Map}(K,X)_{*}$ as a reminder that they take values in a based category.  We obtain tensors and cotensors $X\sma K$ and $F(K,X)$ for based spaces $K$ as the pushouts and pullbacks displayed in the respective diagrams
\[\xymatrix{X\ltimes {*} \ar[r] \ar[d] & X\ltimes K \ar[d]\\
 {*} \ar[r] & X\sma K}
\qquad\text{and}\qquad 
\xymatrix{F(K,X) \ar[r] \ar[d] & \text{Map}(K,X)_{*} \ar[d]\\
 {*} \ar[r] & \text{Map}({*},X)_{*}.}\]
Conversely, given tensors and cotensors $X\sma K$ and $F(K,X)$ for based spaces $K$, we obtain tensors and cotensors $X\ltimes K$ and $\text{Map}(K,X)_{*}$ for unbased spaces $K$ by setting 
$$ X\ltimes K = X\sma K_+\ \ \ \text{and}\ \ \ \text{Map}(K,X)_{*} = F(K_+,X),$$
where $K_{+}$ is the union of $K$ and a disjoint basepoint.
\end{proof}

As usual, for any category $\sC$ and object $B$ in $\sC$, we let $\sC/B$\noteindex{CBa@$\sC/B$} denote the category of objects over $B$.\index{category!of objects over B@of objects over $B$} An object $X=(X,p)$ of $\sC/B$ consists of a total object $X$ together with a projection map
$p\colon X\rtarr B$ to the base object $B$. The morphisms of $\sC/B$ are the maps of total objects that commute with the projections.

\begin{prop}\mylabel{topbicomp}
If $\sC$ is a topologically bicomplete category, then so is $\sC/B$.
\end{prop}

\begin{proof} The product of objects $Y_i$ over $B$, denoted $\times_B Y_i$, is constructed by taking 
the pullback of the product of the projections $Y_i\rtarr B$ along the 
diagonal $B\rtarr \times_i B$.  Pullbacks and arbitrary colimits of objects 
over $B$ are constructed by taking pullbacks and colimits on total objects
and giving them the induced projections. General limits are constructed as
usual from products and pullbacks.  If $X$ is an object over $B$ and 
$K$ is a space, then the tensor
$X\times_B K$ is just $X\times K$ together with the projection
$X\times K \rtarr B\times *\cong B$ induced by the projection 
of $X$ and the projection of $K$ to a point.  Note that this makes
sense even though the tensor $\times$ need have nothing to do with
cartesian products in general; see \myref{mislead} below.
The cotensor $\text{Map}_B(K,X)$ is the pullback of the diagram
\[\xymatrix{B\ar[r]^-\iota & \text{Map}(K,B) & \text{Map}(K,X)\ar[l]}\]
where $\iota$ is the adjoint of $B\times K \rtarr B\times * \cong B$.
\end{proof}

The terminal object in $\sC/B$ is $(B,\text{id})$. Let $\sC_B$\noteindex{CB@$\sC_B$} denote the category of based objects in $\sC/B$,\index{category!of ex-objects} that is, the category of objects under 
$(B,\text{id})$ in $\sC/B$. An object $X=(X,p,s)$ in $\sC_B$, which we call an \emph{ex-object over $B$}, consists of on object $(X,p)$ over $B$ together with a section $s\colon B\rtarr X$. We can therefore think of the ex-objects as retract diagrams
\[\xymatrix{B\ar[r]^s & X \ar[r]^p & B.}\]
The terminal object in $\sC_B$ is $(B,\text{id},\text{id})$, which we denote
by $*_B$; it is also an initial object. The morphisms in $\sC_B$ are the maps 
of total objects $X$ that commute with the projections and sections.  

\begin{prop}\mylabel{btopbicomp}
If $\sC$ is a topologically bicomplete category, then the category $\sC_B$ is based topologically bicomplete.
\end{prop}
\begin{proof} The coproduct of objects $Y_i\in \sC_B$, which we shall refer to as the ``wedge over $B$'' of the $Y_i$ and denote by $\wed_B Y_i$, is constructed by taking the pushout of the coproduct $\amalg B\rtarr \amalg Y_i$ of the sections along the codiagonal $\amalg_i B\rtarr B$. Pushouts and arbitrary limits of objects in $\sC_B$ are constructed by taking pushouts and limits on total objects and giving them the evident induced sections and projections. 
The tensor $X\sma_B K$ of $X=(X,p,s)$ and a based space $K$ is the pushout 
of the diagram  
\[\xymatrix{
B & (X\times {*}) \cup_B (B \times K) \ar[r] \ar[l] 
& X\times K,\\}\]
where the right map is induced by the basepoint of $K$
and the section of $X$. The cotensor $F_B(K,X)$ is the pullback 
of the diagram
\[\xymatrix{
B \ar[r]^-{s} & X &\text{Map}_B(K,X), \ar[l]_-{\epz}\\}\] 
where $\epz$ is evaluation at the basepoint of $K$, that is, the 
adjoint of the evident map $X\times K\rtarr X$ over $B$.
\end{proof}

\begin{rem}\mylabel{mislead}
Notationally, it may be misleading to write $X\times K$ and $X\sma K$ for unbased and based tensors.  It conjures up associations that are appropriate for the examples on hand but that are inappropriate in general. The tensors in a topologically bicomplete category $\sC$ may bear very little relationship to cartesian products or smash products. The standard uniform notation would be $X\otimes K$.  However, we have too many relevant examples to want a uniform notation. In particular, we later use the notations $X\times_B K$ and $X\sma_B K$ in the parametrized context, where a notation such as $X\otimes_B K$ would conjure up its own misleading associations.
\end{rem}

\section{Convenient categories of ex-spaces}

We need a convenient topologically bicomplete category of ex-spaces\footnote{Presumably the prefix ``ex'' 
stands for ``cross'', as in ``cross section''.
The unlovely term ``ex-space'' has been replaced in some
recent literature by ``fiberwise pointed space''. Used repetitively, that is
not much of an improvement. The term ``retractive space'' has also been used.}\index{ex-space}
over a space $B$, where ``convenient'' requires that we have smash 
product and function ex-space functors $\sma_B$ and $F_B$ under 
which our category is closed symmetric monoidal.
Denoting the unit $B\times S^0$ of $\sma_B$ by $S^0_B$, a 
formal argument shows that we will then have isomorphisms
\begin{equation}\label{spsp}
X\sma_B K\iso X\sma_B (S^0_B\sma_B K) \quad \text{and}\quad  F_B(K,Y)\iso F_B(S^0_B\sma_B K,Y)
\end{equation}
relating tensors and cotensors to the smash product and function ex-space functors. In particular, $S_B^0\sma_B K$ is just the product ex-space $B\times K$ with section determined by the basepoint of $K$. 

The point-set topology leading to such a convenient category is delicate, 
and there are quite a few papers devoted to this subject. They do not give exactly what we need, but they come close enough that 
we shall content ourselves with a summary. It is based on the papers 
\cite{B1, B2, B3, BB1, BB2, Lewis} of Booth, Booth and Brown, and Lewis; 
see also James \cite{James, James2}. 

We assume once and for all that our base spaces $B$ are in $\sU$.
We allow the total spaces $X$ of spaces over $B$ to be in $\sK$.
We let $\sK/B$\noteindex{K/B@$\sK/B$}\index{category!of spaces 
over B@of spaces over $B$} and $\sU/B$\noteindex{U/B@$\sU/B$} denote the categories of spaces over $B$ with total spaces in $\sK$ or $\sU$. Similarly, 
we let $\sK_B$\noteindex{KB@$\sK_B$} and $\sU_B$\noteindex{UB@$\sU_B$}\index{category!of ex-spaces} denote the respective categories of ex-spaces over $B$. 

\begin{rem}\mylabel{reasonable}
The section of an ex-space in $\sU_B$ is closed, by \myref{coflemma}.
Quite reasonably, references such as \cite{CJ, James} make the blanket
assumption that sections of ex-spaces must be closed. We have not done 
so since we have not checked that all constructions in sight preserve
this property.
\end{rem}

Both the  separation property on $B$ and the lack of a separation 
property on $X$ are dictated by consideration of the function spaces $\text{Map}_B(X,Y)$ over $B$ that we shall define shortly. These are 
only known to exist when $B$ is weak Hausdorff. However, even when $B$, 
$X$ and $Y$ are weak Hausdorff, $\text{Map}_B(X,Y)$ 
is generally not weak Hausdorff unless the projection $p\colon  X\rtarr B$ 
is an open map.  Categorically, this means that the cartesian monoidal 
category $\sU/B$ is not closed cartesian monoidal. Wishing to retain the separation property, Lewis \cite{Lewis} proposed the following as convenient 
categories of spaces and ex-spaces over a compactly generated space $B$.

\begin{defn} 
Let $\sO(B)$\noteindex{OB@$\sO(B)$} and $\sO_*(B)$\noteindex{O*B@$\sO_*(B)$} be the categories of those compactly generated
spaces and ex-spaces over $B$ whose projection maps are open.
\end{defn}

\begin{rem} 
Bundle projections over $B$ are open maps. Hurewicz fibrations 
over $B$ are open maps if the diagonal $B\rtarr B\times B$ is a
Hurewicz cofibration \cite[2.3]{Lewis}; this holds, for example, 
if $B$ is a CW complex.
\end{rem}

However, the categories $\sO(B)$ and $\sO_*(B)$ are insufficient for our purposes.  Working in these categories, we only have the base change 
adjunction $(f^*,f_*)$ of \S2.1 below for open maps 
$f\colon A\rtarr B$, which is unduly restrictive. For example, we need the adjunction $(\DE^*,\DE_*)$, where $\DE\colon B\rtarr B\times B$ is the 
diagonal map. Moreover, the generating cofibrations of our $q$-type model 
structures do not have open projection maps. This motivates 
us to drop the weak Hausdorff condition on total spaces and to focus 
on $\sK_B$ as our preferred convenient category of ex-spaces over $B$.  
The cofibrant ex-spaces in our $q$-type model structures are weak Hausdorff, hence this separation property is recovered upon cofibrant approximation. Therefore, use of $\sK$ can be viewed as scaffolding in the foundations that can be removed when doing homotopical work.

We topologize the set of ex-maps $X\rtarr Y$ as a 
subspace of the space $\sK(X,Y)$ of maps of total spaces. It is based,
with basepoint the unique map that factors through $*_B$. Therefore the 
category $\sK_B$ is enriched over $\sK_*$.  It is based topologically bicomplete by \myref{topbicomp}.  Recall that we write $\times_B Y_i$ and $\wed_B Y_i$ for
products and wedges over $B$. We also write $Y/\!_BX$ for quotients, which are
understood to be pushouts of diagrams $*_B \longleftarrow X \rtarr Y$. We give 
a more concrete description of the tensors and cotensors in $\sK/B$ and $\sK_B$ given by \myref{topbicomp} and \myref{btopbicomp}. For a space $X$ over $B$, 
we let $X_b$ denote the fiber $p^{-1}(b)$. 
If $X$ is an ex-space, then $X_b$ has the basepoint $s(b)$. 

\begin{defn}\mylabel{exctp}
Let $X$ be a space over $B$ and $K$ be a space. Define $X\times_{B}K$\noteindex{XxBK@$X\times_B K$}\index{tensor!for spaces over B@for spaces over $B$} to be the space $X\times K$ with projection the product of the
projections $X\rtarr B$ and $K\rtarr *$.  Define $\text{Map}_B(K,X)$\noteindex{MapBKX@$\text{Map}_B(K,X)$}\index{cotensor!for spaces over B@for spaces over $B$} to be the subspace of $\text{Map}(K,X)$ 
consisting of those maps $f\colon K\rtarr X$ that factor through some 
fiber $X_b$; the projection sends such a map $f$ to $b$.
\end{defn}

\begin{defn}\mylabel{exct} 
Let $X$ be an ex-space over $B$ and $K$ be a based space. Define $X\sma_B K$\noteindex{XBK@$X\sma_B K$}\index{tensor!for ex-spaces} 
to be the quotient of $X\times_B K$ obtained by taking fiberwise smash products, so that $(X\sma_B K)_b = X_b\sma K$; the basepoints of fibers prescribe the section. Define $F_B(K,X)$\noteindex{FBKX@$F_B(K,X)$}\index{cotensor!for ex-spaces} to be the subspace of $\text{Map}_B(K,X)$ consisting of the based maps $K\rtarr X_b\subset X$ for some $b\in B$, so that 
$F_B(K,X)_b = F(K,X_b)$; the section sends $b$ to the constant map at $s(b)$.
\end{defn}
\begin{rem}
As observed by Lewis \cite[p.\,85]{Lewis}, if $p$ is an open map,
then so are the projections of $X\sma_B K$ and $F_B(K,Y)$. Therefore
$\sO_*(B)$ is tensored and cotensored over $\sT$. 
\end{rem}

The category $\sK/B$ is closed cartesian monoidal under the fiberwise cartesian
product $X\times_B Y$ and the function space $\text{Map}_B(X,Y)$ over $B$. The category $\sK_B$ is closed symmetric monoidal under the fiberwise smash product $X\sma_B Y$ and the function ex-space $F_B(X,Y)$. We recall the definitions. 

\begin{defn} For spaces $X$ and $Y$ over $B$, 
$X\times_B Y$\noteindex{XxBY@$X\times_B Y$}\index{fiberwise product} is the pullback of the projections $p\colon X\rtarr B$ and $q\colon Y\rtarr B$, with the evident projection $X\times_B Y\rtarr B$. When $X$ and $Y$ have sections $s$ and $t$, their pushout $X\vee_B Y$ specifies the coproduct, or wedge, of $X$ and $Y$ in $\sK_B$, and $s$ and $t$ induce a map $X\vee_B Y\rtarr X\times_B Y$ over $B$ that sends $x$ and $y$ to $(x,tp(x))$ and $(sq(y),y)$.  Then $X\sma_B Y$\noteindex{XBY@$X\sma_B Y$}\index{fiberwise smash product} is the pushout in $\sK/B$ displayed in the diagram
$$\xymatrix{
X\vee_B Y \ar[r] \ar[d] & X\times_B Y \ar[d]\\ 
{*}_B  \ar[r] & X\sma_B Y.}$$
This arranges that $(X\sma_B Y)_b = X_b\sma Y_b$, and the section and projection are evident. 
\end{defn}

The following result is \cite[8.3]{BB2}.

\begin{prop}\mylabel{prop:smaB} If $X$ and $Y$ are weak Hausdorff ex-spaces over $B$, then so is $X\sma_BY$. That is, $\sU_B$ is closed under $\sma_B$. 
\end{prop}

Function objects are considerably more subtle, and we need a preliminary definition in order to give the cleanest description.

\begin{defn}\mylabel{partialtilde} 
For a space $Y\in \sK$, define the \emph{partial map classifier}\index{partial map classifier}
$\tilde Y$\noteindex{tY@$\tilde{Y}$} to be the union of $Y$ and a disjoint point $\om$, with the topology
whose closed subspaces are $\tilde Y$ and the closed subspaces of $Y$. The point $\om$ is not a closed subset, and $\tilde{Y}$ is not weak Hausdorff. The name
``partial map classifier'' comes from the observation that, for any space $X$, pairs $(A,f)$ consisting of a closed subset $A$ of $X$ and a continuous map $f\colon  A\rtarr Y$ are in bijective correspondence with continuous maps $\tilde{f}\colon  X\rtarr \tilde Y$. Given $(A,f)$, $\tilde f$ restricts
to $f$ on $A$ and sends $X-A$ to $\om$; given $\tilde{f}$, $(A,f)$ is 
$\tilde{f}^{-1}(Y)$ and the restriction of $\tilde{f}$.
\end{defn}

\begin{defn}\mylabel{MapB}
Let $p\colon X\rtarr B$ and $q\colon Y\rtarr B$ be spaces over $B$. 
Define $\text{Map}_B(X,Y)$\noteindex{MapBXY@$\text{Map}_B(X,Y)$}\index{mapping space!of spaces over B@of spaces over $B$} to be the pullback displayed in the diagram
$$\xymatrix{
\text{Map}_B(X,Y) \ar[r] \ar[d] &  \text{Map}(X,\tilde{Y}) \ar[d]^{\text{Map}(\text{id},\tilde{q})}  \\
B \ar[r]_-{\la}&  \text{Map}(X, \tilde{B}).\\}$$
Here $\la$ is the adjoint of the map $X\times B\rtarr \tilde{B}$
that corresponds to the composite of the inclusion $\text{Graph}(p)\subset X\times B$
and the projection $X\times B\rtarr B$ to the second coordinate. The graph of $p$ is the 
inverse image of the diagonal under $p\times \text{id}\colon  X\times B\rtarr B\times B$, and
the assumption that $B$ is weak Hausdorff ensures that it is a closed subset of 
$X\times B$, as is needed for the definition to make sense. Explicitly, $\la(b)$ sends
$X_b$ to $b$ and sends $X-X_b$ to the point $\om\in \tilde{B}$.
\end{defn}

This definition gives one reason that we require the base spaces of ex-spaces
to be weak Hausdorff. On fibers, $\text{Map}_B(X,Y)_b = \text{Map}(X_b,Y_b)$. The space of 
sections of $\text{Map}_B(X,Y)$ is $\sK/B(X,Y)$. We have (categorically equivalent) 
adjunctions
\begin{gather}\label{maptimes1}
\text{Map}_B(X\times_B Y, Z)\iso \text{Map}_B(X,\text{Map}_B(Y,Z)),\\[1ex]
\label{maptimes2}
\sK/B\,(X\times_B Y, Z)\iso \,\sK/B\,(X,\text{Map}_B(Y,Z)).
\end{gather}
These results are due to Booth \cite{B1, B2, B3}; see also \cite[\S7]{BB1}, 
\cite[\S8]{BB2}, \cite[II\S9]{James}, \cite{Lewis}. 

Examples in \cite[5.3]{BB1} and \cite[1.7]{Lewis} show that $\text{Map}_B(X,Y)$
need not be weak Hausdorff even when $X$ and $Y$ are. The question of when
$\text{Map}_B(X,Y)$ is Hausdorff or weak Hausdorff was studied in \cite[\S5]{BB1}
and later in \cite{James, James2}, but the definitive criterion was given by Lewis 
\cite[1.5]{Lewis}.

\begin{prop}\mylabel{open} Consider a fixed map $p\colon X\rtarr B$ and varying maps \linebreak
$q\colon Y\rtarr B$, where $X$ and the $Y$ are weak Hausdorff. The map $p$ is 
open if and only if the space $\text{Map}_B(X,Y)$ is weak Hausdorff for all $q$.
\end{prop}

\begin{prop}\mylabel{Mapfib} 
If $p\colon X\rtarr B$ and $q\colon Y\rtarr B$ are Hurewicz fibrations,
then the projections $X\times_B Y\rtarr B$ and $\text{Map}_B(X,Y)\rtarr B$ 
are Hurewicz fibrations. The second statement is false with Hurewicz 
fibrations replaced by Serre fibrations.
\end{prop}\begin{proof}
The statement about $X\times_BY$ is clear. The statements about \linebreak
$\text{Map}_B(X,Y)$
are due to Booth \cite[6.1]{B1} or, in the present formulation \cite[3.4]{B2}; 
see also \cite[23.17]{James}.\end{proof}

\begin{defn} 
For ex-spaces $X$ and $Y$ over $B$, define $F_B(X,Y)$\noteindex{FBXY@$F_B(X,Y)$}\index{mapping space!of ex-spaces} to be the
subspace of $\text{Map}_B(X,Y)$ that consists of the points that restrict to based maps
$X_b\rtarr Y_b$ for each $b\in B$; the section sends $b$ to the constant map from 
$X_b$ to the basepoint of $Y_b$. Formally, $F_B(X,Y)$ is the pullback
displayed in the diagram 
$$\xymatrix{
F_B(X,Y) \ar[r] \ar[d] & \text{Map}_B(X,Y) \ar[d]^{\text{Map}_B(s,\text{id})} \\
B\ar[r]_-{t} & Y\iso \text{Map}_B(B,Y),}\\$$
where $s$ and $t$ are the sections of $X$ and $Y$.
\end{defn}

The space of maps $S^0_B\rtarr F_B(X,Y)$ is $\sK_B(X,Y)$, and we have adjunctions
\begin{gather}\label{maptimes1*}
F_B(X\sma_B Y, Z)\iso F_B(X,F_B(Y,Z)),\\[1ex]
\label{maptimes2*}
\sK_B\,(X\sma_B Y, Z)\iso \,\sK_B\,(X,F_B(Y,Z)).
\end{gather}

\myref{open} implies the following analogue of  \myref{prop:smaB}. 

\begin{prop} If $X$ and $Y$ are weak Hausdorff ex-spaces over $B$ and
$X\rtarr B$ is an open map, then  $F_B(X,Y)$ is weak Hausdorff.
\end{prop}

We record the following analogue of \myref{Mapfib}. The second part is again 
due to Booth, who sent us a detailed write-up. The argument is similar to his proofs
in \cite[6.1(i)]{B1} or \cite[3.4]{B2}, but a little more complicated, and a general
result of the same form is given by Morgan \cite{Morgan}.

\begin{prop}\mylabel{Ffib}
If $X$ and $Y$ are ex-spaces over $B$ whose sections are Hurewicz 
cofibrations and whose projections are Hurewicz fibrations, then the 
projections of $X\sma_B Y$ and $F_B(X,Y)$ are Hurewicz fibrations.
\end{prop}

\section{Convenient categories of ex-$G$-spaces}

The discussion just given generalizes readily to the equivariant context. Let $G$ be a compactly generated topological group. Subgroups of $G$ are understood to be closed. Let $B$ be a compactly generated $G$-space (with $G$ acting from the left). We consider $G$-spaces over $B$ and ex-$G$-spaces $(X,p,s)$. The total space $X$ is a $G$-space in $\sK$, and the section and projection are $G$-maps. The fiber $X_b$ is a based $G_b$-space with $G_b$-fixed basepoint $s(b)$, where $G_b$ is the isotropy group of $b$. 

Recall from \cite[II\S1]{MM} the distinction between the category $\sK_G$\noteindex{KG@$\sK_G$}\index{category!of G spaces@of $G$-spaces} of $G$-spaces and non\-equi\-var\-iant maps and the category $G\sK$\noteindex{GK@$G\sK$} of $G$-spaces and equivariant maps; the former is enriched over $G\sK$, the latter over $\sK$. We have a similar dichotomy on the ex-space level. Here we have a conflict of notation with our notation for categories of ex-spaces, and we agree to let $\sK_{G,B}$ denote the category whose objects are the ex-$G$-spaces over $B$ and whose morphisms are the maps of underlying ex-spaces over $B$, that is, the maps $f\colon X\rtarr Y$ such that $f\com s = t$ and $q\com f = p$. Henceforward, we call these maps ``arrows'' to distinguish them from $G$-maps, which we often abbreviate to maps. For $g\in G$, $gf$ is also an arrow of ex-spaces over $B$, so that $\sK_{G,B}(X,Y)$\noteindex{KGBXY@$\sK_{G,B}(X,Y)$} is a $G$-space.  Moreover, composition is given by $G$-maps 
$$\sK_{G,B}(Y,Z)\times \sK_{G,B}(X,Y)\rtarr \sK_{G,B}(X,Z).$$ 
We obtain the category $G\sK_B$\noteindex{GKB@$G\sK_B$}\index{category!of exGspaces@of ex-$G$-spaces} by restricting to $G$-maps $f$, and we may view it as the $G$-fixed point category of $\sK_{G,B}$. Of course, $G\sK_B(X,Y)$ is a space and not a $G$-space. The pair $(\sK_{G,B}, G\sK_B)$ is an example of a \emph{$G$-category}\index{category!G-@$G$- ---}, a structure that we shall recall formally in \S10.2.

Since $*_B$ is an initial and terminal object in both $\sK_{G,B}$ and $G\sK_B$, their morphism spaces are based. Thus $\sK_{G,B}$ is enriched over the category $G\sK_*$ of based $G$-spaces and $G\sK_B$ is enriched over $\sK_*$. As discussed in \cite[II.1.3]{MM}, if we were to think exclusively in enriched category terms, we would resolutely ignore the fact that the $G$-spaces $\sK_{G,B}(X,Y)$ have elements (arrows), thinking of these $G$-spaces as enriched hom objects. From that point of view, $G\sK_B$ is the ``underlying category'' of our enriched $G$-category.  While we prefer to think of $\sK_{G,B}$ as a category, it must be kept in mind that it is not a very well-behaved one. For example, because its arrows are not equivariant, it fails to have limits or colimits. 

In contrast, the category $G\sK_B$ is bicomplete.  Its limits and colimits are constructed in $\sK_B$ and then given induced $G$-actions. The category $\sK_{G,B}$, although not bicomplete, is tensored and cotensored over $\sK_{G,*}$.  The tensors $X\sma_B K$ and cotensors $F_B(K,X)$ are constructed in $\sK_B$ and then given induced $G$-actions. They satisfy the adjunctions
\begin{gather}\label{tensored1}
\sK_{G,B}(X\sma_B K, Y)\iso \sK_{G,*}(K, \sK_{G,B}(X,Y))\iso \sK_{G,B}(X, F_B(K,Y))\\
\intertext{and, by passage to fixed points,}
\label{tensored2}
G\sK_B(X\sma_B K, Y)\iso G\sK_*(K,\sK_{G,B}(X,Y))\iso G\sK_B(X, F_B(K,Y)).
\end{gather} 
It follows that $G\sK_{B}$ is tensored and cotensored over $G\sK_*$ and, in particular, is topologically bicomplete.

The category $\sK_{G,B}$ is closed symmetric monoidal via the fiberwise smash products $X\sma_B Y$ and function objects $F_B(X,Y)$. Again, these are defined in $\sK_B$ and then given induced $G$-actions. The unit is the ex-$G$-space $S^0_B=B\times S^0$. The category $G\sK_B$ inherits a structure of closed symmetric monoidal category. We have homeomorphisms of based $G$-spaces
\begin{gather}\label{ad1}
\sK_{G,B}(X\sma_B Y, Z)\iso \sK_{G,B}(X, F_B(Y,Z))
\intertext{and, by passage to $G$-fixed points, homeomorphisms of based spaces}
\label{ad2}
G\sK_{B}(X\sma_B Y, Z)\iso G\sK_{B}(X, F_B(Y,Z)).
\intertext{The first of these implies an associated homeomorphism of ex-$G$-spaces}
\label{ad1'}
F_B(X\sma_B Y, Z)\iso F_B(X, F_B(Y,Z)).
\end{gather}

Nonequivariantly, the functor that sends an ex-space $X$ over $B$ to the fiber $X_b$ has a left adjoint, denoted $(-)^b$. It sends a based space $K$ to the wedge $K^b = B\wed K$, where $B$ is given the basepoint $b$; the section and projection are evident.  Nonobviously, the same set $B\wed K$ admits a quite different topology under which it gives a {\em right}\, adjoint to the fiber functor $X\mapsto X_b$. We shall prove the equivariant analogue conceptually in \myref{Johann0}, but we describe the left adjoint to the fiber functor explicitly here.

\begin{con}\mylabel{Fibad} 
Let $b\in B$. Then the functor 
$G\sK_{B} \rtarr G_b\sK_{*}$ that sends $Y$ to $Y_b$ has a
left adjoint. It sends a based $G_b$-space $K$ to the
ex-$G$-space $K^b$ given by the pushout
$$K^b = (G\times_{G_b} K)\cup_G B.$$
Here $G\times_{G_b}K$ is the (left) $G$-space $(G\times K)/\sim$, 
where $(gh,k)\sim (g,hk)$ for $g\in G$, $h\in G_b$, and $k\in K$.
The pushout is defined with respect to the map $G\rtarr B$
that sends $g$ to $gb$ and the map $G\rtarr G\times_{G_b}K$ that
sends $g$ to $(g,k_0)$, where $k_0$ is the ($G_b$-fixed) 
basepoint of $K$. 
The section is given by the evident inclusion of $B$ and
the projection is obtained by passage to pushouts from the identity
map of $B$ and the $G$-map $\pi_b\colon G\times_{G_b} K\rtarr B$
given by $\pi_b(g,k) = gb$. Thus we first extend the group action 
on $K$ from $G_b$ to $G$ and then glue the orbit of the basepoint
of $K$ to the orbit of $b$. If $K$ is an unbased $G_b$-space, then 
$(K_+)^b = (G\times_{G_b}K)\amalg B$.
\end{con}

\begin{rem}
There is an alternative parametrized view of equivariance that is important in torsor theory but that we shall not study. It focuses on ``topological groups $G_B$ over $B$'' and ``$G_B$-spaces $E$ over $B$'',  
where $G_B$ is a space over a nonequivariant space $B$ with a product 
$G_B\times _B G_B\rtarr G_B$ that restricts on fibers to the products of topological groups $G_b$ and $E$ is a space over $B$ with 
an action $G_B\times_B E\rtarr E$ that restricts on fibers to actions
$G_b\times E_b\rtarr E_b$.   That theory intersects ours in the special case 
$G_B = G\times B$ for a topological group $G$. Since, at least implicitly, all 
of our homotopy theory is done fiberwise, our work adapts without essential difficulty to give a development of parametrized equivariant homotopy theory 
in that context.
\end{rem}

\section{Appendix: nonassociativity of smash products in ${\sT}op_*$}
\label{sec:topass}

In a 1958 paper \cite{Puppe}, Puppe asserted the following result, 
but he did not give a proof. It was the subject of a series of e-mails 
among Mike Cole, Tony Elmendorf, Gaunce Lewis and the first author. Since we know of no published source that gives the details of this 
or any other counterexample to the associativity of the smash product 
in ${\sT}op_*$, we include the following proof. It is due 
to Kathleen Lewis.

Let $\bQ$ and $\bN$ be the rational numbers and the nonnegative integers, 
topologized as subspaces of $\bR$ and given the basepoint zero. Consider
smash products as quotient spaces, without applying the $k$-ification 
functor. Then we have the following counterexample to associativity. 

\begin{thm} $(\bQ\sma \bQ)\sma \bN$ is not homeomorphic to $\bQ\sma(\bQ\sma \bN)$.
\end{thm}\begin{proof}
Consider the following diagram.
$$\xymatrix{
& \bQ\times \bQ\times \bN \ar[dr]^{p\times\text{id}} \ar[dl]_{\text{id}\times p'} \ar[dd]^{q} & \\
\bQ\times (\bQ\wedge \bN) \ar[d]_{s} &
 & (\bQ\wedge \bQ) \times \bN \ar[d]^-{r} \\
\bQ\wedge(\bQ\wedge \bN) 
  & \bQ\wedge \bQ\wedge \bN \ar[l]_-{t} \ar[r]^-{\cong} 
     & (\bQ\wedge \bQ) \wedge \bN}$$
Here $\bQ\sma \bQ\sma \bN$ denotes the evident quotient space of $\bQ\times \bQ\times \bN$.
The maps $p$, $p'$, $q$, $r$, and $s$ are quotient maps. Since $\bN$ is locally compact,
$p\times \text{id}$ is also a quotient map, hence so is $r\com(p\times \text{id})$. The universal 
property of quotient spaces then gives the bottom right homeomorphism. Since $\bQ$
is not locally compact, $\text{id}\times p'$ need not be a quotient map, and in fact it
is not. The map $t$ is a continuous bijection given by the universal property of
the quotient map $q$, and we claim that $t$ is not a homeomorphism. To show this,
we display an open subset of $\bQ\sma\bQ\sma N$ whose image under $t$ is not open.

Let $\be$ be an irrational number, $0<\be<1$, and let $\ga = (1-\be)/2$. Define
$V'(\be)$ to be the open subset of $\bR\times \bR$ that is the union of the
following four sets.
\begin{enumerate}[(1)]
\item The open ball of radius $\be$ about the origin
\item The tubes $[1,\infty)\times (-\ga,\ga)$, $(-\infty,-1]\times (-\ga,\ga)$,
$(-\ga,\ga)\times [1,\infty)$, and $(-\ga,\ga)\times (-\infty,-1]$ of width $2\ga$
about the axes.
\item The open balls of radius $\ga$ about the four points $(\pm 1,0)$, $(0,\pm 1)$.
\item For each $n\geq 1$, the open ball of radius $\ga/2^n$ about the four points
$(\pm \ga_n,0)$, $(0,\pm \ga_n)$, where $\ga_n = 1 - \sum_{k=0}^{k=n-1}\ga/2^k$.
\end{enumerate}

To visualize this set, it is best to draw a picture. It is symmetric with respect
to $90$ degree rotation. Consider the part lying along the positive $x$-axis. A tube
of width $2\ga$ covers the part of the $x$-axis to the right of $(1,0)$. A ball of
radius $\be$ centers at the origin. A ball of radius $\ga$ centers at $(1,0)$. Its
vertical diagonal is the edge of the tube going off to the right. On the left, by 
the choice of $\ga$, this ball reaches halfway from its center $(1,0)$ to the point 
$(\be,0)$ at the right edge of the ball centered at the origin. The point $(1-\ga,0)$ 
at the left edge of the ball centered at $(1,0)$ is the center of another ball, which 
reaches half the distance from $(1-\ga,0)$ to $(\be,0)$. And so on: 
the point where the left edge of the $n$th ball crosses the $x$-axis is the
center point of the $(n+1)$st ball, which reaches half the distance from its 
center to the edge of the ball centered at the origin.

Define $V(\be)=V'(\be)\cap (\bQ\times \bQ)$. Note that the only points of
the coordinate axes of $\bR\times \bR$ that are not in  $V'(\be)$ are $(\pm \be,0)$ 
and $(0,\pm \be)$. Since $\be$ is irrational, $V(\be)$ contains the coordinate axes 
of $\bQ\times\bQ$.  Because the radii of the balls in the sequence are decreasing, 
for each $\epz > \be$, there is no $\de >0$ such that 
$((-\epz,\epz)\times (-\de,\de))\cap (\bQ\times\bQ)$ is contained in $V(\be)$.

Now let $\al$ be an irrational number, $0<\al<1$. Let $\bullet$ be the basepoint of
$\bQ\sma \bN$ and $*$ be the basepoint of $\bQ\sma\bQ\sma\bN$. Let $U$ be the 
union of $\{*\}$ and the image under $q$ of $\cup_{n\geq 1} V(\al/n)\times\{n\}$. 
This is an open subspace of $\bQ\sma\bQ\sma \bN$ since 
$$q^{-1}(U) = \bQ\times \bQ\times \{0\} \cup (\cup_{n\geq 1} V(\al/n)\times \{n\})$$
is an open subset of $\bQ\times\bQ\times \bN$. We claim that $t(U)$ is not open in 
$\bQ\sma(\bQ\sma\bN)$. Assume that $t(U)$ is open. Then
$$s^{-1}(t(U)) = (\text{id}\times p')(q^{-1}(U))$$
is an open subset of $\bQ\times (\bQ\sma \bN)$, hence it contains an open neighborhood 
$V$ of $(0,\bullet)$. Now $V$ must contain $((-\epz,\epz)\cap \bQ) \times W$ for some $\epz>0$ 
and some open neighborhood $W$ of $\bullet$ in $\bQ\sma \bN$. Since $\bQ\sma\bN$ is 
homeomorphic to the wedge over $n\geq 1$ of the spaces $\bQ\times \{n\}$,
$W$ must contain the wedge over $n\geq 1$ of subsets $((-\de_n,\de_n)\cap \bQ)\times\{n\}$, 
where $\de_n>0$. By the definition of $U$, this implies that 
$$((-\epz,\epz)\times (-\de_n,\de_n))\cap(\bQ\times \bQ) \subset V(\al/n).$$
However, for $n$ large enough that $\epz > \al/n$, there is no $\de_n$ for which this holds.
\end{proof}

\chapter{Change functors and compatibility relations}

\section*{Introduction}

In the previous chapter, we developed the internal properties of the category $G\sK_B$ of ex-$G$-spaces over $B$. As $B$ and $G$ vary, these categories are related by various functors, such as base change functors, change of groups functors, orbit and fixed point functors, external smash product and function space functors, and so forth.  We define these ``change functors'' 
and discuss various compatibility relations among them in this chapter.  

We particularly emphasize base change functors. We give a general 
categorical discussion of such functors in \S2.1, illustrating 
the general constructions with topological examples. In \S2.2, we 
discuss various compatibility relations that relate 
these functors to smash products and function objects.

In \S2.3 and \S2.4 we turn to equivariant phenomena and study restriction of group actions along homomorphisms. As usual, we break this into the study of
restriction along inclusions and pullback along quotient homomorphisms. 

In \S2.3, we discuss restrictions of group actions to subgroups, together 
with the associated induction and coinduction functors. We also consider their compatibilities with base change functors.  In particular, this gives us a convenient way of thinking about passage to fibers and allows us to reinterpret restriction to subgroups in terms of base change and coinduction.  That is the starting point of our generalization of the Wirthm\"uller isomorphism in 
Part IV.

In \S2.4, we consider pullbacks of group actions from a quotient group
$G/N$ to $G$, together with the associated quotient and fixed point functors.
Again, we also consider compatibilities with base change functors. For an $N$-free base space $E$, we find a relation between the quotient functor $(-)/N$ and the fixed point functor $(-)^N$ that involves base change along the quotient map $E\rtarr E/N$. The good properties of the bundle construction in Part IV
can be traced back to this relation, and it is at the heart of the Adams isomorphism in equivariant stable homotopy theory.

In \S2.5, we describe a different categorical framework, one appropriate to ex-spaces with varying base spaces. We show that the relevant category of retracts over varying base spaces is closed symmetric monoidal under external smash product and function ex-space functors. The internal smash product and function ex-space functors are obtained from these by use of base change along diagonal maps.  The external smash products are much better behaved homotopically than the internal ones, and homotopical analysis of base change functors will therefore play a central role in the homotopical analysis of smash products.

In much of this chapter, we work in a general categorical framework. In some places where we restrict to spaces, more general categorical formulations are undoubtedly possible.  When we talk about group actions, all groups are 
assumed to be compactly generated spaces but are otherwise unrestricted.

\section{The base change functors $f_{!}$, $f^*$, and $f_*$}

Let $f\colon  A\rtarr B$ be a map in a bicomplete subcategory $\sB$ of a 
bicomplete category $\sC$. We are thinking of $\sU\subset \sK$ 
or $G\sU\subset G\sK$. We wish to define functors 
$$ f_{!}\colon  \sC_A\rtarr \sC_B,\qquad f^*\colon \sC_B\rtarr \sC_A,  \qquad f_*\colon \sC_A\rtarr \sC_B,$$
such that $f_{!}$ is left adjoint and $f_{*}$ is right adjoint to $f^*$. The definitions
of $f^*$ and $f_{!}$ are dual and require no further hypotheses. The definition of $f_*$
does not work in full generality, but it only requires the further hypothesis that $\sC/B$ 
be cartesian closed. Thus 
we assume given internal hom objects $\text{Map}_B(Y,Z)$ in $\sC/B$ that satisfy the usual
adjunction, as in (\ref{maptimes2}). One reason to work in this generality is to
emphasize that no further point-set topology is needed to construct these base
change functors in the context of ex-spaces. This point is not clear from the literature, where the functor $f_*$ 
is often given an apparently different, but naturally isomorphic, description. 
We work with generic ex-objects 
\[\xymatrix{A\ar[r]^s & X \ar[r]^p & A}
\qquad\text{and}\qquad
\xymatrix{B\ar[r]^t & Y \ar[r]^q & B}\]
in this section.

\begin{defn}\mylabel{retract1}
Define $f_{!}X$\noteindex{fXl@$f_{"!}X$} and its structure maps $q$ and $t$ by means of the map of retracts in the following diagram on the left, where the top square is a pushout and the bottom square is defined by the universal property of pushouts and the requirement that $q\com t= \text{id}$. Define $f^*Y$\noteindex{fYm@$f^*Y$} and its structure maps $p$ and $s$ by means of the map of retracts in the following middle diagram, where the bottom square is a pullback and the top square is defined by the universal property of pullbacks and the requirement that $p\com s = \text{id}$. 
\[\xymatrix{
A\ar[d]_{s} \ar[r]^-{f}   &  B\ar[d]^{t}\\
X\ar[d]_{p} \ar[r]  &   f_{!}X \ar[d]^{q}\\ 
A\ar[r]_-{f}  &  B}
\qquad\qquad
\xymatrix{
A \ar[d]_{s}  \ar[r]^-{f} & B\ar[d]^{t} \\
f^*Y \ar[r] \ar[d]_{p} & Y \ar[d]^{q} \\
A \ar[r]_-{f} & B}
\qquad\qquad
\xymatrix{
B\ar[d]_{t} \ar[r]^-{\io}  & \text{Map}_B(A,A)\ar[d]^{\text{Map}(\text{id},s)}\\
f_*X \ar[d]_{q} \ar[r]  &  \text{Map}_B(A,X)\ar[d]^{\text{Map}(\text{id},p)}\\ 
B \ar[r]_-{\io} & \text{Map}_B(A,A)}\]
Thinking of $X$ and $A$ as objects over $B$ via $f\com p$ and $f$ and observing that the adjoint of the identity map of $A$ gives a map $\io\colon  B\rtarr \text{Map}_B(A,A)$, define $f_*X$\noteindex{fXr@$f_*X$} and its structure maps $q$ and $t$ by means of the map of retracts in the above diagram on the right, where the bottom square is a pullback and the top square is defined by the universal property of pullbacks and the requirement that $q\com t =\text{id}$.
\end{defn}

\begin{prop}\mylabel{first} $(f_{!},f^*)$ is an adjoint pair of functors:
$$\sC_B(f_{!}X,Y)\iso \sC_A(X,f^*Y).$$
\end{prop}\begin{proof}
Maps in both hom sets are specified by maps $k\colon X\rtarr Y$ in $\sC$ 
such that $q\com k = f\com p$ and $k\com s = t\com f$.\end{proof}

\begin{prop}\mylabel{second} $(f^*,f_*)$ is an adjoint pair of functors:
$$\sC_A(f^*Y,X)\iso \sC_B(Y,f_*X).$$
\end{prop}\begin{proof}
A map $k\colon f^*Y = Y\times_B A \rtarr X$ such that $p\com k = p$ and $k\com s = s$
has adjoint $\tilde{k} \colon Y\rtarr \text{Map}_B(A,X)$ such that
$\text{Map}(\text{id},p)\com \tilde k= \io\com q$ and $\tilde{k}\com t = \text{Map}(\text{id},s)\com \io$.
The conclusion follows directly.\end{proof}

\begin{rem} 
Writing these proofs diagrammatically, we see that the adjunction 
isomorphisms are given by homeomorphisms in our context of topological categories.
\end{rem}

We specialize to ex-spaces (or ex-$G$-spaces), in the rest of the section. 
Observe that the fiber $(f_*X)_b$ is the space of sections $A_b\rtarr X_b$ of $p\colon X_b\rtarr A_b$. 

\begin{rem} 
If $f\colon A\rtarr B$ is an open map and $X$ 
is in $\sU$, then $f_*X$ is in $\sU$ and $\sU_A(f^*Y,X)\iso \sU_B(Y,f_*X)$ 
for $Y\in \sU$, by \cite[1.5]{Lewis}.
\end{rem}

\begin{exmp}\mylabel{incsting}
Let $f\colon A\rtarr B$ be an inclusion. Then $f^*Y$
is the restriction of $Y$ to $A$ and $f_!X = B\cup_A X$. The ex-space
$f_*X$ over $B$ is analogous to the prolongation by zero of a sheaf
over $A$. The fiber $(f_*X)_b$ is $X_a$ if $a\in A$ and a point 
$\{b\}$ otherwise. To see this from the definition, recall that 
$\text{Map}(\emptyset, K)$ is a  point for any space $K$ and that 
$\text{Map}_B(A,X)_b = \text{Map}(A_b,X_b)$. As a set, $f_*X \iso B\cup_A X$,
but the topology is quite different. It is devised so that the 
map $Y\rtarr f_*f^*Y$ that restricts to the identity on $Y_a$ for 
$a\in A$ but sends $Y_b$ to $\{b\}$ for $b\notin A$ is continuous.
\end{exmp}

\begin{exmp}\mylabel{r!ex}
Let $r\colon B\rtarr *$ be the unique map. For
a based space $X$ and an ex-space $E = (E,p,s)$ over $B$, we have
$$r^*X = B\times X, \qquad r_{!}E = E/s(B), \qquad \text{and}\qquad r_*E= \text{Sec}(B,E),$$ 
where $\text{Sec}(B,E)$ is the space of maps $t\colon B\rtarr E$ such that 
$p\com t=\text{id}$, 
with basepoint the section $s$.  These elementary base change functors are the key to 
using parametrized homotopy theory to obtain information in ordinary homotopy theory.
Let $\epz\colon r_!r^*\rtarr \text{id}$ and $\et\colon \text{id}\rtarr r^*r_!$ be the counit
and unit of the adjunction $(r_!,r^*)$. Then $r_!r^*X\iso B_+\sma X$ and $\epz$ 
is $r_+\sma\text{id}$.  Similarly, $r_!r^*r_! E\iso B_+\sma E/B$, and 
$r_!\et\colon r_!E\rtarr r_!r^*r_!E$ is the ``Thom diagonal''\index{Thom diagonal} $E/B\rtarr B_+\sma E/B$. 
If $p\colon E\rtarr B$ is a spherical fibration with section, such as the fiberwise 
one-point compactification 
of a vector bundle, then $r_{!}E$ is the Thom complex\index{Thom complex} of $p$.
\end{exmp}

\section{Compatibility relations}

The term ``compatibility relation'' has been used in algebraic geometry in
the context of Grothendieck's six functor formalism that relates base change 
functors to tensor product and internal hom functors in sheaf theory. We describe how
the analogous, but simpler, formalism appears in our categories of ex-objects.

We recall some language. We are especially interested in the behavior of base change functors with respect to closed symmetric monoidal structures that, in
our topological context, are given by smash products and function objects.  Relevant categorical observations are given in \cite{FHM}.  We say that a functor $T\colon\sB\rtarr \sA$ between closed symmetric monoidal categories 
is {\em closed symmetric monoidal} if 
\[TS_{\sB}\iso S_{\sA},\quad T(X\sma_{\sB} Y)\iso TX\sma_{\sA} TY, 
\quad\text{and}\quad 
TF_\sB(X,Y)\iso F_{\sA}(TX,TY),\] 
where $S_{\sB}$, $\sma_{\sB}$ and $F_{\sB}$  denote the unit object, product, and internal hom of $\sB$, and similarly for $\sA$. These isomorphisms must satisfy appropriate coherence conditions. In the language of \cite{FHM}, the following result states that any map $f$ of base spaces gives rise to a ``Wirthm\"uller context",\index{Wirthmuller context@Wirthm\"uller context} 
which means that the functor $f^*$ is closed symmetric monoidal and has both a left adjoint and a right adjoint.

\begin{prop}\mylabel{Wirth0}
If $f\colon A\rtarr B$ is a map of base $G$-spaces, then the functor $f^*\colon G\sK_B\rtarr G\sK_A$ is closed symmetric monoidal. Therefore, by definition and implication, $f^*S^0_B\iso S^0_A$ and there are natural isomorphisms
\begin{gather}\label{oneo}
f^*(Y\sma_B Z)\iso f^*Y\sma_A f^*Z,\\[1ex]
\label{twoo}
F_B(Y,f_*X) \iso f_*F_A(f^*Y,X),\\[1ex]
\label{three0}
f^*F_B(Y,Z)\iso F_A(f^*Y,f^*Z),\\[1ex]
\label{four0}
f_{!}(f^*Y\sma_A X)\iso Y\sma_B f_{!}X,\\[1ex]
\label{five0} 
F_B(f_{!}X,Y)\iso f_*F_A(X,f^*Y),
\end{gather}
where $X$ is an ex-$G$-space over $A$ and $Y$ and $Z$ are ex-$G$-spaces over $B$.
\end{prop}

\begin{proof}
The isomorphism $f^*S^0_B\iso S^0_A$ is evident since $f^*(B\times K)\iso A\times K$ for based $G$-spaces $K$.  The isomorphism (\ref{oneo}) is obtained by passage to quotients from the evident homeomorphism
$$(Y\times_B A)\times _A (Z\times_B A) \iso (Y\times_BZ)\times_B A $$
As explained in \cite[\S\S2, 3]{FHM}, the isomorphism (\ref{oneo}) is equivalent to the
isomorphism (\ref{twoo}), and it determines natural maps from left to right in (\ref{three0}),
(\ref{four0}), and (\ref{five0}) such that all three are isomorphisms if any one is. By a 
comparison of definitions, we see that the categorically defined map in (\ref{three0}),
which is denoted $\al$ in \cite[3.3]{FHM}, coincides in the present situation with the map,
also denoted $\al$, on \cite[p.\,167]{BB2}. As explained on \cite[p.\,178]{BB2}, in the point-set
topological framework that we have adopted, that map $\al$ is a homeomorphism.\end{proof}

\begin{rem} 
Only the very last statement refers to topology.
The categorically defined map $\al$ should quite generally
be an isomorphism in analogous contexts, but we have not pursued
this question in detail. An alternative self-contained proof of the previous proposition is given in \myref{fgext} below by using \myref{Mackey0} to prove (\ref{four0}) instead of (\ref{three0}). In that argument, the only non-formal ingredient is the fact that the functor $D\times_B(-)$ commutes with pushouts.
\end{rem}

We shall later need a purely categorical coherence observation
about the categorically defined map $\al$ of (\ref{three0}).  In fact, it will play a key role in the proof of the fiberwise duality theorem of \S15.1. It is 
convenient to insert it here. 

\begin{rem}\mylabel{coherence}
Let $T\colon \sB\rtarr \sA$ be a symmetric monoidal functor. 
We are thinking of $T$ as, for example, a base change functor $f^*$. 
The map 
$$\alpha\colon TF_\sB(X,Y)\rtarr F_\sA(TX,TY)$$ 
is defined to be the adjoint of 
\[\xymatrix{TF_\sB(X,Y)\sma_\sA TX \iso T(F_\sB(X,Y)\sma_\sB X)\ar[r]^-{T\text{ev}} & TY.}\]
The dual of $X$ is $D_{\sB}X = F_{\sB}(X,S_{\sB})$, where
$S_{\sB}$ is the unit of $\sB$. Taking $Y=S_{\sB}$, the definition of 
$\al$ implies that the top triangle commutes in the diagram
\[\xymatrix{
TD_\sB X\sma_\sA TX \ar[r]^\iso\ar[d]_{\alpha\sma_\sA \text{id}} & 
T(D_\sB X\sma_\sB X) \ar[r]^-{T\text{ev}} & TS_\sB\ar[d]^\iso\\
F_\sA(TX,TS_\sB) \sma_\sA TX \ar[urr]_{\text{ev}}\ar[r]_-\iso & D_\sA f^*X\sma_\sA f^*X \ar[r]_-{\text{ev}} & S_\sA.}\]
The bottom triangle is a naturality diagram. The outer 
rectangle is \cite[3.7]{FHM}, but its commutativity in general was not 
observed there. However, it was observed in \cite[3.8]{FHM} that its commutativity implies the commutativity of the diagram
\[\xymatrix{TD_\sB X \sma_\sA TY \ar[d]_{\alpha\sma_\sA TY}\ar[r]^\iso 
& T(D_\sB X\sma_\sB Y) \ar[r]^-{T\nu} 
& TF_\sB(X,Y)\ar[d]^\alpha\\
D_\sA TX \sma_\sA TY \ar[rr]_-\nu && F_\sA(TX,TY),}\]
where $\nu\colon D_{\sB}X\sma_{\sB}Y\rtarr F_{\sB}(X,Y)$ is the adjoint of
\[\xymatrix@1{
D_{\sB}X\sma_{\sB} Y\sma_{\sB} X \iso D_{\sB}X\sma_{\sB} X\sma_{\sB} Y \ar[r]^-{\text{ev}\sma\text{id}}
& S_{\sB}\sma_{\sB} Y\iso Y.}\]
\end{rem}

In other contexts, the analogue of (\ref{four0}) is called the ``projection formula'',\index{projection formula} and we shall also use that term. The following base change commutation relations with respect to pullbacks are also familiar from other contexts. We state the result for spaces but, apart from use of the fact that the functor $D\times_B (-)$ commutes with pushouts, 
the proof is formal.

\begin{prop}\mylabel{Mackey0} 
Suppose given a pullback diagram of base spaces 
$$\xymatrix{
C \ar[r]^-{g} \ar[d]_{i} & D \ar[d]^{j} \\
A \ar[r]_{f} & B.}$$
Then there are natural isomorphisms of functors 
\begin{equation}\label{bases0}
j^*f_{!} \iso g_{!}i^*, \qquad f^*j_* \iso i_*g^*, \qquad f^*j_{!}\iso i_!g^*, \qquad j^*f_*\iso g_*i^*.
\end{equation}
\end{prop}\begin{proof} The first isomorphism is one of left adjoints, and the second is the corresponding
``conjugate'' isomorphism of right adjoints. Similarly for the third and fourth isomorphisms.
By symmetry, it suffices to prove the first isomorphism. 
The functor $j^* = D\times_B(-)$ commutes with pushouts. For a space $X$ over $A$ regarded by
composition with $f$ as a space over $B$, $C\times_AX\iso D\times_BX$.  This gives
\[ j^*f_{!} X = D\times _B(B\cup_A X) \iso D\cup_C(C\times_A X) = g_{!}i^*X.\qedhere\]
\end{proof}

\section{Change of group and restriction to fibers}

This section begins the study of equivariant phenomena that have no 
non\-equivariant counterparts.  In particular, using a conceptual reinterpretation of the adjoints of the fiber functors $(-)_b$, 
we relate restriction to subgroups to restriction to fibers. Recall that subgroups of $G$ are understood to be closed and fix an inclusion $\io\colon H\subset G$ throughout this section. Parametrized theory gives a convenient way of studying restriction along $\io$ without changing the ambient group from $G$ to $H$.

\begin{prop}\mylabel{homog} The category $G\sK_{G/H}$ of ex-$G$-spaces over 
$G/H$ is equivalent to the category $H\sK_*$ of based $H$-spaces.
\end{prop}\begin{proof}
The equivalence sends an ex-$G$-space $(Y,p,s)$ over $G/H$
to the $H$-space $p^{-1}(eH)$ with basepoint the $H$-fixed point $s(eH)$. 
Its inverse sends a based $H$-space $X$ to the induced 
$G$-space $G\times_H X$, with the evident structure maps.\end{proof}

More formally, recall that there are ``induction''\index{induction} and ``coinduction''\index{coinduction} functors 
$\io_!$ and $\io_*$ from $H$-spaces to $G$-spaces that are left and right adjoint
to the forgetful functor $\io^*$ that sends a $G$-space $Y$ to $Y$ regarded 
as an $H$-space. Explicitly, for an $H$-space $X$,
\begin{equation}\label{GH1}
\io_! X = G\times_H X \qquad \text{and}\qquad \io_*X = \text{Map}_H(G,X).
\end{equation}
The latter is the space of maps of (left) $H$-spaces, with (left)
action of $G$ induced by the right action of $G$ on itself. Similarly,
when $X$ is a based $H$-space, we have the based analogues
\begin{equation}\label{GH2}
\io_! X = G_+\sma_H X \qquad \text{and}\qquad \io_*X = F_H(G,X).
\end{equation}
With this notation, some familiar natural isomorphisms take the forms
\begin{gather}\label{GH3}
\io_!(\io^*Y\times X)\iso Y\times \io_!X \qquad \text{and} \qquad
\io_*\text{Map}(\io^*Y,X)\iso \text{Map}(Y,\io_*X)
\intertext{and, in the based case,}
\label{GH4}
\io_!(\io^*Y\sma X)\iso Y\sma \io_!X \qquad\text{and} \qquad
\io_* F(\io^*Y,X)\iso F(Y,\io_*X).
\end{gather}
By the uniqueness of adjoints, or inspection of definitions, we see that these
familiar change of groups functors are change of base functors along $r\colon G/H\rtarr *$. 

\begin{cor} The change of group and change of base functors associated to
$\io$ and $r$ agree under the equivalence of categories between $H\sK_*$ and $G\sK_{G/H}$: 
$$\io^*\iso r^*, \qquad  \io_!\iso r_{!},  \qquad \text{and} \qquad \io_*\iso r_*.$$
\end{cor}

We can generalize this equivalence of categories, using the following definitions.  We have a forgetful functor 
$\io^*\colon G\sK_{B} \rtarr H\sK_{\io^*B}$.\noteindex{io@$\iota^*$}
It doesn't have an obvious left or right adjoint, but we have obvious 
analogues of induction and coinduction that involve changes of base spaces.
The first will lead to a description of $\io^*$ as a base change functor and
thus as a functor with a left and right adjoint.

\begin{defn}\mylabel{changes0} 
Let $A$ be an $H$-space and $X$ be an $H$-space over $A$. 
Define $\io_!\colon H\sK_A\rtarr G\sK_{\io_! A}$\noteindex{iol@$\iota_{"!}$} by letting $\io_!X$ be the $G$-space
$G\times_H X$ over $\io_!A = G\times_H A$. Define $\io_*\colon H\sK_A\rtarr G\sK_{\io_* A}$\noteindex{ior@$\iota_*$} 
by letting $\io_*X$ be the $G$-space $\text{Map}_H(G,X)$ over $\io_*A = \text{Map}_H(G,A)$.
\end{defn}

For an $H$-space $A$ and a $G$-space $B$, let
\begin{equation}\mylabel{munuin}
\mu\colon G\times_H \io^*B = \io_!\io^*B\rtarr B
\ \, \text{and}\ \, \nu\colon A\rtarr \io^*\io_!A = \io^*(G\times_H A)
\end{equation}
be the counit and unit of the $(\io_!,\io^*)$ adjunction.  
The following result says that ex-$H$-spaces over an $H$-space 
$A$ are equivalent to ex-$G$-spaces over the $G$-space $\io_!A$. 

\begin{prop}\mylabel{ishriek}
The functor $\io_!\colon H\sK_A\rtarr G\sK_{\io_!A}$ is a 
closed symmetric monoidal equivalence of categories with inverse 
the composite 
$$G\sK_{\io_!A}\stackrel{\io^*}{\rtarr} 
H\sK_{\io^*\io_!A}\stackrel{{\nu}^*}{\rtarr} H\sK_A.$$
\end{prop}

Applied to $A =\io^*B$, this equivalence leads to the promised description 
of $\io^*\colon G\sK_B\rtarr H\sK_{\io^*B}$ as a base change functor.

\begin{prop}\mylabel{ishriekb}
The functor $\io^*\colon G\sK_B\rtarr H\sK_{\io^*B}$ is the composite
\[\xymatrix@1
{G\sK_B \ar[r]^-{\mu^*} & G\sK_{\io_!\io^*B} \iso H\sK_{\io^*B}\\}\]
\end{prop}

Change of base and change of groups are related by various further
consistency relations. The following result gives two of them.

\begin{prop}\mylabel{changerel}
Let $f\colon A\rtarr \io^*B$ be a map of $H$-spaces and $\tilde{f}\colon
\io_! A\rtarr B$ be its adjoint map of $G$-spaces.  Then the following
diagrams commute up to natural isomorphism.
\[\xymatrix{
G\sK_{\io_!A} \ar[r]^-{\tilde{f}_!} & G\sK_B \\
H\sK_{A} \ar[r]_-{f_!}\ar[u]^{\io_!}  &
H\sK_{\io^*B}\ar[u]_{\mu_!\com\io_!}}
\quad \ \
\xymatrix{
G\sK_B \ar[r]^-{\tilde{f}^*} \ar[d]_{\io^*}
& G\sK_{\io_!A} \ar[d]^{{\nu}^*\com\io^*}\\
H\sK_{\io^*B} \ar[r]_-{f^*} & H\sK_{A}}\]
\end{prop}

\begin{proof} Since $\tilde{f} = \mu\com \io_!f$, we have
\[\tilde{f}_!\com \iota_!
\cong (\mu\com \iota_!f)_!\com \iota_!
\cong \mu_!\com (\iota_!f)_!\com \iota_!
\cong \mu_!\com \iota_!\com f_!,\]
where the last isomorphism holds because $G\times_H(-)$ commutes with
pushouts. Since $f = \io^*\tilde{f}\com \nu$, we have
\[ f^*\com \iota^* \cong (\io^*\tilde{f}\com \nu)^*\com \io^*
\cong \nu^*\com (\io^*\tilde{f})^*\com \io^*
\cong \nu^*\com \io^*\com \tilde{f}^*, \]
where the last isomorphism holds because pulling the $G$ action back
to an $H$-action commutes with pullbacks.
\end{proof}

The reader may find it illuminating to work out these isomorphisms in the 
context of \myref{homog}.  That result leads to the promised conceptual
reinterpretation of  \myref{Fibad}.

\begin{exmp}\mylabel{Johann0}
For $b\in B$, we also write $b\colon *\rtarr B$ for the map that sends $*$ to $b$, and we write $\tilde{b}\colon G/G_b\rtarr B$ for the induced inclusion of orbits. Thus $b$ is a $G_b$-map and $\tilde{b}$ is a $G$-map. Under the equivalence 
$G\sK_{G/G_b}\iso G_b\sK_*$ of \myref{homog}, ${\tilde{b}}^*$ may be 
interpreted as the fiber functor $G\sK_B\rtarr G_b\sK_*$ that sends $X$ to 
$X_b$, ${\tilde{b}}_{!}$ may be interpreted as the left adjoint of  \myref{Fibad} that sends $K$ to $K^b$, and ${\tilde{b}}_*$ specifies a right adjoint to the fiber functor, 
which we denote by ${^{b}}K$. With these notations, the isomorphisms of 
\myref{Wirth0} specialize to the following natural isomorphisms, 
where $Y$ and $Z$ are in $G\sK_B$ and $K$ is in $G_b\sK_*$. 
\begin{gather*}
(Y\sma_B Z)_b\iso Y_b\sma Z_b,\\[1ex]
F_B(Y,\, ^bK) \iso\, {^{b}}F(Y_b,K),\\[1ex]
F_B(Y,Z)_b\iso F(Y_b,Z_b),\\[1ex]
(Y_b\sma K)^b\iso Y\sma_B K^b,\\[1ex]
F_B(K^b,Y)\iso \, {^{b}}F(K,Y_b).
\end{gather*}
\end{exmp}

\begin{exmp}\mylabel{Johann1}
Several earlier results come together in the following 
situation. Let $f\colon A\rtarr B$ be a $G$-map. For 
$b\in B$, let $b\colon \{b\}\rtarr B$ and $i_b\colon A_b\rtarr A$ 
denote the evident inclusions of $G_b$-spaces. We have the following compatible 
pullback squares, the first of $G_b$-spaces and the second of $G$-spaces.
$$\xymatrix{
A_b \ar[r]^-{f_b} \ar[d]_{i_b} & \{b\} \ar[d]^{b} \\
A \ar[r]_{f} & B}
\qquad 
\xymatrix{
G\times_{G_b} A_b \ar[r]^-{G\times_{G_b} f_b} \ar[d]_{\tilde{\imath}_b} & G/G_b \ar[d]^{\tilde{b}} \\
A \ar[r]_-{f} & B}$$
Applying \myref{Mackey0} to the right-hand square and 
interpreting the conclusion in terms of fibers by \myref{changes0},
we obtain canonical isomorphisms of $G_b$-spaces
$$(f_!X)_b \iso {f_{b}}_!i_b^*X \qquad\text{and}\qquad(f_*X)_b \iso {f_{b}}_*i_b^*X,$$
where $X$ is an ex-$G$-space over $A$, regarded on the right-hand sides as an ex-$G_b$-space over $A$ by pullback along $\io\colon G_b\rtarr G$.
\end{exmp}

\section{Normal subgroups and quotient groups}

Observe that any homomorphism $\tha\colon G\rtarr G'$ factors as the composite of a quotient homomorphism $\epz$, an isomorphism, and an inclusion $\io$. We
studied change of groups along inclusions in the previous section.  Here we consider a quotient homomorphism $\epsilon\colon G\rtarr J$ of $G$ by a normal subgroup $N$. We still have a restriction functor
\[\epsilon^*\colon J\sK_A\rtarr G\sK_{\epsilon^*A},\]
and we also have the functors
\[(-)/N\colon G\sK_B\rtarr J\sK_{B/N}\qquad\text{and}\qquad
(-)^N\colon G\sK_B\rtarr J\sK_{B^N}\]
obtained by passing to orbits over $N$ and to $N$-fixed points. When $B$ is a point, these last two functors are left and right adjoint to $\epsilon^*$, but in general change of base must enter in order to obtain such adjunctions. The following observation follows directly by inspection of the definitions.

\begin{prop}\mylabel{factor0} 
Let $j\colon B^N\rtarr B$ be the inclusion and $p\colon B\rtarr B/N$ be the quotient map.  Then the following factorization diagrams commute.
\[\xymatrix{
G\sK_B \ar[d]_{p_!}\ar[r]^{(-)/N} & J\sK_{B/N} \\
G\sK_{B/N} \ar[ur]_{(-)/N} }
\qquad\text{and}\qquad
\xymatrix{G\sK_B \ar[d]_{j^*}\ar[r]^{(-)^N} & J\sK_{B^N} \\
G\sK_{B^N} \ar[ur]_{(-)^N}}\]
It follows that $((-)/N,p^*\epsilon^*)$ and $(j_!\epsilon^*,(-)^N)$ are adjoint pairs.
\end{prop}

We have the following analogue of \myref{changerel}. 

\begin{prop}\mylabel{fixorbbase}
Let $f\colon A\rtarr B$ be a map of $G$-spaces. Then the following diagrams commute up to natural isomorphisms.
{\[\xymatrix{
G\sK_A \ar[r]^-{f_!} \ar[d]_{(-)/N} & G\sK_B \ar[d]^{(-)/N}\\
J\sK_{A/N} \ar[r]_-{(f/N)_!}  & J\sK_{B/N}}
\quad
\xymatrix{
G\sK_B \ar[r]^-{f^*} \ar[d]_{(-)^N} & G\sK_A \ar[d]^{(-)^N}\\
J\sK_{B^N} \ar[r]_-{(f^N)^*} & J\sK_{A^N}}
\quad
\xymatrix{
G\sK_A \ar[r]^-{f_!} \ar[d]_{(-)^N} & G\sK_B \ar[d]^{(-)^N}\\
J\sK_{A^N} \ar[r]_-{(f^N)_!}  & J\sK_{B^N}}\]}
\end{prop}
\begin{proof}
For ex-$G$-spaces $X$ over $A$ and $Y$ over $B$, these isomorphisms
are given by the homeomorphisms $$(X\cup_A B)/N\iso X/N\cup_{A/N}B/N,$$ 
$$(Y\times_B A)^N\iso Y^N\times_{B^N}A^N,$$
and
$$(X\cup_A B)^N\iso X^N\cup_{A^N}B^N.$$
As a quibble, the last requires $A\rtarr X$ to be a closed inclusion, but this will hold for the sections of compactly generated ex-$G$-spaces over $A$ by \myref{coflemma}(i).
\end{proof}

Specializing to $N$-free $G$-spaces, we obtain a factorization result that is analogous to those in \myref{factor0}, but is less obvious. It is a precursor 
of the Adams isomorphism, which we will derive in \S16.4.

\begin{prop}\mylabel{ouch0}
Let $E$ be an $N$-free $G$-space, let $B = E/N$, and let $p\colon E\rtarr B$ be the quotient map. Then the diagram 
\[\xymatrix{
G\sK_E \ar[r]^{(-)/N} \ar[d]_{p_*} & J\sK_B\\
G\sK_B \ar[ur]_{(-)^N} }\]
commutes up to natural isomorphism. Therefore the left adjoint $(-)/N$ of the functor $p^*\epz^*$ is also its right adjoint.
\end{prop}

\begin{proof}
Let $X$ be an ex-$G$-space over $E$ with projection $q$. Comparing the 
pullbacks that are used to define the functors $p_*$ and $\text{Map}_B$ in Definitions \ref{retract1} and \ref{MapB}, we find that $p_*X$ fits into a pullback diagram
\[\xymatrix{
p_*X \ar[r] \ar[d] & \text{Map}(E,\tilde{X})\ar[d]^{\tilde{q}}\\
B\ar[r]_-{\nu} & \text{Map}(E,\tilde{E}).}\]
Here $\nu(b)$, $b = Ne$, corresponds as in \myref{partialtilde} to the inclusion of the closed subset $Ne$ in $E$.
Passing to $N$-fixed points, we see that it suffices to prove that the 
following commutative diagram is a pullback.
\[\xymatrix{
X/N \ar[r]^-{\mu} \ar[d]_{q/N} & \text{Map}_N(E,\tilde{X})\ar[d]^{\tilde{q}}\\
E/N = B\ar[r]_-{\nu} & \text{Map}_N(E, \tilde{E})}\]
Here $\mu$ is induced from the adjoint of the map $X\times E\rtarr \tilde{X}$ that sends $(x,e)$ to $nx$ if $e = nq(x)$ and sends $(x,e)$ to $\om$ otherwise.
With this description, $\mu$ is well-defined since $E$ is $N$-free. It suffices
to give a continuous inverse to the induced map 
$$\phi\colon X/N \rtarr \text{Map}_N(E,\tilde{X})\times_{\text{Map}_N(E,\tilde{E})} E/N.$$
If $(f,Ne)$ is a point in the pullback, then $f$ corresponds to a map 
$Ne\rtarr X$, and $\ph^{-1}(f,Ne) = Nf(e)$ in $X/N$. For continuity, note that
$\ph^{-1}$ is obtained from the evaluation map 
$\text{Map}(E,\tilde{X})\times E\rtarr \tilde{X}$ by passage to subquotient spaces.
\end{proof}

\begin{rem}\mylabel{iotaalt}
This leads to a useful alternative description of the functor $\iota_!\colon H\sK_A\rtarr G\sK_{\io_!A}$, where $A$ is an $H$-space and 
$\io_!A = G\times_H A$. We have the projection $\pi\colon G\times A \rtarr A$ 
of $(G\times H)$-spaces, where the 
$G\times H$ actions on the source and target are given by
$$(g,h)(g',a) = (gg'h^{-1},ha) \qquad \text{and}\qquad (g,h)a = ha.$$
Consider ex-$H$-spaces $X$ over $A$ as $(G\times H)$-spaces with $G$ 
acting trivially and let $\epsilon\colon G\times H\rtarr H$ be the 
projection. We see from the definition that $\iota_!X=({\pi}^*\epz^*X)/H$. 
Since $G\times A$ is an $H$-free $(G\times H)$-space, we conclude from the 
previous result that $\iota_!X \iso (p_*{\pi}^*\epz^*X)^H$, where 
$p\colon G\times A \rtarr G\times_H A=\iota_!A$ is the quotient map.
\end{rem}

\section{The closed symmetric monoidal category of retracts}

Let $\sB$ be a topologically bicomplete full subcategory of a topologically bicomplete
category $\sC$. We are thinking of $\sU\subset \sK$ or $G\sU\subset G\sK$. 
We have the category of retracts $\sC_{\sB}$.\noteindex{CB@$\sC_{\sB}$}\index{category!of retracts} The objects of $\sC_{\sB}$ are the 
retractions $B\stackrel{s}{\rtarr}X\stackrel{p}{\rtarr}B$
with $B\in\sB$ and $X\in\sC$, abbreviated $(X,p,s)$ or just $X$. The morphisms of 
$\sC_{\sB}$ are the evident commutative diagrams. When $\sB=\sC$, this is just a diagram 
category for the evident two object domain category. 

The importance of the category $\sC_{\sB}$ is apparent from its role in 
\myref{retract1}:  focus on this category is natural when we consider 
base change functors. In our examples, $\sB$ and $\sC$ are enriched and topologically bicomplete over the appropriate category of spaces, $\sU$ for $\sB$ and $\sK$ for $\sC$.  For a space $K\in \sK$, 
the tensors $-\times K$ and cotensors $\text{Map}(K,-)$ applied to retractions give retractions, and we have the adjunction homeomorphisms
\begin{equation}\label{Bsilly}
\sC_{\sB}(X\times K,Y)\iso \sK(K,\sC_{\sB}(X,Y))\iso \sC_{\sB}(X,\text{Map}(K,Y)).
\end{equation} 

The category $G\sK_{G\sU}$ is closed symmetric 
monoidal under an external smash product functor, denoted $X\barwedge Y$,\noteindex{XYe@$X\barwedge Y$}\index{smash product!external} and an external \index{mapping space!external}
function ex-space functor, denoted $\bar{F}(Y,Z)$.\noteindex{FYZe@$\bar{F}(Y,Z)$} If $X$, $Y$, and $Z$ are ex-spaces over $A$, $B$, and $A\times B$, respectively, then $X\barwedge Y$ is an 
ex-space over $A\times B$ and $\bar{F}(Y,Z)$ is an ex-space over $A$. We have
\begin{equation}\label{exad}
G\sK_{A\times B}(X\barwedge Y, Z)\iso G\sK_A(X, \bar{F}(Y,Z)),
\end{equation}
which gives the required adjunction in $G\sK_{G\sU}$. It specializes to parts of (\ref{tensored2}) when
$A$ or $B$ is a point. The ex-space $X\barwedge Y$ is the evident fiberwise smash product, with  
$(X\barwedge Y)_{(a,b)} = X_a\sma Y_b$. The fiber $\bar{F}(Y,Z)_a$ is
$F_B(Y,Z_a)$, where $Z_a$ is the ex-space over $B$ whose fiber  $Z_{a,b}$ over $b$ 
is the inverse image of $(a,b)$ under the projection $Z\rtarr A\times B$.
Rather than describe the topology of the ex-space $\bar{F}(Y,Z)$ directly, 
we give alternative descriptions of $X\barwedge Y$ and $\bar{F}(Y,Z)$
in terms of internal smash products and internal function ex-spaces. 
Let $\pi_A$ and $\pi_B$ be the projections of $A\times B$ on $A$ and $B$ and 
observe that $\pi_A^*X \iso X\times B$ and $\pi_B^*Y \iso A\times Y$.
If one likes, the following results can be taken as a definition of the external operations
and a characterization of the internal operations, or vice versa.

\begin{lem}\mylabel{exin} The external smash product and function ex-space functors 
are determined by the internal functors via natural isomorphisms
$$ X\barwedge Y \iso \pi_A^*X \sma_{A\times B} \pi_B^*Y
\qquad \text{and} \qquad \bar{F}(Y,Z) \iso {\pi_{A}}_{*}F_{A\times B}(\pi_{B}^{*}Y,Z),$$
where $X$, $Y$, and $Z$ are ex-spaces over $A$, $B$, and $A\times B$, respectively.
\end{lem}

With these isomorphisms taken as definitions, the adjunction (\ref{exad}) follows 
from the adjunctions $(\pi^*_A,{\pi_{A}}_{*})$, $(\pi^*_B,{\pi_{B}}_{*})$, and 
$(\sma_{A\times B}, F_{A\times B})$.

\begin{lem}\mylabel{internalize} The internal smash product and function ex-space functors 
are determined by the external functors via natural isomorphisms
$$ X\sma_B Y \iso \DE^*(X\barwedge Y) \qquad \text{and} \qquad 
F_B(X,Y) \iso \bar{F}(X,\DE_*Y),$$
where $X$ and $Y$ are ex-spaces over $B$ and $\DE\colon B\rtarr B\times B$ is the diagonal map.
\end{lem}

With these isomorphisms taken as definitions, the adjunction $(\sma_{B}, F_{B})$  
follows from the adjunctions $(\DE^*,\DE_*)$ and (\ref{exad}).  Since $\DE^*$ is 
symmetric monoidal and the composite of either projection 
$\pi_i\colon B\times B\rtarr B$ with $\DE$ is the identity map of $B$, we see that, 
if we have constructed both internal and external smash products, then they 
must be
related by natural isomorphisms as in Lemmas \ref{exin} and \ref{internalize}. 

\begin{rem} 
The first referee suggests that we point out another consistency check.  
The fiber $(\DE_*Y)_{(b,c)}$ is a point if $b\neq c$ and is
$Y_b$ if $b=c$. Therefore the fiber over $b$ of the restriction 
$(\DE_*Y)_{b}$ of $\DE_* Y$ to $\{b\}\times B$ is $Y_b\cup (B-\{b\})$,
suitably topologized, and 
$$\bar{F}(X,\DE_*Y)_b = F_B(X,(\DE_*Y)_{b})_b \iso F(X_b,Y_b) = F_B(X,Y)_b.$$
\end{rem}

\begin{rem}\mylabel{fgext} 
The description of the internal smash product in terms
of the external smash product sheds light on the basic compatibility
isomorphisms (\ref{oneo}) and (\ref{four0}). For maps $f\colon A\rtarr B$
and $g\colon A'\rtarr B'$ and for ex-spaces $X$ over $B$ and $Y$ over
$B'$, it is easily checked that 
\begin{equation}\label{fgext1}
 f^*Y\barwedge g^*Z\iso (f\times g)^*(Y\barwedge Z).
\end{equation}
Similarly, for ex-spaces $W$ over $A$ and $X$ over $A'$,
\begin{equation}\label{fgext2}
f_!W\barwedge g_!X\iso (f\times g)_!(W\barwedge X).
\end{equation}
Now take $A=A'$, $B=B'$ and $f=g$. For ex-spaces $Y$ and $Z$ over $B$,
$$f^*(Y\sma_B Z)\iso f^*\DE_B^*(Y\barwedge Z) \iso (\DE_B\com f)^*(Y\barwedge Z).$$
On the other hand, using (\ref{fgext1}),
$$f^*Y\sma_A f^*Z\iso \DE_A^*(f\times f)^*(Y\barwedge Z) 
\iso ((f\times f)\com\DE_A)^*(Y\barwedge Z).$$
The right sides are the same since $\DE_B\com f = (f\times f)\com\DE_A$.
Similarly, 
$$f_!(f^*Y\sma_A X)\iso f_!\DE_A^*(f\times\text{id})^*(Y\barwedge X)
\iso f_!((f\times\text{id})\com \DE_A)^*(Y\barwedge X),$$
while
$$Y\sma_B f_!X \iso \DE_B^*(\text{id}\times f)_!(Y\barwedge X).$$
Since the diagram
$$\xymatrix{
A\ar[d]_f \ar[r]^-{\DE_A} & A\times A \ar[r]^-{f\times \text{id}} 
& B\times A \ar[d]^{\text{id}\times f}\\
B \ar[rr]_{\DE_B} & & B\times B\\}$$
is a pullback, the right sides are isomorphic by \myref{Mackey0}.
\end{rem}

It is illuminating conceptually to go further and consider group
actions from an external point of view. For groups $H$ and $G$,
an $H$-space $A$, and a $G$-space $B$, we have an evident 
external smash product  
\begin{equation}\label{GBexex}
\barwedge\colon H{\sK}_A \times G{\sK}_B \to (H\times G)\mathcal{K}_{A\times B}.
\end{equation}
For an ex-$H$-space $X$ over $A$ and an ex-$G$-space $Y$ over $B$, 
$X\barwedge Y$ is just the internal smash product over the 
$(H\times G)$-space $A\times B$
of $\pi_H^*\pi_A^* X$ and $\pi_G^*\pi_B^* Y$, where the $\pi's$
are the projections from $H\times G$ and $A\times B$ to their
coordinates.  It is easily seen that this definition leads to 
another $(\barwedge,\bar{F})$ adjunction. 

When $H = G$, the diagonal 
$\DE\colon G\rtarr G\times G$ is a 
closed inclusion since $G$ is compactly generated. We can pull back
along $\DE$, and then our earlier external smash product $X\barwedge Y$ 
over the $G$-space $\DE^*(A\times B)$ is given in terms of (\ref{GBexex})
as the pullback $\DE^*(X\barwedge Y)$.  Note that, by \myref{ishriekb},
$\DE^*$ here can be viewed as a base change functor.

\chapter{Proper actions, equivariant bundles and fibrations}

\section*{Introduction}

Much of the work in equivariant homotopy theory has focused on compact Lie groups. However, as was already observed by Palais \cite{Palais}, many results can be generalized to arbitrary Lie groups provided that one restricts to proper actions.  These are well-behaved actions whose isotropy groups are compact, and all actions by compact Lie groups are proper.  The classical definition of a Lie group \cite[p. 129]{Chev} includes all discrete groups (even though they need not be second countable) and, for discrete groups, the proper actions are the properly discontinuous ones.

In the parametrized world, the homotopy theory is captured on fibers.  When we restrict to proper actions on base spaces, the fibers have actions by the compact isotropy groups of the base space. So even though our primary interest is still in compact Lie groups of equivariance, proper actions on the base space provide the right natural level of generality. We set the stage for such a theory in this chapter by generalizing various classical results about equivariant bundles and fibrations to a setting focused on proper actions by Lie groups. The reader interested primarily in the nonequivariant theory should skip this chapter since only some very standard material in it is relevant nonequivariantly.

In \S3.1, we recall some basic results about proper actions of locally compact groups. We use this discussion to generalize some results about equivariant bundles in \S3.2. We generalize Waner's equivariant versions of Milnor's results on spaces of the homotopy types of CW complexes in \S3.3.  In \S3.4, we recall and generalize classical theorems of Dold and Stasheff about Hurewicz fibrations. We also recall an important but little known result of Steinberger and West that relates Serre and Hurewicz fibrations. We recall the definition of equivariant quasifibrations in \S3.5.

\section{Proper actions of locally compact groups}

We recall relevant definitions and basic results about proper actions in this section. For appropriate generality and technical convenience, we let $G$ be a locally compact topological group whose underlying topological space is compactly generated. Local compactness means that the identity element, hence any point, has a compact neighborhood. We see from \myref{HauswHaus} that $G$ is Hausdorff and, since all compact subsets are closed, it follows that each neighborhood of any point contains a compact neighborhood.

\begin{rem}
We comment on the assumptions we make for $G$. If $G$ is any topological group whose underlying space is in $\sK$, then an action of $G$ on $X$ in $\sK$ may not come from an action in ${\sT}op$. The point is that the product $G\times X$ in $\sK$ is defined by applying the $k$-ification functor to the product $G\times_c X$ in ${\sT}op$, and not every action $G\times X\rtarr X$ need be continuous when viewed as a function $G\times_c X\rtarr X$. However, when $G$ is locally compact, $G\times_c X$ is already in $\sK$ by \myref{HauswHaus}, and $k$-ification is not needed. There is then no ambiguity about what we mean by a $G$-space, and we need not worry about refining the topology on products with $G$.

Another reason for restricting to locally compact groups is that many useful properties of proper actions only hold in that case.  In the literature, such results are usually derived for actions on Hausdorff spaces, but we shall see that weak Hausdorff generally suffices.
\end{rem}

We begin with some standard equivariant terminology. 

\begin{defn} Let $X$ be a $G$-space and let $H\subset G$.
\begin{enumerate}[(i)]
\item An \emph{$H$-tube} $U$ in $X$ is an open $G$-invariant 
subset of $X$ together with a $G$-map $\pi\colon U\rtarr G/H$. 
If $x\in U$ and $H = G_x$, then $U$ is a \emph{tube around $x$}. 
A tube is \emph{contractible} if $\pi$ is a $G$-homotopy equivalence.
\item An \emph{$H$-slice} $S$ in $X$ is an $H$-invariant subset such 
that the canonical $G$-map $G\times_H S\rtarr GS\subset X$ is an embedding onto 
an open subset. Then $GS$ is an $H$-tube with $S=\pi^{-1}(eH)$. Conversely, if $(U,\pi)$ is an $H$-tube in $X$, then $S=\pi^{-1}(eH)$ is an $H$-slice and $U=GS$.
On isotropy subgroups, we then have $G_y = H_y \subset H$ for all 
$y\in S$, but equality need not hold. If $x\in S$ and $H = G_x$, then 
$S$ is a \emph{slice through $x$}. 
\item We say that $X$ \emph{has enough slices} 
if every point $x\in X$ is contained in an $H$-slice for some \emph{compact} subgroup $H$. This implies that every point $x$ has compact isotropy group,
but in general it does not imply that there must be a slice \emph{through} 
every point $x$.
\item A \emph{$G$-numerable cover of $X$} is a cover $\{U_j\}$ by tubes 
such that there exists a locally finite partition of unity by $G$-maps $\lambda_j\colon X\rtarr[0,1]$ with support $U_j$. 
\end{enumerate}
\end{defn}

The following is the equivariant generalization of \cite[6.7]{Dold0}.

\begin{prop}\mylabel{doldcover}
Any $G$-CW complex admits a $G$-numerable 
cover by contractible tubes.
\end{prop}

\begin{proof}
The proof given by Dold \cite{Dold0} in the nonequivariant case goes through with only a minor change in the initial construction, which we sketch. From there, the technical details are unchanged. Let $X^n$ be the $n$-th skeletal filtration of a $G$-CW complex $X$. Let $\dot{X}^n$ denote the subspace obtained by deleting the centers $G/H\times 0$ of all $n$-cells in $X^n$ and let $r_n\colon \dot{X}^n\rtarr X^{n-1}$ denote the obvious retract. Starting
from the interior $e_n=G/H\times (D^n-S^{n-1})$ of an $n$-cell $c_n$, define $V_n^m$ inductively for $m\geq n$ by setting $V_n^n=e_n$ and $V_n^{m+1}=r_{m+1}^{-1}(V_n^m)$. Then the union $V_n^\infty=\bigcup_{m\geq n}V^m_n$ is a contractible tube, where the projection to $G/H$ is induced by 
the projection of $e_n$ to $G/H\times 0$.
\end{proof}

We now give the definition of a proper group action in $\sK$. We shall see
that the definition could equivalently be made in $\sU$. For further details, but in ${\sT}op$, see for example \cite{Bour, tomDieck}.  Recall that a continuous map is {\em proper} if it is a closed map with compact fibers.

\begin{defn} A $G$-space $X$ in $G\sK$ is \emph{proper} (or \emph{$G$-proper}) 
if the map
\[\theta\colon G\times X\rtarr X\times X\]
specified by $\theta(g,x) = (x,gx)$ is proper.
\end{defn}

We warn the reader that the definition is not quite the standard one. We are working in the category $\sK$, and the product $X\times X$ on the right hand side is the $k$-space obtained by $k$-ifying the standard product topology on $X\times_c X$. In ${\sT}op$ there are various other notions of a proper group action; see \cite{Biller} for a careful discussion. They all agree for actions of locally compact groups on completely regular spaces. If $X$ is proper, then the isotropy groups $G_x$ are compact since they are the fibers 
$\theta^{-1}(x,x)$. Moreover, since points are closed subsets of $G$, the diagonal $\Delta_X =\theta(\{e\}\times X)$ must be a closed subset of $X\times X$ and thus $X$ must be weak Hausdorff. This means that proper $G$-spaces must be in $\sU$.  Since $G$ is locally compact, we have the following useful characterizations.

\begin{prop}\mylabel{properchar}
For a $G$-space $X$ in $G\sK$ the following are equivalent.
\begin{enumerate}[(i)]
\item The action of $G$ on $X$ is proper.
\item The isotropy groups $G_x$ are compact and for all $(x,y )\in X\times X$ 
and all neighborhood $U$ of $\theta^{-1}(x,y)$ in $G\times X$, there is a neighborhood $V$ of $(x,y)$ in $X\times X$ such that $\theta^{-1}(V)\subset U$.
\item The isotropy groups $G_x$ are compact and for all $(x,y)\in X\times X$ and all neighborhoods $U$ of $\{ g \mid gx=y \}$ in $G$, there is a neighborhood $V$ of $(x,y)$ such that 
\[\{g\in G\mid \text{$ga=b$ for some $(a,b)\in V$}\}\subset U.\]
\item The space $X$ is weak Hausdorff and every point $(x,y)\in X\times X$ has a neighborhood $V$ such that 
\[\{g\in G\mid \text{$ga=b$ for some $(a,b)\in V$}\}\]
has compact closure in $G$.
\end{enumerate}
\end{prop}

\begin{proof}
This holds by essentially the same proof as \cite[1.6(b)]{Biller}. One must only keep in mind that we are now working in $\sK$ rather than in ${\sT}op$ and adjust the argument accordingly. 
\end{proof}

\begin{cor} If $G$ is discrete, then a $G$-space $X$ is proper if and 
only if any point $(x,y)\in X\times X$ has a neighborhood $V$ such that
\[\{g\in G\mid \text{$ga=b$ for some $(a,b)\in V$}\}\]
is finite.
\end{cor}

\begin{cor} 
If $G$ is compact, then any $G$-space in $G\sU$ is proper.
\end{cor}

\begin{rem}\mylabel{Hausprop}
There is an alternative description of the set displayed in \myref{properchar}
that may clarify the characterization.  Define
\[\phi\colon G\times X\times X\rtarr X\times X\]
by $\phi(g,x,y)=(gx,y)$.  For $V\subset X\times X$, let $\phi_V$ be the restriction of $\ph$ to $G\times V$ and let $\pi\colon G\times V\rtarr G$ 
be the projection, which is an open map since $G\times V$ has the
product topology.  Then the displayed set is $\pi\phi_V^{-1}(\Delta_X)$.
If $X\times X = X\times_c X$, then the condition in \myref{properchar} is equivalent to the more familiar one that any two points $x$ and $y$ 
in $X$ have neighborhoods $V_x$ and $V_y$ such that
\[\{g\in G\mid gV_x\cap V_y \neq \emptyset\}\]
has compact closure in $G$.
\end{rem}

\begin{prop}\mylabel{propproper}
Proper actions satisfy the following closure properties. 
\begin{enumerate}[(i)]
\item The restriction of a proper action to a closed subgroup is proper.
\item An invariant subspace of a proper $G$-space is also proper.
\item Products of proper $G$-spaces are proper.
\item If $X$ is a proper Hausdorff $G$-space in $G\sK$ and $C$ is a compact Hausdorff $G$-space, then the $G$-space $\text{Map}(C,X)$ is proper.
\item An $H$-space $S$ is $H$-proper if and only if $G\times_H S$ is 
$G$-proper. 
\end{enumerate}
\end{prop}

\begin{proof}
The first three are standard and elementary; see for example \cite[I.5.10]{tomDieck}. The 
fifth is \cite[2.3]{Biller}. We prove (iv). We must show that the map
\[\theta\colon G\times \text{Map}(C,X)\rtarr \text{Map}(C,X)\times \text{Map}(C,X)\]
is proper, which amounts to showing that it is closed and that the isotropy groups $G_f$ are compact for $f\in\text{Map}(C,X)$. For the latter, let $\{g_i\}$ be a net in $G_f$ and fix $c\in C$. Note that $f(g_ic)=g_if(c)$.
Since $C$ is compact, we can assume by passing to a subnet that $\{g_ic\}$ converges to some $\bar{c}\in C$. Let $V$ be a neighborhood of $(f(c),f(\bar{c}))$ such that
\[B=\{g\in G\mid \text{$ga=b$ for some $(a,b)\in V$}\}\]
has compact closure. Since $C$ is compact, $C\times C\times \text{Map}(C,X)$ has the usual product topology. 
Since the map
\[C\times C\times \text{Map}(C,X)\rtarr X\times X\]
that sends $(c,d,f)$ to $(f(c),f(d))$ is continuous and the net $\{c,g_ic,f\}$ converges to $(c,\bar{c},f)$, the net $\{(f(c),f(g_ic))\}=\{(f(c),g_if(c))\}$ must converge to $(f(c),f(\bar{c}))$. It follows that a subnet of $\{g_i\}$ 
lies in $B$ and therefore has a converging sub-subnet.

To show that $\theta$ is closed, let $A$ be a closed subset of $G\times \text{Map}(C,X)$ and let $\{(f_i,g_if_i)\}$ be a net in $\theta(A)$ that converges to $(f, F)$. We must show that $(f, F)$ is in $\theta(A)$. For 
$c\in C$, the net $\{g_i^{-1}c\}$ has a subnet that converges to some $\bar{c}$, by the compactness of $C$, so we may as well assume that the original net 
converges to $\bar{c}$. Let $V$ be a neighborhood of $(f(\bar{c}),F(c))$ 
such that
\[B'=\{g\in G\mid \text{$ga=b$ for some $(a,b)\in V$}\}\]
has compact closure. By continuity and the compactness of $C$, there is a compact neighborhood $K_1\times K_2$ of $(\bar{c},c)$ that $(f, F)$ maps into $V$. Since $\{(f_i, g_if_i)\}$ converges to $(f, F)$, there is an $h$ such that $(f_i, g_if_i)(K_1\times K_2)\subset V$ for $i\geq h$. It follows that there is a $k\geq h$ such that $(f_i(g_i^{-1}c), g_if_i(g_i^{-1}c))\in V$ for all $i\geq k$. Then the subnet $\{g_i\}_{i\geq k}$ is contained in $B'$ and therefore has a sub-subnet that converges to some $g\in G$.
We have now seen that our original net $\{(g_i, f_i)\}$ in $A$ has a subnet $\{(g_{i_j}, f_{i_j})\}$ that converges to $(g,f)$, and $(g,f)\in A$ since $A$ is closed.  By the continuity of $\theta$, $\{\theta(g_i,f_i)\}$ must converge to $(f, F)=\theta(g,f)\in\theta(A)$. In this last statement, we are using the
uniqueness of limits, which we ensure by requiring $X$ and $C$ to be Hausdorff.
\end{proof}

The following theorem of Palais \cite{Palais}, as generalized by
Biller \cite{Biller}, is fundamental. Those sources work in ${\sT}op$, 
but the arguments work just as well in $\sU$.

\begin{thm}[Palais]\mylabel{Palais} Let $X$ be a $G$-space in $G\sU$.
\begin{enumerate}[(i)]
\item If $X$ has enough slices, then it is proper. 
\item Conversely, if $X$ is completely regular and proper, then it has enough
slices.
\item If $G$ is a Lie group and $X$ is completely regular and proper, then there  is a slice through each point of $X$.
\end{enumerate}
\end{thm}
\begin{proof} Part (i) is given by \cite[2.4]{Biller}. Part (iii) is 
given by \cite[2.3.3]{Palais}. Part (ii) is deduced from part (iii) 
in \cite[2.5]{Biller}.
\end{proof}

\section{Proper actions and equivariant bundles}

We introduce here the equivariant bundles to which we will apply our basic foundational results in Part IV. As we explain, \myref{Palais} allows us to generalize some basic results about such bundles from actions of compact Lie groups to proper actions of Lie groups. 

Let $\PI$ be a normal subgroup of a Lie group $\GA$ such that 
$\GA/\PI = G$ and let $q\colon \GA\rtarr G$ be the quotient homomorphism.  By a \emph{principal $(\Pi;\Gamma)$-bundle} we mean the quotient map $p\colon P\rtarr P/\Pi$ where $P$ is a $\Pi$-free $\Gamma$-space such that $\Gamma$ acts properly on $P$. It follows that the induced $G$-action on $B=P/\Pi$ is proper. If $F$ is a $\Gamma$-space, then we have the associated $G$-map $E=P\times_\Pi F\rtarr P\times_\Pi *\cong P/\Pi$, which we say is a \emph{$\Gamma$-bundle with structure group $\Pi$ and fiber $F$}. For compact Lie groups, bundles of this general form are studied in \cite{LM}, which generalizes the study of the classical case $\GA = G\times \PI$ given in \cite{L}. A summary and further references are given in \cite[Chapter VII]{EHCT}.
We recall an observation about such bundles.

\begin{lem}\mylabel{rho}
For $b\in B$, the action of $\GA$ on $F$ induces an action of the isotropy group $G_b$ on the fiber $E_b$ through a homomorphism $\rh_b\colon G_b\rtarr \GA$ such that $q\com \rh_b$ is the inclusion $G_b\rtarr G$ and $E_b\iso \rh_b^*F$.
\end{lem}

\begin{proof}
Choose $z\in P$ such that $\pi(z)=b$.  The isotropy group $\GA_z$ intersects $\PI$ in the trivial group, and $q$ maps $\GA_z$ isomorphically onto $G_b$. Let $\rh_b$ be the composite of $q^{-1}\colon G_b\rtarr \GA_z$ and the inclusion $\GA_z\rtarr \GA$. Since the subspace $\{z\}\times F$ of $P\times F$ is $\GA_z$--invariant and maps homeomorphically onto $E_b$ on passage to orbits over $\PI$, the conclusion follows. Note that changing the choice of $z$ changes $\rh_b$ by conjugation by an element of $\PI$ and changes the identification of $E_b$ with $F$ correspondingly.
\end{proof}

Bundles should be locally trivial. When $P$ is completely regular, local triviality is a consequence of \myref{Palais}(iii), just as in the 
case when $\GA$ is a compact Lie group \cite[Lemma 3]{LM}, and this 
justifies our bundle-theoretic terminology. Note that if $P$ is 
completely regular, then so is $B=P/\Pi$.

\begin{lem}
A completely regular principal $(\Pi;\Gamma)$-bundle $P$ is locally trivial. 
That is, for each $b\in B$, there is a slice $S_b$ through $b$ and a homeomorphism 
\[\xymatrix{
\Gamma\times_{\Lambda} S_b\ar[r]^\cong\ar[d]_{q\times 1} & p^{-1}(GS_b)\ar[d]^{p} \\
G\times_{G_b} S_b \ar[r]^\cong & GS_b}\]
where $\Lambda\subset \Gamma$ only intersects $\Pi$ in the identity element and is mapped isomorphically to $G_b$ by $q$. The $\Lambda$-action on $S_b$ is given by pulling back the $G_b$-action along $q$.
\end{lem}

\section{Spaces of the homotopy types of $G$-CW complexes}

In this section, we recall and generalize the equivariant version of Milnor's results \cite{Milnor} about spaces of the homotopy types of CW complexes. For compact Lie groups, Waner formulated and proved such results 
in \cite[\S4]{Waner1}.  With a few observations, his proofs generalize to 
deal with proper actions by general Lie groups.  We first note the
following immediate consequence of \myref{doldcover} and \myref{Palais}.

\begin{thm} For any locally compact group $G$, a $G$-CW complex is proper if 
and only if it is constructed from cells of the form $G/K\times D^n$, where $K$ is compact.
\end{thm}

We also note the following recent ``triangulation theorem'' of Illman \cite[Theorem II]{Ill}. It is this result that led us to try to generalize some of our results from compact Lie groups to general Lie groups. 

\begin{thm}[Illman]\mylabel{Illman}
If $G$ is a Lie group that acts smoothly and properly on a smooth manifold $M$, then $M$ has a $G$-CW structure.
\end{thm}

Many of our applications of this result are based on the following observation.

\begin{lem}\mylabel{prodproper}
If $H$ and $K$ are closed subgroups of a topological group $G$ and $K$ is
compact, then the diagonal action of $G$ on $G/H\times G/K$ is proper. 
\end{lem}

\begin{proof}
The proof given in \cite[I.5.16]{tomDieck} that $G$ acts properly on $G/K$ generalizes directly. Set $X=G/H\times G$. Let $G$ act diagonally from the left
and let $K$ act on the second factor from the right.  Note that these actions commute. It suffices to show that $\theta\colon G\times X\rtarr X\times X$ is proper.  Indeed, consider the commutative square
\[\xymatrix{G\times X \ar[d]\ar[r]^-\theta & X\times X\ar[d]\\
G\times X/K \ar[r]^-{\bar{\theta}} & X/K\times X/K.}\]
The right vertical map is proper and the left vertical map is surjective.
Therefore, by \cite[VI.2.13]{tomDieck}, the bottom horizontal map
is proper if the top horizontal map is proper.  Since $X$ is a free $G$-space, $\theta$ is proper if and only if the image $\text{Im}(\tha)$ is a closed subspace of $X\times X$ and the map $\phi\colon \text{Im}(\tha)\rtarr G$ specified by $\phi(x,gx)=g$ is continuous. The diagonal subspace of $G/H\times G/H$ is closed, and its preimage under the map
$\zeta\colon X\times X\rtarr G/H\times G/H $
specified by 
$$\zeta((xH,y),(\bar{x}H,\bar{y}))=(\bar{y}y^{-1}xH,\bar{x}H)$$
is precisely $\text{Im}(\tha)$, which is therefore closed. The function 
$\phi$ is the restriction to $\text{Im}(\tha)$ of the continuous map
$\PH\colon X\times X\rtarr G$ specified by
$$\PH((xH,y),(\bar{x}H,\bar{y}))= \bar{y}y^{-1}$$
and is therefore continuous.
\end{proof}

We shall also make essential use of the following corollary of \myref{Illman}.

\begin{cor}\mylabel{Illcor} If $X$ is a proper $G$-CW complex, then, viewed as an $H$-space for any closed subgroup $H$ of $G$, $X$ has the structure of an $H$-cell complex.
\end{cor}
\begin{proof}
Each cell $G/K\times D^n$ has $K$ compact. Since $G$ acts smoothly and properly on the smooth manifold $G/K$, the closed subgroup $H$ also acts smoothly and properly. We use the resulting $H$-CW structure on all of the cells to obtain an $H$-cell structure. It is homotopy equivalent to an $H$-CW complex obtained by ``sliding down'' cells that are attached to higher dimensional ones, but we
shall not need to use that. 
\end{proof} 

\begin{thm}[Milnor, Waner]\mylabel{Milnor}
Let $G$ be a Lie group and $(X;X_i)$ be an $n$-ad 
of closed sub-$G$-spaces of a proper $G$-space $X$. If $(X;X_i)$ has 
the homotopy type of a $G$-CW $n$-ad and $(C;C_i)$ is an $n$-ad of 
compact $G$-spaces, then $(X;X_i)^{(C;C_i)}$ has the homotopy type 
of a $G$-CW $n$-ad.
\end{thm}

\begin{proof}
We only remark how the proof of Waner for the case of actions by a compact Lie group generalizes to the case of proper actions by a Lie group. Define a $G$-simplicial complex to be a $G$-CW complex such that $X/G$ with the induced cell structure is a simplicial complex. In \cite[\S5]{Waner1}, Waner proves that any $G$-CW complex is $G$-homotopy equivalent to a colimit of finite dimensional $G$-simplicial complexes and cellular inclusions and that a $G$-space dominated by a $G$-CW complex is $G$-homotopy equivalent to a $G$-CW complex. The arguments apply verbatim to any topological group $G$. 

The rest of the argument requires two key lemmas. In \cite[4.2]{Waner1}, Waner defines the notion of a $G$-equilocally convex, or $G$-ELC, $G$-space. The first lemma says that every finite dimensional $G$-simplicial complex is $G$-ELC. The essential starting point is that orbits are $G$-ELC, the proof of which uses the Lie group structure just as in \cite[p.358]{Waner1} in the compact case.  From there, Waner's proof \cite[\S6]{Waner1} goes through unchanged. The second says that any completely regular, $G$-paracompact, $G$-ELC, proper $G$-space is dominated by a $G$-CW complex.  When $G$ is compact Lie, this is proven in \cite[\S7]{Waner1}.  However, the hypothesis on $G$ is only used to guarantee the existence of enough slices, hence the proof holds without change for 
proper actions of Lie groups, indeed of locally compact groups.

The rest of the proof goes as in \cite[Theorem 3]{Milnor}.
One only needs to make two small additional observations. First, if a $G$-simplicial complex $K$ has the homotopy type of a proper $G$-space $X$, 
then it is proper. This holds since if $f\colon K\rtarr X$ is a homotopy equivalence, then $G_k\subset G_{f(k)}$ is compact. Second, for an
$n$-ad $(K;K_i)$ of $G$-simplicial complexes and a compact $n$-ad $(C;C_i)$,
$(X;X_i)^{(C;C_i)}$ is proper since it is a subspace of the proper $G$-space $X^C$; see (i) and (iv) of \myref{propproper}. Since it is also completely regular, $G$-paracompact, and $G$-ELC, it is dominated by a $G$-CW complex, 
and the result follows from the steps above. 
\end{proof}

\section{Some classical theorems about fibrations}

A basic principle of parametrized homotopy theory is that homotopical information is given on fibers.  We recall two relevant classical theorems about Hurewicz fibrations and a comparison theorem relating Serre and Hurewicz fibrations.  We begin with Dold's theorem \cite[6.3]{Dold0}. The nonequivariant proof in \cite[2.6]{May} is generalized to the equivariant case in Waner \cite[1.11]{Waner2}. Waner assumes throughout \cite{Waner2} that $G$ is a compact Lie group, but that assumption is not used in the cited proof. 

\begin{thm}[Dold]\mylabel{Dold}
Let $G$ be any topological group and let 
$B$ be a $G$-space that has a $G$-numerable cover by contractible tubes. Let  $X\rtarr B$ and $Y\rtarr B$ be Hurewicz fibrations.  Then a map $X\rtarr Y$ 
over $B$ is a fiberwise $G$-homotopy equivalence if and only if each fiber restriction $X_b\rtarr Y_b$ is a $G_b$-homotopy equivalence.
\end{thm}

We next recall and generalize a classical result that relates the homotopy types of fibers to the homotopy types of total spaces. Nonequivariantly, it is due to Stasheff \cite{Stash} and, with a much simpler proof, Sch\"on \cite{Schon}. The generalization to the equivariant case, for compact Lie groups, is given by Waner \cite[6.1]{Waner2}. With Theorems \ref{Dold}, \ref{Milnor} and \ref{Illman} in place, Sch{\"o}n's argument generalizes directly to give the following version. Since the result plays an important role in our work and the argument is so pretty, we can't resist repeating it in full.

\begin{thm}[Stasheff, Sch\"on]\mylabel{ss}
Let $G$ be a Lie group and $B$ be a proper $G$-space that has the homotopy type of a $G$-CW complex. Let $p\colon X\rtarr B$ be a Hurewicz fibration. Then $X$ has the homotopy type of a $G$-CW complex if and only if each fiber $X_b$ has the homotopy type of a $G_b$-CW complex.
\end{thm}
\begin{proof} 
First assume that $X$ has the homotopy type of a $G$-CW complex. For $b\in B$, let $\iota\colon G_b\rtarr G$ be the inclusion and consider the $G_b$-map $\iota^*p\colon \iota^*X\rtarr \iota^*B$ of $G_b$-spaces. It is still a 
Hurewicz fibration, as we see by using the left adjoint $G\times_{G_b}(-)$ of $\iota^*$. By \myref{Illcor}, $\iota^*X$ and $\iota^*B$ have the homotopy types of $G_b$-CW complexes. Factor $\iota^*p$ through the inclusion into its mapping cylinder $i\colon \iota^*X\rtarr M\iota^*p$. Since $G_b$ is compact, it follows from \myref{Milnor} that the homotopy fiber $F_bi=(M\iota^*p;\{b\},\iota^*X)^{(I;0,1)}$ has the homotopy type of a 
$G_b$-CW complex. Since $F_b i$ is homotopy equivalent to $F_b \iota^*p$, by the gluing lemma, and $F_b \iota^*p$ is homotopy equivalent to the fiber $X_b$, this proves the forward implication.

For the converse, assume that each fiber $X_b$ has the homotopy type of a 
$G_b$-CW complex. Let $\gamma\colon \Gamma X\rtarr X$ be a $G$-CW approximation of $X$. The mapping path fibration of $\gamma$ gives us a factorization of 
$\ga$ as the composite of a $G$-homotopy equivalence 
$\nu\colon \GA X\rtarr N\ga$ and a Hurewicz fibration 
$q\colon N\ga\rtarr X$. We may view $q$ as a map of fibrations over $B$.
\[\xymatrix{N\ga \ar[rr]^-{q} \ar[dr]_{p\com q}& & X \ar[dl]^p\\
& B & \\}\]
The fibers of $p\com q$ have the homotopy types of $G_b$-CW complexes by the first part of the proof, since $\Gamma X$ is a $G$-CW complex, and the fibers
of $p$ have the homotopy types of $G_b$-CW complexes by hypothesis.  Comparison of the long exact sequences associated to $p\com q$ and $p$ gives that $q$ restricts to a $G_b$-homotopy equivalence on each fiber.  Noting that we can pull back a numerable cover by contractible tubes along a homotopy equivalence $B\rtarr B'$, where $B'$ is a $G$-CW complex, it follows from \myref{Dold} that $q$ is a homotopy equivalence.
\end{proof}

Although it no longer plays a role in our theory, the following little known result played a central role in our thinking. It shows that the dichotomy between Serre and Hurewicz fibrations diminishes greatly over CW base spaces. 
It is due to Steinberger and West \cite{SW}, with a correction by Cauty \cite{Cauty}.

\begin{thm}[Steinberger and West; Cauty]\mylabel{qh}
A Serre fibration whose base and total spaces are CW complexes 
is a Hurewicz fibration.
\end{thm}

We believe that this remains true equivariantly for compact Lie groups, 
and it certainly remains true for finite groups. Before we understood 
the limitations of the $q$-model structure, we planned to use this result 
to relate our model theoretic homotopy category of ex-spaces over a
CW complex $B$ to a classical homotopy category defined in terms of 
Hurewicz fibrations and thereby overcome the problems illustrated in \myref{noway}.  Such a comparison is still central to our theory, 
and it is this result that convinced us that such a comparison must hold.

\section{Quasifibrations}

For later reference, we recall the definition of quasifibrations.
Here $G$ can be any topological group. 

\begin{defn}\mylabel{quasifib}
A map $p\colon E\rtarr Y$ in $\sK$ is a \emph{quasifibration} if the map of pairs $p\colon (E,E_y)\rtarr (Y,y)$ is a weak equivalence for all $y$ in $Y$.
A map $p\colon E\rtarr Y$ in $\sK/B$ or $\sK_B$ is a \emph{quasifibration} 
if it is a quasifibration on total spaces. A $G$-map $p\colon E\rtarr Y$ 
is a quasifibration if each of its fixed point maps $p^H\colon E^H\rtarr Y^H$
is a nonequivariant quasifibration.
\end{defn}

The condition that $p\colon (E,E_y)\rtarr (Y,y)$ is a weak equivalence means that for all $e\in E_y$ the following two conditions hold.
\begin{enumerate}[(i)]
\item $p_*\colon \pi_n(E,E_y,e)\rtarr \pi_n(Y,y)$ is an isomorphism for all $n\geq 1$.
\item For any $x\in E$, $p(x)$ is in the path component of $y$ precisely when the path component of $x$ in $E$ intersects $E_y$. In other words, the sequence
$$\pi_0(E_y,e)\rtarr \pi_0(E,e)\rtarr \pi_0(Y,y)$$
of pointed sets is exact.
\end{enumerate}

\begin{warn}
In contrast to the usual treatments in the literature, we do not require
$p$ to be surjective and therefore $\pi_0(E,e)\rtarr \pi_0(Y,y)$ need not
be surjective.  Hurewicz and, more generally, Serre fibrations are examples 
of quasifibrations, and they are not always surjective, as the trivial example $\{0\} \rtarr \{0,1\}$ illustrates.  Model categorically, one point is that
the initial map $\emptyset \rtarr Y$ is always a Serre fibration since the
empty lifting problem always has a solution. 
\end{warn}

The definition of a quasifibration is arranged so that the long exact sequence of homotopy groups associated to the triple $(E,E_y,e)$ is isomorphic to a long exact sequence
\[\cdots \rtarr \pi_{n+1}(Y,y)\rtarr \pi_n(E_y,e) \rtarr \pi_n(E,e)\rtarr \pi_n(Y,y)\rtarr\cdots\rtarr \pi_0(Y,y).\]

\part{Model categories and parametrized spaces}

\chapter*{Introduction}

In Part III, we shall develop foundations for para\-me\-trized equivariant stable homotopy theory. In making that theory rigorous, it became apparent to us that substantial foundational work was already needed on the level of ex-spaces. That work is of considerable interest for its own sake, and it involves general points about the use of model categories that should be of independent interest. Therefore, rather than rush through the space level theory as just a precursor of the spectrum level theory, we have separated it out in this more leisurely and discursive exposition. 

In Chapter 4, which is entirely independent of our parametrized theory, we give general model theoretic background, philosophy, and results. In contrast to the simplicial world, we often have both a classical $h$-type and a derived $q$-type model structure in topologically enriched categories, with respective weak equivalences the homotopy equivalences and the weak homotopy equivalences. We describe what is involved in verifying the model axioms for these two types of model structures.

In Chapter 5, we describe how the parametrized world fits into this general framework. There are several different $h$-type model structures on our categories of parametrized $G$-spaces, with different homotopy equivalences based on different choices of cylinders. These mesh in unexpected ways. Understanding of this particular case leads us to a conceptual axiomatic description of how the classical $h$-type homotopy theory and the $q$-type model structure must be related in order to be able to do homotopy theory satisfactorily in a topologically enriched category.

In Chapter 6, we work nonequivariantly and develop our preferred ``$q$-type'' model category structure, the ``$qf$-model structure'', on the categories 
$\sK/B$ and $\sK_B$. This chapter is taken directly from the second author's thesis \cite{Sig}. 
 
In Chapter 7, we give the equivariant generalization of the $qf$-model
structure and begin the study of the resulting homotopy categories by
discussing those adjunctions that are given by Quillen pairs. There
is another new twist here in that we need to use many Quillen
equivalent $qf$-type model structures. In fact, this 
is already needed nonequivariantly in the study of base change
along bundles $f\colon A\rtarr B$.

In Chapter 8, we discuss ex-fibrations and an ex-fibrant approximation
functor that better serves our purposes than model theoretic fibrant approximation in studying those adjunctions that are not given by
Quillen pairs. In Chapter 9, we describe our parametrized homotopy 
categories in terms of classical homotopy categories of ex-fibrations 
and use this description to resolve the issues concerning base change 
functors and smash products that are discussed in the Prologue.

\chapter{Topologically bicomplete model categories}

\section*{Introduction} 

In \S4.1, we describe a general philosophy about the role of different model structures on a given category $\sC$. It is natural and important in many contexts, and it helps to clarify our thinking about topological categories of parametrized objects. In particular, we advertise a remarkable unpublished insight of Mike Cole.  It is a pleasure to thank him for keeping us informed of his ideas.  We describe how a classical ``$h$-type'' model structure and a suitably related Quillen ``$q$-type'' model structure, can be mixed together to give an ``$m$-type'' model structure such that the $m$-equivalences are the $q$-equivalences and the $m$-fibrations are the $h$-fibrations. This is a completely general phenomenon, not restricted to topological contexts.

In \S\S4.2 and 4.3, we describe classical structure that is present in any topologically bicomplete category $\sC$. Here we follow up a very illuminating paper of Schw\"anzl and Vogt \cite{SVogt}. There are two classes of (Hurewicz) $h$-fibrations and two classes of $h$-cofibrations, ordinary and strong.  Taking weak equivalences to be homotopy equivalences, the ordinary $h$-fibrations pair with the strong $h$-cofibrations and the strong $h$-fibrations pair with the ordinary $h$-cofibrations to give two interrelated model like structures.  For each choice, all of the axioms for a proper topological model category are satisfied except for the factorization axioms, which hold in a weakened form. To prove that $\sC$ is a model category, it suffices to prove one of the factorization axioms since the other will follow.  Again, the theory can easily be adapted to other contexts than our topological one.

We signal an ambiguity of nomenclature.  In the model category literature, the term ``simplicial model structure'' is clear and unambiguous, since there is only one model structure on simplicial sets in common use.  In the topological context, we understand ``topological model structures'' to refer implicitly to the $h$-model structure on spaces for model structures of $h$-type and to the $q$-model structure on spaces for model structures of $q$-type.  The meaning should always be clear from context.

In \S4.4, we give another insight of Cole, which gains power from the work of Schw\"anzl and Vogt.  Cole provides a simple hypothesis that implies the missing factorization axioms for an $h$-model structure of either type on a topologically bicomplete category $\sC$.  When we restrict to compactly generated spaces, the hypothesis applies to give an $h$-model structure on $\sU$. In $\sK$, this seems to fail, and we give a streamlined version of Str{\o}m's original proof \cite{Strom}, together with his proof that the strong $h$-cofibrations in $\sK$ are just the closed ordinary $h$-cofibrations. This works in exactly the same way for the categories $G\sK$ and $G\sU$, where $G$ is any (compactly generated) topological group.

In \S4.5, we describe how to construct compactly generated $q$-type model structures, giving a slight variant of standard treatments. In particular, $G\sK$ and $G\sU$ have the usual $q$-model structures in which the $q$-equivalences are the weak equivalences and the $q$-fibrations are the Serre fibrations. Again, $G$ can be any topological group. However we only know
that the model structure is $G$-topological when $G$ is a compact Lie group.

\section{Model theoretic philosophy: $h$, $q$, and $m$-model
 structures}\label{Sphil}

The point of model categories is to systematize ``homotopy theory''. The homotopy theory present in many categories of interest comes in two flavors. There is a ``classical'' homotopy theory based on homotopy equivalences, and there is a more fundamental ``derived'' homotopy theory based on a weaker notion of equivalence than that of homotopy equivalence. This dichotomy pervades the applications, regardless of field. It is perhaps well understood that both homotopy theories can be expressed in terms of model structures on the underlying category, but this aspect of the classical homotopy theory has usually been ignored in the model theoretical literature, a tradition that goes back to Quillen's original paper \cite{Q}. The ``classical'' model structure on spaces was introduced by Str{\o}m \cite{Strom}, well after Quillen's paper, and the ``classical'' model structure on chain complexes was only introduced explicitly quite recently, by  Cole \cite{Cole2} and Schw\"anzl and Vogt \cite{SVogt}.

Perhaps for this historical reason, it may not be widely understood that these two model structures can profitably be used in tandem, with the $h$-model structure used as a tool for proving things about the $q$-model structure. This point of view is implicit in \cite{EKMM, MM, MMSS}, and a variant of this point of view will be essential to our work. In the cited papers, the terms ``$q$-fibration'' and ``$q$-cofibration'' were used for the fibrations and cofibrations in the Quillen model structures, and the term ``$h$-cofibration'' was used for the classical notion of a Hurewicz cofibration specified in terms of the homotopy extension property (HEP). The corresponding notion of an ``$h$-fibration'' defined in terms of the covering homotopy property (CHP) is fortuitously appropriate\footnote{However, the notation conflicts with the notation often used for Dold's notion of a weak or  ``halb''-fibration. We shall make no use of that notion, despite its real importance in the theory of fibrations. We do not know whether or not it has a model theoretic role to play.}.  Just as the ``$q$'' is meant to suggest Quillen, the ``$h$'' is meant to suggest Hurewicz, as well as homotopy. It is logical to follow this idea further (as was not done in \cite{EKMM, MM, MMSS}) by writing $q$-fibrant, $q$-cofibrant, $h$-fibrant, and $h$-cofibrant for clarity.  Following this still further, we should also write ``$h$-equivalence'' for homotopy equivalence and ``$q$-equivalence'' for (Quillen) weak equivalence. The relations among these notions are as follows in all of the relevant categories $\sC$:

\vspace{1mm}
\begin{center}
\begin{tabular}{|r c l|} \hline
$h$-equivalence & $\Longrightarrow$ &  $q$-equivalence\\
$h$-cofibration & $\Longleftarrow$ &  $q$-cofibration\\
$h$-cofibrant &  $\Longleftarrow$ &   $q$-cofibrant\\
$h$-fibration &  $\Longrightarrow$ & $q$-fibration\\
$h$-fibrant & $\Longrightarrow$ &  $q$-fibrant\\\hline
\end{tabular}
\end{center}
\vspace{1mm}

Therefore, {\em the identity functor is the right adjoint of a
Quillen adjoint pair from $\sC$ with its $h$-model structure
to $\sC$ with its $q$-model structure.}\/  It follows
that we have an adjoint pair relating the classical homotopy
category, $h\sC$\noteindex{hC@$h\sC$} say, to the derived homotopy category
$q\sC = \Ho \sC$\noteindex{HoC@$\Ho \sC$}. This formulation
packages standard information. For example, the Whitehead theorem that a
weak equivalence between cell complexes is a homotopy equivalence,
or its analogue that a quasi-isomorphism between projective complexes
is a homotopy equivalence, is a formal consequence of this adjunction
between homotopy categories.

Recently, Cole \cite{Cole4} discovered a profound new way of thinking
about the dichotomy between the kinds of model structures that we have been
discussing. He proved the following formal model theoretic result.

\begin{thm}[Cole]\mylabel{Colemix}
\index{cofibration!mixed}\index{fibration!mixed}\index{weak equivalence!mixed}
Let  $(\sW_h,\sF\!ib_h,\sC\!of_h)$ and $(\sW_q,\sF\!ib_q,\sC\!of_q)$ be two
mo\-del structures on the same category $\sC$. Suppose that $\sW_h\subset
\sW_q$ and $\sF\!ib_h\subset \sF\!ib_q$. Then there is a \emph{mixed model
structure}\index{model structure!mixed} $(\sW_q,\sF\!ib_h,\sC\!of_m)$ on
$\sC$. The mixed cofibrations $\sC\!of_m$ are the maps in $\sC\!of_h$ that
factor as the composite of a map in $\sW_h$ and a map in $\sC\!of_q$. An
object is $m$-cofibrant if and only it is $h$-cofibrant and of the
$h$-homotopy type of a $q$-cofibrant object. If the $h$ and $q$-model
structures are left or right proper, then so is the $m$-model structure.
\end{thm}

By duality, the analogue with the inclusion $\sF\!ib_h\subset \sF\!ib_q$
replaced by an inclusion $\sC\!of_h\subset \sC\!of _q$ also holds. In the
category of spaces with the $h$ and $q$-model structures discussed above,
the theorem gives a mixed model structure whose $m$-cofibrant spaces are
the spaces of the homotopy types of CW-complexes. This $m$-model structure
combines weak equivalences with Hurewicz fibrations, and it might
conceivably turn out to be as important and convenient as the Quillen model
structure. It is startling that this model structure was not discovered
earlier.

The pragmatic point is two-fold. On the one-hand,
there are many basic results that apply to $h$-cofibrations and not just
$q$-cofibrations. Use of $h$-cofibrations limits the need for $q$-cofibrant
approximation and often clarifies proofs by focusing
attention on what is relevant.  Many examples appear in 
\cite{EKMM, MMSS, MM}, where properties of $h$-cofibrations serve as 
scaffolding in the proof that $q$-model structures are in fact model structures. We shall formalize and generalize this idea in the next chapter. 

On the other hand, there are many vital results that apply only to
$h$-fibrations (Hure\-wicz fibrations), not to $q$-fibrations (Serre
fibrations). For example, a local Hurewicz fibration is a Hurewicz
fibration, but that is not true for Serre fibrations. The mixed
model structure provides a natural framework in which to make use
of Hurewicz fibrations in conjunction with weak equivalences. While
we shall make no formal use of this model structure, it has provided
a helpful guide to our thinking. The philosophy here applies in
algebraic as well as topological contexts, but we shall focus on
the latter.

 \section{Strong Hurewicz cofibrations and fibrations}

Fix a topologically bicomplete category $\sC$ throughout this section
and the next.
With no further hypotheses on $\sC$, we show that it satisfies most of the
axioms for not one but two generally different proper topological $h$-type
model structures. We alert the reader to the fact that we are here using the
term ``$h$-model structure'' in a generic sense.  When we restrict
attention to parametrized spaces, we will use the term in a different
specific sense derived from the $h$-model structure on underlying total
spaces. The material of these
sections follows and extends material in Schw\"anzl and Vogt \cite{SVogt}.

We have cylinders $X\times I$\noteindex{XI@$X\times I$} and cocylinders
$\text{Map}(I,X)$.\noteindex{MapIX@$\text{Map}(I,X)$}  When $\sC$ is based,
we focus on the based cylinders $X\sma I_+$\noteindex{XIa@$X\sma I_+$}
and cocylinders $F(I_+,X)$.\noteindex{FIX@$F(I_+,X)$}  In either case,
these define equivalent notions of homotopy, which we shall sometimes
call \emph{$h$-homotopy}.  We will later use these and cognate notations,
but, for the moment, it is convenient to introduce the common notations
$\text{Cyl}(X)$\noteindex{CylX@$\text{Cyl}(X)$}\index{Cylinder object} and
$\text{Cocyl}(X)$\noteindex{CocylX@$\text{Cocyl}(X)$}\index{Cocylinder
object} for these objects.  There are obvious classes of maps that one
might hope would specify a model structure.

\begin{defn}\mylabel{hmodel}
Let $f$ be a map in $\sC$.
\begin{enumerate}[(i)]
\item $f$ is an \emph{$h$-equivalence}\index{equivalence!h@$h$-}
if it is a homotopy equivalence in $\sC$.
\item $f$ is a \emph{Hurewicz fibration},\index{Hurewicz fibration}
 abbreviated
\emph{$h$-fibration},\index{fibration!h@$h$-} if it
satisfies the CHP\index{CHP}\index{Covering homotopy property} in $\sC$,
 that is,
if it has the right lifting property (RLP) with respect
to the maps $i_0: X\rtarr \text{Cyl}(X)$ for $X\in\sC$.
\item $f$ is a \emph{Hurewicz cofibration},\index{Hurewicz cofibration}
 abbreviated
\emph{$h$-cofibration},\index{cofibration!h@$h$-}
if it satisfies the HEP\index{HEP}\index{Homotopy extension property} in
 $\sC$, that is, if it has
the left lifting property (LLP) with respect to the maps 
$p_0\colon \text{Cocyl}(X)\rtarr X$.
\end{enumerate}
\end{defn}

These sometimes do give a model structure, but then the $h$-cofibrations 
must be exactly the maps that satisfy the LLP with respect to the $h$-acyclic $h$-fibrations, and dually. In general, that does not hold.  We shall characterize the maps
in $\sC$ that do satisfy the LLP with respect to the $h$-acyclic
$h$-fibrations and, dually, the maps that satify the RLP with
respect to the $h$-acyclic $h$-fibrations. For this, we need the
following relative version of the above notions.

\begin{defn}\mylabel{barhmodel} We define strong Hurewicz fibrations and
cofibrations.
\begin{enumerate}[(i)]
\item A map $p\colon E\rtarr Y$ is a \emph{strong Hurewicz
 fibration},\index{Hurewicz fibration!strong} abbreviated
 \emph{$\bar{h}$-fibration},\index{fibration!h@$\bar{h}$-}
if it satisfies the 
\emph{relative CHP}\index{relative CHP}\index{CHP!relative} with 
respect to all $h$-cofibrations $i:A\rtarr X$, in the sense that a 
lift exists in any diagram
\[\xymatrix{A\ar[r]^i\ar[d]_{i_0} & X \ar[r] \ar[d] & E\ar[d]^p\\
\text{Cyl}(A)\ar[r]\ar[urr]|(.55)\hole & \text{Cyl}(X)\ar[r]\ar@{-->}[ur] &
 Y.}\]
 \item A map $i:A\to X$ is a \emph{strong Hurewicz cofibration},\index{Hurewicz cofibration!strong}
abbreviated \emph{$\bar{h}$-cofibration},\index{cofibration!h@$\bar{h}$-}
if it satisfies the \emph{relative HEP}\index{relative HEP}\index{HEP!relative} with respect to all
$h$-fibrations $p:E\to Y$, in the sense that a lift exists in any diagram
\[\xymatrix{A\ar[r]\ar[d]_i & \text{Cocyl}(E) \ar[r] \ar[d] &
 \text{Cocyl}(Y)\ar[d]^{p_0}\\
X\ar[r]\ar@{-->}[ur]\ar[urr]|(.44)\hole & E\ar[r]_p & Y.}\]
\end{enumerate}
\end{defn}

We recall the standard criteria for maps to be $h$-fibrations or $h$-cofibrations. Define the \emph{mapping cylinder}\index{mapping cylinder} $Mf$\noteindex{Mf@$Mf$} and \emph{mapping path fibration}\index{mapping path fibration} $Nf$\noteindex{Nf@$Nf$} by the usual pushout and pullback diagrams
$$\xymatrix{
X \ar[d]_{i_0} \ar[r]^-f & Y \ar[d] \\
\text{Cyl}(X) \ar[r] &  Mf}\\
\ \ \ \ \text{and} \ \ \ \
\xymatrix{Nf \ar[d] \ar[r] & \text{Cocyl}(Y) \ar[d]^{p_0} \\
X \ar[r]_{f} & Y.\\}$$

\begin{lem}\mylabel{babyish} Let $f$ be a map in $\sC$.
\begin{enumerate}[(i)]
\item $f$ is an $h$-fibration if and only if 
it has the RLP with respect to the map $i_0: Nf\rtarr \text{Cyl}(Nf)$.
\item $f$ is an $h$-cofibration if and only if it has
the LLP with respect to the map $p_0: \text{Cocyl}(Mf)\rtarr Mf$.
\end{enumerate}
\end{lem}

The $\bar{h}$-fibrations and $\bar{h}$-cofibrations admit similar characterizations.  These were taken as definitions in \cite[2.4]{SVogt}.

\begin{lem}\mylabel{h-char} (i) A map $p\colon E\rtarr Y$ is an $\bar{h}$-fibration if and only if it has the RLP with respect to the canonical map
$Mi\rtarr \text{Cyl}(X)$ for any $h$-cofibration $i:A\to X$;
this holds if and only if the canonical map $\text{Cocyl}(E)\to Np$
has the RLP with respect to all $h$-cofibrations.\\
(ii) A map $i:A\to X$ is an $\bar{h}$-cofibration if and only if 
it has the LLP with respect to the canonical map
$\text{Cocyl}(E)\rtarr Np$ for any $h$-fibration $p:E\to Y$;
this holds if and only if the canonical map $Mi\to \text{Cyl}(X)$
has the LLP with respect to all $h$-fibrations.
\end{lem}

Observe that the map $i_0:X\rtarr\text{Cyl}(X)$ is an $\bar{h}$-cofibration and the map $p_0:\text{Cocyl}(X)\rtarr X$ is an $\bar{h}$-fibration. Since the cylinder objects associated to initial objects are initial objects, $\bar{h}$-fibrations are in particular $h$-fibrations. Similarly, $\bar{h}$-cofibrations are $h$-cofibrations.  Observe too that every object is both $\bar{h}$-cofibrant and $\bar{h}$-fibrant, hence both $h$-cofibrant and $h$-fibrant.

We shall see in \S4.4 that these distinctions are necessary in $\sK$ but disappear in $\sU$, where the $h$ and $\bar{h}$ notions coincide. Even there, however, the conceptual distinction sheds light on classical arguments.  

The results of this section and the next are quite formal. Amusingly, the
main non-formal ingredient is just the use in the following proof of the 
standard fact that $\{0,1\}\to I$ has the LLP with respect to all 
$h$-acyclic $h$-fibrations.

\begin{lem}\mylabel{h-retract} Let $i\colon A\rtarr X$ and $p\colon E\rtarr
 B$ be maps
in $\sC$.
\begin{enumerate}[(i)]
\item If $i$ is an $h$-acyclic $h$-cofibration, then $i$ is the inclusion
of a strong deformation retraction $r\colon X\rtarr A$.
\item If $i$ is the inclusion of a strong deformation retraction $r:X\to A$, then
$i$ is a retract of $Mi\to \text{Cyl}(X)$.
\item If $p$ is an $h$-acyclic $h$-fibration, then $p$ is a strong
deformation retraction.
\item If $p$ is a strong deformation retraction, then $p$ is a retract of
$\text{Cocyl}(E)\rtarr Np$.
\end{enumerate}
\end{lem}
\begin{proof} The last two statements are dual to the first two. For (i),
since the $h$-equivalence $i$ is an $h$-cofibration, application of the HEP shows that $i$ has a homotopy inverse $r:X\to A$ such that $ri=\text{id}_A$.
Since $\{0,1\}\rtarr I$ has the LLP with respect to $h$-acyclic
$h$-fibrations, an adjunction argument shows that $p_{(0,1)}$ has the RLP
with respect to $h$-cofibrations. Thus a lift
exists in the diagram on the left, which means that $r$ is a strong
deformation retraction with inclusion $i$.
\[\begin{aligned}\xymatrix{
A\ar[d]_i\ar[r]^-{c} & \text{Cocyl}(A) \ar[r]^{\text{Cocyl}(i)} &
\text{Cocyl}(X)\ar[d]^{p_{(0,1)}}\\
X\ar@{-->}[urr]_\beta \ar[rr]_{(i\com r,\text{id}_X)} && X\times X}
\end{aligned}
\quad
\begin{aligned}\xymatrix{
& A \ar[d]_{i_0}\ar[r]^i & X \ar@{-}[d]_{i_0}\ar[dr]^r\\
A\ar[r]^-{i_1}\ar[d]_i &
\text{Cyl}(A)\ar[dr]_{\text{Cyl}(i)}\ar[rr]^(.35){\text{pr}} & \ar[d] &
A\ar[d]^i\\
X \ar[rr]_{i_1} && \text{Cyl}(X) \ar[r]_-\beta & X}
\end{aligned}\]
For (ii), we are given $\beta$ in the diagram on the left displaying $r$ as
a strong deformation retraction with inclusion $i$. Then the diagram on the
right commutes, where the composites displayed in the lower two rows are
identity maps. Using the universal property of $Mi$ to factor the crossing
arrows $i_0$ and $\text{pr}$ through $Mi$, we see that $i$ is a retract
of the canonical map $Mi\to \text{Cyl}(X)$.
\end{proof}

\section{Towards classical model structures in topological
categories}\label{sec:towardh}

We now have two candidates for a classical model structure on $\sC$ based
on the $h$-equivalences.  We can either take the $h$-fibrations and the
$\bar{h}$-cofibrations or the $h$-cofibrations and the $\bar{h}$-fibrations.
The following result shows that all of the axioms for a proper topological
model category are satisfied except that, in general, only a weakened form
of the factorization axioms holds. 

\begin{thm}\mylabel{h-structure}\index{model structure!generich@generic
$h$-structures} The following versions of the axioms for a proper 
topological model category hold.
\begin{enumerate}[(i)]
\item The classes of $h$-cofibrations, $\bar{h}$-cofibrations,
$h$-fibrations and $\bar{h}$-fibrations are closed 
under retracts.
\item Let $i$ be an $h$-cofibration and $p$ be an $h$-fibration. The pair
$(i,p)$ has the lifting property if $i$ is strong and $p$ is $h$-acyclic
or if $p$ is strong and $i$ is $h$-acyclic.
\item Any map $f:X\to Y$ factors as
\[\xymatrix{X \ar[r]^-i & Mf \ar[r]^-r & Y}\]
where $i$ is an $\bar{h}$-cofibration and $r$ has a section that is an
$h$-acyclic $\bar{h}$-cofibration and as
\[\xymatrix{X \ar[r]^-s & Nf \ar[r]^-p & Y}\]
where $p$ is an $\bar{h}$-fibration and $s$ has a retraction that is an
$h$-acyclic $\bar{h}$-fibration.
\item Let $i:A\to X$ be an $h$-cofibration and $p:E\to B$ be an
$h$-fibration, where $i$ or $p$ is strong. Then the map
\[\sC^\Box(i,p): \sC(X,E)\to \sC(A,E)\times_{\sC(A,B)}\sC(X,B)\]
induced by $i$ and $p$ is an $h$-fibration of spaces. It is $h$-acyclic if
$i$ or $p$ is acyclic and it is an $\bar{h}$-fibration if both $i$ and $p$
are strong.
\item The $h$-equivalences are preserved under pushouts along
$h$-cofibrations and pullbacks along $h$-fibrations.
\end{enumerate}
\end{thm}
\begin{proof}
Part (i) is clear since all classes are defined in terms of lifting
properties. Part (ii) follows directly from \myref{h-char} and
\myref{h-retract}.  The factorizations of part (iii)
are the standard ones. We consider the first.  The evident
section $j\colon Y\rtarr Mf$ is an $h$-acyclic $\bar{h}$-cofibration since
it is the pushout of one. Consider the lifting problem in the left diagram
below, in which the middle vertical composite is $i$. Here $p$ is an
$h$-acyclic $h$-fibration, and we choose a section $s$ of $p$.
\[\begin{aligned}\xymatrix{ & X\ar[d]_{i_1}\ar[dr]^{\alpha} \\
X\ar[r]^-{i_0}\ar[d]_f & \text{Cyl}(X) \ar@{-->}[r]^-{\la'}\ar[d]|(.5)\hole
 & E \ar[d]^p\\
Y \ar[r]_j\ar[urr]^(.3){s\beta j} & Mf\ar@{-->}[ur]_{\la}\ar[r]_{\beta} & B}
\end{aligned}
\qquad
\begin{aligned}\xymatrix{X\amalg X \ar[rr]^-{s\com \beta\com j\com
 f\amalg\alpha}\ar[d]_{i_{(0,1)}} && E\ar[d]^p\\
\text{Cyl}(X) \ar[r]\ar@{-->}[urr]^{\la'} & Mf\ar[r]_\beta & B}
\end{aligned}\]
We have a lift $\la'$ in the diagram on the right that makes the diagram on
the left commute, and the universal property of $Mf$ then gives us the lift
$\la$. Part (iv) is a consequence of the ``pairing theorem'' \cite{SVogt}, which we will state below.  Finally we prove the first half of (v). The second half
follows by duality. Assume that $i$ is an $h$-cofibration and
$f$ is an $h$-equivalence in the pushout diagram on the left.
\[\xymatrix{A\ar[r]^f \ar[d]_i & B\ar[d]^j\\
X\ar[r]_g & Y}\qquad\qquad
\xymatrix{B\ar[r]^s\ar[d]_{is} & A\ar[rr]^f\ar[d]\ar[ddr]^(.3)i &&
 B\ar[d]\ar[ddr]^j\\
X\ar[r]^{s'}\ar@{=}[drr] & P\ar[rr]|(.25)\hole^{f'}\ar[dr]|(.4)p && 
X \ar[dr]|(.4)q\\
&& X\ar[rr]^g & & Y}\]
We must prove that $g$ is an $h$-equivalence. By (ii), we can factor $f$ as
a composite of an $h$-acyclic $h$-cofibration and a map that has a section
which is an $h$-acyclic $h$-cofibration. Since pushouts preserve $h$-acyclic
$h$-cofibrations, we may assume that $f$ has a section $s\colon B\rtarr A$
that is an $h$-acyclic $h$-cofibration. We then obtain the diagram on the
right.  Its left back rectangle is a pushout, as is the outer back
rectangle, and therefore the right back rectangle is also a pushout.  This
implies that the bottom square is a pushout. The map $s'$ is an $h$-acyclic
$h$-cofibration since $s$ is one, and therefore $p$ is an $h$-equivalence.
The map $f'$ is also $h$-acyclic since it has the $h$-acyclic section $s'$.
Just as we could assume that $f$ has a section that is an $h$-acyclic
$h$-cofibration, we find that we may assume that $p$ has a section $t$ that
is an $h$-acyclic $h$-cofibration and is a map under $A$.  Chasing pushout
diagrams, we find that $g$ is a retract of $f'$ and is therefore an
$h$-equivalence.
\end{proof}

The following result is the pairing theorem\index{pairing theorem} of
\cite[2.7 and 3.6]{SVogt}. We shall not repeat the proof, which consists 
of careful but formal adjunction arguments. Its general statement is framed 
so as to apply to cartesian products in the unbased situation, smash products 
in the based situation, and tensors in either situation.

\begin{thm}[Schw\"anzl and Vogt]\mylabel{SVogt} Let $\sA$, $\sB$, and $\sC$
 be topologically
bicomplete categories and let
$$ T\colon\sA\times \sB\rtarr \sC,\ \ U\colon\sA^{\text{op}}\times \sC\rtarr
 \sB,\ \ \text{and} \ \ V\colon\sB^{\text{op}}\times \sC\rtarr \sA$$
be continuous functors that satisfy adjunctions
$$ \sC(T(A,B),C) \iso \sB(B, U(A,C))\iso \sA(A,V(B,C)).$$
Let $i\colon A\rtarr X$ be an $h$-cofibration in $\sA$,
$j\colon B\rtarr Y$ be an $h$-cofibration in $\sB$,
and $p\colon E\rtarr Z$ be an $h$-fibration in $\sC$.
\begin{enumerate}[(i)]
\item Assume that $i$ or $j$ is strong.  Then the map
$$T(A,Y)\cup_{T(A,B)}T(X,B)\rtarr T(X,Y)$$
induced by $i$ and $j$ is an $h$-cofibration in $\sC$. It is
$h$-acyclic if $i$ or $j$ is $h$-acyclic and it is strong if
both $i$ and $j$ are strong.
\item Assume that $j$ or $p$ is strong. Then the map
$$V(Y,E)\rtarr V(B,E)\times_{V(B,Z)}V(Y,Z)$$
induced by $j$ and $p$ is an $h$-fibration in $\sA$. It is
$h$-acyclic if $j$ or $p$ is $h$-acyclic and it is strong if
both $j$ and $p$
are strong.
\end{enumerate}
\end{thm}

As Schw\"anzl and Vogt observe, these results imply that the canonical map
$Mi\rtarr\text{Cyl}(X)$ is an $h$-acyclic $\bar{h}$-cofibration for any
$h$-cofibration $i:A\rtarr X$ and, dually, the canonical map
$\text{Cocyl}(X)\rtarr Np$ is an $h$-acyclic $\bar{h}$-fibration for any
$h$-fibration $p:E\rtarr Y$. Together with \myref{h-retract} and the
retract and factorization axioms of \myref{h-structure}, this implies that
all of the various classes of maps are characterized by the expected
lifting properties, just as if we had actual model categories.

\begin{prop}\mylabel{charstrong} The following characterizations hold.
\begin{enumerate}[(i)]
\item The $h$-fibrations are the maps that have the RLP
with respect to the $h$-acyclic $\bar{h}$-cofibrations and the
$h$-acyclic $\bar{h}$-cofibrations are the maps that have the LLP
with respect to the $h$-fibrations.
\item The $h$-cofibrations are the maps that have the LLP
with respect to the $h$-acyclic $\bar{h}$-fibrations and the
$h$-acyclic $\bar{h}$-fibrations are the maps that have the RLP
with respect to the $h$-cofibrations.
\item The $\bar{h}$-fibrations are the maps that have the RLP
with respect to the $h$-acyclic $h$-cofibrations and the
$h$-acyclic $h$-cofibrations are the maps that have the LLP
with respect to the $\bar{h}$-fibrations.
\item The $\bar{h}$-cofibrations are the maps that have the LLP
with respect to the $h$-acyclic $h$-fibrations and the
$h$-acyclic $h$-fibrations are the maps that have the RLP
with respect to the $\bar{h}$-cofibrations.
\end{enumerate}
\end{prop}

To show that $\sC$ has an $h$-type model structure, it suffices
to prove the factorization axioms, and it is unnecessary to prove
them both.

\begin{lem}\mylabel{halffact} For either proposed $h$-model structure,
if one of the factorization axioms holds, then so does the other.
\end{lem}
\begin{proof} For definiteness, consider the case of $h$-fibrations and
$\bar{h}$-cofibrations. By \myref{h-structure}(ii), we can factor
any map $f\colon X\rtarr Y$ as the composite of an
$\bar{h}$-cofibration $i\colon X\rtarr  Mf$ and an $h$-equivalence
$r\colon Mf\rtarr Y$. Suppose that we can factor $r$ as the composite
of an $h$-acyclic $\bar{h}$-cofibration $j\colon Mf\rtarr Z$
and an $h$-fibration $q\colon Z\rtarr Y$. Then $q$ must be
$h$-acyclic, hence $f = q\com (j\com i)$ factors $f$ as the
composite of an $\bar{h}$-cofibration and an $h$-acyclic $h$-fibration.
\end{proof}

A homotopy $X\rtarr Y$ in $\sC$ can be specified by a path $h\colon I\rtarr \sC(X,Y)$. If $i\colon A\rtarr X$ and $p\colon Y\rtarr B$ are maps in $\sC$, then we say that $h$ is a homotopy \emph{relative to} $i$ or 
\emph{corelative to} $p$ if the composite
\[\xymatrix{I\ar[r]^-h & \sC(X,Y) \ar[r]^-{\sC(i,Y)} & \sC(A,Y)}
\qquad\text{or}\qquad
\xymatrix{I\ar[r]^-h & \sC(X,Y) \ar[r]^-{\sC(X,p)} & \sC(X,B)}\]
is constant. When $i$ or $p$ is understood, we also refer to these as homotopies under $A$ or over $B$. The following result is well known and holds in any (based) topologically bicomplete category.  Although we preferred to give a direct proof, we could have derived \myref{h-retract} from this result.  

\begin{prop}\mylabel{relhtpy}
Let $f\colon X\rtarr Y$ be an $h$-equivalence.
\begin{enumerate}[(i)]
\item If $i\colon A\rtarr X$ and $j\colon A\rtarr Y$ are $h$-cofibrations such that $j=f\circ i$, then $f$ is an $h$-equivalence under $A$.
\item If $p\colon Y\rtarr B$ and $q\colon X\rtarr B$ are $h$-fibrations such that $q=p\circ f$, then $f$ is an $h$-equivalence over $B$.
\end{enumerate}
\end{prop}

\begin{proof}
For (i), see for example \cite[p.44]{Concise}. The proof there, although written for spaces, goes through without change. Part (ii) follows by a dual proof.
\end{proof}

\begin{rem}
The current section, as well as the previous and the following one, applies verbatim to the $G$-topologically bicomplete $G$-categories of \S10.2, where $G$ is any topological group.  Of course, $(\sK_{G,B}, G\sK_B)$ is an example. The only changes occur in \myref{h-structure}(iv), where one must take the arrow $G$-spaces $\sC_G(-,-)$ rather than the non-equivariant spaces $G\sC(-,-)$, and in \myref{SVogt}, where the adjunction hypothesis requires a similar equivariant interpretation.  See \S10.3 for a discussion of 
$G$-topological model $G$-categories.
\end{rem}

\section{Classical model structures in general and in $\sK$ and
 $\sU$}\label{sec:classmod}

Again, fix a topologically bicomplete category $\sC$. Independent of the
work of Schw\"anzl and Vogt \cite{SVogt}, Cole \cite{Cole3} proved a general
result concerning when $\sC$ has an $h$-type model structure. As we now
see is inevitable, the core of his argument concerns the verification
of one of the factorization axioms. That requires a hypothesis.

\begin{hyp}\mylabel{hyp}
Let $j_n: Z_n\rtarr Z_{n+1}$ and $q_n: Z_n\rtarr Y$ be maps in $\sC$ such
that $q_{n+1}\com j_n = q_{n}$ and the $j_n$ are $h$-acyclic
$h$-cofibrations.
Let $Z=\text{colim} Z_n$ and let $q\colon Z\rtarr Y$ be obtained by
passage to colimits. Then the canonical map $\text{colim}\, Nq_n \rtarr Nq$ is an isomorphism in $\sC$.
\end{hyp}

\begin{thm}[Cole]\mylabel{Cole}\index{Cole's hypothesis}
If $\sC$ is a topologically bicomplete category which
satisfies \myref{hyp}, then the $h$-equivalences, $h$-fibrations,
and $\bar{h}$-cofibrations specify a proper topological $h$-model
 structure\index{hmodel structure@$h$-model structure} on $\sC$.
\end{thm}
\begin{proof}
It suffices to show that a map $f\colon X\rtarr Y$ factors
as the composite of an $h$-acyclic $\bar{h}$-cofibration $j\colon X\rtarr Z$
and an $h$-fibration $q\colon Z\rtarr Y$. Let $Z_0=X$ and $q_0=f$.
Inductively, given $q_n\colon Z_n\rtarr Y$, construct the following
diagram, in which $Z_{n+1}$ is the displayed pushout.
$$\xymatrix{
Nq_n\ar[d]_{i_0}\ar[r] & Z_n\ar[ddrr]^{q_n}\ar[d]^{j_n} & &\\
\text{Cyl}(Nq_n) \ar[r]^-{\la_n} \ar[rrrd] &
 Z_{n+1}\ar@{-->}[drr]^(0.3){q_{n+1}} & &\\
&&& Y \\}$$
The map $\text{Cyl}(Nq_n)\rtarr Y$ is the adjoint of the projection
$Nq_n\rtarr \text{Cocyl}(Y)$ given by the definition of $Nq_n$, and
$q_{n+1}$ is the induced map. The maps $j_n$ are $h$-acyclic 
$\bar{h}$-cofibrations
since they are pushouts of such maps. Let $Z$ be the colimit of the $Z_n$
and $j$ and $q$ be the colimits of the $j_n$ and $q_n$.  Certainly
$f = q\com j$ and $j$ is an $h$-acyclic $\bar{h}$-cofibration.
By \myref{hyp}, $Nq$ is the colimit of the $Nq_n$. Since the cylinder
functor preserves colimits, we see by \myref{babyish} that $q$ is an
$h$-fibration since the $\la_n$ give a lift $\text{Cyl}(Nq)\rtarr Z$
by passage to colimits.
\end{proof}

The dual version of \myref{Cole} admits a dual proof.

\begin{thm}[Cole]\mylabel{Coledual}
If $\sC$ is a topologically bicomplete category which
satisfies the dual of \myref{hyp}, then the $h$-equivalences,
$\bar{h}$-fibrations, and $h$-cofibrations specify a proper topological
$h$-model structure\index{hmodel structure@$h$-model structure}
on $\sC$.
\end{thm}

From now on, we break the symmetry by focusing on $h$-fibrations and
$\bar{h}$-cofibrations. These give model structures in $\sK$ and $\sU$.
Everything in the rest of the section works equally in $G\sK$ and $G\sU$.
The following theorem combines several results of Str{\o}m
\cite{Strom1, Strom2, Strom}.

\begin{thm}[Str{\o}m]\mylabel{hmodelis}\index{model structure!on spaces}
The following statements hold.
\begin{enumerate}[(i)]
\item  The $h$-equivalences, $h$-fibrations, and $\bar{h}$-cofibrations
give $\sK$ a proper topological $h$-model structure. Moreover, a map in
$\sK$ is an $\bar{h}$-cofibration if and only if it is a closed
$h$-cofibration.
\item The $h$-equivalences, $h$-fibrations, and $\bar{h}$-cofibrations
give $\sU$ a proper topological $h$-model structure. Moreover, a map in
$\sU$ is an $\bar{h}$-cofibration if and only if it is an $h$-cofibration.
\end{enumerate}
\end{thm}
\begin{proof}
\myref{Cole} applies to prove the first statement in $(ii)$,  but it
does not seem to apply to prove the first statement in $(i)$.
The reasons are explained in \myref{Umod}. Taking $Z = Y^I$ and $p=p_0$ there,
the comparison map $\al$ specializes to the map
$\text{colim}\, Nf_n \rtarr Nf$ of \myref{hyp}.
It may be that $\al$ is a homeomorphism in this special case, but
we do not have a proof. It is a homeomorphism when we work in $\sU$.
The characterization of the $\bar{h}$-cofibrations in $\sU$ follows from
\myref{coflemma} and their characterization in $\sK$.

For (i), we give a streamlined version of Str{\o}m's original
arguments that uses the material of the previous section to
prove both statements together. We proceed in four steps.
The first step is Str{\o}m's key observation, the second and
third steps give the second statement, and the fourth step proves
the needed factorization axiom. Consider an inclusion $i\colon A\rtarr X$.

Step 1.  By Str{\o}m's \cite[Thm.\,3]{Strom1}, if $i$ is the inclusion of a
strong deformation retract and there is a map $\ps\colon X\rtarr I$ such
that $\ps^{-1}(0) = A$, then $i$ has the LLP with respect to all
$h$-fibrations.  By \myref{charstrong}(i), this means that $i$ is an 
$h$-acyclic $\bar{h}$-cofibration.

Step 2.  If $i$ is an $h$-cofibration, then the canonical map
$j\colon Mi\rtarr X\times I$ is an $h$-acyclic $h$-cofibration and
therefore, by \myref{h-retract}, the inclusion of a strong deformation 
retract. If $i$ is closed, then $(X,A)$ is an NDR-pair and there exists $\ph\colon X\rtarr  I$
such that $\ph^{-1}(0) = A$. Define $\ps\colon X\times I\rtarr I$ by
$\ps(x,t) = t\ph(x)$. Then $\ps^{-1}(0) = Mi$. Applying Step 1, we
conclude that $j$ has the LLP with respect to all $h$-fibrations. By
\myref{h-char}, this means that $i$ is an $\bar{h}$-cofibration.

Step 3. We can factor any inclusion $i$ as the composite
$$\xymatrix@1{
A\ar[r]^-{i_0} & E \ar[r]^-{\pi} & X,}$$
where $E$ is the subspace $X\times (0,1]\cup A\times I$ of $X\times I$ and
$\pi$ is the projection. Note that $A = \ps^{-1}(0)$, where
$\ps\colon E\rtarr I$ is the projection on the second coordinate.
By direct verification of the CHP \cite[p.\,436]{Strom},
$\pi$ is an $h$-fibration. If $i$ is an $\bar{h}$-cofibration,
then it has the LLP with respect to $\pi$, hence we can lift the identity
map of $X$ to a map $\la\colon X\rtarr E$ such that $\la\com i = i_0$. It
follows that $i(A)$ is closed in $X$ since $i_0(A)$ is closed in $E$.

Step 4.  Let $f\colon X\rtarr Y$ be a map. Use \myref{h-structure}(ii) to
factor $f$ as $p\com s$, where $s\colon X\rtarr Nf$ is the inclusion of a
strong deformation retract and $p$ is an $\bar{h}$-fibration. Use Step 3
to factor $s$ as
$$\xymatrix@1{
X\ar[r]^-{i_0} & Nf\times (0,1]\cup X\times I \ar[r]^-{\pi} & Nf.\\}$$
Here $i_0$ is the inclusion of a strong deformation retract and
$X = \ps^{-1}(0)$, as in Step 3. By Step 1, $i_0$ is an $h$-acyclic
$\bar{h}$-cofibration. By Step 3, $p\com \pi$ is an $h$-fibration.
\end{proof}

There are several further results of Str{\o}m about
$h$-cofibrations that deserve to be highlighted. In order, the following
results are \cite[Theorem 12]{Strom2}, \cite[Lemma 5]{Strom},
and \cite[Corollary 5]{Strom2}.

\begin{prop}\mylabel{StrPull} 
If $p\colon E\rtarr Y$ is an $h$-fibration and the inclusion
$X\subset Y$ is an $\bar{h}$-cofibration, then the induced map
$p^{-1}(X)\rtarr E$ is an $\bar{h}$-cofibration.
\end{prop}

\begin{prop}\mylabel{factorher} If $i\colon A\rtarr B$ and $j\colon B\rtarr
X$ are maps in $\sK$ such that $j$ and
$j\com i$ are $h$-cofibrations, then $i$ is an $h$-cofibration.
\end{prop}

\begin{prop}\mylabel{closure} If an inclusion $A\subset X$ is an
$h$-cofibration, then so is the induced inclusion $\bar{A}\subset X$.
\end{prop}

In view of the characterization of $\bar{h}$-cofibrations in
\myref{hmodelis}, it is natural to ask if there is an analogous
characterization of $\bar{h}$-fibrations. Only the following
sufficient condition is known. It is stated without proof in \cite[4.1.1]{SVogt}, and it gives another reason for requiring
the base spaces of ex-spaces to be in $\sU$.

\begin{prop}\mylabel{hfibshfib}
An $h$-fibration $p\colon E\rtarr Y$ with $Y\in\sU$ is an
$\bar{h}$-fibration.
\end{prop}
\begin{proof}
Let $k\colon A\rtarr X$ be an $h$-acyclic $h$-cofibration and let
$j\colon \overline{A}\rtarr X$ be the induced inclusion.  By Propositions \ref{closure} and \ref{factorher}, $j$ and the inclusion $i\colon A\subset \overline{A}$
are $h$-cofibrations. By \myref{h-retract}(i), $k$ is the inclusion of a deformation retraction $r\colon X\rtarr A$ and the deformation restricts to a homotopy from $(i\circ r)\circ j$ to the identity on $\overline{A}$. 
It follows that $j$ and hence also
$i$ are $h$-acyclic.  Since $j$ is an $h$-acyclic $\bar{h}$-cofibration,
it has the
LLP with respect to $p$, and we see by a little diagram chase that it
suffices to verify that $i$ has the LLP with respect to $p$. Factor $p$
as the composite of $s\colon E\rtarr Np$ and $q\colon Np\colon \rtarr Y$, as
usual. Since $q$ is an
$\bar{h}$-fibration, $(i,q)$ has the lifting property, and it suffices to
show that $(i,s)$ has the lifting property. Suppose given a lifting problem
$f\colon A\rtarr E$ and $g\colon \overline{A}\rtarr Np$ such that
$s\com f = g\com i$. Note that $s(e) = (e,c_{p(e)})$ for $e\in E$,
where $c_y$ denotes the constant path at $y$. Since $Y$ is weak Hausdorff,
the constant paths give a closed subset of $Y^I$ and $Np = Y^I\times_Y E$
is a closed subset of $Y^I\times E$. Therefore $s(E)$ is closed in $Np$. 
We conclude that
$$g(\overline{A})\subset \overline{g(A)} = \overline {s(f(A)} \subset
\overline{s(E)} = s(E),$$
which means that there is a lift $\overline{A}\rtarr E$.
\end{proof}

\section{Compactly generated $q$-type model structures}

We give a variant of the standard procedure for constructing $q$-type model
structures. The exposition prepares the way for a new variant that we will
explain in \S5.4 and which is crucial to our work. Although our discussion is adapted to topological examples, $\sC$ need not be topological until otherwise specified. We first recall the small object argument in settings where compactness allows use of sequential colimits.

\begin{defn}\mylabel{cofhyp} Let $I$ be a set of maps in $\sC$.
\begin{enumerate}[(i)]
\item A \emph{relative $I$-cell complex}\index{relative cell
complex}\index{cell complex!relative} is a map $Z_0\rtarr Z$, where 
$Z$ is the colimit of a sequence of maps $Z_n\rtarr Z_{n+1}$ such that
$Z_{n+1}$ is the pushout $Y\cup_X Z_n$ of a coproduct $X\rtarr Y$ of 
maps in $I$ along a map $X\rtarr Z_n$. 
\item $I$ is \emph{compact}\index{compact!set of maps} if for every domain
object $X$ of a map in $I$ and every relative $I$-complex $Z_0\rtarr Z$,
the map $\text{colim}\,\sC(X,Z_n)\rtarr \sC(X,Z)$ is a bijection.
\item An \emph{$I$-cofibration}\index{cofibration!I@$I$-} is a map that
satisfies the LLP with respect to any map that satisfies the RLP with
respect to $I$.
\end{enumerate}
\end{defn}

\begin{lem}[Small object argument]\mylabel{small}\index{small object
argument}
Let $I$ be a compact set of maps in $\sC$, where $\sC$ is co\-com\-plete.
Then any map $f:X\rtarr Y$ in $\sC$ factors functorially as a composite
\[\xymatrix{X \ar[r]^i & W \ar[r]^p & Y}\]
such that $p$ satisfies the RLP with respect to $I$ and $i$ is a relative
$I$-cell complex and therefore an $I$-cofibration.
\end{lem}

\begin{defn}\mylabel{compgendef}
A model structure on $\sC$ is \emph{compactly generated}\index{compactly
generated!model structure}\index{model structure!compactly generated} if
there are compact sets $I$ and $J$ of maps in $\sC$ such that the following
characterizations hold.
\begin{enumerate}[(i)]
\item  The fibrations are the maps that satisfy the RLP with respect to $J$,
or equivalently, with respect to retracts of relative $J$-cell complexes.
\item The acyclic fibrations are the maps that satisfy the RLP with respect
to $I$, or equivalently, with respect to retracts of relative $I$-cell
 complexes.
\item The cofibrations are the retracts of relative $I$-cell complexes.
\item The acyclic cofibrations are the retracts of relative $J$-cell
 complexes.
\end{enumerate}
The maps in $I$ are called the \emph{generating
cofibrations}\index{generating cofibration}\index{cofibration!generating}
and the maps in $J$ are
called the \emph{generating acyclic cofibrations}. 
\end{defn}

We find it convenient to separate out properties of classes of maps in a
model category, starting with the weak equivalences.

\begin{defn}\mylabel{weakeqsub} A subcategory of $\sC$ is a
\emph{subcategory of weak equivalences}\index{weak equivalence!subcategory
 of --s} if it satisfies the following closure properties.
\begin{enumerate}[(i)]
\item All isomorphisms in $\sC$ are weak equivalences.
\item A retract of a weak equivalence is a weak equivalence.
\item If two out of three maps $f$, $g$, $g\com f$ are weak
equivalences, so is the third.
\end{enumerate}
\end{defn}

\begin{thm}\mylabel{compgen}
Let $\sC$ be a bicomplete category with a subcategory of weak equivalences.
Let $I$ and $J$ be compact sets of maps in $\sC$. Then $\sC$ is a compactly
generated model category with generating cofibrations $I$ and generating
acyclic cofibrations $J$ if the following two conditions hold:
\begin{enumerate}[(i)]
\item (Acyclicity condition) Every relative $J$-cell complex is a weak equivalence.
\item (Compatibility condition) A map has the RLP with respect to $I$ if
and only if it is a weak equivalence and has the RLP
with respect to $J$.
\end{enumerate}
\end{thm}

\begin{proof} This is the formal part of Quillen's original proof of the
 $q$-model structure on topological spaces and is a variant of
 \cite[2.1.19]{Hovey} or \cite[11.3.1]{Hirschhorn}.  The fibrations are
defined to be the maps that satisfy the RLP with respect to $J$. The cofibrations are defined to be the $I$-cofibrations and turn out to be the retracts of relative
$I$-cell complexes. The retract axioms clearly hold and, by (ii), the
cofibrations are the maps that satisfy the LLP with respect to the acyclic
fibrations, which gives one of the lifting axioms. The maps in $J$ satisfy
the LLP with respect to the fibrations and are therefore cofibrations,
which verifies something that is taken as a hypothesis in the versions in
the cited sources. Applying the small object argument to $I$, we factor a
map $f$ as a composite of an $I$-cofibration followed by a map that
satisfies the RLP with respect to $I$; by (ii), the latter is an acyclic
fibration. Applying the small object argument to $J$, we factor $f$ as a
composite of a relative $J$-cell complex that is a $J$-cofibration followed
by a fibration.  By (i), the first map is acyclic, and it is a cofibration
because it satisfies the LLP with respect to all fibrations, in particular
the acyclic ones.  Finally, for the second lifting axiom, if we are given a
lifting problem with an acyclic cofibration $f$ and a fibration $p$, then a
standard retract argument shows that $f$ is a retract of an acyclic
cofibration that satisfies the LLP with respect to all fibrations.
\end{proof}

Using the following companion to \myref{weakeqsub}, we codify the usual
pattern for verifying the acyclicity condition. 

\begin{defn}\mylabel{cofsub} A subcategory of a cocomplete category
$\sC$ is a \emph{subcategory of cofibrations}\index{cofibration!subcategory
 of --s}
if it satisfies the following closure properties.
\begin{enumerate}[(i)]
\item All isomorphisms in $\sC$ are cofibrations.
\item All coproducts of cofibrations are cofibrations.
\item If $i\colon X \rtarr Y$ is a cofibration and $f\colon X\rtarr Z$ is
any map, then the pushout $j\colon Y\rtarr Y\cup_X Z$ of $f$ along $i$ is a cofibration.
\item If $X$ is the colimit of a sequence of cofibrations
$i_n\colon X_n\rtarr X_{n+1}$, then the induced map $i\colon X_0\rtarr X$ is
a cofibration.
\item A retract of a cofibration is a cofibration.
\end{enumerate}
\end{defn}

In more general contexts, (iv) should be given a transfinite
generalization, but we shall not have need of that. Note that if a 
subcategory of cofibrations is defined in terms of a left
lifting property, then all of the conditions hold automatically.

\begin{lem}\mylabel{cofcat}
Let $\sC$ be a cocomplete category together with a subcategory of
cofibrations, denoted \emph{$g$-cofibrations}, and a subcategory 
of weak equivalences, satisfying the following properties.
\begin{enumerate}[(i)]
\item A coproduct of weak equivalences is a weak equivalence.
\item If $i\colon X \rtarr Y$ is an acyclic $g$-cofibration and $f\colon
X\rtarr Z$ is any map, then the pushout $j\colon Y\rtarr Y\cup_X Z$ of $f$
along $i$ is a weak equivalence.
\item If $X$ is the colimit of a sequence of acyclic $g$-cofibrations
$i_n\colon X_n\rtarr X_{n+1}$, then the induced
map $i\colon X_0\rtarr X$ is a weak equivalence.
\end{enumerate}
If every map in a set $J$ is an acyclic $g$-cofibration, then every 
relative $J$-cell complex is a weak equivalence.
\end{lem}

We emphasize that the $g$-cofibrations are not the model category
cofibrations and may or may not be the intrinsic $h$-cofibrations or
$\bar{h}$-cofibrations. They serve as a convenient scaffolding for 
proving the model axioms.

\begin{rem} The properties listed in \myref{cofcat} include some of the
axioms for a ``cofibration category'' given by Baues \cite[pp 6, 182]{Baues}.
However, our purpose is to describe features of categories that
are more richly structured than model categories, often with several
relevant subcategories of cofibrations, rather than to describe
deductions from axiom systems for less richly structured categories, 
which is his focus. The $g$-cofibrations in \myref{cofcat} need not be 
the cofibrations of any cofibration category or model category.
\end{rem}

The $q$-model structures on $\sK$ and $\sU$ are obtained by \myref{compgen},
taking the $q$-equivalences to be the weak equivalences, that is, the maps 
that induce isomorphisms on all homotopy groups, and the $q$-fibrations to
be the Serre fibrations. We also have the equivariant generalization, which
applies to any topological group $G$. We introduce the following notations, which will be used throughout. 

\begin{defn}\mylabel{UrIJ}
Nonequivariantly, let $I$ and $J$ denote the set of inclusions $i\colon S^{n-1} \rtarr D^n$ (where $S^{-1}$ is empty) and the set of maps $i_0\colon D^n\rtarr D^n\times I$. Equivariantly, let $I$ and $J$ denote the set of all maps of the form $G/H\times i$, where $H$ is a (closed) subgroup of $G$ and $i$ runs through the maps in the nonequivariant sets $I$ and $J$.  In the based categories $\sK_*$ and 
$G\sK_*$ we continue to write $I$ and $J$ for the sets obtained 
by adjoining disjoint base points to the specified maps.
\end{defn}

A map $f\colon X\rtarr Y$ of $G$-spaces is said to be a weak equivalence or Serre fibration if all fixed point maps $f^H\colon X^H\rtarr Y^H$ are weak equivalences or Serre fibrations.  Just as nonequivariantly, we also call these $q$-equivalences and $q$-fibrations. Observe that $q$-equivalences are defined in terms of the equivariant homotopy groups $\pi_n^H(X,x) = \pi_n(X^H,x)$ for $H\subset G$ and $x\in X^H$ and that $q$-fibrations are defined in terms of the RLP with respect to the cells in $J$.

If $X_0 \rtarr X$ is a relative $I$ or $J$-cell complex, then $X/X_0$ is in $G\sU$ and \myref{little} gives all that is needed to verify the compactness hypothesis in \myref{cofhyp}(ii).  Taking the $g$-cofibrations to be the $h$-cofibrations, \myref{cofcat} applies to verify the acyclicity condition of \myref{compgen}. With considerable simplification, our verification of the compatibility condition for the $qf$-model structure in Chapter 6 specializes to verify it here.  Nonequivariantly, the $q$-model structure is discussed in \cite[\S8]{DS} and, with somewhat different details, in \cite[2.4]{Hovey} 
(where the details on transfinite sequences are unnecessary). 

Equivariantly, a detailed proof of the following result is given in \cite[III\S1]{MM}. The argument there is given for based $G$-spaces, in $G\sT$, but it works equally well for unbased $G$-spaces, in $G\sK$. 

\begin{thm}\mylabel{Gold}
For any $G$, $G\sK$ is a compactly generated proper mo\-del category whose 
$q$-equivalences, $q$-fibrations, and $q$-co\-fib\-rations are the weak equivalences, the Serre fibrations, and the retracts of relative $G$-cell complexes. The sets $I$ and $J$ are the generating $q$-cofibrations and the generating acyclic $q$-cofibrations, and all $q$-cofibrations are $\bar{h}$-cofibrations. If $G$ is a compact Lie group, then the model structure is
$G$-topological.
\end{thm}

The notion of a $G$-topological model category is defined 
in the same way as the notion of a simplicial or topological 
model category and is discussed formally in \S10.3 below.  
The point of the last statement is that if $H$ and $K$ are
subgroups of a compact Lie group $G$, then  $G/H\times G/K$ 
has the structure of a $G$-CW complex.  By \myref{Illman}, 
this remains true when $G$ is a Lie group and $H$ and $K$ 
are compact subgroups. We shall see how to use this fact 
model theoretically in Chapter 7.

\chapter{Well-grounded topological model categories}

\section*{Introduction}

It is essential to our theory
to understand the interrelationships among the various model structures that appear naturally in the parametrized context, both in topology and in general.
This understanding leads us more generally to an axiomatization of the properties that are required of a good $q$-type model structure in order 
that it relate well to the classical homotopy theory on a topological category.  The obvious $q$-model structure on ex-spaces over $B$ does not satisfy the axioms, and in the next chapter we will introduce a new model structure, the $qf$-model structure, that does satisfy the axioms.  
 
As we recall in \S5.1, any model structure on a category $\sC$ induces
a model structure on the category of objects over, under, or over and under
a given object $B$.  When $\sC$ is topologically bicomplete, so are these over and under categories.  They then have their own intrinsic $h$-type model structures, which differ from the one inherited from $\sC$.  This leads to 
quite a few different model structures on the category $\sC_B$ of objects over and under $B$, each with its own advantages and disadvantages. Letting $B$ vary, we also obtain a model structure on the category of retracts. We shall only 
be using most of these structures informally, but the 
plethora of model structures is eye opening.

In \S5.2, we focus on spaces and compare the various classical notions of fibrations and cofibrations that are present in our over and under categories. Although elementary, this material is subtle, and it is nowhere presented accurately in the literature.  In particular, we discuss $h$-type, $f$-type and $fp$-type model structures, where $f$ and $fp$ stand for ``fiberwise'' and ``fiberwise pointed''.  For simplicity, we discuss this material nonequivariantly, but it applies verbatim equivariantly. 

The comparisons among the $q$, $h$, $f$, and $fp$ classes of maps and model structures guide our development of parametrized homotopy theory. We think of the $f$-notions as playing a transitional role, connecting the $fp$ and 
$h$-notions.  In the rest of the chapter, we work in a general topologically bicomplete category $\sC$, and we sort out this structure and its 
relationship to a desired $q$-type model structure axiomatically.

Here we shift our point of view.  We focus on three basic types of cofibrations that are in play in the general context, namely the Hurewicz cofibrations determined by the cylinders in $\sC$, the ground cofibrations that come in practice from a given forgetful functor to underlying spaces, and the $q$-type model cofibrations. The first two are intrinsic, but we think of the $q$-type cofibrations as subject to negotiation. In $\sK_B$, the Hurewicz cofibrations are the $fp$-cofibrations and the ground cofibrations are the $h$-cofibrations, which is in notational conflict with the point of view taken in the previous chapter.

In \S\S5.3 and 5.4, we ignore model theoretic considerations entirely. We describe how the two intrinsic types of cofibrations relate to each other 
and to colimits and tensors, and we explain how this structure relates to weak equivalences. 

We define the notion of a  ``well-grounded model structure'' in \S5.5.  We believe that this notion captures exactly the right blend of classical and model categorical homotopical structure in topological situations. It describes what is needed for a $q$-type model structure in a topologically bicomplete category to be compatible with its intrinsic $h$-type model structure and its ground structure. Crucially, the $q$-type cofibrations should be ``bicofibrations'', meaning that they are {\em both} Hurewicz cofibrations and ground cofibrations.  To illustrate the usefulness of the axiomatization, and for later reference, we derive the long exact sequences associated to cofiber sequences and the $\text{lim}^1$ exact sequences associated to colimits in \S5.6.

A clear understanding of the desiderata for a good $q$-type model structure reveals that the obvious over and under $q$-model structure is essentially worthless for serious work in parametrized homotopy theory.  This will lead us to introduce the new $qf$-model structure, with better behaved $q$-type cofibrations, in the next chapter.  The formalization given in \S\S5.3--5.6 might seem overly pedantic were it only to serve as motivation for the definition of the $qf$-model structure.  However, we will encounter exactly the same structure in Part III when we construct the level and stable model structures on parametrized spectra.  We hope that the formalization will help guide the reader through the rougher terrain there.

We note parenthetically that there is still another interesting model structure on the category of ex-spaces over $B$, one based on local considerations. It is due to Michelle Intermont and Mark Johnson \cite{IJ}. We shall not discuss their model structure here, but we are indebted to them for illuminating discussions. It is conceivable that their model structure could be used in an alternative development of the stable theory, but that has not been worked out. Their structure suffers the defects that it is not known to be left proper and that, with their definition of weak equivalences, homotopy equivalences of base spaces need not induce equivalences of homotopy categories.

We focus mainly on the nonequivariant context in this chapter, but $G$ 
can be any topological group in all places where equivariance is considered.

\section{Over and under model structures}

Recall from \S1.2 that, for any category $\sC$ and object $B$ in $\sC$, we let $\sC/B$\noteindex{CBa@$\sC/B$} and $\sC_B$ denote the categories of objects over $B$ and of ex-objects over $B$. We also have the category $B\backslash \sC$ of objects under $B$. If $\sC$ is bicomplete, then so are $\sC/B$, $B\backslash \sC$ and $\sC_B$.  We begin with some general observations about over and under model categories before returning to topological categories.

We have forgetful functors $U\colon \sC/B\rtarr \sC$\noteindex{U@$U$} and $V\colon \sC_B\rtarr \sC/B$.\noteindex{V@$V$} The first is left adjoint to the functor that sends an object $Y$ to the object $B\times Y$ over $B$: 
\begin{equation}\label{silly}
\sC(UX,Y) \iso \,\sC/B\,(X, B\times Y).
\end{equation}
The second is right adjoint to the functor that sends an object $X$ over $B$
to the object $X\amalg B$ over and under $B$:
\begin{equation}\label{CexC}
\sC_B(X\amalg B,Y) \iso \,\sC/B\,(X,VY).
\end{equation}
As a composite of a left and a right adjoint, the total object functor
$UV\colon \sC_B\rtarr \sC$ does not enjoy good formal properties. This
obvious fact plays a significant role in our work. For example, it limits 
the value of
the model structures on $\sC_B$ that are given by the following result.
\begin{prop}\mylabel{under}\index{model structure!over and under}
Let $\sC$ be a model category. Then $\sC/B$, $B\backslash\sC$, and $\sC_B$
are model categories in which the weak equivalences, cofibrations, and fibrations are the maps over $B$, under $B$, or over and under $B$ which 
are weak equivalences, fibrations, or cofibrations
in $\sC$. If $\sC$ is left or right proper, then so are $\sC/B$,
 $B\backslash\sC$, and $\sC_B$.
\end{prop}\begin{proof}
As observed in \cite[p.\,5]{Hovey} and \cite[3.10]{DS}, the
statement about $\sC/B$ is a direct verification from the definition of a
model category. By the self-dual nature of the axioms, the statement about
$B\backslash\sC$ is equivalent. The statement about $\sC_B$ follows since
it is the category of objects under $(B,\text{id})$ in
$\sC/B$. The last statement holds since pushouts and pullbacks in these
over and under categories are constructed in $\sC$. \end{proof}

When considering $q$-type model structures, we start with a compactly
generated model category $\sC$. Using the adjunctions (\ref{silly})
and (\ref{CexC}), we then obtain the following addendum to \myref{under}.

\begin{prop}\mylabel{cg} 
If $\sC$ is a compactly generated model category, then $\sC/B$ and $\sC_B$ are compactly generated. The generating (acyclic) cofibrations in $\sC/B$ are the maps $i$ such that $Ui$ is a generating (acyclic) cofibration in $\sC$. The generating (acyclic) cofibrations in $\sC_B$ are the maps $i\amalg B$ where $i$ is a generating (acyclic) cofibration in $\sC/B$.
\end{prop}

We now return to the case when $\sC$ is topologically bicomplete. Then it has the resulting ``classical'', or $h$-type, structure that was discussed in \S\ref{sec:towardh} and \S\ref{sec:classmod}. If our philosophy in \S\ref{Sphil} applies to $\sC$, then it also has $q$ and $m$-structures and the categories $\sC/B$ and $\sC_B$ both inherit over and under model structures that are related as we discussed there. However, since $\sC$ is topologically bicomplete, so is $\sC/B$ by \myref{topbicomp}, and $\sC_B$ is based topologically bicomplete by \myref{btopbicomp}. These categories therefore have classical 
$h$-type structures when they are regarded in their own right as topologically bicomplete categories. To fix notation and avoid confusion we give an overview of all of these structures.

We start with the $h$-classes of maps in $\sC$ that are given in \myref{hmodel} and \myref{h-char}.  As in our discussion of spaces, we work assymmetrically, ignoring the $\bar{h}$-fibrations and focusing on the candidates for $h$-type model structures given by the $h$-fibrations and $\bar{h}$-cofibrations. We agree to use the letter $h$ for the inherited classes of maps in $\sC/B$ and $\sC_B$, although that contradicts our previous use of $h$ for the classical classes of maps in an arbitrary topologically bicomplete category, such as $\sC/B$ or $\sC_B$.  We shall resolve that ambiguity shortly by introducing new names for the classes of ``classical'' maps in those categories.

\begin{defn}
A map $g$
in $\sC/B$ is an $h$-equivalence,\index{equivalence!h@$h$-}
$h$-fibration,\index{fibration!h@$h$-}
$h$-co\-fi\-bra\-tion,\index{cofibration!h@$h$-} or
$\bar{h}$-cofibration\index{cofibration!hs@$\bar{h}$-} if $Ug$ is 
such a map in $\sC$. A map $g$ in $\sC_B$
is an $h$-equivalence, $h$-fibration, $h$-cofibration, or
$\bar{h}$-cofibration if $Vg$ is such a map in $\sC/B$ or, equivalently,
$UVg$ is such a map in $\sC$.
\end{defn}

The $\bar{h}$-cofibrations are $h$-cofibrations, but not conversely in
general. Since the object $*_B=(B,\text{id},\text{id})$ is initial and
terminal in $\sC_B$, an object of $\sC_B$ is $h$-cofibrant (or $\bar{h}$-cofibrant) if its section is an $h$-cofibration (or $\bar{h}$-cofibration) in $\sC$. It is $h$-fibrant if its projection is an $h$-fibration in $\sC$.

In $\sC/B$, we have the notion of a homotopy over $B$, defined in terms of
$X\times_B I$ or, equivalently, $\text{Map}_B(I,X)$. The adjective
``fiberwise'' is generally used in the literature to describe these
homotopies. See, for example, the books \cite{CJ, James} on fiberwise
homotopy theory. To distinguish from the $h$-model structure, we agree to 
write $f$ rather than $h$ for the fiberwise specializations of \myref{hmodel} and \myref{h-char}.  To avoid any possible confusion, we formalize this, 
making use of \myref{charstrong}.

\begin{defn}\mylabel{fmodel}
Let $g$ be a map in $\sC/B$.
\begin{enumerate}[(i)]
\item $g$ is an \emph{$f$-equivalence}\index{equivalence!f@$f$-} if it is a
fiberwise homotopy equivalence.
\item $g$ is an \emph{$f$-fibration}\index{fibration!f@$f$-} if it satisfies
 the fiberwise CHP, that is, if it has the RLP with respect to the maps
$i_0\colon X\rtarr X\times_B I$ for $X\in\sC/B$.
\item $g$ is an \emph{$f$-cofibration}\index{cofibration!f@$f$-} if it
satisfies the fiberwise HEP, that is, if it has the LLP with
respect to the maps $p_0\colon \text{Map}_B(I,X) \rtarr X$.
\item $g$ is an \emph{$\bar{f}$-cofibration}\index{equivalence!fs@$\bar{f}$-}
if it has the LLP with respect to the $f$-acyclic $f$-fibrations.
\end{enumerate}
A map $g$ in $\sC_B$ is an $f$-equivalence, $f$-fibration, $f$-cofibration, 
or $\bar{f}$-cofibration if $Vg$ is one in $\sC/B$.
\end{defn}

Again, $\bar{f}$-cofibrations are $f$-cofibrations, but not conversely in general. \myref{Cole} often applies to show that the $f$-fibrations and $\bar{f}$-cofibrations define an $f$-model structure on $\sC/B$ and therefore, by \myref{under}, on $\sC_B$.  As is always the case for an intrinsic classical model structure, every object of $\sC/B$ is both $f$-cofibrant and $\bar{f}$-cofibrant as well as $f$-fibrant. While this is obvious from the definitions, it may seem counterintuitive. It does not follow that every object of $\sC_B$ is $f$-cofibrant since the two categories have different initial objects.

In $\sC_B$, we also have the notion of a homotopy over and under $B$,
defined in terms of $X\sma_B I_+$ or, equivalently, $F_B(I_+,X)$.
The adjective ``fiberwise pointed'' is used in \cite{CJ, James} to
describe these homotopies. Again, for notational clarity, we agree to
write $fp$ rather than $h$ for the fiberwise pointed specializations of 
\myref{hmodel} and \myref{h-char}, and we formalize this to avoid any
possible confusion.

\begin{defn}\mylabel{fpmodel}
Let $g$ be a map in $\sC_B$.
\begin{enumerate}[(i)]
\item $g$ is an \emph{$fp$-equivalence} if it is a fiberwise pointed
homotopy equivalence.
\item $g$ is an \emph{$fp$-fibration} if it satisfies the fiberwise
point\-ed CHP, that is,
if it has the RLP with respect to the maps $i_0\colon X\rtarr X\sma_B I_+$.
\item $g$ is a \emph{$fp$-cofibration} if it satisfies the fiberwise
point\-ed HEP, that is, if it has the LLP with
respect to the maps $p_0\colon F_B(I_+,X) \rtarr X$.
\item $g$ is an \emph{$\overline{fp}$-cofibration} if it has the LLP
with respect to the $fp$-acyclic $fp$-fibrations.
\end{enumerate}
\end{defn}

Again, $\overline{fp}$-cofibrations are $fp$-cofibrations, but not conversely
in general, and \myref{Cole} often applies to show that the $fp$-fibrations
and $\overline{fp}$-cofibrations define an $fp$-model structure on $\sC_B$. 
We summarize some general formal implications relating our classes of maps.

\begin{prop}\mylabel{compare} Let $\sC$, $\sC/B$ and $\sC_B$ be
topologically bicomplete
categories with $h$, $f$, and $fp$-classes of maps defined as above.
Then the following implications hold for maps in $\sC_B$.

\vspace{1mm}

\begin{center}
\begin{tabular}{|r c c c l|} \hline
$fp$-equivalence & $\Longrightarrow$ &  $f$-equivalence & $\Longrightarrow$
 & $h$-equivalence\\
$fp$-cofibration & $\Longleftarrow$ & $f$-cofibration &$\Longrightarrow$
& $h$-cofibration\\
$\Uparrow$ \ \ \ \ \ \ \ \ & & $\Uparrow$ & & $\ \ \ \ \ \ \ \ \Uparrow$\\
$\overline{fp}$-cofibration & $\Longleftarrow$ &  $\bar{f}$-cofibration &
$\Longrightarrow$ & $\bar{h}$-cofibration\\
$fp$-fibration & $\Longrightarrow$ &  $f$-fibration & $\Longleftarrow$ &
 $h$-fibration\\\hline
\end{tabular}
\end{center}

\vspace{1mm}

\noindent
Moreover, every object of $\sC_B$ is both $fp$-fibrant and $fp$-cofibrant.
\end{prop}\begin{proof}
Trivial inspections of lifting diagrams show that an $h$-fibration is an
$f$-fibration, an $f$-cofibration is an $fp$-cofibration, and an
$\bar{f}$-co\-fi\-bra\-tion is an $\overline{fp}$-cofibration. Use of the
adjunctions (\ref{silly}) and (\ref{CexC}) shows that
an $f$-cofibration is an $h$-cofibration, an $\bar{f}$-cofibration is an
$\bar{h}$-cofibration, and an $fp$-fibration is an $f$-fibration. The last
statement holds since fiberwise pointed homotopies with domain or target
$B$ are constant at the section or projection of the target or source.
\end{proof}

\begin{rem} Assume that these classes of maps define model structures.
Then the implications in \myref{compare} lead via \myref{Colemix} and its
dual version to two new mixed model structures on $\sC_B$, one with weak
equivalences the $f$-equivalences and fibrations the $fp$-fibrations and
one with weak equivalences the $h$-equivalences and cofibrations the
$\bar{f}$-cofibrations.
\end{rem}

The category $\sC_{\sB}$ of retracts introduced in \S2.5 suggests an alternative model theoretic point of view. We give the basic definitions, but we shall not pursue this idea in any detail.  Again, \myref{Cole} often applies to verify the model category axioms.  Note that the intrinsic homotopies are given by homotopies of total objects over and under homotopies of base objects.

\begin{defn}\mylabel{rmodel}
Assume that $\sC_{\sB}$ is topologically bicomplete and let $g$ be a map in
$\sC_{\sB}$.
\begin{enumerate}[(i)]
\item $g$ is an \emph{$r$-equivalence}\index{equivalence!r@$r$-} if it is a
homotopy equivalence of retractions.
\item $g$ is an \emph{$r$-fibration}\index{fibration!r@$r$-} if it satisfies
the retraction CHP, that is, if it has the RLP with respect to the maps
$i_0\colon X\rtarr X\times I$ for $X\in\sC_{\sB}$.
\item $g$ is an \emph{$r$-cofibration}\index{cofibration!r@$r$-} if it
satisfies the retraction HEP,
that is, if it has the LLP with respect to the maps
$p_0\colon \text{Map}(I,X) \rtarr X$.
\item $g$ is an \emph{$\bar{r}$-cofibration}\index{cofibration!r@$\bar{r}$-}
if it has the LLP
with respect to the $r$-acyclic $r$-fibrations.
\end{enumerate}
\end{defn}

\begin{rem}
The initial and terminal object of $\sC_{\sB}$ are the identity retractions of the initial and terminal objects of $\sB$ and every object is both $r$-cofibrant and $r$-fibrant. It might be of interest to characterize the retractions for which the map $*_B\rtarr (X,p,s)$ induced by $s$ is an $r$-cofibration or for which the map $(X,p,s)\rtarr *_B$ induced by $p$ is an $r$-fibration. By specialization of the lifting properties, an ex-map over $B$ that is an $r$-cofibration or $r$-fibration is an $fp$-cofibration or $fp$-fibration in $\sC_B$, but we have not pursued this question further.
\end{rem}

\section{The specialization to over and under categories of spaces}

Now we take $\sC$ to be $\sK$ or $\sU$. We discuss the relationships
among our various classes of fibrations and cofibrations in this
special case, and we consider when the $f$ and $fp$ classes of maps
give model structures. Everything in this section applies equally 
well equivariantly.

We first say a bit about based spaces, which are ex-spaces over $B=\{*\}$.
Here the fact that $*$ is a terminal object greatly simplifies matters.
All of the $f$-notions coincide with the corresponding $h$-notions, and
our trichotomy reduces to the familiar dichotomy between free (or $h$)
notions and based (or $fp$) notions.  Recall that a based space is \emph{well-based},\index{space!well-based --} or \emph{nondegenerately based}, if the inclusion of the 
basepoint is an $h$-cofibration. Every based space is $fp$-cofibrant,
and an $fp$-cofibration between well-based spaces is an $h$-cofibration \cite[Prop.\,9]{Strom}.  Every based space is $fp$-fibrant, and an 
$h$-fibration of based spaces satisfies the based CHP
with respect to well-based source spaces. Of course, the over and 
under $h$-model structure differs from the intrinsic $fp$-model structure.

None of the reverse implications in \myref{compare} holds in general.
We gave details of that result since it is easy to get confused and
think that more is true than we stated.
\begin{sch}
On \cite[p.\,66]{CJ}, it is stated that a fiberwise
pointed cofibration which is a closed inclusion is a fiberwise
cofibration. That is false even when $B$ is a point, since
it would imply that every point of a $T_1$-space is a nondegenerate
basepoint. On \cite[p.\,69]{CJ}, it is stated that a fiberwise pointed
map (= ex-map) is a fiberwise pointed fibration if and only if it is
a fiberwise fibration. That is also false when $B$ is a point, since
the unbased CHP does not imply the based CHP.
\end{sch}

However, as for based spaces, the reverse implications in parts of
\myref{compare} often do hold under appropriate additional hypotheses.
\begin{prop}\mylabel{reverse}
The following implications hold for an arbitrary topologically bicomplete category $\sC$.
\begin{enumerate}[(i)]
\item A map in $\sC/B$ between $h$-fibrant objects over $B$ is an
$h$-equivalence if and only if it is an $f$-equivalence.
\item An ex-map between $f$-cofibrant ex-objects over $B$ is an
$f$-equivalence if and only if it is an $fp$-equivalence.
\end{enumerate}
\end{prop}

\begin{proof}
The first part follows from \myref{relhtpy}(ii) since an $f$-equivalence in $\sC/B$ is the same as an $h$-equivalence over $B$ in $\sC$. The second part follows similarly from \myref{relhtpy}(i) since an $fp$-equivalence in $\sC_B$ is the same as an $f$-equivalence under $B$ in $\sC/B$.
\end{proof}

The following results hold for spaces. We are doubtful that they hold in general.

\begin{prop}\mylabel{reverse2}
The following implications hold in both $G\sK$ and $G\sU$.
\begin{enumerate}[(i)]
\item An ex-map between $\bar{f}$-cofibrant ex-spaces is an $f$-cofibration
if and only if it is an $fp$-cofibration.
\item An ex-map whose source is $\bar{f}$-cofibrant is an $f$-fibration if and
only if it is an $fp$-fibration.
\end{enumerate}
\end{prop}\begin{proof}
Part (ii) is \cite[16.3]{CJ}. 
Part (i) is stated on \cite[p.\,441]{Strom} and the proof given there for based spaces generalizes using the following lemma.
\end{proof}

It is easy to detect $f$-cofibrations by means of the following result,
whose proof is the same as that of the standard characterization of
Hurewicz cofibrations (e.g.\ \cite[p.\,43]{Concise}, see also
\cite[Thm.\,2]{Strom1}, \cite[Lem.\,4]{Strom2} and \cite[4.3]{CJ}).

\begin{lem}\mylabel{fNDR}\index{NDR}
An inclusion $i\colon X\rtarr Y$ in $\sK/B$ is an $f$-cofibration if and
only if $(Y,X)$ is a fiberwise NDR-pair in the sense that there is a map
$u\colon Y\rtarr I$ such that $X \subset u^{-1}(0)$ and a homotopy $h\colon
Y\times_B I \rtarr Y$ over $B$ such that $h_0 = \text{id}$, $h_t = \text{id}$ 
on $X$ for $0\leq t\leq 1$, and $h_1(y)\in X$ if $u(y)<1$.
A closed inclusion $i: X\rtarr Y$ in $\sK/B$ is an $\bar{f}$-cofibration if
and only if the map $u$ above can be chosen so that $X=u^{-1}(0)$.
\end{lem}

We introduce the following names here, but we defer a full discussion
to \S8.1.

\begin{defn}\mylabel{names}  An ex-space is said to be \emph{well-sectioned}\index{well-sectioned!ex-space}\index{ex-space!well-sectioned} 
if it is $\bar{f}$-cofibrant. An ex-space is said to be {\em ex-fibrant}
or, synonomously, to be an \emph{ex-fibration}\index{ex-fibration} if it is both $\bar{f}$-cofibrant and $h$-fibrant. Thus an ex-fibration is a well-sectioned ex-space whose projection is an $h$-fibration.
\end{defn}

The term ex-fibrant is more logical than ex-fibration, since we are defining 
a type of object rather than a type of morphism of $\sK_B$, but the term 
ex-fibration goes better with Serre and Hurewicz fibration and is standard 
in the literature.  We have the following implication of Propositions \ref{compare} and \ref{reverse}. It helps explain the usefulness of ex-fibrations.

\begin{cor}\mylabel{fpequiv} Let $g$ be an ex-map between ex-fibrations 
over $B$.
\begin{enumerate}[(i)]
\item $g$ is an $h$-equivalence if and only if $g$ is an $f$-equivalence,
and this hold if and only if $g$ is an $fp$-equivalence.
\item $g$ is an $f$-cofibration if and only if $g$ is an $fp$-cofibration,
and then $g$ is an $h$-cofibration.
\item $g$ is an $f$-fibration if and only if $g$ is an $fp$-fibration,
and this holds if $g$ is an $h$-fibration.
\end{enumerate}
\end{cor}

\begin{rem}\mylabel{guess}  The model theoretic significance of ex-fibrations
over $B$ is unclear. They are fibrant and cofibrant objects in the mixed
model structure on ex-spaces over $B$ whose weak equivalences are the
$h$-equivalences
and whose cofibrations are the $\bar{f}$-cofibrations.  However, the
converse fails since there are well-sectioned $f$-fibrant ex-spaces that
are $f$-equivalent to $h$-fibrant ex-spaces, hence are mixed fibrant, but
are not themselves $h$-fibrant.
\end{rem}

The previous remark anticipated the following result on over and under
model structures in the categories of spaces and ex-spaces over $B$. 
Note that \myref{coflemma} applies to $\sK/B$ and $\sK_B$ as well 
as to $\sK$ to show that both $f$-cofibrations and $fp$-cofibrations are 
inclusions which are closed when the total spaces are in $\sU$.

\begin{thm}\mylabel{ffpmodel}\index{model structure!fstructure@$f$-structure
on $\sK/B$, $\sU/B$ and $\sU_B$} The following statements hold.
\begin{enumerate}[(i)]
\item The $f$-equivalences, $f$-fibrations, and $\bar{f}$-cofibrations
give $\sK/B$ a proper topological model structure. Moreover, a map in
$\sK/B$ is an $\bar{f}$-cofibration if and only if it is a closed
$f$-cofibration.
\item The $f$-equivalences, $f$-fibrations, and $\bar{f}$-cofibrations
give $\sU/B$ a proper topological model structure. Moreover, a map in
$\sU/B$ is an $\bar{f}$-cofibration if and only if it is an $f$-cofibration.
\item  The $fp$-equivalences, $fp$-fibrations, and
$\overline{fp}$-cofibrations
give $\sU_B$ an $fp$-model structure.
\item The $r$-classes of maps give the category $\sU_{\sU}$ of retracts
a proper topological $r$-model structure.
\end{enumerate}
\end{thm}
\begin{proof}
Apart from the factorization axioms, the model structures follow from the
discussion in \ref{sec:towardh}. In particular, the lifting axioms, the
properness, and the topological property of all of these model structures 
are given by \myref{h-structure}.  In (ii), (iii), and (iv), the 
factorization axioms follow from \myref{Cole} since the argument in 
\myref{Umod} verifies \myref{hyp}. The rest of (i) can be proven by direct mimicry of the proof of \myref{hmodelis}, using \myref{fNDR}, and the
characterization of the $\bar{f}$-cofibrations in (ii) follows.
\end{proof}

\begin{rem}\mylabel{fpmodel?}
We do not know whether or not $\sK_B$ is an $fp$-model category
or whether the $\overline{fp}$-cofibrations in $\sK_B$ are 
characterized as the closed $fp$-cofibrations. We also 
do not know whether or not $\sK_{\sU}$ is an $r$-model category.
The problem here is related to the fact that, while the sections 
of ex-spaces are always inclusions, they need not be closed inclusions
unless the total spaces are in $\sU$. Steps 1 and 3 of the proof
of \myref{hmodelis} fail in $\sK_B$, and we also do not see how
to carry over Str{\o}m's original proofs in \cite{Strom2, Strom}.
Theorem \ref{h-structure} still applies, giving much
of the information carried by a model structure.
Observe too that if $i\colon A \rtarr X$ is a map of well-sectioned
ex-spaces over $B$, then $i$ is an $fp$-cofibration if and only if it
is an $f$-cofibration, by \myref{reverse}(iii).  For ex-spaces that
are not well-sectioned, we have little understanding of $fp$-cofibrations,
even when $B$ is a point. We have little understanding of
$\overline{fp}$-cofibrations that are not $\overline{f}$-cofibrations 
in any case.  
\end{rem}

There is a certain tension between the $fp$ and $h$-notions, with the $f$-notions serving as a bridge between the two. Fiberwise pointed homotopy is the intrinsically right notion of homotopy in $\sK_B$, hence the $fp$-structure is the philosopically right classical $h$-type model structure on $\sK_B$, or at least on $\sU_B$.  It is the one that is naturally related to fiber and cofiber sequences, the theory of which works formally in any based topologically bicomplete category in exactly the same way as for based spaces, as we will
recall in \S5.6. A detailed exposition in the case of ex-spaces is given in \cite{CJ, James, James2}.

However, with $h$ replaced by $fp$, we do not have the implications that 
we emphasized in the general philosophy of \S\ref{Sphil}. In particular, 
with the over and under $q$-model structure,
$q$-cofibrations need not be $fp$-cofibrations and $fp$-fibrations need 
not be $q$-fibrations, let alone $h$-fibrations.  The $q$-model structure
is still related to the $h$-model structure as in \S4.1, but this does not serve to relate the $q$-model structure to parametrized fiber and cofiber sequences in the way that we are familiar with in the nonparametrized context.
This already suggests that the $q$-model structure might not be appropriate
in parametrized homotopy theory.  In the following four sections, we explore conceptually what is required of a $q$-type model structure to connect it up with the intrinsic homotopy theory in a topologically bicomplete category.

\section{Well-grounded topologically bicomplete categories}\label{sec:backstr}

Let $\sC$ be a topologically bicomplete category in either the based or the unbased sense; we use the notations of the based context. In our work here, and in other topological contexts, $\sC$ is \emph{topologically concrete} in the sense that there is a faithful and continuous forgetful functor from $\sC$ to spaces.  In practice, appropriate ``ground cofibrations'' can then be specified in terms of underlying spaces. These cofibrations should be thought 
of as helpful background structure in our category $\sC$.

To avoid ambiguity, we use the term ``Hurewicz cofibration'', abbreviated
notationally to \emph{$cyl$-cofibration},\index{cofibration!cyl- --@$cyl$- --}  for the maps that satisfy the HEP with respect to the cylinders in $\sC$. We
also have the notion of a strong Hurewicz cofibration, which we abbreviate
notationally to \emph{$\overline{cyl}$-cofibration}.\index{cofibration!cyl- --@$\overline{cyl}$- --} 
For example, the $cyl$-cofibrations in $\sK$, $\sK/B$, and $\sK_B$ are the $h$-cofibrations, the $f$-cofibrations, and the $fp$-cofibrations, respectively, and similarly for $\overline{cyl}$-cofibrations. As we have seen, it often happens that $cyl$-cofibrations between suitably nice objects of $\sC$, which we shall call ``well-grounded'', are also ground cofibrations.  We introduce language to describe this situation. The following definitions codify the behavior of the well-grounded objects with respect to the $cyl$-cofibrations, colimits, and tensors in $\sC$. It is convenient to build in the appropriate equivariant generalizations of our notions, although we defer a formal discussion of $G$-topologically bicomplete $G$-categories to \S10.2; see \myref{defn:enrichBG}. The examples in \S1.4 give the idea.

\begin{defn}\mylabel{spaceback} An unbased space is {\em well-grounded} 
if it is compactly generated. A based space is {\em well-grounded} 
if it is compactly generated and well-based. The same definitions
apply to $G$-spaces for a topological group $G$.
\end{defn}

Let $\sC$ be a topologically bicomplete category. 

\begin{defn}\mylabel{back}
A full subcategory of $\sC$ is said to be a subcategory of \emph{well-grounded objects}\index{well-grounded!object}\index{category!sub-- of 
well-grounded objects} if the following properties hold.
\begin{enumerate}[(i)]
\item The initial object of $\sC$ is well-grounded.
\item All coproducts of well-grounded objects are well-grounded.
\item If $i\colon X \rtarr Y$ is a $cyl$-cofibration and $f\colon X\rtarr Z$ is any map, where $X$, $Y$, and $Z$ are well-grounded, then the pushout 
$Y\cup_X Z$ is well-grounded.
\item The colimit of a sequence of $cyl$-cofibrations between well-grounded
objects is well-grounded.
\item A retract of a well-grounded object is well-grounded.
\item If $X$ is a well-grounded object and $K$ is a well-grounded space, then 
$X \sma K$ ($X\times K$ in the unbased context) is well-grounded.
\end{enumerate}
When $\sC$ is $G$-topologically bicomplete, we replace spaces by 
$G$-spaces in (vi). 
\end{defn}

\begin{defn}\mylabel{moreback}
A \emph{ground structure}\index{ground structure} on $\sC$ is a (full) subcategory of well-grounded objects together with a subcategory of cofibrations, called the \emph{ground cofibrations}\index{ground cofibrations}\index{cofibration!ground --s} and denoted \emph{$g$-cofibrations}\index{gcofibration!$g$-cofibration}\index{cofibration!g@$g$-}, such that every $cyl$-cofibration between well-grounded objects is a $g$-cofibration.  A map that is both a $g$-cofibration and a $cyl$-cofibration 
is called a \emph{bicofibration}.\index{bicofibration!bicofibration} 
\end{defn}

Thus a $cyl$-cofibration between well-grounded objects is a bicofibration. 
The need for focusing on bicofibrations and the force of the definition 
come from the following fact.

\begin{warn} In practice, (iii) often fails if $i$ is a $g$-cofibration
between well-grounded objects that is {\em not}\, a $cyl$-cofibration, 
as we shall illustrate in \S6.1.  In particular, in $G\sK_B$ with the
canonical ground structure described below, it can already fail for an 
inclusion $i$ of $I$-cell complexes, where $I$ is the standard set of 
generators for the $q$-cofibrations.
\end{warn}

In the next chapter, we will construct a $q$-type model structure for 
$G\sK_B$ with a set of generating cofibrations to which the following
implication of Definitions \ref{cofsub} and \ref{back} applies.

\begin{lem}\mylabel{cofgcof}
Let $I$ be a set of $cyl$-cofibrations between well-grounded objects and let $f\colon X\rtarr Y$ be a retract of a relative $I$-cell complex $W\rtarr Z$. Then $f$ is a bicofibration.  If $W$ is well-grounded, then so are $X$, $Y$,
and $Z$.
\end{lem}

Our categories of equivariant parametrized spaces have canonical ground structures. Recall that the classes of $f$ and $\bar{f}$-cofibrations in 
$G\sU/B$ and $G\sU_B$ coincide. 

\begin{defn}\mylabel{exbackdef}
A space over $B$ is \emph{well-grounded}\index{space!well-grounded} if 
its total space is compactly generated. An ex-space over $B$ is \emph{well-grounded}\index{ex-space!well-grounded} if
it is well-sectioned and its total space is compactly generated. In both
$G\sK/B$ and $G\sK_B$, define the $g$-cofibrations to be the $h$-cofibrations.
\end{defn}

Note that the only distinction between well-sectioned and well-grounded 
ex-spaces is the condition on total spaces.  The distinction is relevant 
when we consider relative $I$-cell complexes $X_0\rtarr X$
in $G\sK_B$. If $X_0$ is well-sectioned, then so is $X$, whereas $X/X_0$ 
is an $I$-cell complex and is therefore well-grounded for any $X_0$.

\begin{prop}\mylabel{exback} These definitions specify 
ground structures on $G\sK/B$ and on $G\sK_B$.
\end{prop}
\begin{proof}
For $G\sK/B$, the Hurewicz cofibrations are the $f$-cofibrations, and these
are $h$-cofibrations. It is standard that $G\sU/B$ has the closure properties 
specified in \myref{back}.  For $G\sK_B$, the Hurewicz cofibrations are the
$fp$-cofibrations.  Between well-sectioned ex-spaces, these are $f$-cofibrations and therefore $h$-cofibrations by \myref{reverse2}(i).  
Parts (i)--(v) of \myref{back} are clear since well-sectioned means 
$\bar{f}$-cofibrant, which is a lifting property. Finally we consider 
part (vi). Recall that $X\sma_B K$ can be constructed as the pushout of
\[\xymatrix{{*}_B & X\amalg (B\times K) \ar[l]\ar[r] & X\times K}\]
in the category of spaces over $B$. By the equivariant version of
the NDR-pair characterization of $f$-cofibrations in \myref{fNDR}, 
these spaces are $f$-cofibrant and the map on the right is an 
$f$-cofibration. This implies that $X\sma_B K$ is $f$-cofibrant.
\end{proof}

\section{Well-grounded categories of weak equivalences}

The following definition describes how the weak equivalences and the ground structure are related in practice.

\begin{defn}\mylabel{hproper}
Let $\sC$ be a topologically bicomplete category with a given ground structure.
A subcategory of weak equivalences in $\sC$ is \emph{well-grounded}\index{well-grounded!weak equivalence}\index{weak-equivalence!well-grounded}
if the following properties hold  (where acyclicity refers to the weak equivalences).
\begin{enumerate}[(i)]
\item A homotopy equivalence is a weak equivalence.
\item A coproduct of weak equivalences between
well-grounded objects is a weak equivalence.
\item (Gluing lemma)\mylabel{glue}\index{Gluing lemma}
Assume that the maps $i$ and $i'$ are bicofibrations and the vertical
arrows are weak equivalences in the following diagram.
\[\xymatrix{Y \ar[d]  & X \ar[d] \ar[l]_-{i} \ar[r]^f & Z \ar[d] \\
Y' &  X' \ar[l]^-{i'}  \ar[r]_{f'}  &  Z'}\]
Then the induced map of pushouts is a weak equivalence.  In particular, 
pushouts of weak equivalences along bicofibrations are weak equivalences.
\item (Colimit lemma) Let $X$ and $Y$ be the colimits of sequences of bicofibrations $i_n\colon X_n\rtarr X_{n+1}$ and $j_n\colon Y_n\rtarr Y_{n+1}$ such that both $X/X_0$ and $Y/Y_0$ are well-grounded. If $f\colon X\rtarr Y$ is the colimit of a sequence of compatible weak equivalences $f_n\colon X_n\rtarr Y_n$, then $f$ is a weak equivalence. In particular, if each $i_n$ is a weak equivalence, then the induced map $i\colon X_0\rtarr X$ is a weak equivalence. 
\item For a map $i\colon X\rtarr Y$ of well-grounded objects in $\sC$ and a 
map $j\colon K\rtarr L$ of well-grounded spaces, $i\Box j$ is a weak 
equivalence if $i$ is a weak equivalence or $j$ is a $q$-equivalence.
\end{enumerate}
\end{defn}

Here, in the based context, $i\Box j$ is the evident induced map
$$ (X\sma L) \cup_{X\sma K} (Y\sma K) \rtarr Y\sma L.$$
The gluing lemma implies that acyclic bicofibrations are preserved 
under push\-outs, as of course holds for pushouts of acyclic cofibrations in
model categories.  The special case mentioned in (iii) corresponds to the
left proper axiom in model categories. As there, it can be used to prove 
the general case of the gluing lemma provided that we have suitable 
factorizations.

\begin{lem}\mylabel{gluederiv} Assume the following hypotheses.
\begin{enumerate}[(i)]
\item Weak equivalences are preserved under pushouts along bicofibrations.
\item Every map factors as the composite of a bicofibration and a weak equivalence.
\end{enumerate}
Then the gluing lemma holds.
\end{lem}
\begin{proof} We use the notations of \myref{hproper}(iii) and proceed in
three cases. 

If $f$ and $f'$ are both weak equivalences, then, by (i), so are the horizontal arrows in the commutative diagram
$$ \xymatrix{
Y \ar[d] \ar[r] &  Y\cup_X Z \ar[d] \\
Y' \ar[r] & Y'\cup_{X'} Z'.\\}$$
Since $Y\rtarr Y'$ is a weak equivalence, the right arrow is a weak equivalence by the two out of three property of weak equivalences.  

If $f$ and $f'$ are both bicofibrations, consider the commutative diagram
\[\xymatrix@=.6cm{
&& X \ar[rr]^-{i} \ar[dll]_{f} \ar@{-}[d]  & & Y \ar[dd] \ar[ddr] \ar[dl] \\
Z \ar[rrr] \ar[dd]  && \ar[d] & Y\cup_XZ \ar[dd] \ar[ddr] \\
&& X' \ar@{-}[r] \ar[dll]_{f'} & \ar[r]|(.5)\hole & Y\cup_X X' \ar[dl]|(.35)\hole \ar[r] & Y' \ar[dl] \\
Z' \ar[rrr] && & Y\cup_X Z' \ar[r] & Y'\cup_{X'} Z'.}\]
The back, front, top, and two bottom squares are pushouts, and the
middle composite $X'\rtarr Y'$ is $i'$.  Since $f$ and $f'$ are 
bicofibrations, so are the remaining three arrows from the back to the front.  Similarly, $i$ and its pushouts are bicofibrations.
Since $X\rtarr X'$, $Y\rtarr Y'$, and $Z\rtarr Z'$ are weak equivalences, 
(i) and the two out of three property imply that 
$Y\rtarr Y\cup_X X'$, $Y\cup_X X'\rtarr Y'$, 
$Y\cup_X Z\rtarr Y\cup_X Z'$, and $Y\cup_X Z'\rtarr Y'\cup_{X'}Z'$ are weak equivalences. Composing the last two, $Y\cup_X Z\rtarr Y'\cup_{X'}Z'$ is a 
weak equivalence.

To prove the general case, construct the following commutative diagram.
\[\xymatrix@=.6cm{
Y \ar[dd] & & X \ar[dd] \ar[ll]_-{i} \ar[rr]^(.35){f} \ar[dr] && Z \ar[dd] \\
& & & W \ar[ur]_{\bar{f}} \ar[dd] & \\
Y' & & \ar[ll]_-{i'} X' \ar@{-}[r]^(.7){f'} \ar[dr]  & \ar[r] & Z' \\
& & & X'\cup_X W \ar[ur]_{\bar{f}'} }\]
Here we first factor $f$ as the composite of a bicofibration and a weak equivalence $\bar{f}$ and then define a map $\bar{f}'$ by the universal property of pushouts. By hypothesis (i), $W\rtarr X'\cup_X W$ is a weak equivalence, and by the two out of three property, so is $\bar{f}'$. By the second case,
$$Y\cup_X W \rtarr Y'\cup_{X'} (X'\cup_X W) \iso Y'\cup_X W $$
is a weak equivalence and by the first case, so is
\[Y\cup_X Z\iso (Y\cup_X W)\cup_W Z \rtarr (Y'\cup_X W)\cup_{(X'\cup_X W)} Z'
\iso Y'\cup_{X'} Z'. \qedhere\]
\end{proof}

\begin{rem} Clearly the previous result applies to any categories of
weak equivalences and cofibrations that satisfy (i) and (ii).  The 
essential point is that, in practice, we often need bicofibrations 
in order to verify (i).  
\end{rem}

Similarly, but more simply, the following observation reduces the verification of \myref{hproper}(v) to special cases. Here we assume that $\sC$ is based.

\begin{lem}\mylabel{boxacy} Let $i\colon X\rtarr Y$ be a map in $\sC$ and
$j\colon K\rtarr L$ be a map of based spaces. Display $i\Box j$ in the diagram 
$$\xymatrix@=.6cm{
X\sma K \ar[rr]^{\text{id}\sma j} \ar[dd]_{i\sma\text{id}} & & X\sma L \ar[dd]^{i\sma\text{id}} \ar[dl]_{k}\\
& (X\sma  L)\cup_{X\sma K}(Y\sma K)\ar[dr]_-{i\Box j} & \\
Y\sma K \ar[rr]_{\text{id}\sma j} \ar[ur]
& & Y\sma L.}$$
If the maps $i\sma\text{id}$ and the pushout $k$ of $i\sma\text{id}$ along $\text{id}\sma j$ are weak equivalences, then so is $i\Box j$, and 
similarly with the roles of $i$ and $j$ reversed.
\end{lem}

Together with \myref{cofgcof}, the notion of a well-grounded category of weak equivalences encodes a variant of \myref{cofcat} that often applies when the latter does not.

\begin{lem}\mylabel{veriiii}
If $J$ is a set of acyclic $cyl$-cofibrations between well-grounded objects, then all relative $J$-cell complexes are weak equivalences.
\end{lem}

\begin{proof} This follows from (ii), (iii), and (iv) of \myref{hproper}, together with the observation that if $X_0\rtarr X$ is a relative $J$-cell complex, then $X/X_0$ is a $J$-cell complex and is therefore well-grounded, 
so that (iv) applies.
\end{proof}

There is an analogous reduction of the problem of determining when a 
functor preserves weak equivalences.

\begin{lem}\mylabel{reducts}
Let $F\colon \sC\rtarr \sD$ be a functor between topologically bicomplete categories that come equipped with subcategories of well-grounded weak equivalences with respect to given ground structures. Let $J$ be a set of acyclic $cyl$-cofibrations between well-grounded objects in $\sC$. Assume 
that $F$ has a continuous right adjoint and that $F$ takes maps in $J$ to 
weak equivalences between well-grounded objects. Then $F$ takes a retract 
of a relative $J$-cell complex to an acyclic map in $\sD$.
\end{lem}

\begin{proof}
The functor $F$ preserves $cyl$-cofibrations since it has a continuous right adjoint and hence $FJ$ consists of acyclic $cyl$-cofibrations between well-grounded objects. The conclusion follows from \myref{veriiii} and the fact 
that left adjoints commute with colimits and therefore the construction of cell complexes.
\end{proof}

Similarly, cell complexes are relevant to the verification of \myref{hproper}(v). Recall that the $cyl$-cofibrations in $\sK_*$ 
are the $fp$-cofibrations, that is, the based cofibrations.

\begin{lem}\mylabel{topmod}
Let $I$ be a set of $cyl$-cofibrations between well-grounded objects of $\sC$ and let $J$ be a set of $fp$-cofibrations between well-based spaces. If $i$ is a retract of a relative $I$-cell complex, $j$ is a retract of a relative $J$-cell complex, and either $I$ or $J$ consists of weak equivalences, then $i\Box j$ is a weak equivalence.
\end{lem}
\begin{proof}
Assume that $I$ consists of weak equivalences; the proof of the other case is symmetric. Since the functor $-\sma K$ commutes with coproducts, pushouts, sequential colimits, and retracts, we can construct $j\sma K$ by first applying $-\sma K$ to the generators, then construct the cell complex, and finally pass to retracts. Since $-\sma K$ preserves $cyl$-cofibrations and well-grounded objects by \myref{back}(vi), it takes maps in $I$ to $cyl$-cofibrations between well-grounded objects. By \myref{veriiii}, the resulting cell complex is acyclic and therefore so also is any retract of it. Thus $j\sma K$ is an acyclic bicofibration. Since such maps are preserved under pushouts, Lemma \ref{boxacy} applies to give the conclusion.
\end{proof}

The following classical example is implicit in the literature.

\begin{prop}
The $q$-equivalences in $G\sK$ are well-grounded with respect to the
ground structure whose well-grounded objects are the compactly
generated spaces and whose $g$-cofibrations are the $h$-cofibrations.
\end{prop}
\begin{proof}  Parts (i), (ii), and, here in the unbased case, (v) 
of \myref{hproper} are clear, and (iv) follows easily from
\myref{little}.  The essential point is the gluing lemma of (iii). By
passage to fixed point spaces, it suffices to prove this nonequivariantly. 
Using the gluing lemma for the proper $h$-model structure on $\sK$, we 
see that $f$ and $f'$ can be replaced by their mapping cylinders.  Then
the induced map of pushouts is the map of double mapping cylinders
induced by the original diagram. This map is equivalent to a map of excisive
triads, and in that case the result is \cite[1.3]{weak}, whose proof is
corrected in \cite{Wit}. 
\end{proof}

\begin{prop}\mylabel{exwellgr}
The $q$-equivalences in $G\sK/B$ and $G\sK_B$ are well-ground\-ed with respect to the ground structures of \myref{exback}. In these cases, one need only assume that the relevant maps in the gluing and colimit lemmas are ground cofibrations (= $h$-cofibrations), not both ground and Hurewicz cofibrations.
\end{prop}
\begin{proof}
We verify this for $G\sK_B$. Part (i) of \myref{hproper} holds since any $fp$-equivalence is a $q$-equivalence and part (iii) follows directly from the gluing lemma in $G\sK$. For part (ii), the total space of $\vee_B X_i$ is the pushout in $G\sK$ of 
\[\xymatrix{{*}_B & \amalg {*}_B \ar[l]\ar[r] & \amalg X_i.}\]
Since the $X_i$ are well-grounded, the map on the right is an $h$-cofibration,
hence (ii) also follows from the gluing lemma in $G\sK$. In part (iv), the relevant quotient in $G\sK_B$ is given by the pushout, $X/\!_BX_0$, of the diagram $*_B \longleftarrow X_0 \rtarr X$. Since $X/\!_BX_0$ is well-grounded, the quotient total space is in $\sU$ and one can apply \myref{little} just as on the space level. Finally consider (v). As in the proof of \myref{exback}(vi), $X\sma_B K$ can be constructed as the pushout of the following diagram of $f$-cofibrant spaces over $B$.
\[\xymatrix{{*}_B & X\amalg (B\times K) \ar[l]\ar[r] & X\times K}\]
The map on the right is an $f$-cofibration. By the gluing lemma in $G\sK$, 
it suffices to observe that $X\times K$ preserves $q$-equivalences in both variables since homotopy groups commute with products.
\end{proof}

\section{Well-grounded compactly generated model structures}\label{sec:wellgr}

Let $\sC$ be a topologically bicomplete category or, equivariantly, a 
$G$-topologically bicomplete $G$-category. In the notion of a ``well-grounded 
model structure'', we formulate the properties that a compactly generated model structure on $\sC$ should have in order to mesh well with the intrinsic $h$-structure on $\sC$ described in \S\ref{sec:towardh}.  When $\sC$ has such a model structure, and when the classical $h$-structure actually is a model structure, the identity functor on $\sC$ is a Quillen left adjoint from the well-grounded model structure to the $h$-model structure.  Thus this notion gives a precise axiomatization for the implementaton of the philosophy that 
we advertised in \S4.1. We begin with a variant of \myref{compgen}.

\begin{thm}\mylabel{Newcompgen} 
Let $\sC$ be a topologically bicomplete category with a ground structure, a subcategory of well-grounded weak equivalences, and compact sets $I$ and $J$ of maps
that satisfy the following conditions.
\begin{enumerate}[(i)]
\item (Acyclicity condition) Every map in $J$ is a weak equivalence.
\item (Compatibility condition) A map has the RLP with respect to $I$ if
and only if it is a weak equivalence and has the RLP
with respect to $J$.
\item Every map in $I$ and $J$ is a $\overline{cyl}$-co\-fib\-ration
between well-grounded objects.
\end{enumerate}
Then $\sC$ is a compactly generated model category with generating sets 
$I$ and $J$ of cofibrations and acyclic cofibrations. Every cofibration 
is a bicofibration and every cofibrant object is well-grounded.  A
pushout of a weak equivalence along a bicofibration is a weak 
equivalence and, in particular, the model structure is left proper.
The model structure is topological or, equivariantly, $G$-topological 
if the following condition holds.
\begin{enumerate}
\item[(iv)] $i\Box j$ is an $I$-cell complex if $i\colon X\rtarr Y$ is a map 
in $I$ and $j\colon K\rtarr L$ is a map of spaces (or $G$-spaces) in $I$. \end{enumerate}
\end{thm}
\begin{proof} By \myref{veriiii}, \myref{compgen} applies to verify the model
axioms. Condition (iii) implies the statements about cofibrations and
cofibrant objects by \myref{cofgcof}, and the gluing lemma implies the statement about pushouts of weak equivalences. In the last statement, the set $I$ of 
generating cofibrations in the relevant category of (based or unbased) 
spaces is as specified in \myref{UrIJ}. By passage to coproducts, pushouts, sequential colimits, and retracts, (iv) implies that $i\Box j$ is a 
cofibration if $i\colon X\rtarr Y$ is a cofibration in $\sC$ and 
$j\colon K\rtarr L$ is a $q$-cofibration of spaces (or $G$-spaces).
Together with \myref{topmod}, this implies that the model structure is topological. 
\end{proof}

We emphasize the difference between the acyclicity conditions stated in
\myref{compgen} and in \myref{Newcompgen}.  In the applications of the 
former, it is the verification of the acyclicity of $J$-cell complexes 
that is problemmatic, but in the latter our axiomatization has built 
in that verification. Similarly, our axiomatization has built in the verification of the acyclicity condition required for the model 
structure to be topological.

\begin{defn}\mylabel{wellmodel} 
A compactly generated model structure on $\sC$ is said to be 
\emph{well-ground\-ed}\index{well-grounded!model category}\index{model category!well-grounded} if it is right proper and satisfies all of the hypotheses of the preceding theorem. It is therefore proper and topological 
or, equivariantly, $G$-topological.
\end{defn}

\section{Properties of well-grounded model categories}

Assume that $\sC$ is a well-grounded model category. To derive properties
of its homotopy category $\text{Ho}\sC$, we must sort out the relationship between homotopies defined in terms of cylinders and homotopies in the model theoretic sense, which we call ``model homotopies''. Of course, the cylinder objects $\text{Cyl}(X)$ in $\sC$ have maps $i_0$, $i_1\colon X\rtarr \text{Cyl}(X)$ and $p\colon \text{Cyl}(X) \rtarr X$, and $i_0$ (or $i_1$) and $p$ are inverse homotopy equivalences since tensors with spaces preserve homotopies in the space variable. \myref{hproper}(i) ensures that $p$ is therefore a weak equivalence.  This means that $\text{Cyl}(X)$ is a model theoretic cylinder object in $\sC$, provided that we adopt the non-standard definition of \cite{DS}. With the language there, it need not be a {\em good} cylinder object since $i_0\amalg i_1\colon X\amalg X\rtarr \text{Cyl}\,(X)$ need not be a cofibration. As pointed out in \cite[p. 90]{DS}, this already fails for spaces, where the inclusion $X\amalg X\rtarr X\times I$ is not a $q$-cofibration unless $X$ is $q$-cofibrant. With the standard definition given
in \cite{Hirschhorn, Hovey, Q}, cylinder objects are required to have this cofibration property.  Under that interpretation, the cylinder objects $\text{Cyl}(X)$ would not qualify as model theoretic cylinder objects
in general.  (We note parenthetically that ``good cylinders''
are defined in \cite{SVogt} in such a way as to include all 
standard cylinders in the category of spaces).  We record the following observations.

\begin{lem}\mylabel{comphty}
Consider maps $f,g\colon X\rtarr Y$ in $\sC$.
\begin{enumerate}[(i)]
\item If $f$ is homotopic to $g$, then $f$ is left model homotopic to $g$. 
\item If $X$ is cofibrant, then $\text{Cyl}(X)$ is a good cylinder object.
\item If $X$ is cofibrant and $Y$ is fibrant, then $f$ is homotopic to $g$ if and only if $f$ is left and right model homotopic to $g$.
\end{enumerate}
\end{lem}
\begin{proof}
Part (i) is \cite[4.6]{DS}, part (ii) follows from  \myref{back}(iii), and part (iii) follows from \cite[4.23]{DS}.
\end{proof}

Let $[X,Y]$\noteindex{XY@$[X,Y]$} denote the set of morphisms $X\rtarr Y$ in $\text{Ho} \sC$ and let $\pi(X,Y)$\noteindex{piXY@$\pi(X,Y)$} denote the set of homotopy classes of maps $X\rtarr Y$. We shall only use the latter notation when homotopy and model homotopy coincide.

\begin{lem}[Cofiber sequence lemma]\mylabel{modelcofiber}\index{cofiber sequence}
Assume that $\sC$ is based. Consider the cofiber sequence
\[X\rtarr Y \rtarr Cf \rtarr \Sigma X \rtarr \Sigma Y \rtarr \Sigma Cf \rtarr \Sigma^2X \rtarr \cdots\]
of a well-grounded map $f\colon X\rtarr Y$. For any object $Z$, the 
induced sequence
\[\cdots \rtarr [\Sigma^{n+1} X, Z] \rtarr [\Sigma^n Cf, Z] \rtarr [\Sigma^n Y, Z] \rtarr [\Sigma^n X, Z] \rtarr \cdots \rtarr [X, Z]\]
of pointed sets (groups left of $[\Sigma X,Z]$, Abelian groups left of $[\Sigma^2 X, Z]$) is exact.
\end{lem}
\begin{proof} As usual, giving $I$ the basepoint $1$, we define
$$CX = X\sma I, \qquad \SI X = X\sma S^1, \qquad \text{and} \qquad Cf = Y\cup_f CX.$$
If $X$ is cofibrant, then $X$ is well-grounded and $X\rtarr CX$ is a cofibration and therefore a bicofibration. If $X$ and $Y$ are cofibrant, then so is $Cf$, as one sees by solving the relevant lifting problem by first using that $Y$ is cofibrant, then using that $X\rtarr CX$ is a cofibration, and finally using that $Cf$ is a pushout. Thus, taking $Z$ to be fibrant, the conclusion follows in this case from the sequence of homotopy classes of maps
\[\cdots \rtarr \pi(\Sigma X, Z)\rtarr \pi(Cf, Z) \rtarr \pi(Y, Z) \rtarr \pi(X, Z),\]
which is proven to be exact in the same way as on the space level. If $X$ and $Y$ are not cofibrant, let $Qf\colon QX\rtarr QY$ be a cofibrant approximation to $f$. The gluing lemma applies to give that the canonical map $CQf\rtarr Cf$ is a weak equivalence. Therefore the conclusion follows in general from the special case of cofibrant objects.
\end{proof}

\begin{warn}
While the proof just given is very simple, it hides substantial subtleties.  
It is crucial that cofibrant objects $X$ be well-grounded, so that the $cyl$-cofibration $X\rtarr CX$ is a bicofibration and the gluing lemma applies.
\end{warn}

Of course, the group structures are defined just as classically. 
The pinch maps 
\[S^1\cong I/\{0,1\} \rtarr I/\{0,\tfrac12,1\}\cong S^1\vee S^1
\qquad\text{and}\qquad
I \rtarr I/\{\tfrac12,1\}\cong I\vee S^1\]
induce pinch maps 
\[\Sigma X \rtarr \Sigma X \vee \Sigma X \qquad\text{and}\qquad Cf\rtarr Cf\vee \Sigma X\]
that give $\Sigma X$ the structure of a cogroup object in $\text{Ho}\sC$ and $Cf$ a coaction by $\Sigma X$; $\Sigma^2 X$ is an abelian cogroup object for the same reason that higher homotopy groups are abelian. Therefore $[\Sigma X,Z]$ is a group, $[Cf,Z]$ is a $[\Sigma X,Z]$-set, and $[\Sigma X,Z]\rtarr [Cf,Z]$ is a $[\Sigma X,Z]$-map.  

\begin{lem}[Wedge lemma]\mylabel{modelwedges}\index{wedge lemma}
For any $X_i$ and $Y$ in $\sC$, $[\amalg X_i, Y]\iso \PI [X_i,Y]$. 
\end{lem}
\begin{proof} This is standard, using that a coproduct of cofibrant
approximations is a cofibrant approximation.
\end{proof}

\begin{lem}[$\text{Lim}^1$ lemma]\mylabel{modellim1}\index{lim1lemma @$\text{lim}^1$ lemma}
Assume that $\sC$ is based.
Let $X$ be the colimit of a sequence of well-ground\-ed $cyl$-cofibrations $i_n\colon  X_n \rtarr X_{n+1}$. Then, for any object $Y$, there is a ${\text{lim}}^1$ exact sequence of pointed sets
$$* \rtarr {\text{lim}}^1\,[\SI X_n,Y] \rtarr [X,Y] \rtarr \text{lim}\,[X_n,Y]\rtarr *.$$
\end{lem}
\begin{proof}
The telescope $\text{Tel}\,X_n$ is defined to be $\text{colim}\,T_n$, 
where the $T_n$ and a ladder of weak equivalences 
$j_n\colon X_n\rtarr T_n$ and $r_n\colon T_n\rtarr X_n$ 
are constructed inductively by setting $T_0=X_0$ and letting
$j_{n+1}$ and  $r_{n+1}$ be the maps of pushouts induced by 
the following diagram.
\[\xymatrix{
X_n\ar[d]_{i_1} \ar@{=}[r] & X_n \ar[d]^{\nu_2} \ar[r]^-{i_n} & X_{n+1}\ar[d]^{\nu_2}\\
\text{Cyl}\,X_n\ar[d]_{p} & X_n\amalg X_n \ar[l]_{i_{(0,1)}}\ar[r]^{j_n\amalg i_n}\ar@{=}[d] & T_n \amalg X_{n+1}\ar[d]^{r_n\amalg\text{id}}\\
X_n & X_n\amalg X_n \ar[l]^-{\nabla}\ar[r]_-{\text{id}\amalg i_n} & X_n \amalg X_{n+1}}\]
Since $j_{n+1}$ is a pushout of the bicofibration 
$i_1\colon X_n\rtarr \text{Cyl}(X_n)$, the gluing lemma and colimit lemma
specified in \myref{hproper}(iii) and (iv) apply to show that the induced 
maps $\text{Tel}\, X_n \rtarr \text{colim}\, X_n = X$ are weak equivalences.

As in the cofiber sequence lemma, we can use cofibrant approximation to
reduce to a question about $\pi(-,-)$. Then the telescope admits an alternative description from which the ${\text{lim}}^1$ exact sequence is immediate. It would take us too far afield to go into full details of what should be a standard argument, but we give a sketch since we cannot find our preferred 
argument in the literature.  

Recall that the classical homotopy pushout, or double mapping cylinder, of 
\[\xymatrix{Y & X \ar[l]_-f\ar[r]^-{f'} & Y'}\] 
is the ordinary pushout $M(f,f')$ of 
\[\xymatrix{\text{Cyl}\, X & X\amalg X \ar[l]_-{i_{0,1}} \ar[r]^{f\amalg f'} & Y\amalg Y'.}\]
It fits into a cofiber sequence 
\[Y\amalg Y'\rtarr M(f,g)\rtarr \SI X.\]
There results a surjection from $\pi(M(f,g),Z)$ to the evident pullback,
the kernel of which is the set of orbits of the right action
of $\pi(\SI Y,Z)\times \pi(\SI Y', Z)$ on $\pi(\SI X, Z)$ given by
$x(y,y') = (\SI f)^*(y)^{-1}x(\SI f')^*(y')$.

The classical homotopy coequalizer $C(f,g)$ of parallel maps $f,g\colon X\rtarr Y$ is the homotopy pushout of the coproduct $f\amalg g\colon  X\amalg X\rtarr Y\amalg Y$ and the codiagonal $\nabla\colon X\amalg X \rtarr X$.  Using a little algebra, we see that $\pi(C(f,g),Z)$ maps surjectively to the equalizer of $f^*$ and $g^*$ with kernel isomorphic to the set of orbits of $\pi(\SI X,Z)$ under the right action of $\pi(\SI Y,Z)$
specified by $xy = (\SI f)^*(y)^{-1}x(\SI g)^*(y)$. 

In this language, $\text{Tel}\,X_n$ is the classical homotopy coequalizer of the identity and the coproduct of the $i_n$, both being self maps of the coproduct of the $X_n$.  By algebraic inspection, the $\lim^{1}$ exact 
sequence follows directly.  A quicker, less conceptual, argument is 
possible, as in \cite[p.\,146]{Concise} for example. 
\end{proof} 

\begin{rem}
Let $\sC$ be an arbitrary pointed model category with (for simplicity) a functorial cylinder construction $\text{Cyl}$. If $X$ is cofibrant, 
let $\Sigma X$ denote the quotient $\text{Cyl}(X)/(X\vee X)$. Quillen \cite{Q} constructed a natural cogroup structure on $\Sigma X$ in $\text{Ho}\sC$. For
a cofibration $X\rtarr Y$ between cofibrant objects, he also constructed a natural coaction of $\Sigma X$ on the quotient $Y/X$. One can then define cofiber sequences in $\text{Ho}\sC$ just as in the homotopy category of a topological model category, and one can define fiber sequences dually.

The cofiber sequences and fiber sequences each give $\text{Ho}\sC$ a 
suitably weakened form of the notion of a triangulation, called a 
``pretriangulation'' \cite{Hovey, Q}, and they are suitably compatible. 
If $\text{Ho}\sC$ is closed symmetric monoidal one can take this a step further and formulate what it means for the pretriangulation to be compatible with that structure, as was done in \cite{Tri} for triangulated categories. However, proving the compatibility axioms from this general  
point of view would at best be exceedingly laborious, if it could be done 
at all. 

These purely model theoretic constructions of the suspension and looping functors $\Sigma$ and $\Omega$ are more closely related to the familiar topological constructions than might appear. The homotopy category of any model category is enriched and bitensored over the homotopy category of spaces (obtained from the $q$-model structure) \cite{DK, Hovey}, and the suspension and loop functors are given by the (derived) tensor and cotensor with the unit circle. That is, $\Sigma X \simeq X\sma S^1$ and $\Omega X \simeq F(S^1,X)$.

This general point of view is not one that we wish to emphasize. For topological model categories, the structure described in this section is far easier to define and work with directly, as in classical homotopy theory, and we have
axiomatized what is required of a model structure in order to allow the use of such standard and elementary classical methods.  In our topological context, the homotopy category $\text{Ho}\sC$ is automatically enriched over $\text{Ho}\sK_*$ and $(\Sigma, \Omega)$ is a Quillen adjoint pair that descends to an adjoint pair on homotopy categories that agrees with the purely model theoretic adjoint pair just described.

The crucial point for our stable work is that a large part of this structure exists \emph{before} one constructs the desired model structure. It can therefore be used as a tool for carrying out that construction. This is 
in fact how stable model categories were constructed in \cite{EKMM, MM, MMSS},
but there the compatibility between $q$-type and $h$-type structures was too
evident to require much comment.  The key step in our construction of the stable model structure on parametrized spectra in Chapter 12 is to show that cofiber sequences induce long exact sequences on stable homotopy groups.  That will allow us to verify that the stable equivalences are suitably well-grounded,
and from there the model axioms follow as in the earlier work just cited.
\end{rem}
\chapter{The $qf$-model structure on $\sK_B$}

\section*{Introduction}

In this chapter, we introduce and develop our preferred $q$-type model structure, namely the $qf$-model structure. It is a Quillen equivalent 
variant of the $q$-model structure that has fewer, and better structured, cofibrations. For clarity of exposition, we work nonequivariantly 
in this chapter, which is taken from \cite{Sig}. 

We begin by comparing the homotopy theory of spaces and the homotopy theory of ex-spaces over $B$, starting with a comparison of the $q$-model structures that we have on both. In the category $\sK$ of spaces, we have the familiar situation described in \S4.1. The homotopy category $\text{Ho}\sK$ that we care about is defined in terms of $q$-equivalences, the intrinsic notion of homotopy is given by the classical cylinders, and, since all spaces are $q$-fibrant, the category $\text{Ho}\sK$ is equivalent to the classical homotopy category $h\sK_c$ of $q$-cofibrant spaces (or CW complexes). Since the $q$-cofibrations are $h$-cofibrations, the $q$-model structure and the $h$-model structure on $\sK$ mesh smoothly.  Indeed, the classical and model theoretic homotopy theory have been used in tandem for so long that this meshing of structures goes without notice.  In particular, although cofiber and fiber sequences are defined in terms of the $h$-model structure while the homotopy category is defined in terms of the $q$-model structure, the compatibility seems automatic.

Now consider the category $\sK_B$.  The homotopy category $\text{Ho}\sK_B$ that we care about is defined in terms of $q$-equivalences of total spaces, but we need some justification for making that statement.  A map of $q$-fibrant ex-spaces is a $q$-equivalence of total spaces if and only if all of its maps on fibers are $q$-equivalences. This reformulation captures the idea that the homotopical information in parametrized homotopy theory should be encoded on the fibers, and it is such fiberwise $q$-equivalences that we really care about. It is only for $q$-fibrant ex-spaces, or ex-spaces whose projections are at least quasifibrations, that the homotopy groups of total spaces give the ``right answer''. There are three notions of homotopy in sight, $h$, $f$, and $fp$. The last of these is the intrinsic one defined in terms of the relevant cylinders in $\sK_B$, and $\text{Ho}\sK_B$ is equivalent to the classical homotopy category $h{\sK_{B}}_{cf}$ of $q$-cofibrant and $q$-fibrant objects, defined with respect to $fp$-homotopy. It is still true that $q$-cofibrations are $h$-cofibrations. However, it is {\em not} true that $q$-cofibrations are $fp$-cofibrations, and it is the latter that are intrinsic to cofiber sequences. The classical and model theoretic homotopy theory no longer mesh.

Succinctly, the problem is that the $q$-model structure is not an example
of a well-grounded compactly generated model category. The task that lies 
before us is to find a model structure which does satisfy 
the axioms that we set out in \S5.5 and therefore can be used in tandem with 
the $fp$-structure to do parametrized homotopy theory.  Before embarking on this, we point out the limitations of the $q$-model structure more explicitly 
in \S\ref{sec:danger}.  There are two kinds of problems, those that we are focusing on in our development of the model category theory, and the more intrinsic ones that account for \myref{noway} and which cannot be overcome 
model theoretically. 

Ideally, to define the $qf$-model structure, we would like to take the 
$qf$-co\-fi\-bra\-tions to be those $q$-cofibrations that are also $f$-cofibrations.  However, with that choice, we would not know how to prove 
the model category axioms. We get closer if we try to take as generating 
sets of cofibrations and acyclic cofibrations those generators in the 
$q$-model structure that are $f$-cofibrations, but with that choice we 
still would not be able to prove the compatibility condition \myref{Newcompgen}(ii).  However, using this generating set of cofibrations 
and a subtler choice of a generating set of acyclic cofibrations,  we obtain 
a precise enough homotopical relationship to the $q$-equivalences that we 
can prove the cited compatibility. The construction of the $qf$-model 
structure is given in \S\ref{sec:qfstr}, but all proofs are deferred 
to the following three sections.  

\section{Some of the dangers in the parametrized world}\label{sec:danger}

We introduce notation for the generating (acyclic) cofibrations for the 
$q$-model structures on $\sK/B$ and $\sK_B$. These maps are identified in \myref{cg}, starting from the sets $I$ and $J$ in $\sK$ specified in \myref{UrIJ}.  We then make some comments about these maps that 
help explain the structure of our theory.  

\begin{defn}\mylabel{IJB}  For maps $i\colon  C\rtarr D$ and $d\colon D\rtarr B$ of (unbased) spaces, we have the restriction $d\com i\colon  C\rtarr B$ and may view $i$ as a map over $B$. We agree to write $i(d)$ for either the map $i$ viewed as a map over $B$ or the map $i\amalg\text{id}\colon C\amalg B\rtarr D\amalg B$ of ex-spaces over $B$ that is obtained by taking the coproduct with $B$ to adjoin a section.  In either $\sK/B$ or $\sK_B$, define $I_B$ to be the set of all such maps $i(d)$ with $i\in I$, and define $J_B$ to be the set of all such maps $j(d)$ with $j\in J$. Observe that in $\sK_B$, each map in $J_B$ is the inclusion of a deformation retract of spaces under, but not over, $B$.
\end{defn}

\begin{warn}  
We cannot restrict the maps $d$ to be open here. That is one of the 
reasons we chose $\sK_B$ over $\sO_*(B)$ in \S1.3.
\end{warn}

\begin{warn}\mylabel{quien}
The maps in $I_B$ and $J_B$ are clearly not $f$-cofibrations, only $h$-cofibrations. Looking at the NDR-pair characterization of $f$-cofibrations 
given in \myref{fNDR}, we see that, with our arbitrary projections $d$, 
there is in general no way to carry out the required deformation over $B$. 
Since the maps in $I_B$ and $J_B$ are maps between well-sectioned spaces, 
they cannot be $fp$-cofibrations in general, by \myref{reverse2}(i). 
\end{warn}

\begin{rem}\mylabel{mildhelp} Observe that the maps $i$ in $I_B$ or $J_B$ 
are closed inclusions in $\sU$, so that those maps in $I_B$ or $J_B$ which 
are $f$-cofibrations are necessarily $\bar{f}$-cofibrations and 
therefore both $\bar{fp}$-cofibrations and $\bar{h}$-cofibrations, 
by \myref{compare} and  \myref{ffpmodel}. 
\end{rem}

\myref{quien} shows that the $q$-model structure is \emph{not} well-grounded since its generating (acyclic) cofibrations are not ${fp}$-cofibrations. This may sound like a minor technicality, but that is far from the case. We record 
an elementary example.

\begin{ouch0}  Let $B = I$ and define an ex-map $i\colon X\rtarr Y$ over $I$ by letting $X = \{0\}\amalg I$, $Y = I \amalg I$, and $i$ be the inclusion.  The second copies of $I$ give the sections, and the projections are given by the identity map on each copy of $I$. This is a typical generating acyclic $q$-cofibration, and it is not an $fp$-cofibration. Let $Z$ be the pushout of $i$ and $p\colon X\rtarr I$, where the latter is viewed as a map of ex-spaces over $I$. Then  $Z$ is the one-point union $I\vee I$ obtained by identifying the points $0$. The section $I\rtarr Z$ is not an $f$-cofibration, so that $Z$ is not well-sectioned. The same is true if we replace $Y$ by $Y'=\{1/(n+1)\mid n\in\bN\}\amalg I$ and obtain $Z'$. The map $Z'\rtarr C_I Z'$ of $Z'$ into its cone over $I$ is not an $h$-cofibration (and therefore not a $q$-cofibration).
\end{ouch0}

Thus we cannot apply the classical gluing lemma to develop cofiber sequences, as we did in \S5.6.  This and related problems prevent use of the $q$-model structure in a rigorous development of parametrized stable homotopy theory. For example, consider $q$-fibrant approximation.  If we have a map $f\colon X\rtarr Y$ with $q$-fibrant approximation $Rf\colon RX\rtarr RY$, there is no reason to believe that $C_BRf$ is $q$-equivalent to $RC_Bf$.

We are about to overcome model-theoretically the problems pointed out in 
the warnings above.  Turning to the intrinsic problems that must hold in any 
$q$-type model structure, we explain why the base change functor $f^*$ and 
the internal smash product cannot be Quillen left adjoints.

\begin{warn}
If $f\colon A\rtarr B$ is a map and $d\colon D\rtarr B$ is a disk over $B$,
we have no homotopical control over the pullback $A\times_B D \rtarr A$ in
general.
\end{warn}

\begin{warn}\mylabel{ouchtoo} 
In sharp contrast to the nonparametrized case, the generating sets do not behave well with respect to internal smash products, although they do behave well with respect to external smash products. We have
$$(D\amalg A)\barwedge (E\amalg B) \iso (D\times E)\amalg (A\times B).$$
If the projections of $D$ and $E$ are $d$ and $e$, then the projection of $D\times E$ is $d\times e$. However, if $A=B$, then
$$ (D\amalg B)\sma_B (E\amalg B) \iso (d\times e)^{-1}(\DE B)\amalg (A\times B).$$
We have no homotopical control over the space $(d\times e)^{-1}(\DE B)$ in general.
\end{warn}

This has the unfortunate consequence that, when we go on to parametrized spectra in Part III, we will not be able to develop a homotopically well-behaved theory of point-set level parametrized ring spectra. However, we will be able to develop a satisfactory point-set level theory of parametrized module spectra 
over nonparametrized ring spectra.

\section{The $qf$ model structure on the category $\sK/B$}\label{sec:qfstr}

Rather than start with a model structure on $\sK$ to obtain model structures on $\sK/B$ and $\sK_B$, we can start with a model structure on $\sK/B$ and then apply \myref{under} to obtain a model structure on $\sK_B$. This gives us the opportunity to restrict the classes of generating (acyclic) cofibrations present in the $q$-model structure on $\sK/B$ to ones that are $f$-cofibrations, retaining enough of them that we do not lose homotopical information. This has the effect that the generating (acyclic) cofibrations are $f$-cofibrations between well-grounded spaces over $B$, as is required of a well-grounded model structure.  Such maps have closed images, hence are $\bar{f}$-cofibrations, and therefore all of the cofibrations in the resulting model structure on $\sK/B$ are $\bar{f}$-cofibrations.

We call the resulting model structure the ``$qf$-model structure'', where $f$ refers to the fiberwise cofibrations that are used and $q$ refers to the weak equivalences. The latter are the same as in the $q$-model structure, namely the weak equivalences on total spaces, or $q$-equivalences. This model structure restores us to a situation in which the philosophy advertised in \S\ref{Sphil} applies, with the $q$ and $h$-model structures on spaces replaced by the $qf$ and $f$-model structures on spaces over $B$. Since $f$-cofibrations in $\sK_B$ are $fp$-cofibrations, by \myref{compare}, the philosophy also applies to the $qf$ and $fp$-model structures on $\sK_B$, or at least on $\sU_B$ (see \myref{ffpmodel} and \myref{fpmodel?}). 

We need some notations and recollections in order to describe the generating (acyclic) $qf$-cofibrations and the $qf$-fibrations.

\begin{notn}\mylabel{euclid}
For each $n\geq 1$, embed $\bR^{n-1}$ in $\bR^n = \bR^{n-1}\times \bR$ by sending $x$ to $(x,0)$. Let $e_n=(0,1)\in \bR^n$. For
$n\geq 0$, define the following subspaces of $\bR^n$.
\begin{alignat*}{2}
\bR^n_+ &= \{(x,t)\in\bR^n \mid t\geq 0\} &
\bR^n_- &= \{(x,t)\in\bR^n \mid t\leq 0\}\\
D^n &= \{(x,t)\in\bR^n \mid |x|^2+t^2\leq 1\} &\qquad
S^{n-1} &= \{(x,t)\in\bR^n \mid |x|^2+t^2=1\}\\
S^{n-1}_+ &= S^{n-1}\cap \bR^n_+ &
S^{n-1}_- &= S^{n-1}\cap \bR^n_-
\end{alignat*}
Here $\bR^{0} = \{0\}$ and $S^{-1}=\emptyset$.  We think of $S^n\subset \bR^{n+1}$ as having equator $S^{n-1}$, upper hemisphere $S^n_+$ with north pole $e_{n+1}$ and lower hemisphere $S^n_-$.
\end{notn}

We recall a characterization of Serre fibrations.

\begin{prop}
The following conditions on a map $p\colon E\rtarr Y$ in $\sK$ are equivalent; $p$ 
is called a Serre fibration, or $q$-fibration, if they are satisfied.
\begin{enumerate}[(i)]
\item The map $p$ satisfies the covering homotopy property with respect to disks $D^n$; that is, there is a lift in the diagram
\[\xymatrix{D^n \ar[r]^\alpha\ar[d] & E \ar[d]^p\\
D^n\times I \ar[r]_-h\ar@{-->}[ur] & Y.}\]
\item If $h$ is a homotopy relative to the boundary $S^{n-1}$ in the diagram above, then there is a lift that is a homotopy relative to the boundary.
\item The map $p$ has the RLP with respect to the inclusion $S^n_+\rtarr D^{n+1}$ of the upper hemisphere into the boundary $S^n$ of $D^{n+1}$; that is, 
there is a lift in the diagram
\[\xymatrix{S^n_+ \ar[r]^\alpha\ar[d] & E \ar[d]^p\\
D^{n+1} \ar[r]_-{\bar h}\ar@{-->}[ur] & Y.}\]
\end{enumerate}
\end{prop}
\begin{proof}
Serre fibrations $p\colon E\rtarr Y$ are usually characterized by the first condition. Since the pairs $(D^n\times I, D^n)$ and 
$(D^n\times I, D^n \cup (S^{n-1}\times I))$ are homeomorphic, one easily obtains that the first condition implies the second. Similarly a homeomorphism of the pairs $(D^{n+1},S^n_+)$ and $(D^n\times I, D^n)$ gives that the first and third conditions are equivalent. A homotopy $h\colon D^n\times I\rtarr Y$ relative to the boundary $S^{n-1}$ factors through the quotient map 
$D^n\times I\rtarr D^{n+1}$ that sends $(x,t)$ to $(x, (2t-1)\sqrt{1-|x|^2})$.
Conversely, any map $\bar{h}\colon D^{n+1} \rtarr Y$ gives rise to a homotopy $h\colon D^n\times I\rtarr Y$ relative to the boundary $S^{n-1}$. It follows that the second condition implies the third.
\end{proof}

Property (ii) states that Serre fibrations are the maps that satisfy the ``disk lifting property''\index{disk lifting property} and that is the way we shall think about the $qf$-fibrations. In view of property (iii), we sometimes abuse language 
by calling a map $h\colon D^{n+1}\rtarr Y$ a disk homotopy. The restriction to the upper hemisphere $S^n_+$ gives the ``initial disk'' and the restriction to the lower hemisphere $S^n_-$ gives the ``terminal disk''.

\begin{defn}
A disk $d\colon D^n\rtarr B$ in $\sK/B$ is said to be an \emph{$f$-disk} if 
$i(d)\colon S^{n-1}\rtarr D^n$ is an $f$-cofibration. An $f$-disk
$d\colon D^{n+1}\rtarr B$ is said to be a \emph{relative $f$-disk} if 
the lower hemisphere $S^n_-$ is also an $f$-disk, so that the restriction $i(d)\colon S^{n-1}\rtarr S^n_-$ is an $f$-cofibration; the upper hemisphere  $i(d)\colon S^{n-1}\rtarr S^n_+$ need not be an $f$-cofibration. 
\end{defn}

\begin{defn}\mylabel{IJBf}
Define $I^f_B$\noteindex{IBf@$I^f_B$} to be the set of inclusions $i(d)\colon S^{n-1} \rtarr D^n$ in $\sK/B$, where $d\colon D^n\rtarr B$ is an $f$-disk. Define $J^f_B$\noteindex{JBf@$J^f_B$} to be the set of inclusions $i(d)\colon S^n_+\rtarr D^{n+1}$ of the upper hemisphere into a relative $f$-disk $d\colon D^{n+1}\rtarr B$; note that these initial disks are not assumed to be $f$-disks. A map in $\sK/B$ is said to be
\begin{enumerate}[(i)]
\item a \emph{$qf$-fibration}\index{fibration!qf-@$qf$- --} if it has the RLP with respect to $J^f_B$ and
\item a \emph{$qf$-cofibration}\index{cofibration!qf-@$qf$- --} if it has the LLP with respect to all $q$-acyclic $qf$-fibrations, that is, with respect to those $qf$-fibrations that are $q$-equivalences.
\end{enumerate}
Note that $J^f_B$ consists of relative $I^f_B$-cell complexes and that a map is a $qf$-fibration if and only if it has the ``relative $f$-disk lifting property.'' \end{defn}

With these definitions in place, we have the following theorem. Recall
the definition of a well-grounded model category from \myref{wellmodel}
and recall from Propositions \ref{exback} and \ref{exwellgr} that we have ground structures on $\sK/B$ and $\sK_B$ with respect to which the $q$-equivalences
are well-grounded. Also recall the definition of a quasifibration from \myref{quasifib}.

\begin{thm}\mylabel{Thesis}\index{model structure!qf@$qf$- --}
The category $\sK/B$ of spaces over $B$ is a well-grounded model category with respect to the $q$-equivalences, $qf$-fibrations and $qf$-cofibrations. The sets $I^f_B$ and $J^f_B$ are the generating $qf$-cofibrations and the generating acyclic $qf$-cofibrations. All $qf$-cofibrations are also $\bar{f}$-cofibrations and all $qf$-fibrations are quasifibrations.
\end{thm}

Using \myref{under} and \myref{cg}, we obtain the $qf$-model structure on $\sK_B$. We define a $qf$-fibration in $\sK_B$ to be a map which is a $qf$-fibration when regarded as a map in $\sK/B$, and similarly for 
$qf$-cofibrations.

\begin{thm}
The category $\sK_B$ of ex-spaces over $B$ is a well-grounded model category with respect to the $q$-equivalences, $qf$-fibrations, and $qf$-cofibrations. The sets $I^f_B$ and $J^f_B$ of generating $qf$-cofibrations and generating acyclic $qf$-cofibrations are obtained by adjoining disjoint sections to the corresponding sets of maps in $\sK/B$. All $qf$-cofibrations are $\bar{f}$-cofibrations and all $qf$-fibrations are quasifibrations.
\end{thm}

Since the $qf$-model structures are well-grounded, they are in particular proper and topological. Furthermore, the $qf$-cofibrant spaces over $B$ are well-grounded and the $qf$-fibrant spaces over $B$ are quasifibrant. Since $qf$-cofibrations are $q$-cofibrations, we have an obvious comparison.

\begin{thm}
The identity functor is a left Quillen equivalence from $\sK/B$ with the $qf$-model structure to $\sK/B$ with the $q$-model structure, and similarly for
the identity functor on $\sK_B$.
\end{thm}

We state and prove two technical lemmas in \S\ref{sec:qflemmas}, prove that $\sK/B$ is a compactly generated model category in \S\ref{sec:qfproof}, and prove that the $qf$-fibrations are quasifibrations and the model structure is right proper in \S\ref{sec:qqfuasi}. The $\Box$-product condition of \myref{Newcompgen}(iv) follows as usual by inspection of what happens on generating (acyclic) cofibrations and, as in the case $A=*$ of \myref{ouchtoo},
the projections cause no problems here. 

\section{Statements and proofs of the thickening lemmas}\label{sec:qflemmas}

We need two technical ``thickening lemmas''.  They encapsulate the 
idea that no
information about homotopy groups is lost if we restrict from the
general disks and cells used in the $q$-model structure to the $f$-disks
and $f$-cells that we use in the $qf$-model structure.

\begin{lem}\mylabel{thicken1}
Let $(S^m,q)$ be a sphere over $B$. Then there is an $h$-equivalence $\mu\colon (S^m,\bar{q})\rtarr (S^m,q)$ in $\sK/B$ such that $(S^m,\bar{q})$ is an $I^f_B$-cell complex with two cells in each dimension.
\end{lem}

\begin{lem}\mylabel{thicken2}
Let $(D^n,q)$ be a disk over $B$. Then there is an $h$-equivalence $\nu\colon (D^n, \bar{q})\rtarr (D^n, q)$ relative to the upper hemisphere 
$S^{n-1}_+$ such that $(D^n,\bar{q})$ is a relative $f$-disk.
\end{lem}


The rest of the section is devoted to the proofs of these lemmas.
The reader may prefer to skip ahead to \S\ref{sec:qfproof} to see 
how they are used to prove \myref{Thesis}.

\begin{proof}[Proof of \myref{thicken1}]
To define the map $\mu\colon (S^m,\bar{q})\rtarr (S^m,q)$, we begin by defining some auxiliary maps for each natural number $n\leq m$. They will in fact be continuous families of maps, defined for each $s\in [\tfrac12, 1]$. The parameter $s$ will show that $\mu$ is an $h$-equivalence.

First we define the map
\[\phi^n_+\colon D^n\cap \bR^n_+\rtarr A_s\cup s\cdot S^{n-1}_+\]
from the upper half of the disk $D^n$ to the union of the equatorial annulus 
\[A_s=\overline{D^{n-1}-s\cdot D^{n-1}}=\{(x,0)\in \bR^{n}\colon  s\leq |x|\leq 1\}\]
and the upper hemisphere 
\[s\cdot S^{n-1}_+=\{(x,t)\in \bR^n\colon  \text{$t\geq 0$ and $|(x,t)|=s$}\}\]
to be the projection from the south pole $-e_n$. Similarly, we define
\[ \phi^n_-\colon D^n\cap \bR^n_- \rtarr A_s\cup s\cdot S^{n-1}_-\]
to be the projection from the north pole $e_n$.  The map $\phi^n_+$ is drawn schematically in the following picture. Each point in the upper half of the larger disk lies on a unique ray from $-e_n$. The map $\phi^n_+$ sends it to the intersection of that ray with $A_s\cup s\cdot S^{n-1}$; two such points of
intersection are marked with dots in the picture.
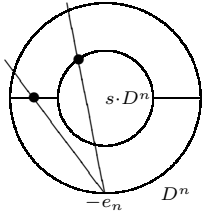
\[\begin{xy}
0;/r3pc/:(1,0),{\ellipse<>{-}}
,(.5,0);(1,0),{\ellipse<>{-}}
,(0,0);(.5,0) **@{-}
,(1.5,0);(2,0) **@{-}
,(1,-1);(.6,1) **@{-}
,(.72,.4) *{\bullet}
,(1,-1);(-.05,.4) **@{-}
,(.25,0) *{\bullet} 
,(1.25,0)*{\sb{s\cdot D^n}}
,(1,-1.1)*{\sb{-e_n}}
,(1.75,-1)*{\sb{D^n}}
\end{xy}\]

Next we use the maps $\phi^n_\pm$ to define a continuous family of maps $f^n_s\colon  D^n\rtarr D^n$ for $s\in[\tfrac12,1]$ by induction on $n$. We let $f^0_s\colon D^0\rtarr D^0$ be the unique map and we define $f^1_s\colon D^1\rtarr D^1$ 
by 
\[f^1_s(t) = \begin{cases}
t/s & \text{if $|t|\leq s$},\\
1 & \text{if $t\geq s$},\\
-1 & \text{if $t\leq -s$};
\end{cases}\]
it maps $[-s,s]$ homeomorphically to $[-1,1]$. We define
$f^n_s\colon D^n\rtarr D^n$ by
\[f^n_s(x,t)=\begin{cases}
s^{-1}\cdot (x,t) & \text{if $|(x,t)|\leq s$},\\
s^{-1}\cdot \phi^n_+(x,t) & \text{if $|(x,t)|\geq s$, $t\geq 0$ and $|\phi^n_+(x,t)|= s$},\\
f^{n-1}_s(\phi^n_+(x,t)) & \text{if $|(x,t)|\geq s$, $t\geq 0$ and $|\phi^n_+(x,t)|\geq s$},\\
s^{-1}\cdot \phi^n_-(x,t) & \text{if $|(x,t)|\geq s$, $t\leq 0$ and $|\phi^n_-(x,t)|= s$},\\
f^{n-1}_s(\phi^n_-(x,t)) & \text{if $|(x,t)|\geq s$, $t\leq 0$ and $|\phi^n_-(x,t)|\geq s$}.
\end{cases}\]
The map $f^n_s$ is drawn in the following picture. The smaller ball $s\cdot D^n$ is mapped homeomorphically to $D^n$ by radial expansion from the origin. Next comes the region in the upper half of the larger ball that is inside the cone and outside the smaller ball. Each segment of a ray from the south pole $-e_n$ that lies in that region is mapped to a point which is determined by where we mapped the intersection of that ray-segment with the smaller ball (which was radially from the origin to the boundary of $D^n$). Third is the region in the upper half of the larger ball that is outside the cone. Each segment of a ray from the south pole $-e_n$ that lies in that region is first projected to the annulus in the equatorial plane of the two balls; we then apply the previously defined map $f^{n-1}_s$ to map the projected points to the equator of $D^n$. The lower half of the ball is mapped similarly.
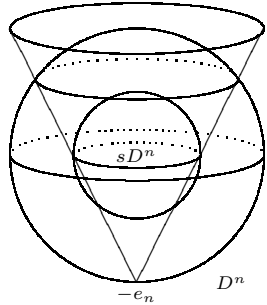
\[\begin{xy}
0;/r4pc/:(1,0),{\ellipse<>{-}},{\ellipse(,.2) d,^u{-}},{\ellipse(,.2) u,^d{.}}
,(.5,0);(1,0),{\ellipse<>{-}},{\ellipse(,.2) d,^u{-}},{\ellipse(,.2) u,^d{.}}
,(0,1);(1,-1) **@{-};(2,1) **@{-}
,(.2,.6);(1,.6),{\ellipse(,.2) d,^u{-}},{\ellipse(,.2) u,^d{.}}
,(0,1);(1,1),{\ellipse(,.2){-}}
,(1,0)*{\sb{sD^n}}
,(1,-1.1) *{\sb{-e_n}}
,(1.75,-1)*{\sb{D^n}}
\end{xy}\]
It is clear that $f^n_s$ gives a homotopy from $f^n_{1/2}$ to the identity and, given any disk $(D^n,q)$ in $\sK/B$, the map $f^n_s$ induces an $h$-equivalence from the $f$-disk $(D^n, q\circ f^n_{1/2})$ to the disk $(D^n,q)$. 

Finally we define the required cell structure on the domain of the desired map $\mu\colon (S^m,\bar{q})\rtarr (S^m,q)$. For each $n\leq m$, the boundary sphere $(S^n, q\circ f^{n+1}_{1/2}|S^n)$ is constructed from two copies of the $f$-disk $(D^n, q\circ f^n_{1/2})$ by gluing them along their boundary. The inclusions $(D^n, q\circ f^n_{1/2})\rtarr (S^n,q\circ f^{n+1}_{1/2}|S^n)$ of the two cells are given by projecting $D^n$ to the upper hemisphere from the south pole $-e_{n+1}$ and, similarly, by projecting $D^n$ to the lower hemisphere from the north pole $e_{n+1}$. The map 
\[\mu=f^{m+1}_{1/2}|S^m\colon  (S^m,q\circ f^{m+1}_{1/2}|S^m) \rtarr (S^m, q).\]
is then the required $f$-cell sphere approximation.
\end{proof}

\begin{proof}[Proof of \myref{thicken2}]
Define $\nu_s\colon D^n\rtarr D^n$ for $s\in[\tfrac12,1]$ by
\[\nu_s(x,t)=\begin{cases}
s^{-1}\cdot (x,t) & \text{if $|(x,t)|\leq s$},\\
|(x,t)|^{-1}\cdot (x,t) & \text{if $|(x,t)|\geq s$, $t\geq 0$ and $|x|\geq s$},\\
s^{-1}\cdot \phi^{n+1}_-(x,t) & \text{if $|(x,t)|\geq s$, $t\leq 0$ and $|\phi^{n+1}_-(x,t)|= s$},\\
|\phi^{n+1}_-(x,t)|^{-1}\cdot \phi^{n+1}_-(x,t) & \text{if $|(x,t)|\geq s$, $t\leq 0$ and $|\phi^{n+1}_-(x,t)|\geq s$},
\end{cases}\]
where $\phi^n_-$ is the projection as in the previous proof. Then $\nu_s$ maps $s\cdot D^n$ homeomorphically to $D^n$, it is radially constant on the region in the upper half space between the disks $D^n$ and $s\cdot D^n$ with respect to projection from the origin, and it is radially constant on the region in the lower half space between the two disks with respect to projection from the north pole.
\end{proof}


\section{The compatibility condition for the $qf$-model structure}\label{sec:qfproof}

This section is devoted to the proof that $\sK/B$ is a compactly generated model category. Since our generating sets $I^f_B$ and $J^f_B$ certainly satisfy conditions (i) and (iii) of \myref{Newcompgen}, it only remains to verify the compatibility condition (ii).  That is, we must show that a map has the RLP with respect to $I^f_B$ if and only if it is a $q$-equivalence and has the RLP with respect to $J^f_B$. Let $p\colon E\rtarr Y$ have the RLP with respect to $I^f_B$. Since all maps in $J^f_B$ are relative $I^f_B$-cell complexes, $p$ has the RLP with respect to $J^f_B$. To show that $\pi_n(p)$ is injective, let $\alpha\colon S^n\rtarr E$ represent an element in $\pi_n(E)$ such that $p\circ \alpha\colon S^n\rtarr Y$ is null-homotopic. Then there is a nullhomotopy $\beta\colon CS^n\rtarr Y$ that gives rise to a lifting problem
\[\xymatrix{S^n\ar[rr]^\alpha\ar[d]_i && E \ar[d]^{p}\\
D^{n+1}\ar[r]_-\nu & D^{n+1}\cong CS^n \ar[r]_-\beta & Y}\]
where $\nu\colon D^{n+1}\rtarr D^{n+1}$ is defined by
$$\nu(x) = \begin{cases}
2x & \text{if $|x|\leq \tfrac12$},\\
|x|^{-1}\cdot x & \text{if $|x|\geq \tfrac12$.}
\end{cases}$$
Then $i$ is an $f$-disk and we are entitled to a lift, which can be viewed as a nullhomotopy of $\alpha$ after we identify $D^{n+1}$ with $CS^n$.

To show that $\pi_n(p)$ is surjective, choose a representative $\alpha\colon S^n\rtarr Y$ of an element in $\pi_n(Y)$. The projection of $Y$ induces a projection $q\colon S^n\rtarr B$ and by \myref{thicken1} there is an $h$-equivalence $\mu\colon (S^n,\bar{q})\rtarr (S^n,q)$ such that $(S^n,\bar{q})$ is an $I^f_B$-complex with two cells in each dimension. We may therefore assume that the source of $\alpha$ is an $I^f_B$-cell complex. Inductively, we can then solve the lifting problems for the diagrams
\[\xymatrix{S^{k-1}\ar[rr]\ar[d]\ar[dr] && E\ar[d]^{p}\\
S^k_\pm \ar[r]_-{i_\pm} & S^k \ar[r]_-{\alpha|S^k} & Y,}\]
where $S^{k-1}\rtarr S^k$ is the inclusion of the equator and $i_\pm\colon S^k_\pm\rtarr S^k$ are the inclusions of the upper and lower hemispheres. We obtain a lift $S^n\rtarr E$. 

Conversely, assume that $p\colon E\rtarr Y$ is an acyclic $qf$-fibration. We must show that $p$ has the RLP with respect to any cell $i$ in $I^f_B$. We are therefore faced with a lifting problem
\[\xymatrix{S^n\ar[r]^\alpha\ar[d]_i & E \ar[d]^{p}\\
D^{n+1}\ar[r]_\beta & Y.}\]
Identifying $D^{n+1}$ with $CS^n$ we see that $\beta$ gives a nullhomotopy of $p\circ \alpha$. Since $\pi_n(p)$ is injective there is a nullhomotopy $\gamma\colon CS^n\rtarr E$ such that $\alpha=\gamma\circ i$, but it may not cover $\beta$. Gluing $\beta$ and $p\circ\gamma$ along $p\circ \alpha$ gives $\delta\colon S^{n+1}\rtarr Y$ such that $\delta|S^{n+1}_+ = \beta$ and $\delta|S^{n+1}_-=p\circ\gamma$. Surjectivity of $\pi_{n+1}(p)$ gives a map $\Delta\colon S^{n+1}\rtarr E$ and a homotopy $h\colon S^{n+1}\wedge I_+\rtarr Y$ from $p\circ \Delta$ to $\delta$. We now construct a diagram
\[\xymatrix{&S^{n+1}_+\ar[d]\ar[r]\ar[dl]_{j} & S^{n+1}_+\cup H \ar[r]^-{(-)/S^n}\ar[d]& S^{n+1}\times 0 \cup S^{n+1}_-\times 1\ar[r]^-{\Delta\cup \gamma}\ar[d] & E \ar[d]^{p}\\
D^{n+2}\ar[r]_\nu & D^{n+2} \ar[r]_\xi & D^{n+2}\ar[r]_-\phi & S^{n+1}\wedge I_+\ar[r]_-h & Y}\]
where the downward maps, except $p$, are inclusions. Part of the bottom row of the diagram is drawn schematically below. Let $H$ be the region on $S^{n+1}_-$ between the equator $S^n$ and the circle through $e_1$ and $-e_{n+2}$ with center on the line $\bR\cdot(e_1-e_{n+2})$. Let $\xi$ be a homeomorphism whose restriction to $S^{n+1}_+$ maps it homeomorphically to $S^{n+1}_+\cup H$. Define $\phi\colon D^{n+2}\rtarr D^{n+2}/S^n\cong S^{n+1}\wedge I_+$ as the composite of the quotient map that identifies the equator $S^n$ of $D^{n+2}$ to a point and a homeomorphism that maps the upper hemisphere $S^{n+1}_+$ to $S^{n+1}\times 0$, maps $H$ to $S^{n+1}_-\times 1$, and is such that 
$(h\circ \phi\circ \xi)|S^{n+1}_-=\beta$. The map $\nu$ is defined as above. 
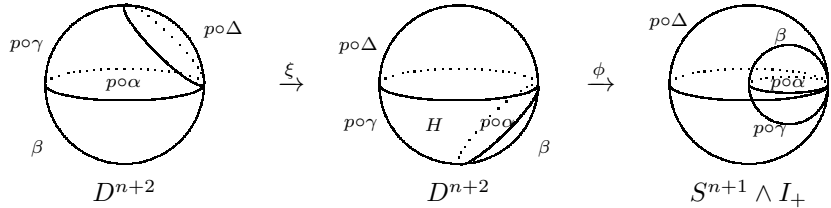
\[
\begin{xy}
0;/r5pc/:(.5,0),{\ellipse<>{-}},{\ellipse(,.2) d,^u{-}},{\ellipse(,.2) u,^d{.}}
,(.5,-.6)*!U{D^{n+2}}
,(.5,.5);(.75,.25),{\ellipse(,.2) d,_r{}},{\ellipse(,.2) u,^l{.}}
,(0,.25)*!R{\sb{p\circ\gamma}}
,(0,-.35)*!UR{\sb\beta}
,(1,.35)*!L{\sb{p\circ\Delta}}
,(.5,0)*{\sb{p\circ\alpha}}
\end{xy}
\quad\xrightarrow{\xi}\quad
\begin{xy}
0;/r5pc/:(.5,0),{\ellipse<>{-}},{\ellipse(,.2) d,^u{-}},{\ellipse(,.2) u,^d{.}}
,(.5,-.6)*!U{D^{n+2}}
,(.5,-.5);(.75,-.25),{\ellipse(,.2) l,_d{.}},{\ellipse(,.2) r,^u{-}}
,(0,-.25)*!R{\sb{p\circ\gamma}}
,(1,-.35)*!UL{\sb\beta}
,(0,.25)*!R{\sb{p\circ\Delta}}
,(.75,-.25)*{\sb{p\circ\alpha}}
,(.35,-.25)*{\sb{H}}
\end{xy}
\quad\xrightarrow{\phi}\quad
\begin{xy}
0;/r5pc/:(-.5,0),{\ellipse<>{-}},{\ellipse(,.2) d,^u{-}},{\ellipse(,.2) u,^d{.}}
;(-.25,0),{\ellipse<>{-}},{\ellipse(,.2) d,^u{-}},{\ellipse(,.2) u,^d{.}}
,(-.5,-.6)*!U{S^{n+1}\wedge I_+}
,(-.25,-.25)*!UR{\sb{p\circ\gamma}}
,(-.25,.25)*!DR{\sb\beta}
,(-1,.35)*!D{\sb{p\circ\Delta}}
,(-.25,0)*{\sb{p\circ\alpha}}
\end{xy}
\]
Since the restriction $S^n\rtarr S^{n+1}_-\cong D^{n+1}$ of $j$ agrees with the $f$-cofibration $i$ in our original lifting problem, we see that $j$ is a $J^f_B$-cell.  Since $p$ is a $qf$-fibration we get a lift in the outer trapezoid. Denote its restriction to $S^{n+1}_-\cong D^{n+1}$ by $k\colon D^{n+1}\rtarr E$. Then $k$ solves our original lifting problem.

\section{The quasifibration and right properness properties}\label{sec:qqfuasi}

We have now established the $qf$-model structures on both $\sK/B$ and $\sK_B$. 
We will derive the right properness of $\sK/B$, and therefore of $\sK_B$, from the fact that every $qf$-fibration is a quasifibration.

\begin{prop}\mylabel{qfles}
If $p\colon E\rtarr Y$ is a $qf$-fibration in $\sK/B$, then $p$ is a quasifibration. Therefore, for any choice of $e\in E$, there results a long exact sequence of homotopy groups
\[\cdots \rtarr \pi_{n+1}(Y,y)\rtarr \pi_n(E_y,e) \rtarr \pi_n(E,e)\rtarr \pi_n(Y,y)\rtarr\cdots\rtarr \pi_0(Y,y),\]
where $y=p(e)$ and $E_y=p^{-1}(y)$.
\end{prop}

\begin{proof}
We must prove that $p$ induces an isomorphism
\[\pi_n(p)\colon \pi_n(E,E_y,e)\rtarr \pi_n(Y,y)\]
for all $n\geq 1$ and verify exactness at $\pi_0(E,e)$. We begin with the latter.  Let $e'\in E$ and suppose that $p(e')$ is in the component of $y'$.
Let $\gamma\colon I\rtarr Y$ be a path in $Y$ from $p(e')$ to $y'$ such that $\gamma$ is constant at $p(e')$ for time $t\leq \tfrac12$. Let $q$ be the projection of $Y$. Then $(I,q\circ \gamma)$ is a relative $f$-disk, and we obtain a lift $\bar\gamma\colon I\rtarr E$ such that $\gamma=p\circ \bar\gamma$. But then $e'$ is in the same component as the endpoint of $\bar\gamma$, which lies in $E_y$.

Now assume that $n\geq 1$. Recall that an element of $\pi_n(X,A,*)$ can be represented by a map of triples $(D^n, S^{n-1}, S^{n-1}_+)\rtarr (X,A,*)$. We begin by showing surjectivity. Let $\alpha\colon (D^n,S^{n-1})\rtarr(Y,y)$ represent an element of $\pi_n(Y,y)$. We can view $D^n$ as a disk over $B$, and \myref{thicken2} gives an approximation $\nu\colon D^n\rtarr D^n$ by a relative $f$-disk. Then we can solve the lifting problem
\[\xymatrix{S^{n-1}_+ \ar[d]\ar[r]^-{c_e} & E\ar[d]^p\\
D^n\ar@{-->}[ur]_{\bar\alpha}\ar[r]_-{\alpha\circ \nu} & Y,}\]
where the top map is the constant map at $e\in E$. A lift is a map of triples $\bar\alpha \colon (D^n, S^{n-1},S^{n-1}_+)\rtarr(E,E_y,e)$ such that $p_*([\bar\alpha])=[\alpha]$.

For injectivity, let $\alpha\colon (D^n,S^{n-1},S^{n-1}_+)\rtarr (E,E_y,e)$ represent an element of $\pi_n(E,E_y,e)$ such that $p_*([\alpha])=0$. Then there is a homotopy $h\colon D^n\times I\rtarr Y$ rel $S^{n-1}$ such that $h|D^n\times 0=p\circ \alpha$ and $h$ maps the rest of the boundary of $D^n\times I$ to $y$. Let $A=D^n\times \{0,1\} \cup S^{n-1}_+\times I\subset \partial (D^n\times I)$ and define $\beta\colon A\rtarr E$ by setting $\beta(x)= \alpha(x)$ if $x\in D^n\times 0$ and $\beta(x)=e$ otherwise. We then have a homeomorphism of pairs $\phi\colon (D^n\times I, A)\rtarr(D^{n+1},S^n_+)$ and an approximation $\nu\colon D^{n+1}\rtarr D^{n+1}$ by an $f$-disk by \myref{thicken2}. We can now solve the lifting problem
\[\xymatrix{S^n_+\ar[d]\ar[r]^-{\beta\circ(\phi|_A)^{-1}} & E \ar[d]\\
D^{n+1}\ar@{-->}[ur]_{\bar\alpha}\ar[r]_-{h\circ\phi^{-1}\circ\nu} & Y,}\]
and this shows that $[\alpha]=0$ in $\pi_n(E,E_y,e)$.
\end{proof}

\begin{cor}\mylabel{fdrp}
The $qf$-model structure on $\sK/B$ is right proper.
\end{cor}

\begin{proof}
Since $qf$-fibrations are preserved under pullbacks, this is a
five lemma comparison of long exact sequences as in \myref{qfles}.
\end{proof}

\chapter{Equivariant $qf$-type model structures}

\section*{Introduction}

We return to the equivariant context in this chapter, letting $G$ be a Lie group throughout.  Actually, our definitions of the $q$ and $qf$-model structures work for arbitrary topological groups $G$, but we must restrict to Lie groups to obtain structures that are $G$-topological and behave well with respect to change of groups and smash products.  A discussion of details special to the
non-compact Lie case is given in \S7.1, but after that the generalization
from compact to non-compact Lie groups requires no extra work.  However,
we alert the reader that passage to {\em stable} equivariant homotopy theory raises new problems in the case of non-compact Lie groups that will not be 
dealt with in this book; see \S11.6.

The equivariant $q$-model structure on $G\sK_B$ is just the evident over and under $q$-model structure. However, the equivariant generalization of the 
$qf$-model structure is subtle. In fact, the subtlety is already relevant nonequivariantly when we study base change along the projection of a bundle.  The problem is that there are so few generating $qf$-cofibrations that many functors that take generating $q$-cofibrations to $q$-cofibrations do not take generating $qf$-cofibrations to $qf$-cofibrations.  We show how to get around this in \S7.2.  For each such functor that we encounter, we find an enlargement of the obvious sets of (acyclic) generating $qf$-cofibrations on the target of the functor so that it is still a model category, but now the functor does send generating (acyclic) $qf$-cofibrations to (acyclic) $gf$-cofibrations.  

The point is that there are many different useful choices of Quillen equivalent $qf$-type model structures, and they can be used in tandem.  For all of our choices, the weak equivalences are the $\sG$-equivalences and all cofibrations are both $q$-cofibrations and $f$-cofibrations.  Given a finite number of adjoint pairs with composable left adjoints such that each is a Quillen adjunction with its own choice of $qf$-type model structure, we can successively expand generating sets in target categories of the left 
adjoints to arrange that the composite be one of Quillen left adjoints
with respect to well chosen $qf$-type model structures. 

In \S7.2, we describe the $qf(\sC)$-model structure associated to a 
``generating set'' $\sC$ of $\sG$-complexes. Each such model structure
is $G$-topological.  In \S\ref{sec:qfadj}, we show that external smash products
are Quillen adjunctions when $\sC$ is a ``closed'' generating set, as can always
be arranged, and we show that all base change adjunctions $(f_!,f^*)$ are Quillen adjunctions.  We show further that there are generating sets for which $(f^*,f_*)$ is a Quillen adjunction when $f$ is a bundle with cellular fibers.  In \S\ref{sec:JHadj}, we show similarly that various change of group functors are given by Quillen adjunctions when the generating sets are well chosen.
In \S7.5, we show that $\text{Ho}G\sK_B$ has the properties required for application of the Brown representability theorem.  Those adjunctions between our basic functors that are not given by Quillen adjoint pairs in any choice of $qf$-model structure are studied in Chapter 9.

\section{Families and non-compact Lie groups}

Two sources of problems in the equivariant homotopy theory of general 
topology groups $G$ are that we only know that orbit types $G/K$ are 
$H$-CW complexes for $H\subset G$ when $G$ is a Lie group and $K$ is 
a compact subgroup and we only know that a product of orbits 
$G/H\times G/K$ is a $G$-CW complex when $G$ is a Lie group and $K$ 
(or $H$) is a compact subgroup.  This motivates us to restrict to Lie 
groups, for which these conclusions are ensured by \myref{Illman} and \myref{prodproper}. 

The compactness requirements force us to restrict orbit types when we prove
properties of our model structures, and the family $\sG$ of all compact subgroups of our Lie group $G$ plays an important role. We recall the relevant definitions, which apply to any topological group $G$ and are familiar and important in a variety of contexts. They provide a context that allows us to work with non-compact Lie groups with no more technical work than is required for compact Lie groups.  

A {\em family} $\sF$ in $G$ is a set of subgroups that is closed under passage to subgroups and conjugates.  An {\em $\sF$-space} is a $G$-space all of whose isotropy groups are in $\sF$. An {\em $\sF$-equivalence} is a $G$-map $f$ such that $f^H$ is a weak equivalence for all $H\in \sF$. If $X$ is an $\sF$-space,
then the only non-empty fixed point sets $X^H$ are those for groups $H\in\sF$. In particular, an $\sF$-equivalence between $\sF$-spaces is the same as a $q$-equivalence. For based $G$-spaces, the definition of an $\sF$-space must be altered to require that all isotropy groups except that of the $G$-fixed 
base point must be in $\sF$. The notion of an $\sF$-equivalence remains unchanged.

A map in $G\sK/B$ or $G\sK_B$ is an $\sF$-equivalence if its map of total $G$-spaces is an $\sF$-equivalence.  If $B$ is an $\sF$-space, then so is any $G$-space $X$ over $B$ and any fiber $X_b$. The only orbits that can 
then appear in our parametrized theory are of the form $G/H$ with $H\in \sF$ 
and the only non-empty fixed point sets $X^H$ are those for groups $H\in\sF$.
In particular, $H$ must be subconjugate to some $G_b$. An $\sF$-equivalence 
of $G$-spaces over an $\sF$-space $B$ is the same as a $q$-equivalence.

It is well-known that equivariant $q$-type model structures generalize naturally to families. One takes the weak equivalences to be the $\sF$-equivalences, and one restricts the orbits $G/H$ that appear as factors in the generating (acyclic) cofibrations to be those such that $H\in \sF$. The resulting cell complexes are called $\sF$-cell complexes.  Restricting tensors from $G$-spaces 
to $\sF$-spaces, we obtain a restriction of the notion of a $G$-topological model category to an $\sF$-topological model category that applies here; see \myref{Gtopfamily}.

Proper $G$-spaces are particularly well-behaved $\sG$-spaces, where $\sG$ is the family of compact subgroups of our Lie group $G$, and $\sG$-cell complexes are proper $G$-spaces. Restricting base $G$-spaces to be proper, or more generally to be $\sG$-spaces, has the effect of restricting all relevant orbit types 
$G/H$ to ones where $H$ is compact. However, this is too restrictive for some purposes. For example, we are interested in developing nonparametrized equivariant homotopy theory for non--compact Lie groups $G$. Here $B=*$ is a 
$G$-space which, in the unbased sense, is not a $\sG$-space.  

We therefore do not make the blanket assumption that $B$ is a $\sG$-space.  
We give the $q$-model structure in complete generality, in \myref{qoverB}, but after that we restrict to $\sG$-model structures throughout. That is, our weak
equivalences will be the $\sG$-equivalences. This ensures that, after cofibrant approximation, our total $G$-spaces are $\sG$-spaces. This convention 
enables us to arrange that all of our model categories are $G$-topological.  Everything in this chapter applies more generally to the study of parametrized $\sF$-homotopy theory for any family $\sF$; see \myref{Fremark}.

The reader may prefer to think in terms of either the case when $B=*$
or the case when $B$ is proper.  Indeed, in order to resolve the problems intrinsic to the parametrized context that are described in the Prologue, 
which we do in Chapter 9, it seems essential that we restrict to proper 
actions on base spaces. The reason is that Stasheff's \myref{ss} relating 
the equivariant homotopy types of fibers and total spaces plays a 
fundamental role in the solution.  Alternatively, the reader may prefer
to focus just on compact Lie groups, reading $q$-equivalence 
instead of $\sG$-equivalence and $G$-space instead of $\sG$-space. 

\section{The equivariant $q$ and $qf$-model structures}\label{sec:qfeq}

Recall from \myref{UrIJ} that the sets $I$ and $J$ of generating cofibrations and generating acyclic cofibrations of $G$-spaces are defined as the sets of all maps of the form $G/H\times i$, where $i$ is in the corresponding set $I$ or $J$ of maps of spaces. 

\begin{defn}\mylabel{IJBG}
Starting from the sets $I$ and $J$ of maps of $G$-spaces, define sets $I_B$ and $J_B$ of maps of ex-$G$-spaces over $B$ in exactly the same way that their nonequivariant counterparts were defined in terms of the sets $I$ and $J$ of maps of spaces in \myref{IJB}. Note that if $B$ is a $\sG$-space, then only orbits $G/H$ with $H$ compact appear in the sets $I_B$ and $J_B$.
\end{defn}

Taking $Y = B$ in the usual composite adjunction
\begin{equation}\label{adadad}
G\sK(G/H\times T , Y) \iso H\sK(T,Y) \iso \sK(T,Y^H)
\end{equation}
for non-equivariant spaces $T$ and $G$-spaces $Y$, we can translate back and forth between equivariant homotopy groups and cells for $G$-spaces over $B$ on the one hand and nonequivariant homotopy groups and cells for spaces over $B^H$ on the other. Maps in each of the equivariant sets specified in \myref{IJBG} correspond by adjunction to maps in the nonequivariant set with the same name. Systematically using this translation, it is easy to use \myref{compgen} to generalize the $q$-model structures on $\sK/B$ and $\sK_B$ to corresponding model structures on $G\sK/B$ and $G\sK_B$.  We obtain the following theorem.

\begin{thm}[$q$-model structure]\mylabel{qoverB}
The categories $G\sK/B$ and $G\sK_B$ are compactly generated proper $\sG$-top\-o\-lo\-gi\-cal model categories whose $q$-equivalences, $q$-fibrations, and $q$-co\-fibrations are the maps whose underlying maps of total $G$-spaces are $q$-equivalences, $q$-fibrations, and $q$-cofibrations. The sets $I_B$ and $J_B$ are the generating $q$-cofibrations and generating acyclic $q$-cofibrations, and all $q$-cofibrations are $\bar{h}$-cofibrations. If $B$ is a $\sG$-space, then the model structure is $G$-topological.
\end{thm}

To show that the $q$-model structures are $\sG$-topological, and $G$-topological if $B$ is a $\sG$-space, we must inspect the maps $i\Box j$ in $G\sK/B$, where $i$ is a generating $q$-cofibration in $G\sK/B$ and $j$ is a generating cofibration in $G\sK$. They have the form
\[i\Box j\colon G/H\times G/K \times \partial(D^m\times D^n) 
\rtarr G/H\times G/K \times D^m\times D^n\]
given by the product of $G/H\times G/K$ with the inclusion of the boundary of $D^m\times D^n$. By \myref{prodproper}, $G/H\times G/K$ is a proper $G$-space
if $H$ or $K$ is compact.  Since we are assuming that $G$ is a Lie group, we can then triangulate $G/H\times G/K$ as a $\sG$-CW complex and use the triangulation to write $i\Box j$ as a relative $I_B$-cell complex. The case when either $i$ or $j$ is acyclic works in the same way. As explained in \myref{ouchtoo}, there is no problem with projection maps in this external context. Moreover, if $i$ is an $f$-cofibration, then so is $i\Box j$, as we see from the fiberwise NDR characterization. 

One might be tempted to generalize the $qf$-model structure to the equivariant context in exactly the same way as we just did for the $q$-model structure. This certainly works to give a model structure.  However, there is no reason to think that it is either $G$ or $\sG$-topological.  The problem is that we need 
$i\Box j$ above to be a $qf$-cofibration when $i$ is a generating $qf$-cofibration, and triangulations into $f$-cells are hard to come by. Therefore the $G$-CW structure on $G/H\times G/K$ will rarely produce a relative $I^f_B$-cell complex. This means that we must be careful when selecting the generating (acyclic) $qf$-cofibrations if we want the resulting model structure to be $G$-topological. We will build the solution into our definition of $qf$-type model structures, but we need a few preliminaries.

We shall make repeated use of the adjunction 
\begin{equation}\label{keyadj}
G\sK(C \times T , Y)\cong \sK(T,\text{Map}_G(C, Y))
\end{equation}
for non-equivariant spaces $T$ and $G$-spaces $C$ and $Y$. This is a generalization of (\ref{adadad}).
Taking $Y=B$, we note in particular that it gives a correspondence 
between maps $f\colon T\rtarr T'$ over $\text{Map}_G(C,B)$ and 
$G$-maps $\text{id}\times f\colon C\times T\rtarr C\times T'$ over $B$.

\begin{lem}\mylabel{mapgood}
If $C$ is a $\sG$-cell complex, then the functor $\text{Map}_G(C,-)\colon G\sK \rtarr \sK$ preserves all $q$-equivalences.
\end{lem}

\begin{proof} The functor $\text{Map}(C,-)$ is a Quillen right adjoint since the $q$-model structure on $G\sK$ is $\sG$-topological. The $G$-fixed point functor is also a Quillen right adjoint, for example by \myref{fixedptrQa0} below. The composite $\text{Map}_G(C,-)$ therefore preserves $q$-equivalences between $q$-fibrant $G$-spaces.  However, every $G$-space is $q$-fibrant.
\end{proof}

Observe that \myref{prodproper} gives that the collection of $\sG$-cell complexes is closed under products with arbitrary orbits $G/H$ of $G$.

\begin{defn}\mylabel{IJBG2}
Let $\sO_G$ denote the set of all orbits $G/H$ of $G$. Any set $\sC$ of $\sG$-cell complexes in $G\sK$ that contains all orbits $G/K$ with $K\in \sG$ and is closed under products with arbitrary orbits in $\sO_G$ is called a \emph{generating set}. It is a {\em closed} generating set if it is 
closed under finite products. The {\em closure} of a generating set
$\sC$ is the generating set consisting of the finite products of the 
$\sG$-cell complexes in $\sC$.
We define sets of generating $qf(\sC)$-cofibrations and acyclic $qf(\sC)$-cofibrations in $G\sK/B$ associated to any generating set $\sC$ as follows.
\begin{enumerate}[(i)]
\item Let $I^f_B(\sC)$ consist of the maps
\[(\text{id}\times i)(d)\colon C\times S^{n-1} \rtarr C\times D^n\]
such that $C\in \sC$, $d\colon C\times D^n\rtarr B$ is a $G$-map, $i$ is the boundary inclusion, and the corresponding map $i$ over $\text{Map}_G(C,B)$ is a generating $qf$-cofibration in $\sK/\text{Map}_G(C,B)$; that is, $i$ must be an $f$-cofibration. 
\item Similarly let $J^f_B(\sC)$ consist of the maps 
$$(\text{id}\times i)(d)\colon C\times S^n_+\rtarr C\times D^{n+1}$$
such that $C\in\sC$, $d\colon C\times D^{n+1}\rtarr B$ is a $G$-map, $i$ is the inclusion, and the corresponding map $i$ over $\text{Map}_G(C,B)$ is a generating acyclic $qf$-cofibration in $\sK/\text{Map}_G(C,B)$.
\end{enumerate}
Adjoining disjoint sections, we obtain the corresponding sets $I^f_B(\sC)$ and $J^f_B(\sC)$ in $G\sK_B$.
\end{defn}

Fix a generating set $\sC$.  We define a $qf$-type model structure based on $\sC$, called the $qf(\sC)$-model structure. Its weak equivalences are the $\sG$-equivalences, which are the same as the $q$-equivalences when $B$ is a $\sG$-space. We define the $qf(\sC)$-fibrations.

\begin{defn}\mylabel{qffibdef}
A map $f$ in $G\sK/B$ is a \emph{$qf(\sC)$-fibration} if $\text{Map}_G(C,f)$ is a $qf$-fibration in $\sK/\text{Map}_G(C,B)$ for all $C\in\sC$. A map in $G\sK_B$ is a \emph{$qf(\sC)$-fibration} if the underlying map in $G\sK/B$ is one. In
either category, a map $f$ is a {\em $\sG$-quasifibration} if $f^H$ is a quasifibration for $H\in \sG$.
\end{defn}

\begin{thm}[$qf$-model structure]\mylabel{Gqfstr}
For any generating set $\sC$, the categories $G\sK/B$ and $G\sK_B$ are well-grounded (hence $G$-topological) model categories. The weak equivalences and fibrations are the $\sG$-equi\-valences and the $qf(\sC)$-fibrations. The sets $I^f_B(\sC)$ and $J^f_B(\sC)$ are the generating $qf(\sC)$-cofib\-rations and the generating acyclic $qf(\sC)$-cofibrations. All $qf(\sC)$-cofibrations are both $q$-cofibrations and $\bar{f}$-cofib\-rations, and all $qf(\sC)$-fibrations are $\sG$-quasifibrations.
\end{thm}

\begin{proof} Recall from \myref{exwellgr} that the $q$-equivalences in $G\sK/B$ and $G\sK_B$ are well-grounded with respect to the ground structure given in \myref{exbackdef} and \myref{exback}. It follows that the $\sG$-equivalences
are also well-grounded. It suffices to verify conditions (i)--(iv) of \myref{Newcompgen}.  The acyclicity condition (i) is obvious.

Consider the compatibility condition (ii). By the adjunction (\ref{keyadj}), a map $f$ has the RLP with respect to $I^f_B(\sC)$ if and only if $\text{Map}_G(C,f)$ has the RLP with respect to $I^f_{\text{Map}_G(C,B)}$ for all $C\in\sC$. By the compatibility condition for the nonequivariant $qf$-model structure, that holds if and only if $\text{Map}_G(C,f)$ is a $q$-equivalence and has the LLP with respect to $J^f_{\text{Map}_G(C,B)}$ for all $C\in\sC$. By \myref{mapgood}, $\text{Map}_G(C,f)$ is a $q$-equivalence if $f$ is one. Conversely, if $\text{Map}_G(C,f)$ is a $q$-equivalence for all $C\in\sC$, then the case $C=G/K$ shows that $f^K$ is a $q$-equivalence for every compact $K$ and thus $f$ is a $\sG$-equivalence. By the adjunction again, we see that $f$ has the RLP with respect to $I^f_B(\sC)$ if and only if $f$ is a $\sG$-equivalence which has the RLP with respect to $J^f_B(\sC)$.

The fiberwise NDR characterization of $\bar{f}$-cofibrations given in \myref{fNDR} shows that $I^f_B(\sC)$ and $J^f_B(\sC)$ consist of $\bar{f}$-cofibrations, as stipulated in (iii).  More precisely, if $(u,h)$, $u\colon D^n\rtarr I$ and $h\colon D^n\times I \rtarr D^n$, represents $(D^n,S^{n-1})$ as a fiberwise NDR-pair over $\text{Map}_G(C,B)$, then the map $v=u\circ\pi\colon C\times D^n\rtarr D^n\rtarr I$ and the homotopy given by the maps $\text{id}\times h_t$ over $B$ corresponding to the $h_t$ represent $(C\times D^n,C\times S^{n-1})$ as a fiberwise NDR pair over $B$.

Since $\text{Map}_G(G/K,f)\cong f^K$ is a nonequivariant $qf$-fibration for any $qf(\sC)$-fibration $f$, every $qf(\sC)$-fibration is a $\sG$-quasifibration by \myref{qfles}. That the model structure is right proper follows as in \myref{fdrp}.

Finally, we must verify the $\Box$-product condition (iv). The relevant maps $i\Box j$,
\[i\colon C\times S^{m-1}\rtarr C\times D^m\quad \text{and}\quad
j\colon G/H\times S^{n-1}\rtarr G/H\times D^n,\]
are of the form
\[C\times G/H \times k\colon C\times G/H \times \partial(D^m\times D^n)\rtarr 
C\times G/H \times D^m\times D^n,\]
where $k$ is the boundary inclusion. Now $C\times G/H\in \sC$ by the closure property of the generating set, so we don't need to triangulate. The projection of the target factors through the projection of the target $C\times D^m$ of $i$.  To see that the corresponding map $k$ over $\text{Map}_G(C\times G/H,B)$ is an $\bar{f}$-cofibration, let $(u,h)$ represent $(D^m,S^{m-1})$ as a fiberwise NDR-pair over $\text{Map}_G(C,B)$ and let $(v,j)$ represent $(D^n,S^{n-1})$ as an NDR-pair; we can think of the latter as a fiberwise NDR-pair over $* = \text{Map}_G(G/H,*)$. Then the usual product pair representation (for example, \cite[p. 43]{Concise}) exhibits $k$ as a fiberwise NDR over $\text{Map}_G(C,B)\times \text{Map}_G(G/H,*)$ and thus, by the factorization of the projection of $i\Box j$, also over $\text{Map}_G(C\times G/H,B\times *)$.
\end{proof}

\begin{thm} If $\sC\subset \sC'$ is an inclusion of generating sets, then the identity functor is a left Quillen equivalence from $G\sK/B$ with the $qf(\sC)$-model structure to $G\sK/B$ with the $qf(\sC')$-model structure. The identity functor is also a left Quillen equivalence from $G\sK/B$ with the $qf(\sC)$-model structure to $G\sK/B$ with the $q$-model structure.  Both statements also hold for the identity functor on $G\sK_B$.
\end{thm}

\begin{proof}
The first statement is obvious. For the second, if $\text{id}_C\times i$ is a generating $qf(\sC)$-cofibration, then $C$ is a $\sG$-cell complex and we can use the triangulation to write $\text{id}_C\times i$ as a relative $I_B$-cell complex.
\end{proof}

\begin{thm}
For any $\sC$, the identity functor is a left Quillen adjoint from $G\sK/B$ with the $qf(\sC)$-model structure to $G\sK/B$ with the $f$-model structure. Similarly, the identity functor is a left Quillen adjoint from $G\sK_B$ with the $qf(\sC)$-model structure to $G\sK_B$ with the $fp$-model structure.
\end{thm}

The last result, which implements the philosophy of \S4.1, is false for the 
$q$-model structures.

\begin{rem}\mylabel{Theqf}
The smallest generating set $\sC$ is the set of all (non-empty) finite products of orbits $G/H$ of $G$ such that at least one of the factors has $H$ compact. 
Clearly it is a closed generating set. 
Henceforward, by {\em the} $qf$-model structure, we mean the $qf(\sC)$-model structure associated to this choice of $\sC$. In the nonequivariant 
case, this is the $qf$-model structure of the previous chapter.
\end{rem}

\begin{rem}
In the nonparametrized setting, the $\sG$-model structure associated to the
$q$-model structure and the $qf(\sC)$-model structures on $G\sK = G\sK/*$ 
coincide, and similarly for $G\sK_*$. This holds since the $f$-cofibrations and $h$-cofibrations over a point coincide and since the $C\in\sC$ for any choice
of $\sC$ are $\sG$-cell complexes. Of course, the $qf(\sC)$-model structures have more generating (acyclic) cofibrations.
\end{rem}

\begin{rem}
It might be useful to combine the various $qf(\sC)$-model structures by taking the union of the $qf(\sC)$-cofibrations over some suitable collection of generating sets $\sC$ and so obtain a ``closure'' of the $qf$-model structure whose cofibrations are as close as possible to being the intersection of the $q$-cofibrations with the $\bar{f}$-cofibrations. We do not know whether or not 
that can be done.
\end{rem}

\begin{rem}\mylabel{Fremark}
As noted in the introduction, we can generalize the $q$ and $qf(\sC)$-model structures to the context of families $\sF$. We generalize the $q$-model structure to the $\sF$-model structure by taking the $\sF$-equivalences and $\sF$-fibrations and by restricting the sets $I_B$ and $J_B$ to be constructed from orbits $G/H$ with $H\in \sF$. The resulting model structure will then be $(\sF\cap \sG)$-topological and $\sF$-topological if the base space $B$ is a $\sG$-space.

To generalize the $qf(\sC)$ model structure, we take the weak equivalences to be the $\sF\cap\sG$-equivalences and we require the generating set $\sC$ to consist of $\sF\cap\sG$-cell complexes, to contain the orbits $G/K$ for $K\in \sF\cap \sG$, and to be closed under products with orbits $G/K$ where $K\in \sF$.  With that modification, everything else above goes through unchanged.
\end{rem}

\section{External smash product and base change adjunctions}\label{sec:qfadj}

The following results relate the $q$ and $qf(\sC)$-model structures to smash products and base change functors and show that various of our adjunctions are given by Quillen adjoint pairs and therefore induce adjunctions on passage to homotopy categories.  For uniformity, we must understand the $q$-model structure to mean the associated $\sG$-model structure, although many of the results do apply to the full $q$-model structure. Those results that refer to $q$-equivalences by name work equally well for $\sG$-equivalences. Most of the results in this section and the next apply both to the $\sG$-model structure 
and to the $qf(\sC)$-model structure for any generating set $\sC$. We agree to omit the $q$ or $qf(\sC)$ from the notations in those cases. In other cases, we will have to restrict to well chosen generating sets $\sC$.

With these conventions, our first result is clear from the fact that our 
model structures are $G$-topological.

\begin{prop}\mylabel{smaB} 
For a based $G$-CW complex $K$, the functor $(-)\sma_B K$ preserves cofibrations and acyclic cofibrations, hence the functor $F_B(K,-)$ preserves fibrations and acyclic fibrations. Thus $((-)\sma_B K, F_B(K,-))$ is a Quillen adjoint pair of endofunctors of $G\sK_B$. 
\end{prop}

For the rest of our results, recall from \myref{reducts} that a left adjoint 
that takes generating acyclic cofibrations to acyclic cofibrations preserves acyclic cofibrations.  The following two results apply to the $qf(\sC)$-model structure for any closed generating set $\sC$.

\begin{prop}\mylabel{Boxcof20}
If $i\colon X\rtarr Y$ and $j\colon W\rtarr Z$ are cofibrations over base $G$-spaces $A$ and $B$, then 
$$i\Box j\colon (Y\barwedge  W)\cup_{X\barwedge  W}(X\barwedge Z)\rtarr Y\barwedge Z$$
is a cofibration over $A\times B$ which is acyclic if either $i$ or $j$ is acyclic.
\end{prop}

\begin{proof} 
It suffices to inspect $i\Box j$ for generating (acyclic) cofibrations as was done for the case $A=*$ in the proof of \myref{Gqfstr}. For generating 
cofibrations, the argument there generalizes without change to this setting.
The assumption that $\sC$ is closed avoids the need for triangulations here. 
For the acyclicity, it suffices to work in the $q$-model structure, for which
the conclusion is both more general and easier to prove. There it is easily
checked using triangulations of products of $\sG$-cell complexes that if $i$ is 
a generating cofibration and $j$ is a generating acyclic cofibration, then 
$i\Box j$ is an acyclic cofibration.
\end{proof}

Of course, by \myref{ouchtoo}, the analogue for internal smash products fails. Taking $W =*_B$ and changing notations, we obtain the following special case.

\begin{cor}\mylabel{smaAB}
Let $Y$ be a cofibrant ex-space over $B$. Then the functor 
$(-)\barwedge Y$
from ex-spaces over $A$ to ex-spaces over $A\times B$ preserves cofibrations 
and acyclic cofibrations, hence the functor $\bar{F}(Y,-)$ 
from ex-spaces over $A\times B$ to ex-spaces over $A$ preserves 
fibrations and acyclic fibrations. Thus
$((-)\barwedge Y,\, \bar{F}(Y,-))$ is a 
Quillen adjoint pair of functors between $G\sK_A$ and $G\sK_{A\times B}$.
\end{cor}

The next two results apply to the $qf(\sC)$-model structures for any $\sC$, provided that we use the same generating set $\sC$ for both $G\sK_A$ and $G\sK_B$. 

\begin{prop}\mylabel{Qad10}
Let $f\colon A\rtarr B$ be a $G$-map. Then the functor $f_!$ preserves cofibrations and acyclic cofibrations, hence $(f_{!},f^*)$ is a Quillen adjoint pair. The functor $f_!$ also preserves $q$-equivalences between well-sectioned ex-spaces. If $f$ is a $q$-fibration, then the functor $f^*$ preserves all 
$q$-equivalences.
\end{prop}

\begin{proof}
If $(D,p)$ is a space over $A$, then $f_!((D,p)\amalg A)=(D,f\circ p)\amalg B$. Therefore $f_!$ takes generating (acyclic) $q$-cofibrations over $A$ to such maps over $B$. If $(u, h)$ represents $(D^n, S^{n-1})$ as a fiberwise NDR-pair over $\text{Map}_G(C,A)$, then, after composing the projection maps with
$\text{Map}_G(C,A)\rtarr \text{Map}_G(C,B)$,
it also represents $(D^n, S^{n-1})$ as a fiberwise NDR-pair over $\text{Map}_G(C,B)$. It follows that $f_!$ also preserves the generating (acyclic) $qf$-cofibrations. Recall that the well-sectioned ex-spaces are those that are $\bar{f}$-cofibrant and that $f$-cofibrations are $h$-cofibrations. Since $f_!X$ is defined by a pushout in $G\sK$, the gluing lemma in $G\sK$ implies that $f_!$ preserves $q$-equivalences between well-sectioned ex-spaces.

If $f$ is a $q$-fibration and $k\colon  Y\rtarr Z$ is a $q$-equivalence of ex-spaces over $B$, consider the diagram
\[\xymatrix@=.4cm{
&& f^*Z \ar[rrr]\ar[dddl]|(.4)\hole &&& Z\ar[dddl]\\
f^*Y \ar[urr]^-{f^*k}\ar[ddr]\ar[rrr] &&& Y \ar[urr]^-k\ar[ddr]\\
\\
& A \ar[rrr]_f &&& B.}\]
The relation $(A\times_B Z)\times _Z Y\iso A\times_BY$ shows that the top square is a pullback, and the pullback $f^*Z\rtarr Z$ of $f$ is a $q$-fibration. Since the $q$-model structure on the category of $G$-spaces is right proper, it follows that $f^*k$ is a $q$-equivalence.
\end{proof}

\begin{prop}\mylabel{ffequiv0}
If $f\colon A\rtarr B$ is a $q$-equivalence, then $(f_{!},f^*)$ is a Quillen equivalence.
\end{prop}

\begin{proof} The conclusion holds if and only if the induced adjunction
on homotopy categories is an adjoint equivalence \cite[1.3.3]{Hovey}, so
it suffices to verify the usual defining condition for a Quillen
adjunction in either model structure. The condition for the other
model structure follows formally.  We choose the $q$-model structure. 
Let $X$ be a $q$-cofibrant ex-space over $A$ and $Y$ be a $q$-fibrant 
ex-space over $B$, so that $A\rtarr X$ is a $q$-cofibration and $Y\rtarr B$ 
is a $q$-fibration of $G$-spaces. Since the model structure on the category 
of $G$-spaces is left and right proper, inspection of the defining diagrams 
in \myref{retract1} shows that the canonical maps  $X\rtarr f_{!}X$ and $f^*Y\rtarr Y$ of total spaces are $q$-equivalences.  For an ex-map
$k\colon f_!X\rtarr Y$ with adjoint $\tilde k\colon X\rtarr f^*Y$, the commutative diagram
$$\xymatrix{
X\ar[r] \ar[d]_{\tilde{k}} & f_!X \ar[d]^{k} \\
f^*Y \ar[r] & Y}$$
of total spaces then implies that $k$ is a $q$-equivalence if and only if $\tilde{k}$ is a $q$-equivalence.
\end{proof}

In view of \myref{noway}, we can at best expect only a partial and restricted analogue of \myref{Qad10} for $(f^*,f_*)$. We first give a result for the
$q$-model structure and then show how to obtain the analogue for the
$qf(\sC)$-model structures using well chosen generating sets $\sC$.

\begin{prop}\mylabel{Qad20}
Let $f\colon A\rtarr B$ be a $G$-bundle such that $B$ is a $\sG$-space and
each fiber $A_b$ is a $G_b$-cell complex. Then $(f^*,f_*)$ is a Quillen adjoint pair with respect to the $q$-model structures. Moreover, if the total space of an ex-$G$-space $Y$ over $B$ is a $\sG$-cell complex, then so is the total space of $f^*Y$.
\end{prop}

\begin{proof}
Since $f$ is a $q$-fibration, $f^*$ preserves $q$-equivalences. It therefore suffices to show that $f^*$ takes generating cofibrations in $I_B$ to relative $I_A$-cell complexes. Observe first that if $\ph\colon G/H \rtarr B$ is a $G$-map with $\ph(eH)= b$, then $H\subset G_b$ and the pullback $G$-bundle $\phi^*f\colon f^*(G/H,\ph)\rtarr G/H$ of $f$ along $\ph$ is $G$-homeomorphic to $G\times_H A_b\rtarr G/H$. We can triangulate orbits in a $G_b$-cell decomposition of $A_b$ as $H$-CW complexes, by \myref{Illman}, and so give $A_b$ the structure of an $H$-cell complex.  Then $G\times_H A_b$ has an induced structure of a $\sG$-cell complex and thus so does $f^*(G/H,\ph)$.

For a space $d\colon E\rtarr B$ over $B$ with associated ex-space $E\amalg B$ over $B$, we have $f^*(E\amalg B) = f^*E\amalg A$. Let $E = G/H\times D^n$ and let $i\colon G/H\rtarr G/H\times D^n$ be the inclusion $i(gH)=(gH,0)$. The composite $d\circ i$ is a map $\phi$ as above. Since the identity map on 
$G/H\times D^n$ is homotopic to the composite $i\circ \pi\colon G/H\times D^n\rtarr G/H\times D^n$, where $\pi$ is the projection, the pullback $G$-bundle $d^*f\colon f^*(E,d)\rtarr E$ is equivalent to the pullback bundle $(\phi\circ\pi)^*f\colon f^*(E,\phi\circ\pi)\rtarr E$. But the latter is the product of  $\phi^*f\colon f^*(G/H,\ph)\rtarr G/H$ and the identity map of $D^n$ as we see from the following composite of pullbacks
\[\xymatrix{f^*(G/H\times D^n,\ph\com\pi) \ar[r]\ar[d]_{(\phi\circ\pi)^*f} 
& f^*(G/H,\ph) \ar[r]\ar[d]_{\phi^*f} & f^*(G/H\times D^n,d) \ar[r]\ar[d]_{d^*f} & A\ar[d]^f\\
G/H\times D^n\ar[r]_-{\pi} & G/H \ar[r]_-i & G/H\times D^n\ar[r]_-d & B.}\]
The $\sG$-cell structure on $f^*(G/H,\ph)$ gives a canonical decomposition of the inclusion $f^*(G/H,\ph) \times S^{n-1}\rtarr f^*(G/H,\ph) \times D^n$ as a relative $\sG$-cell complex. The last statement follows by applying this analysis inductively to the cells of $Y$. 
\end{proof}

The previous result fails for the $qf$-model structure. In fact, it already fails nonequivariantly for the unique map $f\colon A\rtarr *$, where $A$ is a 
CW-complex.  The proof breaks down when we try to use a 
cell decomposition of $A$ (the fiber over $*$) to decompose cells 
$A\times S^{n-1}\rtarr A\times D^n$ over $A$ as relative $I^f_A$-cell 
complexes. Similarly, the equivariant proof above breaks down when we try to use the $G$-cell structure of $f^*(G/H,\ph)$ to obtain a relative $I^f_A$-cell complex. Note, however, that there is no problem when the fibers are 
homogeneous spaces $G/H$; the nonequivariant analogue is just the trivial 
case when $f$ is a homeomorphism, but principal bundles and projections
$G/H\times B\rtarr B$ give interesting equivariant examples.
For the general equivariant case, we choose a closed generating set $\sC(f)$ that depends on the $G$-bundle $f$ and a given closed generating set $\sC$. Using the $qf(\sC)$-model structures on $G\sK_A$ and $G\sK_B$, we then 
recover the Quillen adjunction.
 
\begin{con}
Let $f\colon A\rtarr B$ be a $G$-bundle such that $B$ is a $\sG$-space and
each fiber $A_b$ is a $G_b$-cell complex and let $\sC$ be a closed generating set.  We construct the set $\sC(f)$ inductively.  Let $\sC(f)_0=\sC$ and suppose that we have constructed a set $\sC(f)_n$ of $\sG$-cell complexes in $G\sK$ that is closed under both finite products and products with arbitrary orbits $G/H$ of $G$. Let
\[\sA_n=\{f^*(C,\ph)\mid \text{$C\in \sC(f)_n$ and $\phi\in G\sK(C,B)$}\}.\]
Then let $\sC(f)_{n+1}$ consist of all finite products of spaces in $\sC(f)_n\cup \sA_n$. Note that $\sC(f)_{n+1}$ contains $\sC(f)_n$ and that the $f^*(C,\phi)$ are $\sG$-cell complexes by the last statement of \myref{Qad20}. Finally, let $\sC(f)=\bigcup \sC(f)_n$.  Clearly $\sC(f)\supset \sC$ is a closed 
generating set that contains $f^*(C,\ph)$ for all $C\in \sC(f)$ and all 
$G$-maps $\phi\colon C\rtarr B$.
\end{con}

\begin{prop}\mylabel{Qad202}
Let $f\colon A\rtarr B$ be a $G$-bundle such that $B$ is a $\sG$-space
and all fibers $A_b$ are $G_b$-cell complexes. Then $(f^*,f_*)$ is a 
Quillen adjoint pair with
respect to the $qf(\sC(f))$-model structures on $G\sK_A$ and $G\sK_B$.
\end{prop}

\begin{proof}
Reexamining the proof of \myref{Qad20}, but starting with a map $d\colon E=C\times D^n\rtarr B$ where $C\in \sC(f)$, we see that 
$$f^*E\iso f^*(C,\ph)\times D^n$$
where $\phi=d\circ i$. Since $f^*(C,\ph)$ is a $\sG$-cell complex in $\sC(f)$, it remains only to show that $f^*(C,\ph)\times S^{n-1}\rtarr f^*(C,\ph)\times D^n$ is an $f$-cofibration. Let $(u,h)$ represent $(D^n,S^{n-1})$ as a fiberwise NDR-pair over $\text{Map}_G(C,B)$. Applying $f^*$ to the corresponding maps $h_t\colon C\times D^n\rtarr C\times D^n$ over $B$, we obtain maps $f^*h_t\colon f^*E\rtarr f^*E$ over $A$. Under the displayed isomorphism, these maps give a homotopy $f^*h\colon D^n\times I\rtarr D^n$ that, together with $u$, represents $(D^n,S^{n-1})$ as a fiberwise NDR-pair over $\text{Map}_G(f^*(C,\phi),A)$.
\end{proof}

\section{Change of group adjunctions}\mylabel{sec:JHadj}

We consider change of groups in the $q$ and the $qf$-model structures, starting with the former. The context of the following results is given in \S\S2.3 and 2.4. 

\begin{prop}\mylabel{grpres}
Let $\tha\colon G\rtarr G'$ be a homomorphism of Lie groups. The restriction of action functor
\[\tha^*\colon G'\sK_B \rtarr G\sK_{\tha^*B}\]
preserves $q$-equivalences and $q$-fibrations. If $B$ is a $\sG'$-space, then it also preserves $q$-cofibrations.
\end{prop}

\begin{proof}
Since $(\tha^*A)^H=A^{\tha(H)}$ for any subgroup $H$ of $G$ and a map 
$f:X\rtarr Y$ of $G$-spaces is a $q$-equivalence or $q$-fibration if and only if each $f^H$ is a $q$-equivalence or $q$-fibration, it is clear that $\tha^*$ preserves $q$-equivalences and $q$-fibrations. To study $q$-cofibrations, recall that $\tha$ factors as the composite of a quotient homomorphism, an isomorphism, and an inclusion.  If $\tha$ is an inclusion and $H'$ is a compact subgroup of $G'$, then we can triangulate $G'/H'$ as a $G$-CW complex by \myref{Illman}. If $\tha$ is a quotient homomorphism with kernel $N$ and $H'$ is a subgroup of $G'$, then $H' = H/N$ for a subgroup $H$ of $G$ and $\tha^*(G'/H') = G/H$ so that no triangulations are required. Thus in both of these cases, $\tha^*$ takes generating $q$-cofibrations to $q$-cofibrations. Since $\theta^*$ is also a left adjoint in both cases, it preserves $q$-cofibrations in general.
\end{proof}

\begin{rem}
We did not require $\theta^*B$ to be a $\sG$-space in \myref{grpres}. 
However, if the kernel of $\tha$ is compact and $B$ is a $\sG'$-space, 
then $\theta^*B$ is a $\sG$-space. Indeed, $\tha$ is then
a proper map and $G_b=\theta^{-1}(G'_b)$ is compact since $G'_b$ is compact.
The restriction to compact kernels is the price we must 
pay in order to stay in the context of compact isotropy groups. 
We might instead consider $G'$-spaces $B$ such that the isotropy groups 
of both $B$ as a $G'$-space and $\tha^*B$ as a $G$-space are compact, 
but the assumption on $\tha^*B$ would be unnatural.  Note however that
one of the main reasons for restricting to compact isotropy groups is to 
obtain $G$-CW structures.  If $X$ is a $G'$-CW complex where $G' = G/N$ 
is a quotient group of $G$, then $\tha^*X$ is a $G$-CW complex with the same cells since the relevant orbits $G'/H'$ can be identified with $G/H$, 
where $H' = H/N$.
\end{rem}

For the $qf$-model structures, and to study adjunctions, it is convenient to consider quotient homomorphisms and inclusions separately. For the 
former, we consider the adjunctions of \myref{factor0}.

\begin{prop}\mylabel{fixedptrQa0}
Let $\epsilon\colon G\rtarr J$ be a quotient homomorphism of $G$ by a 
normal subgroup $N$. For a $G$-space $B$, 
consider the functors
\[(-)/N\colon G\sK_B \rtarr J\sK_{B/N} \ \ \text{and} \ \
(-)^N\colon G\sK_B \rtarr J\sK_{B^N}.\]
Let $j\colon B^N\rtarr B$ be the inclusion and $p\colon B\rtarr B/N$ be the quotient map. Then $((-)/N,p^*\epsilon^*)$ and $(j_!\epsilon^*,(-)^N)$ are Quillen adjoint pairs with respect to the $q$-model structures on 
both $G\sK_B$ and $J\sK_{B/N}$.  Let $\sC_G$ and $\sC_J$ be generating
sets of $G$-cell complexes and $J$-cell complexes. Consider $G\sK_B$
with the $qf(\sC_G)$-model structure and $J\sK_{B/N}$ and $J\sK_{B^N}$
with the $qf(\sC_J)$-model structure. Then 
\begin{enumerate}[(i)]
\item $((-)/N,p^*\epsilon^*)$ is a Quillen adjunction if $C/N\in\sC_J$ for $C\in \sC_G$. 
\vspace{1mm}
\item $(j_!\epsilon^*,(-)^N)$ is a Quillen adjunction if $\epz^*C\in \sC_G$ for $C\in \sC_J$. 
\end{enumerate}
\end{prop}

\begin{proof} Since $(j_!,j^*)$ and $(p_!,p^*)$ are Quillen adjoint pairs
in both the $q$ and the $qf$ contexts, it suffices to consider the case 
when $N$ acts trivially on $B$, so that $j$ and $p$ are identity maps. 
Then $\epz^*$ is right 
adjoint to $(-)/N$ and left adjoint to $(-)^N$. The properties of $\epz^*$
in the previous result give the conclusion for the $q$-model structures.
The functors $\epz^*$ and $(-)^N$ preserve $q$-equivalences. Since 
$$ \text{Map}_G(C,\epz^*f')\cong\text{Map}_G(C/N,f') \ \ 
\text{and} \ \ \text{Map}_J(C',f^N)\cong\text{Map}_G(\epsilon^*C',f)$$
for a $J$-map $f'$ and a $G$-map $f$, the conditions on generating sets in
(i) and (ii) ensure that $\epz^*$ and $(-)^N$ preserve the relevant $qf$-fibrations.
\end{proof} 

\begin{rem} In (i), we can take $\sC_J$ to consist of all finite products of quotients $C/N$ with $C\in \sC_G$ and orbits $J/H$ to arrange that $\sC_J$
be closed and contain these $C/N$. In (ii), we can take $\sC_G$ to 
consist of all products of pullbacks $\epz^*C$ for $C\in\sC_J$ with finite
products of orbits $G/H$. This set will be closed if $\sC_J$ is closed
since $\epz^*$ preserves products.
\end{rem}

Using \myref{fixedptrQa0} in conjunction with the additional change of
group relations of Propositions \ref{fixorbbase} and \ref{ouch0}, we
obtain the following compendium of equivalences in homotopy categories.

\begin{prop}\mylabel{orbfixdescend}
Let $A$ and $B$ be $G$-spaces. Let $j\colon B^N\rtarr B$ be the inclusion and $p\colon B\rtarr B/N$ be the quotient map, and let $f\colon A\rtarr B$ be 
a $G$-map.  Then, for ex-$G$-spaces $X$ over $A$ and $Y$ over $B$,
\begin{alignat*}{3}
&(p_!Y)/N \simeq Y/N, &\qquad\qquad
&(f_!X)/N \simeq (f/N)_!(X/N),\\
&(j^*Y)^N \simeq Y^N, &\qquad\qquad
&(f^*Y)^N \simeq (f^N)^*(Y^N), \\
&(p_*Y)^N \simeq Y/N, &\qquad\qquad
&(f_!X)^N \simeq (f^N)_!(X^N),
\end{alignat*}
where, for the last equivalence on the left, $B$ must be an $N$-free $G$-space.
\end{prop}

\begin{proof}
The equivalences displayed in the first line come from isomorphisms between Quillen left adjoints and are therefore clear. Similarly the equivalences in the second line come from isomorphisms between Quillen right adjoints. The first equivalence in the third line (in which we have changed notations from
\myref{ouch0}) comes from an isomorphism between a Quillen right adjoint on the left hand side, by \myref{Qad20}, and a Quillen left adjoint on the right hand side and therefore also descends directly to an equivalence on homotopy categories. For the last equivalence, note that $(-)^N$ preserves all $q$-equivalences and also preserves well-grounded ex-spaces and that $(f^N)_!$ preserves $q$-equivalences between well-grounded ex-spaces. Letting $Q$ and $R$ denote cofibrant and fibrant replacement functors, as usual, it follows that the maps
\[(R(f_!X))^N \ltarr (f_!X)^N \cong (f^N)_!(X^N) \ltarr (f^N)_!(Q(X^N)) \]
are $q$-equivalences on ex-spaces $X$ that are $qf$-fibrant and $qf$-cofibrant.
As noted in the proof of \myref{fixorbbase}, the point set level isomorphism $(f_!X)^N\iso (f^N)_!(X^N)$ is only valid for an ex-space $X$ whose section is a closed inclusion. However, if $X$ is $qf$-cofibrant, then it is compactly generated and this holds by \myref{coflemma}(i). Thus the equivalence holds 
in general in the homotopy category.
\end{proof}

The context for the next result is given in \myref{changes0} and \myref{ishriek}.

\begin{prop}\mylabel{Lishriek}
Let $\iota\colon H\rtarr G$ be the inclusion of a subgroup and let $A$ be 
an $H$-space. The adjoint equivalence $(\iota_!,\nu^*\iota^*)$ relating $H\sK_A$ and $G\sK_{\iota_!A}$ is a Quillen equivalence in the $q$-model structures and
also in the $qf(\sC_H)$ and $qf(\sC_G)$-model structures for any generating
sets $\sC_H$ and $\sC_G$ of $H$-cell complexes and $G$-cell complexes such
that $\io_!C = G\times_H C \in \sC_G$ for $C\in\sC_H$.
If $A$ is proper and completely regular, then the functor $\iota_!$ is 
also a Quillen right adjoint with respect to the $q$ and $qf$-model 
structures.
\end{prop}

\begin{proof}
Recall that $\nu\colon A\rtarr \iota^*\iota_!A= G\times_H A$ is the natural inclusion of $H$-spaces and that $(\nu_!,\nu^*)$ is a Quillen adjunction in both
the $q$ and $qf$ contexts. The functor $\io^*$ preserves $q$-equivalences and 
$q$-fibrations. It takes $qf(\sC_G)$-fibrations to $qf(\sC_H)$-fibrations
when $\io_!C \in \sC_G$ for $C\in\sC_H$ since
$$\text{Map}_H(C,\io^*f)\cong\text{Map}_G(\io_!C,f).$$

To show that $(\iota_!,\nu^*\iota^*)$ is a Quillen equivalence, we may
as well check the defining condition in the $q$-model structure.
Let $X$ be a $q$-cofibrant ex-$H$-space over $A$ and $Y$ be a $q$-fibrant ex-$G$-space over $\iota_!A$. Consider a $G$-map $f\colon \iota_!X\rtarr Y$. We must show that $f$ is a $q$-equivalence if and only if its adjoint $H$-map 
$\tilde{f}\colon X\rtarr \nu^*\iota^*Y$ is a $q$-equivalence. Since $\iota_!$ preserves acyclic $q$-cofibrations, we can extend $f$ to $f'\colon \iota_!RX\rtarr Y$, where $RX$ is a $q$-fibrant approximation. Since $f'$ is a $q$-equivalence if and only if $f$ is one, and similarly for their adjoints, we may assume without loss of generality that $X$ is $q$-fibrant. Recall from \myref{ishriek} that $\io_!$ and $\nu^*\io^*$ are inverse equivalences of categories and observe that  $\nu^*\iota^*Y$ can be viewed as the restriction, $Y|_A$, of $Y$ along the inclusion of $H$-spaces $\nu\colon A\rtarr G\times_H A$.  From that point of view, $\tilde{f}\colon X \rtarr \nu^*\iota^*Y$ is just the map $X\rtarr Y|A$ of ex-$H$-spaces over $A$ obtained by restriction of $\io^*f$ along $\nu$.

Now $f$ is a $q$-equivalence if and only if $f$ restricts to a $q$-equivalence $f_{[g,a]}$ on each fiber, meaning that this restriction is a weak equivalence after passage to fixed points under all subgroups of the isotropy group of $[g,a]$. For $a\in A$, the isotropy subgroup $H_a\subset H$ of $a$ coincides with the isotropy subgroup $G_{[e,a]}\subset G$ of $[e,a]\in G\times_H A$. For $g\in G$, the isotropy subgroup of $[g,a]$ is $gH_a g^{-1}$. Since the action by $g\in G$ induces a homeomorphism between the fibers over $[e,a]$ and over $[g,a]$, we see that $f$ is a $q$-equivalence if and only if each of the restrictions $f_{[e,a]}$ is a $q$-equivalence. But that holds if and only if $\tilde{f}$ is a $q$-equivalence.

For the last statement, recall the description of $\iota_!$ in 
\myref{iotaalt} as the composite $(p_*\pi^*\epz^*(-))^H$, where
$\epz\colon G\times H\rtarr H$ and  $\pi\colon G\times A\rtarr A$ are the projections and $p\colon G\times A\rtarr G\times_H A$ is the quotient map. 
Since $G\times A$ is completely regular, $p$ is a bundle with fiber 
$G/H_a$ over $[g,a]$, and $H_a$ is compact since $A$ is
proper. Therefore, by Propositions \ref{Qad20} and \ref{Qad202}, 
$p_*$ is a Quillen right adjoint with respect to the $q$ and $qf$-model structures.  In view of \myref{fixedptrQa0}, this displays $\iota_!$ as a composite of Quillen right adjoints.
\end{proof}

\begin{rem} We can take $\sC_G$ to consist of all finite products of 
the $\io_!C$ with $C\in \sC_H$ and orbits $G/K$ to arrange that 
$\sC_G$ be closed and contain these $\io_!C$. 
\end{rem}

We shall prove that $(\iota_!,\nu^*\iota^*)$ descends to a closed symmetric
monoidal equivalence of homotopy categories in \myref{imonoidaldescends} below.
The first statement of \myref{Lishriek} implies that the description of $\io^*$ in terms of base change that is given in \myref{ishriekb} descends to homotopy categories.

\begin{cor}\mylabel{LishriekCor}
The functor $\io^*\colon \text{Ho}G\sK_B\rtarr \text{Ho}H\sK_{\io^*B}$ is the composite
\[\xymatrix@1
{\text{Ho} G\sK_B \ar[r]^-{\mu^*} & \text{Ho}G\sK_{\io_!\io^*B} \htp \text{Ho} H\sK_{\io^*B}\\}\]
\end{cor}

\section{Fiber adjunctions and Brown representability}

For a point $b$ in $B$, we combine the special case 
$\tilde{b}\colon G/G_b\rtarr B$ of \myref{Qad10} with the special 
case $\io\colon G_b\rtarr G$ and $A = *$, hence $\nu\colon *\rtarr G/G_b$, 
of \myref{Lishriek} to obtain the following result concerning passage
to fibers. Recall from
\myref{Johann0} that the fiber functor $(-)_b\colon G\sK_B\rtarr G_b\sK_*$ is given by $\nu^*\iota^*\tilde{b}^*=b^*\iota^*$. By conjugation, its left adjoint $(-)^b$ therefore agrees with $\tilde{b}_!\iota_!$.

\begin{prop}\mylabel{FibadQ0}
For $b\in B$, the pair of functors $((-)^b,(-)_b)$ relating $G_b\sK_*$ and $G\sK_B$ is a Quillen adjoint pair.
\end{prop}

We use this result to verify the formal hypotheses of Brown's 
representability theorem \cite{Brown} for the category $\text{Ho}G\sK_B$. 
Of course, this verification is independent of the choice of model structure.
The category $G\sK_B$ has coproducts and homotopy pushouts, hence 
homotopy colimits of directed sequences. The usual constructions of homotopy
pushouts as double mapping cylinders and of directed homotopy colimits as 
telescopes makes clear that if the total spaces of their inputs are compactly 
generated, as they are after $q$-cofibrant approximation, then so are the 
total spaces of their outputs.  We need a few preliminaries.

\begin{defn}\mylabel{detectset}
For $n\geq 0$, $b\in B$, and $H\subset G_b$, let  $S^{n,b}_H$\noteindex{SHnb@$S^{n,b}_H$} be the ex-$G$-space 
$((G_b/H\times S^n)_+)^b$ over $B$. Explicitly, by \myref{Fibad}, $S^{n,b}_H = (G/H\times S^n)\amalg B$, with the obvious section and with the projection that maps $G/H\times S^n$ to the point $b$ and maps $B$ by the identity map.  Equivalently, taking $d$ to be the constant map at $b$, $S^{n,b}_H$ is the quotient ex-$G$-space associated to the generating cofibration $i(d)$,  $i\colon G/H\times S^{n-1}\rtarr G/H\times D^n$. Therefore, $S^{n,b}_H$ is cofibrant in both the $q$ and the $qf$-model structures. Let $\sD_B$\noteindex{DB@$\sD_B$} be the ``detecting set'' of all such ex-$G$-spaces $S^{n,b}_H$.
\end{defn}

Let $[X,Y]_{G,B}$ denote the set of maps $X\rtarr Y$ in $\text{Ho}G\sK_B$. 

\begin{lem} Each $X$ in $\sD_B$ is compact, in the sense 
that
$$\text{colim}\, [X, Y_n]_{G,B}\iso [X, \text{hocolim}\,Y_n]_{G,B}$$
for any sequence of maps $Y_n\rtarr Y_{n+1}$ in $G\sK_B$.  
\end{lem}
\begin{proof} If $X = S^{n,b}_H$, then 
$[X,Y]_{G,B} \iso [G_b\times S^n)_+,Y_b]_{G_b}$. 
In $G_b\sK_*$, every object is fibrant and the target is the set of 
homotopy classes of $G_b$-maps $(G_b\times S^n)_+ \rtarr Y_b$, which
is the set of unbased nonequivariant homotopy classes of maps $S^n\rtarr Y_b$.
Using cofibrant replacement, we can arrange that the $(Y_n)_b$ have total spaces 
in $\sU$. Then the conclusion follows from \myref{little}.
\end{proof}

The following result says that the set $\sD_B$ detects $q$-equivalences.

\begin{prop} A map $\xi\colon Y\rtarr Z$ in $G\sK_B$ is a $q$-equivalence 
if and only if the induced map $\xi_*\colon [X,Y]_{G,B}\rtarr [X,Z]_{G,B}$ is a bijection for all $X\in \sD_B$. 
\end{prop} 
\begin{proof} We may assume that $Y$ and $Z$ are fibrant. By the evident long exact sequences of homotopy groups and the five lemma, $\xi$ 
is a $q$-equivalence if and only if each $Y_b\rtarr Z_b$ is a $q$-equivalence.
This is detected by the based $G_b$-spaces $(G_b/H\times S^n)_+$ and the 
conclusion follows by adjunction.
\end{proof}

\begin{thm}[Brown]\mylabel{brown0} A contravariant set-valued functor on the
category $\text{Ho}G\sK_B$ is
representable if and only if it satisfies the wedge and Mayer-Vietoris
axioms.
\end{thm}

 \chapter{Ex-fibrations and ex-quasifibrations}

To complete the foundations of parametrized homotopy theory,
we are faced with two problems that were discussed in the 
Prologue. In our preferred $qf$-model structure, the base change 
adjunction $(f_!, f^*)$ is a Quillen pair for any map $f$ and is a
Quillen equivalence if $f$ is an equivalence. As shown by \myref{noway}, 
this implies that the base change adjunction $(f^*,f_*)$ cannot be a 
Quillen adjoint pair. Some such defect must hold for any model structure. 
Therefore, we cannot turn to model theory to construct the functor $f_*$ 
on the level of homotopy categories.  The same counterexample illustrates that passage to derived functors is not functorial in general, so that a relation between composites of functors that holds on the point-set level need not 
imply a corresponding relation on homotopy categories.

In any attempt to solve those two problems, one runs into a third one that concerns a basic foundational problem in ex-space theory. Model theoretical considerations lead to the use of Serre fibrations as projections, or to the even weaker class of $qf$-fibrations. However, only Hurewicz fibrations are considered in most of the literature.  There is good reason for that. 
Fiberwise smash products, suspensions, cofibers, function spaces, and other fundamental constructions in ex-space theory do not preserve Serre fibrations.

The solutions to all three problems are obtained by the use of ex-fibrations.
Recall that these are the well-sectioned $h$-fibrant ex-spaces. We study
their properties in \S8.1. They seem to give the definitively right kind of ``fibrant ex-space'' from the point of view of classical homotopy theory, and they behave much better under the cited constructions than do Serre fibrations,
as we show in \S8.2.  Many variants of this notion appear in the literature.  Precisely this variant, with this name, appears in Monica Clapp's paper \cite{Clapp}, and we are indepted to her work for an understanding of the centrality of the notion. Perversely, as we noted in \myref{guess}, it is unclear how it fits into the model categorical framework.

We construct an elementary ex-fibrant approximation functor in \S8.3. It plays a key role in bridging the gap between the model theoretic and classical worlds. In a different context, the classification of sectioned fibrations, the first author introduced this construction in \cite[\S5]{May}.  We record
some its properties in \S8.4. 

We define quasifibrant ex-spaces and ex-quasifibrations and show that they inherit some of the good properties of ex-fibrations in \S8.5. They will 
play a key role in the stable theory.

Everything in this chapter works just as well equivariantly as nonequivariantly
for any topological group $G$ of equivariance.

\section{Ex-fibrations}

Under various names, the following notions were in common use in the 1970's.
We shall see shortly that these definitions agree with those given in \myref{names}. 

\begin{defn}\mylabel{well} 
Let $(X,p,s)$ be an ex-space over $B$. 
\begin{enumerate}[(i)]
\item $(X,p,s)$ is \emph{well-sectioned}\index{ex-space!well-sectioned}\index{well-sectioned ex-space} if $s$ is a closed inclusion and there is a retraction 
$$\rh\colon X\times I\rtarr X\cup_B(B\times I) = Ms$$
over $B$.
\item $(X,p,s)$ is \emph{well-fibered}\index{ex-space!well-fibered}\index{well-fibered ex-space} if there is a coretraction, or \emph{path-lifting function},\index{path-lifting function}
$$\io\colon Np = X\times_B B^I \rtarr X^I$$
under $B^I$, where $B^I$ maps to $Np$ via $\al \rtarr (s\al(0),\al)$. 
\item $(X,p,s)$ is an \emph{ex-fibration}\index{ex-fibration} if it is both well-sectioned 
and well-fibered.
\end{enumerate}
\end{defn}

The requirement in (i) that the retraction $\rh$ be a map over $B$ ensures that it restricts on fibers to a retraction that exhibits the nondegeneracy of the basepoint $s(b)$ in $X_b$ for each $b\in B$.  In view of \myref{ffpmodel}(i), 
we have the following characterization of well-sectioned ex-spaces, in 
agreement with \myref{names}. 

\begin{lem}\mylabel{charcof} 
An ex-space $X$ is well-sectioned if and only if $X$ is $\bar{f}$-cofibrant.
\end{lem}

We use the term ``well-sectioned'' since it goes well with ``well-based''.
The category of well-sectioned ex-spaces is the appropriate para\-me\-trized generalization of the category of well-based spaces, and restricting to well-sectioned ex-spaces is analogous to restricting to well-based spaces. 

Note that the section of $X$ provides a canonical way of lifting a path in $B$
that starts at $b$ to a path in $X$ that starts at $s(b)$. The requirement in \myref{well}(ii) that the path-lifting function $\io$ be a map under $B^I$ says that  $\io(s\al(0),\al)(t) = s(\al(t))$ for all $\al\in B^I$ and $t\in I$. That is, $\iota$ is required to restrict to the canonical lifts provided by the section, so that paths in $X$ that start in $s(B)$ remain in $s(B)$. In 
contrast with \myref{charcof}, the well-fibered condition does not by itself fit naturally into the model theoretic context of Chapter 5.  However, we have the following characterization of ex-fibrations, which again is in agreement with the original definition we gave in \myref{names}.

\begin{lem}\mylabel{charexfib}
If $X$ is well-fibered, then $X$ is $h$-fibrant.
If $X$ is well-sectioned, then $X$ is an ex-fibration if and only 
if $X$ is $h$-fibrant.
\end{lem}
\begin{proof}
The first statement is clear since the coretraction $\io$ is a path-lifting function. This gives the forward implication of the second statement, and the
converse is a special case of the following result of Eggar \cite[3.2]{Eggar}.
\end{proof}

\begin{lem}\mylabel{Eggar}
Let $i\colon X\rtarr Y$ be an $\bar{f}$-cofibration of 
ex-spaces over $B$, where $Y$ is $h$-fibrant. Then any map 
$\io\colon X\times_B B^I\rtarr Y^I$ such that the composite
$$\xymatrix@1{
X\times_B B^I \ar[r]^-{\io} & Y^I \ar[r] & Y\times_B B^I\\}$$
is the inclusion can be extended to a coretraction $Y\times_B B^I\rtarr Y^I$.
\end{lem}
\begin{proof} 
The inclusion $X\times_B B^I \rtarr Y\times_B B^I$ is an 
$\bar{h}$-cofibration by \myref{StrPull}. Therefore 
there is a lift ${\nu}$ in the diagram
$$\xymatrix{
(Y\times_B B^I)\times\{0\} \cup (X\times_B B^I)\times I \ar[r]^-{f} \ar[d]
& Y \ar[d]\\
(Y\times_B B^I)\times I \ar[r]_-{g} \ar@{-->}[ur]^{\nu} & B,\\}$$
where $f(y,\om,0) = y$, $f(x,\om,t) = \io(x,\om)(t)$, and $g(y,\om,t)= \om(t)$.
The adjoint $Y\times_B B^I\rtarr Y^I$ of $\nu$ is the required extension to a 
coretraction.
\end{proof}

\begin{rem}
We comment on the terminology.

(1) We are following \cite{CJ, James} and others in saying that an 
$\bar{f}$-cofibrant ex-space is well-sectioned; the term ``fiberwise well-pointed'' is also used.  For a based space, the terms ``nondegenerately based'' and ``well-based'' or ``well-pointed'' are used interchangeably to mean that the inclusion of the basepoint is an $h$-cofibration. In contrast, for an ex-space, the term ``fiberwise nondegenerately pointed'' is used in \cite{CJ, James} to indicate a somewhat weaker condition than well-sectioned.

(2) The term ``well-fibered'' is new but goes naturally with well-sec\-tion\-ed. The concept itself is old. We believe that it is due to Eggar \cite[3.3]{Eggar}, who calls a coretraction under $B^I$ a {\em special lifting function}.

(3) Becker and Gottlieb \cite{BG1} may have been the first to use the term ``ex-fibration'', but for a slightly different notion with sensible CW restrictions. As noted in the introduction, precisely our notion is used by Clapp \cite{Clapp}.  Earlier, in \cite[\S5]{May} and \cite{May1}, the first author called ex-fibrations ``$\sT$-fibrations'', and he studied their classification and their fiberwise localizations and completions. The equivariant generalization appears in Waner \cite{Waner}. A more recent treatment of the classification of ex-fibrations has been given by Booth \cite{Boothbrag}. 
\end{rem}

\section{Preservation properties of ex-fibrations}

We have a series of results that show that ex-fibrations behave well with respect to standard constructions. In some of them, one must use the equivariant version of \myref{fNDR} to verify that the given construction preserves well-sectioned objects. In all of them, if we only assume that the input ex-spaces are well-sectioned, then we can conclude that the output ex-spaces are well-sectioned.  It is the fact that the given constructions preserve well-fibered objects that is crucial. Few if any of these results hold with Serre rather 
than Hurewicz fibrations as projections. 

\begin{prop}\mylabel{pres1}
Ex-fibrations satisfy the following properties.
\begin{enumerate}[(i)]
\item A wedge over $B$ of ex-fibrations is an ex-fibration.
\item If $X$, $Y$ and $Z$ are ex-fibrations and $i$ is an $\bar{f}$-cofibration in the following pushout diagram of ex-spaces over $B$, then $Y\cup_X Z$ is an ex-fibration.
$$\xymatrix{
X \ar[r]^{i} \ar[d] & Y \ar[d] \\
Z \ar[r] & Y\cup_X Z}$$
\item The colimit of a sequence of $\bar{f}$-co\-fi\-bra\-tions $X_i\rtarr X_{i+1}$  between ex-\-fi\-bra\-tions is an ex-fibration.
\end{enumerate}
If the input ex-spaces are only assumed to be well-sectioned, then the output ex-spaces are well-sectioned. 
\end{prop}
\begin{proof}  The last statement is clear. Using it, we see that the
colimits in (i), (ii), and (iii) are well-sectioned, hence it suffices
to prove that they are $h$-fibrant.  This is done by constructing path 
lifting functions for the colimits from path lifting functions for their
inputs.  In (i), we start with path lifting functions under $B^I$ for 
the wedge summands and see that they glue together to define a path 
lifting function under $B^I$ for the wedge.  Part (ii) is due to Clapp \cite[1.3]{Clapp}, and we omit full details. She starts with a path lifting function for $X$ and uses \myref{Eggar} to extend it to a path lifting
function for $Y$. She also starts with a path lifting function for $Z$.
She then uses a representation $(h,u)$ of $(Y,X)$ as a fiberwise NDR pair
to build a path lifting function for the pushout from the given path
lifting function for $Z$ and a suitably deformed version of the path
lifting function for $Y$.  In (iii), \myref{Eggar} shows that we can
extend a path lifting function for $X_i$ to a path lifting function
for $X_{i+1}$.  Inductively, this allows the construction of compatible
path lifting functions for the $X_i$ that glue together to give a path 
lifting function for their colimit.
\end{proof}

Although of little use to us, since the $f$-homotopy category is not the right one for our purposes, many of our adjunctions give Quillen adjoint pairs with respect to the $f$-model structure.  For example, the following result, which should be compared with \myref{Qad10}, implies that $(f_!,f^*)$ is a Quillen adjoint pair in the $f$-model structures and that it is a Quillen equivalence 
if $f$ is an $h$-equivalence.

\begin{prop}\mylabel{fexpres}
Let $f\colon A\rtarr B$ be a map, let $X$ be a well-sectioned ex-space over $A$, and let $Y$ be a well-sectioned ex-space over $B$. Then $f_!X$ and $f^*Y$ are well-sectioned. If $Y$ is an ex-fibration, then so is $f^*Y$, and the functor $f^*$ preserves $f$-equivalences. If $f$ is an $h$-equivalence, then $(f_!,f^*)$ induces an equivalence of $f$-homotopy categories.
\end{prop}
\begin{proof}
It is easy to check that representations of $(X,A)$ and $(Y,B)$ as fiberwise NDR-pairs induce representations of $(f_!X,B)$ and $(f^*Y,A)$ as fiberwise NDR-pairs.  As a pullback, the functor $f^*$ preserves both $f$-fibrant and $h$-fibrant ex-spaces, and $f^*$ preserve $f$-equivalences since it preserves $f$-homotopies. For the last statement, if $f$ is a homotopy equivalence with homotopy inverse $g$, then standard arguments with the CHP imply that $f^*$ induces an equivalence of $f$-homotopy categories with inverse $g^*$; see, for example, \cite[2.5]{May}. It follows that $g^*$ is equivalent to $f_!$ and that $(f_!,f^*)$ is a Quillen equivalence.
\end{proof}

The following result appears in \cite{Eggar} and \cite[3.6]{May}. It also leads to a Quillen adjoint pair with respect to the $f$-model structure; compare \myref{smaAB}.

\begin{prop}\mylabel{Hursma}
Let $X$ and $Y$ be well-sectioned ex-spaces over $A$ and $B$. Then 
$X\barwedge Y$ is a well-sectioned ex-space over $A\times B$. If $X$ 
and $Y$ are ex-fibrations, then $X\barwedge Y$ is an ex-fibration.
\end{prop}
\begin{proof}
Representations of $(X,A)$ and $(Y,B)$ as fiberwise NDR-pairs determine a representation of $(X\barwedge Y, A\times B)$ as a fiberwise NDR-pair, by standard formulas \cite[p.\,43]{Concise}. Similarly, path lifting functions for $X$ and $Y$ can be used to write down a path lifting function for 
$X\barwedge Y$.
\end{proof}

\begin{cor}
If $X$ and $Y$ are ex-fibrations over $B$, then so is $X\sma_B Y$.
\end{cor}

\begin{cor}\mylabel{HursmaK}
If $X$ is an ex-fibration over $B$ and $K$ is a well-based space, then $X\sma_B K$ is an ex-fibration over $B$.
\end{cor}

\begin{prop}\mylabel{savior}
Let $X$ and $Y$ be well-sectioned and let $f\colon X\rtarr Y$ be an ex-map that is an $h$-equivalence. Then $f\sma_B\text{id}\colon X\sma_B Z\rtarr Y\sma_B Z$ is an $h$-equivalence for any ex-fibration $Z$. In particular, $f\sma_B\text{id}\colon X\sma_B K\rtarr Y\sma_B K$ is an $h$-equivalence for any well-based space $K$. 
\end{prop}
\begin{proof}
As observed by Clapp \cite[2.7]{Clapp}, this follows from the gluing lemma by comparing the defining pushouts.
\end{proof}

As in ordinary topology, function objects work less well, but we do have the following analogue of \myref{HursmaK}.

\begin{prop}
If $X$ is an ex-fibration over $B$ and $K$ is a compact well-based space, then $F_B(K,X)$ is an ex-fibration over $B$.
\end{prop}
\begin{proof}
Let $(h,u)$ represent $(X,B)$ as a fiberwise NDR-pair. Then $(j,v)$ represents $(F_B(K,X),B)$ as a fiberwise NDR-pair, where 
\[v(f) = \text{sup}_{k\in K}u(f(k)) \qquad\text{and}\qquad 
j_t(f)(k) = h_t(f(k))\]
for $f\in F_B(K,X)$. Note for this that $F_B(K,B)=B$ and that, by \myref{Ffib}, $F_B(K,X)$ is $h$-fibrant.
\end{proof}

\section{The ex-fibrant approximation functor}

We describe an elementary ex-fibrant replacement functor ${P}$. It is just 
the composite of a whiskering functor $W$ with a version of the mapping path fibration functor $N$. The functor $P$ replaces ex-spaces by naturally $h$-equivalent ex-fibrations. From the point of view of model theory, ${P}$ can be thought of as 
a kind of $q$-fibrant replacement functor that gives Hurewicz fibrations rather than just Serre fibrations as projections.  The nonequivariant version of $P$ appears in \cite[5.3, 5.6]{May}, and the equivariant version appears in \cite[\S3]{Waner}. With motivation from the theory of transports in fibrations, those sources work with Moore paths of varying length.  Surprisingly, that choice turns out to be essential for the construction to work.

We therefore begin by recalling that the space of \emph{Moore paths}\index{Moore paths} in $B$ is given by\noteindex{LB@$\Lambda B$}
\[\Lambda B = \{(\la, l)\in B^{[0,\infty]}\times [0,\infty) \mid \text{$\la(r)=\la(l)$ for $r\geq l$}\}\]
with the subspace topology. We write $\la$ for $(\la,l)$ and
$l_{\la}$ for $l$, which is the length of $\la$. 
Let $e\colon \Lambda B \rtarr B$ be the endpoint projection $e(\lambda)=\lambda(l_{\la})$. The composite of Moore paths $\mu$ and $\lambda$ such that $\lambda(l_\lambda)=\mu(0)$ is defined by $l_{\mu\lambda}=l_\mu+l_\lambda$ and
\[(\mu\lambda)(r)=\begin{cases}
\lambda(r) &\text{if $r\leq l_\lambda$},\\
\mu(r-l_\lambda) &\text{if $r\geq l_\lambda$}.
\end{cases}\]
Embed $B$ and $B^I$ in $\Lambda B$ as the paths of length $0$ and $1$. 
For a Moore path $\lambda$ in $B$ and real numbers $u$ and $v$ such that
$0\leq u\leq v$, let $\lambda|_u^v$ denote the Moore path $r\mapsto \lambda(u+r)$ of length $v-u$.

\begin{defn}
Consider an ex-space $X = (X,p,s)$ over $B$.
\begin{enumerate}[(i)]
\item Define the \emph{whiskering functor}\index{whiskering functor}\index{functor!whiskering --} $W$\noteindex{WX@$WX$} by letting
\[WX=(X\cup_B (B\times I),q,t),\]
where the pushout is defined with respect to $i_0\colon B\rtarr B\times I$. 
The projection $q$ is given by the projection $p$ of $X$ and the projection 
$B\times I\rtarr B$, and the section $t$ is given by $t(b)=(b,1)$.
\item Define the \emph{Moore mapping path fibration functor}\index{Moore mapping path fibration}\index{functor!Moore mapping path fibration --} $L$\noteindex{LX@$LX$} by letting
\[LX=(X\times_B \Lambda B, q,t),\]
where the pullback is defined with respect to the map $\LA B\rtarr B$
given by evaluation at $0$. The projection $q$ is given by $q(x,\lambda)=e(\lambda)$ and the section $t$ is given by 
$t(b)=(s(b), b)$, where $b$ is viewed as a path of length $0$.
\end{enumerate}
\end{defn}

Thus $WX$ is obtained by growing a whisker on each point in the section of $X$, and the endpoints of the whiskers are used to give $WX$ a section. Similarly, $LX$ is obtained by attaching to $x\in X$ all Moore paths in $B$ starting at $p(x)$. The endpoints of the paths give the projection.  In the language of \S\ref{sec:towardh}, $WX$ is the standard mapping cylinder construction of the section of $X$, thought of as a map in $G\sK/B$. The section $t$ of $WX$ is just the $f$-cofibration in the standard factorization $\rho \circ t$ of $s$ through its mapping cylinder. In particular, $WX$ is well-sectioned. Similarly, $LX$ is a modification of the mapping path fibration $Np$ in $G\sK$. The projection $p$ of $X$ factors through the projection $q$ of $LX$, which is an $h$-fibration; a path lifting function $\xi\colon LX\times_B B^I \rtarr (LX)^I$ is given by
$\xi((x,\lambda),\gamma)(t)=(x,\gamma|_0^t\lambda)$. Thus $LX$ is $h$-fibrant, but it need not be well-fibered. 

We can display all of this conveniently in the following diagram. The third
square on the top is a pushout and the second square on the bottom is a pullback. That defines the maps $\ph$ and $\pi$, and the maps $\rho$ and 
$\iota$ are induced by the universal properties from the identity map of $X$.
\[\xymatrix@=.8cm{
&&&B\ar[d]_{i_1}\ar@{=}[dr]\\
B \ar@{=}[r]\ar[d] & B \ar@{=}[r]\ar[d] & B \ar[d]^s \ar[r]_-{i_0} & B\times I \ar[d]\ar[r]_-{\text{pr}} & B\ar[d] \\
X \ar[d]\ar@{-->}[r]^{\iota} & LX \ar[r]^-\pi\ar[d] & X \ar[d]^p \ar[r]^-\phi & WX \ar[d]\ar@{-->}[r]^-{\rho} & X \ar[d] \\
B\ar@{=}[dr]\ar[r] & \Lambda B \ar[d]^{e}\ar[r]^-{p_0} & B\ar@{=}[r] & B \ar@{=}[r] & B\\
& B}\]
Thus $\rho$ projects whiskers on fibers to the original basepoints and $\io$
is the inclusion $x\mapsto (x,p(x))$, where $p(x)$ is the path of length zero. Note that $\phi$ is not a map under $B$ and $\pi$ is not a map over $B$. They give an inverse $f$-equivalence to $\rho$ and an inverse $h$-equivalence to $\iota$.

\begin{prop}\mylabel{fpequiv2}
The map $\rho\colon  WX\rtarr X$ is a natural $f$-equivalence of ex-spaces and $WX$ is well-sectioned. The map $\iota\colon X\rtarr LX$ is a natural $h$-equivalence of ex-spaces and $LX$ is $h$-fibrant. Therefore $W$ takes $f$-equivalences to $fp$-equivalences and $L$ takes $h$-equivalences to $f$-equivalences.
\end{prop}

The last statement follows from \myref{reverse}. We think of $\rho$ and 
$\iota$ as giving a \emph{well-sectioned approximation}\index{well-sectioned approximation}\index{approximation!well-sectioned --} and an 
\emph{$h$-fibrant approximation}\index{fibrant approximation@$h$-fibrant approximation}\index{approximation!fibrant@$h$-fibrant --} in the category of ex-spaces. We will combine them to obtain the promised ex-fibrant approximation, but we first insert a technical lemma. 

\begin{lem}\mylabel{closedP}
If $X$ is an ex-space with a closed section, then $WLX$ is an ex-fibration. 
If $X$ is well-fibered, then $WX$ is an ex-fibration.
\end{lem}
\begin{proof}
A path lifting function $\xi\colon NWLX = WLX\times_B B^I \rtarr (WLX)^I$ for $WLX$ is obtained by letting
\[\xi(z,\gamma)(t) = 
\begin{cases}
(x,\gamma|_0^t\lambda)\in LX & \text{if $z=(x,\lambda)\in LX$},\\
(\gamma(t),u-t)\in B\times I & \text{if $z=(b,u)$ and $t\leq u$},\\
(s(\gamma(u)), \gamma|_u^t)\in LX & \text{if $z=(b,u)$ and $t\geq u$.}
\end{cases}\]
It is easy to verify that, as a map of sets, $\xi$ gives a well-defined section of the canonical retraction $\pi\colon (WLX)^I\rtarr WLX\times_B B^I$. Continuity is a bit more delicate, but if the section of $X$ is closed, then one verifies that
\[\PH=\{(z,\gamma)  \mid \text{$z$ is the equivalence class of $(s(b),b)\sim (b,0)$}\}\]
is a closed subset of WLX and hence $N\PH$ is a closed subset of $NWLX$.
To see the implication, note that $(-)\times B^I$ preserves closed inclusions and $Z\times_B B^I\subset Z\times B^I$ is a closed  inclusion because $B$ is in $\sU$ (see \myref{Umod}). 
Continuity follows since we are then piecing together continuous 
functions on closed subsets.

If $X$ is well-fibered and $\xi\colon X\times_B B^I\rtarr X^I$ is a path-lifting function under $B^I$, we can define a path lifting function $\bar\xi\colon  WX\times_B B^I \rtarr (WX)^I$ for $WX$ by
\[\bar\xi(x,\gamma)=\begin{cases}
\xi(x,\gamma) & \text{if $x\in X$},\\
(\gamma,u) & \text{if $x=(b,u)$}.
\end{cases}\]
To check that $\bar\xi$ is continuous, we use the fact that the functor
$N(-)=B^I\times_B (-)$ commutes with pushouts to write $NWX$ as a pushout.
We then see that $\bar{\xi}$ is the map obtained by passage to pushouts
from a pair of continuous maps.
\end{proof}

Recall that the sections of ex-spaces in $G\sU_B$ are closed, by \myref{coflemma}.  Since we shall only need to apply the constructions 
of this section to ex-spaces in $G\sU_B$, the closed section 
hypothesis need not concern us. 

\begin{defn}\mylabel{exfibapp}
Define the \emph{ex-fibrant approximation functor}\index{approximation!ex-fibrant --} $P$\noteindex{PX@$PX$} by the natural 
zig-zag of $h$-equivalences $\ph = (\rh,W\io)$ displayed in the diagram
$$\xymatrix{ 
X & WX \ar[l]_-{\rh} \ar[r]^-{W\io} & WLX = PX.}
$$
By \myref{fpequiv2}, $P$ takes $h$-equivalences between arbitrary ex-spaces to $fp$-equivalences. If $X$ has a closed section, then $PX$ is an ex-fibration. If $X$ is an ex-fibration, then it has a closed section, and the above display is a natural zig-zag of $fp$-equivalences between ex-fibrations.
\end{defn}

\section{Preservation properties of ex-fibrant approximation}

One advantage of ex-fibrant approximation over $q$ or $qf$-fibrant approximation is that there are explicit commutation natural transformations relating it to many constructions of interest. The following result is an elementary illustrative example. 

\begin{lem}
Let $\sD$ be a small category, $X\colon \sD\rtarr G\sK_B$ be a functor, and
$$ \om\colon \text{colim}\, WX_d \rtarr W\text{colim}\, X_d  \ \ \ \text{and}\ \ \
\nu\colon \text{colim}\, LX_d\rtarr L\text{colim}\, X_d $$
be the evident natural maps. Then $\om$ is a map over $\text{colim}\, X_d$ and $\nu$ is a map under $\text{colim}\, X_d$, so that the following diagrams commute. All maps in these diagrams are $h$-equivalences.
$$\xymatrix{
\text{colim}\, WX_d \ar[rr]^{\om} \ar[dr]_{\text{colim}\,\rh}\, & & W\text{colim}\, X_d \ar[dl]^{\rh} \\
& \text{colim}\, X_d  & }$$
$$\xymatrix{
& \text{colim}\, X_d  \ar[dl]_{\text{colim}\,\io}   \ar[dr]^{\io}  &  \\
\text{colim}\, LX_d \ar[rr]_{\nu} & & L\text{colim}\, X_d}$$
Let $\mu = W\nu\com \om\colon \text{colim}\, PX_d\rtarr P\text{colim}\, X_d$.  Then the following diagram of $h$-equivalences commutes.
$$\xymatrix{
\text{colim}\, X_d\ar@{=}[d] & \text{colim}\, WX_d \ar[l]_-{\text{colim}\, \rh} 
\ar[r]^-{\text{colim}\, W\io} \ar[d]^{\om} & \text{colim}\, PX_d \ar[d]^{\mu}\\
\text{colim}\, X_d   &  W\text{colim}\, X_d \ar[l]^-{\rh} \ar[r]_-{W\io} & P\text{colim}\, X_d}$$
The analogous statements for limits also hold. 
\end{lem}
\begin{proof}
This is clear from the construction of limits and colimits in
\myref{btopbicomp}. The relevant $h$-equivalences of total 
spaces are natural and piece together to pass to limits and colimits.
\end{proof}

\begin{warn}\mylabel{DoubleTrouble}
We would like an analogue of the previous result for tensors. In particular, we would like a natural map $(LX)\sma K\rtarr L(X\sma K)$ under $X\sma K$ for ex-spaces $X$ over $B$ and based spaces $K$.  Inspection of definitions makes clear that there is no such map. The obvious map that one might write down, as in the erroneous \cite[5.6]{May}, is not well-defined.  In Part III, this complicates the extension of $P$ to a functor on spectra over $B$.
\end{warn}

\begin{lem}\mylabel{munu} 
Let $f\colon A\rtarr B$ be a map.
\begin{enumerate}[(i)]
\item Let $X$ be an ex-space over $A$. Then there are natural maps
$$\om\colon f_!WX \rtarr Wf_!X \ \ \ \text{and} \ \ \ \nu\colon f_!LX\rtarr Lf_!X$$
of ex-spaces over $B$ such that $\om$ is a map over $f_!X$ and $\nu$ is a map under $f_!X$. Let $\mu = W\nu\com \om \colon f_!PX\rtarr Pf_!X$. Then the following diagram commutes.
$$\xymatrix{
f_!X\ar@{=}[d] & f_!WX \ar[l]_-{f_!\rh} 
\ar[r]^-{f_!W\io} \ar[d]^{\om} & f_!PX \ar[d]^{\mu}\\
f_!X   &  Wf_!X \ar[l]^-{\rh} \ar[r]_-{W\io} & Pf_!X }$$

\item Let $Y$ be an ex-space over $B$. Then there are natural maps
$$\om\colon Wf^*Y \rtarr f^*WY \ \ \text{and} \ \ \nu\colon  Lf^*Y\rtarr f^*LY$$ of ex-spaces over $A$, the first an isomorphism, such that $\om$ is a map over $f^*Y$ and $\nu$ is a map under $f^*Y$. Let $\mu = \om\com W\nu \colon Pf^*Y\rtarr f^*PY$. Then the following diagram commutes.
$$\xymatrix{
f^*Y\ar@{=}[d] & Wf^*Y \ar[l]_-{\rh} 
\ar[r]^-{W\io} \ar[d]^{\om} & Pf^*Y \ar[d]^{\mu}\\
f^*Y   & f^*WY \ar[l]^-{f^*\rh} \ar[r]_-{f^*W\io} & f^*PY}$$
If $Y$ is an ex-fibration, then $\mu$ is an $fp$-equivalence.

\item Let $X$ be an ex-space over $A$. Then there are natural maps
$$\om\colon Wf_*X \rtarr f_*WX \ \ \ \text{and} \ \ \ \nu\colon Lf_*X\rtarr f_*LX$$
of ex-spaces over $B$ such that $\om$ is a map over $f_*X$ and $\nu$ is a map under $f_*X$. Let $\mu = \om\com W\nu \colon Pf_*X\rtarr f_*PX$. Then the following diagram commutes.
$$\xymatrix{
f_*X\ar@{=}[d] & Wf_*X \ar[l]_-{\rh} 
\ar[r]^-{W\io} \ar[d]^{\om} & Pf_*X \ar[d]^{\mu}\\
f_*X   &  f_*WX \ar[l]^-{f_*\rh} \ar[r]_-{f_*W\io} & f_*PX}$$
\end{enumerate}
\end{lem}
\begin{proof} Again, the proof is by inspection of definitions.
Since $f_!$ does not preserve ex-fibrations, we do not have an 
analogue for $f_!$ of the last statement about $f^*$ in (ii).
\end{proof}

\begin{warn}\mylabel{PffP}
We offer another example of the technical dangers lurking in this subject. The maps $\mu$ in the previous proposition are {\em not} $h$-equi\-va\-lences in general, the problem in (ii), say, being that $f^*$ does not preserve $h$-equivalences in general.  If $\mu\colon Pf^*Y\rtarr f^*PY$ were always an $h$-equivalence, then one could prove by the methods in \S9.3 below that the relations (\ref{bases0}) descend to homotopy categories for all pullbacks of the form displayed in \myref{Mackey0}. In view of \myref{noway}, that conclusion is false. This is another pitfall we fell into, and it invalidated much work in an earlier draft.
\end{warn}

\section{Quasifibrant ex-spaces and ex-quasifibrations}

By analogy with the fact that an ex-fibration is a well-sectioned 
$h$-fibrant ex-space, we adopt the following terminology.

\begin{defn}
An ex-space $X$ is \emph{quasifibrant}\index{ex-space!quasifibrant} if its projection $p$ is a quasifibration. An \emph{ex-quasifibration}\index{ex-quasifibration} is a well-sectioned
quasifibrant ex-space. 
\end{defn}

If $X$ is quasifibrant, there is a long exact sequence of homotopy groups
\[\cdots\rtarr \pi^H_{q+1}(B,b)\rtarr \pi^H_q(X_b,x)\rtarr \pi^H_q(X,x) \rtarr \pi^H_q(B,b) \rtarr \cdots \rtarr \pi^H_0(B,b)\]
for any $b\in B$, $x\in X_b$ and $H\subset G_b$.  Using this and the long exact sequences of 
the pairs $(X,X_b)$, five lemma comparisons give the following observations.

\begin{lem}\mylabel{quasichar} Let $f\colon X\rtarr Y$ be a 
$q$-equivalence of ex-spaces over $B$. Then each map of fibers 
$f\colon X_b\rtarr Y_b$ is a $q$-equivalence if and only if each map
of pairs $f\colon (X,X_b)\rtarr (Y,Y_b)$ is a $q$-equivalence. If $X$ 
and $Y$ are quasifibrant, then these maps of pairs are $q$-equivalences. Conversely, if these maps of pairs are $q$-equivalences and either $X$ or 
$Y$ is quasifibrant, then so is the other.
\end{lem}

Working in $G\sU_B$, we obtain the following result. The same pattern of 
proof gives many other results of the same nature that we leave to the reader.

\begin{prop}\mylabel{quasicof}
The following statements hold.
\begin{enumerate}[(i)]
\item A wedge over $B$ of ex-quasifibrations is an ex-quasifibration.
\item If $f:X\rtarr Y$ is a map such that $X$ is an ex-quasifibration 
and $Y$ is quasifibrant, then the cofiber $C_B f$ is quasifibrant.
\item If $X$ is an ex-quasifibration and $K$ is a well-based space, 
then $X\sma_B K$ is an ex-quasifibration.
\end{enumerate}
\end{prop}
\begin{proof}
This follows from \myref{quasichar}, the natural zig-zag
\[\xymatrix{X & WX \ar[l]\ar[r] & PX }\]
of $h$-equivalences, the corresponding preservation properties for 
ex-fibrations, and the properties of $q$-equivalences given by the 
statement that they are well-grounded; see \myref{hproper} and 
\myref{exwellgr}. It is also relevant that in each case passage to 
fibers gives the nonparametrized analogue of the construction under consideration. Since this result plays a vital role in our work, we 
give more complete details of (ii) and (iii); (i) works the same way.

The cofiber $C_Bf$ is the pushout of the diagram
\[ \xymatrix{C_B X & X \ar[l]\ar[r]^-f & Y.}\]
If $X$ is well-sectioned, then the left arrow is an $h$-cofibration
and $WX$ and $PX$ are well-sectioned. Replacing $f$
by $Wf$ and $Pf$ we obtain three such cofiber 
diagrams. Together with our original zig-zag this gives a 
$3\times 3$-diagram. Applying the gluing lemma, \myref{hproper}(iii), 
we obtain a zig-zag of $q$-equivalences
\[\xymatrix{C_B f & C_B Wf \ar[l]\ar[r] & C_B Pf.}\]
Similarly, on fibers we obtain zig-zags of $q$-equivalences
\[\xymatrix{Cf_b & C(Wf)_b \ar[l]\ar[r] & CW(Lf)_b.}\]
There results a zig-zag of $q$-equivalences of pairs 
\[\xymatrix{(C_B f, Cf_b) & (C_B Wf, CWf_b)\ar[l]\ar[r] & 
(C_B Pf, CW(Lf)_b).}\]
Since $C_BPf$ is ex-fibrant and in particular quasifibrant, $C_Bf$ is quasifibrant.

Similarly, by \myref{hproper}(v), we have natural zig-zags of $q$-equivalences
\[\xymatrix{X\sma_B K & WX\sma_B K \ar[l]\ar[r] & PX\sma_B K}\]
and 
\[\xymatrix{X_b\sma K & WX_b\sma K \ar[l]\ar[r] & W(LX)_b\sma K.}\]
We therefore have a zig-zag of $q$-equivalences of pairs
\[\xymatrix@=.6cm{(X\sma_B K, X_b\sma K) & (WX\sma_B K, WX_b\sma K)\ar[l]\ar[r] 
& (PX \sma_B K, W(LX)_b\sma K).}\]
Since $PX\sma_B K$ is ex-fibrant and in particular quasifibrant, $X\sma_B K$ is quasifibrant.
\end{proof}

\chapter{The equivalence between $\Ho G\sK_B$ and $hG\sW_B$}

\section*{Introduction}

We developed the point-set level properties of the category $G\sK_B$
of ex-$G$-spaces over $B$ in Chapters 1 and 2, and we developed those homotopical properties that are accessible to model theoretic techniques 
in Chapter 4 -- 7. In this chapter, we use ex-fibrations to prove that 
certain structure on the point-set level that seems inaccessible from the 
point of view of model category theory nevertheless descends to homotopy categories. In particular, we prove that $\Ho G\sK_B$ is closed symmetric monoidal and that the right derived functor $f^*$ of the Quillen adjunction $(f_!,f^*)$ in the $qf$-model structure is closed symmetric monoidal and 
has a right adjoint $f_*$.

In \S9.1 we use the ex-fibrant approximation functor to prove that our model theoretic homotopy category of ex-$G$-spaces over $B$ is equivalent to the classical homotopy category of ex-$G$-fibrations over $B$. In \S9.2, we discuss how to pass to derived functors on either side of that equivalence in certain general cases. Replacing the model-theoretic method of constructing derived functors by a more classical method given in terms of ex-fibrant approximation, we construct the functors $f_*$ and $F_B$ on homotopy categories in \S9.3. By a combination of methods, we prove that $\Ho G\sK_B$ is a symmetric monoidal category and that the base change functor $f^*$ descends to a closed symmetric monoidal functor on homotopy categories in \S9.4. We also obtain such descent 
to homotopy categories results for change of group adjunctions and for passage to fibers in that section.  These results are central to the theory, and there 
seem to be no shortcuts to their proofs.

Everything is understood to be equivariant in this chapter, and we abbreviate ex-$G$-fibration and ex-$G$-space to ex-fibration and ex-space throughout. 
We shall retreat just a bit from all--embracing generality. We assume that $G$ is a Lie group and that all given base $G$-spaces $B$ are proper and 
are of the homotopy types of $G$-CW complexes. The reader may prefer to 
assume that $G$ is compact, but there is no gain in simplicity. In view of 
the properties of the base change adjunction $(f_!,f^*)$ given in 
\myref{Qad10}, there would be no real loss of generality if we restricted further to base spaces that are actual $G$-CW complexes, but that would be inconveniently restrictive.

\section{The equivalence of $\Ho G\sK_B$ and $hG\sW_B$}

Recall that $X\sma_B I_+$  is a cylinder object in the sense of the $qf$-model structure. When we restrict to $qf$-fibrant and $qf$-cofibrant objects, homotopies in the $qf$-model sense are the same as $fp$-homotopies, by \myref{comphty}. The morphism set $[X,Y]_{G,B}$ in $\Ho G\sK_B$ is naturally isomorphic to $[RQX,RQY]_{G,B}$, and this is the set of $fp$-homotopy classes of maps $RQX\rtarr RQY$. Here $R$ and $Q$ denote the functorial $qf$-fibrant and
$qf$-cofibrant approximation functors obtained from the small object argument. The total space of $RQX$ has the homotopy type of a $G$-CW complex since $B$
does. This leads us to introduce the following categories. 

\begin{defn}
Define $G\sV_B$ to be the full subcategory of $G\sK_B$ whose objects are  
well-grounded and $qf$-fibrant with total spaces of the homotopy types
of $G$-CW complexes. Define $G\sW_B$ to be the full subcategory of $G\sV_B$ whose objects are the ex-fibrations over $B$. Let $hG\sW_B$ denote the category obtained from $G\sW_B$ by passage to $fp$-homotopy classes of maps.
\end{defn}

Note that the definition of $G\sW_B$ makes no reference to model category theory. Recall that well-grounded means well-sectioned and compactly generated. When $B=*$, $G\sW_*$ is just the category of well-based compactly generated $G$-spaces of the homotopy types of $G$-CW complexes, and it is standard that its classical homotopy category is equivalent to the homotopy category of based $G$-spaces with respect to the $q$-model structure.  We shall prove a parametrized generalization.

We think of $G\sV_B$ as a convenient half way house between $G\sK_B$ and $G\sW_B$. It turns out to be close enough to the category of $qf$-cofibrant and $qf$-fibrant objects in $G\sK_B$ to serve as such for our purposes, while already having some of the properties of $G\sW_B$. The following crucial theorem fails for the $q$-model structure. It is essential for this result that we  allow the objects of $\sV_B$ to be well-sectioned rather than requiring them
to be $qf$-cofibrant. This will force an assymmetry when we deal with left and right derived functors in \myref{cderiv} below.

\begin{thm}\mylabel{CWfix}
The $qf$-cofibrant and $qf$-fibrant approximation functor $RQ$ and the ex-fibrant approximation functor $P$, together with the forgetful functors
$I$ and $J$, induce the following equivalences of homotopy categories.
\[\xymatrix{\Ho G\sK_B \ar@<.5ex>[r]^-{RQ} 
& \Ho G\sV_B \ar@<.5ex>[r]^-P\ar@<.5ex>[l]^-I 
& hG\sW_B \ar@<.5ex>[l]^-J }\]
\end{thm}
\begin{proof}
For $X$ in $G\sK_B$, we have a natural zig-zag of $q$-equivalences in $G\sK_B$ 
\[\xymatrix{X & QX\ar[l]\ar[r] & RQX.}\]
Therefore $X$ and $IRQX$ are naturally $q$-equivalent in $G\sK_B$. If $X$ is in $G\sV_B$, then it is $qf$-fibrant and therefore so is $QX$. Then the above zig-zag is in $G\sV_B$ and thus $X$ and $RQIX$ are naturally $q$-equivalent in $G\sV_B$.

Since $q$-equivalences in $G\sV_B$ are $h$-equivalences, and $P$ takes $h$-equivalences to $fp$-equivalences, it is clear that $P$ induces a functor on homotopy categories. Conversely, since $fp$-equivalences are in particular 
$q$-equivalences, the forgetful functor $J$ induces a functor in the other direction. For $X$ in $G\sV_B$ we have the natural zig-zag of $h$-equivalences 
\[\xymatrix{X & WX\ar[r]^{W\iota}\ar[l]_\rho & PX}\]
of \myref{exfibapp}. However $WX$ may not be in $G\sV_B$ since it may not be
$qf$-fibrant. Applying $qf$-fibrant approximation, we get a natural zig-zag of $q$-equivalences in $G\sV_B$ connecting $X$ and $PX$. It follows that $X$ and $JPX$ are naturally $q$-equivalent in $G\sV_B$. Starting with $X$ in $G\sW_B$, the above display is a zig-zag of $fp$-equivalences in $G\sW_B$, by \myref{fpequiv2}. It follows that $X$ and $JPX$ are naturally 
$fp$-equivalent in $G\sW_B$.
\end{proof}

\section{Derived functors on homotopy categories}

Model category theory tells us how Quillen functors $V\colon G\sK_A\rtarr G\sK_B$ induce derived functors on the homotopy categories on the left hand side of the equivalence displayed in \myref{CWfix}. We now seek an equivalent way of passing to derived functors on the right hand side. We begin with an informal discussion. We focus on functors of one variable, but functors of several
variables work the same way.

Following the custom in algebraic topology, we have been abusing notation by using the same notation for point-set level functors and for derived homotopy category level functors. We will continue to do so. However, the more accurate notations of algebraic geometry, $LV$ and $RV$ for left and right derived functors, might clarify the discussion. As we have already seen in \myref{noway}, passage to derived functors is not functorial in general, so that a relation between composites of functors that holds on the point-set level need not imply a corresponding relation on passage to homotopy categories.

Recall that, model theoretically, if $V$ is a Quillen right adjoint, then the right derived functor of $V$ is obtained by first applying fibrant approximation $R$ and then applying $V$ on homotopy categories, which makes sense since $V$ preserves weak equivalences between fibrant objects. The left derived functor of a Quillen left adjoint $V$ is defined dually, via $VQ$. Problems arise when one tries to compose left and right derived functors, which is what we must do to prove some of our compatibility relations.

The equivalence of categories proven in \myref{CWfix} gives us a way of putting the relevant left and right adjoints on the same footing, giving a ``straight'' passage to derived functors that is neither ``left'' nor ``right''.  We need
mild good behavior for this to work.

\begin{defn}
A functor $V\colon G\sK_A\rtarr G\sK_B$ is \emph{good}\index{functor!good}\index{good functor} if it is continuous,
takes well-grounded ex-spaces to well-grounded ex-spaces, and takes ex-spaces
whose total spaces are of the homotopy types of $G$-CW complexes to ex-spaces with that property. Since $V$ is continuous, it preserves $fp$-homotopies.
\end{defn}

\begin{prop}\mylabel{cderiv}
Let $V\colon G\sK_A\rtarr G\sK_B$ be a good functor that is a left or
right Quillen adjoint.  If $V$ is a Quillen left adjoint, assume further 
that it preserves $q$-equivalences between well-grounded ex-spaces. Then,
under the equivalence of categories in \myref{CWfix}, the derived functor 
$\Ho G\sK_A\rtarr\Ho G\sK_B$ induced by $VQ$ or $VR$ 
is equivalent to the functor
$PVJ\colon hG\sW_A\rtarr hG\sW_B$ obtained
by passage to homotopy classes of maps. \end{prop}
\begin{proof}
If $V$ is a Quillen right adjoint, then it preserves $q$-equivalences between $qf$-fibrant objects. If $V$ is a Quillen left adjoint, then we are assuming that it preserves $q$-equivalences between well-grounded objects. Since $G\sV_A$ consists of well-sectioned $qf$-fibrant objects, it follows in both cases that  $V\colon  G\sV_A \rtarr G\sV_B$ passes straight to homotopy categories to give  $V\colon  \text{Ho}G\sV_A \rtarr \text{Ho}G\sV_B$. Since $V$ preserves $G$-CW homotopy types on total spaces,  $V$ takes $q$-equivalences to $h$-equivalences. Therefore $PV$ takes $q$-equivalences
to $fp$-equivalences and induces a functor $\Ho G\sV_A\rtarr hG\sW_B$. To show that $PVJ$ and either $VQ$ or $VR$ agree under the equivalence of categories, 
it suffices to verify that the following diagram commutes.
\[\xymatrix{
\Ho G\sK_A\ar[d]_{RQ} \ar[rr]^-{VQ \ \ \text{or} \ \ VR} 
& &  \Ho G\sK_B \ar[d]^{PRQ}\\
\Ho G\sV_A \ar[rr]_{PV} & &  hG\sW_B}\]
We have a natural acyclic $qf$-fibration $QX\rtarr X$ and a natural 
acyclic $qf$-cofibration $X\rtarr RX$. If $V$ is a Quillen left adjoint, 
then we have a zig-zag of natural $q$-equivalences
\[\xymatrix{RQVQ \ar[r] & RVQ & VQ \ar[l]\ar[r] & VRQ}\]
because $V$ preserves acyclic $qf$-cofibrations.
If $V$ is a Quillen right adjoint, then we have a zig-zag of natural $q$-equivalences
\[\xymatrix{RQVR  & RQVRQ \ar[r]\ar[l] & RVRQ & VRQ\ar[l]}\]
because $V$ preserves $q$-equivalences between $qf$-fibrant objects. In both cases, all objects have total spaces of the homotopy types of $G$-CW complexes, so in fact we have zig-zags of $h$-equivalences.  Therefore, applying $P$ gives
us zig-zags of $fp$-equivalences in $G\sW_B$, by \myref{fpequiv2}.  
\end{proof}

\begin{rem}
When $V$ preserves ex-fibrations, $PV$ is naturally $fp$-equivalent to $V$ on ex-fibrations, by \myref{fpequiv2}. The derived functor of $V$ can then be obtained directly by applying $V$ and passing to equivalence classes of maps under $fp$-homotopy. 
\end{rem}

\section{The functors $f_*$ and $F_B$ on homotopy categories}

We use the equivalence between $\Ho G\sK_B$ and $hG\sW_B$ to prove that, for any map $f\colon A\rtarr B$ between spaces of the homotopy types of $G$-CW complexes, the $(f^*,f_*)$ adjunction descends to homotopy categories. 
We begin by verifying that $f^*$ satisfies the hypotheses of \ref{cderiv}. 

\begin{prop}\mylabel{fstarderiv}
Let $f\colon A\rtarr B$ be a map of base spaces. Then the base change functor $f^*$ restricts to a functor $f^*\colon G\sW_B\rtarr G\sW_A$.
\end{prop}
\begin{proof}
Consider $Y$ in $G\sW_B$. Since the total space of $Y$ is of the homotopy type of a $G$-CW complex, the fibers $Y_b$ are of the homotopy types of $G_b$-CW complexes by \myref{ss}. The fiber $(f^*Y)_a$ is a copy of $Y_{f(a)}$, and $G_a$ acts through the evident inclusion $G_a\subset G_{f(a)}$. Therefore $(f^*Y)_a$ is of the homotopy type of a $G_a$-CW complex. The total space of $f^*Y$ is therefore of the homotopy type of a $G$-CW complex, again by \myref{ss}.  Moreover, $f^*Y$ is an ex-fibration by \myref{fexpres}.  Thus $f^*$ restricts to a functor $f^*\colon G\sW_B\rtarr G\sW_A$.
\end{proof}

\begin{thm}\mylabel{descendf0}
For any map $f\colon A\rtarr B$ of base spaces, the right derived functor $f^*\colon \Ho G\sK_B\rtarr \Ho G\sK_A$ has a right adjoint $f_*$, so that $$ [f^*Y,X]_{G,A} \iso [Y,f_*X]_{G,B}$$ for $X$ in $G\sK_A$ and $Y$ in $G\sK_B$.
\end{thm}
\begin{proof}
In view of the equivalence of categories in \myref{CWfix} and the fact that $f^*$ descends directly to a functor $f^*\colon hG\sW_B\rtarr hG\sW_B$
on homotopy categories, by Propositions \ref{cderiv} and \ref{fstarderiv}, it suffices to construct a right adjoint $f_*\colon hG\sW_A \rtarr hG\sW_B$. We do that using the Brown representability theorem. By \myref{brown0}, $\Ho G\sK_B$ satisfies the formal hypotheses for Brown representability, and therefore so does $hG\sW_B$. In fact $G\sW_B$ has all of the relevant wedges and homotopy colimits since these constructions preserve ex-fibrations by \myref{pres1} and \myref{HursmaK} and since they clearly preserve $G$-CW homotopy types on the total space level and stay within $G\sU_B$.  The objects in the detecting set $\sD_B$ of \myref{detectset} are not in $G\sW_B$, but we can apply the ex-fibrant approximation functor $P$ to them to obtain a detecting set of objects in $hG\sW_B$. Therefore a contravariant set-valued functor on $hG\sW_B$ is representable if and only if it satisfies the wedge and Mayer-Vietoris axoms.

For a fixed ex-fibrant space $X$ over $A$, consider the functor $\pi(f^*Y,X)_{G,A}$ on $Y$ in $G\sW_B$, where $\pi$ denotes $fp$-homotopy classes of maps. Since the functor $\pi(W,X)_{G,A}$ on $W$ in $G\sW_A$ is represented and is computed using homotopy classes of maps, it clearly satisfies the wedge and Mayer-Vietoris axioms. It therefore suffices to show that the functor $f^*$ preserves wedges and homotopy pushouts, since that will imply 
that the functor $\pi(f^*Y,X)_{G,A}$ of $Y$ also satisfies the wedge and Mayer-Vietoris axioms. We can then conclude that there is an object $f_*X\in G\sW_B$ that represents this functor. It follows formally that $f_*$ is a functor of $X$ and that the required adjunction holds.

Because $f^*\colon G\sK_B\rtarr G\sK_A$ is a left adjoint, it preserves colimits, and this implies that $f^*\colon G\sW_B\rtarr G\sW_A$  preserves the relevant homotopy colimits. Moreover, $f^*$ preserves $fp$-homotopies and so induces a functor on homotopy categories that still preserves these homotopy colimits.
\end{proof}

We agree to write $\simeq$ for natural equivalences on homotopy categories.

\begin{rem} 
For composable maps $f$ and $g$, $g_*\com f_* \simeq (g\com f)_*$ on homotopy categories since $f^*\com g^*\simeq (g\com f)^*$ on homotopy categories. The latter equivalence is clear since $f^*$ and $g^*$ are derived from Quillen right adjoints. More sophisticated commutation laws are proven in the next section.
\end{rem}

Applying \myref{descendf0} to diagonal maps and composing with the homotopy category level adjunction between the external smash product and function ex-space functors, we obtain the following basic result; compare \myref{internalize}.

\begin{thm}\mylabel{smashing0}
Define $\sma_B$ and $F_B$ on $\Ho G\sK_B$ as the composite (derived) functors
$$X\sma_B Y = \DE^*(X\barwedge Y) \qquad\text{and}\qquad F_B(X,Y) = \bar{F}(X,\DE_*Y).$$
Then 
$$ [X\sma_B Y, Z]_{G,B}\iso [X,F_B(Y,Z)]_{G,B}$$
for $X$, $Y$, and $Z$ in $\Ho G\sK_B$.
\end{thm}
\begin{proof}
The displayed adjunction is the composite of adjunctions for the (derived) external smash and function ex-space functors and for the (derived) adjoint pair $(\DE^*,\DE_*)$.
\end{proof}

\begin{rem} 
The referee points out that the ex-space analogue of \cite[7.2]{BB2} shows that we can work directly with the point-set topology to show that the $(\sma_B,F_B)$ adjunction on the original category $G\sK_B$ is continuous and so descends to (classical) $fp$-homotopy categories to give the adjunction
$$ hG\sK_B(X\sma_B Y,Z)\iso hG\sK_B(X,F_B(Y,Z)).$$
Presumably similar point-set topological arguments work to show that, for a map $f\colon A\rtarr B$, we have an adjunction
$$hG\sK_A(f^*X,Y)\iso hG\sK_B(X,f_*Y).$$
These adjunctions do {\em not}\, imply our Theorems \ref{descendf0} and \ref{smashing0}. By definition, our category $hG\sW_B$ is a full subcategory of $hG\sK_B$, but it is not an {\em equivalent}\, full subcategory. The objects of $G\sW_B$ are very restricted, and general function ex-spaces $F_B(Y,Z)$ are not $fp$-homotopy equivalent to such objects. The force of our theorems is that, after restricting to our subcategories $hG\sW_B$, we still have right adjoints {\em in these categories}.  It is this fact that we need to obtain right adjoints in our preferred homotopy categories $\Ho G\sK_B$.
\end{rem}

\section{Compatibility relations for smash products and base change}

We first prove that $\Ho G\sK_B$ satisfies the associativity, commutativity and unity conditions required of a symmetric monoidal category. We then show that all of the isomorphisms of functors in \myref{Wirth0} and some of the isomorphisms of functors in \myref{Mackey0} still hold after passage to homotopy categories. Finally, we relate change of groups and passage to fibers to the
symmetric monoidal structure on homotopy categories.  In some of our arguments, it is natural to work in $\Ho G\sK_{B}$. In others, it is natural to work in the equivalent category $hG\sW_{B}$.

\begin{prop}\mylabel{fgext'}
For maps $f\colon A\rtarr B$ and $g\colon A'\rtarr B'$ of base spaces and for ex-spaces $X$ over $B$ and $Y$ over $B'$, 
\begin{equation}\label{fgext1'}
 (f^*Y\barwedge g^*Z)\simeq (f\times g)^*(Y\barwedge Z)
\end{equation}
in $\Ho G\sK_A$. For ex-spaces $W$ over $A$ and $X$ over $A'$,
\begin{equation}\label{fgext2'}
(f_!W\barwedge g_!X)\simeq (f\times g)_!(W\barwedge X)
\end{equation}
in $\Ho G\sK_B$.
\end{prop}
\begin{proof} 
For (\ref{fgext1'}), we work with ex-fibrations, starting in $hG\sW_{B\times B'}$. By Propositions \ref{fexpres} and \ref{Hursma}, the functors we are dealing with preserve ex-fibrations and therefore descend straight to homotopy categories. The conclusion is thus immediate from its point-set level analogue. For (\ref{fgext2'}), we work with model theoretic homotopy categories, starting in $\Ho G\sK_{A\times A'}$. Since $(f\times g)_!\simeq (f\times\text{id})_!\com (\text{id}\times g)_!$, we can proceed in two steps and so assume that $g=\text{id}$. By \myref{smaAB} and \myref{Qad10}, we are then composing Quillen left adjoints. Starting with $qf$-cofibrant objects, we do not need to apply $qf$-cofibrant approximation, and the conclusion follows directly from its point-set level analogue.
\end{proof}

We use this to complete the proof that $\Ho G\sK_B$ is symmetric monoidal.

\begin{thm}\mylabel{clsymmon}
The category $\Ho G\sK_B$ is closed symmetric monoidal under the functors $\sma_B$ and $F_B$. 
\end{thm}
\begin{proof} 
In view of \myref{smashing0}, we need only prove the associativity, commutativity, and unity of $\sma_B$ up to coherent natural isomorphism. The external smash product has evident associativity, commutativity, and unity isomorphisms, and these descend directly to homotopy categories since the external smash product of $qf$-cofibrant ex-spaces over $A$ and $B$ is $qf$-cofibrant over $A\times B$. To see that these isomorphisms are inherited after internalization along $\DE^*$, we use (\ref{fgext1'}).  For the associativity of $\sma_B$, we have
\begin{multline*}
\DE^*(\DE^*(X\barwedge Y)\barwedge Z)
\simeq \DE^*(\DE\times\text{id})^*((X\barwedge Y)\barwedge Z)
\simeq ((\DE\times\text{id})\DE)^*((X\barwedge Y)\barwedge Z)\\
\quad \quad \simeq ((\text{id}\times \DE)\DE)^*(X\barwedge(Y\barwedge Z)) 
\simeq \DE^*(\text{id}\times \DE)^*(X\barwedge(Y\barwedge Z))
\simeq \DE^*(X\barwedge \DE^*(Y\barwedge Z)).
\end{multline*}

\noindent
The commutativity of $\sma_B$ is similar but simpler. For the unit, we observe that $S^0_B\simeq {r^*S^0}$,
$r\colon B\rtarr *$. Therefore, since $(\text{id}\times r)\DE = \text{id}$, 
\[X\sma_B S^0_B \simeq \DE^*(X\barwedge r^*S^0) 
\simeq \DE^*(\text{id}\times r)^*(X\barwedge S^0)
\simeq ((\text{id}\times r)\DE)^* (X) = X.\qedhere\]
\end{proof}

We turn next to the derived versions of the base change compatibilities of Propositions \ref{Wirth0} and \ref{Mackey0}.  Observe that the functor $f_!$ is
good since the section of a well-sectioned ex-space is an $h$-cofibration and
since $G$-CW homotopy types are preserved under pushouts, one leg of which is 
an $h$-cofibration. Moreover, $f_!$ preserves $q$-equivalences between well-sectioned ex-spaces by \myref{Qad10}. Therefore \myref{cderiv} applies to $f_!$.

\begin{thm}\mylabel{fclsymmon}
For a $G$-map $f\colon A\rtarr B$, $f^*\colon \Ho G\sK_B\rtarr \Ho G\sK_A$ is a closed symmetric monoidal functor. 
\end{thm}
\begin{proof}
Since $f^*S^0_B\iso S^0_A$ in $G\sK_A$ and $S^0_B$ is $qf$-fibrant, $f^*S^0_B\simeq S^0_A$ in $\Ho G\sK_A$. We must prove that the isomorphisms (\ref{oneo}) through (\ref{five0}) descend to equivalences on homotopy categories.  Categorical arguments in \cite[\S\S2, 3]{FHM} show that it suffices to show that the two isomorphisms (\ref{oneo}) and (\ref{four0}) descend to equivalences on homotopy categories. These two isomorphisms do not involve the right adjoints $f_*$ or $\DE_*$ and are therefore more tractable than the other three. First consider 
(\ref{oneo}):
\[f^*(Y\sma_B Z)\iso f^*Y\sma_A f^*Z.\]
If $Y$ and $Z$ are in $G\sW_B$, then the two sides of this isomorphism are both in $G\sW_A$, by \myref{fexpres} and \myref{Hursma}. Therefore the point-set level isomorphism descends directly to the desired homotopy category level equivalence. Next, consider (\ref{four0}):
\[f_{!}(f^*Y\sma_A X)\iso Y\sma_B f_{!}X.\]
Assume that $X$ is in $G\sW_A$ and $Y$ is in $G\sW_B$. The functor $f_!$ does not preserve ex-fibrations so, to pass to derived categories, we must replace 
it by $Pf_!$ on both sides.  By \myref{savior}, the functor $Y\sma_B(-)$ preserves $h$-equivalences between well-sectioned ex-spaces. Since $P$ sends $h$-equivalences to $fp$-equivalences, we therefore have $fp$-equivalences, natural up to $fp$-homotopy,
\[\xymatrix{
Pf_{!}(f^*Y\sma_A X)\iso P(Y\sma_B f_{!}X)\ \ar[r]^-{P(\text{id}\sma_B\ph)} & 
\ P(Y\sma_BPf_!X) & Y\sma_BPf_!X, \ar[l]_-{\ph}}\]
where $\ph = (\rh,W\io)$ is the zigzag of $h$-equivalences of \myref{exfibapp}. This implies the desired equivalence in the homotopy category.
\end{proof}

The reader is invited to try to prove directly that the projection formula holds in the homotopy category. Even the simple case of $f\colon *\rtarr B$, the inclusion of a point, should demonstrate the usefulness of \myref{cderiv}.

\begin{thm}\mylabel{pullbackfix} 
Suppose given a pullback diagram of $G$-spaces 
$$\xymatrix{
C \ar[r]^-{g} \ar[d]_{i} & D \ar[d]^{j} \\
A \ar[r]_{f} & B}$$
in which $f$ (or $j$) is a $q$-fibration. Then there are natural equivalences of functors on homotopy categories
\begin{equation}\label{basesmore0}
j^*f_{!} \simeq g_{!}i^*, \quad f^*j_* \simeq i_*g^*, \quad f^*j_{!}\simeq i_!g^*, \quad j^*f_*\simeq g_*i^*.
\end{equation} 
\end{thm}
\begin{proof}
As in \myref{Mackey0} the second and fourth equivalences are conjugate to the first and third. However, since the situation is no longer symmetric, we must prove both the first and third equivalences, assuming $f$ is a $q$-fibration.

First consider the desired equivalence $f^*j_{!}\simeq i_!g^*$. We work with ex-fibrations, starting with $X\in hG\sW_D$.  We must replace $j_!$ and $i_!$ by $Pj_!$ and $Pi_!$ before passing to homotopy categories. By \myref{Qad10}, $f^*$ preserves $q$-equivalences since $f$ is a $q$-fibration. Moreover, our $q$-equivalences are $h$-equivalences since we are dealing with total spaces of the homotopy types of $G$-CW complexes.  By the diagram in \myref{munu}(ii), we see that $\mu\colon Pf^*\rtarr f^*P$ is a natural $h$-equivalence here. This would be false for arbitrary maps $f$, as observed in \myref{PffP}. Since $\mu$ is an $h$-equivalence between ex-fibrations, it is an $fp$-equivalence. Therefore
$$f^*Pj_{!}X\htp Pf^*j_!X\iso  Pi_!g^*X.$$

Now consider the desired equivalence $j^*f_{!}X \simeq g_{!}i^*X$ in $\Ho G\sK_D$.  Our assumption that $f$ is a $q$-fibration gives us no direct help with this.  However, we may factor $j$ as the composite of a homotopy equivalence and an $h$-fibration. Expanding our pullback diagram as a composite of pullbacks, we see that it suffices to prove our commutation relation when $j$ is an $h$-fibration and when $j$ is a homotopy equivalence. The first case is immediate by symmetry from the first part. Thus assume that $j$ is a homotopy equivalence. Then $i$ is also a homotopy equivalence. By \myref{Qad10}, $(i_!,i^*)$ and $(j_!,j^*)$ are adjoint equivalences of homotopy categories. Therefore
\[ j^*f_! \simeq j^*f_!i_!i^* \simeq j^*j_!g_!i^*\simeq g_!i^*. \qedhere\]
\end{proof}

Finally, we turn to a promised compatibility relationship between products
and change of groups. We observed in \myref{Lishriek} that the point-set level closed symmetric monoidal equivalence of \myref{ishriek} is given by a Quillen equivalence. The following addendem shows that the resulting equivalence on homotopy categories is again closed symmetric monoidal.

\begin{prop}\mylabel{imonoidaldescends}
Let $\iota\colon H\rtarr G$ be the inclusion of a subgroup and $A$ be an $H$-space. The Quillen equivalence $(\iota_!, \nu^*\iota^*)$ descends to a closed symmetric monoidal equivalence between $\text{Ho}H\sK_A$ and 
$\text{Ho}G\sK_{\io_!A}$.
\end{prop}

\begin{proof}
Let $\Delta\colon A\rtarr A\times A$ be the diagonal map. The isomorphisms
\[\iota^*\Delta^*(X\barwedge Y)\cong \Delta^*\iota^* (X\barwedge Y)
\cong \Delta^* (\iota^*X\barwedge \iota^*Y)\]
descend to equivalences on homotopy categories, the first since it is between Quillen right adjoints, the second since $\iota^*$ preserves all $q$-equivalences. It follows that $\nu^*\iota^*$ is a symmetric monoidal 
functor on homotopy categories. Since it is also an equivalence, it follows formally that it is closed symmetric monoidal.
\end{proof}

Combined with \myref{fclsymmon} applied to the inclusion 
$\tilde{b}\colon G/G_b \rtarr B$, this last observation 
gives us the following conclusion. 

\begin{thm}\mylabel{fiberfun}
The derived fiber functor $(-)_b\colon \Ho G\sK_B\rtarr \Ho G_b\sK_b$ is closed symmetric monoidal, and it has a left adjoint $(-)^b$ and a right adjoint ${^b}(-)$.
\end{thm} 

We emphasize that this innocent looking result packages highly non-trivial and important information. It gives in particular that, for ex-$G$-spaces $X$ 
and $Y$, the (derived) fiber $F_B(X,Y)_b$ of the (derived) function space $F_B(X,Y)$ is equivalent in $\Ho G_b\sK_b$ to the (derived) function space $F(X_b,Y_b)$ of the (derived) fibers $X_b$ and $Y_b$. On the point set level, that is what motivated the definition of the internal function ex-space.  
That it still holds on the level of homotopy categories is a reassuring consistency result.

\part{Parametrized equivariant stable homotopy theory}

\chapter*{Introduction}

We develop rigorous foundations for parametrized equivariant
stable homotopy theory. The idea is to start with a fixed base
$G$-space $B$ and to build a good category, 
here denoted $G\sS_{B}$, of $G$-spectra over $B$. We assume once and 
for all that our base spaces $B$ must be compactly generated and must 
have the homotopy types of $G$-CW complexes.  By ``good'' we mean 
that $G\sS_{B}$ is a closed symmetric monoidal topological model category 
whose associated homotopy category has properties analogous to those 
of the ordinary equivariant stable homotopy category.

Informally, the homotopy theory of $G\sS_{B}$ is specified by the 
homotopy theory seen on the fibers of $G$-spectra over $B$.
One compelling reason for taking the parametrized stable homotopy 
category seriously, even nonequivariantly, is to build a natural home 
in which one can do stable homotopy theory while still keeping track of fundamental groups and groupoids.  Stable homotopy theory has tended to 
ignore such intrinsically unstable data. This has the effect of losing 
contact with more geometric branches of mathematics in which the 
fundamental group cannot be ignored. 

For example, one basic motivation for the 
equivariant theory is that it gives a context in which to better 
understand equivariant orientations, Thom isomorphisms, and 
Poincar\'e duality. There is no problem for $G$-simply connected 
manifolds $M$ \cite[III\S6]{LMS}, but restriction to such $M$ is 
clearly inadequate for applications to transformation group theory. 
Despite a great deal of work on the subject by Costenoble and Waner, 
and some by May, \cite{CMW, CW1, CW2, CW3, MayR}, this circle of 
ideas is not yet fully understood.  Costenoble and Waner \cite{CWNew} 
use our work to study this problem for ordinary equivariant theories,
and for general theories this is work in progress by the second 
author.

There are many problems that 
make the development far less than an obvious generalization of the 
nonparametrized theory. Problems on the space level were dealt
with in Parts I and II, and we deal with the analogous spectrum level 
problems here.  We give some categorical preliminaries on enriched
equivariant categories in Chapter 10.  We define and develop the basic properties of our preferred category of parametrized $G$-spectra in 
Chapter 11, study its model structures in Chapter 12, and study 
adjunctions and compatibility relations in Chapter 13.  All of the 
problems that we faced on the space level are still there, but their 
solutions are considerably more difficult. In Chapter 14, we go on to 
study further such compatibilities that more fundamentally involve 
equivariance. 
 
The theory of highly structured spectra is highly cumulative.
We build on the theory of equivariant orthogonal
spectra of Mandell and May \cite{MM}.  In turn, that 
theory builds on the theory of
nonequivariant orthogonal spectra. A self-contained treatment of 
nonequivariant diagram spectra, including orthogonal spectra, 
is given by Mandell, May, Schwede, Shipley in \cite{MMSS}.  
The treatments of \cite{MM} and \cite{MMSS}, like this one, are topological 
as opposed to simplicial. That seems to be essential when dealing with 
infinite groups of equivariance. It also allows use of orthogonal 
spectra rather than symmetric spectra.  These are much more natural 
equivariantly and, even nonequivariantly, they have the major 
convenience that their weak equivalences are exactly the maps that 
induce isomorphisms of homotopy groups. 

The theory of equivariant parametrized spectra can be thought
of as the pushout over the theory of spectra of the theories of 
equivariant spectra and of nonequivariant parametrized spectra. 
However, there is no nonequivariant precursor of the present treatment
of parametrized spectra in the literature. There are preliminary forms 
of such a theory \cite{BG1, BG2, Clapp, CP, CJ}, but these either do not
go beyond suspension spectra or are based on obsolescent technology.
None of them go nearly far enough into the theory for the purposes we 
have in mind, although the early first approximation of 
Monica Clapp \cite{Clapp}, written up in more detail with Dieter Puppe 
\cite{CP}, deserves considerable credit. Clapp gave the strongest previous 
version of our fiberwise duality theorem, and her emphasis on ex-fibrations, 
together with some key technical results about them, have been very helpful. 
The reader primarily interested in classical homotopy theory should ignore all 
details of equivariance in reading Chapters 11--13. In fact, 
given \cite{MM}, the equivariance adds few serious difficulties to the passage from spectra to parametrized spectra, although it does add many interesting new features.

There are at least two possible alternative cumulative approaches. Rather 
than building on the theory of orthogonal $G$-spectra of \cite{MM, MMSS}, 
one can build on the theory of $G$-spectra of \cite{LMS}, the theory of 
$S$-modules of \cite{EKMM}, and the pushout of these, the theory of 
$S_G$-modules of \cite{MM}.  Po Hu \cite{Hu} began work on the first 
stage of a treatment along these lines, using parametrized $G$-spectra, 
but she did not address the foundational issues concerning smash products, function spectra, base change functors, and compatibility relations considered
here.  Moreover, following the first author's misleading unpublished notes \cite{May0}, she took the $q$-model structure on ex-$G$-spaces 
as her starting point, and the stable model structure cannot be made rigorous from there. It appears to us that resolving all of these issues in that framework is likely to be more difficult than in the framework that we have adopted. In particular, homotopical control of the parametrized spectrification functor and of cofiber sequences seems problematic.  

Alternatively, for finite groups $G$, one can build on the theory of symmetric 
spectra of Hovey, Smith, and Shipley \cite{HSS} and its equivariant 
generalization due to Mandell \cite{Mandell}.  Such an approach would avoid 
the point-set topological technicalities of the present approach and would
presumably lead to rather different looking problems with fibrations and 
cofibrations.  The problems with the stable homotopy category level adjunctions
that involve base change functors, smash products, and function spectra are
intrinsic and would remain. Our solutions to these problems do not seem to 
carry over to the simplicial context in an obvious way, and an alternative simplicial treatment could prove to be quite illuminating.
 
In view of the understanding of unstable equivariant homotopy 
theory for proper actions of non-compact Lie groups that was 
obtained in Part II, it might seem that there should be no real 
difficulty in obtaining a good stable theory along the same lines 
as the theory for compact Lie groups. However, in contrast with the 
rest of this book, equivariant stable homotopy theory for non-compact Lie 
groups is in preliminary and incomplete form, with still unresolved 
technical problems. We leave its study to future work, explaining
in \S11.6 where some of the problems lie. Except in that section, 
$G$ is asssumed to be a compact Lie group from Chapter 11 onwards.
A few other notes on terminology may be helpful.  We shall not use the term ``ex-spectrum over $B$'' since, stably, there is no meaningful unsectioned theory. Instead, we shall use the term ``spectrum over $B$''. This is especially convenient when considering base change.  We write out ``orthogonal $G$-spectrum over $B$'' until \S11.4. However, since we consider no other kinds 
of $G$-spectra and work equivariantly throughout, we later abbreviate this 
to ``spectrum over $B$'' when there is no danger of confusion.  That is,
we work equivariantly throughout, but we only draw attention to this fact
when it plays a significant mathematical role. 

\chapter{Enriched categories and $G$-categories}

\section*{Introduction}

To give context for the structure enjoyed by the categories of 
parametrized orthogonal $G$-spectra that we shall define, we first 
describe the kind of equivariant para\-metrized enrichments that we shall encounter.  In fact, our categories have several layers of enrichment, and 
it is helpful to have a consistent language, somewhat non-standard from a categorical point of view, to keep track of them. 
In \S\S10.1 and 10.2, we give some preliminaries on enriched categories,
working non\-equi\-variant\-ly in \S10.1 and adding considerations of
equivariance in \S10.2. We discuss the role of the several enrichments 
in sight in our $G$-topological model $G$-categories in \S10.3. In
this chapter, $G$ can be any topological group.

\section{Parametrized enriched categories}

As discussed in \S1.2, all of our categories $\sC$ are topological, 
meaning that they are enriched over the category $\sK_*$ of based spaces 
(= $k$-spaces). In contrast with general enriched category theory and our further enrichments, the topological enrichment is given just by a topology on the underlying set of morphisms, and we denote the space of morphisms $X\rtarr Y$ by $\sC(X,Y)$.  We say that a topological category $\sC$ is {\em topologically bicomplete}\, if it is bicomplete and bitensored over
$\sK_*$.  In fact, we shall have enrichments and bitensorings over the 
category $\sK_B$ of ex-spaces over $B$ that imply the topological 
enrichment and bitensoring by restriction to ex-spaces $B\times T$ for $T\in\sK_*$. 

Recall from \S1.3 that $\sK_B$ is topologically bicomplete, with tensors 
and cotensors denoted by $K\sma_B T$ and $F_B(T,K)$ for $T\in\sK_*$ and $K\in \sK_B$.  (Since we shall use letters like $X$, $Y$, and $Z$ for spectra, we 
have changed the letters that we use generically for spaces and ex-spaces 
from those that we used earlier). It is also closed symmetric monoidal under 
its fiberwise smash 
product and function space functors, which are also denoted by $\sma_B$ and $F_B$; its unit object is $S^0_B = B\times S^0$. It is therefore enriched and bitensored over itself. The two enrichments are related by natural based homeomorphisms
\begin{equation}\label{eqn:enrich}
\sK_B(K,L) \iso \sK_B(S^0_B, F_B(K,L)).
\end{equation}
This is the case $T = S^0$ of the more general based homeomorphism
\begin{equation}\label{eqn:autoenrich}
\sK_*(T,\sK_B(K,L)) \iso \sK_B(S^0_B\sma_B T, F_B(K,L))
\end{equation}
for $T\in \sK_*$ and $K$, $L\in \sK_B$. The Yoneda lemma, (\ref{eqn:enrich}),
and the bitensoring adjunctions imply that the two bitensorings are related by 
the equivalent natural isomorphisms of ex-spaces
\begin{equation}\label{eqn:compatible}
K\wedge_BT \cong K\wedge_B  (S^0_B\wedge_B T)
\quad\text{and}\quad  F_B(T,K)\cong F_B(S^0_B\wedge_B T, K).
\end{equation}
These in turn imply the equivalent generalizations
\begin{equation}\label{eqn:trans}
K\sma_B (L\sma_B T)\iso (K\sma_B L)\sma_B T 
\quad\text{and}\quad F_B(T, F_B(K,L))\iso F_B(K\sma_B T,L).
\end{equation}
Formally, rather than defining the enrichments and bitensorings over $\sK_*$ independently,
we can take (\ref{eqn:autoenrich}) and (\ref{eqn:compatible}) as definitions of these 
structures in terms of the enrichment and bitensoring over $\sK_B$. Then (\ref{eqn:trans}) 
and the bitensoring adjunction homeomorphisms
\begin{equation}
 \sK_B(K\sma_B T, L)\iso \sK_*(T,\sK_B(K,L))\iso \sK_B(K,F_B(T,L))
\end{equation}
follow directly.  

\begin{rem}\mylabel{confusion}
We shall be making much use of the functor 
$S^0_B\sma_B(-)$, and we henceforward abbreviate 
notation by setting
$$  T_B = B\times T = S^0_B\sma_B T$$
for a based space $T$, and similarly for maps. Observe that
$K\sma_B T$ and $K\sma_B T_B$ are two names for the same
ex-space over $B$. When working on a formal conceptual level, 
it is often best to think in terms of tensors over $\sK_*$ and 
use the first name. However, on a pragmatic level, to avoid
confusion, it is best to view based spaces as embedded 
in ex-spaces via $S^0_B\sma_B(-)$ and to use the second
notation, working only with tensors over $\sK_B$.
\end{rem}

We generalize and formalize several aspects of the discussion above.

\begin{defn}
A topological category $\sC$ is \emph{topological over $B$}\index{category!topological over B@topological over $B$} if it is enriched and bitensored over $\sK_B$. It is \emph{topologically bicomplete over $B$}\index{category!topologically bicomplete over B@topologically bicomplete over $B$} if it is 
also bicomplete. We write $P_B(X,Y)$\noteindex{PBXY@$P_B(X,Y)$} for
the hom ex-space over $B$, and we write $X\sma_B K$\noteindex{XBK@$X\sma_B K$} and $F_B(K,X)$\noteindex{FBKX@$F_B(K,X)$} for the
tensor and cotensor in $\sC$, where $X$, $Y\in\sC$ and $K\in\sK_B$.  
Explicitly, we 
require bitensoring adjunction homeomorphisms of based spaces
\begin{equation}\label{eqn:adj}
\sC(X\wedge_B K, Y)\cong \sK_B(K, P_B(X,Y)) \cong \sC(X, F_B(K,Y)).
\end{equation}
By Yoneda lemma arguments, these imply unit and transitivity isomorphisms in $\sC$
\begin{equation}\label{eqn:transag}
X\iso X\sma_B S^0_B\quad\text{and}\quad X\sma_B (K\sma_B L) \iso (X\sma_B K)\sma_B L.
\end{equation}
and also bitensoring adjunction isomorphisms of ex-spaces
\begin{equation}\label{eqn:intadj}
P_B(X\wedge_B K, Y)\cong F_B(K, P_B(X,Y)) \cong P_B(X,F_B(K,Y)).
\end{equation}
Conversely, there is a natural homeomorphism
\begin{equation}\label{eqn:under}
\sC(X,Y)\cong \sK_{B}(S^0_B, P_B(X,Y)),
\end{equation}
and the isomorphisms (\ref{eqn:adj}) follow from (\ref{eqn:intadj}) by 
applying $\sK_B(S^0_B,-)$.
\end{defn}

If we do not require $\sC$ to be topological to begin with, we can take
(\ref{eqn:under}) as the definition of the space $\sC(X,Y)$ and so recover 
the topological enrichment. With the
notation of \myref{confusion}, we obtain tensors and cotensors
with based spaces $T$ by setting
\begin{equation}\label{eqn:biten3}
X\sma_B T = X\sma_B T_{B} \quad\text{and}\quad  F_B(T,X) = F_B(T_{B}, X).
\end{equation}
The adjunction homeomorphisms
\begin{equation} 
\sC(X\wedge_B T, Y)\cong \sK_*(T, \sC(X,Y)) \cong \sC(X, F_B(T,Y))
\end{equation}
are obtained by replacing $K$ by $T_{B}$ in (\ref{eqn:adj}) and using
(\ref{eqn:autoenrich}) and (\ref{eqn:under}).

In the cases of interest, $\sC$ is closed symmetric monoidal, and then the 
hom ex-spaces $P_B(X,Y)$ can be understood in terms of the internal hom in $\sC$ by the following definition and result.

\begin{defn}\mylabel{topsym}
Let $\sC$ be a topological category over $B$ with a closed symmetric
monoidal structure given by a product $\sma_B$ and function object 
functor $F_B$, with unit object $S_B$. We say that $\sC$ is a \emph{topological
closed symmetric monoidal category over $B$}\index{category!topological closed symmetric monoidal over B@topological closed symmetric monoidal over $B$} if the tensors and products
are related by a natural isomorphism
\[X\wedge_B K \cong X\wedge_B(S_B\wedge_B K)\]
in $\sC$ for $K\in \sK_B$ and $X\in \sC$.
\end{defn}

\begin{prop}\mylabel{prop:PvsF} 
Let $\sC$ be a topological closed symmetric monoidal category over $B$.
Then, for $K\in\sK_B$ and $X$, $Y$, $Z\in \sC$, there are natural isomorphisms
\begin{gather*}
F_B(K, Y) \cong F_B(S_B\wedge_B K, Y),\\
P_B(X, Y)\iso P_B(S_B,F_B(X,Y)),\\
P_B(X\wedge_B Y, Z)\cong P_B(X,F_B(Y,Z))
\end{gather*}
in $\sC$ and a natural homeomorphism of based spaces
\[\sK_B(K,P_B(X,Y)) \cong \sC(S_B\wedge_B K, F_B(X,Y)).\]
\end{prop}

\section{Equivariant parametrized enriched categories}

Turning to the equivariant generalization, we give details of the context 
of topological $G$-categories, continuous $G$-functors, and natural $G$-maps
that we first alluded to in \S1.4. The discussion elaborates on that given
in \cite[II\S1]{MM}. Generically, we use notations of the form $\sC_G$ and $G\sC$ to denote a category $\sC_G$ enriched over the category $G\sK_*$ of 
based $G$-spaces and its associated ``$G$-fixed category'' $G\sC$ with the same objects and the $G$-maps between them; $G\sC$ is enriched over $\sK_*$. We shall write $(\sC_G,G\sC)$ for such a pair, and we shall refer to the pair as a ``$G$-category''.

In the terminology of enriched category theory, $G\sC$ is the 
underlying topological category of $\sC_G$. The hom objects of
$\sC_G$ are $G$-spaces $\sC_G(X,Y)$; $G$-functors and natural $G$-maps just 
mean functors and natural transformations enriched over  $G\sK_*$. 
Consistently with enriched category theory, the space 
$G\sC(X,Y) = \sC_G(X,Y)^G$ can be identified with the space of $G$-maps 
$S^0\rtarr \sC_G(X,Y)$. We call the points of $\sC_G(X,Y)$ ``arrows'' to distinguish them from the points of $G\sC(X,Y)$, which we call ``$G$-maps'', 
or often just ``maps'', with the equivariance understood. 

We cannot expect $\sC_G$ to have limits and colimits, but $G\sC$ is usually 
bicomplete. In many of our examples, both $\sC_G$ and $G\sC$ are closed symmetric
monoidal under functors $\sma_B$ and $F_B$.  For example, we have the closed 
symmetric monoidal $G$-category $(\sK_{G,B},G\sK_B)$ of ex-$G$-spaces over a 
$G$-space $B$ described in \S1.4.  

\begin{defn}\mylabel{defn:enrichBG}
A $G$-category $(\sC_G,G\sC)$\noteindex{CG@$\sC_B$}\noteindex{GC@$G\sC$} is \emph{$G$-topological over $B$}\index{category!G-topological over B@$G$-topological over $B$} if $\sC_G$ 
is enriched over $G\sK_B$ and bitensored over $\sK_{G,B}$. It follows that
$G\sC$ is enriched over $\sK_B$ and bitensored over $G\sK_B$. We say
that $(\sC_G,G\sC)$ is \emph{$G$-topologically bicomplete over $B$}\index{category!G-topologically bicomplete over B@$G$-topologically bicomplete over $B$} if, in addition, $G\sC$ is bicomplete. We write $P_B(X,Y)$ for the hom ex-$G$-space 
over $B$, and we write $X\sma_B K$ and $F_B(K,X)$ for the
tensor and cotensor in $\sC_G$, where $X$, $Y\in\sC_G$ and $K\in\sK_{G,B}$.  
Explicitly, we require bitensoring adjunction homeomorphisms of based $G$-spaces
\begin{equation}
\sC_G(X\wedge_B K, Y)\cong \sK_{G,B}(K, P_B(X,Y)) \cong \sC_G(X, F_B(K,Y)).
\mylabel{eqn:adjG}
\end{equation}
There result coherent unit and transitivity isomorphisms in $G\sC$
\begin{equation}\label{eqn:transagG}
X\iso X\sma_B S^0_B\quad \text{and}\quad X\sma_B (K\sma_B L) \iso (X\sma_B K)\sma_B L
\end{equation}
and also bitensoring adjunction isomorphisms of ex-$G$-spaces
\begin{equation}\label{eqn:intadjG}
P_B(X\wedge_B K, Y)\cong F_B(K, P_B(X,Y)) \cong P_B(X,F_B(K,Y)).
\end{equation}
Conversely, there is a natural homeomorphism of based $G$-spaces
\begin{equation}\label{eqn:underG}
\sC_G(X,Y)\cong \sK_{G,B}(S^0_B, P_B(X,Y)),
\end{equation}
and the isomorphisms (\ref{eqn:adjG}) follow from (\ref{eqn:intadjG}) by 
applying $\sK_{G,B}(S^0_B,-)$. Passage to $G$-fixed points from 
(\ref{eqn:adjG}) 
gives the bitensoring adjunction homeomorphisms of based spaces
\begin{equation}\label{eqn:adjG2}
G\sC(X\wedge_B K, Y)\cong G\sK_B(K, P_B(X,Y)) \cong G\sC(X, F_B(K,Y)).
\end{equation}
\end{defn}

We warn the reader that we shall not always adhere strictly to the notational 
pattern of \myref{defn:enrichBG} for our several layers of enrichment.
In particular, in the domain categories for our equivariant diagram spaces
and diagram spectra, only $\sC_G$ is of interest, not $G\sC$, and our notations will reflect that.  On the other hand, when studying model categories, it is always the bicomplete category $G\sC$ that is of fundamental interest.

If $(\sC_G,G\sC)$ is $G$-topological over $B$, then it is automatically
$G$-topo\-lo\-gic\-al (over $*$). This enrichment is recovered by taking (\ref{eqn:under}),
read equivariantly, as the definition of the based $G$-space $\sC_G(X,Y)$. 
Just as in the nonequivariant case, for based $G$-spaces $T$ and objects
$X$ of $\sC_G$, the tensors and cotensors in $\sC_G$ and $G\sC$ are given 
on objects by
\begin{equation}\label{SOB}
X\wedge_B T = X\wedge_B T_{B} \quad\text{and}\quad F_B(T,X) = F_B(T_{B}, X),
\end{equation}
using the notation of \myref{confusion} equivariantly. The required $G$-homeomorphisms 
\begin{equation}\label{SOB2}
\sC_G(X\wedge_B T, Y)\cong \sK_{G,*}(T, \sC_G(X,Y)) \cong \sC_G(X, F_B(T,Y))
\end{equation}
follow directly.

We have equivariant analogues of \myref{topsym} and \myref{prop:PvsF}.

\begin{defn}\mylabel{topsymG}
Let $(\sC_G,G\sC)$ be a $G$-topological $G$-category over $B$ with a closed 
symmetric monoidal structure given by a product $G$-functor $\sma_B$ and a
function object $G$-functor $F_B$, with unit object $S_B$. We say that $(\sC_G,G\sC)$ 
is a \emph{$G$-topological closed symmetric monoidal $G$-category over $B$}\index{category!G-toplogical closed symmetric@$G$-topological closed symmetric monoidal over $B$} if the tensors and 
products are related by a natural isomorphism
\[X\wedge_B K \cong X\wedge_B(S_B\wedge_B K)\]
in $G\sC$ for $K\in G\sK_B$ and $X\in G\sC$.
\end{defn}

\begin{prop}\mylabel{prop:PvsFG} 
Let $(\sC_G,G\sC)$ be a $G$-topological closed symmetric mon\-oid\-al $G$-category over $B$.
Then, for $K\in\sK_B$ and $X$, $Y$, $Z\in \sC$, there are natural isomorphisms
\begin{gather*}
F_B(K, Y) \cong F_B(S_B\wedge_B K, Y),\\
P_B(X, Y)\iso P_B(S_B,F_B(X,Y)),\\
P_B(X\wedge_B Y, Z)\cong P_B(X,F_B(Y,Z))
\end{gather*}
in $G\sC$ and there is a natural homeomorphism of based $G$-spaces
\[\sK_{G,B}(K,P_B(X,Y))\cong \sC_G(S_B\wedge_B K, F_B(X,Y)).\]
\end{prop}

\section{$G$-topological model $G$-categories}

We explain what it means for a $G$-topological $G$-category $(\sC_G, G\sC)$ 
over $B$ to have a $G$-topological model structure. This structure implies 
in particular that the homotopy category $\text{Ho}G\sC$ is bitensored over
the homotopy category $\text{Ho}G\sK$.  We need some notation.
Throughout this section, we consider maps  
$$i\colon W\rtarr X, \ j\colon V\rtarr Z,\  \text{and} \ p\colon E\rtarr Y$$
in $G\sC$ and a map $k\colon K\rtarr L$ in either  $G\sK_B$ or $G\sK_*$; in
the latter case we apply the functor $(-)_B = B\times (-)$ to $k$ and so regard it as a map in $G\sK_B$, as suggested in \myref{confusion}. We shall define the notion of a $G$-topological model category in terms of the induced map
\begin{equation}\label{gbox3}
\sC_G^\Box(i,p)\colon  \sC_G(X,E) \rtarr \sC_G(W,E)\times_{\sC_G(W,Y)}\sC_G(X,Y)
\end{equation}
of based $G$-spaces.  Passing to $G$-fixed points, this gives rise to a map
\begin{equation}\label{gbox3bis}
G\sC^\Box(i,p)\colon  G\sC(X,E) \rtarr G\sC(W,E)\times_{G\sC(W,Y)}G\sC(X,Y)
\end{equation}
of based spaces, and we have the following motivating observation.

\begin{lem}\mylabel{lemma:lppair}
The pair $(i,p)$ has the lifting property if and only if the function $G\sC^\Box(i,p)$ is surjective.
\end{lem}

\begin{defn}\mylabel{Gtopmodel}
Let $(\sC_G,G\sC)$ be a $G$-topological $G$-category over $B$ such that 
$G\sC$ is a model category. We say that the model structure is 
\emph{$G$-topological}\index{G-topological@$G$-topological}\index{model category!G-top@$G$-topological} if $\sC_G^\Box(i,p)$ is a fibration in 
$G\sK_*$ when $i$ is a cofibration and $p$ is a fibration and is acyclic 
when, further, either $i$ or $p$ is acyclic.
\end{defn}

\begin{rem}\mylabel{Gtopfamily}
The definition must refer consistently to either $h$-type or $q$-type model
structures. The resulting notions are quite different. We usually have in
mind a $q$-type model structure. In that case, the weak equivalences and
fibrations are often characterized by conditions on the $H$-fixed point maps
$f^H$ of a map $f$. If $\sF$ is a family of subgroups of $G$, such as the
family $\sG$ of compact subgroups, we can restrict attention to those $H\in
\sF$. The resulting $\sF$-equivalences and $\sF$-fibrations usually specify
another model structure on $G\sC$. In particular, we have the $\sF$-model
structure on $G\sK_*$. For the $qf$-type model structures of \S7.2, we must
start with a generating set $\sC$ that contains the orbits $G/H$ with $H\in
\sF\cap\sG$ and consists of $\sF\cap\sG$-cell complexes. We say that an
$\sF$-model structure on $G\sC$ is {\em $\sF$-topological} if the condition
of the previous definition holds when we use the $\sF$-notions of fibration,
cofibration and weak equivalence throughout. The observations of this
section generalize to $\sF$-topological model categories for any family
$\sF$.
\end{rem}

In addition to the map of $G$-spaces displayed in (\ref{gbox3}), we have a map
\begin{equation}\label{PBoxmap}
P^\Box_B(i,p)\colon  P_B(X,E)\rtarr P_B(W,E)\times_{P_B(W,Y)} P_B(X,Y)
\end{equation} 
of ex-$G$-spaces over $B$.  

\begin{warn}
We can define what it means for $(\sC_G,G\sC)$ to be $G$-topological 
{\em over $B$}, using the map $P^{\Box}_B(i,p)$ of ex-spaces rather than the 
map $\sC_G^\Box(i,p)$ of spaces. However, we know of no examples where this 
condition is satisfied. For example, $(\sK_{G,B},G\sK_B)$ is $G$-topological, by Theorems \ref{qoverB} and \ref{Gqfstr}, but, as \myref{ouchtoo} makes clear by
adjunction, we cannot expect it to be $G$-topological over $B$.
\end{warn}

Just as in the classical theory of simplicial or topological model categories, there are various equivalent reformulations of what it means for $G\sC$ to be $G$-topological.  To explain them, observe that the tensors and cotensors with
ex-$G$-spaces over $B$ give rise to induced maps
\begin{equation}\label{gbox1}
\ \ \ \ i\Box_B k\colon  (X\wedge_B K) \cup_{W\wedge_B K} (W\wedge_B L) 
\rtarr X\wedge_B L
\end{equation}
and
\begin{equation}\label{gbox2}
F^\Box_B(k,p)\colon F_B(L,E)\rtarr F_B(K,E)\times_{F_B(K,Y)} F_B(L,Y)
\end{equation}
of ex-$G$-spaces over $B$. If $(\sC_G,G\sC)$ is closed symmetric monoidal, 
then we also have the induced maps
\begin{equation}
i\Box_B j\colon  (X\wedge_B V) \cup_{W\wedge_B V} (W\wedge_B Z) 
\rtarr X\wedge_B Z
\end{equation}
and
\begin{equation}
F^\Box_B(j,p)\colon F_B(Z,E)\rtarr F_B(V,E)\times_{F_B(V,Y)} F_B(Z,Y)
\end{equation}
in $G\sC$. We have various adjunction isomorphisms relating these various $\Box$-product maps and $\Box$-function object maps. 

\begin{prop}\mylabel{prop:topveri}
If $k$ is a map of ex-$G$-spaces over $B$, then there are adjunction isomorphisms
\begin{equation}
P^\Box_B(i\Box_B k, p)\cong F^\Box_B(k,P^\Box_B(i,p)) 
\cong P^\Box_B(i,F^\Box_B(k,p))
\end{equation}
of maps of ex-$G$-spaces over $B$ and
\begin{equation}
\sC^\Box_G(i\Box_B k, p)\cong \sK^\Box_{G,B}(k,P^\Box_B(i,p))
\cong \sC^\Box_G(i,F^\Box_B(k,p))
\end{equation}
of maps of based $G$-spaces. If $k$ is a map of based $G$-spaces, then 
the last pair of isomorphisms can be rewritten as
\begin{equation}\label{gtopequiv}
\sC^\Box_G(i\Box_B k, p)\cong \sK^\Box_{G,*}(k,\sC^\Box_G(i,p))
\cong \sC^\Box_G(i,F^\Box_B(k,p)).
\end{equation}
When $(\sC_G,G\sC)$ is closed symmetric monoidal there are  
adjunction isomorphisms
\begin{equation}
P^\Box_B(i\Box_B k, p)\cong P^\Box_B(i,F^\Box_B(k,p))
\end{equation}
of maps of ex-$G$-spaces over $B$ and
\begin{equation}
\sC^\Box_G(i\Box_B k, p)\cong \sC^\Box_G(i,F^\Box_B(k,p))
\end{equation}
of maps of based $G$-spaces.
\end{prop}

Together with \myref{lemma:lppair}, this implies the promised alternative equivalent conditions that describe when a model category is $G$-topological.

\begin{prop}\mylabel{Gtopchar}
Let $(\sC_G,G\sC)$ be a $G$-topological $G$-category over $B$ such that $G\sC$ has a model structure. Then the following conditions are equivalent.
\begin{enumerate}[(i)]
\item The map $i\Box_B k$ of (\ref{gbox1}) is a cofibration in $G\sC$ if $i$ is a cofibration in $G\sC$ and $k$ is a cofibration in $G\sK_*$. It is acyclic if either $i$ or $k$ is acyclic.
\item The map $F^\Box_B(k,p)$ of (\ref{gbox2}) is a fibration in $G\sC$ if $p$ is a fibration in $G\sC$ and $k$ is a cofibration in $G\sK_*$. It is acyclic if either $p$ or $k$ is acyclic.
\item The map $\sC^\Box_G(i,p)$ of (\ref{gbox3}) is a fibration in $G\sK_*$ if $i$ is a cofibration in $G\sC$ and $p$ is a fibration in $G\sC$. It is acyclic if either $i$ or $p$ is acyclic.
\end{enumerate}
\end{prop}

\begin{proof}
The third condition is our definition of the model structure being 
$G$-topological. We prove that the first condition is equivalent to the 
third. A similar argument shows that the second condition is also equivalent
to the third. The map $\sC^\Box_G(i,p)$ is a fibration if and only if $(k,\sC^\Box_G(i,p))$ has the lifting property with respect to all acyclic cofibrations $k$ in $G\sK_*$. By \myref{lemma:lppair} and the first adjunction isomorphism in (\ref{gtopequiv}), that holds if and only if $(i\Box_B k, p)$ 
has the lifting property, that is, if and only if $i\Box_B k$ is an acyclic cofibration. If either $i$ or $p$ is acyclic, then we take $k$ to be a 
cofibration in $G\sK_*$ and argue similarly.
\end{proof} 

\chapter{The category of orthogonal $G$-spectra over $B$}

\section*{Introduction}

Intuitively, an orthogonal spectrum $X$ over $B$ consists of ex-spaces
$X(V)$ over $B$ and ex-maps $\si\colon X(V)\sma_B S^W\rtarr X(V\oplus W)$ for 
suitable inner product spaces $V$ and $W$. The orthogonal group $O(V)$ 
must act on $X(V)$, and $\si$ must be $(O(V)\times O(W))$-equivariant. The 
orthogonal group actions enable the definition of a good external smash
product. Moreover, they will later allow us to define stable weak equivalences
in terms of homotopy groups, as would not be possible if we only had 
actions by symmetric groups.  

Similarly, use of general inner product
spaces allows us to build in actions by a compact Lie group $G$ without 
difficulty. For non-compact Lie groups, we should ignore inner products
and use linear isomorphisms, replacing the compact
orthogonal group $O(V)$ by the general linear group $GL(V)$. 
However, as we explain in \S11.6, there are more serious problems
in generalizing to non-compact Lie groups; except in that section, 
we require $G$ to be a compact Lie group. 

Working equivariantly, we first describe $X$ as a suitable diagram of
ex-$G$-spaces in \S11.1. The domain category for our diagrams is denoted
$\sI_G$ and is independent of $B$. We then build in the structure maps $\si$ 
in \S11.2, where we define the category of orthogonal $G$-spectra over $B$.
In \S11.3, we show that it too can be described as a category of diagrams 
of ex-$G$-spaces. The domain category here is denoted $\sJ_{G,B}$. It does depend on $B$, as indicated by the notation. The formal properties of the category of ex-$G$-spaces 
over $B$ carry over to the category of orthogonal $G$-spectra over $B$, 
but there are some new twists. For example, our category of $G$-spectra 
over $B$ is enriched not just over based $G$-spaces, but more generally 
over ex-$G$-spaces over $B$. We discussed the relevant formalities in
the previous chapter.  This enhanced enrichment is essential to the 
definition of function $G$-spectra over $B$. 

We show in \S11.4 that the base change functors and their properties also 
carry over to these categories of parametrized $G$-spectra, and we discuss change of group functors and restriction to fibers in \S11.5.

\section{The category of $\sI_G$-spaces over $B$}

We recall the $G$-category $(\sI_G,G\sI)$ from \cite[II.2.1]{MM}. The objects and arrows of $\sI_G$ are finite dimensional $G$-inner product spaces and linear isometric isomorphisms. The maps of $G\sI$ are $G$-linear isometries. More precisely, as dictated by the general theory 
of \cite{MM, MMSS}, we take $\sI_G(V,W)$ as based with basepoint disjoint from the space of linear isometric isomorphisms $V\rtarr W$.  As in \cite[II.1.1]{MM}, the objects $V$ run over the collection $\sV$ of all representations that embed up to isomorphism in a given ``$G$-universe'' $U$, where a $G$-universe is a sum of countably many copies of representations in a set of representations that includes the trivial representation. We think in terms of a ``complete $G$-universe'', one that contains 
all representations of $G$, but the choice is irrelevant until otherwise 
stated. 
As in \cite[II.2.2]{MM}, we can restrict from $\sV$ to any cofinal 
subcollection $\sW$ that is closed under direct sums. 

Based $G$-spaces are ex-$G$-spaces over $*$, and $\sI_G$-spaces are defined in 
\cite[II.2.3]{MM} as $G$-functors $\sI_G\rtarr \sT_G$, where $\sT_G$ is the
$G$-category of compactly generated based $G$-spaces.  One can just as well drop the weak Hausdorff condition, which plays no necessary mathematical role in 
\cite{MM, MMSS}, and allow general $k$-spaces. With the notations of Part II, 
we can thus change the target $G$-category to 
$\sK_{G,*}$.  Then we generalize the definition to the parametrized context simply 
by changing the target $G$-category to the category $\sK_{G,B}$ of ex-$G$-spaces over a $G$-space $B$.  Thus we define an $\sI_G$-space $X$ over (and under) $B$ 
to be a $G$-functor $X\colon \sI_G\rtarr \sK_{G,B}$. Using nonequivariant arrows and equivariant maps, we obtain the $G$-category $(\sI_G\sK_B,G\sI\sK_B)$ of $\sI_G$-spaces.

To unravel definitions, for each representation $V\in\sV$ we are given an ex-$G$-space $X(V)$
over $B$, for each arrow (linear isometric isomorphism) $f\colon V\rtarr W$ we are given an arrow
(non-equivariant map) 
\[X(f)\colon X(V)\rtarr X(W)\]
of ex-$G$-spaces over $B$, and the continuous function
\[X\colon \sI_G(V,W)\rtarr \sK_{G,B}(X(V),X(W))\]
is a based $G$-map. An arrow $\al\colon  X\rtarr Y$ is just a natural transformation,
and a $G$-map is a $G$-natural transformation, for which each $\al_V\colon  X(V)\rtarr Y(V)$
is a $G$-map. For both arrows and $G$-maps, the naturality diagrams
$$\xymatrix{
X(V)\ar[r]^-{\al_V} \ar[d]_{X(f)} & Y(V) \ar[d]^{Y(f)}\\
X(W)\ar[r]_-{\al_W} & Y(W)}$$
must commute for all arrows $f\colon V\rtarr W$. The group $G$ acts on the space
$\sI_G\sK_B(X,Y)$ of arrows by levelwise conjugation.  The $G$-fixed category
is denoted by $G\sI\sK_B$. It has objects the $\sI_G$-spaces $X$ and maps the $G$-maps.

To study the parametrized enrichment of the $G$-category of orthogonal 
$G$-spectra over $B$, it is convenient to extend the domain category $\sI_G$, 
which is enriched over $\sK_{G,*}$, to a new domain category $\sI_{G,B}$ that 
is enriched over $\sK_{G,B}$.  Departing from the notational pattern of 
\myref{defn:enrichBG} and using \myref{confusion}, we define the hom ex-$G$-spaces over $B$ of $\sI_{G,B}$ by
\begin{equation}\label{IGB}
\sI_{G,B}(V,W) = \sI_G(V,W)_B \equiv B\times \sI_G(V,W).
\end{equation} 
If $X\colon \sI_G\rtarr \sK_{G,B}$ is an $\sI_G$-space, 
then the given based $G$-maps 
$$X\colon  \sI_G(V,W) \rtarr \sK_{G,B}(X(V),X(W))$$
correspond by adjunction (see (\ref{SOB}) and (\ref{SOB2})) to ex-$G$-maps
$$X(V)\sma_B \sI_{G,B}(V,W)\rtarr X(W).$$
In turn, these correspond by the internal hom adjunction to ex-$G$-maps
$$X\colon  \sI_{G,B}(V,W)\rtarr F_B(X(V),X(W)).$$
These give an equivalent version of the original $G$-functor $X$, but
now in terms of categories enriched over the category $G\sK_B$.

\begin{lem}\mylabel{omnilem}
The $G$-category $(\sI_G\sK_B,G\sI\sK_B)$ of $\sI_G$-spaces is equivalent 
to the $G$-category of $\sI_{G,B}$-spaces, where an $\sI_{G,B}$-space is a
$G$-functor $X\colon \sI_{G,B}\rtarr \sK_{G,B}$ enriched over $G\sK_B$.
\end{lem}

\begin{prop}\mylabel{omnispace} 
The $G$-category $(\sI_G\sK_B,G\sI\sK_B)$ is $G$-top\-o\-log\-i\-cal
over $B$ and thus also $G$-topological. Therefore the category
$G\sI\sK_B$ is topologically bicomplete over $B$.
\end{prop}
\begin{proof} 
We define tensor and cotensor $\sI_G$-spaces over $B$
$$ X\sma_B K \qquad\text{and}\qquad  F_B(K,X)$$
levelwise, where $K$ is an ex-$G$-space and $X$ is an $\sI_G$-space. 
For $\sI_G$-spaces $X$ and $Y$, we must define a parametrized morphism
ex-$G$-space $P_B(X,Y)$ over $B$. Parallelling a standard formal description 
of the $G$-space $\sI_G\sK_B(X,Y)$, we define $P_B(X,Y)$ to be the end
\begin{equation}\label{P}
P_B(X,Y) = \int_{\sI_{G,B}} F_B(X(V),Y(V)).
\end{equation}
Explicitly, it is the equalizer displayed in the following diagram of ex-$G$-spaces.
$$
\xymatrix{P_B(X,Y) \ar[d]\\
\prod_{V}F_B(X(V),Y(V))
\ar@<1ex>[d]^-{\tilde{\nu}}
\ar@<-1ex>[d]_-{\tilde{\mu}}\\
\prod_{V,W} F_B(\sI_{G,B}(V,W),F_B(X(V),Y(W))).}
$$
The products run over the objects and pairs of objects of a skeleton $sk\sI_G$
of $\sI_G$. The $(V,W)$th coordinate of $\tilde{\mu}$ is given by the composite of the
projection to $F_B(X(W),Y(W))$ and the $G$-map
$$F_{B}(X(W),Y(W))
\rtarr F_B(\sI_{G,B}(V,W),F_{B}(X(V),Y(W)))$$
adjoint to the composite ex-$G$-map
\[\xymatrix{
F_{B}(X(W),Y(W))\sma_B \sI_{G,B}(V,W) \ar[d]^-{\text{id}\sma_B X}\\
F_{B}(X(W),Y(W))\sma_B F_{B}(X(V),X(W)) \ar[d]^-{\com}\\
F_B(X(V),Y(W)).}\]
The $(V,W)$th coordinate of $\tilde{\nu}$ is the composite of the
projection to $F_B(X(V),Y(V))$ and the $G$-map
$$\tilde{\nu}_{V,W}\colon F_B(X(V),Y(V))
\rtarr F_B(\sI_{G,B}(V,W), F_B(X(V),Y(W))$$
adjoint to the composite ex-$G$-map
$$\xymatrix{
\sI_{G,B}(V,W)\sma_B F_B(X(V),Y(V)) \ar[d]^-{Y\sma_B \text{id}}\\
F_{B}(Y(V),Y(W))\sma_B F_B(X(V),Y(V)) \ar[d]^-{\com}\\
F_B(X(V),Y(W)).\\}$$
Passage to ends from the isomorphisms of ex-$G$-spaces
\[F_B(X(V) \wedge_B K, Y(V))
\cong F_B(K, F_B(X(V),Y(V)))
\cong F_B(X(V), F_B(K,Y(V)))\]
gives natural isomorphisms of ex-$G$-spaces
\begin{equation}\label{Pad0}
P_B(X \wedge_B K, Y)
\cong F_B(K, P_B(X,Y))
\cong P_B(X, F_B(K,Y)).
\end{equation}
With these constructions, we see that $(\sI_G\sK_B,G\sI\sK_B)$ 
is $G$-top\-o\-log\-i\-cal over $B$; compare \myref{defn:enrichBG} and the 
discussion following it. The last statement follows since $G\sI\sK_B$ 
is complete and cocomplete, with limits and colimits constructed levelwise 
from the limits and colimits in $G\sK_{B}$. 
\end{proof}

We have several kinds of smash products and function objects in this context.
For $\sI_G$-spaces $X$ and $Y$ over $B$, define the \emph{``external'' smash 
product}\index{external smash product} $X \barwedge_B Y$\noteindex{XbwY@$X\barwedge_B Y$} by
$$X \barwedge_B Y = \sma_B\com (X\times Y)\colon  \sI_G\times \sI_G \rtarr \sK_{G,B}.$$
Thus $(X\barwedge_B Y)(V,W) = X(V)\sma_B Y(W)$.  Here we have used the word ``external''
to refer to the use of pairs of representations, as is usual in the theory of diagram
spectra.  It is standard category theory \cite{Day, MMSS} to use left Kan extension to 
internalize this external smash product over $B$. This gives the internal smash product 
$X\sma_B Y$ of
$\sI_G$-spaces over $B$, which is again an $\sI_G$-space over $B$. For an $\sI_G$-space $Y$ 
over $B$ and an $(\sI_G\times\sI_G)$-space $Z$ over $B$, define the {\em external function 
$\sI_G$-space over $B$}, denoted $\bar F_B(Y,Z)$, by
$$\bar F_B(Y,Z)(V) = P_B(Y,Z\langle V\rangle),$$
where $Z\langle V\rangle (W) = Z(V,W)$. It is mainly to allow this definition that we 
need the morphism ex-$G$-spaces $P_B(-,-)$. It is also formal to obtain an internal function
$\sI_G$-space functor $F_B$ on $\sI_G$-spaces over $B$ by use of right Kan extension
\cite{Day, MMSS}. Using these internal smash product and function $\sI_G$-space functors,
we obtain the following result. Recall \myref{topsymG} and \myref{prop:PvsFG}. 

\begin{thm} $(\sI_G\sK_B,G\sI\sK_B)$ is a $G$-top\-o\-lo\-gi\-cal 
closed symmetric mon\-oid\-al $G$-category over $B$.
\end{thm}

\begin{rem}\mylabel{extsmash1} 
In the theory of ex-spaces, we also have the ``external smash product'' of ex-spaces over different base spaces defined in \S2.5.  Using the two different notions of ``external'' together, we obtain the definition of the ``external external smash product'' of an $\sI_G$-space $X$ over $A$ and an $\sI_G$-space $Y$ over $B$; it is an  $(\sI_G\times \sI_G)$-space over $A\times B$.  We write $X\barwedge Y$ for the left Kan extension internalization of this smash product. Thus $X\barwedge Y$ is an $\sI_G$-space over $A\times B$.  Similarly, using the external function ex-space construction $\bar{F}$ of \S2.5, for an $\sI_G$-space $Y$ over $B$ and an $\sI_G$-space $Z$ over $A\times B$, we obtain the ``internalized external function $\sI_G$-space'' $\bar{F}(Y,Z)$ over $A$. Notationally, use of $\barwedge$ and $\bar{F}$ without an ensuing subscript always denotes these internalized external operations with respect to varying base spaces. We shall return to these functors in \myref{extsmash2}.  

Similarly, but more simply, we have the ``external tensor'' $K\barwedge Y$ of an ex-$G$-space $K$ over $A$ and an $\sI_G$-space $Y$ over $B$, which again is an $\sI_G$-space over $A\times B$. When $A = *$, this is just the tensor of based $G$-spaces with $\sI_G$-spaces over $B$. The case $B=*$ shows how to construct an
$\sI_G$-space over $A$ from an ex-$G$-space over $A$ and an $\sI_G$-space.
Since these external tensors can be view as special cases of external smash products, via variants of \myref{topsymG} and (\ref{smasma}) below, we shall not treat them formally and shall not repeat the definitions on the $G$-spectrum level. However, we shall find several uses for them. 
\end{rem}

\section{The category of orthogonal $G$-spectra over $B$}

For a representation $V$ of $G$ and an $\sI_G$-space $X$, we define
\begin{equation}\label{SIOM}
\SI^V_B X = X\sma_B S^V_B \qquad\text{and}\qquad \OM^V_B X = F_B(S^V_B,X),
\end{equation}
where $S^V$ is the one-point compactification of $V$.  

\begin{defn}\mylabel{spheres} 
Define the $G$-sphere $S_B$, written $S_{G,B}$ when necessary for clarity, 
to be the $\sI_G$-space over $B$ that sends $V$ to $S_B^V$.
\end{defn}
 
Clearly $S_B^V\sma_B S_B^W\iso S_B^{V\oplus W}$, and the functor $S_{B}$ is strong symmetric monoidal, where the monoidal structure on $\sI_G$ is given by direct sums. It follows that $S_{B}$ is a commutative monoid in the symmetric monoidal category $G\sI\sK_B$, and we can define $S_{B}$-modules $X$ in terms of (right) actions $X\sma_B S_{B}\rtarr X$. These $S_{B}$-modules are our orthogonal $G$-spectra over $B$, but it is more convenient to give the definition using the equivalent reformulation in terms of the external 
smash product.

\begin{defn} 
An \emph{$\sI_G$-spectrum, or orthogonal $G$-spectrum, over $B$}\index{spectrum!over B@over $B$} is an $\sI_G$-space $X$ over $B$ together with a structure $G$-map 
$$\si\colon  X\barwedge_B S_{B}\rtarr X\com \oplus$$ 
such that the
evident unit and associativity diagrams commute. Thus we have compatible 
equivariant structure maps
$$\si\colon  \SI_B^W X(V) = X(V)\sma_B S_B^W\rtarr X(V\oplus W).$$
Let $\sS_{G,B}$\noteindex{SGB@$\sS_{G,B}$} denote the topological $G$-category of $\sI_G$-spectra over $B$
and arrows $f\colon  X\rtarr Y$ that commute with the structure maps, with $G$ acting
by conjugation on arrows.  Let $G\sS_{B}$\noteindex{GSB@$G\sS_B$} denote the topological category of 
$\sI_G$-spectra over $B$ and $G$-maps (equivariant arrows) between them.
\end{defn}

\begin{defn} 
Define the suspension orthogonal $G$-spectrum functor\noteindex{SiB@$\SI^\infty_B$}
and the $0$th ex-$G$-space functor\noteindex{OmB@$\OM^\infty_B$} 
$$\SI^{\infty}_B\colon \sK_{G,B}\rtarr \sS_{G,B} \ \  \text{and}\ \
\OM^{\infty}_B\colon \sS_{G,B}\rtarr \sK_{G,B}$$
by 
$(\SI^{\infty}_B K)(V) = \SI_B^VK$, with the 
evident isomorphisms as structure maps, and $\OM^{\infty}_B X = X(0)$. Then
$\SI^{\infty}_B$ and $\OM^{\infty}_B$ give left and right adjoints between 
$\sK_{G,B}$ and $\sS_{G,B}$ and, on passage to $G$-fixed points, between $G\sK_B$ and $G\sS_B$.
\end{defn}

The category $G\sS_{B}$ is our candidate for a good category of para\-me\-trized $G$-spectra 
over $B$. It inherits all of the properties of the category $G\sI\sK_B$ of $\sI_G$-spaces 
that were discussed in the previous section and, in the case $B=*$, it is exactly the 
category $G\sS$ of orthogonal $G$-spectra that is studied in \cite{MM}.  We summarize its formal 
properties in the following omnibus theorem. In the language of \S10.2, much of it can
be summarized by the assertion that the $G$-category $(\sS_{G,B},G\sS_B)$ is a 
$G$-topological closed symmetric monoidal $G$-category over $B$, but we prefer to be 
more explicit than that.

\begin{thm}\mylabel{omnigood}
The $G$-category $\sS_{G,B}$ is enriched over $G\sK_{B}$ and is tensored and cotensored over $\sK_{G,B}$.  The category $G\sS_{B}$ is enriched over $\sK_{B}$ and is tensored and cotensored over $G\sK_B$. The $G$-category $\sS_{G,B}$ and the category $G\sS_{B}$ admit smash product and function spectrum functors $\sma_B$ and $F_B$ under which they are closed symmetric monoidal with unit object $S_{B}$. Let $X$ and $Y$ be orthogonal $G$-spectra over $B$ and $K$ be an ex-$G$-space over $B$. The morphism ex-$G$-spaces $P_B(X,Y)$ can be specified by
$$P_B(X,Y) = \OM^{\infty}_B F_B(X,Y),$$
and there are natural isomorphisms
$$\SI^{\infty}_B K\iso S_B\sma_B K \qquad\text{and}\qquad \OM^{\infty}_BX\iso P_B(S_B,X).$$
The tensors and cotensors are related to smash products and function $G$-spectra by natural isomorphisms
\begin{equation}\label{smasma}
X\sma_B K \iso X\sma_B \SI^{\infty}_B K 
\qquad\text{and}\qquad
F_B(K,X)\iso F_B(\SI^{\infty}_B K,X)
\end{equation}
of orthogonal $G$-spectra.  There are natural isomorphisms
\begin{equation}\label{Pad!}
P_B(\SI^{\infty}_BK,X) \iso F_B(K,\OM^{\infty}_BX)
\end{equation}
and
\begin{equation}\label{Pad0S}
P_B(X \wedge_B K, Y)
\cong F_B(K, P_B(X,Y))
\cong P_B(X, F_B(K,Y))
\end{equation}
of ex-$G$-spaces,
\begin{equation}\label{Pad1S}
\sS_{G,B}(X \wedge_B K, Y)
\cong \sK_{G,B}(K,P_B(X,Y))
\cong \sS_{G,B}(X, F_B(K,Y))
\end{equation}
of based $G$-spaces, and
\begin{equation}\label{Pad2S}
G\sS_B(X \wedge_B K, Y)
\cong G\sK_B(K, P_B(X,Y))
\cong G\sS_B(X, F_B(K,Y))
\end{equation}
of based spaces. Moreover, $G\sS_{B}$ is $G$-topologically bicomplete over $B$.
\end{thm}\begin{proof} For the enrichment, the $G$-space $\sS_{G,B}(X,Y)$ is the evident sub $G$-space
of $\sI_G\sK_B(X,Y)$, and the space $G\sS_B(X,Y)$ is the evident sub space of $G\sI\sK_B(X,Y)$.
The tensors and cotensors in $\sS_{G,B}$ are constructed in $\sI_G\sK_B$ and given induced 
structure maps. The limits and colimits in $G\sS_B$ are constructed in the same way.  
As in \cite[II\S3]{MM}, we think of orthogonal $G$-spectra over $B$ as 
$S_{B}$-modules, and we construct the smash product and function spectra functors 
by passage to coequalizers and equalizers from the smash product and function $\sI_G$-space 
functors, exactly as in the definition of tensor products and hom functors in algebra.  We
have defined $P_B(X,Y)$ in the statement, but we shall give a more intrinsic alternative
description later. The first isomorphism of (\ref{smasma}) is given by unit and associativity 
relations
$$X\sma_B K \iso (X\sma_B S_{B})\sma_B K \iso X \sma_B \SI^{\infty}_B K.$$
The second follows from the Yoneda lemma since
\begin{align*}
G\sS_B(X,F_B(K,Y)) & \iso G\sS_B(X\sma_B K,Y)\\
 & \cong  G\sS_B(X\sma_B\SI^{\infty}_B K,Y) \\
 & \iso   G\sS_B(X,F_B(\SI^{\infty}_B K,Y)).
\end{align*}
Now (\ref{Pad!}) and (\ref{Pad0S}) follow from already established adjunctions. For part of the latter, we apply $\OM^{\infty}_B$ to the composite isomorphism
\begin{align*}
F_B(X\sma_B K,Y) & \iso F_B(X\sma_B \SI^{\infty}_B K,Y) \\
& \cong  F_B(X,F_B(\SI^{\infty}_B K,Y))\\
& \cong  F_B(X,F_B(K,Y)).
\end{align*}
Comparisons of definitions, seen more easily from (\ref{altPB}) below, give
\begin{equation}
\sS_{G,B}(X,Y) = \sK_{G,B}(S_B^0,P_B(X,Y))
\end{equation}
and
\begin{equation}
G\sS_B(X,Y)\iso G\sK_B(S^0_B,P_B(X,Y)).
\end{equation}
Therefore the isomorphisms (\ref{Pad1S}) and (\ref{Pad2S}) follow from (\ref{Pad0S}). 
\end{proof}

As noted in \S10.1, we obtain the following corollary by replacing 
$K$ with  $T_B$ for a based $G$-space $T$ in the 
tensors and cotensors of the theorem.  Of course, these tensors 
and cotensors with $G$-spaces could just as well be defined directly.  
It will be important in our discussion of model category structures 
to keep separately in mind the tensors and cotensors over ex-$G$-spaces 
over $B$ and over based $G$-spaces.

\begin{cor}\mylabel{omnigoodcor} The $G$-category $\sS_{G,B}$ is enriched over $G\sK_{*}$ and is 
tensored and cotensored over $\sK_{G,*}$. The category $G\sS_{B}$ is enriched over
$\sK_{G,*}$ and is tensored and cotensored over $G\sK_*$. Thus, for orthogonal $G$-spectra $X$ 
and $Y$ and based $G$-spaces $T$,
\begin{equation}\label{tencoten2}
\sS_{G,B}(X\sma_B T,Y)\iso \sK_{G,*}(T,\sS_{G,B}(X,Y))\iso \sS_{G,B}(X,F_B(T,Y))
\end{equation}
and
\begin{equation}
\label{tencoten2G}
G\sS_{B}(X\sma_B T,Y)\iso G\sK_{*}(T,\sS_{G,B}(X,Y))\iso G\sS_B(X,F_B(T,Y)).
\end{equation}
\end{cor}

We have the parallel definition of $G$-prespectra over $B$.

\begin{defn}  
A \emph{$G$-prespectrum $X$ over $B$}\index{prespectrum over B@prespectrum over $B$} consists of ex-$G$-spaces
$X(V)$ over $B$ for $V\in \sV$ together with structure $G$-maps 
$\si\colon \SI_B^WX(V)\rtarr X(V\oplus W)$ such that $\si$ is the identity
if $W=0$ and the following diagrams commute.
$$\xymatrix{
\SI^Z_B\SI^W_B X(V)\ar[d]_{\SI^Z_B\si} \ar[r]^-{\iso} 
& \SI^{W\oplus Z}_B X(V) \ar[d]^{\si}\\
\SI^Z_BX(V\oplus W) \ar[r]_-{\si} & X(V\oplus W\oplus Z)\\}$$
Let $\sP_{G,B}$\noteindex{PGB@$\sP_{G,B}$} denote the
$G$-category of $G$-prespectra and nonequivariant arrows, and let $G\sP_B$\noteindex{GPB@$G\sP_B$} denote
its $G$-fixed category of $G$-prespectra and $G$-maps. There result 
forgetful functors 
$$\bU\colon  \sS_{G,B}\rtarr \sP_G \qquad\text{and}\qquad  \bU\colon  G\sS_{B}\rtarr G\sP_{B}.$$
\end{defn}

The categories $\sP_{G,B}$ and $G\sP_{B}$ enjoy the same properties that were specified 
for $\sS_{G,B}$ and $G\sS_{B}$ in \myref{omnigood} and \myref{omnigoodcor}, 
except for the statements about smash product and function spectra.  Here,
since we do not have the internal hom functor $F_B$, we must give
an alternative direct description of $P_B(X,Y)$, as in (\ref{altPB}) below.  

\section{Orthogonal $G$-spectra as diagram ex-$G$-spaces}

Arguing as in \cite[\S2]{MMSS} and \cite[II\S4]{MM}, we
construct a new domain category $\sJ_{G,B}$\noteindex{JGB@$\sJ_{G,B}$} which
has the same object set $\sV$ as $\sI_G$ and, like $\sI_{G,B}$, is enriched 
over $G\sK_{B}$. It builds in spheres in such a way that the category of $\sI_G$-{\em spectra}\, over $B$ is equivalent to the category 
of  $\sJ_{G,B}$-{\em spaces}\, over $B$. Here, just as for $\sI_{G,B}$ 
in \myref{omnilem}, we understand a $\sJ_{G,B}$-space to be an enriched 
$G$-functor $X\colon \sJ_{G,B}\rtarr \sK_{G,B}$. Thus it is specified 
by ex-$G$-spaces $X(V)$ and ex-$G$-maps
$$X\colon \sJ_{G,B}(V,W)\rtarr F_B(X(V), X(W)).$$
To construct $\sJ_{G,B}$, recall from \cite[II\S4]{MM} that we have a topological $G$-category 
$\sJ_G$ with object set $\sV$ such that the category of $\sI_G$-spectra is equivalent to the category of $\sJ_G$-spaces. We define 
\begin{equation}\label{JG1}
\sJ_{G,B}(V,W) = \sJ_G(V,W)_B,
\end{equation}
just as we defined $\sI_{G,B}$ in (\ref{IGB}), and the desired equivalence
of categories follows. Rather than repeat either of the different constructions of $\sJ_G$ given in \cite{MMSS} and \cite{MM}, we shall shortly give a direct
description of $\sJ_{G,B}$. The intuition is that an extension of an $\sI_{G,B}$-space to a $\sJ_{G,B}$-space builds in an action by $S_{B}$.

The alternative description of $G\sS_{B}$ as the category of enriched $G$-functors 
$\sJ_{G,B}\rtarr\sK_{G,B}$ and enriched $G$-natural transformations 
leads to a more 
conceptual proof of \myref{omnigood}: it is a specialization of 
general results 
about diagram categories of enriched functors.  In analogy with (\ref{P}) we could have 
defined $P_B(X,Y)$ to be the end
\begin{equation}\label{altPB}
P_B(X,Y) = \int_{\sJ_{G,B}} F_B(X(V),Y(V))
\end{equation}
and derived the isomorphism (\ref{Pad0S}) just as we derived (\ref{Pad0})
in the previous section. By the Yoneda lemma, the two definitions of $P_B(X,Y)$ 
agree.  With this description of $P_B$, some of the adjunctions in 
\myref{omnigood} become more transparent. 

This leads to an alternative description of $\sJ_{G,B}$ in terms of 
$\sI_{G,B}$, following \cite[2.1]{MMSS}. We have the represented functors 
$V^*\colon \sI_{G}\rtarr \sK_{G,B}$ specified by $V^*(W) = \sI_{G,B}(V,W)$.  
If $X$ is an $\sI_G$-space, such as $V^*$, then the smash product
$X\sma_B S_{B}$ in the category of $\sI_G$-spaces is a ``free'' orthogonal 
$G$-spectrum over $B$. Let
\begin{equation}\label{JG2}
\sJ_{G,B}(V,W) = P_{B}(W^*\sma_B S_{B}, V^*\sma_B S_{B}),
\end{equation}
with the evident composition. Then we can mimic the arguments of \cite[\S\S2, 23]{MMSS}
to check that the category of $\sJ_{G,B}$-spaces is equivalent to the category of
$\sI_{G}$-spectra over $B$. An enriched Yoneda lemma argument \cite[2.4]{Ke} shows that this 
description of $\sJ_{G,B}$ coincides up to isomorphism with our original one.  

Although we will not have occasion to quote it formally, we record the following consequence 
of the identification of $\sI_G$-spectra over $B$ with $\sJ_{G,B}$-spaces.

\begin{lem}\mylabel{save} For any enriched $G$-functor $T\colon\sK_{G,B}\rtarr \sK_{G,B}$
and orthogonal $G$-spectrum $X$ over $B$, the composite functor $T\com X$ is an orthogonal 
$G$-spectrum over $B$. Similarly, an enriched natural transformation 
$\xi\colon T\rtarr T'$ induces a natural $G$-map 
$\xi\colon T\com X\rtarr T'\com X$. 
\end{lem}\begin{proof} The enriched functor $T$ is given by maps
$$T\colon F_B(K,L)\rtarr F_B(T(K),T(L)).$$ 
Composing levelwise with $X$ gives maps
$$\sJ_{G,B}(V,W)\rtarr F_B(T(X(V)),T(X(W)))$$
that specify $T\com X$. It is a direct categorical 
implication of the fact that $T$ is an enriched functor that there are natural maps of 
ex-$G$-spaces
$$T(K)\sma_B L\rtarr T(K\sma_B L) \qquad\text{and}\qquad TF_B(K,L)\rtarr F_B(K,T(L))$$
for ex-$G$-spaces $K$ and $L$. This explains more concretely why the 
structure maps of $X$ induce structure maps for $T\com X$. Similarly, since $\xi$ is enriched, 
it is given by maps from the unit ex-$G$-space $S^0_B$ to $F_B(T(K),T'(K))$ such that the 
appropriate diagrams commute. We specialize to $K=X(V)$ to obtain 
$\xi\colon T\com X\rtarr T'\com X$.
\end{proof}

The following functors relating ex-$G$-spaces to orthogonal $G$-spectra over $B$ play a central role in our theory. In particular, they give ``negative 
dimensional'' spheres $\SI^{\infty}_VS^0_B = S^{-V}_B$. 

\begin{defn}\mylabel{FVs}  
Let $V^* = V^*_B$ denote the represented $\sJ_{G,B}$-space 
specified by $V^*(W) = \sJ_{G,B}(V,W)$.  Define the \emph{shift desuspension 
functor}\index{functor!shift desuspension --}\index{shift desuspension functor}\noteindex{FV@$F_V$}
$$F_V \colon \sK_{G,B}\rtarr \sS_{G,B}$$ 
by letting
$F_V K = V^*\sma_B K$ for an ex-$G$-space $K$.  Let 
$\text{Ev}_V\colon  \sS_{G,B}\rtarr \sK_{G,B}$ be the functor given by evaluation at $V$.\index{functor!evaluation --}\index{evaluation functor}\noteindex{EvV@$\text{Ev}_V$}
The alternative notations
$$\SI^{\infty}_V K = F_VK \qquad\text{and}\qquad \OM^{\infty}_V K = \text{Ev}_V$$ 
are often used. In particular,  $F_0 = \SI^{\infty}_0= \SI^{\infty}_B$ and 
$\text{Ev}_0 = \OM^{\infty}_0 = \OM^{\infty}_B$.
\end{defn}

\begin{lem}\mylabel{FVEV} The functors $F_V$ and $\text{Ev}_V$ are left and 
right adjoint, and there is a natural isomorphism
$$ F_V K\sma_B F_W L \iso F_{V\oplus W}(K\sma_B L).$$ 
\end{lem}\begin{proof} The first statement is clear, and the verification of the
second statement is formal, as in \cite[\S1]{MMSS}. 
\end{proof}

\section{The base change functors $f^*$, $f_!$, and $f_*$}

From now on, we drop the adjective ``orthogonal'' (or prefix $\sI_G$), and we 
generally take the equivariance for granted, referring to orthogonal $G$-spectra
over $B$ just as spectra over $B$.  We return $G$ to the notations when
considering change of groups, or for emphasis, but otherwise $G$-actions 
are tacitly assumed throughout. 
 
We first show that the results on base change functors proven for
ex-spaces in \S2.2 extend to parametrized spectra. We then show
that the results in \S2.5 relating external and internal smash
product and function ex-spaces also extend to parametrized spectra. Let 
$A$ and $B$ be base $G$-spaces.

\begin{thm}\mylabel{Wirth} Let $f\colon A\rtarr B$ be a $G$-map. Let $X$ be in
$\sS_{G,A}$ and let $Y$ and $Z$ be in $\sS_{G,B}$. There are $G$-functors 
\[f_!\colon \sS_{G,A} \rtarr \sS_{G,B},\qquad
f^*\colon \sS_{G,B} \rtarr \sS_{G,A}, \qquad
f_*\colon \sS_{G,A} \rtarr \sS_{G,B}\]
and $G$-adjunctions
\[\sS_{G,B}(f_!X,Y)\iso \sS_{G,A}(X,f^*Y) \qquad\text{and}\qquad
\sS_{G,A}(f^*Y,X)\iso \sS_{G,B}(Y,f_*X).\]
On passage to $G$-fixed points levelwise, there result functors
\[f_!\colon G\sS_{A} \rtarr G\sS_{B},\qquad
f^*\colon G\sS_{B} \rtarr G\sS_{A}, \qquad
f_*\colon G\sS_{A} \rtarr G\sS_{B}\]
and adjunctions
\[G\sS_{B}(f_!X,Y)\iso G\sS_{A}(X,f^*Y)\qquad\text{and}\qquad
G\sS_{A}(f^*Y,X)\iso G\sS_{B}(Y,f_*X).\]
The functor $f^*$ is closed symmetric monoidal. Therefore, 
by definition and implication, $f^*S_{B}\iso S_{A}$ and there 
are natural isomorphisms
\begin{gather}\label{one}
f^*(Y\sma_B Z)\iso f^*Y\sma_A f^*Z,\\[1ex]
\label{two}
F_B(Y,f_*X) \iso f_*F_A(f^*Y,X),\\[1ex]
\label{three}
f^*F_B(Y,Z)\iso F_A(f^*Y,f^*Z),\\[1ex]
\label{four}
f_{!}(f^*Y\sma_A X)\iso Y\sma_B f_{!}X,\\[1ex]
\label{five}
F_B(f_{!}X,Y)\iso f_*F_A(X,f^*Y).
\end{gather}
\end{thm}\begin{proof} We define the functors $f^*$, $f_!$, and $f_*$ 
levelwise.  This certainly gives well-defined functors on $\sI_G$-spaces
that satisfy the appropriate adjunctions there. We shall show shortly 
that these functors preserve $\sI_G$-spectra. For a based $G$-space $T$, $f^*(T_B) \iso T_A$, and this implies $f^*S_{B}\iso S_{A}$. If we replace $\sI_G$-spectra by 
$\sI_G$-spaces and replace the internal smash product and function object 
functors ($\sma$ and $F$) by their 
external precursors ($\barwedge$ and $\bar F$), then everything is 
immediate by levelwise application of the corresponding results for
ex-spaces.  Still working with $\sI_G$-spaces, we first show how to
internalize the isomorphisms (\ref{one}) and (\ref{four}) by use of the 
universal property of left Kan extension. Indeed, noting that
$(f_*X)\circ \oplus \iso f_*(X\circ \oplus)$, and similarly for $f^*$ and $f_!$, 
we have
\begin{align*}
\sI_G\sK_{A}(f^*(Y\wedge_B Z), X)
&\cong  \sI_G\sK_{B}(Y\wedge_B Z, f_* X) \\
&\cong  (\sI_G\times\sI_G)\sK_{B}(Y\barwedge_B Z, f_*X\circ \oplus)\\
&\cong  (\sI_G\times\sI_G)\sK_{A}(f^*(Y\barwedge_B Z), X\circ \oplus )\\
&\cong  (\sI_G\times\sI_G)\sK_{A}(f^*Y\barwedge_A f^*Z, X\circ \oplus)\\
&\cong  \sI_G\sK_{A}(f^*Y \wedge_A f^*Z, X)
\end{align*}
and
\begin{align*}
\sI_G\sK_{B}(f_!X\wedge_B Y,Z)
&\cong  (\sI_G\times\sI_G)\sK_{B}(f_!X\barwedge_B Y,Z\circ \oplus)\\
&\cong  (\sI_G\times\sI_G)\sK_{B}(f_!(X\barwedge_A f^*Y),Z\circ \oplus)\\
&\cong  (\sI_G\times\sI_G)\sK_{A}(X\barwedge_A f^*Y,f^*Z\circ \oplus) \\
&\cong  \sI_G\sK_{A}(X\wedge_A f^*Y,f^*Z)\\
&\cong  \sI_G\sK_{A}(f_!(X\wedge_A f^*Y),Z).
\end{align*}
As explained in \cite[\S\S2--3]{FHM}, the remaining isomorphisms on the 
$\sI_G$-space level follow formally.  

We must show that our functors
on $\sI_G$-spaces preserve $\sI_G$-spectra. The given structure map
$\si\colon Y\barwedge_B S_{B}\rtarr Y\com \oplus$ gives rise
via the external version of (\ref{one}) to the required structure map 
$$f^*Y\barwedge_A S_{A}\iso f^*(Y\barwedge_B S_{B})\rtarr f^*Y\com \oplus.$$
Similarly, the given structure map 
$\si\colon X\barwedge S_{A}\rtarr X\com \oplus$
gives rise to the required structure map 
$$ f_!X\barwedge_B S_{B} \iso f_!(X\barwedge_A S_{A})\rtarr f_!X\com \oplus.$$
As in \cite[(3.6)]{FHM}, there is a canonical natural map, not usually an isomorphism,
$$\pi\colon f_*X\barwedge_B Y\rtarr f_*(X\barwedge_A f^*Y).$$
Taking $Y=S_{B}$, we see that $\si$ also induces the required structure map
$$ f_*X\barwedge_B S_{B} \rtarr f_*(X\barwedge_A S_{A})\rtarr f_*X\com \oplus.$$
Now the spectrum level adjunctions follow directly from their $\sI_G$-space analogues.
The spectrum level isomorphisms (\ref{one}) and (\ref{four}) follow from their 
$\sI_G$-space analogues by comparisons of coequalizer diagrams, and the remaining isomorphisms again follow formally.
\end{proof}

\begin{rem}\mylabel{FVvsf*} 
Since the base change functors are defined levelwise, they 
commute with the evaluation functors $\text{Ev}_V$. These commutation
relations for the right adjoints $f_*$ and $f^*$ imply conjugate 
commutation isomorphisms
$$f^*F_V\iso F_Vf^* \qquad\text{and}\qquad f_!F_V \iso F_V f_!$$ 
of left adjoints. In particular, 
$$f^*\SI^{\infty}_B\iso \SI^{\infty}_A f^* \qquad\text{and}\qquad
f_! \SI^{\infty}_A \iso \SI^{\infty}_B f_!.$$ 
Via (\ref{smasma}), these isomorphisms 
and the isomorphisms of the theorem imply isomorphisms relating base change functors 
to tensors and cotensors.  For example (\ref{four}) implies isomorphisms
$$ f_!(f^*Y\sma_A K) \iso Y\sma_B f_!K \qquad\text{and}\qquad
f_!(f^*L\sma_A X) \iso L\sma_B f_!X.$$
Here $K$ and $L$ are ex-spaces over $A$ and $B$ and $X$ and $Y$ are spectra over $A$ and $B$.
\end{rem}

The following result is immediate from its precursor \myref{Mackey0} for ex-spaces.

\begin{prop}\mylabel{Mackey} 
Suppose given a pullback diagram of $G$-spaces 
$$\xymatrix{
C \ar[r]^-{g} \ar[d]_{i} & D \ar[d]^{j} \\
A \ar[r]_{f} & B.\\}$$
Then there are natural isomorphisms of functors 
\begin{equation}\label{bases}
j^*f_{!} \iso g_{!}i^*, \qquad f^*j_* \iso i_*g^*, 
\qquad f^*j_{!}\iso i_!g^*, \qquad j^*f_*\iso g_*i^*.
\end{equation}
\end{prop}

Returning to \myref{extsmash1}, we have the following important results
on external smash product and function spectra and their internalization by
means of base change along diagonal maps.

\begin{prop}\mylabel{extsmash2}  Let $X$ be a spectrum over $A$, $Y$ be 
a spectrum over $B$, and $Z$ be a spectrum over $A\times B$. There
is an external smash product functor that assigns a spectrum $X\barwedge Y$
over $A\times B$ to $X$ and $Y$ and an external function spectrum functor
that assigns a spectrum $\bar{F}(Y,Z)$ over $A$ to $Y$ and $Z$, and 
there is a natural isomorphism
\[G\sS_{A\times B}(X\barwedge Y, Z)\iso G\sS_A(X, \bar{F}(Y,Z)).\]
The internal smash products are determined from the external ones via
$$ X\sma_B Y \iso \DE^*(X\barwedge Y) \qquad\text{and}\qquad
F_B(X,Y) \iso \bar{F}(X,\DE_*Y),$$
where $X$ and $Y$ are spectra over $B$ and $\DE\colon B\rtarr B\times B$ is the 
diagonal map.
\end{prop} \begin{proof}
It is not hard to start from \myref{extsmash1} and construct
these functors directly. We instead follow \myref{exin} and observe that 
the spectrum level external functors can and, up to isomorphism, must be 
defined in terms of the internal functors as
$$ X\barwedge Y \iso \pi_A^*X \sma_{A\times B} \pi_B^*Y
\qquad\text{and}\qquad \bar{F}(Y,Z) \iso \pi_{A\, *}F_{A\times B}(\pi_{B}^{*}Y,Z),$$
where $\pi_A\colon A\times B\rtarr A$ and $\pi_B\colon A\times B\rtarr B$ are the
projections. The displayed adjunction is immediate from the adjunctions 
$(\pi^*_A,\pi_{A\, *})$, $(\pi^*_B,\pi_{B\, *})$, and 
$(\sma_{A\times B}, F_{A\times B})$. The second statement follows formally,
as in \myref{internalize}.
\end{proof}

\begin{prop}\mylabel{SISISI} For ex-spaces $K$ over $A$ and $L$
over $B$, there is a natural isomorphism
$$\SI^{\infty}_{A\times B} (K\barwedge L) 
\iso \SI^{\infty}_{A}K\barwedge \SI^{\infty}_{B}L. $$
\end{prop}\begin{proof}
This is most easily seen using adjunction and the Yoneda lemma.
Using external function objects, we see that 
$\bar{F}(\SI^{\infty}_BL,Z)\iso \bar{F}(L,Z)$ for $Z\in G\sS_{A\times B}$. 
This has zeroth ex-space $\bar{F}(L,Z(0))$ over $A$.
\end{proof}

\section{Change of groups and restriction to fibers}

We give the analogues for parametrized spectra of the results concerning 
change of groups and restriction to fibers that were given for parametrized
ex-spaces in \S2.3. We shall say more about change of groups in Chapter 14.  
Fix an inclusion $\io\colon H\rtarr G$ of
a (closed) subgroup $H$ of $G$ and let $A$ be an $H$-space and $B$ be a 
$G$-space. We index $H$-spectra over $A$ on the collection $\io^*{\sV}$ 
of $H$-representations $\io^*V$ with $V\in \sV$. As we discuss in \S\S14.2
and 14.3, when $\sV$ is the collection of all representations of $G$, we can change indexing to the collection of all representations of $H$ since our assumption that $G$ is compact ensures that every representation of $H$ 
is a direct 
summand of a representation $\io^*V$. We have an evident forgetful functor
\begin{equation}
\io^*\colon G\sS_{B} \rtarr H\sS_{\io^*B}.
\end{equation}
On the space level, we write $\io_!$ ambiguously for both the based and 
unbased induction functors $G_+\sma_H(-)$ and $G\times_H(-)$, and similarly 
for coinduction $\io_*$. Context should make clear which is intended. 
Applying the unbased versions to retracts, we defined induction and 
coinduction functors $\io_!$ and $\io_*$ on ex-spaces in \myref{changes0}. 
These functors extend to the spectrum level.  Recall that $S_{G,B}$ 
denotes the $G$-sphere spectrum over $B$.

\begin{prop}\mylabel{eyeeye} Levelwise application of $\io_!$ and 
$\io_*$ gives functors
\[\io_!\colon H\sS_A\rtarr G\sS_{\io_! A} \qquad\text{and}\qquad
\io_*\colon H\sS_A\rtarr G\sS_{\io_* A}.\]
\end{prop}\begin{proof}
We must show that the structure $H$-maps
$\si\colon X\barwedge S_{H,A}\rtarr X\com \oplus$ 
of an $H$-spectrum $X$ over $A$ induce structure $G$-maps 
for the $\sI_G$-spaces $\io_!X$ and $\io_*X$.
It is clear that $\io_!(X\com \oplus)\iso \io_!X\com \oplus$
and $\io_*(X\com \oplus)\iso \io_*X\com \oplus$.
Using (\ref{GH3}), we see that $S_{G,\io_!A} \iso \io_!S_{H,A}$. 
Since the functor $\io_!$ on the ex-space level is symmetric monoidal by
\myref{ishriek}, its levelwise $\sI_G$-space analogue commutes
up to isomorphism with the external smash product $\barwedge$. Thus
$\si$ induces a structure $G$-map 
$$\io_!X \barwedge_{\io_!A} S_{G,\io_!A}
\iso \io_!(X\barwedge_A S_{H,A})\rtarr \io_!(X\com \oplus)\iso \io_!X\com \oplus.$$
For $\iota_*$, let $\mu:\iota^*\iota_*\rtarr \text{Id}$ be the counit 
of the space
level adjunction $(\iota_*,\iota^*)$ (see (\ref{GH1})). For
an $H$-space $A$, $\mu$ is the $H$-map $\text{Map}_H(G,A)\rtarr A$ 
given by evaluation at the identity element of $G$. Applied to an
ex-space $K$ over $A$, thought of as a retract,
$\mu$ gives a map $\iota^*\iota_*K\rtarr K$ of total spaces over and 
under the map $\mu:\iota^*\iota_*A\rtarr A$ of base spaces in the category 
of retracts of \S2.5. We can apply this to $X$ levelwise. 
We also have the projection 
$\text{pr}:\mu^*S_{H,A}\rtarr S_{H,A}$ over $\mu$. Together, these maps give 
$$
\xymatrix{
\iota^*(\iota_*X\barwedge_{\iota_*A}S_{G,\iota_*A})
\cong \iota^*\iota_*X\barwedge_{\iota^*\iota_*A} \mu^*S_{H,A} 
\ar[r]^-{\mu\barwedge \text{pr}} & X\barwedge_A S_{H,A}.\\}
$$
For the isomorphism, we have used the facts that $\iota^*$ is strong monoidal and 
that $\iota^*S_{G,\iota_*A}\cong S_{H,\iota^*\iota_*A}\cong \mu^*S_{H,A}$. 
The adjoint of the composite of this map with the structure map
$\sigma:X\barwedge_A S_{H,A}\rtarr X\circ \oplus$ gives the required
structure map 
$\iota_*X\barwedge_{\iota_*A} S_{G,\iota_*A}\rtarr \iota_*X\circ\oplus$.
\end{proof}

As on the ex-space level, the categories $H\sS_A$ and $G\sS_{G\times_H A}$ 
can be used interchangeably. The following result is immediate from \myref{ishriek}.

\begin{prop}\mylabel{changes}
Let $\nu\colon A\rtarr \io^*\io_!A$ be the natural inclusion of $H$-spaces.  
Then $\io_!\colon H\sS_A\rtarr G\sS_{\io_{!}A}$ is a closed symmetric monoidal equivalence of categories with inverse the composite 
${\nu}^*\com \io^*\colon G\sS_{\io_!A}\rtarr H\sS_{\io^*\io_!A}\rtarr H\sS_A$.
\end{prop}

In particular, if $A = *$ then $\nu$ maps $*$ to the identity coset 
$eH\in G/H$ and we see that $H\sS$ and $G\sS_{G/H}$ can be used 
interchangeably. Arguing as in \myref{homog}, we could more easily
prove this directly.

\begin{cor}\mylabel{changestoo}
The category $H\sS$ is equivalent as a closed symmetric mon\-oid\-al 
category to $G\sS_{G/H}$. Under this equivalence,
$$\io^*\iso r^*, \qquad \io_!\iso r_{!},  \qquad\text{and}\qquad \io_*\iso r_*,$$
where $r\colon G/H \rtarr *$.
\end{cor}

Looking at the fiber $X_b(V) = X(V)_b$ over $b$ of a $G$-spectrum $X$ over $B$, we see a $G_b$-spectrum $X_b$ of the sort that has been studied in \cite{MM}, 
where $G_b$ is the isotropy group of $b$.  Our homotopical analysis of
parametrized $G$-spectra will be based on the idea of applying the results 
of \cite{MM} fiberwise. By the previous result, we can think of this fiber
as a $G$-spectrum over $G/G_b$. The following spectrum level analogues of 
\myref{Johann0} and \myref{Johann1} analyze the relationships among passage to 
fibers, base change, and change of groups. 

\begin{exmp}\mylabel{Johann}
For $b\in B$, we write $b\colon *\rtarr B$ for the $G_b$-map that sends $*$ to $b$ and $\tilde{b}\colon G/G_b \rtarr B$ for the induced inclusion of orbits. Under the equivalence $G\sS_{G/G_b}\iso G_b\sS$, $\tilde{b}^*$ may be interpreted as the fiber functor $G\sS_B\rtarr G_b\sS$ that sends $Y$ to $Y_b$.  Its left and right adjoints $\tilde{b}_{!}$ and $\tilde{b}_*$ may be interpreted as the functors that send a $G_b$-spectrum $X$ to the $G$-spectra $X^b$ and $^bX$ over $B$ obtained by levelwise application of the corresponding ex-space level adjoints of \myref{Fibad} and \myref{Johann0}. With these notations, the isomorphisms of \myref{Wirth} specialize to the following natural isomorphisms, where $Y$ and $Z$ are in $G\sS_B$ and $X$ is in $G_b\sS$. 
\begin{gather*}
(Y\sma_B Z)_b\iso Y_b\sma Z_b,\\[1ex]
F_B(Y,\, ^bX) \iso\, {^{b}}F(Y_b,X),\\[1ex]
F_B(Y,Z)_b\iso F(Y_b,Z_b),\\[1ex]
(Y_b\sma X)^b\iso Y\sma_B X^b,\\[1ex]
F_B(X^b,Y)\iso \, {^{b}}F(X,Y_b).
\end{gather*}
\end{exmp}

\begin{exmp}\mylabel{Johann2} 
Let $f\colon A\rtarr B$ be a $G$-map and let $i_b\colon A_b\rtarr B$
be the inclusion of the fiber over $b$, which is a $G_b$-map. As in
\myref{Johann1}, we have the compatible pullback squares
$$\xymatrix{
A_b \ar[r]^-{f_b} \ar[d]_{i_b} & \{b\} \ar[d]^{b} \\
A \ar[r]_{f} & B\\}
\qquad\qquad
\xymatrix{
G\times_{G_b} A_b \ar[r]^-{G\times_{G_b} f_b} \ar[d]_{\tilde{\imath}_b} 
& G/G_b \ar[d]^{\tilde{b}} \\
A \ar[r]_-{f} & B.\\}$$
Applying \myref{Mackey} to the right-hand square and 
interpreting the conclusion in terms of fibers, we obtain canonical 
isomorphisms of $G_b$-spectra
$$(f_!X)_b \iso {f_{b}}_!i_b^*X \qquad\text{and}\qquad (f_*X)_b \iso {f_{b}}_*i_b^*X,$$
where $X$ is a $G$-spectrum over $A$, regarded on the right-hand sides as a $G_b$-spectrum 
over $A$ by pullback along $\io\colon G_b\rtarr G$.
\end{exmp}

\section{Some problems concerning non-compact Lie groups}

In equivariant stable homotopy theory, the key idea 
is that the one-point compactification of a representation $V$ of 
dimension $n$ is a $G$-sphere and that smashing with that sphere 
should be a self-equivalence of the equivariant stable homotopy category.
That is, the idea is to invert $G$-spheres in just the way that we
invert spheres when constructing the nonequivariant stable homotopy
category. For compact Lie groups of equivariance, the philosophy
and its implementation and applications are well understood. When
we invert representation spheres, we invert other homotopy spheres
as well, and the relevant Picard group is analyzed in \cite{FLM}.

For non-compact Lie groups, the present work seems to be the first
attempt to consider foundations for equivariant stable 
homotopy theory.  The philosophy is less clear, and its technical 
implementation is problematic. The need for such a theory
is evident, however.  The focus on finite dimensional representations
is intrinsic to the philosophy but fails to come to grips with basic
features of the representation theory of non-compact Lie groups. A theory 
based on finite dimensional representations
should still have its uses, but there are real difficulties to obtaining
even that much.   In particular,
a focus on spheres associated to linear representations, rather than
on less highly structured homotopy spheres, may be misplaced.

A non-compact semi-simple Lie group will generally have no non-trivial
finite dimensional unitary or orthogonal representations, hence our
theory of ``orthogonal'' $G$-spectra is clearly too restrictive.  This is 
easily remedied.
The use of linear isometries in the definition of orthogonal
spectra is a choice dictated more by the
history than by the mathematics.  In the alternative 
approach to equivariant stable homotopy theory based on 
Lewis-May spectra and EKMM \cite{EKMM, LMS, MM}, use of orthogonal 
complements is certainly convenient and perhaps essential.  However,
the diagram orthogonal spectra of \cite{MM, MMSS} could just as well
have been developed in terms of diagram ``general linear spectra''.
In the few places where complements are used, they can by avoided.
For consistency with the previous literature, we have chosen to give
our exposition in the compact case using the word ``orthogonal'' and 
the language from the cited references, but for general
Lie groups of equivariance, we should eliminate all 
considerations of isometries. 

More precisely, for the complete case, we redefine $\sI$ by taking $\sV$ to be the collection of all finite dimensional representations $V$ of $G$.  More generally, we can index on any subcollection that contains the trivial representation and is closed under finite direct sums. Since we are only interested in a skeleton of $\sI$, we may as well restrict to orthogonal representations in $\sV$ when $G$ is compact.
We replace linear isometries by linear isomorphims when defining the $G$-spaces $\sI(V,W)$. Thus we replace orthogonal groups by general linear groups. Otherwise, the formal definitional framework developed in this chapter (or, 
in the nonparametrized case, \cite[II]{MM}) goes through verbatim for general
topological groups $G$.  

However, we emphasize the formality.  When considering change of groups,
for example, the significance changes drastically. As noted at the start 
of the previous section, for an inclusion $\io\colon H\rtarr G$ of
a (closed) subgroup $H$ of $G$, we index $H$-spectra on the collection $\io^*{\sV}$ of $H$-representations $\io^*V$ with $V\in \sV$. We also
pointed out the relevance of the compact case of the following result.

\begin{prop} If $G$ is either a compact Lie group or a matrix group
and $W$ is a representation of a subgroup $H$, then there is a representation $V$ 
of $G$ and an embedding of $W$ as a subrepresentation of $\io^*V$.
\end{prop}

This is clear in the compact case and is given by \cite[3.1]{Palais} 
for matrix groups.  However, the following striking counterexample,
which we learned from Victor Ginzburg, shows just how badly this 
basic result fails in general.

\begin{ouch0}[Ginzburg] Let $\bH$ be the Heisenberg group of 
$3\times 3$ matrices
\[ \left( \begin{array}{ccc}
1 & a & c \\
0 & 1 & b \\
0 & 0 & 1 \\
\end{array} \right) \]
where $a$, $b$, and $c$ are real numbers.  Embed $\bR$ in $\bH$ 
as the subgroup of matrices with $a=b=0$.  Embed $\bZ$ in $\bR$
as usual.  Then $\bR$ is a central subgroup of $\bH$.  Define
$G = \bH/\bZ$. Then $T = \bR/\bZ$ is a circle subgroup of $G$.
Moreover, $T$ is the center of $G$ and coincides with the commutator
subgroup $[G,G]$.  Let $V$ be any finite dimensional (complex linear) representation
of $G$.  Since $T$ is compact, the action of $T$ on $V$ is semisimple,
and since $T$ is central, any weight space of $T$ is a $G$-submodule.
Therefore $V$ is a direct sum of $G$-submodules $V_i$ such that $T$ 
acts on each $V_i$ by scalar matrices.  Since $T = [G,G]$,
this scalar action of $T$ on $V_i$ is trivial: the determinant
of $g$ is $1$ for any $g\in [G,G]$. Therefore no nontrivial 
$1$-dimensional character of $T$ can embed in $V$. Reinterpreting
in terms of real representations, as we may, we conclude that, for 
$\io\colon T\rtarr G$, $\io^*\sV$ is the trivial $T$-universe.
\end{ouch0}

For a compact Lie group $G$ and inclusion $\io\colon H\subset G$, 
$\io^*X$ is a dualizable $H$-spectrum if $X$ is a dualizable $G$-spectrum,
and an $H$-spectrum indexed on the trivial $H$-universe is dualizable if 
and only if it is a retract of a finite $H$-CW spectrum built up from 
trivial orbits. We conclude that duality theory (in the nonparametrized
context) cannot work as one would wish in the context of the previous example.

Looking ahead, much of the theory of the following three chapters also
works formally in the context of non-compact Lie groups.  However, there is at
least one serious technical difficulty.  Our theory is
based on the use of one-point compactifications $S^V$.  If $V$ is a linear
representation of a non-compact Lie group $G$, there is no reason to
think that $G$ acts smoothly and properly on $S^V$, even if the isotropy
groups of $V$ are compact.  In fact, if Illman's \myref{Illman} were 
to apply, then $S^V$ would be a $G$-cell complex, hence it would be built 
up from non-compact orbits $G/H$ given by compact subgroups $H$.  However, as closed
subsets of $S^V$, the closed cells would have to be compact. That is, the
putative $G$-CW structure would contradict the
compactness of $S^V$. Said another way, we see no reason to believe
that the $S^V$ are $q$-cofibrant $G$-spaces.  Therefore, the functors 
$(-)\sma S^V$ need not be Quillen left adjoints and the functors $\OM^V$ 
and $\OM^V_B$ need not preserve fibrant objects in the relevant model
structures.  Compare, for example, \myref{QuillK} and the derivation 
of the long exact sequences (\ref{les2}) and (\ref{les3}) below. What
seems to be needed, for a start, is something like a model structure on
$G$-spaces such that $X$ is cofibrant if $\io^*X$ is $H$-cofibrant for
all inclusions $\io\colon H\to G$ of compact subgroups.

 \chapter{Model structures for parametrized $G$-spectra}

\section*{Introduction}

We define and study two model structures on the category $G\sS_B$ of (orthogonal) $G$-spectra over $B$. We emphasize that, except for the theory of smash products, everything in this chapter applies equally well to the category $G\sP_B$ of $G$-prespectra over $B$.  That fact will become 
important in the next chapter.

We start in \S12.1 by defining a ``level model structure'' on $G\sS_B$, based on the $qf$-model structure on $G\sK_B$.  In \S12.2, we record analogues for this model structure of the results on external smash product and base change functors that were given for $G\sK_B$ in \S7.2. The level model structure serves as a stepping stone to the stable model structure, which we define in \S12.3.  It has the same cofibrations as the level model structure, and we therefore call these ``$s$-cofibrations''. An essential point in our approach is a fiberwise definition of the homotopy groups of a parametrized $G$-spectrum that throws much of our work onto the theory of nonparametrized orthogonal $G$-spectra developed by Mandell and the first author in \cite{MM}.   We define homotopy groups using the level $qf$-fibrant replacement functor provided by the level model structure, and we define stable equivalences to be the $\pi_*$-isomorphisms. It is essential to think in terms of fibers and not total spaces since the total spaces of a parametrized spectrum do not assemble into a spectrum.  We show in \S12.4 that the $\pi_*$-isomorphisms give a well-grounded subcategory of weak equivalences, and we complete the proofs of the model axioms in \S12.5.  We return to the context of \S12.2 in \S12.6, where we prove 
that various Quillen adjoint pairs in the level model structures are also
Quillen adjoint pairs in the stable model structures.

The basic conclusion is that $G\sS_B$ is a well-grounded model category under the stable structure.  Although not very noticeable on the surface, essential use is made of the $qf$-model structure on $G\sK_B$ throughout this chapter. It is possible to obtain a level model structure on $G\sS_B$ from the $q$-model structure on $G\sK_B$, as we explain in \myref{qnogood}. However this model
structure is not well-grounded and therefore does not provide the necessary tools to work out the technical details of \S12.4. The results there are crucial to prove that the relative cell complexes over $B$ defined in terms of the appropriate generating acyclic $s$-cofibrations are acyclic.\footnote{In \cite[3.4]{Hu}, such acyclicity of relative cell complexes is assumed without proof.}  It was our fruitless attempt to obtain a stable model structure starting from the level $q$-model structure that led us to the construction of the $qf$-model structure on $G\sK_B$ and to the notion of a well-grounded model category.

When there are no issues of equivariance, we generally abbreviate 
$G$-spectrum over $B$, ex-$G$-space, and $G$-space to spectrum 
over $B$, ex-space, and space; $G$ is a compact Lie group throughout.

\section{The level model structure on $G\sS_B$}

After changing the base space from $*$ to $B$, the level model structure works in much the same way as in the nonparametrized case of \cite{MM}.

\begin{defn}\mylabel{never}
Let $f\colon X\rtarr Y$ be a map of spectra over $B$.  With one exception, for any type of ex-space and any type of map of ex-spaces, we say that $X$ or $f$ is a \emph{level type of spectrum over $B$}\index{spectrum!level type of --} or a \emph{level type of map of spectra over $B$}\index{fibration!level type of --} if each $X(V)$ or $f(V)\colon X(V) \rtarr Y(V)$ is that type of ex-space or that type of map. Thus, for example, we have level $h$, level $f$ and level $fp$-fibrations, cofibrations and equivalences from \S5.1 together with the corresponding fibrant and cofibrant objects. We have level $q$-equivalences and level $q$ and $qf$-fibrations from \S7.1 and we have level ex-fibrations and level ex-quasifibrations from \S8.1 and \S8.5. The exceptions concern cofibrations and cofibrant objects. We shall \emph{never} be interested in ``level $q$-cofibrations'' or ``level $qf$-cofibrations'', nor in ``level $q$-cofibrant'' or ``level $qf$-cofibrant'' objects, since these do not correspond
to cofibrations and cofibrant objects in the model structures that we consider. Instead we have the following definitions.
\begin{enumerate}[(i)]
\item $f$ is an \emph{$s$-cofibration}\index{cofibration!s-@$s$- --} if it satisfies the LLP with respect to the level acyclic $qf$-fibrations.
\item $f$ is a \emph{level acyclic $s$-cofibration} if it is both a level $q$-equivalence and an $s$-cofibration.
\end{enumerate}
To reiterate, in the phrase ``level acyclic $qf$-fibration'', the adjective ``level'' applies to ``acyclic $qf$-fibration'', but in the phrase ``level 
acyclic $s$-cofibration'' it  applies \emph{only} to ``acyclic''; the cofibrations are not defined levelwise.
\end{defn}

\begin{defn}\mylabel{sillybilly}
A spectrum $X$ over $B$ is \emph{well-sectioned} if it is level well-sectioned,\index{spectrum!well-sectioned}\index{well-sectioned!spectrum} 
so that each ex-space $X(V)$ is $\bar{f}$-cofibrant. It is \emph{well-grounded}\index{well-grounded!spectrum}\index{spectrum!well-grounded} 
if it is level well-grounded, so that each $X(V)$ is well-sectioned 
and compactly generated.
\end{defn}

The discussion of \S\ref{sec:towardh} applies to the category $G\sS_B$ of $G$-spectra over $B$ with homotopies defined in terms of the cylinders $X\sma_B I_+$. In particular, we have the notion of a Hurewicz cofibration in $G\sS_B$, abbreviated $cyl$-cofibration,\index{cofibration!cyl-@$cyl$- --} defined in terms of these cylinders, and we also have the notion of strong Hurewicz cofibration, abbreviated $\overline{cyl}$-cofibration.

\begin{lem}\mylabel{levelcofs}
A $cyl$-cofibration of spectra over $B$ is a level $fp$-cofibration and a $cyl$-fibration of spectra over $B$ is a level $fp$-fibration. A $cyl$-cofibration between well-sectioned spectra over $B$ is a level $f$-cofibration and therefore both a level $h$-cofibration and a level $fp$-cofibration.
\end{lem}

\begin{proof}
By the mapping cylinder retraction criterion of Hurewicz cofibrations, a $cyl$-cofibration of spectra over $B$ is a level $fp$-cofibration. The statement about fibrations follows similarly from the path lifting function characterization of Hurewicz fibrations. An $fp$-cofibration between well-sectioned ex-spaces is an $f$-cofibration by \myref{reverse2}, and all $f$-cofibrations are $h$-cofibrations.
\end{proof}

Recall the notions of a ground structure and of a well-grounded 
subcategory of weak equivalences from Definitions \ref{back},
\ref{moreback}, and \ref{hproper}.

\begin{prop}\mylabel{levelwellgr}
The well-grounded spectra over $B$ give $G\sS_B$ a ground structure whose
ground cofibrations, or $g$-cofibrations, are the level $h$-cofibrations. The level $q$-equivalences specify a well-grounded subcategory of weak equivalences
with respect to this ground structure. In the gluing and colimit lemmas, one need only assume that the relevant maps are level $h$-cofibrations, not necessarily also $cyl$-cofibrations.
\end{prop}

\begin{proof}
That we have a ground structure follows levelwise from the ground structure on ex-spaces in \myref{exback}. That the level $q$-equivalences 
are well-grounded follows levelwise from \myref{exwellgr}.
\end{proof}

We construct the level model structure on $G\sS_{B}$ from the $qf$-model structure on $G\sK_B$ specified in \myref{Theqf}, but all results apply
verbatim starting from the $qf(\sC)$-model structure for any closed
generating set $\sC$ (as defined in \myref{IJBG2}). We shall need the extra generality for the reasons discussed in Chapter 7. Recall that $I^f_B$ and $J^f_B$ denote the sets of generating $qf$-cofibrations and generating acyclic $qf$-cofibrations in $G\sK_B$. We use the shift desuspension functors $F_V$ of \myref{FVs} to obtain corresponding sets on the spectrum level. We need the following observations.

\begin{lem}\mylabel{goodFV} 
The functor $F_V$ enjoys the following properties.
\begin{enumerate}[(i)]
\item If $K$ is a well-grounded ex-space over $B$, then $F_VK$ is well-grounded. If $K$ is an ex-fibration, then $F_VK$ is a level ex-fibration.
\item If $i\colon K\rtarr L$ is an $h$-equivalence between well-grounded ex-spaces over $B$, then $F_Vi$ is a level $h$-equivalence.
\item If $i\colon K\rtarr L$ is an $fp$-cofibration, then $F_Vi$ is a $cyl$-cofibration and therefore a level $fp$-cofibration. If, further,
$K$ and $L$ are well-sectioned, then $F_Vi$ is a level $f$-cofibration 
and therefore a level $h$-cofibration.
\item If $i\colon K\rtarr L$ is an $\overline{fp}$-cofibration, then 
$F_Vi$ is a $\overline{cyl}$-cofibration.
\item If $i\colon K\rtarr L$ is an $\overline{f}$-cofibration
between well-grounded ex-spaces over $B$,
then $F_Vi$ is a $\overline{cyl}$-cofibration which is a level 
$\overline{f}$-cofibration and therefore both a level 
$\overline{fp}$-cofibration and a level $\overline{h}$-cofibration.
\end{enumerate}
\end{lem}

\begin{proof}
By \myref{FVs}, $(F_V K)(W) = \sJ_G(V,W)_B\wedge_B K$, and the $G$-space $\sJ_G(V,W)$ is well-based. Now (i) holds by \myref{HursmaK} and (ii) holds 
by \myref{savior}. Since $F_V$ is left adjoint to the evaluation functor $\text{Ev}_V$ and since $cyl$-fibrations are level $fp$-fibrations, (iv) 
and the first statement of (iii) follow from the definitions by adjunction. 
The second statement of (iii) follows from \myref{reverse2}. The first 
half of (v) follows from (iv) since $\overline{f}$-cofibrations are $\overline{fp}$-cofibrations, and the second half follows from (iii) since $F_Vi$ is a level $f$-cofibration between well-grounded spectra and therefore a level $\overline{f}$-cofibration by \myref{ffpmodel}(ii).
\end{proof}

\begin{defn}\mylabel{FBJB}
Define $FI^f_B$\noteindex{FIBf@$FI^f_B$} to be the set of maps $F_V i$ with $V$ in a skeleton $sk \sI_G$ of $\sI_G$ and $i$ in $I^f_B$. Define $FJ^f_B$\noteindex{FJBf@$FJ^f_B$} to be the set of maps $F_Vj$ with 
$V$ in $sk \sI_G$ and $j$ in $J^f_B$.
\end{defn}

Recall the notion of a well-grounded model structure from \myref{wellmodel}.  
Among other properties, such model structures are compactly generated, proper,
and $G$-topological.

\begin{thm}\mylabel{levelqf}\index{model structure!level -- on parametrized spectra}\index{level model structure}
The category $G\sS_B$ is a well-grounded model category with respect to the level $q$-equivalences, the level $qf$-fibrations and the $s$-cofibrations. The sets $FI^f_B$ and $FJ^f_B$ give the generating $s$-co\-fi\-bra\-tions and the generating level acyclic $s$-co\-fi\-bra\-tions. All $s$-cofibrations are level $\overline{f}$-cofibrations, hence level $\overline{fp}$ and level $\overline{h}$-cofibrations, and all $s$-cofibrant spectra over $B$ are well-grounded.
\end{thm}

\begin{proof}
By \myref{goodFV}, the maps in $FI^f_B$ and $FJ^f_B$ are
$\overline{cyl}$-co\-fi\-bra\-tions between well-grounded 
objects and $\overline{f}$-co\-fi\-bra\-tions. Moreover, the 
maps in $FJ^f_B$ are level acyclic. Therefore, to prove the 
model axioms, we need only verify the compatibility condition (ii) in \myref{Newcompgen}.  Adjunction arguments show that a map is a level $qf$-fibration if and only if it has the RLP with respect to $FJ^f_B$ and that it is a level acyclic $q$-fibration if and only if it has the RLP with respect to $FI^f_B$. This implies that the classes of $s$-cofibrations and of $FI^f_B$-cofibrations (in the sense of \myref{cofhyp}(iii)) coincide. Therefore, if a map has the RLP with respect to $FI^f_B$, then it is a level acyclic $qf$-fibration. The required compatibility condition now follows from its analogue for $G\sK_B$. Condition (iv) in \myref{Newcompgen} holds by its ex-space level analogue and the fact that $(F_V K)\sma_B T\iso F_V(K\sma_B T)$ for an ex-space $K$ over $B$ and a based space $T$. Right properness follows directly from the space level analogue.  
\end{proof}

\begin{rem}\mylabel{qnogood}
Just as in \myref{FBJB}, we can also define sets $FI_B$ and $FJ_B$ based on 
the generating sets $I_B$ and $J_B$ for the $q$-model structure on $G\sK_B$. 
We can then use \myref{compgen} to prove the analogue of \myref{levelqf} 
stating that $G\sS_B$ is a cofibrantly generated model category under the 
level $q$-model structure. Since the compatibility condition holds by the same proof as for the level $qf$-model structure, we need only verify the acyclicity condition to show this. 

For a generating acyclic $q$-cofibration $j\in J_{B}$, we have 
$F_Vj=V^*\sma_B j$, where $V^*(W)=\sJ_{G,B}(V,W)$. This map is a level 
$h$-equivalence by \myref{goodFV}(ii). Although $j$ is an $h$-cofibration, 
it is not immediate that $F_Vj$ is a level $h$-cofibration. (This holds for $j\in J_B^f$ by \myref{goodFV}(iii), since $j$ is then an $fp$-cofibration). Indeed, for general ex-spaces $K$ and $h$-cofibrations $f$, $K\sma_B f$ need 
not be an $h$-cofibration. However, since $\sJ_{G,B}(V,W) = \sJ_G(V,W)_B$, we see directly that $F_Vj$ is indeed a level $h$-cofibration. By inspection of 
the definition of wedges over $B$ in terms of pushouts, the gluing lemma in
$\sK$ then applies to show that wedges over $B$ of maps in $FJ_B$ are level 
acyclic $h$-cofibrations. Since pushouts and colimits in $\sS_B$ are 
constructed levelwise on total spaces, it follows that relative $FJ_B$ 
complexes are acyclic $h$-cofibrations since the $q$-model structure on 
$\sK$ is well-grounded. 
\end{rem}

\begin{rem}\mylabel{plus1} As in the nonparametrized case \cite{MM}, ``positive'' model structures would be needed to obtain a comparison
with the as yet undeveloped alternative approach to parametrized stable 
homotopy theory based on \cite{EKMM, LMS}. Such model structures can be
defined as in \cite[p.\, 44]{MM}, starting from the subsets $(FI^f_B)^+$ 
and $(FJ^f_B)^+$ that are obtained by restricting to those $V$ such that $V^G\neq 0$. One then defines the \emph{positive level} versions of all 
of the types of maps specified in \myref{never} by restricting to those 
levels $V$ such that $V^G\neq 0$. The positive level analogue of 
\myref{levelqf} holds, where the positive $s$-cofibrations are the 
$s$-cofibrations that are isomorphisms at all levels $V$ such that 
$V^G=0$; compare \cite[III.2.10]{MM}.  However, we shall make no use of the positive model structure in this paper, and we will make little further reference to it.
\end{rem}

The same proof as in \cite[I.2.10, II.4.10, III.2.12]{MM} gives the following result. 

\begin{thm}\mylabel{QPU}
The forgetful functor $\bU$ from spectra over $B$ to prespectra over $B$ has a left adjoint $\bP$ such that $(\bP,\bU)$ is a Quillen equivalence.
\end{thm}

\section{Some Quillen adjoint pairs relating level model structures}

This section gives the analogues for the level model structure of some of the ex-space level results in \S\S7.2-7.4.  These results are also analogues of 
results in \cite[III.\S2]{MM}, which in turn have non-equivariant precursors 
in \cite[\S6]{MMSS}.  They admit essentially the same proofs as in Chapter 7 
or in the cited references.  The level $qf$-model structure is understood throughout. More precisely, where a $qf(\sC)$-model structure was used in
Chapter 7, we must use the corresponding level $qf(\sC)$-model structure here.
Since we want our model structures to be $G$-topological, we only use generating
sets $\sC$ that are closed under finite products.

Our first observation is immediate from the fact that equivalences and
fibrations are defined levelwise, the next follows directly from its
ex-space analogue Proposition \ref{smaB}, and the third and fourth
are proven in the same way as their ex-space analogues \ref{Boxcof20} and
\myref{smaAB}. All apply to the level $qf(\sC)$-model structures for any 
choice of $\sC$.

\begin{prop}
The pair of adjoint functors $(F_V,\text{Ev}_V)$ between $G\sK_B$ and $G\sS_B$ is a Quillen adjoint pair.
\end{prop}

\begin{prop}\mylabel{QuillK}
For a based $G$-CW complex $T$, $((-)\sma_B T, F_B(T,-))$ is a Quillen adjoint pair of endofunctors of $G\sS_B$.
\end{prop} 

\begin{prop}\mylabel{Boxcof2}
If $i\colon  X\rtarr Y$ and $j\colon  W\rtarr Z$ are $s$-cofibrations 
of spectra over base spaces $A$ and $B$, then 
\[i\Box j\colon  (Y\barwedge  W)\cup_{X\barwedge  W}(X\barwedge Z)\rtarr Y\barwedge Z\]
is an $s$-cofibration over $A\times B$ which is level 
acyclic if either $i$ or $j$ is acyclic. 
\end{prop}

As in \S7.2, we cannot expect this result to hold for internal smash 
products over $B$. The case $A=*$, which relates spectra to spectra 
over $B$, is particularly important.  As we explain in \S14.1, it 
leads to a fully satisfactory theory of \emph{parametrized} module 
spectra over \emph{nonparametrized} ring spectra.

\begin{cor}\mylabel{ext} 
If $Y$ is $s$-cofibrant over $B$, then the functor $(-)\barwedge Y$ from $G\sS_A$ to $G\sS_{A\times B}$ is a Quillen left adjoint with Quillen right adjoint $\bar{F}(Y,-)$. 
\end{cor}

Again the next result is a direct consequence of its 
ex-space analogue \myref{Qad10} and applies with any 
choice of $\sC$.

\begin{prop}\mylabel{Qad1}
Let $f\colon A\rtarr B$ be a $G$-map. Then $(f_{!},f^*)$ is a Quillen 
adjoint pair. The functor $f_!$ preserves level $q$-equivalences
between well-sectioned $G$-spectra over $B$. If $f$ is a $qf$-fibration, 
then $f^*$ preserves all level $q$-equivalences.
\end{prop}

\begin{prop}\mylabel{ffequiv}
If $f\colon A\rtarr B$ is a $q$-equivalence, then $(f_{!},f^*)$ is a 
Quillen equivalence.
\end{prop}
\begin{proof} We mimic the proof of \myref{ffequiv0}, but with $X$
and $Y$ taken to be an $s$-cofibrant $G$-spectrum over $A$ and 
a level $qf$-fibrant $G$-spectrum over $B$.  It is clear that 
$f^*Y\rtarr Y$ is a level $q$-equivalence since $A\rtarr B$ is 
a $q$-equivalence.  Since $X$ is $s$-cofibrant, $*_A\rtarr X$ is a 
level $h$-cofibration. Note that it is essential for this statement 
that we start from the $qf$ and not the $q$-model structure on ex-spaces. 
Since pushouts along level $h$-cofibrations preserve level $q$-equivalences, $X\rtarr f_!X$ is a level $q$-equivalence. The conclusion follows as in
\myref{ffequiv0}.
\end{proof}

\begin{prop}\mylabel{Qad2}
Let $f\colon A\rtarr B$ be a $G$-bundle whose fibers $A_b$ are $G_b$-CW complexes. Then $f^*$ preserves level $q$-equivalences and 
$s$-cofibrations. Therefore $(f^*,f_*)$ is a Quillen adjoint pair.
\end{prop}

\begin{proof} Here we must use a generating set $\sC(f)$ as specified 
in \myref{Qad202}. The proof that $f^*$ preserves $s$-cofibrations 
reduces to showing that the maps $f^*F_Vi \cong F_Vf^*i$ are 
$s$-cofibrations for generating $s$-cofibrations $i$. Since $F_V$ is a 
Quillen left adjoint it takes $qf$-cofibrations to $s$-cofibrations, 
so we are reduced to the ex-space level, where $f^*i$ is shown to be a $qf$-cofibration in \myref{Qad202}.
\end{proof}

Now consider the change of groups functors of \S11.5. The following result
shows that the equivalence of \myref{changes} descends to homotopy categories.
It is proven by levelwise application of its ex-space analogue
\myref{Lishriek}, together with change of universe considerations 
that are deferred until \S14.2 and \S14.3. 

\begin{prop}\mylabel{Lchanges}
Let $\iota\colon H\rtarr G$ be the inclusion of a subgroup. The pair of functors $(\iota_!,\nu^*\iota^*)$ relating $H\sS_A$ and $G\sS_{\iota_!A}$ give a Quillen equivalence. If $A$ is completely regular, then $\iota_!$ is also a Quillen right adjoint.
\end{prop}

For a point $b$ in $B$, we combine the special case $\tilde{b}\colon G/G_b\rtarr B$ of \myref{Qad1} with \myref{Lchanges}, where $\iota\colon G_b\rtarr G$ and $\nu\colon *\rtarr G/G_b$, to obtain the following analogue of 
\myref{FibadQ0}.  Recall from \myref{Johann} that the fiber functor 
$(-)_b\colon G\sS_B\rtarr G_b\sS$ is given by $\nu^*\iota^*\tilde{b}^*=b^*\iota^*$. Its left adjoint $(-)^b$ 
therefore agrees with $\tilde{b}_!\iota_!$.

\begin{prop}\mylabel{FibadQ}
For $b\in B$, the pair of functors $((-)^b,(-)_b)$ relating $G_b\sS_*$ and $G\sS_B$ is a Quillen adjoint pair.
\end{prop}

\section{The stable model structure on $G\sS_B$}

The essential point in the construction of the stable model structure is to define the appropriate (stable) homotopy groups. The weak equivalences will then be the maps of parametrized spectra that induce isomorphisms on all homotopy groups. We refer to them as the $\pi_*$-isomorphisms or $s$-equivalences, using
these terms interchangeably. There are several motivating observations for our definitions. We return the group $G$ to the notations for the moment.

First, a $G$-spectrum $X$ over $B$ is level $qf$-fibrant if and only if each
projection $X(V)\rtarr *_B(V)=B$ is a $qf$-fibration of ex-$G$-spaces. It is equivalent that each fixed point map $X(V)^H\rtarr B^H$ be a non-equivariant $qf$-fibration, and, by \myref{qfles}, we have resulting long exact sequences of homotopy groups 
\begin{equation}\label{les1}
\cdots \rtarr \pi_{q+1}^H(B)
\rtarr \pi^H_q(X_b(V)) 
\rtarr \pi^H_q(X(V))
\rtarr \pi_q^H(B)
\rtarr \cdots
\end{equation}
for each $b\in B^H$. Here, for a $G$-space $T$, $\pi_q^H(T)$ denotes $\pi_q(T^H)$.

Second, as we have already discussed in \S11.4, the fibers $X_b$ of a $G$-spectrum $X$ are $G_b$-spectra, and our guiding principle is to use these nonparametrized spectra to encode the homotopical information about our parametrized spectra. \myref{FibadQ} allows us to encode levelwise 
information in the level homotopy groups of fibers, and it is plausible that 
we can similarly encode the full structure of our parametrized $G$-spectrum 
$X$ in the spectrum level homotopy groups of the fiber $G_b$-spectra $X_b$.
However, we can only expect to do so when $X$ is level $qf$-fibrant and we have
the long exact sequences (\ref{les1}).  

Recall that the homotopy groups $\pi_q^H(Y)$ of a nonparametrized 
$G$-spectrum $Y$ are defined in \cite[III.3.2]{MM} as the colimits of the groups $\pi_q^H(\OM^VY(V))$, where the maps of the colimit system are induced in the
evident way by the adjoint structure maps $\tilde\sigma\colon  Y(V)\rtarr \Omega^{W-V}Y(W)$ of $Y$. The functor $\OM^V$ on based $G$-spaces preserves 
$q$-fibrations and the functor $\Omega^V_B=F_B(S^V,-)$ on $G$-spectra over $B$ preserves level $qf$-fibrations. Formally, these hold since $S^V$ is a $q$-cofibrant $G$-space and the relevant model structures are $G$-topological. 
This leads to two families
of long exact sequences relating the homotopy groups $\pi_q^H(\OM^VX_b(W)$ of
fibers to the homotopy groups of the base space $B$ and of the total spaces
$X(W)$.  First, if $X$ is a level $q$-fibrant $G$-spectrum over $B$, then, using basepoints determined by a point $b\in B^H$ for any $H\subset G_b$, the $q$-fibrations $\OM^VX(W) \rtarr \OM^VB$ of based $G$-spaces with fibers $\OM^VX_b(W)$ induce long exact sequences 
\begin{equation}\label{les2}
\cdots \rtarr \pi_{q+1}^H(\OM^V B)\rtarr \pi^H_q(\OM^VX_b(W)) 
\rtarr \pi^H_q(\OM^VX(W)) \rtarr \pi_q^H(\OM^V B)\rtarr \cdots.
\end{equation}
Second, if $X$ is level $qf$-fibrant, then the $qf$-fibrations $(\Omega^V_BX)(W)\rtarr *_B$ of ex-$G$-spectra over $B$ with 
fibers $\Omega^VX_b(W)$ induce long exact sequences
\begin{equation}\label{les3}
\cdots \rtarr \pi_{q+1}^H(B)
\rtarr \pi^H_q(\OM^VX_b(W)) 
\rtarr \pi^H_q((\OM^V_BX)(W))
\rtarr \pi_q^H(B)
\rtarr \cdots.
\end{equation}
The first allows us to relate the homotopy groups of the 
$X_b$ to the homotopy groups of the ordinary loops $\OM^VX(W)$ 
on total spaces.  The second
allows us to relate the homotopy groups of the $X_b$ to 
the homotopy groups of the parametrized loop ex-spaces $(\OM_B^VX)(W)$. 
It is the second that is most relevant to our work.

\begin{defn}\mylabel{htygps}
The \emph{homotopy groups}\index{homotopy groups!of parametrized spectra} of a level $qf$-fibrant $G$-spec\-trum over $B$, or of a level $qf$-fibrant $G$-prespectrum $X$, are all of the homotopy groups $\pi_q^H(X_b)$ of all of the fibers $X_b$, where $H\subset G_b$. The homotopy groups of a general $G$-spectrum, or $G$-prespectrum, $X$ over $B$ are the homotopy groups $\pi_q^H((RX)_b)$ of a level $qf$-fibrant approximation $RX$ to $X$. We still denote these homotopy groups by $\pi_q^H(X_b)$.
In either category, a map $f\colon X\rtarr Y$ is said to be a \emph{$\pi_*$-isomorphism}\index{equivalence!s-@$s$- --}\index{equivalence!p-iso@$\pi_*$-isomorphism}\index{p-iso@$\pi_*$-isomorphism} or, synonymously, an \emph{$s$-equi\-va\-lence}, if, after level $qf$-fibrant approximation, it induces an isomorphism on all homotopy groups.
\end{defn} 

There are also homotopy groups specified in terms of maps out of sphere spectra over $B$, but we choose to ignore them in setting up our model theoretic foundations. Our choice captures the intuitive idea that spectra over $B$ should be {\em parametrized spectra}\/:  the fiber spectra should carry all of the homotopy theoretical information. With this choice, a good deal of the work needed to set up the stable model structure reduces to work that has already been done in \cite{MM}. The following observation is a starting point that illustrates the pattern of proof. Now that we have seen how the equivariance
appears in the definition of homotopy groups, we revert to our custom of
generally deleting $G$ from the notations.

\begin{lem}\mylabel{levelpi}
A level $q$-equi\-va\-lence of spectra over $B$ is a $\pi_*$-iso\-mor\-phism. \end{lem}
\begin{proof} 
A level $qf$-fibrant approximation to the given level $q$-equi\-va\-lence is a level acyclic $qf$-fibration, and it induces a level $q$-equi\-va\-lence on fibers over points of $B$ by \myref{FibadQ}. This allows us to apply {\cite[III.3.3]{MM}}, which gives the same conclusion for nonparametrized spectra, one fiber at a time.
\end{proof}

To exploit our definition of homotopy groups, we need the following accompanying definition and proposition. 

\begin{defn}  
An \emph{$\Omega$-prespectrum over $B$}\index{prespectrum!Om@$\Omega$- --} is a level $qf$-fibrant prespectrum $X$ over $B$ such that  each of its adjoint structure maps $\tilde{\si}\colon  X(V)\rtarr \OM_B^{W-V}X(W)$ is a $q$-equi\-va\-lence of ex-spaces over $B$, that is, a $q$-equi\-va\-lence of total spaces. An (orthogonal) \emph{$\OM$-spectrum over $B$}\index{spectrum!Om@$\Omega$- -- over $B$} is a level $qf$-fibrant spectrum over $B$ such that each of its adjoint structure maps is a $q$-equi\-va\-lence; equivalently, its underlying prespectrum must be an $\OM$-prespectrum over $B$. 
\end{defn}

Since we are omitting the adjective ``orthogonal'' from ``orthogonal spectrum over $B$'', we must use the term ``$\OM$-prespectrum over $B$'' on the prespectrum level to avoid confusion; the more standard term ``$\OM$-spectrum'' was used in \cite{MM}.

\begin{prop}\mylabel{OmOm} 
A level fibrant $G$-spectrum $X$ over $B$ is an $\OM$-$G$-spectrum over $B$ if and only if each fiber $X_b$ is an $\OM$-$G_b$-spectrum. The $G$-prespectrum analogue also holds.
\end{prop}

\begin{proof} 
By the five lemma, this is immediate from a comparison of the long exact sequences in (\ref{les1}) and (\ref{les3}).
\end{proof}

This result leads to the following partial converse to \myref{levelpi}.

\begin{thm}\mylabel{bombsaway} 
A $\pi_*$-isomorphism between $\OM$-spectra over $B$ is a level $q$-equi\-va\-lence.
\end{thm}

\begin{proof}
The analogue for nonparametrized $\OM$-spectra is \cite[III.3.4]{MM}. In view of \myref{OmOm}, we can apply that result on fibers and then use that $\OM$-spectra over $B$ are required to be level $qf$-fibrant to deduce the claimed level 
$q$-equi\-va\-lence on total spaces from (\ref{les1}).
\end{proof}

Technically, the real force of our definition of homotopy groups is that this result describing the $\pi_*$-isomorphisms between $\OM$-spectra over $B$ is an immediate consequence of the work in \cite{MM}. Given this relationship between $\OM$-spectra and homotopy groups, many of the arguments of \cite{MM} apply fiberwise to allow the development of the \emph{stable model structure}.  However, as discussed in the next section, careful use of level fibrant approximation is required. We shall use the terms ``stable model structure'' and ``$s$-model structure'' interchangeably.  The $s$-cofibrations are the same as those of the level $qf$-model structure and the $s$-fibrant spectra over $B$ turn out to be the $\Omega$-spectra over $B$. 

\begin{defn}\mylabel{Def5} 
A map of spectra or prespectra over $B$ is 
\begin{enumerate}[(i)]
\item an \emph{acyclic $s$-cofibration} if it is a $\pi_*$-isomorphism and an $s$-cofibration,
\item an \emph{$s$-fibration}\index{fibration!s-@$s$- --} if it satisfies the RLP with respect to the acyclic $s$-co\-fi\-bra\-tions, 
\item an \emph{acyclic $s$-fibration} if it is a $\pi_*$-isomorphism and an $s$-fibration. 
\end{enumerate}
\end{defn}

We shall prove the following basic theorem in the next two sections.  

\begin{thm}\mylabel{modelS} 
The categories $G\sS_{B}$ and $G\sP_B$ are well-grounded model categories with respect to the 
$\pi_*$-iso\-mor\-phisms (= $s$-equi\-va\-lences), 
$s$-fibrations and $s$-cofibrations.  The $s$-fibrant 
objects are the $\Omega$-spectra over $B$.
\end{thm} 

\begin{rem}\mylabel{plus2} Recall \myref{plus1}. We can define positive
$\OM$-prespectra and positive analogues of our $s$-classes of maps,
starting with the positive level $qf$-model structure. 
As in \cite[III\S5]{MM}, the positive analogue of the previous theorem also holds, with the same proof.  The identity functor is the left adjoint of a Quillen equivalence from $G\sS_{B}$ or $G\sP_B$ with its positive stable 
model structure to $G\sS_{B}$ or $G\sP_B$ with its stable model structure. 
\end{rem}

The proof of the following result is virtually the same as the proof of its
nonparametrized precursor \cite[III.4.16 and III.5.7]{MM} and will not be 
repeated. 

\begin{thm}\mylabel{modelPU} 
The adjoint pair $(\bP,\bU)$ relating the categories $G\sP_{B}$ and $G\sS_B$ of prespectra and spectra over $B$ is a Quillen equivalence with respect to either the stable model structures or the positive stable model structures. 
\end{thm}

As in \cite[III.\S6]{MM}, \myref{modelS} leads to the following definition and theorem, whose proof is the same as the proof of \cite[III.6.1]{MM}.

\begin{defn}
Let $[X,Y]^{\ell}$ denote the morphism sets in the homotopy category associated to the level $qf$-model structure on $G\sP_B$ or $G\sS_B$. A map $f\colon  X\rtarr Y$ is a \emph{stable equivalence}\index{equivalence!stable --}\index{stable equivalence} if  $f^*\colon  [Y,E]^{\ell}\rtarr [X,E]^{\ell}$ is an isomorphism for all $\OM$-spectra $E$ over $B$. Define the positive analogues similarly. Let $[X,Y]$ denote the morphism sets in the stable homotopy category $\Ho G\sS_B$ of spectra over $B$.  
\end{defn}

\begin{thm}\mylabel{modelT}
The following are equivalent for a map $f\colon X\rtarr Y$ of spectra or prespectra over $B$.
\begin{enumerate}[(i)]
\item $f$ is a stable equivalence.
\item $f$ is a positive stable equivalence.
\item $f$ is a $\pi_*$-isomorphism. 
\end{enumerate} 
Moreover $[X,E] = [X,E]^{\ell}$ if $E$ is an $\OM$-spectrum.
\end{thm}

\myref{compaR} below should make it clear why the last statement is true.

\section{The $\pi_*$-isomorphisms}

In the main, the proof of \myref{modelS} is obtained by applying the results in \cite{MM} fiberwise. Since total spaces are no longer assumed to be weak Hausdorff, we have to be a little careful:  we are quoting results proven for $\sT$ and using them for $\sK_*$. However, we can just as well interpret \cite{MM} in terms of $\sK_*$.  The total spaces $X(V)$ of an $s$-cofibrant 
spectrum over $B$ are weak Hausdorff, hence $s$-cofibrant approximation
places us in a situation where total spaces are in $\sU$ and therefore fibers are in $\sT$.

There is a more substantial technical problem to overcome in adapting the proofs of \cite{MM, MMSS} to the present setting. In the situations encountered in those references, all objects were level $q$-fibrant, and that simplified matters considerably. Here, level $qf$-fibrant approximation entered into our definition of homotopy groups, and for that reason the results of this section are considerably more subtle than their counterparts in the cited sources.

We begin by noting that any level ex-quasifibrant approximation, not necessarily
a $qf$-fibrant approximation, can be used to calculate the homotopy groups of parametrized spectra.

\begin{lem}\mylabel{zigzag}
A zig-zag of level $q$-equivalences connecting a spectrum $X$ over $B$ to a level ex-quasifibrant spectrum $Y$ over $B$ induces an isomorphism between the homotopy groups of $X$ and of $Y$, and the latter can be computed directly in terms of the fibers of $Y$.
\end{lem}

\begin{proof}
This follows from \myref{levelpi} by applying a level $qf$-fibrant approximation functor to the zig-zag. 
\end{proof}

\begin{thm}\mylabel{exact}
Let $f\colon  X\rtarr Y$ be a map between $G$-spectra over $B$. For any $H\subset G$ and $b\in B^H$, there is a natural long exact sequence
\[ \cdots \rtarr \pi^H_{q+1}(Y_b) \rtarr \pi^H_q((F_B f)_b)\rtarr \pi^H_q(X_b)\rtarr \pi^H_q(Y_b)\rtarr \cdots\]
and, if $X$ is well-sectioned, there is also a natural long exact sequence
\[ \cdots \rtarr \pi^H_q(X_b)\rtarr \pi^H_q(Y_b)\rtarr \pi^H_{q}((C_B f)_b)\rtarr \pi^H_{q-1}(X_b)\rtarr \cdots.\]
\end{thm}

\begin{proof}
For the first long exact sequence, let $R$ be a level $qf$-fibrant approximation functor and consider $Rf$. We claim that the induced map $F_Bf\rtarr F_BRf$ is a level $q$-equivalence and that $F_BRf$ is level $qf$-fibrant. This means that $F_BRf$ is a level $qf$-fibrant approximation to $F_Bf$, so that the homotopy
groups of the fibers $(F_BRf)_b \iso F((Rf)_b)$ are the homotopy groups of $F_Bf$. When restricted to fibers over $b$, the parametrized fiber sequence 
$RX\rtarr RY\rtarr F_BRf$ of spectra over $B$ gives the 
nonparametrized fiber sequence $(RX)_b \rtarr (RY)_b \rtarr F((Rf)_b)$, and the long exact sequence follows from \cite[III.3.5]{MM}. To prove the claim, observe that since $F_B(I,Y)\rtarr Y$ is a Hurewicz fibration, it has a path-lifting function which levelwise shows that $F_B(I,Y)\rtarr Y$ is a level $fp$-fibration and therefore a level $qf$-fibration (since all $qf$-cofibrations are $fp$-cofibrations in $G\sK_B$). The dual gluing lemma (see \myref{hproper}(iii)) then gives that the induced map $F_Bf\rtarr F_BRf$ is a level $q$-equivalence. Since $F_B(I,-)$ preserves level $qf$-fibrant objects and since pullbacks of level $qf$-fibrant objects along a level $qf$-fibration are level $qf$-fibrant, $F_BRf$ is level $qf$-fibrant.

Since the maps $X\rtarr C_B X$ and $RX\rtarr C_B RX$ are $cyl$-cofibrations between well-sectioned spectra and therefore level $h$-cofibrations by \myref{levelcofs}, the gluing lemma gives that $C_B f\rtarr C_B Rf$ is a level $q$-equivalence. Since $RX$ and $RY$ are level well-sectioned and level $qf$-fibrant, they are level ex-quasifibrations. It follows from \myref{quasicof} that $C_B Rf$ is a level ex-quasifibration.  We cannot conclude that $C_BRf$ is
level $qf$-fibrant, but by \myref{zigzag} we can nevertheless use $C_B Rf$ to calculate the homotopy groups of $C_B f$. On fibers over $b$, the cofiber sequence of $Rf$ is just the cofiber sequence of $(Rf)_b$, and the long exact sequence follows from \cite[III.3.5]{MM}.
\end{proof}

Recall \myref{levelwellgr}, which specifies the ground structure
in $G\sS_B$ and shows that the level $q$-equivalences give a well-grounded subcategory of weak equivalences; the $g$-cofibrations are just the level
$h$-cofibrations.  The following result shows that the same 
is true for the $\pi_*$-isomorphisms. However, in contrast to \myref{levelwellgr}, it is crucial to assume that the relevant maps 
in the gluing and colimit lemmas are both $cyl$-cofibrations and 
$g$-cofibrations, as prescribed in \myref{hproper}.

\begin{thm}\mylabel{piwellgr} 
The $\pi_*$-isomorphisms in $G\sS_B$ give a well-grounded 
subcategory of weak equivalences. In detail, the following
statements hold. 
\noindent\begin{enumerate}[(i)]
\item A homotopy equivalence is a $\pi_*$-isomorphism.
\item The homotopy groups of a wedge of well-grounded spectra over $B$ are the direct sums of the homotopy groups of the wedge summands. 
\item The $\pi_*$-isomorphisms are preserved under pushouts along maps that are both $cyl$ and $g$-cofibrations.
\item Let $X$ be the colimit of a sequence $i_n\colon X_n\rtarr X_{n+1}$ of maps that are both $cyl$ and $g$-cofibrations and assume that $X/\!_B X_0$ is well-grounded. Then the homotopy groups of $X$ are the colimits of the homotopy groups of the $X_n$.
\item For a map $i\colon X\rtarr Y$ of well-grounded spectra over $B$ and a map $j\colon K\rtarr L$ of well-based spaces, $i\Box j$ is a $\pi_*$-isomorphism if either $i$ is a $\pi_*$-isomorphism or $j$ is a $q$-equivalence.
\end{enumerate}
\end{thm}

\begin{proof} 
The conclusion that the $\pi_*$-isomorphisms give a well-grounded
subcategory of weak equivalences, as prescribed in \myref{hproper},
follows directly from the listed properties, using \myref{gluederiv} to derive the gluing lemma. Since level $q$-equivalences are $\pi_*$-isomorphisms, 
$s$-cofibrant approximation in the level $qf$-model structure gives the
factorization hypothesis \myref{gluederiv}(ii).

A homotopy equivalence of spectra is a level $fp$-equivalence and hence a level $q$-equivalence, so (i) follows from \myref{levelpi}. For finite wedges, (ii) is immediate from the evident split cofiber sequences and \myref{exact}. For arbitrary wedges of well-grounded spectra over $B$, $\vee_B X_i \rtarr \vee_B RX_i$ is a level $q$-equivalence since the level $q$-equivalences are well-grounded and $\vee_B RX_i$ is level quasifibrant by \myref{quasicof}. 
By \myref{zigzag} we can use $\vee_B RX_i$ to calculate the homotopy groups of $\vee_B X_i$. Over a point $b$ in $B$, $\vee_B RX_i$ is just $\vee (RX_i)_b$ and the result follows from the nonparametrized analogue \cite[III.3.5]{MM}.

Now consider (iii). Let $i\colon X\rtarr Y$ be both a $cyl$-cofibration and
a $g$-cofibration and let $f\colon X\rtarr Z$ be a $\pi_*$-isomorphism. Since
$i$ and its $s$-cofibrant approximation $Qi$ are both $cyl$ and 
$g$-cofibrations and since the level $q$-equivalences give a well-grounded
subcategory of weak equivalences, the gluing lemma shows that we may 
approximate our given pushout diagram by one in which all objects are well-sectioned. Let $j\colon Z\rtarr Y\cup_X Z$ be the pushout of $i$ along $f$. Since $i$ and $j$ are $cyl$-cofibrations and $j$ is the pushout of $i$, their 
cofibers are homotopy equivalent. Comparing the long exact sequences of 
homotopy groups associated to the cofiber sequences of $i$ and $j$ gives 
that the pushout $Y\rtarr Y\cup_X Z$ of $f$ along $i$ is a $\pi_*$-isomorphism.

For (iv), we may use $s$-cofibrant approximation in the level model structure 
to replace our given tower by one in which all objects are well-sectioned.  We note as in the proof of \myref{modellim1} that the natural map $\text{Tel} X_n \rtarr \text{colim} X_n$ is a level $q$-equivalence and therefore a
$\pi_*$-isomorphism. Relating the telescope to a classical homotopy
coequalizer as in the cited proof, we reduce the calculation of the 
homotopy groups of the telescope to an algebraic inspection based on (ii). 
Alternatively, one can commute double colimits to reduce the verification
to its space level analogue.

For (v), it suffices to show that the tensor $X\sma_B T$ preserves 
$\pi_*$-isomorphisms in either variable, by \myref{boxacy}.  That follows 
from \myref{pismash} below.
\end{proof}

\begin{prop}\mylabel{pismash}
Let $f\colon X\rtarr Y$ be a map between well-grounded spectra over $B$.
\begin{enumerate}[(i)]
\item If $f$ is a level $q$-equivalence and $g\colon T\rtarr T'$ is a $q$-equivalence of well-based spaces, then
$$ \text{id}\sma_B g\colon X\sma_B T\rtarr X\sma_B T'$$
is a level $q$-equivalence and therefore a $\pi_*$-isomorphism.
\item If $f$ is a $\pi_*$-iso\-mor\-phism, then
$$ f\sma_B\text{id}\colon X\sma_B T\rtarr Y\sma_B T$$
is a $\pi_*$-isomorphism for any well-based space $T$ and
$$F_B(\text{id},f)\colon F_B(T,X)\rtarr F_B(T,Y)$$
is a $\pi_*$-isomorphism for any finite based CW complex $T$.
\item For a representation $V$ in $\sV$, the map $f$ is a $\pi_*$-iso\-mor\-phism if and only if $\SI^V_Bf$ is a $\pi_*$-isomorphism.
\end{enumerate}

\end{prop}

\begin{proof}
Part (i) holds since the level $q$-equivalences are well-grounded.
Therefore, for the first part of (ii), we may assume by $q$-cofibrant approximation in the space variable that $T$ is a based CW complex.  
Using \myref{quasicof}, it also implies that $-\sma_B T$ preserves approximations of well-grounded spectra over $B$ by level ex-quasifibrations. 
Now the first part of (ii) follows fiberwise from its nonparametrized analogue \cite[III.3.11]{MM} and (iii) follows fiberwise from its nonparametrized analogue \cite[III.3.6]{MM}.  Since $F_B(-,X)$ takes cofiber sequences of based spaces to fiber sequences of spectra over $B$, the second part of (iii) follows from the first exact sequence in \myref{exact}, as in the proof of \cite[III.3.9]{MM}.
\end{proof}

This leads to the following result, which shows that we are in a stable situation. 

\begin{prop}\mylabel{SIOMV}
For all well-grounded spectra $X$ over $B$ and all representations $V$ in $\sI_G$, the unit $\et\colon X\rtarr \OM_B^V\SI_B^VX$ and counit $\epz\colon  \SI_B^V\OM_B^VX\rtarr X$ of the $(\SI_B^V,\OM^V_B)$ adjunction are $\pi_*$-isomorphisms. Therefore, if $f\colon X\rtarr Y$ is a map between well-grounded spectra over $B$, then the natural maps $\eta\colon F_B f\rtarr \Omega_B C_B f$ and $\epsilon\colon  \Sigma_BF_Bf\rtarr C_B f$ are $\pi_*$-isomorphism.
\end{prop}
\begin{proof}
For $\eta$, after approximation of $X$ by an ex-quasifibration, the 
conclusion follows fiberwise from its nonparametrized analogue \cite[III.3.6]{MM}. Using the two out of three property and the
triangle equality for the adjunction, it follows that $\OM^V_B\epz$ 
is a $\pi_*$-isomorphism, hence so is $\epz$. For the last statement,
the maps $\et$ and $\epz$ are the parametrized analogues of the maps 
defined for ordinary loops and suspensions in \cite[p. 61]{Concise},
and they fit into diagrams relating fiber and cofiber sequences
like those displayed there. Now the last statement follows from the 
five lemma and the exact sequences in \myref{exact}.
\end{proof}

\section{Proofs of the model axioms}

We need some $G$-spectrum level recollections from \cite{MM} and their
analogues for $G$-spectra over $B$ to describe the generating acyclic 
$s$-cofibrations.  Let $(\sS_G,G\sS)$ denote the $G$-category of 
$G$-spectra. To keep track of enrichments, we return $G$ to the
notations for the moment.

We have a shift desuspension functor 
$F_V$ from based $G$-spaces to $G$-spectra given by $F_VT = V^*\sma T$, 
where $V^*(W) = \sJ_G(V,W)$ \cite[III.4.6]{MM}. It is left adjoint to 
evaluation at $V$.  For $G$-spectra $X$, the adjoint structure $G$-map 
$$\tilde{\si}\colon X(V)\rtarr \OM^WX(V\oplus W)$$ 
may be viewed by 
adjunction as a $G$-map
$$\tilde{\si}\colon \sS_{G}(F_VS^0,X) \rtarr \sS_{G}(F_{V\oplus W}S^W,X).$$
Passing to $G$-fixed points and taking $X= F_VS^0$, the image of the 
identity map gives a map of $G$-spectra
$$ \la^{V,W}\colon F_{V\oplus W}S^W\rtarr F_VS^0.$$
(The notation $\la_{V,W}$ was used in \cite{MM}, but we need room for a 
subscript). A Yoneda lemma argument then shows that the map of $G$-spaces
$$ \sS_G(\la^{V,W},\text{id})\colon \sS_G(F_VS^0,X)\rtarr \sS_G(F_{V\oplus W}S^W,X)$$
can be identified with $\tilde{\si}\colon X(V)\rtarr \OM^WX(V\oplus W)$.

We need the analogue for $G$-spectra over $B$. Recall from \myref{FVs} that, for an ex-$G$-space $K$ over $B$, 
$(F_V K)(W) = V^*(W)\wedge_B K$, where 
$$V^*(W) = \sJ_{G,B}(V,W) = \sJ_G(V,W)_B = (F_VS^0)(W)\sma_B S^0_B.$$  
It follows that we can identify $F_VK$ with the evident external tensor 
$F_VS^0 \sma_B K$ of the $G$-spectrum $F_VS^0$ and the ex-$G$-space $K$ over $B$; compare \myref{extsmash1}.  We have used the notation $\sma_B$ for this generalized tensor, but viewing it as a special case of the external smash product of spectra over $*$ and over $B$ would suggest the alternative notation $\barwedge$.  

\begin{defn}\mylabel{lambdas} 
For ex-$G$-spaces $K$ over $B$, we define a natural map
$$ \lambda^{V,W}_B\colon F_{V\oplus W}\Sigma^W_B K\rtarr F_VK.$$
Namely, identifying the source and target with external tensor products, define
$$\la_B^{V,W} = \lambda^{V,W}\sma_B \text{id}\colon (F_{V\oplus W}S^W)\sma_B K \rtarr (F_VS^0)\sma_B K.$$ 
\end{defn}

We can describe the adjoint structure maps of $G$-spectra over $B$ in
terms of these maps $\la^{V,W}_B$.

\begin{lem}\mylabel{eqn:lambdaadj}
Under the adjunctions 
\[P_B(F_VS^0_B, X)\cong F_B(S^0_B, X(V))\cong X(V)\]
and 
\[P_B(F_{V\oplus W}S^W_B, X)\cong F_B(S^0_B, \Omega^W_B X(V\oplus W))
\cong \Omega^W_BX(V\oplus W),\]
the map 
$$P_B(\lambda^{V,W}_B, \text{id})\colon  P_B(F_VS^0_B, X)\rtarr
P_B(F_{V\oplus W}S^W_B, X)$$
corresponds to 
$$\tilde\sigma\colon  X(V)\rtarr \Omega_B^W X(V\oplus W).$$
\end{lem}\begin{proof} 
When $X=F_VS^0_B$, the conclusion holds by comparison with the case
of $G$-spectra. The general case follows from the Yoneda lemma of 
enriched category theory.  See, for example, \cite[6.3.5]{Borceux}.
\end{proof}

We could have started off by defining $\la^{V,W}_B$ in a conceptual manner
analogous to our definition of $\la^{V,W}$, but we want the explicit 
description of $\la^{V,W}_B$ in terms of $\la^{V,W}$ in order to deduce homotopical properties in the parametrized context from homotopical properties in the nonparametrized context.  For that and other purposes, we need the following observation. We return to our convention of deleting $G$ from
the notations, on the understanding that everything is equivariant.

\begin{lem}\mylabel{gentensor}
If $\phi\colon X\rtarr Y$ is an $s$-equivalence of level well-based 
non\-par\-a\-me\-trized spectra and $K$ is a well-grounded ex-space with total space of the homotopy type of a $G$-CW complex, then 
$\phi\sma_B \text{id}\colon X\sma_B K\rtarr Y\sma_B K$ is an $s$-equivalence.
\end{lem}

\begin{proof} We use the ex-fibrant approximation functor $P$ of 
\myref{exfibapp}. We have a natural zig-zag of $h$-equivalences between 
$K$ and $PK$. 
By \myref{savior}, it induces a zig-zag of level $h$-equivalences between $X\sma_B K$ and $X\sma_B PK$ and, by \myref{HursmaK},  $X\sma_B PK$ is a level ex-fibration. Therefore, by \myref{zigzag}, it suffices to consider the case when $K$ is an ex-fibration. Since $(X\sma_B K)_b = X\sma K_b$ and $K_b$ is of the homotopy type of a $G_b$-CW complex, by \myref{ss}, each $(\phi\sma_B \text{id})_b$ is an $s$-equivalence by \cite[III.3.11]{MM}.
\end{proof}

The following result is crucial.

\begin{prop}\mylabel{key}
Let $K$ be a well-grounded ex-space with total space of the homotopy 
type of a CW complex. Then
\[\lambda^{V,W}_B\colon F_{V\oplus W} \Sigma^W_B K\rtarr F_VK\]
and 
\[\la^{V,W}\barwedge \text{id}\colon 
F_{V\oplus W}S^W \barwedge F_ZK \rtarr F_VS^0\barwedge F_ZK\]
are $\pi_*$-isomorphisms of spectra over $B$.
\end{prop}

\begin{proof}
Since $\la^{V,W}_B = \la^{V,W}\sma_B \text{id}$, \myref{gentensor} and the corresponding nonparametrized statement \cite[III.4.5]{MM} imply the first statement. For the second statement, observe that for spectra $X$ we have the associativity relation
$$X\barwedge F_ZK \cong X\barwedge (F_Z S^0\sma_B K)\cong 
(X\sma F_Z S^0)\sma_B K.$$
Taking $X = F_V T$ for a based space $T$ and using \myref{FVEV}, we see that 
\[F_VT \barwedge F_ZK \cong F_{V\oplus Z}(T \sma_B K).\]
Using equivalences of this form and checking definitions, we conclude that the map $\la^{V,W}\barwedge \text{id}$ of the statement can be identified with the map
$$\la^{V\oplus Z,W}\sma_B\text{id}\colon
(F_{V\oplus Z\oplus W}S^W)\sma_B K \rtarr 
(F_{V\oplus Z}S^0)\sma_B K.$$
Thus the second $\pi_*$-isomorphism is a special case of the first.
\end{proof}

From here, the proof of \myref{modelS} closely parallels arguments in \cite[III.\S4]{MM}, but simplified a little  by \myref{Newcompgen}. The generating set of $s$-cofibrations is again $FI^f_B$. The generating set $FK^f_B$ of acyclic $s$-cofibrations is given by a variant of the definition in the nonparametrized case \cite[III.4.6]{MM}.

\begin{defn}\mylabel{Def6} 
Recall the factorization of $\lambda^{V,W}$ through the mapping cylinder (in the category of spectra) as 
\[\xymatrix{\lambda^{V,W}\colon  F_{V\oplus W} S^W \ar[r]^-{k^{V,W}} & M\lambda^{V,W}\ar[r]^-{r^{V,W}} & F_VS^0.}\]
Here $k^{V,W}$ is an $s$-cofibration and $r^{V,W}$ is a deformation retraction. For $i\colon C\rtarr D$ in $I^f_B$, the map 
\[i\Box k^{V,W}\colon C\sma_B M\lambda^{V,W} \cup_{C\sma_B F_{V\oplus W} S^W} D\sma_B F_{V\oplus W} S^W \rtarr D\sma_B M\lambda^{V,W}\]
is an $s$-cofibration in $G\sS_B$ by \myref{Boxcof2}, and it is therefore also a $cyl$-cofibration by \myref{levelqf}.  It is a $\pi_*$-isomorphism by \myref{key} and inspection of definitions. The $s$-cofibrations in $FJ^f_B$ are level acyclic and are therefore also $\pi_*$-isomorphisms. Restricting to $V$ and $W$ in $\text{sk}\sI_G$, define the generating set $FK^f_B$ of acyclic $s$-cofibrations to be the union of $FJ^f_B$ and the set of all maps of the form $i\Box k^{V,W}$ with $i\in I^f_B$.
\end{defn}

A fortiori, the following result identifies the $s$-fibrations, but it must be proven a priori as a first step towards the verification of the model axioms.

\begin{prop}\mylabel{RLPL}
A map $f\colon X\rtarr Y$ satisfies the RLP with respect to $FK^f_B$ if and only if $f$ is a level $qf$-fibration and the diagrams 
\begin{equation}\label{OMpb}
\xymatrix{
X(V) \ar[r]^-{\tilde{\si}} \ar[d]_{f(V)} & \OM^W_B X(V\oplus W) 
\ar[d]^{\OM^W_B f(V\oplus W)} \\
Y(V) \ar[r]_-{\tilde{\si}} & \OM^W_B Y(V\oplus W) \\}
\end{equation}
are homotopy pullbacks for all $V$ and $W$.
\end{prop}

\begin{proof}  
As in \cite[III.4.7]{MM}, the homotopy pullback property must be interpreted as requiring a $q$-equi\-va\-lence from $X(V)$ into the pullback in the displayed diagram. Recall that $FJ^f_B$ is contained in $FK^f_B$ and that a map has the RLP with respect to $FJ^f_B$ if and only if it is a level $qf$-fibration. This gives part of both implications. It remains to show that a level $qf$-fibration $f$ has the RLP with respect to $i\Box k^{V,W}$ for all $i\in I^f_B$ if and only if the displayed diagram is a homotopy pullback. This is a formal but not altogether trivial exericise from the fact that the level $qf$-model structure is $G$-topological in the sense characterized in \myref{Gtopchar}. Notice that the map 
$i\Box k^{V,W}$ is isomorphic to the map $i\Box k^{V,W}_B$, where $k^{V,W}_B = k^{V,W}\sma_B S^0_B$. With notation as in (\ref{PBoxmap}), $f$ has the RLP with respect to $i\Box k^{V,W}_B$ for all $i\in I^f_B$ if and only if the pair 
$(i, P_B^\Box(k^{V,W}_B,f))$ has the lifting property for all $i\in I^f_B$, which holds
if and only if the map $P_B^\Box(k^{V,W}_B,f)$ of ex-spaces over $B$ is an acyclic $qf$-fibration.  This map is a $qf$-fibration since, for $j\in J^f_B$, the map $j\Box k^{V,W} \iso j\Box k^{V,W}_B$ is a level acyclic $s$-cofibration of spectra over $B$ by \myref{Boxcof2}. Since $f$ is a level $qf$-fibration, $(j\Box k^{V,W}_B, f)$ has the lifting property, hence, by adjunction, so does $(j, P_B^\Box(k^{V,W}_{B},f))$.   Finally, $P_B^\Box(k^{V,W}_{B},f)$ is homotopy equivalent to $P_B^\Box(\la^{V,W}_B,f)$
so one is a $q$-equi\-va\-lence if and only if the other is.  Under the isomorphisms  in \myref{eqn:lambdaadj}, the map $P_B^\Box(\lambda^{V,W}_{B},f)$ coincides with the map from $X(V)$ into the pullback in the displayed diagram and is thus a $q$-equi\-va\-lence if and only if that diagram is a homotopy pullback.
\end{proof}

Let $*_B$ be the terminal spectrum over $B$, so that each $*_B(V)$ is the terminal ex-space $*_B$. Observe that $*_B$ is an $\OM$-spectrum with trivial homotopy groups.

\begin{cor} 
The terminal map $F\rtarr *_B$ satisfies the RLP with respect to $FK_B$ if and only if $F$ is an $\OM$-spectrum over $B$.
\end{cor}

\begin{cor}\mylabel{stabislevel}
If $f\colon X\rtarr Y$ is a $\pi_*$-isomorphism that satisfies the RLP with respect to $FK_B$, then $f$ is a level acyclic $qf$-fibration. 
\end{cor}

\begin{proof}
Since $f$ is a level $qf$-fibration by \myref{RLPL}, the dual of the gluing lemma applied to the diagram
\[\xymatrix{
{*}_B\ar[r]\ar[d] & Y\ar@{=}[d] & X\ar[l]_f\ar@{=}[d] \\
F_B(I, Y) \ar[r] & Y & X\ar[l]^f}\]
gives that the induced map $F\rtarr F_Bf$ of pullbacks is a level $q$-equivalence. Since $f$ has the RLP with respect to $FK_B$, so does its pullback $F\rtarr *_B$. By the previous corollary, $F$ is thus an $\Omega$-spectrum over $B$. In particular, it is level $qf$-fibrant. We conclude that $F$ is a level $qf$-fibrant approximation for $F_Bf$. Since $f$ is a $\pi_*$-isomorphism, \myref{exact} gives that $F$ is acyclic. By \myref{bombsaway}, this implies that $F\rtarr *_B$ is a level $q$-equi\-va\-lence. Thus the fibers $F(V)_b$ all have trivial homotopy groups. We conclude (with a bit of extra argument as in \cite[9.8]{MMSS} to handle $\pi_0$) that each map of fibers $f(V)_b$ induces an isomorphism on homotopy groups. Therefore, since each $f(V)$ is a $qf$-fibration, each $f(V)$ induces an isomorphism on homotopy groups. 
\end{proof}

The proof of the model axioms for the stable model structure is now immediate.

\begin{proof}[Proof of \myref{modelS}]
The $\pi_*$-isomorphisms give a well-grounded subcategory of weak
equivalences, by \myref{piwellgr}.  Conditions (i), (iii), and (iv) in \myref{Newcompgen} are clear from our specification of the generating
acyclic $s$-cofibrations and the result for the level $qf$-model 
structure. For condition (ii), a $\pi_*$-isomorphism that satisfies 
the RLP with respect to $FK_B$ has the RLP with respect to $FI_B$ by \myref{stabislevel}. Conversely, a map that has the RLP with respect to $FI_B$ is a level acyclic $qf$-fibration and therefore has the RLP with respect to $FK_B$ by \myref{RLPL}. It is a $\pi_*$-isomorphism since it is level acyclic.
Since all $s$-fibrations are level $qf$-fibrations, right properness follows from the slightly stronger observation in the following result.
\end{proof}

\begin{prop}\mylabel{s-rightproper}
The $\pi_*$-isomorphisms in $G\sS_B$ are preserved under pullbacks along level $qf$-fibrations.
\end{prop}

\begin{proof}
Let $g$ be the pullback of a level $qf$-fibration $f$ along a $\pi_*$-isomorphism. Then $g$ is a level $qf$-fibration and the fibers of $g(V)$ 
are isomorphic to the fibers of $f(V)$. Therefore the homotopy fibers 
$F_Bg$ are level $q$-equivalent to the homotopy fibers $F_Bf$. 
The result follows by comparison of the first long exact sequence in \myref{exact} for $f$ and $g$.
\end{proof}

\section{Some Quillen adjoint pairs relating stable model structures}

We prove here that all of the adjoint pairs that were shown to be Quillen adjoints with respect to the level model structure in \S12.2 are still Quillen adjoints with respect to the stable model structure.  In view of the role played by level $qf$-fibrant approximation in our definition of homotopy groups, it is helpful to first understand the relationship between $s$-fibrant approximation and level $qf$-fibrant approximation. Now that the model structures have been established, we henceforward use the term $s$-equivalence rather than the synonymous term $\pi_*$-isomorphism. 

\begin{lem}\mylabel{compaR}
Let $\nu\colon X\rtarr RX$ and $\nu_{\ell}\colon X\rtarr R_{\ell}X$ be an $s$-fibrant approximation of $X$  and a level $qf$-fibrant approximation of $X$. Then there is an $s$-equivalence $\xi\colon R_{\ell} X\rtarr RX$ under $X$. 
\end{lem}

\begin{proof} 
Since $\nu_{\ell}$ is a level acyclic $s$-cofibration, it is an acyclic $s$-cofibration by \myref{levelpi}. Since $RX$ is $s$-fibrant, the RLP gives a map $\xi$ under $X$, and it is an $s$-equivalence since $\nu$ and $\nu_{\ell}$ are $s$-equivalences.
\end{proof}

We have the following relationship between the homotopy categories of ex-spaces over $B$ and of spectra over $B$.

\begin{prop}\mylabel{stablepair} 
The pair $(\SI^{\infty}_B,\OM^{\infty}_B)$ is a Quillen adjunction relating $G\sS_B$ and $G\sK_B$.  More generally, $(\SI^{\infty}_V,\OM^{\infty}_V) = (F_V,Ev_V)$ is a Quillen adjunction for any representation $V\in \sV$.
\end{prop}

\begin{proof}
The maps $\SI^{\infty}_V i$, where $i\in I^f_B$ is a generating cofibration for the $qf$-model structure on $G\sK_B$, are among the generating cofibrations of the $s$-model structure on $G\sS_B$, and it follows that $\SI^{\infty}_V$ preserves cofibrations. Since $\SI^{\infty}_V$ takes acyclic $qf$-cofibrations to level acyclic $qf$-cofibrations, and these are acyclic by \myref{levelpi}, $\SI^{\infty}_V$ also preserves acyclic cofibrations.
\end{proof}

Now consider an adjoint pair $(F,V)$ between categories of parametrized spectra that is a Quillen adjunction with respect to the level model structures.  Since the cofibrations are the same in the level model structure and in the stable model structure, the left adjoint $F$ certainly preserves cofibrations. Thus, 
to show that $(F,V)$ is also a Quillen adjunction with respect to the stable model structures, we need only show that $F$ carries acyclic $s$-cofibrations 
to $s$-equivalences. When $F$ preserves all $s$-equivalences, this is obvious;
otherwise, by \myref{reducts}, it suffices to verify this for the generating acyclic $s$-cofibrations. The cited result applies in general to subcategories of well-grounded weak equivalences, and in our context it applies to both the
level $q$-equivalences and the $s$-equivalences.  Recall that a Quillen left adjoint in any model structure preserves weak equivalences between cofibrant objects, by Ken Brown's lemma \cite[1.1.12]{Hovey}.  The following parenthetical observation applies to give a stronger conclusion for the Quillen left adjoints that we shall encounter.  It will play a crucial role in exploiting 
the equivalence of homotopy categories that we will establish in 
the next chapter. Note that the $s$-cofibrant spectra are the cofibrant 
objects in both the level and the stable model structures, and they are 
well-grounded.

\begin{prop}\mylabel{neat}
Let $F$ be a Quillen left adjoint between categories of parametrized 
spectra with their stable model structures and suppose that $F$ preserves 
level $q$-equivalences between well-grounded spectra. Then $F$ preserves 
$s$-equivalences between well-grounded spectra. 
\end{prop}

\begin{proof}
If $g\colon X\rtarr Y$ is an $s$-equivalence, where $X$ and $Y$ are well-grounded, factor $g$  in the level model structure as 
\[\xymatrix{X\ar[r]^{g'} & W \ar[r]^{g''} & Y,}\]
where $g'$ is an $s$-cofibration and $g''$ is a level acyclic $qf$-fibration. Then $W$ is well-grounded and $Fg''$ is a level $q$-equivalence by assumption. Since $F$ is a Quillen left adjoint in the $s$-model structures, $Fg'$ is an $s$-equivalence. Since level $q$-equivalences are $s$-equivalences it follows that $Fg=Fg''\circ Fg'$ is an $s$-equivalence.
\end{proof}

The following sequence of results consists of analogues for the stable
model structures of results proven for the level model structures in \S12.2.
Recall that we actually have well-grounded stable model structures $s(\sC)$ 
for any closed generating set $\sC$. 
As in \S12.2, wherever a $qf(\sC)$-model structure was used in Chapter 7 
for some particularly well chosen $\sC$, we must use the corresponding 
$s(\sC)$-model structure here.
 
\begin{prop}\mylabel{spacesmashpair}
Let $T$ be a based $G$-CW complex. Then $(-\sma_B T, F_B(T,-))$ is a Quillen adjunction on $G\sS_B$. When $T = S^V$, it is a Quillen equivalence.
\end{prop}

\begin{proof}
This is immediate from the fact that the stable model structure is 
$G$-topological, together with Propositions \ref{pismash} and \ref{SIOMV}.
\end{proof}

\begin{prop}\mylabel{Boxcof2too}
If $i\colon  X\rtarr Y$ and $j\colon  W\rtarr Z$ are $s$-cofibrations 
of spectra over base spaces $A$ and $B$, then 
\[i\Box j\colon  (Y\barwedge  W)\cup_{X\barwedge  W}(X\barwedge Z)\rtarr Y\barwedge Z\]
is an $s$-cofibration over $A\times B$ which is $s$-acyclic if either $i$ or $j$ is $s$-acyclic. 
\end{prop}
\begin{proof}
The statement about $s$-cofibrations is part of the analogue, \myref{Boxcof2},
for the level model structure.  
As usual, it suffices to show that $i\Box j$ is an $s$-equivalence if 
$i\in FI^f_B$ and $j\in FK^f_B$, where $FK^f_B$ is the set of generating 
acyclic $s$-cofibrations specified in \myref{Def6}. Arguing
as in \myref{boxacy} and using properness, this will hold if smashing the 
source and the target of $i$ with $j$ give $s$-equivalences. The reduction so far would work just as well for internal smash products. The required last step 
reduces via inspection of \myref{Def6} to an application of \myref{key}, with base space taken to be $A\times B$.  The reason that this last step works for external smash products but fails for internal smash products is made clear in \myref{ouchtoo}.
\end{proof}

\begin{cor}\mylabel{exttoo}
If $Y$ is an $s$-cofibrant spectrum over $B$, then the functor $(-)\barwedge Y$ from $G\sS_A$ to $G\sS_{A\times B}$ is a Quillen left adjoint with Quillen right adjoint $\bar{F}(Y,-)$. 
\end{cor}

\begin{prop}\mylabel{Qad1too} 
Let $f\colon A\rtarr B$ be a $G$-map. Then $(f_{!},f^*)$ is a Quillen adjoint pair.  If $f$ is a $q$-equivalence, then $(f_{!},f^*)$ is a Quillen equivalence.
\end{prop}

\begin{proof} 
We must show that $f_!$ takes acyclic $s$-co\-fi\-bra\-tions to $s$-equi\-va\-lences. Since $f_!$ preserves well-grounded objects and level $q$-equivalences between well-grounded objects by \myref{Qad1}, 
it suffices by \myref{reducts} to prove that $f_!k$ is an $s$-equivalence for each map $k$ in $FK^f_A$. This follows from the corresponding Quillen adjunction with respect to the level model structure if $k\in FJ^f_A$, so  assume that $k$ is of the form $i\Box k^{V,W}\iso i\Box k^{V,W}_A$. We claim that $f_!k$ is a map in $FK^f_B$ and is therefore an $s$-equivalence.  Observe that $k^{V,W}_A\iso f^*k^{V,W}_B$. Using (\ref{four}) and the fact that $f_!$ preserves pushouts, we see from the definition of the $\Box$-product that $f_!(i\Box f^*k^{V,W}_B)\iso (f_!i)\Box k^{V,W}_B$.  Since $i$ is obtained from a map over $A$ by adjoining a disjoint section, $f_!i$ is obtained from a map over $B$ by adjoining a disjoint section and is thus in $I^f_B$.

Now assume that $f$ is a $q$-equivalence.  By \cite[1.3.16]{Hovey}, $(f_!,f^*)$ is a Quillen equivalence if and only if $f^*$ reflects $s$-equivalences between $s$-fibrant objects and the composite $X\rtarr f^*f_!X\rtarr f^*Rf_!X$ given by the unit of the adjunction and $s$-fibrant approximation is an $s$-equivalence for all $s$-cofibrant $X$. Since the $s$-fibrant objects are the $\OM$-spectra over $B$ and the $s$-equivalences between $\OM$-spectra over $B$ are the level $q$-equivalences, the reflection property follows directly from the corresponding Quillen equivalence with respect to the level model structure. That result also gives that the composite $X\rtarr f^*f_!X\rtarr f^*R_{\ell}f_!X$ is a level $q$-equivalence and hence an $s$-equivalence. Applying \myref{compaR} with $X$ replaced by $f_!X$ and observing that $f^*$ preserves $s$-equivalences between level $qf$-fibrant $G$-spectra over $B$ since $(f^*Y)_a\iso Y_{f(a)}$, a little diagram chase shows that the composite $X\rtarr f^*f_!X\rtarr f^*Rf_!X$ is an $s$-equivalence.
\end{proof}

Observe that \myref{neat} applies to $f_!$.

\begin{prop}\mylabel{Qad2too}
Let $f\colon A\rtarr B$ be a $G$-bundle whose fibers $A_b$ are $G_b$-CW complexes. Then $(f^*,f_*)$ is a Quillen adjoint pair.
\end{prop}

\begin{proof} 
We must show that $f^*$ preserves acyclic $s$-cofibrations. Again it suffices by \myref{reducts} to prove that $f^*k$ is an $s$-equivalence between well-grounded spectra for each map $k\in FK^f_B$. That $f^*k$ is a map between well-grounded spectra follows from the fact that if $K\amalg B$ is a space over $B$ with a disjoint section, then $f^* F_V (K\amalg B) = F_V f^*K \amalg A$ is well-grounded. To see that $f^*k$ is an $s$-equivalence, it is enough, as in the proof of \myref{Qad1too}, to consider $k = i\Box k^{V,W}_B$ with $i\in I^f_B$. We have that $f^*k^{V,W}_B = k^{V,W}_A$ and, since $f^*$ preserves pushouts, smash products, and factorizations through mapping cylinders, we see as in the cited proof that $f^*k\iso f^*i\Box k^{V,W}_A$, which is an acyclic $s$-cofibration.
\end{proof}

\begin{prop}\mylabel{Lchanges2}
Let $\iota\colon H\rtarr G$ be the inclusion of a subgroup. The pair of functors $(\iota_!,\nu^*\iota^*)$ relating $H\sS_A$ and $G\sS_{\iota_!A}$ gives a Quillen equivalence. If $A$ is completely regular, then $\iota_!$ is also a Quillen right adjoint.
\end{prop}

\begin{proof}
By \myref{grprestrrQa} below, $(\iota_!,\nu^*\iota^*)$ is a Quillen adjoint pair. The proof that it is a Quillen equivalence is the same as the proof of 
the ex-space level analogue in \myref{Lishriek}. The last statement is less
obvious. As in the proof of the corresponding statement in
\myref{Lishriek}, it follows from the spectrum level analogue of
\myref{iotaalt}, which in turn requires the spectrum level analogue
of \myref{ouch0}, and the analogue in the stable model structure of \myref{fixedptrQa0}. The required analogues are proven in \S14.4 below.
\end{proof}

We shall see that $(\iota_!,\nu^*\iota^*)$ descends to a closed symmetric
monoidal equivalence of homotopy categories in \myref{Symmoni} below.

\begin{cor}\mylabel{LishriekCor2}
The functor $\io^*\colon \text{Ho}G\sS_B\rtarr \text{Ho}H\sS_{\io^*B}$ is the composite
\[\xymatrix@1
{\text{Ho}G\sS_B \ar[r]^-{\mu^*} & \text{Ho} G\sK_{\io_!\io^*B} \htp \text{Ho} H\sK_{\io^*B}\\}\]
\end{cor}

Using \myref{Johann} as in \myref{FibadQ}, the
following result is now a special case of 
Propositions \ref{Lchanges2} and \ref{Qad1too}. 

\begin{prop}\mylabel{FibadQtoo} 
For $b\in B$, the pair of functors $((-)^b,(-)_b)$ relating $G_b\sS$ and $G\sS_B$ is a Quillen adjoint pair. 
\end{prop}

\chapter{Adjunctions and compatibility relations}

\section*{Introduction}

The utility of the stable homotopy category $\Ho G\sS_B$ depends on the fact that the usual functors and adjunctions descend to it and still satisfy appropriate commutation relations. We consider such matters in this chapter. 
Many of our basic adjunctions are Quillen adjunctions in the stable model structure.  We recorded those in \S12.6.  The crucial adjunction missing from \S12.6 is $(f^*,f_*)$ for a general map $f$ of base spaces.  This cannot be a Quillen adjoint pair by the argument in \myref{noway}. We used Brown representability to construct the right adjoint $f_*$ between homotopy categories of ex-spaces in \myref{descendf0}. Analogously, in \S13.1 we use Brown representability to construct $f_*$ between homotopy categories of
parametrized spectra, and we use base change along diagonal maps to internalize smash products and function spectra.  There is an interesting twist here. It is not easy to verify the Mayer-Vietoris axiom directly.  Rather, we use the triangulated category variant of the Brown representability theorem, whose hypotheses turn out to be easier to check.

In \S13.7, we complete the proof that our stable homotopy categories are symmetric monoidal and prove some basic compatibility relations among smash products and base change functors.  These results involve commutation of Quillen left and right adjoints, and we would not know how to prove them using only model theoretic fibrant and cofibrant replacement functors.
Rather, their proofs depend on an equivalence between our model theoretic 
stable homotopy category of parametrized $G$-prespectra and a classical 
homotopy category of what we call ``excellent'' parametrized $G$-prespectra.
We used an analogous, but more elementary, equivalence of categories in Chapter 9.  It is essential to use parametrized $G$-prespectra rather than parametrized $G$-spectra to make the comparison since the relevant constructions do not all preserve functoriality on linear isometries; that is, they do not preserve $\sI_G$-spaces.  Results proven using the comparison are then translated to parametrized $G$-spectra along the Quillen equivalence between 
parametrized $G$-prespectra and parametrized $G$-spectra.

These equivalences of categories allow us to use a prespectrum level analogue $T$ of the ex-fibrant approximation functor $P$ to study derived functors.  
We define excellent parametrized $G$-prespectra in \S13.2.  We lift the ex-fibrant approximation functor $P$ from ex-$G$-spaces to parametrized $G$-spectra
in \S13.3.  There are several further twists here. First, the functor $P$ on 
ex-$G$-spaces does not behave well with respect to tensors, so extending it to 
a functor on parametrized $G$-prespectra is subtle. Second, with the extension, the zig-zag of $h$-equivalences connecting $P$ to the identity functor is no longer given by honest maps of parametrized $G$-prespectra, only weak maps. Third, the functor $P$ does not take parametrized $G$-prespectra to excellent ones. To remedy this, we introduce two auxiliary functors $K$ and $E$ in \S13.4. The composite $T=KEP$ does land in excellent parametrized $G$-prespectra, and $K$ converts weak maps to honest maps.  In \S\S13.5 and 13.6 we use $T$ to 
prove the promised equivalence of homotopy categories and show how to study derived funtors in this context.

There are few issues of equivariance in this chapter, and we generally continue to omit the (compact Lie) group $G$ from the notations.  We adopt the convention of calling isomorphisms in homotopy categories \emph{equivalences} and we denote them by $\simeq$ rather than $\cong$.

\section{Brown representability and the functors $f_*$ and $F_B$}

We need some preliminaries about the two versions of Brown representability
that are applicable in stable situations.  Recall \myref{Johann}. 

\begin{defn}\mylabel{detecting}
For $n\in \bZ$ and $H\subset G$, we have an $s$-cofibrant sphere $G$-spectrum $S^n_H$ such that $\pi_n^H(X) = [S^n_H,X]_G$ for all $G$-spectra $X$.  Explicitly, 
$$S^n_H = \begin{cases}
\SI^{\infty}(G/H_+ \sma S^n) & \text{if $n\geq 0$,}\\
F_{-n}(G/H_+\sma S^0) & \text{if $n<0$},
\end{cases}$$ 
as in \cite[II.4.7]{MM}, where $F_{-n}$ is the shift desuspension by $\bR^n$. 
We may allow the ambient group to vary.  Replacing $G$ by $G_b$ for $b\in B$ and letting $H\subset G_b$, define $S^{n,b}_H$ to be the $G$-spectrum $(S^n_H)^b$ over $B$. Note that $S^{n,b}_H$ is $s$-cofibrant, by \myref{FibadQ}.  By adjunction, for $G$-spectra $X$ over $B$, $\pi_n^H(X_b)$ is isomorphic to $[S^{n,b}_H,X]_{G,B}$.  Let $\sD_B$  be the set of all such $G$-spectra $S^{n,b}_H$ over $B$.
\end{defn}

From here, the following three results work in exactly the same way as their ex-space analogues in \S7.4.  Observe that the category $\Ho G\sK_B$ has coproducts and homotopy pushouts, hence homotopy colimits of directed sequences. 

\begin{lem}\mylabel{compact} 
Each $X\in \sD_B$ is compact, in the sense that 
$$\text{colim}\, [X, Y_n]_{G,B}\iso [X, \text{hocolim}\, Y_n]_{G,B}$$ 
for any sequence of maps $Y_n\rtarr Y_{n+1}$ in $G\sS_B$.
\end{lem}

\begin{prop}\mylabel{detect} 
A map $\xi\colon Y\rtarr Z$ in $G\sS_B$ is an $s$-equivalence if and only if the induced map $\xi_*\colon [X,Y]_{G,B}\rtarr [X,Z]_{G,B}$ is a bijection for all $X\in \sD_B$. 
\end{prop} 

\begin{proof} 
This is a tautology since as $X$ ranges through the $S^{n,b}_H$,
$[X,Y]_{G,B}$ ranges through the homotopy groups $\pi_n^H(Y_b)$
that define the $s$-equivalences. 
\end{proof}

\begin{thm}[Brown]\mylabel{brown} 
A contravariant set-valued functor on the cat\-egory $\Ho G\sS_B$ is representable if and only if it satisfies the wedge and Mayer-Vietoris axioms.
\end{thm}

Since we have the Quillen adjoint pair $(f_!,f^*)$, we have the right derived
functor $f^*\colon \Ho G\sS_B\rtarr \Ho G\sS_A$.   As in the proof of the analogous result on the level of ex-spaces, \myref{descendf0}, we can obtain the desired right adjoint $f_*$ to $f^*$ by use of Brown's theorem provided that we can show that $f^*$ preserves the relevant homotopy colimits. However, since $f^*\colon G\sS_B\rtarr G\sS_A$ does not preserve $s$-cofibrant objects, this is not obvious.  We will later give results that would allow us to carry out the proof in a manner analogous to the proof of \myref{descendf0}, but it is instructive to switch gears and give a more direct proof. It is based on the use of triangulated categories and would not have applied on the ex-space level.  

\begin{lem}\mylabel{yestrian}
The category $\Ho G\sS_B$ is triangulated.
\end{lem}

\begin{proof} 
The treatment of triangulated categories in \cite{Tri} gives a general pattern of proof for showing that homotopy categories associated to appropriate model categories are triangulated. It applies here. The distinguished triangles are those equivalent in $\Ho G\sS_B$ to cofiber sequences that start with a well-grounded spectrum or, equivalently by \myref{SIOMV}, those equivalent to the negatives of fiber sequences. Note that, by the proof of \myref{exact}, every cofiber sequence is equivalent in 
$\Ho G\sS_B$ to a cofiber sequence of level ex-quasifibrations.
\end{proof}

In triangulated categories, there is an alternative version of Brown's representability theorem due to Neeman \cite{Nee2}. It requires a ``detecting set of compact objects''. In triangulated categories with coproducts (or sums), an object $X$ is said to be compact if 
$\bigoplus [X, Y_i] \iso [X,\bigoplus Y_i]$ for any set of objects $Y_i$. In our topological situations, this reduces to the compactness of spheres, exactly as the proof of \myref{compact}.  A {\em detecting} set of objects is one that detects equivalences, in the sense suggested by \myref{detect}.  We have the following result.

\begin{lem} 
$\sD_B$ is a detecting set of compact objects in $\Ho G\sS_B$.
\end{lem}

Recall that an additive functor between triangulated categories is said to be {\em exact} if it commutes with $\SI$ up to a natural equivalence and preserves distinguished triangles.  The following theorems are proven in 
\cite[3.1, 4.1]{Nee2}; they are discussed with an eye to applications
such as ours in \cite[\S8]{FHM}.

\begin{thm}\mylabel{BR} 
Let $\sA$ be a compactly detected triangulated category. A functor $H\colon \sA^{op}\rtarr \sA b$ that takes distinguished triangles to long exact sequences and converts coproducts to products is representable.
\end{thm}

\begin{thm}\mylabel{TAFT} 
Let $\sA$ be a compactly detected triangulated category and $\sB$ be any triangulated category. An exact functor $F\colon \sA \rtarr \sB$ that 
preserves coproducts has a right adjoint $G$.
\end{thm}

\begin{thm}\mylabel{descendf}
For any $G$-map $f\colon  A\rtarr B$, there is a right adjoint 
$f_*$ to the functor $f^*\colon \Ho G\sS_B\rtarr \Ho G\sS_A$, 
so that
$$[f^*Y,X]_{G,A} \iso [Y,f_*X]_{G,B}$$
for $X$ in $G\sS_A$ and $Y$ in $G\sS_B$. 
\end{thm}

\begin{proof} The left adjoint $f_!$ commutes with $\SI$ and preserves
cofiber sequences, and this remains true after passage to derived homotopy
categories. Therefore the derived functor $f_!$ is exact.  Since $f^*$ is Quillen right adjoint to $f_!$, the derived functor $f^*$ is right adjoint to $f_!$ and is therefore also exact; see, for example, \cite[3.9]{Nee1}.  
If $X$ is in $\sD_A$, then $f_!X$ is compact in $\Ho G\sS_B$, as we 
see from commutation relations between relevant Quillen left adjoints
given in \myref{FVvsf*}. It follows formally that $f^*$ preserves 
coproducts, by \cite[5.1]{Nee2} or \cite[7.4]{FHM}. 
\end{proof}

\begin{rem} 
For composable maps $f$ and $g$, there is a natural equivalence $g_*\com f_* \simeq (g\com f)_*$ on homotopy categories since $f^*\com g^*\simeq (g\com f)^*$.
\end{rem}

Exactly as for ex-spaces in \myref{smashing0}, we apply change of base along the diagonal map $\DE\colon B\rtarr B\times B$ to obtain internal smash product and function spectra functors in $\Ho G\sS_B$.  

\begin{thm}\mylabel{smashing}  
Define $\sma_B$ and $F_B$ on $\Ho G\sS_B$ to be the composite (derived) functors
$$X\sma_B Y = \DE^*(X\barwedge Y) \quad\text{and}\quad F_B(X,Y) = \bar{F}(X,\DE_*Y).$$
Then 
$$ [X\sma_B Y, Z]_{G,B}\iso [X,F_B(Y,Z)]_{G,B}$$
for $X$, $Y$ and $Z$ in $\Ho G\sS_B$.
\end{thm}

\begin{proof}
The displayed adjunction is the composite of the adjunction for the external smash product and function spectra functors given by \myref{exttoo} and the adjunction $(\DE^*,\DE_*)$.
\end{proof}

\section{The category $G\sE_B$ of excellent prespectra over $B$}

We must still prove that $\Ho G\sS_B$ is a closed symmetric monoidal category under the derived internal smash product, that the derived functor $f^*$ is closed symmetric monoidal, and that various compatibility relations that hold on the point-set level descend to homotopy categories. In particular, since our right adjoints $f_*$, $\DE_*$, and therefore $F_B$ come from Brown's representability theorem, it is not at all obvious how to prove that they are well-behaved homotopically. In Chapter 9, we solved the corresponding ex-space level problems by proving that $\Ho G\sK_B$ is equivalent to the more classical and elementary homotopy category $hG\sW_B$. Here $G\sW_B$ is the category of ex-fibrations over $B$ whose total spaces are compactly generated and of the homotopy types of $G$-CW complexes, and $hG\sW_B$ is obtained from $G\sW_B$ simply by passage to homotopy classes of maps. This equivalence allowed us to exploit the ex-fibrant approximation functor $P$ of \S8.3 to resolve the 
cited problems. 

We shall resolve our spectrum level problems similarly, and the
following definitions give the appropriate analogues of $G\sW_B$ and
$hG\sW_B$.  However, to keep closer to the ex-space level, it is essential 
to work with parametrized prespectra rather than parametrized spectra.  It is safe to do so in view of the Quillen equivalence $(\bP,\bU)$ of \myref{modelPU} relating $G\sP_B$ and $G\sS_B$.

\begin{defn}\mylabel{Sigma}
Let $X$ be a $G$-prespectrum over $B$.
\begin{enumerate}[(i)]
\item $X$ is \emph{well-structured} if each level $X(V)$ is in $G\sW_B$.
\item $X$ is \emph{$\SI$-cofibrant} if it is well-grounded and each structure map \[\si\colon \SI^W_BX(V)\rtarr X(V\oplus W)\] is an $fp$-cofibration.
\end{enumerate}
\end{defn}

We can now give the definition of excellent $G$-prespectra over $B$ and of the associated classical homotopy category.  Working with classical nonequivariant
and nonparametrized coordinatized prespectra $\{E_n\}$, it has been known 
since the 1960's that the following definition gives the simplest 
quick and dirty rigorous construction of the stable homotopy category.

\begin{defn}\mylabel{excel}
The category $G\sE_B$ of \emph{excellent} $G$-prespectra over $B$ is the full subcategory of $G\sP_B$ whose objects are the well-structured $\Sigma$-cofibrant $\Omega$-$G$-prespectra over $B$. Let $hG\sE_B$ denote the classical homotopy category obtained from $G\sE_B$ by passage to homotopy classes of maps.
\end{defn}

We comment on the conditions we require of excellent prespectra over $B$. We require that they be well-structured so that we can exploit levelwise our equivalence of homotopy categories on the ex-space level. We require that
they be $\Sigma$-cofibrant since that provides ``homotopical glue'' that is necessary for the transition from the known equivalence on the ex-space level 
to the desired equivalence on the prespectrum level. We shall make this idea
precise shortly, in \myref{CPhelp}. We require that they be $\Omega$-prespectra over $B$ since it is clearly sensible to restrict attention to $s$-fibrant objects in $G\sS_B$ if we hope to compare homotopy categories.
Recall that $X$ is an $\Omega$-prespectrum if it is a level $qf$-fibrant prespectrum over $B$ whose adjoint structure maps
$$\tilde\sigma\colon  X(V) \rtarr \Omega_B^{W-V}X(W)$$ 
are $q$-equivalences. Since excellent prespectra over $B$ are required to be
level ex-fibrations, they are automatically level $qf$-fibrant. The condition 
on the adjoint structure maps is stronger than it appears on the surface.

\begin{lem}\mylabel{fpOM}
For excellent $G$-prespectra $X$ over $B$, the adjoint structure maps
$$\tilde{\si}\colon X(V)\rtarr \OM^W_BX(V\oplus W)$$
are $fp$-equivalences.
\end{lem}\begin{proof} The $\tilde{\si}$ are $q$-equivalences between
$G$-CW homotopy types and are therefore $h$-equivalences.
Since they are maps between ex-fibrations, they are
$fp$-equivalences by \myref{reverse}.
\end{proof}

This implies, for example, that homotopy-preserving functors $G\sE_B\rtarr G\sP_B$ that may not preserve level $q$-equivalences nevertheless do 
preserve the equivalence property required of the adjoint structure maps.

\begin{rem} Our definition of excellent parametrized prespectra is close to 
that used by Clapp and Puppe \cite{Clapp, CP}, who in turn were influenced by definitions in \cite{MQR}.  Curiously, while Clapp \cite{Clapp} focuses on ex-fibrations, Clapp and Puppe \cite{CP} never mention fibration conditions. These papers are nonequivariant, but the second is written in terms of what the authors call ``coordinate-free spectra'' over $B$. These are the same as our nonequivariant prespectra over $B$, except that their adjoint structure maps $\tilde{\si}$ are required to be closed inclusions, which holds automatically for $\SI$-cofibrant prespectra. Clapp and Puppe \cite{CP} use the term ``cofibrant'' for our notion of $\SI$-cofibrant.
\end{rem}

A crucial result of Clapp and Puppe makes the idea
of homotopical glue precise.  It is stated nonequivariantly in 
\cite[6.1]{CP}, but it works just as well equivariantly. 
Translated to our language, it reads as follows.

\begin{prop}[Clapp-Puppe]\mylabel{CPhelp}
If $f\colon X\rtarr Y$ is a level $fp$-equivalence between $\SI$-cofibrant prespectra over $B$, then $f$ is a homotopy equivalence of prespectra 
over $B$.  Therefore, if $f\colon X\rtarr Y$ is a level $h$-equivalence 
between well-structured $\SI$-cofibrant prespectra over $B$, then $f$ is a homotopy equivalence of prespectra over $B$.
\end{prop}

\begin{proof}[Sketch proof] 
The proof is analogous to the proof that a ladder of homotopy equivalences
connecting sequences of cofibrations induces a homotopy equivalence on passage
to colimits. The point is that, for $\SI$-cofibrant parametrized prespectra $Y$, we can carry out inductive arguments just as if $Y$ were just such a colimit.
Using standard cofibration arguments, carried over to the parametrized case, we can extend an $fp$-homotopy inverse of $\SI_B^{W_i}X(V_i)\rtarr \SI_B^{W_i}Y(V_i)$ to an $fp$-homotopy inverse of $X(V_{i+1})\rtarr Y(V_{i+1})$ and proceed inductively. The last statement follows by \myref{fpequiv}(i),
which shows that a level $h$-equivalence between well-structured prespectra
over $B$ is a level $fp$-equivalence.
\end{proof}

\section{The level ex-fibrant approximation functor $P$ on prespectra} 

We seek an approximation functor to play the role on the parametrized prespectrum level that the functor $P$ played on the ex-space level functor. 
We shall introduce three approximation functors, $P$, $E$ and $K$, that successively build in the properties of being well-structured, being an $\Omega$-prespectrum, and being $\Sigma$-cofibrant, each preserving the properties already obtained.  We define $P$ in this section and $E$ and
$K$ in the next.

Lifting the ex-space level functor $P$ of \S8.3 to the prespectrum level requires care. Recall that $P$ is the composite of the whiskering functor 
$W$ and the Moore mapping path space functor $L$, together with the 
natural zig-zag of $h$-equivalences 
\begin{equation}\label{repeat}
\xymatrix{K & WK \ar[l]_{\rho}\ar[r]^-{W\iota} & WLK=PK}
\end{equation}
of \myref{exfibapp} for ex-spaces $K$ over $B$. The functors $W$ and $L$ do not commute with tensors with based spaces, hence cannot be enriched over $G\sK_B$, by \myref{save}. There is therefore no canonical way of inducing structure maps after applying $P$ levelwise to a prespectrum, as one might at first hope. We shall resolve this by constructing by hand certain non-canonical but natural maps 
\begin{equation}\label{alphamap}
\alpha_V\colon  WK\sma_B S^V\rtarr W(K\sma_B S^V)\end{equation}
and
\begin{equation}\label{betamap}
\beta_V\colon  LK\sma_B S^V\rtarr L(K\sma_B S^V)
\end{equation}
such that 
$\al_0 = \text{id}$, $\beta_0=\text{id}$ and the following associativity 
diagram commutes, where $(F,f_V)$ stands for either $(W,\alpha_V)$ or $(L,\beta_V)$.
\begin{equation}\label{eqn:assoc}
\xymatrix{FK \sma_B S^V\sma_B S^{V'}\ar[r]^{f_V\sma \text{id}}\ar[d]_{\iso} & F(K\sma_B S^V)\sma_B S^{V'} \ar[r]^{f_{V'}} & F(K\sma_B S^V \sma_B S^{V'})\ar[d]^{\iso}\\
FK \sma_B S^{V\oplus V'} \ar[rr]^{f_{V\oplus V'}} && F(K\sma_B S^{V\oplus V'})}
\end{equation}

The definitions of these maps and the proofs that these diagrams commute depend
on chosen decompositions of $V$ and $V'$ as direct sums of indecomposable
representations, and we cannot choose compatible decompositions for all
representations $V$ and $V'$ at once. For this reason, and for other reasons
that will become apparent later, we must switch gears and work with sequentially indexed prespectra.

Thus, to be precise about the constructions in this section and the next, we restrict our original collection $\sV$ of indexing representations to a countable cofinal sequence $\sW$ of expanding representations in our given universe $U$.  More precisely, $\sW$ consists of representations $V_i$ for $i\geq 0$ such that $V_0=0$ and $V_i\subset V_{i+1}$. We set $W_i=V_{i+1}-V_{i}$. Such a sequence can be chosen in any universe.  We could just as well start with representations $W_i$ and
define $V_i$ inductively by $V_{i+1} = V_i\oplus W_i$.  There is no need to 
use orthogonal complements. We shall write in terms of complements, but on
the understanding that that is just a notational convenience.

\begin{rem}\mylabel{indexingreps}
There is a small quibble here since we originally defined our categories of parametrized prespectra only on collections of representations that are closed under finite direct sums, which $\sW$ clearly is not. However, if we let ${\sW'}$ consist of all finite sums of the $W_i$, then we recover such a collection.  As in \S11.3 (or \cite[\S2]{MMSS}), we can interpret $G\sP_B^{\sW'}$ as a diagram category indexed on a certain small category,
say $\sD^{\sW'}_G$, with object set $\sW'$, and we can interpret $G\sP_B^{\sW}$ 
as a diagram category indexed on the full subcategory $\sD^{\sW}_G$ of $\sD^{\sW'}_G$ whose object set is $\sW$. This gives a restriction functor
$\bU\colon G\sP_B^{\sW'}\rtarr G\sP_B^{\sW}$ that is right adjoint to a prolongation functor $\bP$ \cite[\S3]{MMSS}, and $(\bP,\bU)$ induces an adjoint equivalence of homotopy categories. We shall study such ``change of universe'' adjunctions in \S14.2. They allow us to lift all results we prove about the categories of parametrized prespectra indexed on cofinal sequences to our 
usual ones indexed on collections of representations closed under direct sums.
\end{rem}

\begin{defn}\mylabel{wierddef} Let $X$ be a prespectrum over $B$ indexed on the countable cofinal sequence $\sW = \{V_i\}$, where $V_{0} = 0$ and 
$V_{i+1} = V_i\oplus W_i$.  Let $X$ have structure maps 
$\si_i\colon \SI_B^{W_i}X(V_i)\rtarr X(V_{i+1})$. Then the maps 
$$W\sigma_i \circ \alpha\colon WX(V_i)\sma_B S^{W_i}\rtarr WX(V_{i+1})$$ 
and 
$$L\sigma_i \circ \beta\colon LX(V_i)\sma_B S^{W_i}\rtarr LX(V_{i+1})$$ 
specify structure maps for prespectra $WX$ and $LX$ over $B$. Therefore 
$PX = WLX$ is a prespectrum over $B$.
\end{defn}

Unfortunately, as will be clear from the following construction, the maps in the zig-zag (\ref{repeat}) do not lift to the prespectrum level. They only induce \emph{weak} maps of prespectra, that is, levelwise maps that only commute with the structure maps up to (canonical) 
$fp$-homotopy. Fortunately, the last approximation functor $K$, which arranges $\Sigma$-cofibrancy and will be discussed in the next section, turns weak maps 
into honest ones. 

\begin{con}\mylabel{wierdcon} 
We define $\al_V$ and $\be_V$. Fix a decomposition of $V$ into irreducible representations 
and let $\sP_V$ be the set of the projections from $V$ to the irreducible subrepresentations in this fixed decomposition. Define three equivariant 
maps from $V$ to the real numbers by setting
\[\|v\|_V=\max_{\pi\in \sP_V}|\pi v|, \quad
\mu_V(v)=\prod_{\pi\in \sP_V}(1-|\pi v|), \quad
\nu_V(v)=\prod_{\pi\in \sP_V}\max(1,|\pi v|).\]
Applying the same definitions to another representation $V'$ and to 
$V\oplus V'$ with its induced decomposition as a sum of irreducible representations, we see that the following equations hold.
\begin{gather*}
\|v\oplus v'\|_{V\oplus V'}=\max\{\|v\|_V, \|v'\|_{V'}\},\\
\mu_{V\oplus V'}(v\oplus v')=\mu_V(v)\mu_{V'}(v'),\\
\nu_{V\oplus V'}(v\oplus v')=\nu_V(v)\nu_{V'}(v').
\end{gather*}
Define a natural map
\[h_V\colon  WK \wedge_B S^V\sma_B [1,\infty)_+ \rtarr W(K \wedge_B S^V),\]
by setting
\begin{gather*}
h_V(x\sma v\sma t) = 
\begin{cases}
\ x\sma \mu(t^{-1}v)^{-1}\cdot v & \  \ \text{if $\|v\| \leq t$,}\\
(p(x),1-\nu(t^{-1}v)^{-1}) & \ \ \text{if $\|v\| \geq t$,}
\end{cases} \\
h_V((b,s)\sma v\sma t) = 
\begin{cases}
(b,s) & \text{if $\|v\| \leq t$,}\\
(b,1-(1-s)\nu(t^{-1}v)^{-1}) & \text{if $\|v\| \geq t$.}
\end{cases}
\end{gather*}
At time $t=1$ this specifies $\alpha_V$ and it is easy to verify that the associativity diagram (\ref{eqn:assoc}) commutes. Further, 
the map $\rho\circ h_V$ extends to $t=\infty$ to give an $fp$-homotopy from $\rho\circ\alpha_V$ to $\rho\sma_B \text{id}$. It follows that $\rho$ induces levelwise a weak map of prespectra $WX\rtarr X$.

Similarly define
\[k_V\colon  LK \sma_B S^V \sma_B [1,\infty)_+ \rtarr L(K \wedge_B S^V),\]
by setting
\[k_V((x,\lambda)\sma v\sma t) =
\begin{cases}
(x\sma  v,\lambda) & \text{if $\|v\| \leq t$,}\\
(x\sma v,\nu(t^{-1}v)\lambda & \text{if $\|v\| \geq t$.}
\end{cases}\]
Here, if $1\leq a <\infty$, and $\lambda\in \Lambda B$, then $a\lambda$ denotes the Moore path of length $l_\lambda/a$ given by $a\lambda (t)=\lambda(at)$. At time $t=1$ this specifies $\beta_V$, and it is again easy to check the
required associativity. The map $k_V\circ (\iota\sma\text{id})$ extends to an $fp$-homotopy from $\beta_V\circ (\iota\sma_B \text{id})$ to $\iota$, hence
$\io$ induces levelwise a weak map of prespectra $X\rtarr LX$, to which we
can apply $W$ to obtain a weak map $WX\rtarr WLX = PX$.
\end{con}

In view of \myref{exfibapp}, naturality arguments from
\myref{wierddef} and \myref{wierdcon} prove the following theorem.

\begin{thm}\mylabel{PPP} There are functors $L$, $W$, and 
$P = WL$ on $G\sP_B$ that are given levelwise by the 
functors $L$, $W$, and $P$ on $G\sK_B$. There are natural weak maps 
$\rh\colon WX \rtarr X$ and $\io\colon X\rtarr LX$ that are given
levelwise by the ex-space maps $\rh$ and $\io$. Therefore,
there is a natural zig-zag of weak maps $\ph = (\rh,W\io)$ as displayed 
in the diagram
$$\xymatrix{ 
X & WX \ar[l]_-{\rh} \ar[r]^-{W\io} & WLX = PX.}
$$
These maps are level $h$-equivalences, and $P$ converts level $h$-equivalences
to level $fp$-equivalences. If each $X(V)$ is compactly generated and of the 
homotopy type of a $G$-CW complex, then $PX$ is well-structured.  If $X$ is well-structured, then the weak maps in the above display are level $fp$-equivalences between well-structured $G$-prespectra over $B$. If, further,
the adjoint structure maps of $X$ are $h$-equivalences or $q$-equivalences, 
then so are the adjoint structure maps of $LX$, $WX$, and $PX$.
\end{thm}
\begin{proof}
The only point that may need elaboration is the last clause.  For a
weak map $f\colon X\rtarr Y$, we have a homotopy commutative diagram
$$\xymatrix{
X(V) \ar[r]^-{\tilde{\sigma}} \ar[d]_f & \Omega^W_B X(V\oplus W) \ar[d]^{\OM_B^Wf}\\
Y(V) \ar[r]_-{\tilde{\sigma}}          & \Omega^W_B Y(V\oplus W). \\}
$$
The functor $\Omega^W_B$ preserves $fp$-equivalences. Therefore, if
$f$ is an $fp$-equivalence, then the $\tilde{\si}$ for $X$ are
$h$-equivalences or $q$-equivalences if and only if the $\tilde{\si}$
for $Y$ are so. We apply this to $f=\rh$ and $f=W\io$.
\end{proof}

\section{The auxiliary approximation functors $K$ and $E$}\label{sec:K}

We begin with the parametrized $\Omega$-prespectrum approximation functor $E$. 
This is a folklore construction when $B$ is a point. In the parametrized 
context, the proof of the following result makes essential use of Stasheff's
theorem, \myref{ss}, and therefore depends on our standing assumption that 
$G$ acts properly on $B$.

\begin{prop}\mylabel{EEE}
There is a functor $E\colon G\sP_B\rtarr G\sP_B$ and a natural map $\al\colon X\rtarr EX$ with the following properties.
\begin{enumerate}[(i)]
\item The functor $E$ preserves level $fp$-equivalences and well-grounded prespectra. \item If $X$ is well-structured, then $EX$ is a well-structured $\OM$-prespectrum and the map $\alpha\colon X\rtarr EX$ is an $s$-equivalence.
\end{enumerate}
\end{prop}

\begin{proof}
Define $EX$ by letting $EX(V_i)$ be the telescope over $j\geq i$ of the ex-spaces $\OM^{V_j-V_i}_BX(V_j)$ with respect to the adjoint structure maps
\[ \OM^{V_j-V_i}_B\tilde{\si}\colon \OM^{V_j-V_i}_BX(V_j)
\rtarr \OM^{V_j-V_i}_B\OM^{W_j}_BX(V_{j+1})
\iso \OM^{V_{j+1}-V_i}_BX(V_{j+1}).\]
Since the functor $\OM^{W_i}_B$ commutes with telescopes,  $\OM_B^{W_i}EX(V_{i+1})$ is isomorphic to the telescope over $j\geq i+1$ of the ex-spaces $\OM^{V_j-V_{i+1}}_BX(V_j)$.  The adjoint structure map $EX(V_i)\rtarr \OM^{W_i}_B EX(V_{i+1})$ is induced by the maps $\Omega_B^{V_j-V_i}\tilde\sigma_j$ for $j\geq i$. The map $\al\colon X\rtarr EX$ is given by the inclusion of the bases of the telescopes. If $f\colon X\rtarr Y$ is a level $fp$-equivalence, then $Ef\colon EX\rtarr EY$ is a level $fp$-equivalence since a standard inductive argument (applicable in any topologically bicomplete category) shows that the telescope of a ladder of $fp$-equivalences is an $fp$-equivalence.

If $X$ is well-grounded or level ex-fibrant, then so is $EX$ since the construction clearly stays in the category of compactly generated spaces and since it preserves the conditions of being well-sectioned or level ex-fibrant by results in \S8.2. To show that $E$ preserves well-structured prespectra, it remains to show that if $X$ has total spaces of the homotopy types of $G$-CW complexes, then so does $EX$. By Stasheff's theorem (\myref{ss}), the fibers $X(V)_b = X_b(V)$ have the homotopy types of $G_b$-CW complexes.  We have the analogous construction $E$ in the category of $G_b$-prespectra and, by
Milnor's theorem (\myref{Milnor}) and standard facts about telescopes, the
$(E(X_b))(V)$ have the homotopy types of $G_b$-CW complexes. It is clear from the definition of $E$ that $(E(X_b))(V) = ((EX)(V))_b$. That is, the $G_b$-prespectrum $E(X_b)$ is the fiber $(EX)_b$ of the $G$-prespectrum $EX$ over 
$B$. By Stasheff's theorem again, it follows that the $(EX)(V)$ have the homotopy types of $G$-CW complexes.

To check that the adjoint structure maps are $q$-equivalences when $X$ is well-structured, it suffices to check that they induce $q$-equivalences on the fibers over $b$ for all $b\in B$. That holds by inspection of the homotopy groups of the colimits that define $(EX)_b \iso E(X_b)$. Similarly, we see that $\al$ is a $\pi_*$-equivalence when $X$ is well-structured by fiberwise comparison of the colimits of homotopy groups of fibers that define the homotopy groups of $X$ 
and $EX$.
\end{proof}

To approximate parametrized prespectra by level $fp$-equivalent $\Sigma$-cofibrant prespectra, we use the elementary cylinder construction $K$ that
was first defined in \cite{May69} and has been used in various papers since.  
We recall the construction and its main properties from \cite[6.8]{LMS}, which carries over verbatim to the parametrized context. A more sophisticated but less convenient treatment is given in \cite{EKMM}. 

\begin{prop}\mylabel{KKK} 
There is a functor $K\colon G\sP_B\rtarr G\sP_B$ and a natural level $fp$-equivalence $\pi\colon KX\rtarr X$.  Therefore $K$ preserves level $fp$-equivalences. If $X$ is well-grounded, then $KX$ is $\SI$-cofibrant. If $X$ is well-structured, then $KX$ is well-structured. If $X$ is a well-structured $\OM$-prespectrum, then so is $KX$ and thus $KX$ is excellent. There is a natural weak map $\iota\colon X\rtarr KX$ that is a right inverse of $\pi$, and $K$ takes weak maps $f$ to honest maps $Kf$ such that $\iota\circ f = Kf \circ \iota$. 
\end{prop}

\begin{proof}
Define $KX$, a level inclusion $\io\colon X\rtarr KX$, and a level $fp$-deformation retraction $\pi\colon KX\rtarr X$ right inverse to $\io$ as follows. Let $KX(0) = X(0)$ and $\io(0) = \pi(0) = \text{id}$. Inductively, suppose given $KX(V_i)$, an inclusion $\io(V_i)\colon X(V_i)\rtarr KX(V_i)$ and an inverse $fp$-deformation retraction $\pi(V_i)\colon KX(V_i)\rtarr X(V_i)$. Let $KX(V_{i+1})$ be the double mapping cylinder in $G\sK_B$ of the pair of maps
\[\xymatrix@=.6cm{\Sigma^{W_i}_B KX(V_i) && \Sigma^{W_i}_B X(V_i)\ar[ll]_-{\Sigma^{W_i}_B \iota(V_i)} \ar[rr]^-\sigma && X(V_{i+1})}\]
in $G\sK_B$. Let $\si\colon \SI^{W_i}_BKX(V_i)\rtarr KX(V_{i+1})$ be the inclusion of the left base of the double mapping cylinder, which is an
$fp$-cofibration and let $\io(V_{i+1})\colon X(V_{i+1})\rtarr KX(V_{i+1})$ be the inclusion of the right base. Let $\pi(V_{i+1})\colon KX(V_{i+1})\rtarr X(V_{i+1})$ be the map obtained by first using the $fp$-equivalence $\Sigma^{W_i}_B\pi(V_i)$ on the left base to map to the mapping cylinder of $\si$ and then using the evident deformation retraction to the right base. There is an equivalent description as a finite telescope. Certainly $\pi$ is a map of prespectra over $B$ and a level $fp$-deformation retraction with level inverse the weak map $\io$. The functoriality of the construction is clear.

If $X$ is well-grounded, then $KX$ is clearly also well-grounded and thus $KX$ is $\SI$-cofibrant.  If $X$ is well-structured, then so is $KX$ by Propositions \ref{pres1} and \ref{Hursma}. If, further, the adjoint structure maps of $X$ are $q$-equivalences, then they are $fp$-equivalences since $X$ is well-structured. Since $K$ preserves $fp$-homotopies, it follows that $KX$ is also an $\Omega$-prespectrum. Alternatively, since $\Omega_B^V$ is a Quillen right adjoint in the $qf$-model structure, it preserves $q$-equivalences between $qf$-fibrant ex-spaces. In particular, the maps $\Omega^W_B\pi(V_i)$ are $q$-equivalences.

If $f\colon X\rtarr Y$ is a weak map with $fp$-homotopies 
$$h_i\colon \Sigma_B^{W_i}X(V_i)\sma_B I_+ \rtarr Y(V_{i+1})$$ 
from $\sigma_Y\circ \Sigma^{W_i}f(V_i)$ to $f(V_{i+1})\circ \sigma_X$, define $Kf$ inductively by setting $Kf(0)=f(0)$ and letting $Kf(V_{i+1})$ be $\Sigma^{W_j}_BKf(V_i)$ on the left end of the mapping cylinder, $f(V_{i+1})$ on the right end and as follows on the cylinder itself:
\[Kf(V_{i+1})[x,t]=\begin{cases}
[\Sigma^{W_i}_B f(V_i)(x),2t] & \text{if $0\leq t\leq \tfrac12$},\\
h_i(x,2t-1) & \text{if $\tfrac12 \leq t\leq 1$.}
\end{cases}\]
Then $Kf$ is a map of prespectra over $B$ and $\iota\circ f = Kf\circ \iota$.
\end{proof}

The composite approximation functor $T = KEP$ has various good preservation properties. The ex-space level properties of $P$ recorded in \S8.4 are 
inherited on the prespectrum level, and we have the following sample result 
for $E$ and $K$. 

\begin{lem}
For a $G$-map $f\colon A\rtarr B$, a prespectrum $Y$ over $B$ and a prespectrum $X$ over $A$, there are natural isomorphisms
\[f^*EY \cong Ef^*Y,\quad f^*KY\cong Kf^*Y \quad\text{and}\quad Kf_!X\cong f_!KX.\]
\end{lem}

\begin{proof}
The relevant telescopes commute with $f^*$ since it is a symmetric monoidal left
adjoint and with $f_!$ since it is a left adjoint and the projection formula (\ref{four0}) holds.
\end{proof}

\section{The equivalence between $\Ho G\sP_B$ and $h G\sE_B$}

We can now extend the results of \S9.1 to parametrized prespectra. As in the previous section, our parametrized prespectra are indexed on a countable 
cofinal sequence of expanding representations in our given universe.  We
begin by collating the results of the previous two sections. 

\begin{thm}\mylabel{Tzigzag} 
Let $X$ be a well-grounded $G$-prespectrum over $B$ whose total spaces are of the homotopy types of $G$-CW complexes and define $TX = KEP X$. 
\begin{enumerate}[(i)]
\item $TX$ is an excellent $G$-prespectrum.
\item $T$ takes level $q$-equivalences between $G$-prespectra over $B$
that satisfy the hypotheses on $X$ to homotopy equivalences of $G$-prespectra.
\item There is a zig-zag of $s$-equivalences between $X$ and $TX$. 
\item If $X$ is an excellent $G$-prespectrum over $B$, then the zig-zag 
consists of level $fp$-equivalences, and it gives rise to a zig-zag of 
homotopy equivalences of $G$-prespectra over $B$ connecting $X$ and $TX$.
\end{enumerate}
\end{thm}

\begin{proof}
We have that $PX$ is well-structured by \myref{PPP}, $EPX$ is a well-structured $\OM$-prespectrum by \myref{EEE}, and $TX$ is excellent by \myref{KKK}. In (ii), a level $q$-equivalence is a level $h$-equivalence. By the results just quoted, $P$ takes level $h$-equivalences to level $fp$-equivalences, which are preserved by $E$, and $K$ takes level $fp$-equivalences to homotopy equivalences. Since $K$ converts weak maps to genuine maps, we have the following diagram of maps of $G$-presepectra over $B$.
\begin{equation}\label{GROWL}
\xymatrix{KX\ar[d]_\pi & KWX\ar[l]_{K\rho}\ar[d]^\pi & WKX\ar[d]_{W\pi}\ar[r]^{WK\iota} & WKLX\ar[d]^{W\pi} & KEPX\ar[d]^\pi\\
X & WX\ar@{=}[r] & WX & WLX\ar[r]_{\al} & EPX}
\end{equation}
The vertical maps $\pi$, hence also the vertical maps $W\pi$, are level $fp$-equivalences. The map $\rho$ is a level $f$-equivalence. The map $\io$ is a level $h$-equivalence, hence so is $WK\io$. The map $\al$ is an $s$-equivalence because $PX$ is well-structured. Since level $q$-equivalences are also
$s$-equivalences, the diagram displays a zig-zag of $s$-equivalences 
between $X$ and $TX$.

For the last statement, observe that all prespectra in the diagram are well-structured $\OM$-prespectra over $B$. Moreover, $\al$ is a level $q$-equivalence by \myref{bombsaway}. It is therefore a level $h$-equivalence since our total spaces have the homotopy types of $G$-CW complexes. Since all prespectra in our diagram are well-structured, our level $h$-equivalences are level $fp$-equivalences, by \myref{fpequiv2}. Applying $K$ where needed, we can expand the diagram to a zig-zag of level $fp$-equivalences between $\SI$-cofibrant prespectra.  By \myref{CPhelp}, this gives a zig-zag of homotopy equivalences connecting $X$ and $TX$.
\end{proof}

We introduce a category that is intermediate between $G\sP_B$ and $G\sE_B$.

\begin{defn}
Define $G\sQ_B$ to be the full subcategory of $G\sP_B$ consisting of the well-grounded $\Omega$-prespectra over $B$ whose total spaces are of the homotopy types of $G$-CW complexes. Define $\text{Ho} G\sQ_B$ to be the homotopy category obtained by inverting the $s$-equivalences in $G\sQ_B$; by the proof of the next theorem, there are no set-theoretic problems in defining $\text{Ho} G\sQ_B$.  Define $T = KEP\colon G\sQ_B \rtarr G\sE_B$. 
\end{defn}

Since the $\Omega$-prespectra over $B$ are the $s$-fibrant prespectra over $B$ and since $s$-cofibrant spectra are well-grounded and have total spaces of the homotopy types of $G$-CW complexes, all $G$-prespectra over $B$ that are 
$s$-cofibrant and $s$-fibrant are in $G\sQ_B$.  We prove that $\text{Ho}G\sP_B$ is equivalent to $hG\sE_B$ by proving that these categories are both equivalent to $\text{Ho} G\sQ_B$.

\begin{thm}\mylabel{exceq}
The canonical $s$-cofibrant and $s$-fibrant approximation functor $RQ$ and the composite approximation functor $T=KEP$, together with the forgetful functors, induce the following equivalences of homotopy categories.
\[\xymatrix{\Ho G\sP_B \ar@<.5ex>[r]^-{RQ} 
& \Ho G\sQ_B \ar@<.5ex>[r]^-{T}\ar@<.5ex>[l]^-I 
& hG\sE_B \ar@<.5ex>[l]^-J }\]
\end{thm}

\begin{proof}
For $X$ in $G\sP_B$, we have a natural zig-zag of $s$-equivalences in $G\sP_B$ 
\[\xymatrix{X & QX\ar[l]\ar[r] & RQX.}\]
Therefore $X$ and $IRQX$ are naturally $s$-equivalent in $G\sP_B$. If $X$ is in $G\sQ_B$, then it is $s$-fibrant and therefore so is $QX$. Then the above zig-zag is in $G\sQ_B$ so $X$ and $RQIX$ are naturally $s$-equivalent in $G\sQ_B$.

By \myref{bombsaway}, $s$-equivalences in $G\sQ_B$ are level $q$-equivalences, and $T$ takes level $q$-equivalences to homotopy equivalences by \myref{Tzigzag}.  Conversely, since homotopy equivalences are $s$-equivalences, the forgetful functor $J$ induces a functor in the other direction.

For $X$ in $G\sQ_B$ we have the natural zig-zag of $s$-equivalences displayed in (\ref{GROWL}). Applying $s$-fibrant approximation, we get a natural zig-zag of $s$-equivalences in $G\sQ_B$ so $X$ and $JTX$ are naturally $s$-equivalent in $G\sQ_B$. Starting with $X$ in $G\sE_B$, the last statement of \myref{Tzigzag} shows that $X$ and $TJX$ are naturally homotopy equivalent in $G\sE_B$.
\end{proof}

\section{Derived functors on homotopy categories}

With $P$ replaced by $T$, the discussion of derived functors in \S9.2 carries over from the level of ex-spaces to the level of parametrized prespectra indexed on cofinal sequences. In \S13.7 and \S14.2 we will discuss how to pass from there to conclusions on the level of parametrized spectra indexed on our usual collections of representations closed under direct sums. We must show that if $V$ is a Quillen left or right adjoint, then its model theoretic left or right derived functor agrees under our equivalences of categories with the functor obtained simply by passing to homotopy classes of maps from the composite $TV$. As on the ex-space level, we need some mild good behavior for this to work. 

\begin{defn}\mylabel{goodfun}
A functor $V\colon G\sP_A\rtarr G\sP_B$ is \emph{good} if it is continuous,  preserves well-grounded parametrized prespectra, and takes prespectra over $A$ whose levelwise total spaces are of the homotopy types of $G$-CW complexes to prespectra over $B$ with that property. Since $V$ is continuous, it preserves homotopies. There are evident variants for functors $V$ with source or target
$G\sK_*$: $V$ must be continuous, preserve well-grounded objects, and preserve
$G$-CW homotopy type conditions on objects.
\end{defn}

Note that a good functor $V$ need not take $\OM$-$G$-prespectra to
$\OM$-$G$-prespectra and recall that a Quillen right adjoint must
preserve fibrant objects and thus, in our context, must preserve
$\OM$-$G$-prespectra.

\begin{prop}\mylabel{scderiv}
Let $V\colon G\sP_A\rtarr G\sP_B$ be a good functor that is a part of a Quillen adjoint pair. If $V$ is a Quillen left adjoint, assume further that it preserves level $q$-equivalences between well-grounded objects. Then the derived functor $\Ho G\sP_A\rtarr\Ho G\sP_B$, induced by $VQ$ or $VR$, is equivalent to the functor $TVJ\colon hG\sE_A\rtarr hG\sE_B$ under the equivalence of categories 
in \myref{exceq}
\end{prop}

\begin{proof}
If $V$ is a Quillen right adjoint, then it preserves $s$-equivalences between $s$-fibrant objects. If $V$ is a Quillen left adjoint, then it preserves $s$-equivalences between well-grounded objects by \myref{neat}. Therefore, since $G\sQ_A$ consists of well-sectioned $s$-fibrant objects, the functor $V\colon  G\sQ_A \rtarr G\sP_B$ passes straight to homotopy categories to give $V\colon \text{Ho}G\sQ_A \rtarr \text{Ho}G\sP_B$ in both cases.

If $V$ is a Quillen right adjoint, then it takes an $s$-equivalence $f$ in $G\sQ_A$ to an $s$-equivalence since the objects of $G\sQ_A$ are $s$-fibrant. Then $Vf$ is a level $q$-equivalence by \myref{bombsaway} and, since $V$ is good, it is a level $h$-equivalence. On the other hand, if $V$ is a Quillen left adjoint, then \myref{bombsaway} gives that $f$ is a level $q$-equivalence and, by assumption, $Vf$ is then a level $q$-equivalence. Since $V$ is good, $Vf$ is actually a level $h$-equivalence. In both cases it follows that $V$ takes $s$-equivalences to level $h$-equivalences and therefore $TV$ passes to a functor $\Ho G\sQ_A\rtarr hG\sE_B$.

To show that $TVJ$ and either $VQ$ or $VR$ agree under the equivalence of categories, it suffices to verify that the following diagram commutes.
\[\xymatrix{
\Ho G\sP_A\ar[d]_{RQ} \ar[rr]^-{VQ \ \ \text{or} \ \ VR} 
& &  \Ho G\sP_B \ar[d]^{TRQ}\\
\Ho G\sQ_A \ar[rr]_{TV} & &  h G\sE_B}\]
We have functorial $s$-cofibrant and $s$-fibrant approximation functors $Q$ and $R$, with natural acyclic $s$-fibrations $QX\rtarr X$ and acyclic $s$-cofibrations $X\rtarr RX$. Clearly $Q$ and $R$ preserve $s$-equivalences. If $V$ is a Quillen left adjoint, then we have a zig-zag of natural $s$-equivalences
\[\xymatrix{RQVQ \ar[r] & RVQ & VQ \ar[l]\ar[r] & VRQ}\]
because $V$ preserves acyclic $s$-cofibrations.
If $V$ is a Quillen right adjoint, then we have a zig-zag of natural $s$-equivalences
\[\xymatrix{RQVR  & RQVRQ \ar[r]\ar[l] & RVRQ & VRQ\ar[l]}\]
because $V$ preserves $s$-equivalences between $s$-fibrant objects. In both cases, all objects have total spaces of the homotopy types of $G$-CW complexes, 
hence we have zig-zags of level $h$-equivalences. Applying $T$, we obtain a zig-zag of homotopy equivalences in $G\sE_B$ by \myref{Tzigzag}.
\end{proof}

\begin{rem}
If $V$ preserves excellent parametrized prespectra, then $TV$ is naturally homotopy equivalent to $V$ on excellent parametrized prespectra. The derived functor of $V$ can then be obtained directly by applying $V$ and passing to homotopy classes of maps. 
\end{rem}

\section{Compatibility relations for smash products and base change}

This section is parallel to \S9.3. The main change is just that we must replace the functor $P$ used there with the functor $T =KEP$ that we have here. This gives us results for the categories $G\sP_B^\sW$ of parametrized prespectra indexed on a collection $\sW$ consisting of a cofinal sequence in some universe $U$. In order to obtain statements about $G\sS_B^\sV$, where $\sV=\sV(U)$, we have two pairs of Quillen equivalences, both of which can be viewed as consisting of a prolongation functor left adjoint to a forgetful functor that creates the weak equivalences; see \cite[1.2]{MM}.
\[\xymatrix@1{G\sP_B^{\sW} \ar@<.5ex>[r]^-{j_*} 
& G\sP_B^\sV \ar@<.5ex>[r]^-{\bP}\ar@<.5ex>[l]^-{j^*} 
& G\sS_B^\sV \ar@<.5ex>[l]^{\bU}}\]
We postpone until \S14.2 consideration of the pair $(j_*,j^*)$ and the extension from $G\sP_B^\sW$ to $G\sP_B^\sV$ and deal with the extension from $G\sP_B^\sV$ to $G\sS_B^\sV$ in this section. 

One general remark is in order, though. The forgetful functors $j^*$ and $\bU$ create weak equivalences and therefore pass directly to homotopy categories. If they commute on the point set level with a functor $V$ which is a part of a Quillen adjoint pair, then they will also commute with its derived functor on the level of homotopy categories. It follows formally that the derived prolongation functors $\bP$ and $j_*$ then also commute with the derived functor $V$ and its adjoints. This holds in particular for the base change functor $V=f^*$. Extending commutation results for such functors from $G\sP_B^\sW$ to $G\sS_B^\sV$ is therefore easy. However, some of the functors $V$ that we need to consider only exist on some of the categories in the above display, and such
functors require special care. These include the change of universe functors that we discuss in \S14.2, which don't exist on the level of $G\sP_B^\sW$, and the smash product $\sma_B$, which we have only specified on the spectrum level and which we now discuss on the prespectrum level.

\begin{rem}\mylabel{subtlety}
Because the domain category for the diagram category of (equivariant and parametrized) prespectra is only monoidal, not symmetric mon\-oid\-al, 
we cannot use left Kan extension to internalize ``external'' smash products of prespectra; see \cite[4.1]{MMSS}. Here ``external'' is understood in the sense of indexing on pairs of representations.  Therefore, on the equivariant parametrized prespectrum level, when we write $X\barwedge Y$ for prespectra
$X$ over $A$ and $Y$ over $B$, we should understand the external external smash product, in the sense of \myref{extsmash1}. When passing from prespectrum level arguments to spectrum level conclusions using $(\bP, \bU)$, we are implicitly using composites of the general form $\bP V\bU$, and similarly for functors of several variables involving smash products.  We can carry out the several variable arguments externally on the prespectrum level, only internalizing 
with left Kan extension after passage to spectra, where we have good
homotopical control by \myref{exttoo}.

Alternatively, we can make use of classical ``handicrafted smash products'' of prespectra, which are defined by use of arbitrary choices of sequences of representations.  As discussed on the nonequivariant nonparametrized level in \cite[\S11]{MMSS}, handicrafted smash products of prespectra agree under the adjoint equivalence $(\bP,\bU)$ with the internalized smash products. Provided that we use external parametrized handicrafted smash products over varying base spaces, only internalizing along diagonal maps at the end, the discussion there adapts readily to give the same conclusion for homotopy categories of equivariant parametrized prespectra and spectra. The advantage of handicrafted smash products is that their definition involves only direct use of ex-space level constructions that enjoy good preservation properties with respect to ex-fibrations. This often allows direct transposition of ex-space level arguments in $h G\sW_B$ to parametrized prespectrum level arguments 
in $hG\sE_B$.
\end{rem}

We state the following results in terms of parametrized spectra, and we indicate which parts of the proofs require the use of $hG\sE_B$ and which parts work directly in the stable homotopy category $\text{Ho}G\sS_B$.

\begin{prop}\mylabel{helpme} 
Let $f\colon A\rtarr B$ and $g\colon A'\rtarr B'$ be $G$-maps. If $W$ and $X$ are spectra over $A$ and $A'$, then
$$f_!W\barwedge g_!X \simeq (f\times g)_!(W\barwedge X)$$
in $\Ho G\sS_{B\times B'}$. If $Y$ and $Z$ are spectra over $B$ and $B'$, then
$$f^*Y\barwedge g^*Z \simeq (f\times g)^*(Y\barwedge Z)$$
in $\Ho G\sS_{A\times A'}$.
\end{prop}

\begin{proof}
Working directly in $\Ho G\sS_{B\times B'}$, the first equivalence reduces to its point-set level analogue by consideration of Quillen left adjoints, as in the corresponding proof of \myref{fgext'}. We work in $h G\sE_{A\times A'}$ to prove the second equivalence. Here $f^*$ and $\barwedge$ (understood in the external or handicrafted sense) are both good, and both preserve excellent prespectra. Indeed, they preserve well-structured prespectra by levelwise application of Propositions \ref{fexpres} and \ref{Hursma}, they preserve $\Sigma$-cofibrant prespectra since $f^*$ and $\barwedge$ on ex-spaces preserve $fp$-cofibrations because they are left adjoints that commute with $fp$-homotopies, and they preserve $\Omega$-prespectra by \myref{fpOM} since they preserve $fp$-homotopies. Therefore, using excellent prespectra, we can pass straight to homotopy categories, without use of $T$, as in the corresponding proof of \myref{fgext'}.
\end{proof}

\begin{thm}\mylabel{symmonhtp} 
The category $\Ho G\sS_B$ is closed symmetric monoidal under the functors $\sma_B$ and $F_B$. 
\end{thm}

\begin{proof}
Working in $\Ho G\sS_B$, the associativity, commutativity, and unity of $\sma_B$ follow by pullback along diagonal maps from their easily proven external analogues and the second equivalence in the previous result, exactly as in \myref{clsymmon}.
\end{proof}

We have a commutation relation between change of base and suspension spectrum functors that is analogous to the relation between change of base and smash products recorded in \myref{helpme}.

\begin{prop}\mylabel{SIfSI} 
For a $G$-map $f\colon A\rtarr B$, there are natural equivalences
\[\SI^{\infty}_Bf_!\simeq f_!\SI^{\infty}_A
\quad\text{and}\quad \SI^{\infty}_Af^*\simeq f^*\SI^{\infty}_B\]
of (derived) functors. The same conclusion holds more generally for the shift desuspension functors $F_V = \SI^{\infty}_V$.
\end{prop}

\begin{proof}
Working in $\Ho G\sS_B$, the first equivalence is clear since it is a comparison of Quillen left adjoints that commute on the point-set level. For the second equivalence, we start in $h G\sW_B$ and end in $h G\sE_A$. For $K\in G\sW_B$, the point set level suspension prespectrum $\SI^{\infty}_BK$ is both $\Sigma$-cofibrant and well-structured, by \myref{HursmaK}, but of course it is not an $\OM$-prespectrum over $B$. Since $\Sigma^\infty_B$ is good and takes well-grounded $q$-equivalences to well-grounded level $q$-equivalences, $T\SI^{\infty}_B$ is equivalent to the model theoretic left derived functor of the Quillen left adjoint $\SI^{\infty}_B$. 
Here we may omit $P$ from the composite functor $T$ and, since $f^*$ commutes with both $K$ and $E$, the conclusion follows on passage to homotopy categories.
\end{proof}

Applying this to $\DE\colon B\rtarr B\times B$ and using \myref{SISISI}, we obtain the following consequence.

\begin{prop}\mylabel{SISISI2}
For ex-spaces $K$ and $L$ over $B$,
\[\SI^{\infty}_{B} (K\sma_B L) \simeq \SI^{\infty}_{B}K \sma_B  \SI^{\infty}_{B}L\]
in $\Ho G\sS_B$.
\end{prop}

For $f\colon A\rtarr B$, evident properties of the functor $f_!$ on ex-spaces 
imply that the functor $f_!\colon G\sP_A\rtarr G\sP_B$ is good, and $f_!$ satisfies the other hypotheses of \myref{scderiv} by \myref{Qad1}. We use 
this to prove the following basic result.

\begin{thm}\mylabel{Wirthmore}
For a $G$-map $f\colon A\rtarr B$ between base spaces, the derived functor $f^*\colon \Ho G\sS_B\rtarr \Ho G\sS_A$ is closed symmetric monoidal.
\end{thm}

\begin{proof}
Since $S_B$ is not $s$-fibrant, the isomorphism $f^*S_B\iso S_A$ in $G\sS_B$ does not immediately imply the required equivalence $f^*S_B\simeq S_A$ in $\Ho G\sS_A$, where $f^*S_B$ means $f^*RS_B$. However, \myref{SIfSI} specializes to give this equivalence. For the rest, we must show that the isomorphisms (\ref{one}) through (\ref{five}) descend to equivalences on homotopy categories.  By category theory in \cite{FHM}, it suffices to consider (\ref{one}) and (\ref{four}), and the proofs are similar to those in \myref{fclsymmon}. Since $\barwedge$ and $\Delta^*$ both preserve excellent prespectra, so do the internalized smash products $\sma_A$ and $\sma_B$. For excellent prespectra $Y$ and $Z$ over $B$, it follows that both sides of
\[f^*(Y\sma_B Z) \iso f^*Y \sma_A f^*Z\]
are excellent prespectra over $A$, hence the point-set level isomorphism descends directly to the desired equivalence on the homotopy category level. Next consider
\[f_!(f^*Y\sma_A X)\iso Y\sma_B f_!X,\]
where $X$ is an excellent prespectrum over $A$. Here we must replace $f_!$ by $Tf_!$ on both sides. By \myref{Tzigzag} we have a natural zig-zag $\phi$ of level $h$-equivalences connecting $T$ to the identity functor which, when applied to excellent parametrized prespecra gives rise to a zig-zag $\psi$ of actual homotopy equivalences. We obtain the following zig-zag.
\[\xymatrix{Tf_!(f^*Y\sma_A X)\cong  T(Y\sma_B f_!X) \ar[rr]^-{T(\text{id}\sma_B \phi)} && T(Y\sma_B Tf_!X) \ar[r]^-\psi \ar[ll] & Y\sma_B Tf_!X. \ar[l]}\]
Using handicrafted products with their termwise construction in terms of smash products of ex-spaces, it follows from \myref{savior} that $\text{id}\sma_B -$ preserves level $h$-equivalences between well-sectioned spectra. Thus $\text{id}\sma_B \phi$ is a zig-zag of level $h$-equivalences and $T(\text{id}\sma_B\phi)$ is a zig-zag of actual homotopy equivalences.
\end{proof}

\begin{thm}\mylabel{Mackeymore} 
Suppose given a pullback diagram of $G$-spaces 
\[\xymatrix{
C \ar[r]^-{g} \ar[d]_{i} & D \ar[d]^{j} \\
A \ar[r]_{f} & B}\]
in which $f$ (or $j$) is a $q$-fibration. Then there are natural equivalences of (derived) functors on stable homotopy categories
\begin{equation}\label{basesmore}
j^*f_{!} \simeq g_{!}i^*, \quad f^*j_* \simeq i_*g^*, \quad f^*j_{!}\simeq i_!g^*, \quad j^*f_*\simeq g_*i^*.
\end{equation} 
\end{thm}

\begin{proof}
Working in $hG\sE_B$, the proof is similar to that of \myref{pullbackfix} but with $P$ replaced by $T$. Again it suffices to consider the first equivalence, and, as explained there, since $f$ is a $q$-fibration there is a level $fp$-equivalence $\mu\colon Pf^*\rtarr f^*P$. Since $f^*$ commutes with both $K$ and $E$, we obtain a level $fp$-equivalence $\mu\colon Tf^*\rtarr f^*T$ between $\Sigma$-cofibrant prespectra over $A$ so it is in fact a homotopy equivalence by \myref{CPhelp}. Then $f^*Tj_!X\simeq Tf^*j_!X\cong Ti_!g^*X$. 
\end{proof}

The following observation holds by the same proof as the analogous ex-space level result \myref{imonoidaldescends}.

\begin{prop}\mylabel{Symmoni}
Let $\iota\colon H\rtarr G$ be the inclusion of a subgroup and $A$ be 
an $H$-space. The closed symmetric monoidal Quillen equivalence $(\iota_!, \nu^*\iota^*)$ descends to a closed symmetric monoidal equivalence between
$\text{Ho}H\sS_A$ and $\text{Ho}G\sS_{\io_!A}$.
\end{prop}

Combined with \myref{Wirthmore}, applied to the inclusion 
$\tilde{b}\colon G/G_b \rtarr B$, and \myref{SIfSI}, this last 
observation gives us the following stable analogue of \myref{fiberfun}.

\begin{thm}\mylabel{reassuring}
The derived fiber functor $(-)_b\colon \Ho G\sK_B\rtarr \Ho G_b\sK_*$ is closed symmetric monoidal and it has both a left adjoint $(-)^b$ and a right adjoint ${^b}(-)$. Moreover, the derived fiber functor commutes with the derived 
suspension spectrum functor, $(\SI^{\infty}_B K)_b\simeq \SI^{\infty}(K_b)$
as $G_b$-spectra.
\end{thm} 

For emphasis, we repeat a remark that we made after the analogous ex-space level result.  This innocent looking result packages highly non-trivial and important information. In particular, it gives that
$F_B(X,Y)_b\htp F(X_b,Y_b)$ in $\Ho G_b\sS$
for $X,Y\in \Ho G\sS_B$, where the fiber and function 
object functors are understood in the derived sense. This reassuring 
consistency result is central to our applications in the last two chapters,
where parametrized duality is studied fiberwise.

%
%

\chapter{Module categories, change of universe, and change of groups}

\section*{Introduction}

We first give a discussion of module categories of parametrized
spectra over nonparametrized ring spectra.  Although we shall not
go into these applications here, one basic motivation for our work
is to set up the homotopical foundations for studying the generalized 
homology and cohomology theories of parametrized spectra that are 
represented by such nonparametrized ring spectra.  The good behavior
of the external 
smash product $G\sS\times G\sS_B\rtarr G\sS_B$ makes it easy to do this.
While the mathematics here is evident, it deserves emphasis since
the ideas are likely to be central to future applications.

In the rest of the chapter, we focus on problems that are special to the equivariant context. We give the parametrized generalization of some of the
work in \cite{MM} concerning change of universe, change of groups, and fixed point and orbit spectra.  As usual, an essential point is to determine which 
of the standard adjunctions are given by Quillen adjoint pairs and to prove 
that other adjunctions and compatibilities that are evident on the point set level also descend to homotopy categories.  We discuss change of universe
in \S14.2.  Here the use of prespectra indexed on cofinal sequences in the
previous chapter introduces some minor difficulties that were not studied in
the nonparametrized theory of \cite[V\S1]{MM} and are already relevant
nonequivariantly.  We study subgroups and fixed point spectra in \S14.3. We study quotient groups and orbit spectra in \S14.4.  Aside from some analogues for parametrized spectra 
of earlier results for parametrized spaces, these sections are precisely
parallel to \cite[V\S\S2 and 3]{MM}. We have not written down the parametrized analogue of \cite[V\S4]{MM}, which gives the theory of geometric fixed point
spectra, since it would be tedious to repeat the constructions given there.
It will be apparent to the interested reader that, mutatis mutandis, the definitions and results in \cite[V\S4]{MM} generalize to the parametrized
context.

\section{Parametrized module $G$-spectra}

We can define a parametrized (strict) ring $G$-spectrum $R$ over $B$ to be a monoid in the symmetric monoidal category $G\sS_B$, and we can then define parametrized $R$-modules and $R$-algebras in the usual way, as has become standard in stable homotopy theory \cite{EKMM, HSS, MM, MMSS}. However, even though the 
smash product $\sma_B$ in $G\sS_B$ gives a point-set level symmetric monoidal structure, we cannot expect to obtain Quillen model structures on the categories
of such $R$-modules or $R$-algebras, as was done for orthogonal $G$-spectra in \cite[III\S\S7,8]{MM}. To do that, we would need better homotopical behavior than we can prove here.  We have only set up adequate foundations for the classical style theory of up to homotopy parametrized module spectra over up to homotopy parametrized ring spectra. From that point of view, our
homotopical foundations are entirely satisfactory. The source of the problem is \myref{ouchtoo}, which implies that $X\sma_B (-) $ in $G\sS_B$ cannot be a Quillen functor. 

However, in applications, it is natural to start with a nonparametrized orthogonal ring $G$-spectrum $R$. We are then interested in understanding 
the $R$-homology and $R$-cohomology theories of $G$-spectra over $B$ and their relationships with the $R$-homology and $R$-cohomology of the fibers. For this study, just as in the nonparametrized work of \cite{EKMM, HSS, MM, MMSS}, one 
is interested in the theory of $R$-modules.  The \emph{external} smash product $\barwedge\colon G\sS\times G\sS_B \rtarr G\sS_B$ has enough of the good properties of the nonparametrized smash product $G\sS\times G\sS \rtarr G\sS$ 
to give us homotopical control over parametrized module spectra over nonparametrized ring spectra. We devote this section to 
developing the relevant theory, which is parallel to \cite[III\S7]{MM}. 
Let $R$ be a ring spectrum in $G\sS$ which is well-grounded when viewed 
as a spectrum, meaning that each $R(V)$ is well-based and 
compactly generated.

\begin{defn} A {\em (left) $R$-module over $B$} is a $G$-spectrum $M$ 
over $B$ together with a left action $R\barwedge M\rtarr M$ satisfying 
the usual associativity and unit conditions. The category $GR\sM_B$ of 
left $R$-modules over $B$ consists of the $G$-spectra $M$ over $B$ and
the maps of $G$-spectra over $B$ that preserve the action by $R$.
\end{defn}

Since $(R\barwedge X)_b = R\sma X_b$, a parametrized $R$-module over $B$ is 
precisely that: each $X_b$ is an $R$-module $G_b$-spectrum.  More
formally, we have the $G$-category $(R\sM_{G,B},GR\sM_B)$, as 
discussed in \S\S1.4 and 12.2, and the following result is clear.  

\begin{prop}
The $G$-category $(R\sM_{G,B},GR\sM_B)$ is $G$-topologically bicomplete in the sense of \myref{defn:enrichBG}. All of the required limits, colimits, 
tensors, and cotensors are constructed in the underlying $G$-category $(\sS_{G,B},G\sS_B)$ and then given induced $R$-module structures in 
the evident way.  A $\text{cyl}$-cofibration of $R$-modules 
is a $\text{cyl}$-cofibration of underlying $G$-spectra 
over $B$.
\end{prop}

The last statement holds by the retract of mapping cylinders characterization of $\text{cyl}$-cofibrations. This immediately implies that $GR\sM_B$ inherits a ground structure from $G\sS_B$, in the sense of \myref{back}. Recall that the well-grounded $G$-spectra over $B$ are those that are level well-grounded (well-sectioned and compactly generated) and that the $g$-cofibrations of $G$-spectra over $B$ are the level $h$-cofibrations; see \myref{sillybilly} and \myref{levelwellgr}.

\begin{defn}\mylabel{Rground}
An $R$-module over $B$ is well-grounded if its underlying
$G$-spectrum over $B$ is well-grounded.  A map of $R$-modules over $B$
is a $g$-cofibration, level $q$-equivalence, or $s$-equivalence if its
underlying map of $G$-spectra over $B$ is such a map.
\end{defn}

Also recall the notion of a subcategory of well-grounded weak equivalences from \myref{hproper}. Since colimits and tensors for $R$-modules are defined in terms of the underlying $G$-spectra over $B$, the following theorem is immediate from its counterpart for $G$-spectra over $B$, which is given by \myref{levelwellgr} and \myref{piwellgr}.

\begin{thm}
\myref{Rground} specifies a ground structure on $GR\sM_B$
such that the level $q$-equivalences and the $s$-equivalences 
both give subcategories of well-grounded weak equivalences.
\end{thm}

Finally, recall the definition of a well-grounded model structure from 
\myref{wellmodel}. Such model structures are compactly generated, and we 
must define the generators of $GR\sM_B$.  The free $R$-module 
functor $\bF_R=R\barwedge - \colon G\sS_B \rtarr GR\sM_B$ is left 
adjoint to the forgetful functor $\bU\colon GR\sM_B\rtarr G\sS_B$. 
Adjunction arguments from the definitions show that $\bF_R$ preserves 
$\text{cyl}$-cofibrations and $\overline{\text{cyl}}$-cofibrations. 

\begin{defn}
Define $\bF_R FI^f_B$, $\bF_R FJ^f_B$ and $\bF_R FK^f_B$ by applying the 
free $R$-module functor to the maps in the sets specified in \myref{FBJB} 
and \myref{Def6}.
A map of $R$-modules over $B$ is 
\begin{enumerate}[(i)]
\item a level $qf$-fibration or an $s$-fibration if it is one in $G\sS_B$,
\item an $s$-cofibration if it satisfies the LLP with respect
to the level acyclic $qf$-fibrations,
\end{enumerate}
\end{defn}

\begin{thm}
The category $GR\sM_B$ is a well-grounded model category with respect to the level $q$-equivalences, the level $qf$-fibrations, and the $s$-cofibrations. The sets $\bF_R FI^f_B$ and $\bF_R FJ^f_B$ give the generating $s$-cofibrations and generating level acyclic $s$-cofibrations.  All $s$-cofibrations of $R$-modules over $B$ are $s$-cofibrations of $G$-spectra over $B$.
\end{thm}

We omit the proof since it is virtually the same as the proof of the
following theorem, which gives the starting point for serious work on the homology and cohomology theory of parametrized $G$-spectra.

\begin{thm}
The category $GR\sM_B$ is a well-grounded model category with respect to the $s$-equivalences, the $s$-fibrations, and the $s$-cofibrations; $\bF_RFK^f_B$ gives the generating acyclic $s$-cofibrations.
\end{thm}

\begin{proof}
The compatibility condition is automatic by adjunction from the para\-metrized spectrum level, and we have already observed that the free $R$-module functor $\bF_R$ preserves $\overline{cyl}$-cofibrations. It also preserves the relevant $\Box$-products, and $\bF_R F_VK = (R\sma F_VS^0)\barwedge K$ is well-grounded 
if $K$ is a well-grounded ex-space. Only the acyclicity condition remains.
If $R$ is $s$-cofibrant as a ring spectrum, then $R$ is also $s$-cofibrant 
as a spectrum, by \cite[III.7.6(iv) and (v)]{MM}. In that case, the functor 
$R\barwedge (-)=\bU\bF_R$ is a Quillen left adjoint by \myref{exttoo} and therefore preserves level acyclic $s$-cofibrations. It follows that the 
maps in $\bF_RK^f_B$ are $s$-equivalences. The case of a general
well-grounded $R$ reduces to the cofibrant case by use of the next result;
compare \myref{Rinvar} below.
\end{proof}

\begin{prop}\mylabel{barwcof}
The following statements hold.
\begin{enumerate}[(i)]
\item For an $s$-cofibrant spectrum $X$ over $B$, the functor $-\barwedge X\colon G\sS \rtarr G\sS_B$ preserves $s$-equivalences between well-grounded spectra in $G\sS$.
\item If $Y$ is well-grounded in $G\sS$, $j\colon A\rtarr X$ is an acyclic $s$-cofibration in $G\sS_B$, and $A$ is well-grounded, then $Y\barwedge j\colon Y\barwedge A\rtarr Y\barwedge X$ is an $s$-equivalence.
\end{enumerate}
\end{prop}

\begin{proof} Let $\ph\colon Y\rtarr Z$ be an $s$-equivalence between
well-grounded spectra. By parts (ii)--(iv) of \myref{hproper},  it 
suffices to show that $\ph\barwedge F_V K$ is an $s$-equivalence if
$K$ is the source or target of a map in $I_B^f$. This map is isomorphic
to the map $(\ph\sma F_VS^0)\sma_B K$, where $F_VS^0$ is the shift desuspension
in $G\sS$, not $G\sS_B$. Here $\ph\sma F_VS^0$  is an $s$-equivalence by
the nonparametrized analogue \cite[III.7.3]{MM}, and the 
conclusion follows from \myref{gentensor}. (The comment on the notations $\barwedge$ and $\sma_B$ above \myref{lambdas} is relevant: the former
is an external smash product and the latter is a tensor).

For (ii), we apply an argument from \cite[12.5]{MMSS}. We let $Z = X/\!_BA$,
which is $s$-cofibrant, and we let $QY\rtarr Y$ be an $s$-cofibrant approximation. Since $j$ is an $s$-cofibration, it is a 
$\text{cyl}$-cofibration and $Cj$ is homotopy equivalent to $Z$.
Since $A$ is well-grounded, we can apply the long exact sequence
of homotopy groups of \myref{exact} to conclude that $Z$ is $s$-acyclic.
The map $Z\rtarr *_B$ is then an $s$-equivalence between $s$-cofibrant spectra over $B$.  Since $QY\barwedge -$ is a Quillen left adjoint, by \myref{Boxcof2}, $QY\barwedge Z\rtarr QY\barwedge *_B \cong *_B$ is an $s$-equivalence. Since $QY\barwedge Z \rtarr Y\barwedge Z$ is an $s$-equivalence by part (i),  $Y\barwedge Z$ is $s$-acyclic. Since the functor $Y\barwedge -$ preserves cofiber 
sequences, another application of \myref{exact} shows that $Y\barwedge j$ is an $s$-equivalence.
\end{proof}

\begin{prop}\mylabel{Rinvar}
If $\phi\colon Q\rtarr R$ is an $s$-equivalence of well-grounded ring spectra, then the functors 
$$\phi_*=R\sma_Q(-)\colon GQ\sM_B\rtarr GR\sM_B \ \ \text{and} \ \  
\phi^*\colon GR\sM_B\rtarr GQ\sM_B $$ 
given by extension of scalars and restriction of action define a Quillen equivalence $(\phi_*,\phi^*)$ between the categories of 
$Q$-modules and of $R$-modules over $B$.
\end{prop}

\begin{proof}
Since $s$-fibrations and $s$-equivalences are created in the underlying
category of spectra over $B$, it is clear that they are preserved by
$\ph^*$, so that we have a Quillen pair. If $M$ is an $s$-cofibrant $Q$-module, then, by the previous result, the unit map 
$\phi\sma\text{id}\colon M\cong Q\sma_Q M\rtarr \ph^*(R\sma_Q M)$ 
of the adjunction is an $s$-equivalence of spectra over $B$.  Therefore, if $N$ is an $s$-fibrant $R$-module, then a map $M\rtarr \phi^*N$ of $Q$-modules is an $s$-equivalence if and only if its adjoint map $R\sma_Q M\rtarr N$ of $R$-modules is an $s$-equivalence.
\end{proof}

Implicitly, we have been dealing all along with the case when $R$ is the 
sphere spectrum $S$, and we can 
mimic all of the model theoretic work that we have done in that case. The
results of \S12.6 and \S13.1 carry over directly. For $f\colon A\rtarr B$, base change preserves $R$-modules, $(f_!,f^*)$ gives a Quillen adjoint pair relating the categories of $R$-modules over $A$ and over $B$, and we obtain a Quillen equivalence if $f$ is a  $q$-equivalence. If $f$ is a bundle with CW fibres, we obtain a Quillen pair $(f^*,f_*)$, and we can apply the triangulated category version of Brown representability to construct a right adjoint $f_*$ in general. 
However, we do not know how to generalize the rest of Chapter 13 to the module context since we have not worked out a theory of excellent $R$-modules with an accompanying excellent $R$-module approximation functor.  In view of the retreat to prespectra with their primitive handicrafted smash products in that theory,
it seems unlikely to us that any such construction can be expected.

We also have the notion of a right $R$-module over a nonparametrized ring spectrum $R$. If $M$ and $N$ are right and left $R$-modules over $A$ and
$B$ and $L$ is a left $R$-module over $A\times B$, then we define 
spectra $M\barwedge_R N$ over $A\times B$ and $\bar{F}_R(N,L)$ over $A$
by the usual coequalizer
\[\xymatrix{M\barwedge R\barwedge N \ar@<.5ex>[r]\ar@<-.5ex>[r] & M\barwedge N \ar[r]& M\barwedge_R N}\]
and equalizer
\[\xymatrix{\bar{F}_R(N,L)\ar[r] & \bar{F}(N, L) \ar@<.5ex>[r]\ar@<-.5ex>[r] & \bar{F}(R\barwedge N, L).}\]
If $R$ is commmutative, then $M\sma_R N$ and $F_R(N,L)$ are naturally $R$-modules. 

We have good homotopical control over these external constructions, as in
Propositions \ref{Boxcof2} and \ref{Boxcof2too}.  For example, if we take $A=*$, then we have good homotopical
control over the smash product spectrum $M\sma_R N$ over $B$ and the non-parametrized function spectrum $F_R(N,L)$, where $M$ is a non-parametrized
right $R$-module and $N$ and $L$ are left $R$-modules over $B$. However, if we take $A=B$ and internalize $M\barwedge_R N$ along the diagonal $\Delta\colon B\rtarr B\times B$ by setting $M\sma_R N=\Delta^* M\barwedge_R N$ and $F_R(M,N)=\bar{F}_R(M,\Delta_*N)$, we lose homotopical control.  

Similarly, when $R$ is commutative, $R\sM_B$ has the structure of a closed symmetric monoidal category, and that allows us to define (commutative) $R$-algebras over $B$ to be (commutative) monoids in $R\sM_B$. However, because of the lack of homotopical control, in the absence of the theory of Chapter 13, we cannot give the categories of $R$-algebras and of commutative $R$-algebras over $B$ model structures.

\begin{rem}
Although we have not pursued the idea, it seems highly likely that there are interesting examples of rings and modules that allow varying base spaces and
are defined in terms of the external smash product.  
For example, one might consider $G$-spectra $R_n$ over $B^n$ with products $R_m\barwedge R_n\rtarr R_{m+n}$, or one might consider ``globally defined'' parametrized ring spectra $R$ consisting of spectra $R_B$ over $B$ for all $B$ together with appropriate products $R_A\barwedge R_B\rtarr R_{A\times B}$. The
$R_B$ would in particular be module spectra over the nonparametrized ring
spectrum $R_*$. As in the nonparametrized theory, one must use the positive stable model structures to study such ring objects model theoretically when
$R_*$ is commutative. The essential point is that the external smash product 
is sufficiently well-behaved homotopically that there is no obstacle to 
passage from point-set level constructions to homotopy category level conclusions.
\end{rem}

\section{Change of universe}\mylabel{universesec}

Recall that $G$-spectra over $B$ are defined in terms of a chosen collection $\sV$ of representations of $G$. As usual in equivariant stable homotopy theory, we must introduce functors that allow us to change the collection $\sV$. Such functors are usually referred to as ``change of universe'' functors, since $\sV$ is often given as the collection $\sV(U)$ of all representations that embed up to isomorphism in a given $G$-universe $U$. It is however often convenient to restrict $\sV$ to be a cofinal subcollection of $\sV(U)$ that is closed under direct sums, and when we dealt with excellent prespectra it became essential to restrict $\sV$ further to a countable cofinal sequence of expanding representations in $U$. In both cases it is usual to insist that the trivial representation $\bR$ is included in $\sV$. In order to deal with the change functors in all of the above cases at once, we adopt a slightly different approach from the one that was used in \cite[V.\S1]{MM}. We then explain how it specializes to the more explicit approach given there.

Let $G\sS_B^{\sV}$ denote the category of $G$-spectra over $B$ indexed on $\sV$. If $\sV$ is not closed under direct sums, then we are thinking of $G\sS_B^\sV$ as the restriction of the diagram category corresponding to $G\sS_B^{\sV'}$,
where $\sV'$ is the closure of $\sV$ under sums, as discussed in \myref{indexingreps}. 

Let $i\colon \sV\subset \sV'$ be the inclusion of one collection of representations in another. Thinking of parametrized spectra as 
diagram ex-spaces, we see that the evident forgetful functor
\[i^*\colon G\sS^{\sV'} \rtarr G\sS^\sV\]
has a left adjoint $i_*$ given by the prolongation, or expansion of universe, functor
\[(i_*X)(V') = \sJ^{\sV'}_G(-,V') \otimes_{\sJ^\sV_G} X.\]
Such prolongation functors are discussed in detail in \cite[I\S3]{MMSS} 
and \cite[I\S2]{MM}. By \cite[I.2.4]{MM}, the unit $\text{Id}\rtarr i^*i_*$ of
the adjunction is a natural isomorphism.

We have more concrete descriptions of the functor $i_*$ when $\sV$ consists of a cofinal sequence of representations in some universe $U$. Recall that $\sJ_G^\sV(V,V)$ is the orthogonal group $O(V)$ with a disjoint base point.

\begin{lem}\mylabel{describei*}
If $\sV=\{V_i\}\subset \sV'$ is a countable expanding sequence in some $G$-universe $U$, then
\[(i_*X)(V')\cong \sJ^{\sV'}_G(V_i,V')\sma_{O(V_i)} X(V_i)\]
where $i$ is the largest natural number such that there is a linear isometry $V_i\rtarr V'$.
\end{lem}

 \begin{proof}
  The  forgetful functor $i^*$ is restriction along a functor $\io\colon
   \sJ_G^{\sV}\rtarr \sJ_G^{\sV'}$ and $(i_*X)(V')$ is constructed as the
   coequalizer of the pair of parallel maps
   \[\xymatrix{\bigvee_{j,k}                   \sJ^{\sV'}_G(V_j,V')\sma_B
   \sJ^\sV_G(V_k,V_j)\sma_B    X(V_k)    \ar@<.5ex>[r]\ar@<-.5ex>[r]    &
   \bigvee_j \sJ^{\sV'}_G(V_j,V')\sma_B X(V_j)}\]
   given  by  composition  in  $\sJ_G^{\sV'}$  and by the evaluation maps
   associated  to the diagram $X$. A cofinality argument that is easily
   made precise by use of the explicit description of the category
   $\sJ^{\sV'}_G$  given  in  \cite[II.\S4]{MM} shows that the above
   coequalizer agrees with the coequalizer of the subdiagram
   \[\xymatrix{\sJ^{\sV'}_G(V_i,V')\sma_B \sJ^\sV_G(V_i,V_i)\sma_B X(V_i)
   \ar@<.5ex>[r]\ar@<-.5ex>[r] & \sJ^{\sV'}_G(V_i,V')\sma_B X(V_i).}\]
   This coequalizer is the space that we have denoted by                                  $\sJ^{\sV'}_G(V_i,V')\sma_{O(V_i)} X(V_i)$.
   \end{proof}

\begin{rem}
The argument above works in the same way for prespectra. It gives 
the conclusion that, for parametrized prespectra $X$ in $G\sP_B^\sV$,
\[(i_*X)(V)\cong \Sigma^{V-V_i}_B X(V_i).\]
\end{rem}

\begin{rem}\mylabel{uniMM}
Assume that $\sV$ and $\sV'$ are closed under finite sums and contain the trivial representation. We can then define the change of universe functors
\[I^\sV_{\sV'}=i_*i'^* \colon G\sS_B^{\sV'} \rtarr G\sS_B^\sV\]
where $i\colon \{\bR^n\} \subset \sV$ and $i'\colon \{\bR^n\} \subset \sV'$. Explicitly 
\[(I^\sV_{\sV'} X)(V)\cong \sJ_G^\sV(\bR^n,V)\sma_{O(n)} X(\bR^n).\]
This is the definition given in \cite[V.1.2]{MM}. These change of universe functors $I^\sV_{\sV'}$ are exceptionally well behaved on the point set level
and agree with those we are using when $\sV\subset \sV'$. They are symmetric monoidal equivalences of categories. For collections $\sV$, $\sV'$ and $\sV''$, they satisfy
\[I_{\sV'}^{\sV}\com \SI^{\sV'}_B \cong \SI^{\sV}_B,\qquad
I_{\sV'}^{\sV} \com I_{\sV''}^{\sV'} \cong I_{\sV''}^{\sV},\qquad
I_{\sV}^{\sV} \cong \text{Id}.\]
Moreover, $I^\sV_{\sV'}$ is continuous and commutes with smash products with 
ex-spaces. In particular, it is homotopy preserving and therefore induces equivalences of the classical homotopy categories. Unfortunately, however, the functors $I^\sV_{\sV'}$ are as poorly behaved on the homotopy level as they are
well behaved on the point set level. They do not preserve either level $q$-equivalences or $s$-equivalences in general and the point set level relations above do not descend to the model theoretic homotopy categories that we are interested in. Furthermore, these functors $I^\sV_{\sV'}$ do not exist if $\sV$ is a cofinal expanding sequence. We shall therefore not make much use of them.
\end{rem}

Returning to our full generality, let $i\colon \sV\subset \sV'$. The adjoint pair $(i_*,i^*)$ has good homotopical properties.

\begin{thm}\mylabel{change1} 
Let $i\colon \sV\subset \sV'$. Then $i^*$ preserves level $q$-equivalences, level $qf$-fibrations, $s$-fibrations, and $s$-acyclic $s$-fibrations. Therefore $(i_*,i^*)$ is a Quillen adjoint pair in the level $qf$-model structure and in the $s$-model structure. Moreover, $i_*$ on homotopy categories is symmetric monoidal. If $\sV$ is cofinal in $\sV'$, then $i^*$ creates the weak equivalences and $(i_*,i^*)$ is a Quillen equivalence.
\end{thm}

\begin{proof}  It is clear from its levelwise definition that $i^*$ 
preserves level $q$-equi\-va\-lences and level $qf$-fibrations. It follows
that its left adjoint $i_*$ preserves $s$-cofibrations and level acyclic
$s$-cofibrations. This in turn implies that $i^*$ preserves $s$-acyclic
$s$-fibrations, since those are the maps that satisfy the RLP with respect
to the $s$-cofibrations. The levelwise description of $s$-fibrations in
\myref{RLPL} implies that $i^*$ preserves $s$-fibrations. The last statement follows from the definition of homotopy groups and the fact that the unit $\text{id}\rtarr i^*i_*$ is an isomorphism.  The functor $i_*$
commutes with $\barwedge$ on the point set level, by \cite[I.2.14]{MM}, 
and this commutation relation descends directly to homotopy categories. 
Applying \myref{chvschuni} below to the diagonal map of $B$, it follows 
that the derived functor $i_*$ is symmetric monoidal.
\end{proof}

We have constructed the change of universe functors on both the spectrum and prespectrum level and they are compatible with the restriction functors $\bU$. However, in order to make use of excellent parametrized prespectra, we must restrict to parametrized prespectra indexed on cofinal sequencess $j\colon \sW\subset \sV$ and $j'\colon \sW'\subset \sV'$ of indexing representations in the given universes $U\subset U'$. But then there need not be an induced inclusion $i\colon \sW\subset \sW'$. We therefore also define change of 
universe functors for prespectra indexed on cofinal sequences.

\begin{defn} Let $i\colon \sV\subset \sV'$ and choose cofinal sequences $\sW=\{V_i\}$ and $\sW'=\{V_i'\}$ in $\sV$ and $\sV'$ such that $V_{i+1}=V_i\oplus W_i$ and $V_i'=V_i\oplus Z_i$, where 
$Z_{i+1}= Z_{i}\oplus W_i'$ and thus $V_{i+1}' = V_i'\oplus W_i\oplus W_i'$. Define a pair of adjoint functors
\[\xymatrix{G\sP_B^\sW \ar@<.5ex>[r]^{\bar\imath_*} & G\sP_B^{\sW'} \ar@<.5ex>[l]^{\bar\imath^*}}\]
by setting
\[(\bar\imath_*X)(V_i')=\Sigma_B^{Z_i}X(V_i)\qquad\text{and}\qquad
(\bar\imath^*Y)(V_i)=\Omega_B^{Z_i}Y(V_i').\]
The structure maps are induced from the given structure maps in the evident way.
\end{defn}

\begin{prop}\mylabel{chunicomp}
The pair $(\bar\imath_*, \bar\imath^*)$ is a Quillen adjoint pair with respect to both the level $qf$-model structure and the stable model structure.
The following diagram commutes when the vertical arrows point in the same direction.
\[\xymatrix{\Ho G\sP^\sW_B \ar@<.5ex>[d]^{\bar\imath_*} 
& \Ho G\sP^\sV_B \ar@<.5ex>[d]^{i_*}\ar[l]_-{j^*}\\
\Ho G\sP_B^{\sW'} \ar@<.5ex>[u]^{\bar\imath^*} 
& \Ho G\sP_B^{\sV'} \ar@<.5ex>[u]^{i^*}\ar[l]^-{(j')^*}}\]
\end{prop}

\begin{proof}
This is clearly a Quillen adjunction in the level $qf$-model structure, and to show that it is a Quillen adjunction in the stable model structure it therefore suffices to verify the condition of \myref{RLPL}. The homotopy pullback \ref{OMpb} associated to the pair $(V_i,W_i)$ and an $s$-fibration $f\colon X\rtarr Y$ is still a homotopy pullback after we apply $\Omega_B^{Z_i}$ to it and displays the required diagram \ref{OMpb} for the map $\bar\imath^*f$. We have that 
\[(\bar\imath_*j^*X)(V_i')=\Sigma_B^{Z_i}X(V_i)
\cong \Sigma_B^{V_i'-V_i}X(V_i)=((j')^*i_*X)(V_i')\]
and this point set level isomorphism descends to homotopy categories since the functors $j^*$ and $(j')^*$ preserve all $s$-equivalences. The adjoint structure maps of $X\in G\sP_B^{\sV'}$ induce maps
\[(j^*i^*X)(V_i) = X(V_i)\rtarr \Omega_B^{Z_i}X(V_i') = (\bar\imath^*(j')^*X)(V_i).\]
When $X$ is $s$-fibrant, its adjoint structure maps are level $q$-equivalences, and we thus obtain an equivalence $j^*i^*\simeq \bar\imath^*(j')^*$ on homotopy categories. 
\end{proof}

On the point-set level, we have the following commutation relations
between change of universe functors and change of base functors.

\begin{lem}\mylabel{pointsetchchuni}
Let $i\colon \sV\subset \sV'$ and let $f\colon A\rtarr B$ be a $G$-map. Then $i^*$ commutes up to natural isomorphism with the change of base functors $f_!$, $f^*$, and $f_*$, and $i_*$ commutes up to natural isomorphism with $f_!$ and $f_*$. 
\end{lem}

\begin{proof}
The first statement is clear from the levelwise constructions of the base change functors, and the second statement follows by conjugation since $i_*$, $f_!$,
and $f^*$ are left adjoints of $i^*$, $f^*$, and $f_*$.
\end{proof}

The missing relation, $i_*f^*\iso f^*i_*$, would hold with the alternative
point-set level definitions of \myref{uniMM}, where $i^*$ and $i_*$ are inverse equivalences.  However, these are point-set level relationships that need not descend to model theoretic homotopy categories.  With our preferred definition
of $i_*$ in terms of prolongation, the 
following result shows that $i_*f^*\htp f^*i_*$ on homotopy
categories even though we need not have an isomorphism on the 
point-set level.

\begin{prop}\mylabel{chvschuni}
Let $i\colon \sV\subset \sV'$ and let $f\colon A\rtarr B$ be a $G$-map. Then there are natural equivalences of derived functors 
\[i^*f^* \simeq f^*i^*,\quad
i_*f_! \simeq f_!i_*,\quad
i_*f^* \simeq f^*i_*,\quad
i^*f_* \simeq f_*i^*,\quad
i^*f_! \simeq f_!i^*\]
in the relevant homotopy categories.
\end{prop}

\begin{proof} 
The first two equivalences are clear since we are commuting Quillen right adjoints and their corresponding Quillen left adjoints. The fourth will follow by adjunction from the third. If $f$ is a homotopy equivalence, then $f^* \htp (f_!)^{-1}$ and in this case the third follows from the second and the fifth from the first. Factoring $f$ as the composite of an $h$-fibration and a homotopy equivalence, we see that the third will hold in general if it holds when $f$ is an $h$-fibration. Similarly, factoring $f$ as the composite of an $h$-cofibration and a homotopy equivalence, we see that the fifth will hold in general if it holds when $f$ is an $h$-cofibration.

Further, for the third equivalence, it suffices to show that $\bar\imath_*f^*\simeq f^*\bar\imath_*$ since \myref{chunicomp} then gives that
\[i_*f^*\simeq i_*j_*j^*f^*\simeq (j')_*\bar\imath_*f^*j^*
\simeq (j')_*f^*\bar\imath^*j^*\simeq f^*(j')_*(j')^*i^*\simeq f^*i^*.\]
Similarly, for the fifth equivalence, it suffices to show that $\bar\imath^*f_!\simeq f_!\bar\imath^*$, for then
\[i^*f_!\simeq i^*(j')_*(j')^*f_!\simeq j_*\bar\imath^*f_!(j')^*
\simeq j_*f_!\bar\imath^*(j')^*\simeq f_!(j')_*(j')^*i^*\simeq f_!i^*.\]

We have reduced the proof of the third equivalence to the situation when $f$ is an $h$-fibration and $i_*$ is replaced by $\bar\imath_*$. The functor $f^*$ preserves excellent prespectra over $B$, but we must apply $T$ to $\bar\imath_*$ before passing to homotopy categories. As in the proof of \myref{Mackeymore}, since $f$ is assumed to be an $h$-fibration 
we have a natural homotopy equivalence $\mu\colon Tf^*\rtarr f^*T$
in our categories indexed on $\sW$ or on $\sW'$.  Therefore
\[T\bar\imath_*f^*\cong  Tf^*\bar\imath_*\htp f^*T\bar\imath_*.\]

Similarly, we have reduced the proof of the fifth equivalence to the situation when $f$ is an $h$-cofibration and $i^*$ is replaced by $\bar\imath^*$. Then $f_!$ preserves level $h$-equivalences, and so does $\bar\imath^*$ since it preserves level $q$-equivalences and preserves objects whose total spaces are of the homotopy types of $G$-CW complexes. Since $T$ takes zig-zags of level $h$-equivalences to homotopy equivalences,
\[\xymatrix{Tf_!T\bar\imath^* \ar@{<->}[r]^-{\htp} & Tf_!\bar\imath^*\cong T\bar\imath^*f_! \ar@{<->}[r]^-{\htp} & T\bar\imath^*Tf_!}\]
displays a zig-zag of homotopy equivalences showing that $f_!\bar\imath^*\simeq \bar\imath^*f_!$.
\end{proof}

\section{Restriction to subgroups}

Let $\tha\colon G'\rtarr G$ be a homomorphism and let $\tha^*\sV$ be the collection of $G'$-representations $\tha^*V$ for $V\in \sV$, where $\sV$
is our chosen collection of indexing $G$-representations.  We have
implicitly used the following result in our earlier results on change 
of groups.

\begin{prop}\mylabel{grprestrrQa} 
The functor $\tha^*\colon G\sS_B \rtarr G'\sS^{\tha^*\sV}_{\tha^*B}$ preserves level $q$-equi\-va\-lences, level $qf$-fi\-bra\-tions, $s$-fi\-bra\-tions, and $s$-equivalences provided that the model structures are defined with respect
to generating sets $\sC_{G}$ and $\sC_{G'}$ of $G$-cell complexes and $G'$-cell complexes such that $\tha_!C = G\times_{G'} C \in \sC_G$ for $C\in\sC_{G'}$.
\end{prop}

\begin{proof} Since $(\tha^*A)^H=A^{\tha^*(H)}$ for a $G$-space $A$ and 
a subgroup $H$ of $G'$, this is clear from the definitions of homotopy 
groups and from the characterizations of fibrations given in \myref{qffibdef}
and \myref{RLPL}. Note in particular that $\tha^*$ preserves the level 
$qf$-fibrant approximations that are used in the definition of the stable 
homotopy groups.
\end{proof}

For the remainder of this section fix a subgroup $H$ of $G$ and consider the inclusion $\io\colon H\subset G$.  For an $H$-space $A$, we simplify notation by letting $H\sS_A^{\sV}$ denote the category of $H$-spectra over $A$ indexed on $\io^*\sV$.  Clearly, we then have the restriction of action functor
$$\io^*\colon G\sS_B^{\sV}\rtarr H\sS_{\io^*B}^{\sV}.$$ 
For $i\colon \sV\subset \sV'$, we have $\io^*i^* = i^*\io^*$ since with 
either composite we are just restricting from the representations in $\sV'$ 
to the representations in $\sV$ and viewing all $G$-spaces in sight as $H$-spaces. 
 
When $\sV = \sV(U)$ for a $G$-universe $U$, there is a quibble here (as was
discussed in \cite[V.10]{MM}).  We are using $\io^*\sV$ as the corresponding indexing collection for $H$. However, if $V$ is an irreducible representation 
of $G$, $\io^*V$ is generally not an irreducible representation of $H$ and we should expand $\io^*\sV$ to include all representations that 
embed up to isomorphism in $\io^*U$ to fit the definitions into our usual
framework.  However, there is a change of universe functor associated to
the inclusion $i\colon \iota^*\sV(U)\subset \sV(\iota^*U)$ that fixes this. The
functor $i^*$ preserves all $s$-equivalences and descends to an equivalence 
on homotopy categories. We can and should use these rectifications when restricting to $H$-spectra over $\io^*B$ for a fixed chosen $H$.

\begin{rem} Consider passage to fibers and recall \myref{FibadQtoo}.
\begin{enumerate}[(i)]
\item Applied to inclusions of orbits, \myref{chvschuni}
implies that the functors $i^*$ for $i\colon \sV\subset \sV'$ are 
compatible with passage to fibers, in the sense that 
$$(i^*X)_b \iso i^*(X_b)\ \ \text{for}\ \ b\in B,$$
where $i^*$ on the right is the change of universe functor on 
$G_b$-spectra. 
\item When $\sV = \sV(U)$, we can view the fiber functor
\[(-)_b\colon G\sS_B \rtarr G_b\sS\]
as landing in spectra indexed on $\sV(\iota^* U)$, $\io\colon G_b\rtarr G$ by composing with $i_*$ for $i\colon \iota^*\sV(U)\subset \sV(\iota^*U)$. However, these change of universe functors must be used with caution since they are not compatible as $b$ and therefore $G_b$ vary.
\end{enumerate}
\end{rem}

Recall from Propositions \ref{Lchanges2} and \ref{Symmoni} that the equivalence
of categories $(\iota_!, \nu^*\iota^*)$ between $H\sS_A$ and $G\sS_{\io_!A}$
induces a closed symmetric monoidal equivalence of categories between 
$\text{Ho}H\sS_A$ and $\text{Ho}G\sS_{\io_!A}$. By \myref{LishriekCor2}, 
we can interpret the restriction
functor $\io^*\colon \text{Ho}G\sS_B\rtarr \text{Ho}H\sS_{\io^*B}$ as the
composite of base change $\mu^*$ along $\mu\colon \io_!\io^*B\rtarr B$ and this
equivalence applied to $A = \io^*B$.  The following spectrum level analogue of \myref{changerel} gives compatibility relations between change of
base functors and these results on change of groups. 

\begin{prop}\mylabel{substitute}
Let $f\colon A\rtarr \io^*B$ be a map of $H$-spaces and $\tilde{f}\colon
\io_! A\rtarr B$ be its adjoint map of $G$-spaces.  Then the following
diagrams commute up to natural isomorphism, where 
$\mu\colon \io_!\io^*B\rtarr B$ and $\nu\colon A\rtarr \io^*\io_!A$ 
are the counit and unit of the adjunction $(\io_!,\io^*)$.
\[\xymatrix{
G\sS_{\io_!A} \ar[r]^-{\tilde{f}_!} & G\sS_B \\
H\sS_{A} \ar[r]_-{f_!}\ar[u]^{\io_!}  &
H\sS_{\io^*B}\ar[u]_{\mu_!\com\io_!}}
\quad \ \
\xymatrix{
G\sS_B \ar[r]^-{\tilde{f}^*} \ar[d]_{\io^*}
& G\sS_{\io_!A} \ar[d]^{{\nu}^*\com\io^*}\\
H\sS_{\io^*B} \ar[r]_-{f^*} & H\sS_{A}}\]
These diagrams descend to natural equivalences of composites of 
derived functors on homotopy categories.
\end{prop}
\begin{proof}
The point set level diagrams commute by \myref{changerel}, applied
levelwise. The left diagram is one of Quillen left adjoints and the
right diagram is one of Quillen right adjoints, by Propositions 
\ref{Qad1too} and \ref{Lchanges2} and \myref{LishriekCor2}. 
\end{proof}
 
We now define a parametrized fixed point functor associated to the inclusion $\iota\colon H\rtarr G$. Its target is a category of nonequivariant parametrized spectra. In the next section we will consider a fixed point functor that takes values in a category of parametrized $WH$-spectra, where $WH=NH/H$ is the Weyl group.

Write $G\sS^{\text{triv}}_B$ for $G$-spectra over $B$ indexed on trivial representations.  These are ``naive'' parametrized $G$-spectra. As usual, to define fixed point spectra, we must change to the trivial universe before taking fixed points levelwise. Thus let $\sV^G = \{V^G \mid V\in \sV\}$. It is contained in $\sV$ if $\sV=\sV(U)$ for some universe $U$.

\begin{defn}
The \emph{$G$-fixed point functor} $(-)^G\colon G\sS_B\rtarr \sS_{B^G}$ is the composite of $i^*$, $i\colon \sV^G\subset \sV$, and levelwise passage to fixed points. For a subgroup $H$ of $G$ the \emph{$H$-fixed point functor} 
$(-)^H\colon G\sS_B\rtarr \sS_{B^H}$ is the composite of $\io^*$, $\io\colon H\subset G$, and $(-)^H$.
\end{defn}

Since the homotopy groups of a level $qf$-fibrant $G$-spectrum $X$ over $B$ are the homotopy groups $\pi_q^H(X_b)$, we see from the nonparametrized analogue \cite[V.3.2]{MM} that these are then the homotopy groups of $X^H$.  Recall
in particular that the $s$-fibrant $G$-spectra over $B$ are the $\OM$-$G$-spectra over $B$, which are level $qf$-fibrant.  Therefore, for all subgroups
$H$ of $G$, the homotopy groups of a parametrized $G$-spectrum $X$ are the nonequivariant homotopy groups of the nonequivariant spectra $X^H$, provided
that $(-)^H$ is understood to mean the derived fixed point functor. 

On the point-set level, the functor $(-)^G$ is a right adjoint. Thinking of the homomorphism $\epz\colon G\rtarr e$ to the trivial group, let $\epz^*\colon \sS_A\rtarr G\sS^{\text{triv}}_{\epz^*A}$ be the functor that sends spectra 
over a space $A$ to $G$-trivial $G$-spectra over $A$ regarded as a $G$-trivial $G$-space. The following result is immediate by passage to fibers from its nonparametrized special case \cite[V.3.4]{MM}.  Let $\sA\!\ell\ell$ denote
the collection of all representations of $G$.

\begin{prop}\mylabel{fixad} 
Let $A$ be a space. Let $Y$ be a naive $G$-spectrum over $\epz^*A$ and $X$ be a spectrum over $A$. There is a natural isomorphism
$$G\sS^{\text{triv}}_{\epz^*A}(\epz^*X, Y)\iso \sS_A(X, Y^G).$$
For (genuine) $G$-spectra $Y$ over $\epz^*A$, there is a natural isomorphism
$$G\sS_{\epz^*A}(i_*\epz^*X,Y)\iso \sS_A(X,(i^*Y)^G),$$
where $i\colon \text{triv}\subset \sA\!\ell\ell$.
Both of these adjunctions are given by Quillen adjoint pairs relating the respective level and stable model structures.
\end{prop}

Returning to $G$-spaces $B$ and comparing \myref{FVs} with the proof of \cite[V.3.5-3.6]{MM}, we obtain the following curious results.

\begin{prop}\mylabel{fixcof} 
For a representation $V$ and an ex-$G$-space $K$, we have that $(F_V K)^G = *_{B^G}$ unless $G$ acts trivially on $V$, when $(F_VK)^G \iso F_V(K^G)$ as a nonequivariant spectrum over $B^G$. The functor $(-)^G$ preserves $s$-cofibrations, but it does not preserve acyclic $s$-cofibrations.
\end{prop}

\begin{cor}\mylabel{fixsusp}
For ex-$G$-spaces $K$,
$$(\SI^{\infty}_BK)^G\iso \SI^{\infty}_B(K^G).$$ 
\end{cor}

This isomorphism of spectra over $B^G$ does \emph{not} descend to the homotopy category $\Ho G\sS_{B^G}$. The reader is warned to consult \cite[V\S3]{MM} for the meaning of these results.  There is also an analogue of the comparison between $G$-fixed points and smash products in \cite[V.3.8]{MM}, but only when $B=B^G$ and only with good behavior with respect to cofibrant objects when external smash products are used. We shall not state the result formally.

\section{Normal subgroups and quotient groups}

We now turn to quotient homomorphisms and associated orbit and fixed point functors. The material of this section generalizes a number of results from \S2.4, \S7.3, and \S9.5 to the level of parametrized spectra.

Just as we have been using $\io$ generically for inclusions of subgroups, we shall use $\epz$ generically for quotient homomorphisms. In particular, for an inclusion $\io\colon H\subset G$, we let $WH = NH/H$, where $NH$ is the normalizer of $H$ in $G$, and we have the quotient homomorphism $\epz\colon NH\rtarr WH$. We can study this situation by first restricting from $G$ to $NH$, thus changing the ambient group.  Therefore, there is no loss of generality if we focus attention on a normal subgroup $N$ of $G$ with quotient group $J=G/N$, as we do throughout this section.

\begin{defn}\mylabel{defnchN}
Let $G\sS^{\text{$N$-triv}}_B$ be the category of $G$-spectra over $B$ indexed on the $N$-trivial representations of $G$. Regard representations of $J$ as $N$-trivial representations of $G$ by pullback along $\epz\colon G\rtarr J$.
For a $J$-space $A$, define
\[\epz^*\colon  J\sS_A\rtarr G\sS^{\text{$N$-triv}}_{\epz^*A}\]
levelwise by regarding ex-$J$-spaces over $A$ as $N$-trivial $G$-spaces over $\epz^*A$. For a $G$-space $B$, define 
\[(-)/N\colon G\sS^{\text{$N$-triv}}_B \rtarr J\sS_{B/N}
\qquad\text{and}\qquad
(-)^N\colon  G\sS^{\text{$N$-triv}}_B \rtarr J\sS_{B^N}\]
by levelwise passage to orbits over $N$ and to $N$-fixed points.
\end{defn}

\begin{lem}\mylabel{fixedptrQa}
The $N$-fixed point functor $(-)^N$ preserves level $q$-equivalences, level $qf$-fibrations, $s$-fibrations, and $s$-equivalences, provided that the model structures are defined with respect to generating sets $\sC_{G}$ and $\sC_{J}$ of $G$-cell complexes and $J$-cell complexes such that $C/N \in \sC_J$ for $C\in\sC_{G}$.
\end{lem}

\begin{proof}
This is a special case of \myref{grprestrrQa}; it also follows directly from the ex-space level analogue in \myref{fixedptrQa0}, the characterization of $s$-fibrations in \myref{RLPL}, and inspection of the definition of the $s$-equivalences.  
\end{proof}

\begin{prop}\mylabel{factor} Let $j\colon B^N\rtarr B$ be the inclusion 
and $p\colon B\rtarr B/N$ be the quotient map.  Then the following
factorization diagrams commute
\[\xymatrix{G\sS^{\text{$N$-triv}}_B \ar[d]_{p_!}\ar[r]^{(-)/N} 
& J\sS_{B/N} \\
G\sS^{\text{$N$-triv}}_{B/N} \ar[ur]_{(-)/N} }
\qquad\text{and}\qquad
\xymatrix{G\sS^{\text{$N$-triv}}_B \ar[d]_{j^*}\ar[r]^{(-)^N} 
& J\sS_{B^N} \\
G\sS^{\text{$N$-triv}}_{B^N} \ar[ur]_{(-)^N}}\]
and they descend to natural equivalences on homotopy categories
\[(p_!X)/N \simeq X/N \qquad\text{and}\qquad (j^*X)^N\simeq X^N\]
for $X$ in $\Ho G\sS_B^{\text{$N$-triv}}$.
The following adjunction isomorphisms follow.
\begin{enumerate}[(i)]
\item For $Y\in G\sS^{\text{$N$-triv}}_B$ and $X\in J\sS_{B/N}$,
$$J\sS_{B/N}(Y/N,X)\iso G\sS^{\text{$N$-triv}}_B(Y,p^*\epz^*X).$$
\item For $Y\in G\sS^{\text{$N$-triv}}_B$ and $X\in J\sS_{B^N}$,
$$G\sS^{\text{$N$-triv}}_B(j_!\epz^*X,Y)\iso J\sS_{B^N}(X,Y^N).$$
\item For (genuine) $G$-spectra $Y\in G\sS_B$ and $X\in J\sS_{B^N}$, 
$$G\sS_B(i_*j_!\epz^*X,Y)\iso J\sS_{B^N}(X,(i^*Y)^N),$$
where $i\colon\text{triv}\subset\sA\!\ell\ell$.
\end{enumerate}
All of these adjunctions are Quillen adjoint pairs with respect to both 
the level and the stable model structures and so descend to homotopy categories.
\end{prop}

\begin{proof}
The factorizations follow from the ex-space level analogue \myref{factor0}. The statement about Quillen adjunctions holds since $(-)^N$, $\epsilon^*$ and $i^*$ preserve level $q$-equivalences, level fibrations, $s$-equivalences and level $s$-fibrations, by \myref{fixedptrQa}, \myref{grprestrrQa} and \myref{change1}.
\end{proof}

The behavior of the orbit and fixed point functors with respect to base change is recorded in the following result.

\begin{prop}\mylabel{chvsfixorbit}
Let $f\colon A\rtarr B$ be a map of $G$-spaces. Then the following diagrams commute up to natural isomorphism
{\small\[\xymatrix{
G\sS_A^{\text{$N$-triv}} \ar[r]^-{f_!} \ar[d]_{(-)/N} 
& G\sS^{\text{$N$-triv}}_B \ar[d]^{(-)/N}\\
J\sS_{A/N} \ar[r]_-{(f/N)_!}  & J\sS_{B/N}}
\;\;
\xymatrix{
G\sS^{\text{$N$-triv}}_B \ar[r]^-{f^*} \ar[d]_{(-)^N}
& G\sS^{\text{$N$-triv}}_A \ar[d]^{(-)^N}\\
J\sS_{B^N} \ar[r]_-{(f^N)^*} & J\sS_{A^N}}
\;\;
\xymatrix{
G\sS_{A}^{\text{$N$-triv}} \ar[r]^-{f_!} \ar[d]_{(-)^N} 
& G\sS^{\text{$N$-triv}}_B \ar[d]^{(-)^N}\\
J\sS_{A^N} \ar[r]_-{(f^N)_!}  & J\sS_{B^N}}\]}
and they descend to the following natural equivalences on homotopy categories
\[(f_!X)/N \simeq (f/N)_! (X/N),
\quad
(f^*X)^N \simeq (f^N)^*(Y^N),
\quad
(f_!X)^N \simeq (f^N)_!(X/N)\]
for $X\in\Ho G\sS^\text{$N$-triv}_A$ and $Y\in \Ho G\sS^\text{$N$-triv}_B$.
\end{prop}

\begin{proof}
The first statement follows levelwise from the ex-space level analogue \myref{fixorbbase}. The proof that it descends to equivalences on homotopy categories is the same as for the ex-space level analogue \myref{orbfixdescend}.
\end{proof}

Specializing to $N$-free $G$-spaces, we obtain a factorization result that is analogous to those in \myref{factor}, but is less obvious. It is a precursor of the Adams isomorphism.

\begin{prop}\mylabel{ouch} 
Let $E$ be an $N$-free $G$-space, let $B = E/N$, and let $p\colon E\rtarr B$ be the quotient map. Then the diagram 
\[\xymatrix{
G\sS_E^{\text{$N$-triv}} \ar[d]_{p_*} \ar[r]^-{(-)/N} & J\sS_B\\
G\sS_B^{\text{$N$-triv}} \ar[ur]_-{(-)^N}}\]
commutes up to a natural isomorphism, and it descends to a natural equivalence
$$X/N\simeq (p_*X)^N$$
in $G\sS_E^{\text{$N$-triv}}$ for $X\in \Ho J\sS_B$. Therefore the left adjoint $(-)/N$ of the functor $p^*\epz^*$ is also its right adjoint. 
\end{prop}

\begin{proof}
The point set level result follows levelwise from the ex-space level result \myref{ouch0}. Since it is an isomorphism between a Quillen left adjoint on the left hand side and a composite of Quillen right adjoints on the right hand side, it descends directly to homotopy categories.
\end{proof}

\part{Duality, transfer, and base change isomorphisms}

\chapter{Fiberwise duality and transfer maps}

\section*{Introduction}
We put the foundations of Part III to use in the two chapters of this last part. Unless otherwise stated, we work in the derived homotopy categories, and all functors should be understood in the derived sense.  For example, we have the derived fiber functor
$$(-)_b\colon \Ho G\sS_B \rtarr \Ho G_b\sS.$$
Since passage to fibers is a Quillen right adjoint, this means that we
replace $G$-spectra $X$ over $B$ by $s$-fibrant approximations before taking
point-set level fibers.  For emphasis, and to make the notation $X_b$ clear
and unambiguous, we may assume that $X$ is $s$-fibrant, but there is no loss
of generality.  A map $f$ in $\Ho G\sS_B$ is an equivalence if and only if 
$f_b$ is an equivalence for all $b\in B$, and that allows us to transer
information back and forth between the parametrized and unparametrized 
homotopy categories with impunity.  Here we use the word ``equivalence'' 
to mean an isomorphism in $\Ho G\sS_B$, and we use the notation $\simeq$ for this relation. We reserve the symbol $\cong$ to mean an isomorphism on the 
point set level.

We have proven that the basic structure enjoyed by the category $G\sS_B$ of parametrized spectra descends coherently to the homotopy category $\Ho G\sS_B$.  In particular, $\Ho G\sS_B$ is closed symmetric monoidal, and
the derived fiber functor is closed symmetric monoidal.  In any symmetric monoidal category, we have standard categorical notions of dualizable and invertible objects.  In \S\ref{sec:fibdual}, we prove the fiberwise duality theorem, which says that a $G$-spectrum $X$ over $B$ is dualizable or invertible if and only 
if each fiber $X_b$ is dualizable or invertible. This allows us to recognize dualizable or invertible $G$-spectra over $B$ when we see them.  

In \S\ref{sec:trfr}, we explain how the fiberwise duality theorem leads to a simple and general conceptual definition of trace and transfer maps with good properties.  To define the transfer, we regard a Hurewicz fibration $p\colon E\rtarr B$ with stably dualizable fibers as a space over $B$. We adjoin a copy of $B$ to obtain a section, and we suspend to obtain a $G$-spectrum over $B$.  It is dualizable since its fibers are dualizable, hence it has a transfer map defined by categorical nonsense. Pushing down to $G$-spectra by base change along the map $r\colon B\rtarr *$, we obtain the transfer map of $G$-spectra $\SI^{\infty}B_+\rtarr \SI^{\infty}E_+$. This construction is a generalization of various earlier constructions of the transfer \cite{BG1, BG2, CG, Clapp, Waner}, most of which restrict to finite dimensional base spaces and are nonequivariant. An essential point is that the homotopy category of $G$-spectra over $B$ is closed symmetric monoidal with a ``compatible triangulation'', in the sense specified in \cite{Tri}.  We defer the proof of the required compatibility relations to \S\ref{sec:comptriang}.  This point implies that our traces and transfers satisfy additivity relations as well as the more elementary standard properties. 

Some of the classical constructions of the transfer work only for bundles, 
but have various properties that are inaccessible to the more general 
construction and are important in calculations. These transfers also 
admit a perhaps more satisfying construction. Rather than relying on duality 
on the level of parametrized spectra, they are obtained by inserting duality
maps for fibers fiberwise into bundles. In the literature, the construction 
again usually requires finite dimensional base spaces and is nonequivariant. 
We give a general conceptual version of this alternative construction in \S\ref{sec:fibtrfr}.

As a first preliminary, in \S\ref{sec:bdlconstr} we show how to insert parametrized spectra fiberwise into the standard construction of equivariant bundles associated to principal bundles. The general construction is of considerable interest nonequivariantly.  The construction on the ex-space 
level is easy enough, but even here many of the properties that we describe 
seem to be new. The construction is likely to have many further applications. The idea is to generalize the standard construction of the bundle of tangents along the fibers of a bundle by replacing the tangent bundle of the fiber by any spectrum over the fiber. In more detail, we consider $G$-bundles $p\colon E\rtarr B$ with fibers $F$.  We allow the structure group $\PI$ and ambient group $G$ to be related by an extension $1 \rtarr \PI\rtarr \GA\rtarr G\rtarr 1$, and we take $F$ to be a $\GA$-space. The bundle $p$ has an associated principal $(\PI;\GA)$-bundle $\pi\colon P\rtarr B$, where $P$ is a $\PI$-free $\GA$-space and $B = P/\PI$.  We show how to construct a $G$-spectrum $P\times_{\PI}X$ over $E$ from a $\GA$-spectrum $X$ over $F$.

As a second preliminary, in \S\ref{sec:PIFree} we develop the theory of
$\PI$-free parametrized $\GA$-spectra.  This is a direct generalization
of the nonparametrized theory and is important in many contexts. In
particular, it will play a role in our proof of the Adams isomorphism
in \S\ref{sec:adams}.

The application to transfer maps in \S\ref{sec:fibtrfr} can be described as follows. When $F$ is dualizable, we have a transfer map $\ta\colon S_{\GA}\rtarr \SI^{\infty}_{\GA}F_+$ of $\GA$-spectra. We insert this into the functor $P\times_{\PI}(-)$ to obtain a map
\[P\times_{\PI}\ta\colon P\times_{\PI} S_{\GA}\rtarr P\times_{\PI}\SI^{\infty}_{\GA}F_+\]
of $G$-spectra over $B$.  Again pushing down to a map of $G$-spectra along $r\colon B\rtarr *$, we obtain the transfer $G$-map $\SI^{\infty}_GB_+\rtarr \SI^{\infty}_GE_+$. This description hides a subtlety.  The construction involves a change of universe functor, and the key point is that this functor is a symmetric monoidal equivalence between categories of parametrized $\PI$-free $\GA$-spectra.  This makes it transparent from the naturality of transfer maps with respect to symmetric monoidal functors that the fiberwise transfer map of 
a bundle agrees with its transfer map as a Hurewicz fibration.

We assume throughout that all given groups $G$ are compact Lie 
groups and all given base $G$-spaces are of the homotopy types 
of $G$-CW complexes.

\section{The fiberwise duality theorem}\label{sec:fibdual}

We characterize the dualizable and invertible $G$-spectra over $B$.  A recent exposition of the general theory of duality in closed symmetric monoidal
categories appears in \cite{Pic}, to which we refer the reader for discussion
of the relevant categorical definitions and arguments. The following theorem is a substantial generalization of various early results of the same nature about ex-fibrations. These are due, for example, to Becker and Gottlieb \cite[\S4]{BG1}, Clapp \cite[3.5]{Clapp}, and Waner \cite[4.6]{Waner}. 

\begin{thm}[The fiberwise duality theorem]\mylabel{bingo}
Let $X$ be an ($s$-fibrant) $G$-spec\-trum over $B$. Then $X$ is dualizable 
(respectively, invertible) if and only if $X_b$ is dualizable (respectively, invertible) as a $G_b$-spectrum for each $b\in B$. 
\end{thm}

\begin{proof}  
By definition, $X$ is dualizable if and only if the natural map
$$\nu\colon  D_BX\sma_B X\rtarr F_B(X,X)$$
in $\Ho G\sS_B$ is an equivalence. Passing to (derived) fibers, 
this holds if and only if the resulting map 
\[\xymatrix{DX_b\sma X_b \simeq (D_BX\sma_B X)_b \ar[r]^-{\nu_b} & F_B(X,X)_b\simeq F(X_b,X_b)}\]
in $\Ho G_b\sS$ is an equivalence for all $b\in B$. By the 
categorical coherence observation \myref{coherence}, the latter map is the corresponding natural map $\nu$ in $\Ho G_b\sS$. Again by definition, that
map is an equivalence if and only if $X_b$ is dualizable.

Similarly, $X$ is invertible if and only if
the evaluation map
$$\text{ev}\colon  D_BX\sma_B X\rtarr S_B$$ 
in $\Ho G\sS_B$ is an equivalence. Passing to (derived) fibers,
this holds if and only if the resulting map
\[\xymatrix{DX_b\sma X_b\simeq (D_BX\sma_B X)_b \ar[r]^-{\text{ev}_b} & (S_B)_b\simeq S}\]
in $\Ho G_b\sS$ is an equivalence for all $b\in B$. Again by
\myref{coherence}, the latter map
is the evaluation map for $X_b$ in $\Ho G_b\sS$, and that map
is an equivalence if and only if $X_b$ is invertible.
\end{proof}

Therefore, to recognize parametrized dualizable and invertible $G$-spectra, it
suffices to recognize nonparametrized dualizable and invertible $G$-spectra.
As we now recall from \cite{FLM}, these are well-understood.
 
Recall that a $G$-space $X$ is dominated by a $G$-space $Y$ if $X$ is a 
retract up to homotopy of $Y$, so that the identity map of $X$ is homotopic
to a composite $X\rtarr Y\rtarr X$. If $Y$ has the homotopy type of a $G$-CW
complex, then so does $X$.  We say that $X$ is finitely dominated
if it is dominated by a finite $G$-CW complex.  This does not imply that 
$X$ has the homotopy type of a finite $G$-CW complex, even when $X$ and all
of its fixed point spaces $X^H$ are simply connected and therefore, since 
they are finitely dominated, homotopy equivalent to finite CW complexes. 

For example, a $G$-space $X$ is a $G$-ENR (Euclidean neighborhood retract)\index{Euclidean neighborhood retract} if it can be embedded as a retract of an open subset of some representation $V$.  Such open subsets
are triangulable as $G$-CW complexes, so $X$ has the homotopy type of a
$G$-CW complex.  A compact $G$-ENR is a retract of a finite $G$-CW complex 
and is thus finitely dominated, but it need not have the homotopy type of a finite $G$-CW complex.  Non-smooth topological $G$-manifolds give examples 
of such non-finite compact $G$-ENRs.

The following result is \cite[2.1]{FLM}. 

\begin{thm}\mylabel{dualG}
Up to equivalence, the dualizable $G$-spectra are the $G$-spectra of the
form $\Sigma^{-V}\Sigma^\infty X$ where $X$ is a finitely dominated based 
$G$-CW complex and $V$ is a representation of $G$.
\end{thm}

\begin{defn}\mylabel{hrep}
A \emph{generalized homotopy representation} $X$ is a finitely dominated based $G$-CW complex such that, for each subgroup $H$ of $G$, $X^H$ is equivalent to a sphere $S^{n(H)}$. A \emph{stable homotopy representation} is a $G$-spectrum of the form $\SI^{-V}\SI^{\infty}X$, where $X$ is a generalized homotopy representation and $V$ is a representation of $G$.
\end{defn}

The following result is \cite[0.5]{FLM}.

\begin{thm}\mylabel{inverG}
Up to equivalence, the invertible $G$-spectra are
the stable homotopy representations.
\end{thm}

Combining results, we obtain the following conclusion about ex-$G$-fibrations.

\begin{thm}\mylabel{spaceFDT} 
Let $E$ be an ex-$G$-fibration over $B$. If each fiber $E_b$ is a finitely
dominated $G_b$-space, then $\SI^{\infty}_BE$ is a dualizable $G$-spectrum 
over $B$. If each $E_b$ is a generalized homotopy representation of $G_b$, 
then $\SI^{\infty}_BE$ is an invertible $G$-spectrum over $B$.
\end{thm}

\begin{proof} 
Since the derived suspension spectrum functor commutes
with passage to derived fibers, by \myref{reassuring}, the derived fiber  $(\SI^{\infty}_BE)_b$ is equivalent to $\SI^{\infty}E_b$. The conclusion 
follows directly from Theorems \ref{bingo}, \ref{dualG}, and \ref{inverG}. 
\end{proof}

In particular, sphere $G$-bundles and, more generally, spherical 
$G$-fibrations over $B$, have invertible suspension $G$-spectra over $B$.

\section{Duality and transfer maps}\label{sec:trfr}

Since the stable homotopy category $\Ho G\sS_B$ is closed 
symmetric mon\-oi\-dal, we have the following generalized trace maps
at our disposal. We state the definition and recall its properties
in full generality, and we then specialize to show how it gives
a simple conceptual definition of the transfer maps associated 
to equivariant Hurewicz fibrations. 

\begin{defn}\mylabel{tracemap}
Let $\sC$ be any closed symmetric monoidal category with unit object $S$. 
For a dualizable object $X$ of $\sC$ with a ``coaction'' map 
$\Delta_X\colon  X\rtarr X\wedge C_X$ for some object $C_X\in\sC$, define the 
\emph{trace}\index{trace} $\tau(f)$ of a self map $f$ of $X$ by the diagram
\[\xymatrix{ S\ar[r]^-\eta \ar[d]_{\tau(f)} & X\wedge DX \ar[r]^-\gamma 
& DX\wedge X \ar[d]^{Df\wedge \Delta_X} \\
C_X & S\wedge C_X\ar[l]^-\cong & DX\wedge X\wedge C_X.\ar[l]^-{\epsilon\wedge 1}}\]
\end{defn}

\begin{rem}
Such a categorical description of generalized trace maps was first given by Dold 
and Puppe \cite{DP}, where they showed that it gives the right framework for trace 
maps in algebra, the transfer maps of Becker and Gottlieb \cite{BG1, BG2}, and 
the fixed point theory of Dold \cite{Dold}. These early constructions of 
transfer maps had finiteness conditions that were first eliminated by Clapp 
\cite{Clapp, CP}. Indeed, she gave an early construction of a parametrized 
stable homotopy category and proved a precursor of our fiberwise duality theorem.  
The equivariant analogue of the attractive space level treatment of Spanier-Whitehead duality given by Dold and Puppe was worked out in \cite{LMS}, and a recent categorical 
exposition of duality has been given in \cite{Pic}.
\end{rem}

Two cases are of particular interest. The first is when $C_X=S$ and 
$\Delta_X$ is the unit isomorphism. Then $\tau(f)$ is called the 
\emph{Lefschetz constant}\index{Lefschetz constant} of $f$ and is denoted by $\chi(f)$; in the special 
case when $f=\text{id}$ it is called the \emph{Euler characteristic}\index{Euler characteristic}
of $X$ and denoted by $\chi(X)$. The second is when $C_X=X$. We then
think of $\Delta_X$ as a diagonal map, and $\tau_X=\tau(\text{id})$ 
is called the \emph{transfer map}\index{transfer map} of $X$. 

\begin{rem}
If $C_X$ comes with a ``counit'' map $\xi\colon  C_X\rtarr S$ such that 
the composite 
$$\xymatrix{X\ar[r]^-{\Delta} & X\wedge C_X \ar[r]^-{\text{id}\wedge \xi} & X}$$ 
is the identity, then $\chi(f) = \xi\circ\tau(f)$ by a little diagram chase. The reason for the terminology ``coaction'' and ``counit'' for the maps $\Delta_X$ and $\xi$ is that in many situations $C_X$ is a comonoid and $\Delta_X$ is a coaction of $C_X$ 
on $X$.
\end{rem}

The following basic properties of the trace are proven in \cite[III\S7]{LMS}
and in \cite{Tri}, where more detailed statements are given. Define a map 
$$(f,\alpha)\colon (X,\Delta_X)\rtarr(Y,\Delta_Y)$$
to be a pair of maps $f\colon X\rtarr Y$ and $\alpha\colon C_X\rtarr C_Y$ such 
that the following diagram commutes.
\[\xymatrix{X\ar[r]^-{\Delta_X}\ar[d]_f & X\wedge C_X \ar[d]^{f\wedge \alpha}\\
Y\ar[r]_-{\Delta_Y} & Y\wedge C_Y}\]

\begin{prop}\mylabel{traceprop}
The trace satisfies the following properties, where $X$ and $Y$ are dualizable
and $\DE_X$ and $\DE_Y$ are given. 
\begin{enumerate}[(i)]
\item{\em (Naturality)} If $\sC$ and $\sD$ are closed symmetric monoidal
categories and \linebreak 
$F\colon\sC\rtarr \sD$ is a lax symmetric monoidal functor
such that $FS_{\sC}\iso S_{\sD}$, then
\[\tau(Ff)=F\tau(f),\]
where $C_{FX} = FC_X$ and $\Delta_{FX}=F\Delta_X$.
\item{\em (Unit property)} If $f$ is a self map of the unit object, then 
$\chi(f)=f$.
\item{\em (Fixed point property)}
If $(f,\alpha)$ is a self map of $(X,\Delta_X)$, then 
$$\alpha\circ \tau(f)=\tau(f).$$
\item{\em (Invariance under retracts)} 
If $X\stackrel{i}{\rtarr} Y \stackrel{r}{\rtarr} X$ is a retract, $f$ is a self map of $X$, and 
$(i,\alpha)$ is a map $(X,\Delta_X)\rtarr (Y,\Delta_Y)$, then 
$$\alpha\circ\tau(f)=\tau(ifr).$$
\item{\em (Commutation with $\sma$)}
If $f$ and $g$ are self maps of $X$ and $Y$, then 
\[\tau(f\wedge g)=\tau(f)\wedge \tau(g),\]
where $\Delta_{X\wedge Y} = (\text{id}\sma\ga\sma\text{id})\circ (\DE_X\sma \DE_Y)$ with
$\ga$ the transposition.
\item{\em (Commutation with $\vee$)}  If $\sC$ is additive 
and $h\colon X\vee Y\rtarr X\vee Y$ induces $f\colon X\rtarr X$ and 
$g\colon Y\rtarr Y$ by inclusion and retraction, then
\[\tau(h)=\tau(f) + \tau(g),\]
where $C_X = C_Y = C_{X\vee Y}$ and $\Delta_{X\vee Y} = \DE_X \vee \DE_Y$.
\item{\em (Anticommutation with suspension)} If $\sC$ is triangulated, then
$$\tau(\Sigma f)=-\tau(f)$$
for all self maps $f$, where $\DE_{\SI X} = \SI \DE_X$. 
\end{enumerate}
\end{prop}

In the triangulated context, there is another and very much deeper property.

\begin{thm}[Additivity]\mylabel{addprop} Let $\sC$ be a closed symmetric monoidal
category with a ``compatible triangulation''.  Let $X$ and $Y$ be dualizable
and let  $\DE_X$ and $\DE_Y$ be given, where $C=C_X=C_Y$. Let $(f,\text{id})$ be a 
map $(X,\DE_X)\rtarr (Y,\DE_Y)$ and extend $f$ to a distinguished triangle
$$\xymatrix{X\ar[r]^-f  & Y \ar[r]^{g} & Z \ar[r]^-h &\SI X.\\}$$
Assume given maps $\ph$ and $\ps$ that make the left square commute in the
first of the following two diagrams.
$$\xymatrix{
X \ar[r]^-f \ar[d]_{\ph}  &  Y \ar[r]^-g \ar[d]^{\ps} 
& Z \ar[r]^-h\ar[d]^{\om} & \SI X \ar[d]^{\SI\ph}\\
X \ar[r]_-f  &  Y \ar[r]_-g & Z \ar[r]_-h & \SI X\\}$$
$$\xymatrix{
X \ar[r]^-f\ar[d]_{\DE_X}  &  Y \ar[r]^-g \ar[d]^{\DE_Y} 
& Z \ar[r]^-h \ar[d]^{\DE_Z} & \SI X \ar[d]^{\SI\DE_X}\\
X\sma C \ar[r]_-{f\sma\text{id}}  &  Y\sma C \ar[r]_-{g\sma\text{id}} 
& Z\sma C \ar[r]_-{h\sma\text{id}} & \SI (X\sma C)\\}$$
Then there are maps $\om$ and $\DE_Z$ such that the diagrams commute and
$$\ta(\psi) = \ta(\om) + \ta(\ph).$$
\end{thm}

The most important case starts with only the distinguished triangle $(f,g,h)$ and concludes with the fundamental additivity relation
$$\ch(Y) = \ch(X) + \ch(Z).$$
The  additivity of traces was studied in \cite[III\S7]{LMS} in the equivariant stable homotopy category, but the proof there is incorrect. A thorough investigation of precisely what is needed to prove the additivity of traces is given in \cite{Tri}, where the axioms for a ``compatible triangulation'' are formulated. These axioms hold in all situations previously encountered in algebraic topology and algebraic geometry. However, the model theoretic method of proof described in \cite{Tri} assumes the usual model theoretic compatibilities, such as the pushout-product axiom of \cite{SS}, and these fail to hold in the present context.  Since the proof of the following result only makes sense by close comparison with the proof in \cite{Tri}, we shall defer it to \S\ref{sec:comptriang}.

\begin{thm} The category $\Ho G\sS_B$ is a closed symmetric monoidal category
with a compatible triangulation.
\end{thm}

With these foundations in place, we can now generalize the 
classical construction of transfer maps. The results above specialize 
to give more information about them than is to be found in the literature. 
If $X$ is a dualizable $G$-spectrum over $B$ with a diagonal map
$\Delta_X\colon X\rtarr X\sma_B X$, then we have the transfer map 
$\tau_X\colon  S_{B}\rtarr X$. 
We shall apply this to suspension $G$-spectra associated to
$G$-fibrations $p\colon E\rtarr B$, but we do not assume that $p$ has a section. 
We need some notation. It has been the custom since the beginnings of algebraic 
topology to use the same letter $E$ for a bundle and for its underlying total space.  
It seems to us that this standard abuse of notation seriously obscures the literature of 
parametrized homotopy theory, and for that reason we shall be very pedantic at this 
point.  

\begin{notn} 
For a $G$-space $E$ over $B$, let $(E,p)_+$ denote the 
ex-$G$-space $E\amalg B$ over $B$, with section at the disjoint copy of $B$.
The usual notation is $E_+$, but we shall reserve that notation for the union
of the total $G$-space $E$ with a disjoint basepoint. Observe that if $p$ is
a Hurewicz $G$-fibration, then $(E,p)_+$ is an ex-$G$-fibration.  Except where
otherwise indicated, we agree to write $r$ for the unique map $B\rtarr *$ for 
any based $G$-space $B$.
\end{notn}

Recall the desription of the base change functors associated to $r$ from \myref{r!ex}. The spectrum level versions of these functors are central to the deduction of results in classical stable homotopy theory from results in parametrized stable homotopy theory. The following observation is particularly relevant.

\begin{lem}\mylabel{lesstrivial} 
For a $G$-map $p\colon E\rtarr B$, thought of as a $G$-space over $B$, \[r_!\SI^{\infty}_B (E,p)_+\simeq \Sigma^\infty E_+,\]
where $r\colon B\rtarr *$. In particular, $r_!S_B \simeq \SI^{\infty}B_+$.
\end{lem}

\begin{proof} We have $r_!\SI^{\infty}_B\htp \SI^{\infty}r_!$. This is a 
commutation relation between Quillen left adjoints, and the corresponding
commutation relation for right adjoints holds since 
$$r^*\OM^{\infty}X  = B\times X_0 \iso \OM^{\infty}_B r^* X$$
for a $G$-spectrum $X$.  It therefore suffices to show that $r_!(E,p)_+$ is equivalent to $E_+$, where $r_!$ denotes the functor on derived categories. 
By \myref{Qad10}, $r_!$ preserves $q$-equivalences between well-sectioned ex-spaces and it follows that $r_!Q(E,p)_+\simeq r_!(E,p)_+\cong E_+$ where the first equivalence is induced by $qf$-cofibrant approximation of $(E,p)_+$.
\end{proof}

To be precise about diagonal maps on the parametrized level, 
we consider base change along $\DE\colon B\rtarr B\times B$. 
We have the obvious commutative diagram 
$$\xymatrix{
E \ar[d]_{p} \ar[r]^-{\DE} & E\times E \ar[d]^{p\times p}\\
B \ar[r]_-{\DE} & B\times B.}$$
We consider $E$ as a space over $B\times B$ via this composite. The
diagonal map of $E$ then specifies a natural map
$$\DE_!((E,p)_+) = (E,\DE\com p)_+ \rtarr (E\times E,p\times p)_+\iso (E,p)_+\barwedge(E,p)_+$$
of ex-spaces over $B\times B$. This is a comparison map between 
Quillen left adjoints and therefore descends to a natural map in 
$\Ho G\sK_{B\times B}$. Its adjoint is a natural map $(E,p)_+\rtarr (E,p)_+\sma_B (E,p)_+$ in $\Ho G\sK_B$. Apply the (derived) suspension functor $\SI^{\infty}_B$ to this map and note that the target is equivalent 
to $\SI^{\infty}_{B}(E,p)_+\sma_B \SI^{\infty}_{B}(E,p)_+$, by \myref{SISISI2}.
This gives the required natural diagonal map
$$\DE_{(E,p)_+}\colon \SI^{\infty}_{B}(E,p)_+
\rtarr \SI^{\infty}_{B}(E,p)_+\sma_B \SI^{\infty}_{B}(E,p)_+$$
in $\Ho G\sS_B$.

\begin{defn}[The transfer map]\index{transfer map!for fibrations}\mylabel{fibtransfer}
Let $p\colon E\rtarr B$ be a Hurewicz $G$-fibration over $B$ such that each
fiber $E_b$ is homotopy equivalent to a retract of a finite $G_b$-CW-complex. 
Then $\Sigma_B^\infty (E,p)_+$ is a dualizable $G$-spectrum over $B$ by \myref{spaceFDT} and we obtain the transfer map 
$$\ta_{(E,p)_{+}}\colon S_{B}\rtarr \Sigma_B^\infty (E,p)_+$$
in $\Ho G\sS_B$.  Define the \emph{transfer map} of $E$ to be the map
$$\ta_E = r_!\ta_{(E,p)_{+}} 
\colon \Sigma^{\infty} B_+ \iso r_! S_B
\rtarr r_!\Sigma_B^\infty (E,p)_+\iso \Sigma^{\infty} E_+$$
in $\Ho G\sS$.
\end{defn} 

With this definition, all of the standard properties of transfer maps
are direct consequences of the general categorical results
\myref{traceprop} and \myref{addprop} and the properties of $r_!$.
\section{The bundle construction on parametrized spectra}\label{sec:bdlconstr}

The construction of the transfer in the previous section works ``globally'', starting on the parametrized spectrum level.  We now give a fiberwise construction of ``stable bundles'' that leads to an alternative fiberwise perspective.  However, it is natural to work in greater generality than is needed for the construction of the transfer. The extra generality will be needed in the proof of the Wirthm\"uller isomorphism in \S\ref{sec:fibwirth} and will surely find other applications. The relevant bundle theoretic background was recalled in \S3.2.

Let $\PI$ be a normal subgroup of a compact Lie group $\GA$ such that 
$\GA/\PI = G$ and let $q\colon \GA\rtarr G$ be the quotient homomorphism.  
Let $p\colon E\rtarr B$ be a $(\PI;\GA)$-bundle with fiber a $\GA$-space $F$
and with associated principal $(\PI;\GA)$-bundle $\pi\colon P\rtarr B$.  
Then $P$ is a $\PI$-free $\GA$-space, $\pi$ is the quotient map to 
the orbit $G$-space $B=P/\PI$, and $p$ is the associated $G$-bundle 
$E \iso \T{F}\rtarr B$.  
To simplify the homotopical analysis, we assume for the rest of this section that $F$ and $P$ are $\GA$-CW complexes such that $P$ is $\PI$-free.  We let $E ={\T{F}}$ and $B = \T{*}$. Note that $B$ is a $G$-CW complex.  We are thinking of the cases when $F$ is a point or when $F$ is a smooth $\GA$-manifold. On the ex-space level, application of $P\times_{\PI}(-)$ to retracts gives the functor
\[\U{_F} = \T(-)\colon \GA\sK_F \rtarr G\sK_E.\]\index{bundle construction}
Thus, for an ex-$\GA$-space $K$ over $F$, the ex-$G$-space $P\times_{\PI} K$ over $P\times_{\PI} F$ has section and projection induced by the section and projection of $K$. Observe that if $F$ is a smooth manifold and $S^{\ta}$ is the sphere bundle obtained by fiberwise one-point compactification of the tangent bundle of $F$, then $\T{S^{\ta}}$ is the $G$-bundle of spherical tangents along the fiber associated to $p$. 

We can extend the functor $\U{_F}$ from ex-spaces to ex-spectra.  Change of universe must enter since $\GA$-spectra are indexed on representations of $\GA$ and $G$-spectra are indexed on representations of $G$.  We view representations of $G$ as $\PI$-trivial representations of $\GA$. This gives $i\colon q^*\sV_G\rtarr \sV_{\GA}$.  Implicitly applying the functor $i^*$ to $\GA$-spectra, {\em we agree to index both $G$-spectra and $\GA$-spectra on $\sV_G$ for the rest of the section}.  We are interested in $\GA$-spectra indexed on a complete universe, and we shall return to this point in the next section.  Since $\PI$ acts trivially on our representations $V$, we have
$$ \U{_F{K}}\sma_E S^V \iso \U{_F{(K\sma_F S^V)}}.$$
Therefore, for a $\GA$-spectrum $X$ over $F$, the ex-$G$-spaces $\U{_F{X(V)}}$ over $E$ inherit structure maps from $X$, so that $\U{_FX}$ is a $G$-spectrum over $E$. We have the same definition on the prespectrum level. These functors $\U{_F}$ are exceptionally well-behaved, as the following results show.

\begin{prop}
The functor $\U{_F}\colon \GA\sS^{\text{$\Pi$-triv}}_F \rtarr G\sS_E$ is both a left and a right Quillen adjoint with respect to the level and stable model structures.  Moreover, the functor $\U{_F}\colon \GA\sP^{\text{$\Pi$-triv}}_F \rtarr G\sP_E$ takes excellent $\GA$-prespectra over $F$ to excellent $G$-prespectra over $E = \T{F}$.
\end{prop}

\begin{proof}
Let ${\pi}\colon P\times F\rtarr F$ be the projection. Clearly $\U{_F}$ is the composite of ${\pi}^*\colon \GA\sS_F \rtarr \GA\sS_{P\times F}$ and $(-)/\PI\colon \GA\sS_{P\times F}\rtarr G\sS_E$. By Propositions \ref{Qad1}, \ref{Qad2}, \ref{Qad1too}, and \ref{Qad2too}, ${\pi}^*$ is both a left and a right Quillen adjoint, provided we use appropriate generating sets in our
definitions of the model structures. By \myref{factor}, the functor $(-)/\PI$ is a Quillen left adjoint. By \myref{ouch}, it coincides with the right adjoint 
$(-)^{\PI}\com p_*$, where $p$ here is the quotient map $P\times F\rtarr P\times_{\PI}F = E$.  Using \myref{rho}, we see that $p\colon E\rtarr B$ is a  $G$-bundle with CW fibers. Therefore $p_*$ is a Quillen right adjoint by Propositions \ref{Qad2} and \ref{Qad2too}, and $(-)^{\PI}$ is a Quillen right adjoint by \myref{factor}. The last statement is easily checked from \myref{excel} and \myref{fpOM}.
\end{proof}

We need an observation about the behavior of $\U{_F}$ on fibers.

\begin{lem}\mylabel{tEtB} 
Fix $b\in B$. Let $\io\colon G_b\rtarr G$ and $\rh_b\colon G_b\rtarr \GA$ be the inclusion and the homomorphism of \myref{rho}. Let $b\colon \{b\}\rtarr B$ and  $i_b\colon E_b\rtarr E$ denote the evident inclusions of $G_b$ spaces. The following diagrams commute, and these commutation relations descend to homotopy categories.
\[\xymatrix{
\Gamma\sS_*^{\text{$\Pi$-triv}}\ar[d]_{\U{_{*}}}\ar[r]^-{\rho_b^*} & G_b\sS_{b}\\
G\sS_{B}\ar[r]_{\io^*} & G_b\sS_B \ar[u]_{b^*}}
\quad \text{and} \quad
\xymatrix{
\Gamma\sS^{\text{$\Pi$-triv}}_F\ar[d]_{\U{_F}}\ar[r]^-{\rho_b^*} & G_b\sS_{E_b} \\
G\sS_E\ar[r]_-{\io^*} & G_b\sS_E \ar[u]_{i_b^*}}\]
\end{lem}
\begin{proof} On the level of ex-spaces, this is immediate by inspection.
The diagrams extend levelwise to parametrized spectra, and passage to
homotopy categories is clear from the previous result.
\end{proof}

Writing $\SI^{\infty}$ for suspension spectra functors indexed on complete universes, we have that $i^*\SI^{\infty}_{\GA,F}$, where 
$i\colon q^*\sV_G\subset \sV_\Gamma$, is the suspension $\GA$-spectrum functor indexed on the $\PI$-trivial $\GA$-universe $q^*\sV_G$.

\begin{thm}\mylabel{PSmash}
There is a natural isomorphism of functors
$$\U{_F{i^*\SI^{\infty}_{\GA,F}}}\iso\SI^{\infty}_{G,E}\U{{_F}}\colon 
\GA \sK_F\rtarr G\sS_E,$$
and this isomorphism descends to homotopy categories.  The functor
$$\U{_F}\colon \Ho \GA\sS^{\text{$\Pi$-triv}}_F \rtarr \Ho G\sS_E$$
is closed symmetric monoidal.
\end{thm}

\begin{proof}
Let $K$ be an ex-$\GA$-space over $F$. Since we are indexing on representations $V$ of $G$, we have isomorphisms
\[(\U{_F{i^*\SI^{\infty}_{\GA,F}K}})(V)  =  \T{(K\sma_F S^V_F)}
\iso  (\T{K})\sma_E S^V = (\SI^{\infty}_{G,E}\U{_{F}K})(V).\]
This gives a natural isomorphism of $G$-spectra over $E$, and it descends to homotopy categories since it is a comparison of composites of Quillen left adjoints. Note in particular that $\U{_F}i^*S_{\GA,F}$ is isomorphic to $S_{G,E}$. We must show that the functor $\U{_F}$ commutes up to coherent natural isomorphism with smash products and function objects.  For ex-$\GA$-spaces $K$ and $L$ over $F$, it is easy to check that there is a natural isomorphism
$$ \U{{_F}{(K\sma_F L)}} \rtarr \U{{_F}{K}}\sma_E \U{{_F}{L}}$$
of ex-$G$-spaces over $E$.  This isomorphism extends levelwise to external smash products (external in the sense of pairs of representations). However, since external pairings (in the sense of pairs of base spaces) do not naturally come into play here, to retain homotopical control it seems simplest to just extend levelwise to handicrafted smash products of $\GA$-prespectra; compare \myref{subtlety}. Using excellent prespectra to pass to homotopy categories of prespectra and then using the equivalence $(\bP,\bU)$ to pass to homotopy categories of spectra, we obtain the required natural equivalence
$$ \U{{_F}{(X\sma_F Y)}} \rtarr \U{{_F}{X}}\sma_E \U{{_F}{Y}}$$
in $\Ho G\sS_E$ for $\GA$-spectra $X$ and $Y$ over $F$. The
adjoint of the composite
$$\xymatrix{
\U{{_F}{F_F(X,Y)}}\sma_E \U{{_F}{X}}\simeq
\U{{_F}{(F_F(X,Y)\sma_F X)}} \ar[rr]^-{\U{{_F}{(\text{ev})}}}&&
\U{{_F}{Y}}}$$
is a natural map 
$$ \U{{_F}{F_F(X,Y)}} \rtarr F_E(\U{{_F}{X}},\U{{_F}{Y}})$$
in $\Ho G\sS_E$, and we must show that it is an equivalence. This will hold if it holds when restricted to fibers over points of $E$. Since each point is in some $E_b$, it suffices to show that the restriction to each $E_b$ is an equivalence.  However, using \myref{tEtB}, we see that the restriction to $E_b$ is the adjoint to the $G_b$-map $\text{ev}\colon F_{E_b}(\rh_b^*X,\rh_b^*Y)\sma_{E_b} \rh_b^*X \rtarr \rh_b^*Y$, and is thus the identity map. 
\end{proof}

We have the following relations between $\U{_F}$ and base change functors.

\begin{prop}\mylabel{Pfibers} 
Consider $r\colon F\rtarr *$ and $p = \T{r}\colon E\rtarr B$. 
For $Y \in \GA\sS_*^{\text{$\Pi$-triv}}$ and $X\in \GA\sS^{\text{$\Pi$-triv}}_F$, 
there are natural isomorphisms
\[p_!\U{_F{X}} \rtarr \U{_*{r_!X}}, \quad
\U{_F{r^*Y}}\rtarr p^*\U{_*{Y}},  \quad \text{and} \quad
\U{_*{r_*X}}\rtarr p_*\U{_F{X}},\]
and these isomorphisms induce natural equivalences on homotopy categories.
\end{prop}

\begin{proof}
We first work on the ex-space level. Let $T$ be a based $G$-space and $K$ be an ex-$G$-space over $F$. Applying the functor $\T{(-)}$ to the maps of retracts that define $r_!K$ and $r^*T$ (see \myref{retract1}), we immediately obtain the first two maps.  The first is the natural isomorphism
$$(\T{K})\cup_E B \iso \T{(K/F)}$$
in which the section $F$ is collapsed to a point in $K$ on both sides. The second is the evident natural isomorphism
$$P\times_{\Pi} (F\times T)\iso (P\times_{\Pi} F)\times_B (P\times_{\Pi} T).$$
For the third map, recall that $r_*K =\text{Sec}(F,K)$.  The adjoint of 
$$\xymatrix{
(\T{\text{Sec}(F,K)})\times_B E 
\iso \T{(\text{Sec}(F,K)\times F)} \ar[rr]^-{\T{\text{ev}}} && \T{K}}$$
gives a map
$$ \U{_*{r_*K}} = \T{\text{Sec}(F,K)}\rtarr \text{Map}_B(E,\T{K}).$$ 
Together with the projection of the source to $B$, it induces an isomorphism to $p_*\U{_F{K}}$, which is the pullback along $B\rtarr \text{Map}_B(E,E)$ of the projection of the target induced by the projection $P\times_{\PI}K\rtarr P\times_{\PI}F = E$. Applied levelwise, these point-set level isomorphisms carry over directly to parametrized prespectra and spectra. We must show that they descend to equivalences in homotopy categories. Since \myref{Qad2too} applies to show that both $p_*$ and $r_*$ are Quillen right adjoints (and we have no need to use Brown representability here), the first commutation relation is between composites of left Quillen adjoints, the second is between functors that are both left and right Quillen adjoints, and the third is between Quillen right adjoints, so descent to homotopy categories is immediate.
\end{proof}

\section{$\PI$-free parametrized $\GA$-spectra}\label{sec:PIFree}
We retain the notations of the previous section in this section and the next.
In the next section, we show that the bundle construction on parametrized 
spectra leads to a fiberwise generalization of the restriction to bundles 
of the trace and transfer maps for fibrations that we described in \S\ref{sec:trfr}. The definition depends on a result that is proven by use 
of the theory of $\PI$-free $\GA$-spectra that we present here.  

We first recall what it means to say that a $\GA$-spectrum $X$ (indexed on any 
universe) is $\PI$-free. Let $\sF(\PI;\GA)$ be the family of subgroups $\LA$ of 
$\GA$ such that $\LA\cap \PI = e$. A $\GA$-CW 
complex $T$ is $\PI$-free if and only if the only orbit types $\GA/\LA$ that appear 
in its construction have $\LA\in \sF(\PI;\GA)$. We then say that $T$ is an 
$\sF(\PI;\GA)$-CW complex.  We can make the same definitions for $\GA$-CW spectra, 
and in general we say that a $\GA$-spectrum is $\PI$-free if it is isomorphic in 
$\Ho \GA\sS$ to an $\sF(\PI;\GA)$-CW spectrum. There is a more conceptual 
homotopical reformulation that is the one relevant to the parametrized point of 
view and that does not depend on the theory of $\GA$-CW spectra.

Let $E(\PI;\GA)$ be the universal $\PI$-free $\GA$-space, so that $E(\PI;\GA)^{\LA}$ 
is contractible if $\LA\cap \PI = e$ and is empty otherwise. We may take $E(\PI;\GA)$ 
to be an $\sF(\PI;\GA)$-CW complex. Let $B(\PI;\GA) = E(\PI;\GA)/\PI$ and observe that 
$B$ is a $G$-CW complex and therefore also a $\GA$-CW complex. We note parenthetically
that the quotient map $p: E(\PI;\GA) \rtarr B(\PI;\GA)$ is the universal principal 
$(\PI;\GA)$-bundle. That is, pullback along $p$ gives a bijection
\[[X,B(\PI;\GA)]_G\rtarr \sB{(\PI;\GA)}(X),\]
where $\sB{(\PI;\GA)}(X)$ denotes the set of equivalence classes of 
principal $(\PI;\GA)$-bundles over the $G$-space $X$; see \cite{LM} or
\cite [VII\S2]{EHCT}. 

\begin{defn}\mylabel{Nfree} Let $r\colon E(\PI;\GA)\rtarr *$ be the projection
and let $\si$ be the counit of the (derived) adjunction $(r_!,r^*)$.
A $\GA$-spectrum $X$ is said to be \emph{$\PI$-free} if $\si\colon r_!r^* X\rtarr X$ is an equivalence.
\end{defn}

The definition should seem reasonable since $r_!r^* T\iso E(\PI;\GA)_+\sma T$
for a $\GA$-space $T$. It is equivalent to the original definition in 
terms of an equivalence in $\Ho G\sS$ to an $\sF(\PI;\GA)$-CW spectrum; 
see \cite[II.2.12]{LMS} or \cite[VI\S4]{MM}. The definition generalizes
readily to the parametrized context.
\begin{defn} Let $\pi\colon  E(\Pi;\Gamma)\times F\rtarr F$ be the 
projection and let $\sigma$ be the counit of the (derived) adjunction $(\pi_!,\pi^*)$. An ex-$\GA$-space or $\Gamma$-spectrum $X$ over a 
$\GA$-space $F$ is said to be \emph{$\Pi$-free} if 
$\sigma\colon  \pi_!\pi^* X \rtarr X$ is an equivalence.
\end{defn}

Since the fiber $(\pi_!\pi^* X)_f$ is $E(\PI;\GA)_+\sma X_f$, the definition should seem reasonable. Since equivalences are detected fiberwise, we have
the following results.
\begin{lem}\mylabel{FibPI} 
A $\GA$-spectrum $X$ over $F$ is $\PI$-free if and only if each of its fibers $X_f$ is a $(\Pi\cap\GA_f)$-free $\GA_f$-spectrum.
\end{lem}

\begin{proof}
The fiber of $E(\Pi;\Gamma)\times F\rtarr F$ over $f\in F$ is the $\Gamma$-space $E(\Pi;\Gamma)$ with the action restricted along $\iota\colon \Gamma_f\rtarr \Gamma$. It is a model of the universal $(\Pi\cap\Gamma_f)$-free $\Gamma_f$-space $E(\Pi\cap\Gamma_f,\Gamma_f)$. Applying $(-)_f$ to the counit $\pi_!\pi^* X \rtarr X$ and using \myref{Mackeymore} we obtain the counit $r_!r^* X_f\rtarr X_f$ where $r\colon \iota^* E(\Pi;\Gamma)\rtarr *$.
\end{proof}

\begin{lem}\mylabel{IsPI}  If $P$ is a $\PI$-free $\GA$-space and $X$ 
is any ex-$\GA$-space or $\GA$-spectrum over $F$, then $P\times X$
is a $\PI$-free ex-$\GA$-space or $\GA$-spectrum over $P\times F$.
\end{lem}

A useful slogan asserts that ``$\PI$-free $\GA$-spectra live in the 
$\PI$-trivial universe''. To explain it, consider the inclusion $i\colon q^*\sV_G\rtarr \sV_{\GA}$ of the complete $G$-universe $\sV_G$ as the 
universe of 
$\PI$-trivial representations in the complete $\GA$-universe $\sV_{\GA}$. 
Then the slogan is given meaning by the following result.  In the nonparametrized case $F=*$, it is proven in \cite[II\S2]{LMS} and is
discussed further in \cite[VI\S4]{MM}. Since the parametrized case 
presents no complications and the proof is quite easy, we only give a 
sketch.

\begin{prop}\mylabel{Nfreeii}
The change of universe adjunction $(i_*,i^*)$ descends to a symmetric monoidal equivalence between the homotopy categories of $\PI$-free $\GA$-spectra over $F$ indexed on $\PI$-trivial representations of $\GA$ on the one hand and indexed on all representations of $\GA$ on the other. For $\PI$-free $\GA$-spectra $X$ over $F$ indexed on $\sV_{\GA}$, there is a natural equivalence 
$i_*(E(\PI;\GA)_+\sma i^*X)\simeq X$.
\end{prop}

\begin{proof}[Sketch Proof] If $\LA\cap \PI = e$, then the quotient map
$q\colon \GA\rtarr G$ maps $\LA$ isomorphically onto a subgroup of $G$. Any
representation $V$ of $\LA$ is therefore of the form $q^*W$ for a representation
$W$ of $q(\LA)$. It follows that the restrictions to $\LA$ of the universes 
$\sV_{\GA}$ and $q^*\sV_G$ have the same representations.  This makes clear that,
on $\PI$-free $\GA$-spectra over $F$, the unit and counit of the adjunction 
$(i_*,i^*)$ must be $\sF(\PI;\GA)$-equivalences, in the sense that they are 
$\LA$-equivalences for any $\LA$ in $\sF(\PI;\GA)$. Smashing the unit and counit 
with $E(\PI;\GA)_+$, which has trivial fixed point sets for subgroups not in 
$\sF(\PI;\GA)$, we obtain natural equivalences, and it follows from 
\myref{Nfree} that the unit and counit are themselves equivalences 
when applied to $\PI$-free $\GA$-spectra. Alternatively, restricting to 
$s$-fibrant $\GA$-spectra over $F$, the conclusion follows fiberwise from its
nonparametrized precursor. Since $i_*$ is symmetric monoidal, by \myref{change1}, 
so is the equivalence. The last statement holds since
\[i_*(E(\PI;\GA)_+\sma i^*X)\htp E(\PI;\GA)_+\sma i_*i^*X\htp X. \qedhere\]
\end{proof}

\section{The fiberwise transfer for $(\PI;\GA)$-bundles}\label{sec:fibtrfr}

We consider a fixed given principal $(\PI;\GA)$-bundle $P$, where 
$\PI$ is a normal subgroup of $\GA$ with quotient group $G$ and 
quotient map $q\colon \GA\rtarr G$. We also consider a $\GA$-space 
$F$ and the associated $(\PI;\GA)$-bundle 
$$p\colon E = P\times_{\PI}F\rtarr P\times_{\PI}* = B.$$
We have the inclusion $i\colon q^*\sV_G\rtarr \sV_{\GA}$ of the complete 
$G$-universe $\sV_G$ as the universe of $\PI$-trivial representations 
in the complete $\GA$-universe $\sV_{\GA}$.

The change of universe functor $i^*\colon \GA\sS_{F}\rtarr \GA\sS_{F}^{{\PI}-{\text{triv}}}$ is {\em not}\, symmetric mon\-oid\-al, and it does not preserve dualizable objects. For example, with $F=*$ and $\PI =e$, the orbit spectrum $i^*\SI^{\infty}\GA/\LA$ is not dualizable if $\LA$ is a non-trivial subgroup of $\GA$.  The bundle theoretic study of transfer maps is based on the following
result, whose proof is based on the theory of $\PI$-free $\GA$-spectra
given in the previous section.

\begin{thm}\mylabel{bomb} The composite functor 
$\U{{_F}{i^*}}\colon \Ho \GA\sS_F\rtarr \Ho G\sS_E$
is symmetric mon\-oid\-al.
\end{thm}

\begin{proof}
Let $\pi\colon P\times F\rtarr F$ be the projection and note that
${\pi}^*X = P\times X$. The functor $\U{_F}$ is the composite of the 
symmetric monoidal Quillen left adjoint ${\pi}^*$ and the Quillen left 
adjoint $(-)/\PI$.  By \myref{PSmash}, the functor $\U{_F}$ on homotopy
categories is also symmetric monoidal since the $\GA$-space $P$ is $\PI$-free. We observe first that the composite ${\pi}^*i^*$ is symmetric monoidal.
Indeed, for $\GA$-spectra $X$ and $Y$ over $F$, we have
\begin{align*}
{\pi}^*i^*(X\sma_FY) & \htp i^*{\pi}^*(X\sma_FY) 
&& \text{by \myref{chvschuni}}\\
& \htp  i^*({\pi}^*X\sma_{P\times F}{\pi}^*Y) 
&& \text{by \myref{Wirthmore}}\\
& \htp  i^*{\pi}^*X\sma_{P\times F}i^*{\pi}^*Y 
&& \text{by \myref{IsPI} and \myref{Nfreeii}}\\
& \htp  {\pi}^*i^*X\sma_{P\times F} {\pi}^*i^*Y 
&& \text{by \myref{chvschuni}.}
\end{align*}
It follows directly that $\U{{_F}i^*}$ is symmetric monoidal:
\begin{align*}
\U{{_F}i^*(X\sma_FY)} & =  ({\pi}^*i^*(X\sma_FY))/\PI
&& \text{by definition}\\
& \htp  ({\pi}^*i^*X\sma_{P\times F} {\pi}^*i^*Y)/\PI
&& \text{by the previous display}\\
& \htp  ({\pi}^*(i^*X\sma_{F} i^*Y))/\PI
&& \text{by \myref{Wirthmore}}\\
& =  \U{{_F}(i^*X\sma_{F} i^*Y})
&& \text{by definition}\\
& \htp  \U{{_F}i^*X}\sma_E \U{{_F}i^*Y}
&& \text{by \myref{PSmash}}. \qedhere
\end{align*}
\end{proof}

Now \myref{traceprop}(i) shows that $\U{{_F}{i^*}}$ commutes with trace maps. 

\begin{thm}\mylabel{PtrP}
Let $X\in \Ho\Gamma \sS_F$ be dualizable. Then $\bar{P}_F i^*X \in \Ho G\sS_{E}$ is dualizable. Suppose given a coaction map $\Delta_X:X\to X\wedge_F C_X$ and a self map $\ph\colon X\rtarr X$. Then 

\[\tau(\bar{P}_F i^*\ph)\simeq \bar{P}_F i^* \tau(\ph)\colon S_{E}\rtarr \U{{_F}{i^*C_X}},\]
where $\bar{P}_F i^*X$ is given the coaction map 
$$\bar{P}_Fi^*(\Delta_X):\bar{P}_Fi^* X \rtarr 
\bar{P}_Fi^* (X \wedge_{F} C_X) \htp \bar{P}_Fi^* X \wedge_{E} \bar{P}_Fi^*C_X.$$
\end{thm}

These trace maps are maps of $G$-spectra over $E$, rather than over $B$.
We can apply $r_!$, $r\colon E\rtarr *$, to obtain trace maps of 
nonparametrized spectra.  This kind of trace map can be viewed as a
fiberwise generalization of the kind of nonparametrized trace map 
that is defined bundle theoretically in the literature.  To connect 
up with the latter, we specialize and change our point of view so as 
to arrive at bundle theoretic trace maps over $B$. Specializing further 
to transfer maps, we obtain the promised comparison with the transfer 
maps of \myref{fibtransfer}.

With these goals in mind, we now focus on the case $F=*$, so that $E$ above becomes $B$, with $p$ the identity map, and our trace maps are parametrized 
over $B$.  We study our original fixed given $(\PI;\GA)$-bundle $p\colon E\rtarr B$ in a different fashion. We rename its fiber $M$ to avoid confusion with respect to the role that space is playing. In the theory above, $F$ was 
a base space for paramentrized spectra and there was no need for $F$ to be dualizable.  We now consider the case when $M$ is stably dualizable, so that $\SI^{\infty}M_+$ is dualizable, and we write $\ta_M$ for the transfer map $S\rtarr \SI^{\infty}M_+$ in $\GA\sS$, as defined in and after \myref{tracemap}. We apply \myref{PtrP} with $F=*$ and $X = \SI^{\infty}M_+$ to obtain the following special case. Here we use the diagonal map induced by the diagonal 
map of $M$. Observe that, by \myref{PSmash},
$$\U{{_*}{i^*\SI^{\infty}M_+}} \simeq \SI^{\infty}\U{{_*}{M_+}} 
= \Sigma^\infty_B (E,p)_+.$$

\begin{thm}\mylabel{trantran}
Let $M$ be a compact $\Gamma$-ENR and let $p\colon E\rtarr B$
be a $(\PI;\GA)$-bundle with fiber $M$ and associated principal 
$(\PI;\GA)$-bundle $P$.  Let $\ph$ be a self-map of $\SI^{\infty}M_+$.  
Then
\[\tau(\bar{P}_*i^*\ph)\simeq \bar{P}_*i^*(\tau(\ph)): 
S_B\to \Sigma^\infty_B (E,p)_+.\]
Therefore, taking $\ph = \text{id}$ and applying $r_!$, $r\colon B\rtarr *$,
\[\ta_E \simeq r_!\bar{P}_*i^*\ta_M\colon \SI^{\infty} B_+\rtarr \SI^{\infty} E_+.\]
\end{thm}

This result gives a clear and precise comparison between 
the specialization to bundles of the globally defined transfer map for 
Hurewicz fibrations and the fiberwise transfer map for bundles. Effectively, 
we have inserted the transfer map for $M_+$ fiberwise into $P\times_{\PI}(-)$ 
to obtain an alternative description of the transfer map for the dualizable 
$G$-spectrum $\SI^{\infty}(E,p)_+$ over $B$.

There is a useful reinterpretation of the 
description of transfer maps given by \myref{trantran}. Consider 
${\pi}\colon P\rtarr *$. Observe that, by \myref{chvsfixorbit}, 
instead of applying $r_!$,
$r\colon B\rtarr *$, to orbit spectra under the action of $\PI$, 
we could first apply ${\pi}_!$ and then pass to orbits. For a 
$\GA$-spectrum $X$, we have a natural isomorphism 
$${\pi}_!{\pi}^*i^*X\iso P_+\sma i^*X$$
and a natural equivalence
$$i_*(P_+\sma i^*X) \htp P_+\sma X.$$

\begin{cor}
let $M$ be a compact $\Gamma$-ENR and let $p\colon E\rtarr B$ be a 
$(\PI;\GA)$-bundle with fiber $M$ and associated 
principal $(\PI;\GA)$-bundle $P$.  Then the transfer $\ta_E\colon \SI^{\infty} B_+\rtarr \SI^{\infty} E_+$
is obtained by passage to orbits over $\PI$ from the map
$$\tilde{\ta} = \text{id}\sma i^*\ta_M 
\colon P_+\sma i^*S\rtarr P_+\sma i^*\SI^{\infty} M_+,$$
and $i_*\tilde{\ta}$ can be identified with 
$$\text{id}\sma \ta_M\colon P_+\sma S\rtarr P_+\sma \SI^{\infty}M_+.$$
\end{cor}

\begin{rem} The corollary gives exactly the transfer map as defined
by Lewis and May \cite[IV.3.1]{LMS}.  Working in the nonparametrized
context, they tried in vain to obtain a spectrum level transfer map
for Hurewicz fibrations over general base spaces. The comparison here
also sheds light on the relationship between the two constructions of 
Becker and Gottlieb \cite{BG1, BG2}, both of which require finite
dimensional base spaces. The first is bundle theoretic and is easily
seen to be equivalent to the construction in this section by using
Atiyah duality to interpret $\ta_M$ for a $\GA$-manifold $M$. 
Precisely, by \cite[IV.2.3]{LMS}, if $M$ is embedded in $V$ with normal
bundle $\nu$ and $\ta$ is the tangent bundle of $M$, then the transfer
map $\ta_M$ is homotopic to the map obtained by applying the functor 
$\SI^{-V}\SI^{\infty}$ to the composite of the Pontryagin-Thom map 
$S^V\rtarr T\nu$ and the map $T\nu\rtarr T(\nu\oplus\ta)\iso M_+\sma S^V$ 
induced by the inclusion $\nu\rtarr \nu\oplus \ta$. The second,
which is generalized to the equivariant setting by Waner \cite{Waner}, 
is fibration theoretic and is easily seen to be equivalent to the 
construction of \S\ref{sec:trfr}.  Another approach to the comparison is to show 
that suitable Hurewicz fibrations are equivalent to bundles,
as is done by Casson and Gottlieb in \cite{CG}.
\end{rem}

\begin{rem} Since our definition coincides with that of \cite[IV.3.1]{LMS},
the properties of the transfer catalogued in \cite[IV\S\S3--7]{LMS} 
apply verbatim. Many of these properties generalize directly to the 
parametrized trace and transfer maps of \myref{PtrP}. Actually, the 
definition of 
\cite[IV.3.1]{LMS} works more generally with $P$, or rather $i^*\SI^{\infty}P_+$, 
replaced by a general $\PI$-free $\GA$-spectrum $\bP$ indexed on $\sV_G$. The 
constructions here admit similar generalizations.  One way to achieve this with 
minimal work is to use the case $P=E(\PI;\GA)$ of the construction 
already on hand. Thus, for a $\PI$-free $\GA$-spectrum $\bP$ over $F$
indexed on $\sV_G$, we can define
$$ \bar{\bP}_{F}i^*X = \overline{E(\PI;\GA)}_{F}(\bP\sma_F i^*X)$$
and develop parametrized trace and transfer maps from there. 
We leave the further development of the theory to the interested reader.
\end{rem}

\section{Sketch proofs of the compatible triangulation axioms}\label{sec:comptriang}

We must explain why $\Ho G\sS_B$ is a closed symmetric monoidal 
category with a compatible triangulation, in the sense specified in
\cite{Tri}. We have the closed symmetric monoidal structure and the 
triangulation, the latter by \myref{yestrian}. We must prove the 
compatibility axioms (TC1)--(TC5) of 
\cite[\S4]{Tri}. The essential idea is to verify the axioms using external 
smash products and function 
objects and then pull back along diagonal maps to obtain the conclusions. The axiom 
(TC1) only involves suspension, in our case $\SI_B$, and is thus easily checked 
using \myref{spacesmashpair}. For (TC2), we must show that the functors
$X\sma_B(-)$, $F_B(X,-)$, and $F_B(-,Y)$ preserve distinquished triangles, where
$X$ and $Y$ are $G$-spectra over $B$.
Either model theoretically or by standard topological arguments with cofiber
sequences and fiber sequences, it is easy to see that these conclusions hold with 
$\sma_B$ and $F_B$ replaced by the external functors $\barwedge$ and $\bar{F}$.
Since $\DE^*$ and therefore its right adjoint $\DE_*$ are exact, the conclusion 
internalizes directly. Similarly, the braid axiom (TC3) and additivity axiom (TC4) 
hold for $\barwedge$ by the arguments explained in \cite[\S6]{Tri}, and they pull 
back along $\DE^*$ to give these axioms internally in $\Ho G\sS_B$. 
 
The braid duality axiom (TC5) is more subtle because it involves simultaneous use 
of $\sma_B$ and $F_B$.  Externally, we can work over $B\times B$, using $\barwedge$.
Inspecting the argument in \cite[\S7]{Tri}, we see that the only internal homs used 
in the verification of the braid duality axiom are duals of the form 
$F(-,T)$ for a suitable approximation $T$ of the unit object. In our context, it
turns out that we need to use two analogues of this functor, one to mimic the 
proof of (TC5a) given in \cite[pp\,62-64]{Tri} and another to mimic the proof of
(TC5b) given in \cite[pp\,65-67]{Tri}.  For the first, let $T\in G\sS_{B\times B}$ 
be a fibrant model of the derived $\DE_*S_B$, so that $\bar{F}(X,T)$ is a model for 
$DX = F_B(X,S_B)$ in $\Ho G\sS_B$.  With this replacement for $F(X,T)$,
the cited proof of (TC5a) goes through, first working externally and then 
internalizing along $\DE^*$.  The cited proof of (TC5b) relies on a natural 
point-set level map
\begin{equation}\label{oldpair} 
F(X,T)\sma F(Y,T)\rtarr F(X\sma Y,T),
\end{equation}
and this makes no sense in our external context.  Working internally, in 
$\Ho G\sS_B$, we have such a map 
\begin{equation}\label{intequ}
F_B(X,S_B)\sma_B F(Y,S_B)\rtarr F_B(X\sma_B Y,S_B),
\end{equation}
but we need a point-set level external model for it to carry out the cited argument.
Let $U$ be a fibrant model for $\DE_*S_{B\times B}$ in 
$G\sS_{B\times B\times B\times B}$. Replacing
the functor $D'(-) = F(X,T)$ used in \cite[pp.\, 66-67]{Tri} with the functor 
$$D'(-) =\bar{F}(-,U)\colon G\sS_{B\times B} \rtarr G\sS_{B\times B},$$
we find that the cited argument goes through verbatim on the external level,
working in the category $G\sS_{B\times B}$, once we construct a natural map
\begin{equation}\label{goodpair}
\bar{F}(X,T)\barwedge \bar{F}(Y,T) \rtarr \bar{F}(X\barwedge Y,U)
\end{equation}
in $G\sS_{B\times B}$ to substitute for the pairing (\ref{oldpair}). Starting 
from the $(\barwedge,\bar{F})$ adjunction, we obtain an external pairing
\begin{equation}\label{extpair}
\bar{F}(X,T)\barwedge \bar{F}(Y,T)\rtarr \bar{F}(X\barwedge Y,T\barwedge T).
\end{equation}
We also have the natural map
\[\xymatrix{
\DE_*X\barwedge \DE_* Y \ar[d]^-{\et} \\
\DE_*\DE^*(\DE_*X\barwedge \DE_* Y)\ar[d]^{\iso}\\
\DE_*(\DE^*\DE_*X\barwedge \DE^*\DE_* Y)\ar[d]^-{\DE_*(\epz\barwedge\epz)}\\
\DE_*(X\barwedge Y).}\]
Applying this with $X=Y=S_B$ and using that $S_B\barwedge S_B$ is isomorphic
to $S_{B\times B}$, we obtain a lift $\xi$ in the diagram
$$\xymatrix{
\DE_* S_B\barwedge \DE_*S_B \ar[r] \ar[d]
& \DE_*(S_B\barwedge S_B)\iso \DE_*S_{B\times B}  \ar[r] & U \ar[d]\\
T\barwedge T \ar@{->}[urr]^{\xi} \ar[rr] & & {*}_{B\times B}}$$
Composing $\bar{F}(X\barwedge Y, \xi)$ with the pairing (\ref{extpair}),
we obtain the required pairing (\ref{goodpair}). 
Internalization along $\DE^*$ is then a not altogether trivial exercise which
shows that, on passage to homotopy categories, application of $\DE^*$ to the 
pairing (\ref{goodpair}) gives a model for the pairing (\ref{intequ}). 
The latter pairing can be viewed as a map 
$$ \DE^*(\bar{F}(X,T)\barwedge \bar{F}(Y,T))\rtarr 
\bar{F}(\DE^*(X\barwedge Y),T),$$
and the essential point of the exercise is to verify that 
$\DE^*\bar{F}(X\barwedge Y,U)$
is equivalent to $\bar{F}(\DE^*(X\barwedge Y),T)$. Using that
$\DE^*S_{B\times B} \iso S_{B}$ and looking at represented functors,
we see that a Yoneda lemma argument reduces the verification to the 
proof of a derived analogue of (\ref{four}) that is proven in the 
same way as \myref{Wirthmore}.  

\chapter{The Wirthm\"uller and Adams isomorphisms}

\section*{Introduction}
This chapter consists of variations on a theme. For a $G$-map 
$f\colon A\rtarr B$, the base change functor $f^*$ from $G$-spectra over 
$B$ to $G$-spectra over $A$ has a left adjoint $f_!$ and a right adjoint $f_*$.
We study comparisons between $f_!$ and $f_*$.  As preamble, we show in \S16.1 that there is always a natural map $\ph\colon f_!\rtarr f_*$ that relates the two adjunctions. It is an equivalence when $f$ is a homotopy equivalence, but not in general. This comparison is largely formal and applies to analogous sheaf theoretic contexts.

In the rest of the chapter, we use our foundations together with formal arguments developed in \cite{FHM} to obtain a simple proof of a general version of the Wirthm\"uller isomorphism and to reprove the Adams isomorphism 
as a special case. This material constitutes a considerably simplified 
version of work of Po Hu on the same topic \cite{Hu}. We consider $G$-bundles 
$p\colon E\rtarr B$, as in \S3.2 and \S15.3.  We assume that the fiber $M$ 
is a smooth closed $\GA$-manifold; manifolds with boundary work similarly. 
The generalized Wirthm\"uller isomorphism computes the relatively mysterious right adjoint $p_*$ of the functor $p^*$ as a suitable shift of the relatively familiar left adjoint $p_!$. 

We explain the result in the special case when $E\rtarr B$ is $M\rtarr *$ in
\S16.2, but we defer the proof to \S16.5.  We also show how to relate the Wirthm\"uller 
isomorphisms for $M$ and $N$ when $N$ is smoothly embedded in $M$.  
When $M = G/H$, 
the result specializes under the equivalence between the category of 
$G$-spectra over $G/H$ and the category of $H$-spectra to the Wirthm\"uller 
isomorphism in the form proven by Lewis and May \cite[II\S6]{LMS}.
As explained in \cite{MayW}, the categorical analysis in \cite{FHM} allows 
considerable simplification of that proof.  Our proof for general $M$ follows 
the same pattern, but it is quite different in detail since the special case 
$M=G/H$ has certain simplifying features. For example, when $G$ is finite, that case follows formally from Atiyah duality for $G/H$ and the trivial observation that 
$H/K_+$ is an $H$-retract of $G/K_+$ for $K\subset H\subset G$. 

In \S16.3, we show that the general case of $G$-bundles $p\colon E\rtarr B$ reduces 
fiberwise to the special case $M\rtarr *$.  The proof is an immediate application
of the construction $P\times_{\PI}(-)$ on parametrized $\GA$-spectra that was 
studied in \S15.3.  This allows 
a simple fiberwise construction of the $G$-spectrum over $E$ by which one must shift 
$p_!$ to obtain the desired isomorphism.  With this construction, it is immediate that 
the map of $G$-spectra over $B$ that we wish to prove to be an equivalence coincides 
on the fiber over $b$ with a map that we know to be an equivalence by the case 
$M\rtarr *$.  Since equivalences are detected fiberwise, that proves the result.

In turn, we prove in \S16.4 that the Adams isomorphism relating orbit spectra and fixed point spectra that was proven by Lewis and May in \cite[II\S8]{LMS} 
is a virtually immediate special case of our generalized Wirthm\"uller isomorphism.  These results complete the program originated in \cite{MM} of reproving conceptually all of the basic foundational results that were first proven in a less satisfactory ad hoc way in \cite{LMS}. The pioneering work of Po Hu \cite{Hu} paved the way but, in the absence of adequate foundations, the 
bundle construction of \S15.3, and the simplifying framework of \cite{FHM},
her arguments were very long and difficult. Our work recovers variant 
versions of all of her results. The basic idea that parametrized $G$-spectra should clarify and simplify the Wirthm\"uller and Adams isomorphisms is due 
to Gaunce Lewis \cite{LewisE}.

Again, we assume throughout that all given groups $G$ are compact Lie groups 
and all given base $G$-spaces are of the homotopy types of $G$-CW complexes.

\section{A natural comparison map $f_!\rtarr f_*$}

The Wirthm\"uller isomorphism that is the subject of the next few sections gives an equivalence between $f_*$ and a shift of $f_!$ for certain equivariant bundles $f$. In the course of our work on that, we came upon a curious natural comparison map $f_!\rtarr f_*$ for any map $f$ whatever.  We have no current applications for it, but since the relationships among base change functors are so central to the theory and its applications, we shall describe that map in this digressive section. It works just as well on the level of ex-spaces and indeed quite generally in other contexts where one has analogous base change adjunctions. 

\begin{thm}\mylabel{Please} 
Let $f\colon A\rtarr B$ be a $G$-map and let $X$ be a $G$-spectrum over $A$. Let $\epz\colon f^*f_*\rtarr\text{Id}$ and $\si\colon f_!f^*\rtarr \text{Id}$ denote the counits of the adjunctions $(f^*,f_*)$ and $(f_!,f^*)$ relating $\Ho G\sS_A$ and $\Ho G\sS_B$. There is a natural map $\ph\colon f_!X\rtarr f_*X$ in $\Ho G\sS_B$ such that the following diagram commutes:
\begin{equation}\label{newdia}
\xymatrix{
& f_!f^*f_*X \ar[dl]_{f_!\epz} \ar[dr]^{\si}  & \\
f_!X \ar[rr]_-{\ph} & & f_* X.\\}
\end{equation}
\end{thm}

\begin{proof}
Let $K = (K,p,s)$ be an ex-space over $A$. Then $f_!K = K\cup_A B$ and $f^*f_!K = f_!K\times_B A$. Here points $(f(a),a)$ in $B\times_B A$ are identified with points $(s(a),a)$ in $K\times_B A$, and we see that $f^*f_!K$ can be identified with the pullback $K\times_B A$. The projection to $K$ is then a map $\ps\colon f^*f_!K\rtarr K$ of ex-spaces over $A$. When $K = f^*L$ for an ex-space $L$ over $B$, $\ps = f^*\si\colon f^*f_!f^*L\rtarr f^*L$ since $f^*\si$ is also given by the projection $f^*L\times_B A\rtarr f^*L$. Passing to spectra over $A$ levelwise, we obtain a natural map $\ps\colon f^*f_!X\rtarr X$ of spectra over $A$ such that $\ps = f^*\si$ when $X = f^*Y$. 

To pass to homotopy categories, we take two steps. Factoring $f$ as a composite of a homotopy equivalence and an $h$-fibration, we see that we may assume that
$f$ is either a homotopy equivalence or an $h$-fibration. In the former case, $f_*$ must be inverse to the equivalence $f^*$ and thus equivalent to $f_!$. Here $\epz$ and $\si$ are equivalences and we may as well define $\ph$ by the commutativity of (\ref{newdia}). In the latter case, we may work in $G\sE_A$. Since $f$ is an $h$-fibration, we have a natural homotopy equivalence $\mu\colon Tf^*Y\rtarr f^*TY$ for $Y\in G\sE_B$. The derived functor $f_!$ is induced by $Tf_!$, and $TX$ is naturally homotopy equivalent to $X$ when $X\in G\sE_A$. The composite 
$$\xymatrix{f^*Tf_!X \htp Tf^*f_! X \ar[r]^-{T\ps} & TX\htp X}$$
gives a natural map $\ps\colon f^*f_!X\rtarr X$ in $hG\sE_A$. When $X = f^*Y$, we have the commutative naturality diagram
$$\xymatrix{
Tf^*f_!f^*Y \ar[d]_{\mu} \ar[rr]^-{T\ps = Tf^*\si} & &Tf^*Y \ar[d]^{\mu}\\
f^*Tf_!f^*Y \ar[rr]_-{f^*T\si} & & f^*TY.}$$
The bottom arrow is the derived version of $f^*\si$ and the composite around the top is the derived version of $\ps$. Using the equivalences of categories of \S13.5, we obtain a natural map $\ps\colon  f^*f_!X \rtarr X$ in $\Ho G\sS_A$. Let $\et\colon \text{Id}\rtarr f_*f^*$ be the unit of the (derived) adjunction $(f^*,f_*)$ and define $\ph\colon f_!X\rtarr f_*X$ to be the adjoint of $\ps$ in $\Ho G\sS_A$, so that $\ph = f_*\ps\com \et$. For $Y\in \Ho G\sS_B$, we have $f^*\si = \ps\colon f^*f_!f^*Y\rtarr f^*Y$. It follows formally that (\ref{newdia}) commutes. Indeed, 
$$\epz\com f^*\si = \epz\com \ps = \ps\com f^*f_!\epz.$$
The adjoint of $\epz\com f^*\si$ is $\si$ since
$$f_*(\epz\com f^*\si)\com\et = f_*\epz\com f_*f^*\si\com \et
=f_*\epz\com \et\com \si = \si,$$
while the adjoint of $\ps\com f^*f_!\epz$ is $\ph\com f_!\epz$ since
\[f_*(\ps\com f^*f_!\epz)\com \et = f_*\ps\com f_*f^*f_!\epz\com\et
=f_*\ps \com\et\com f_!\epz = \ph\com f_!\epz. \qedhere\]
\end{proof}

\section{The Wirthm\"{u}ller isomorphism for manifolds}\label{sec:mfldwirth}

The classical Wirthm\"{u}ller isomorphism in the equivariant stable homotopy category relates induction and coinduction, the left and right adjoints of the restriction functor from $G$-spectra to $H$-spectra. More precisely, it says that for $H$-spectra $X$, there is a natural equivalence of $G$-spectra \begin{equation}\label{Wirth1} F_H(G_+, X) \simeq G_+\wedge_H (X\wedge S^{-L}), \end{equation} where $L$ is the tangent representation at the identity coset in $G/H$ and $S^{-L}$ is the inverse of the invertible $H$-spectrum $\Sigma^\infty S^L$. Here again, ``equivalence'' means isomorphism in the relevant stable homotopy category and is denoted by $\htp$.

One can also think of this in terms of base change functors. Recall from \myref{changestoo} that the category of $H$-spectra is equivalent to the category of $G$-spectra over $G/H$.  The equivalence is given in one direction by applying the functor $G\times_H -$, and in the other by taking the fiber over the identity coset. This equivalence preserves all structure in sight, including the symmetric monoidal and model structures. The map $r\colon G/H\rtarr *$ induces a pullback functor $r^*$ from $G$-spectra to $G$-spectra over $G/H$, and it has left and right adjoints $r_!$ and $r_*$. The functor $r^*$ corresponds under the equivalence to the restriction functor and therefore $r_!$ and $r_*$ correspond to the induction and coinduction functors. In this terminology, the Wirthm\"{u}ller isomorphism (\ref{Wirth1}) takes the form
\begin{equation}\label{Wirth1bis}
r_* X \simeq r_!(X\wedge_{G/H} C_{G/H})
\end{equation}
for $G$-spectra $X$ over $G/H$, where $C_{G/H}= \io_!S^{-L}$, $\io:H\subset G$ (see \myref{eyeeye}). 

We think of $G/H\rtarr *$ as the simplest kind of a bundle with a compact manifold as fiber, and we generalize (\ref{Wirth1bis}) to maps $p\colon E\rtarr B$ that are equivariant bundles with a smooth closed manifold $M$ as a fiber.  We discuss the case $B=*$ in this section and prove the general case in the next.  However, it is convenient to begin by describing the form of the map that gives the equivalence in general.  For that, we require a $G$-spectrum $C_p$ over $E$ together with an equivalence 
\begin{equation}\label{Wobj} 
\xymatrix{\al_p\colon p_!C_p \ar[r]^-{\htp} & D(p_! S_{E})}
\end{equation} 
that identifies the dual of $p_!S_E$. 
We call $C_p$,\noteindex{Cp@$C_p$} together with $\al_p$, a {\em Wirthm\"uller object}.\index{Wirthmuller object@Wirthm\"uller object}

In \cite{FHM}, Fausk, Hu, and May give a categorical discussion of equivalences of Wirthm\"uller type, including a simplifying formal analysis that describes the minimal amount of information that is needed to prove such a result. In particular, given a Wirthm\"uller object $C_p$, they define a canonical candidate
$$\om_p\colon p_*X\rtarr p_!(X\sma_E C_p)$$
for an equivalence, namely the composite displayed in the commutative diagram
\begin{equation}\label{omega}
\xymatrix{
p_*X \htp p_*X\sma_B D(S_B)\ar[dd]_{\om_{p}} \ar[rr]^-{\text{id}\sma_BD(\si)}
& & p_*X\sma_B D(p_!S_E) \\
& & p_*X\sma_B p_!C_p  \ar[u]_{\text{id}\sma_B \al_p}^{\htp}\\
p_!(X\sma_E C_p) & & p_!(p^*p_*X\sma_E C_p) \ar[u]_{\htp}\ar[ll]^{p_!(\epz\sma_E\text{id})}.}
\end{equation}
The maps $\si\colon p_!S_E\htp p_!p^*S_B \rtarr S_B$ and $\epz\colon p_*p^*X\rtarr X$
are given by the counits of the adjunctions $(p_!,p^*)$ and $(p^*,p_*)$. The 
arrow labelled $\htp$ is an equivalence given by the derived version of the
projection formula (\ref{four}) that is proven in \myref{Wirthmore}. 

When $M$ is a smooth closed $G$-manifold and $r$ is the map $\colon M\rtarr *$,  we write $C_M$ for a Wirth\-m\"uller object $C_r$ and we write $\om_M$ for $\om_r$. It is easy to describe $C_M$. Let $\ta$ be the tangent $G$-bundle of $M$. Embed $M$ in a $G$-representation $V$ and let $\nu$ be the normal $G$-bundle of the embedding.  By Atiyah duality, the union $M_+$ of $M$ and a
disjoint basepoint is $V$-dual to the Thom $G$-space $T\nu$. A detailed equivariant proof is given in \cite[III\S5]{LMS}, but we require little 
beyond the mere statement.

For a $G$-vector bundle $\xi$ over a $G$-space $B$, let $S^{\xi}$ denote the fiberwise one-point compactification of $\xi$, with section given by the points at infinity.  This ex-$G$-space over $B$ must not be confused with the Thom complex $T\xi$. The latter is obtained by identifying the section to a point and is precisely $r_!S^{\xi}$, $r\colon B\rtarr *$.

\begin{defn}\mylabel{CM}
Define $C_M$ to be the $G$-spectrum $\SI_M^{-V}\SI^{\infty}_M S^{\nu}$ over $M$.
\end{defn}

\begin{rem}\mylabel{tauinv} 
By \myref{spaceFDT}, the suspension $G$-spectrum
$\SI^{\infty}_MS^{\ta}$ is invertible. Visibly, $C_M$ is its inverse.
\end{rem}

\begin{lem}\mylabel{starting} 
There is an equivalence $\al_M\colon r_!C_M \rtarr D(r_!S_{M})$,
$r\colon M\rtarr *$.
\end{lem}

\begin{proof}
Since $S_{M}(V) = M\times S^V$, $r_!S_{M} = \SI^{\infty}M_+$. Since $r_!S^{\nu} = T\nu$ and $r_!$ commutes with shift desuspension functors, $r_!C_M$ is equivalent to $\SI^{-V}\SI^{\infty}T\nu$. There is a canonical evaluation map $\text{ev}\colon T\nu\wedge M_+\rtarr S^V$ of a duality \cite[p.152]{LMS}. Explicitly, using the diagonal of $M$ and the zero section of $\nu$ we obtain an embedding of $M$ in $\nu\times M$ with trivial normal bundle $M\times V$, and $\text{ev}$ is composite of the Pontryagin-Thom map associated to this embedding and the projection $M_+\sma S^V\rtarr S^V$. We apply the functor 
$\SI^{-V}\SI^{\infty}$ to obtain
$$\SI^{-V}\SI^{\infty}T\nu \sma \SI^{\infty}M_+ \htp \SI^{-V}\SI^{\infty}(T\nu\sma M_+) \rtarr \SI^{-V}\SI^{\infty} S^V\htp S.$$
Atiyah duality states that the adjoint of this map is an equivalence from $\SI^{-V}\SI^{\infty}T\nu$ to $D(M_+)$. This is the required map $\al_M$. 
\end{proof}

We shall prove the following result in \S\ref{sec:wirthpf}.

\begin{thm}[The Wirthm\"uller isomorphism for manifolds]\mylabel{thmW1}
For $G$-spectra $X$ over $M$ and $r\colon M\rtarr *$, the map
$$\om_M\colon r_*X\rtarr r_!(X\sma_M C_M)$$
is a natural equivalence in the homotopy category $\Ho G\sS$ of $G$-spectra.
\end{thm}

In an earlier draft of this paper, we thought we could reduce the general case of \myref{thmW1} to the special case $M = G/H$. However, instead of leading to a simplifiction, the argument we had in mind leads to an interesting relative version of the Wirthm\"uller isomorphism. Its starting point is the following observation.

\begin{lem}\mylabel{WMGH} 
Let $i\colon N\rtarr M$ be an embedding of smooth closed $G$-manifolds and let $\nu_{M,N}$ be the normal bundle of $i$. Then $C_N$ is equivalent to $S^{\nu_{M,N}}\sma_N i^*C_M$.  
\end{lem}

\begin{proof}
An embedding of $M$ in a representation $V$ restricts along $i$ to an embedding of $N$ in $V$, and $i^*\nu_M\oplus \nu_{M,N} \iso \nu_{N}$. Commutation relations in \myref{SIfSI} give that
$$i^*\SI_M^{-V}\SI^{\infty}_M S^{\nu_M} 
\htp \SI_{N}^{-V}\SI^{\infty}_{N}i^*S^{\nu_M}.$$
The conclusion follows after smashing with $S^{\nu_{M,N}}$.
\end{proof}

\begin{cor}[The relative Wirthm\"uller isomorphism]
Let $i\colon N\rtarr M$ be a smooth embedding of closed $G$-manifolds. For $G$-spectra $X$ over $N$, there is a natural equivalence
$$\om_{M,N}\colon r_*i_* X\rtarr r_*i_!(X\sma_N S^{\nu_{M,N}}).$$
\end{cor}

\begin{proof}
Here $r\colon M\rtarr *$. Write $q = r\com i\colon N\rtarr *$. Then $q_*\htp r_*i_*$ and $q_!\htp r_!i_!$. Define $\om_{M,N}$ by commutativity of the diagram of equivalences
$$\xymatrix{
q_*X \ar[d]_{\htp} \ar[r]^{\om_N} & q_!(X\sma_N C_N) \ar[d]^{\htp}\\
r_*i_*X  \ar[d]_{\om_{M,N}} 
&  r_!i_! (X\sma_N S^{\nu_{M,N}}\sma_N i^*C_M) \ar[d]^{\htp}\\
r_*i_!(X\sma_N S^{\nu_{M,N}})\ar[r]_-{\om_M} 
& r_!(i_!(X\sma_N S^{\nu_{M,N}})\sma_M C_M).}$$
Here the derived version of the projection formula (\ref{four}) gives the lower right equivalence.\end{proof}

We explain the strategy that we have not implemented for deducing the Wirth\-m\"uller 
isomorphism for $M$ from the Wirthm\"uller isomorphism for orbits. 

\begin{rem} 
One can use relative Atiyah duality to define an intrinsic map $\al_{M,N}$ and use $\al_{M,N}$ to define a map $\om_{M,N}$ directly. One can then obtain the displayed diagram by a chase. If one could prove directly that $\om_{M,N}$ was an equivalence, then, using the invertibility of $\SI^{\infty}_NS^{\nu_{M,N}}$, one could deduce that $\om_M$ is an equivalence on all $i_!X$ if $\om_N$ is an equivalence. By \myref{detect}, $\om_M$ is an isomorphism for all $Y$ if it is an isomorphism for all $Y$ in the detecting set $\sD_M$ of \myref{detecting}.  Those $Y$, namely the $S^{n,b}_H$ for $b\in M$ and $H\subset G_b$, are of the form $\tilde{b}_!X$, where $\tilde{b}\colon G/G_b\rtarr M$ is the inclusion of the orbit of $b$ and $X$ is a $G$-spectrum over $G/G_b$. Thus the Wirthm\"uller isomorphism for orbits would imply the Wirthm\"uller isomorphism for $M$. 
\end{rem}

\section{The fiberwise Wirthm\"uller isomorphism}\label{sec:fibwirth}

As in \S\ref{sec:bdlconstr}, let $G$ be a quotient $\GA/\PI$, where $\PI$ is a normal subgroup of a compact Lie group $\GA$.  Let $M$ be a smooth closed $\GA$-manifold and let $p\colon E\rtarr B$ be a $(\PI;\GA)$-bundle with fiber $M$.  This means that $p$ has an associated principal $(\PI;\GA)$-bundle $\pi\colon  P\rtarr B$ and $p$ is the associated $G$-bundle $E = \T{M}\rtarr P/\PI=B$.  We apply the functor $\U{_M}$ to the Wirthm\"uller object $C_M$ to obtain the Wirthm\"uller object $C_p$, and we apply $\U{_M}$ to $\al_M$ to obtain the required equivalence $\al_p$. This means that the Wirthm\"uller object for $p$ is obtained by inserting the Wirthm\"uller object for $M$ fiberwise into the functor $\T{(-)}$.

\begin{defn} 
Define $C_p$ to be the $G$-spectrum $\U{_{M}{i^*C_{M}}}$ over $E$, where $i^*$ is the change of universe functor associated to the inclusion of the $\PI$-trivial $\GA$-universe in the complete $\GA$-universe.
\end{defn}

\begin{rem}\mylabel{tauinvtoo}  
Recall \myref{tauinv}. By \myref{spaceFDT}, the suspension $G$-spectrum $\SI^{\infty}_E(P\times_{\PI} S^{\ta})$ is invertible. The Wirthm\"uller object $C_p$ is its inverse.
\end{rem}

\begin{lem}\mylabel{startingtoo} 
There is an equivalence $\al_p\colon p_!C_p \rtarr D(p_!S_{E})$.
\end{lem}

\begin{proof}
We define $\al_p$ to be the composite
\begin{equation}\label{alf}
\xymatrix{
p_!\U{_M{i^*C_M}} \ar[r] & \U{_{*}{i^*r_!C_M}} \ar[rr]^-{\U{_{*}{i^*\al_M}}} & & 
\U{_{*}{i^*D(r_!S_M)}} \ar[r] & D(p_!S_{E}).\\}
\end{equation}
The left arrow is given by the first equivalence of \myref{Pfibers} and the last equivalence of \myref{chvschuni}. The middle arrow is an equivalence since $\al_M$ is one. The right arrow is the following composite equivalence,
\begin{align*}
\U{_{*}{i^*D(r_!S_M)}} & \htp  \U{_{*}{D(i^*r_!S_M)}} \\
 & \htp  \U{_{*}{D(r_!i^*S_M}}) && \text{by Propositions \ref{chvschuni} and \ref{Nfreeii}}\\
 & \htp  D(\U{_{*}{r_!i^*S_M}}) && \text{by \myref{PSmash}}\\
 & \htp  D(p_!\U{_{M}{i^*S_M}}) && \text{by \myref{Pfibers}}\\
 & \htp  D(p_!S_E) && \text{by \myref{PSmash}}.
\end{align*}
For the first displayed equivalence, $r_!S_M \simeq \SI^{\infty}M_+$ by \myref{lesstrivial}, hence 
$$D(r_!S_M) \simeq D(\SI^{\infty}M_+)\simeq F(M_+,S).$$
For based $\GA$-spaces $T$, $i^*F(T,S)\iso F(T,i^*S)$ by inspection. If
$T$ is a based $\GA$-CW complex, this is an isomorphism of Quillen right
adjoints and so descends to homotopy categories. Again by inspection,
$$i^*\SI^{\infty}\iso \SI^{\infty}\colon G\sK_*\rtarr G\sS^{\PI-\text{triv}}.$$
This isomorphism passes to homotopy categories since both sides take 
$q$-equi\-va\-lences to level $q$-equi\-va\-lences. Therefore $i^*S\htp S$
and $F(T,i^*S)\htp D(i^*\SI^{\infty}T)$.
\end{proof}

\begin{thm}[The fiberwise Wirthm\"uller isomorphism]\mylabel{thmW2}
For $G$-spec\-tra $X$ over $E$, the map
$$\om_p\colon p_* X \rtarr p_!(X\sma_E C_p)$$
is a natural equivalence of $G$-spectra over $B$.
\end{thm}

\begin{proof}
The action of $G_b$ on the fiber $E_b\iso \rh_b^*M$ of $b\in B$ is smooth, hence the Wirthm\"uller isomorphism for manifolds gives the result for $r\colon E_b\rtarr *$.  We claim that the restriction $\al_b$ of $\al_p$ to the fiber over $b$ is an equivalence 
$$\al_{E_b}\colon r_!C_{E_b}\rtarr D(r_!S_{E_b})$$
of $G_b$-spectra of the form used to prove \myref{thmW1} for $r$. Indeed, with $p_b = r$, $i\colon E_b\subset E$ and $\io\colon G_b\subset G$, the derived version of \myref{Johann2} and \myref{tEtB} give that the source of $\al_b$ is
$$ (p_!\U{_M{C_M}})_b \htp r_!i^*\io^*\U{_M{C_M}}\htp r_!\rh_b^*C_M \iso r_!C_{E_b}.$$
For the last isomorphism we must view the representation $V$ of $\GA$ that appears in the definition of $C_M$ as a representation of $G_b$ by pullback along $\rh_b$. Similarly, using \myref{reassuring} and the derived version of \myref{Johann2}, the target of $\al_b$ is
$$ D(p_!S_{G,E})_b \htp D((p_!S_{G,E})_b) \htp D(r_!i^*\io^*S_{G,E}) \htp D(r_!S_{G_b,E_b}).$$
In view of the role of $\al_M$ in the definition of $\al_p$, diagram chases from the definitions show that $\al_b$ agrees under these equivalences with the $G_b$-equivalence $\al_{E_b}$. 

Now, looking at the definition of $\om_p$ (\ref{omega}), we see that, aside from the equivalence $\al_p$, its constituent maps are just counits of adjunctions and derived isomorphisms coming from the closed symmetric monoidal structures. By \myref{reassuring} and the derived versions of commutation relations in \myref{Johann}, these maps restrict on fibers to maps of the same form. Therefore the restriction
$$\om_b\colon (p_* X)_b \rtarr (p_{!}(X\sma_E C_p))_b$$ 
of $\om_p$ to the fiber over $b$ can be identified with the map of $G_b$-spectra
$$ \om_{E_b}\colon r_* X_b \rtarr r_!(X_b\sma_{E_b} C_{E_b}).$$
This map is an equivalence of $G_b$-spectra by \myref{thmW1}. Since equivalences of $G$-spectra over $B$ are detected fiberwise, this implies that $\om_p$ is an equivalence.
\end{proof}

\begin{rem} 
When $\GA = G\times \PI$ and only $\PI$ acts on $M$, one can think of $p\colon E\rtarr B$ as a topological $G$-bundle with a reduction of its structural group to a suitably large compact subgroup $\PI$ of the group of diffeomorphisms of $M$.  Our fiberwise Wirthm\"uller isomorphism theorem is a variant of the main theorem, \cite[4.8]{Hu}, of a paper of Po Hu. She worked with $\text{Diff}(M)$ itself as an implicit structure group, without use of an auxiliary group $\PI$ and without an ambient group $\GA$. That bundle theoretic framework leads to formidable complications, hence her arguments are very much more difficult than ours. Her result is both more and less general than the specialization of ours to the case $\GA = G\times \PI$: it allows bundles that might not admit a single compact structure group $\PI$, but it requires the base spaces to be $G$-CW
complexes with countably many cells. It does not handle more general group extensions.
\end{rem} 

\section{The Adams isomorphism}\label{sec:adams}

Let $N$ be a normal subgroup of $G$ and let $\epsilon\colon G\rtarr J$ be the quotient by $N$. The conjugation action of $G$ on $N$ induces an action of $G$ on the tangent space of $N$ at the identity element, giving us the adjoint representation $A = A(N;G)$. Let $(i_*,i^*)$ be the change of universe adjunction associated to the inclusion $i\colon q^*\sV_J\rtarr \sV_{G}$ of 
the complete $J$-universe $\sV_J$ as the universe of $N$-trivial 
representations in the complete $G$-universe $\sV_{G}$.

Recall the discussion of $N$-free $G$-spectra from \S\ref{sec:PIFree}, where $\PI$ and $\GA$ played the roles of $N$ and $G$.

\begin{thm}[Adams isomorphism] 
For $N$-free $G$-spectra $X$ in $G\sS^{\text{$N$-triv}}$, there is a natural equivalence
\[ X/N\htp (i^*\Sigma^{-A}i_* X)^N \]
in $\Ho J\sS^{\text{$N$-triv}}.$
\end{thm}

We shall derive this by applying the fiberwise Wirthm\"{u}ller isomorphism to the quotient $G$-map $p\colon  E(N;G) \rtarr B(N;G)$, where $E(N;G)$ is the
universal $N$-free $G$-space and $B(N;G) = E(N;G)/N$.  To place ourselves in the bundle theoretic context of the previous section, we give another description of $p$, following \cite[II\S7]{EHCT}. It is formal and would similarly identify $p\colon E\rtarr E/N$ for any $N$-free $G$-space $E$.  Let $\GA = G\ltimes N$ be the semi-direct product of $G$ and $N$, where $G$ acts by conjugation on $N$. Write $\PI$ for the normal subgroup $\{e\}\ltimes N$ of $\Gamma$. We then have an extension
\[1\rtarr \Pi \rtarr \Gamma \stackrel{\theta}{\rtarr} G \rtarr 1,\]
where $\theta(g,n)=gn$. Give $N$ the $\Gamma$-action $(g,n)\cdot m = gnmg^{-1}$. Then $N\cong \Gamma/G$ as $\Gamma$-spaces, where we view $G$ as the subgroup $G\ltimes\{e\}$ of $\Gamma$. The composite
\[E(N;G)\cong \theta^*E(N;G)\times_\Pi (\Gamma/G) 
\rtarr \theta^* E(N;G)\times_\Pi * \cong B(N;G)\]
induced by $\GA/G\rtarr *$ is $p$.  Since $\theta^*E(N;G)$ is a $\Pi$-free $\Gamma$-space, we see that $p$ is a bundle with fiber $\Gamma/G\cong N$ to which the fiberwise Wirthm\"uller isomorphism applies. We must identify the relevant Wirthm\"uller object. We write $r$ for the map $E(N;G)\rtarr *$.

\begin{prop}\mylabel{WirthAd}
The Wirthm\"{u}ller object $C_{p}$ is $\rE^*S^{-A}$.
\end{prop}

\begin{proof}
The tangent bundle of $\Gamma/G\cong N$ is the trivial bundle $N\times A$ \cite[p.\,99]{LMS}. Indeed, let $\Gamma$ act on $A$ via the projection $\epsilon\colon  \Gamma\rtarr G$,  $\epz(n,g)=g$. We obtain a $\Gamma$-trivialization of the tangent bundle of $\Gamma/G$ by sending $(n,a)\in N\times A$ to $d_eL_n(a)$, where $d_eL_n$ is the differential at $e$ of left translation by $n$. It follows that the tangent bundle along the fibers of $p$ is also trivial: 
\[\theta^*E(N;G)\times_N(\Gamma/G\times A)
\cong (\theta^*E(N;G)\times_N\Gamma/G))\times A\cong E(N;G)\times A.\]
Thus the spherical bundle of tangents along the fiber is $E(N;G)\times S^A=\rE^*S^A$, and the inverse of its suspension $G$-spectrum over $E(N;G)$ is $\rE^*S^{-A}$. In view of \myref{tauinvtoo}, this gives the conclusion.
\end{proof}

\begin{proof}[Proof of the Adams isomorphism]
Let $X\in G\sS^{\text{$N$-triv}}$ be $N$-free. Applying the fiberwise Wirthm\"{u}ller isomorphism to the $G$-spectrum $\rE^*i_*X$ over $E(N;G)$ and using that $C_p$ is $\rE^*S^{-A}$, we obtain a natural equivalence
$$ p_*\rE^*i_* X \htp p_!(\rE^*i_*X\sma_{E(N;G)} \rE^*S^{-A})$$
of $G$-spectra over $B(N;G)$.  Write $\bar{r}$ for the map $B(N;G)\rtarr *$, so that $\bar{r}\com p = r$. Applying the functor $\rB_!((i^*(-))^N)$ to the displayed equivalence, we obtain a natural equivalence
\[\rB_!((i^*p_*\rE^*i_*X)^N) \htp \rB_!((i^*p_!(\rE^*i_*X\wedge_{E(N;G)}\rE^*S^{-A}))^N)\]
in $\Ho J\sS^{\text{$N$-triv}}$. We proceed to identify both sides. The source is
\begin{align*}
\rB_!((i^*p_*\rE^*i_*X)^N)
&\htp  \rB_!((p_*\rE^*i^*i_*X)^N) && \text{by \myref{chvschuni}}\\
&\htp  \rB_!((p_*\rE^*X)^N) && \text{by \myref{Nfreeii}}\\
&\htp  \rB_!((p_!\rE^*X)/N) && \text{by \myref{ouch}}\\
&\htp  (\rB_!p_!\rE^*X)/N && \text{by \myref{chvsfixorbit}}\\
&\htp  (\rE_!\rE^*X)/N && \text{by functoriality}\\
&\htp  X/N. && \text{by \myref{Nfree}}.
\end{align*}
The target is
\begin{align*}
\rB_!((i^*p_!(\rE^*i_*X &\wedge_{E(N;G)}\rE^*S^{-A}))^N)\\
&\htp  \rB_!((i^*p_!\rE^*\Sigma^{-A}i_*X)^N) && \text{by \myref{symmonhtp}}\\
&\htp (\rB_!i^*p_!\rE^*\Sigma^{-A}i_*X)^N && \text{by \myref{chvsfixorbit}}\\
&\htp  (i^*\rB_!p_!\rE^*\Sigma^{-A}i_*X)^N && \text{by Propositions \ref{chvschuni} and \ref{Nfreeii}}\\
&\htp  (i^*\rE_!\rE^*\Sigma^{-A}i_*X)^N && \text{by functoriality}\\
&\htp  (i^*\Sigma^{-A}i_*X)^N && \text{by \myref{Nfree}.}\qedhere
\end{align*}
\end{proof}

\begin{rem} 
In outline, the proof just given is essentially that indicated by 
Po Hu \cite[pp 81--99]{Hu}. However, her argument, although more conceptual, 
is a good deal longer and more complicated than the original 
proof in \cite[pp 96--102]{LMS}.
\end{rem} 

\section{Proof of the Wirthm\"uller isomorphism for manifolds}\label{sec:wirthpf}

We prove \myref{thmW1} here.  Thus consider $r\colon M\rtarr *$ for a smooth compact $G$-manifold $M$. With $C_M = \SI_M^{-V}\SI^{\infty}_M S^{\nu}$, the diagram (\ref{omega}) displays a canonical map 
$$\om = \om_M\colon r_*(X)\rtarr r_!(X\sma_M C_M)$$
of $G$-spectra, where $X$ is a $G$-spectrum over $M$. We must show that $\om$ is an equivalence. In outline, we follow the pattern of proof explained in \cite{FHM} and illustrated in the case $M = G/H$ in \cite{MayW}, but the details are very different from those applicable in that special case. 

We first describe a formal reduction implied by the results of \cite{FHM}. Consider the set $\sD_M$ of detecting objects in $\Ho G\sS_M$ that is specified in \myref{detecting}. The objects in $\sD_M$ are compact, by \myref{compact}, and dualizable.  We have the analogous detecting set 
$\sD_*$ of compact objects in $\Ho G\sS$. For $Y$ in $\sD_*$, $r^*Y$ is dualizable and it follows formally, by \cite[2.1.3(d)]{HPS}, that $r^*Y$ is compact (in the sense of \myref{compact}). Therefore $r_*$, as well as $r_!$, preserves coproducts \cite[7.4]{FHM}. This verifies the hypotheses of the 
formal Wirthm\"uller isomorphism theorem, \cite[8.1]{FHM}, and that result 
shows that $\om$ will be an equivalence for all $G$-spectra $X$ over $M$ if 
it is an equivalence for those $X$ in $\sD_M$.

Such $X$ are of the form $S^{n,m}_H=\widetilde{m}_!\io_!S^n_H$, where $n\in \bZ$, $m\in M$, $H\subset G_m$, and $\io$ is the inclusion of $G_m$ in $G$.  By commutation with suspension, we can assume that $n\geq 0$. Then $X$ is of the form $\SI^{\infty}_MK$ for an ex-$G$-space $K$ over $M$, and $X$ can be any such $G$-spectrum over $M$ in the rest of the proof. By \cite[6.3]{FHM}, it suffices to construct a map ${\xi_X} \colon r^*r_!(X\wedge_M C_M)\rtarr X$ such that certain diagrams commute. To be precise, let $\si$ and $\zeta$ be the counit and unit of the $(r_!,r^*)$ adjunction, note that $r^*S\iso S_M$, and define maps $\tau = \ta_S$ and $\xi = \xi_{S_M}$ by commutativity of the diagrams
\begin{equation}\label{diagdis}
\xymatrix{S\ar[r]^-{\ta} \ar[d]_{\simeq} & r_!C_M\ar[d]^{\alpha_M}\\
DS \ar[r]_-{D\si} & Dr_!r^*S}
\quad\text{and}\quad
\xymatrix{r^*r_!C_M\ar[rr]^-{\xi}\ar[d]_{r^*\alpha_M} & & S_M \ar[d]^\simeq\\
r^*Dr_!S_M \ar[r]_-{\simeq}  & Dr^*r_!S_M\ar[r]_-{D\zeta} & DS_M}
\end{equation}
Then define $\ta_Y$ for a general $G$-spectrum $Y$ to be the composite 
\begin{equation}\label{tau}
\xymatrix{
\ta_Y\colon Y\simeq Y\wedge S \ar[r]^-{\text{id}\wedge \tau} & Y\wedge r_!C_M  
\htp r_!(r^*Y\wedge_M C_M)\\}
\end{equation}
and define $\xi_{r^*Y}$ for the $G$-spectrum $r^*Y$ over $M$ to be the composite
\begin{equation}\label{xi}
\xymatrix{
\xi_{r^*Y}\colon r^*r_!(r^*Y\wedge_M C_M)\simeq
r^*Y \wedge_M r^*r_!C_M \ar[r]^-{\text{id}\wedge \xi} & r^*Y\wedge_M S_M 
\simeq  r^*Y.\\}
\end{equation}
Here the equivalences are given by the derived versions of (\ref{one})
and the projection formula (\ref{four}). With these notations, we shall 
prove the following result.

\begin{prop}\mylabel{prop:spec}
For $X=\Sigma^\infty_M K$, there is a map 
$$\xi_X\colon  r^*r_!(X\wedge_M C_M)\rtarr X$$ 
such that the composite
\begin{equation}\label{eq:trieq}
\xymatrix{
r_!(X\wedge_M C_M) \ar[d]^-{\ta_{r_!(X\sma_M C_M)}} \\
r_!(r^*r_!(X\wedge_M C_M)\sma_M C_M)
\ar[d]^-{r_!(\xi_X\sma_M\text{id})}\\   
r_!(X\wedge_M C_M)\\}
\end{equation}
is the identity map (in $\Ho G\sS$)
and, for any map  $\theta\colon  r^*Y \rtarr X$ of $G$-spectra over $M$, the 
following diagram commutes in $\Ho G\sS_M$.
\begin{equation}\label{eq:nat}
\xymatrix{r^*r_!( r^*Y\wedge_M C_M)\ar[r]^-{\xi_{r^*Y}}
\ar[d]_{r^*r_!(\theta\wedge \text{id})} & r^*Y \ar[d]^\theta\\
r^*r_!(X\wedge_M C_M) \ar[r]_-{\xi_X} & X}
\end{equation}
\end{prop}

This will complete the proof of the theorem by the cited reduction 
from \cite{FHM}.

\begin{cor} For $X$ in $\sD_M$, $\om_M\colon r_*(X)\rtarr r_!(X\sma_M C_M)$
is an equivalence with inverse the adjoint of $\xi_X$.
\end{cor}
\begin{proof} Taking $Y$ to be $r_*X$ and $\theta$ to be the counit 
of the $(r^*,r_*)$ adjunction in (\ref{eq:nat}), the conclusion
is a direct application of \cite[6.3]{FHM}. 
\end{proof}

Thus it suffices to prove \myref{prop:spec}.  We shall construct the map 
$\xi_X$ and prove that it satisfies the stated properties by reducing to 
space level considerations.  We begin with a space level description 
of the maps $\tau$ and $\xi$ displayed in (\ref{diagdis}),
and we need some space level notations.

\begin{notns}\mylabel{notns}
Recall that $\nu$ denotes the normal bundle of $M$ and that we have the 
duality  map $\text{ev}\colon T\nu\wedge M_+\rtarr S^V$ specified in the 
proof of \myref{starting}. Also, recall that 
\[r_!K = K/s(M),  \quad r^*T = T_M = M\times T, 
\quad\text{and}\quad
r_*K = \text{Sec}(M,K)\]
for any based $G$-space $T$ and any ex-$G$-space $(K,p,s)$ over $M$.
In particular, 
$$r_!S^{\nu} = T\nu, \quad r_!S^0_M = M_+, \quad \text{and}\quad 
r_*r^*T \iso F(M_+,T).$$
Therefore the adjoint $\widetilde{\text{ev}}\colon T\nu\rtarr F(M_+,S^V)$ 
is a map $r_!S^{\nu}\rtarr r_*S^V_M$.
Let $t\colon S^V\rtarr r_!S^\nu$ be the Pontryagin-Thom 
construction and $k$ be the composite
$$\xymatrix{k\colon r^*r_!S^\nu \ar[r]^-{r^*\tilde{\text{ev}}}
& r^*r_*S^V_M  \ar[r]^-{\epz} & S^V_M,\\}$$
where $\epz\colon r^*r_*\rtarr \text{id}$ is the counit of the adjunction 
$(r^*,r_*)$; note that, in general, $\epz$ is just the evaluation map 
$M\times \text{Sec}(M,K)\rtarr K$.
\end{notns}

Recall from Propositions \ref{SIfSI} and \ref{SISISI2} that we can commute 
suspension spectrum functors past smash products and base change functors.

\begin{lem}\mylabel{lemma:units} With these definitions of $t$ and $k$,
\[\tau\simeq \Sigma^{-V}\Sigma^\infty t\colon S\iso \Sigma^{-V}\Sigma^\infty S^V
\rtarr  \Sigma^{-V}\Sigma^\infty r_!S^\nu \htp r_!C_M\]
and
\[\xi \simeq \Sigma^{-V}_M\Sigma^\infty_M k\colon 
r^*r_!C_M\htp \Sigma^{-V}_M\Sigma^\infty_M r^*r_!S^{\nu} \rtarr
\Sigma^{-V}_M\Sigma^\infty_M S^V_M \htp S_M.\]
\end{lem}
\begin{proof} By \cite[III.5.2]{LMS},
the dual of $t$ is the projection $\de\colon M_+\rtarr S^0$.
This means that the following diagram is stably homotopy commutative.
$$\xymatrix{
S^V\sma M_+ \ar[r]^-{t\sma\text{id}} \ar[d]_{\text{id}\sma \de}
& T\nu\sma M_+ \ar[d]^{\text{ev}}\\
S^V\sma S^0 \ar@{=}[r] & S^V\\}$$
Here $\de\colon M_+ =r_!r^*S^0\rtarr S^0$ is the counit of the space level 
adjunction $(r_!,r^*)$, and we can identify $\SI^{\infty}\de$ with the 
counit $\si\colon \SI^{\infty}M_+\iso r_!r^*S\rtarr S$ of the spectrum level 
adjunction $(r_!,r^*)$.  Applying $\Sigma^{-V}\Sigma^\infty$ to the 
diagram and passing to adjoints, the right vertical arrow becomes 
$$\al_M\colon r_!C_M \htp \SI^{-V}\SI^{\infty}T\nu\rtarr D(M_+) = Dr_!r^*S,$$ 
by the proof of \myref{starting}. Comparing the resulting diagram with 
the diagram that defines $\ta$, we conclude that 
$\tau\simeq \Sigma^{-V}\Sigma^\infty t$. 

For the identification of $\xi$, we consider the 
composite equivalence $r^*r_!C_M\htp Dr^*r_!S_M$ in the diagram 
that defines $\xi$ to be an identification of the dual of 
$r^*r_!S_M$. To identify the dual of $\ze$
modulo that identification, we observe that the following diagram is 
commutative.  
$$\xymatrix{
r^*r_!S^{\nu}\sma_M S^0_M \ar[rr]^-{k\sma\text{id}} 
\ar[d]_{\text{id}\sma\zeta}
& & S^V_M\sma_M S^0_M \ar@{=}[d] \\
r^*r_!S^{\nu}\sma_M r^*r_!S^0_M \ar[r]_{\cong}
& r^*(r_!S^{\nu}\sma r_!S^0_M) \ar[r]_-{r^*\text{ev}} & r^*S^V\\}$$
Indeed, recalling that $k = \epz\circ r^*\tilde{\text{ev}}$ and rewriting 
the diagram in more familiar notation, it becomes
$$\xymatrix{
M\times T\nu \ar[r]^-{\text{id}\times \tilde{\text{ev}}} \ar[d]_-{\text{id}\times \zeta}
& M\times \text{Sec}(M,M\times S^V) \ar[d]^{\epz}\\
M\times (T\nu\sma M_+) \ar[r]_-{\text{id}\times \text{ev}} 
& M\times S^V,\\}$$
and both composites send $(m,x)$ to $(m,\text{ev}(x\sma m))$.
Applying $\Sigma^{-V}\Sigma^{\infty}$ to the first diagram and 
comparing with the definition of $\xi$, we conclude that 
$\xi \simeq \Sigma^{-V}_M\Sigma^\infty_M k$. 
\end{proof}

The following space level result will imply \myref{prop:spec}. 

\begin{prop}\mylabel{prop:space}
Let $K$ be a $G$-space over $M$.  There is a natural map 
$$u_K\colon  r^*r_!(K\wedge_M S^\nu) \rtarr  \Sigma^V_M K$$ 
in $\Ho G\sK_M$ which satisfies the following properties.
\begin{enumerate}[(i)]
\item When $K=S^0_M$, 
$u_K\simeq k\colon r^*r_!S^{\nu}\rtarr S^V_M$.  
\item For a based $G$-space $T$, the following diagram commutes in 
$\Ho G\sK_M$.
\[\xymatrix{ r^*r_!(T_M\wedge_M K\wedge_M S^\nu) \ar[r]^\simeq 
\ar[d]_{_{u_{T_M\sma_M K}}} & T_M \wedge_M r^*r_!(K\wedge_M S^\nu) 
\ar[d]^{\text{id}\wedge u_K}\\
\Sigma^V_M(T_M\wedge_M K) \ar[r]^\cong & T_M \wedge_M \Sigma^V_M K}\]
Here the top equivalence is given by (\ref{one}) and (\ref{four}) in 
$\Ho G\sK_M$.
\item The following diagram commutes in $\Ho G\sK_*$.
\[\xymatrix{r_!(K\wedge_M S^\nu)\wedge S^V 
\ar[rr]^-{\text{id}\wedge t}\ar[d]_\simeq 
&& r_!(K\wedge_M S^\nu)\wedge r_! S^\nu \ar[d]^\simeq \\
r_!(\Sigma^V_M K\wedge_M S^\nu) 
&& r_!(r^*r_!(K\wedge_M S^\nu)\wedge_M S^\nu) \ar[ll]^-{r_!(u_K\wedge \text{id})}}\]
Here the right vertical equivalence is given by (\ref{four}) in $\Ho G\sK_M$.
The left vertical equivalence is the composite
\[r_!(K\wedge_M S^\nu)\wedge S^V 
\simeq r_!(K\wedge_M S^\nu\wedge_M S^V_M)
\simeq r_!(K\wedge_M S^V_M\wedge_M S^\nu),\]
where the second equivalence is obtained by moving the copy of $S^\tau$ 
from $S^V_M\cong S^\tau\wedge_M S^\nu$ and amalgamating it with the displayed 
copy of $S^\nu$.
\end{enumerate}
\end{prop}

\begin{proof}[Proof of \myref{prop:spec}]
Let $X=\Sigma^\infty_M K$.  Define $\xi_X$ to be the map
\[\xymatrix{
r^*r_!(X\wedge_M C_M) 
\simeq \Sigma^{-V}_M \Sigma^\infty_M r^*r_!(K\wedge_M S^\nu) 
\ar[rr]^-{\Sigma^{-V}_M\Sigma^\infty_M u_K} &&
\Sigma^{-V}_M \Sigma^\infty_M \Sigma^V_M K\iso X.}\]
Using \myref{prop:space}(ii), we see that 
$\xi_{\Sigma^V_MX}$ can be identified with $\Sigma^V_M\xi_X$, 
which in turn can be identified with $\Sigma^{\infty}_Mu_K$. 
To show that the diagram (\ref{eq:nat}) commutes, it suffices to 
show that the diagram obtained from it by applying $\Sigma^V_M$ 
commutes. We have just identified the lower horizontal arrow of 
the resulting diagram in space level terms.  Similarly, the 
definition (\ref{xi}) of its upper horizontal arrow, together with 
\myref{lemma:units} and \myref{prop:space}(i), 
identifies its upper horizontal arrow, with $Y$ serving as a dummy
variable. More explicitly, using the projection formula
(\ref{four}), we see that the diagram can be rewritten as
\[\xymatrix{r^*r_!(r^*Y\wedge_M S^\nu) \ar[r]^\simeq 
\ar[d]_{r^*r_!(\theta\wedge \text{id})} & r^*Y \wedge_M r^*r_! S^\nu 
\ar[r]^-{\text{id}\wedge u_{S^0_M}} & r^*Y\wedge_M S^V_M \ar[d]^{\theta\wedge \text{id}}\\
r^*r_!(X\wedge_M S^\nu)\ar[rr]_{\Sigma^\infty_M u_K} && X\wedge_M S^V_M}\] 
Consider the dummy variable $Y$ levelwise.  We see from the case
$K=S^0_M$ of \myref{prop:space}(ii) that, at level $V$, 
the top row is the map $u_{r^*Y(V)}$. Therefore the diagram commutes levelwise 
by the naturality of $u$.

To prove that the composite (\ref{eq:trieq}) is the identity map, we apply 
$\Sigma^V$ to it. To abbreviate notation, write $Y = r_!(X\wedge_M C_M)$ and 
consider the following diagram.
\[\xymatrix{ 
\SI^V Y \ar[rrr]^-{\SI^V(\text{\text{id}}\sma\ta)} \ar@{=}[d]
\ar[drr]^-{\Sigma^V\tau_{Y}} & & & \SI^V(Y\sma r_!C_M) 
\ar[dl]^{\htp} \ar[d]^{\htp}\\
\SI^V Y \ar[d]_{\htp} & &\SI^V r_!(r^*Y\sma_M C_M)
\ar[ll]_-{\SI^Vr_!(\xi_X\sma_M\text{id})} \ar[dr]^{\htp}  
&  Y\wedge r_! S^\nu \ar[d]^\simeq\\
r_!(X\sma_M S^{\nu}) & & &
r_!(r^*Y\wedge_M S^\nu) 
\ar[lll]^-{r_!(\xi_X\wedge_M \text{id})} 
}\]
We must prove that the triangle at the upper left commutes. The 
arrows marked $\htp$ are given by (\ref{four}) and the evident
equivalence $\SI^V_MC_M\htp \SI^{\infty}S^{\nu}$.  The upper 
triangle commutes by the definition of $\ta_Y$ in terms of $\ta$, 
the bottom trapezoid commutes by naturality, and the triangle at
the right commutes by inspection of projection formula isomorphisms. 
Thus it suffices to prove that traversal of the perimeter gives a 
commutative diagram, and this will hold for $X$ if it holds for 
$\SI^V_MX$. In that case, we see from \myref{lemma:units}, 
\myref{prop:space}(ii), and a diagram chase that the 
perimeter agrees with the diagram that is obtained by applying the 
suspension spectrum functor to the diagram in 
\myref{prop:space}(iii). 
\end{proof}

The proof of \myref{prop:space} is based on the following construction 
of a natural map 
\[w_K\colon r^*r_!K \rtarr \tilde{K}\sma_M S^\tau\]
that will give rise to the required map $u_K$. Here $\tilde{K}$ is a
suitably ``fattened up'' version of $K$.

\begin{con} As in \cite[11.5]{MilS}, identify the tangent bundle of $M$
with the normal bundle of the diagonal embedding $M\rtarr M\times M$.
Let $U$ be a tubular neighborhood of the diagonal.  Let $\text{pr}_1$ and 
$\text{pr}_2$ be the projections $M\times M\rtarr M$ and let 
$\pi_i\colon U\rtarr M$ be their restrictions to $U$.  For an ex-space 
$(K,p_K,s_K)$ over $M$, consider the following diagram of 
retracts, where $\Delta$ is the diagonal and $\iota$ is the inclusion
of $U$ in $M\times M$.
Note that $\pi_i = \text{pr}_i \com \io$ and define 
$\tilde{K} = (\pi_1)_!\pi_2^* K$.
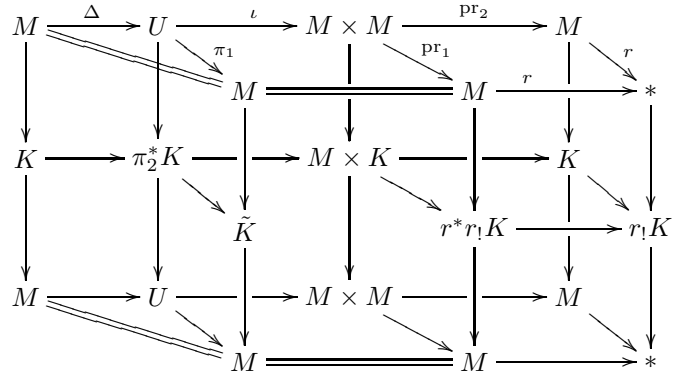
\[\xymatrix@=.4cm{M \ar[rr]^\Delta\ar[dd]\ar@{=}[drrr] 
&& U \ar[rr]^{\iota}\ar[dd]\ar[dr]^(.6){\pi_1} 
&& M\times M \ar[rr]^-{\text{pr}_2}\ar@{-}[d]\ar[dr]^(.6){\text{pr}_1} 
&& M\ar[dr]^(.6)r \ar@{-}[d]\\
&&& M\ar@{=}[rr]\ar[dd] &\ar[d]& M \ar[rr]^(.3){r}\ar[dd] &\ar[d]& {*} \ar[dd]\\
K \ar[rr]\ar[dd] && \pi_2^* K \ar@{-}[r]\ar[dd]\ar[dr] 
&\ar[r]& M\times K \ar@{-}[r]\ar[dd]\ar[dr] &\ar[r]& K\ar[dr]\ar@{-}[d]\\
&&& \tilde{K}\ar[dd] && r^*r_!K \ar[rr]\ar[dd] &\ar[d]& r_!K\ar[dd]\\ 
M \ar[rr]\ar@{=}[drrr] && U \ar@{-}[r]\ar[dr] &\ar[r]
& M\times M \ar@{-}[r]\ar[dr] &\ar[r]& M\ar[dr]\\
&&& M\ar@{=}[rr] && M \ar[rr] && {*}}\]
The floor and ceiling of the diagram are identical. 
The back wall is formed by pulling $K$ back along the maps of base spaces and 
then the front wall is obtained from the back wall by applying lower shriek 
functors. Here we have used the canonical isomorphism 
${\text{pr}_1}_!(M\times K)={\text{pr}_1}_!\text{pr}_2^*K\iso r^*r_!K$ 
associated to the pullback square on the right side of the floor. Since
the $\pi_i$ are homotopy equivalences, the maps 
\begin{equation}\label{kappak}
K\rtarr \pi_2^*K\rtarr \tilde K
\end{equation} 
at the left are equivalences when $K$ is $q$-cofibrant and $q$-fibrant, and
we denote the displayed composite equivalence as $\ka$.

To get a better feeling for the spaces in the diagram, we make the 
following schematic picture.
   
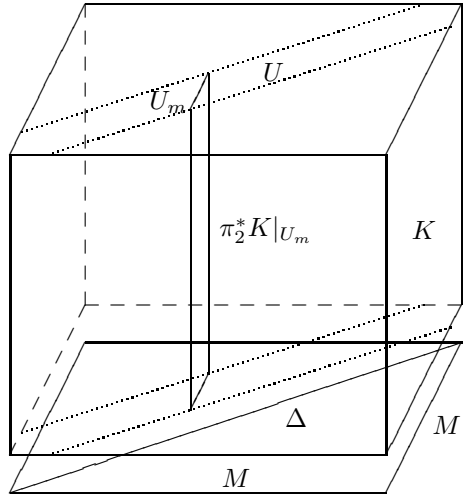
\[\begin{xy}
   0;(50,0) **@{-};(60,20) **@{-};(60,60) **@{-};(10,60) **@{-};(0,40)
   **@{-};(50,40) **@{-};(60,60) **@{-}
   ,(50,40);(50,0) **@{-},(0,40);(0,0) **@{-};(10,20) **@{--};(10,60)
   **@{--},(10,20);(60,20) **@{--}
   ,(5,0);(58.5,17) **@{.},(1.5,3);(55,20) **@{.}
   ,(5,40);(58.5,57) **@{.},(1.5,43);(55,60) **@{.} ,(24,6);(26.4,10.8)
   **@{-};(26.4,50.8) **@{-};(24,46) **@{-};(24,6) **@{-}
   ,(55,30) *{K},(34,30) *{\pi^*_2K|_{U_m}}
   ,(35,51) *{U},(21,47) *{U_m}
    ,(0,-5);(50,-5)**@{-};(60,15)**@{-};(10,15)**@{-};(0,-5)**@{-};(60,15)
   **@{-}
   ,(38,5) *{\Delta},(30,-3) *{M},(58,4) *{M}
   \end{xy}\]
   The cube represents $M\times K$ with $M$ in the horizontal direction
   and $K$ in the other two directions. It is sitting over $M\times M$ 
   with vertical projection and we think of the top face as being the
   image of the section. We can view $\pi_2^* K$ as the subspace of
   $M\times K$ sitting over the neighborhood $U$ of the diagonal in     
   $M\times M$. Passing to the front face of the diagram, the fibers of
   $\tilde{K}$ over points $m$ in $M$ are the slices $\pi_2^*K|_{U_m}$ of
   $\pi_2^* K$ as displayed with basepoint obtained by identifying the
   fiber $U_m=\pi_1^{-1}(m)$ over $m$ in the bundle $U$ to a single
   point. We therefore think of $\tilde{K}$ as a fattening of the space
   $K$ that sits over the diagonal in $M\times M$.

For a sub $G$-space $A$ of a $G$-space $K$ over $M$, not necessarily sectioned,
passage to fiberwise quotients gives an ex-$G$-space $K/\!_M A$ over $M$ with 
total space $K\cup_A M$. Given two such pairs $(K,A)$ and $(L,B)$, we obtain
a product pair by setting
\[(K,A)\times_M (L,B)=(K\times_M L, K\times_M B \cup A\times_M L).\]
Its fiberwise quotient is the ex-$G$-space $(K/\!_M A)\wedge_M (L/\!_M B)$ over $M$. 

The pair $(M\times M, M\times M-U)$ is a model for the Thom
complex $T\ta$, and we can identify $T\ta$ with the 
quotient space $(M\times M)/(M\times M-U)$. More relevantly
for us, the fiberwise quotient $(M\times M)/\!_M(M\times M-U)$
is a model for $S^{\ta}$. View $M\times K=\text{pr}_2^* K$, 
$M\times M$, and $U$ as $G$-spaces over $M$ via projection to the 
first factor and embed $U$ in $M\times K$ by sending $(m,n)$ to 
$(m,s_K(n))$. We have the diagonal map
\[M\times K\rtarr   (M\times K, U) \times_M (M\times M, M\times M - U)\]
of $G$-spaces over $M$ that sends $(m,k)$ to $((m,k), (m,p_K(k)))$ for $m\in M$
and $k\in K$. It induces the top map in the following diagram in $G\sK_M$.
\[\xymatrix{ M\times K \ar[r]\ar[d]\ar@{->}[dr] 
& ((M\times K)/\!_M U)\wedge_M {S^\tau}\\
r^*r_!K \ar@{->}[r]_{w_K} & \tilde{K}\wedge_M {S^\tau} \ar[u]}\]
Here $(M\times K)/\!_M U$ is obtained from $M\times K$ by identifying all points
of the form $(m,s_K(n))$ such that $(m,n)\in U_m$ to a single basepoint in the 
fiber over $m$.  It therefore contains $\tilde K$ as a subspace, and this gives 
the right vertical inclusion.  The image of the top arrow lands in the 
image of the right vertical arrow since if $(m,p_K(k))$ is not in $U_m$, 
then $(m,k)$ maps to the basepoint in $S^{\ta}$ and therefore to the
basepoint in $(M\times K)/\!_M U\sma_M S^{\ta}$.  This gives the diagonal arrow. 
Note that $r^*r_!K = M\times (K/s_K(M))$ and the left vertical arrow is the
obvious quotient map. Since the diagonal arrow maps $(m,x)$ with $x\in s_K(M)$ 
to the base point of the fiber over $m$ in $\tilde K\wedge_M {S^\tau}$, it is 
constant on the fibers of the left vertical arrow. It therefore factors through
a map $w_K$. Explicitly, $w_K$ is specified by 
\begin{equation}\label{doubleu}
w_K(m,[x])=\begin{cases}
\; [m,x]\wedge[m,p_K(x)] &\text{if $(m,p_K(x))\in U_m$,}\\
\; * & \text{otherwise,}
\end{cases}
\end{equation}
where $m\in M$, $x\in K$, and the square brackets denote equivalence classes.
\end{con}

\begin{proof}[Proof of \myref{prop:space}]
Here we are working in homotopy categories, and we may assume that $K$
is $qf$-fibrant and $qf$-cofibrant.  Let $L=K\wedge_M S^\nu$. 
We define the map $u_K$ in $\Ho G\sK_M$ by the natural zig-zag 
\[\xymatrix{
r^*r_! L \ar[r]^-{w_{L}}
& \tilde{L}\wedge_M {S^\tau} 
&&  L \wedge_M S^\tau \ar[ll]_-{\ka\sma_M\text{id}} \ar[r]^-{\mu} 
& \Sigma^V_M K}\]
of arrows in $G\sK_M$, where $\mu$ is induced by the isomorphism
\[\xymatrix{
S^\nu\wedge_MS^\tau \ar[r]^-{\text{twist}}
& S^\tau\wedge_MS^\nu \ar[r]^-{\cong}
& S^V_M}\]
and $\ka\colon L\rtarr \tilde L$ is an equivalence as in ({kappak}).

\noindent
\emph{Proof of (i)}
We must show that $u_{S^0_M}\simeq k$ in $\Ho G\sK_M$. Using our zig-zag definition of $u_K$, we see that it suffices to show that the composite 
\[\xymatrix{
M\times T\nu \ar[r]^-{k} & M\times S^V 
\ar[r]^-{\mu^{-1}} & S^\nu\wedge_M S^\tau 
\ar[r]^-{\kappa\wedge\text{id}} & \widetilde{S^\nu}\wedge_M S^\tau}\]
is homotopic to the map $w_{S^\nu}$ defined in (\ref{doubleu}). As noted in the proof of \myref{lemma:units}, $k(m,[v])=(m,\text{ev}([v]\wedge m))$, where $v\in S^\nu$ and brackets denote equivalence classes. Recall that the map $\text{ev}$ depends on a choice of a tubular neighborhood of $M$ in $\nu\times M$ (as in the proof of \myref{starting}).  We use the obvious choice
\[\{(v,m)\, |\, v\in \nu,\ \, m\in M, \ \, \text{and}\ \  (p_\nu(x),m)\in U\}.\]
Under our identification of the normal bundle of $\Delta:M\to M\times M$ and thus of its tubular neighborhood $U$ with $\ta$, this tubular neighborhood is identified with $M\times V\cong \nu\oplus \tau$. When not in the section, we can view points $[m,n]\in S^\tau=(M\times M)/_M(M\times M-U)$ as vectors $(m,n)$ in the tangent space $U_m\cong\tau_m\subset V_m\cong V$ of $M$ at $m$. We then have that $(m,\text{ev}([v]\wedge m)) = (m,(p_\nu(v),m)+v)$. To identify this point in the image of $\mu$, let $u\in\nu_m$ be such that $(p_\nu(v),m)+v = (m,p_\nu(v))+u$ in $V$ and note that $u$ depends continuously on $m$ and $v$. Since $\mu(u\wedge [m,p_\nu(v)])=(m,(m,p_\nu(v))+ u)$, the composite displayed above is given by
\[(\kappa\wedge_M \text{id})\mu^{-1}k(m,[v])=
\begin{cases}
\; [m,u]\wedge [m,p_\nu(v)] & \text{if $(m,p_\nu(v))\in U_m$,}\\
\; * & \text{otherwise.}
\end{cases}\]
A linear homotopy in the fibers of $\nu$ shows that this map is 
homotopic to $w_{S^\nu}$.

\noindent
\emph{Proof of (ii)}
Inspection of the construction of $w$ gives the following naturality diagram for based $G$-spaces $T$ and ex-$G$-spaces $K$ over $M$.
\[\xymatrix{ r^*r_!(T_M\wedge_M K) \ar[r]^\simeq\ar[d]_{w} 
& T_M \wedge_M r^*r_!K \ar[d]^{\text{id}\wedge w}\\
(\widetilde{T_M\wedge_M K}) \wedge_M {S^\tau}\ar[r]^\simeq 
& T_M \wedge_M \tilde{K}\wedge_M {S^\tau} }.\]
Here, using $r^*r_!\iso (\pi_1)_!\pi_2^*$, the bottom equivalence is the following application of the projection formula.
\begin{align*}
(\pi_1)_!\pi_2^*(r^*T\wedge_M K)
&\simeq  (\pi_1)_!(\pi_2^*r^*T\wedge_M \pi_2^*K)\\
&\simeq  (\pi_1)_!(\pi_1^*r^*T\wedge_M \pi_2^*K)\\
&\simeq  r^*T\wedge_M (\pi_1)_!\pi_2^*K
\end{align*}
This use of the projection formula is compatible with its use for $r_!$ to obtain the equivalence of the top row. Analogous naturality diagrams for the other two maps in the definition of $u_K$ give the conclusion.

\noindent
\emph{Proof of (iii)}
Again let $L=K\wedge_M S^\nu$. Expanding the diagram in the statement of (iii) in terms of the definition of $u_K$, we must prove that the following diagram commutes in $\Ho G\sK_*$, where the equivalences here are the vertical arrows of the diagram in (iii).
\[\xymatrix@=.75cm{
r_!(\Sigma^V_MK\wedge_MS^{\nu}) 
&& r_!L \wedge S^V \ar[rr]^-{\text{id}\wedge t} \ar[ll]_-{\simeq}
&&  r_!L \wedge r_! S^\nu \ar[d]^\simeq\\
r_!(L\wedge_M {S^\tau} \wedge_M S^\nu) \ar[u]^{r_!\mu} 
\ar[rr]_-{r_!(\kappa\wedge\text{id})}
&& r_!(\tilde L\wedge_M {S^\tau} \wedge_M S^\nu)
&& r_!(r^*r_!L\wedge_M S^\nu) \ar[ll]^-{r_!(w_L\wedge \text{id})}}\]
We chase the diagram starting in $r_!(L\wedge_M {S^\tau} \wedge_M S^\nu)$ and mapping to $r_!(\tilde L\wedge_M {S^\tau} \wedge_M S^\nu)$. Let $x = k\sma w \in L = K\sma_M S^{\nu}$, $u\in S^\tau$, and $v\in S^\nu$ be points in fibers over a given $m\in M$. Using square brackets to denote passage to quotient spaces (the lower shriek functors), we see that $r_!(\kappa\sma\text{id})$  sends $[x\sma u\sma v]$ to $[[m,x]\sma u\sma v]$. The definitions of $\mu$ and of the top left equivalence (which is the left vertical equivalence in the diagram of the statement) are arranged in such a way that the composite of $\mu$ and the inverse of the equivalence sends $[x\sma u\sma v]$ to $[x\sma(u+v)]$. Let $t(u+v)=[z]$, $z\in S^\nu$, and let $n=p_\nu(z)$. Chasing $[x\sma u\sma v]$ around the top of the diagram, when we do not arrive at the basepoint we arrive at the point $[[n,x]\sma [n,m]\sma z]$, where $[n,m]$ is an element of $U\cong \tau$. We can identify the target space with $r_!(\tilde{L}) \sma S^V$ using the identification of $S^\tau \sma S^\nu$ with $M \times S^V$ and the projection formula. Then our two maps are homotopic by a homotopy $h$ that can be written in the form
\[h([x\wedge u\wedge v],s)=[m+s[m,n],x]\wedge(u+v+2s[n,m]).\]
Here $m+s[m,n]$ denotes a point on the path from $m$ to $n$ in $M$ that is the image under the exponential map of the line segment from $0$ to $[m,n]$ in the tangent space at $m$.   The Thom map takes $u+v$ in $S^V$ (which in $M \times S^V$ is based at $m$) to the point $z$ in $S^\nu$ based at $n$. In $S^V$, we have $u+v = [m,n]+z$.  Since $[n,m]= -[m,n]$ under the identification of $\tau$ with $U$ (as in \cite[11.5]{MilS}), we see that the homotopy ends at the composite around the top of the diagram, and it clearly begins at $r_!(\kappa\sma\text{id})$. 
\end{proof}

\backmatter


\end{document}